# Exact quasiclassical asymptotics beyond Maslov canonical operator.

J.Foukzon[1], A.A.Potapov[2], S.A. Podosenov[3],.
[1]Israel Institute of Technology.
[2]IRE RAS,
[3]All-Russian Scientific-Research Institute
 for Optical and Physical Measurements.

**Abstract**: The main purpose of this paper is to calculate exact quasiclassical asimptotic of the quantum averages without any reference to the corresponding quasiclassical asimptotic of the Schrödinger wave function $\Psi(x,t,x_0;\hbar)$ given via Maslov canonical operator $K_{\Lambda^n(t)}^{1/\hbar}[\bullet]$. We suggest a new asymptotic representation for the quantum averages with position variable with localized initial data.

# Content

Colombeau-Wiener measure $\widetilde{D}^{\mathbb{C}}[x(\tau), v_\varepsilon]$.

**I.0.7**.Colombeau-Wiener Path Integral and generalized Haba theorem.

**I.1.Colombeau generalized functions**.

**I.1.1**.Point Values of Generalized Functions on $M$ and Generalized Numbers.

**I.1.2**.Colombeau algebra of compactly supported generalized functions and Colombeau algebra of tempered generalized functions.

**I.1.3**.Algebra of generalized functions based on $\mathcal{L}(\mathbb{R}^n)$. Colombeau algebra $G_{\tau,\mathcal{L}}(\mathbb{R}^n)$.

**I.1.4**.Regular Colombeau algebras.

**I.1.5**.Integration in Colombeau theory.

**I.1.6**.Fourier transform in Colombeau theory.

**I.1.7**.Topological $\widetilde{\mathbb{C}}$-modules.

**I.1.8**.Completeness and duality theory for topological $\widetilde{\mathbb{C}}$-modules.

**I.1.9**.Definition and basic properties of $G_E$ and $\widetilde{G}_B$.

**I.1.10**.Definition and basic properties of $G_B$ and $\widetilde{G}_B$.

**I.1.11**.Colombeau Sobolev spaces.Definition and basic properties of
$\left(W_p^{l,k,\varepsilon}(\Omega)\right)_{\varepsilon\in(0,1]}$.

**I.2**.Pseudodifferential operators with coefficients are distributions.

**I.3**.The generalized oscillatory integral.

**I.4**.Pseudodifferential operators with generalized symbol.

**I.5**.Formal series and symbols.

**I.6**.Oscillatory integrals with generalized phase function on $\mathbb{R}^n$.

**I.7**.Generalized symbols and amplitudes of weighted type.

**I.8**.Pseudodifferential operators acting on $G_{\tau,\mathcal{L}}(\mathbb{R}^n)$.

**I.9**.Oscillatory integrals of weighted type.

**II.1.1**.Feynman-Colombeau path integral.

**II.1.2**.Generalized Maslov theorem.

**II.2**.Master Equation for the Feynman-Colombeau path integral.

**II.2.1**.Master Equation for the Feynman-Colombeau path integral corresponding with a Schrödinger equation with a non-Hermitian Hamiltonian.

**II.2.2**.Master Equation for the Feynman-Colombeau path integral corresponding with a Schrödinger equation with Hermitian Hamiltonian.

**II.3**. Quantum anharmonic oscillator with a cubic potential.

**II.4**. Quantum anharmonic oscillator with a duble well potential.

**III**.Comparizon with a perturbation theory.

**III.1**.Path-integrals calculation by using saddle-point approximation.

**IV**. Wave function collapse and weak Colombeau solutions of the Colombeau-Schrödinger equation.

**IV.1**. Weak Colombeau solutions of the Colombeau-Schrödinger equation.

**V**.Schrödinger's cat paradox resolution using GRW collapse model.
**Appendix**. CAT-A.
**Appendix**. CAT-B.
**Appendix**. CAT-C.
**Apendix I.0**.Generalized Colombeau Propagator.
**Apendix I.1**.Derivation of the Feynman-Colombeau path integral for the non-Hermitian generalized Hamiltonian.
**Apendix I.2**. Colombeau-Schrödinger Equation from Feynman-Colombeau path integral.
**Apendix I.3**.Colombeau-Feynman oscillatory path integral.

**Apendix II.1.1**.Upper and lower infinite-dimensional Feynman subintegrals.
**Apendix II.1.2**.Upper and lower infinite-dimensional Colombeau-Feynman subintegrals.
**Apendix II.2.1**.Colombeau Signed Measures.
**Apendix II.2.2**. Colombeau Submeasures.
**Apendix II.2.3**.Infinite-Dimensional Colombeau-Feynman Signed Measures.
**Apendix III**. Colombeau-Albeverio normed path integral.
**Apendix IV**.The Schrödinger equation with singular potential $(V_\varepsilon)_\varepsilon$ in Colombeau module $E_M[\mathcal{F}(\mathbb{R}^d)]$.

# I.0.Introduction

## I.0.1.Colombeau solutions of the Schrödinger equation using Maslov canonical operator.

Let **H** be a complex infinite dimensional separable Hilbert space,with inner product $\langle \cdot, \cdot \rangle$ and norm $|\cdot|$. Let us consider Schrödinger equation:

$$-i\hbar \frac{\partial \Psi(t)}{\partial t} + \widehat{H}(t)\Psi(t) = 0,$$

$$\textbf{(I.0.1.1)}$$

$$\Psi(0) = \Psi_0(x), x \in \mathbb{R}^n,$$

where operator $H(t)$ given via formula

$$H(t) = -\frac{\hbar^2}{2m}\Delta + V(x,t) \qquad \textbf{(I.0.1.2)}$$

and $H$ is essentially self-adjoint, $\widehat{H}$ is the closure of $H, \Psi(t) : \mathbb{R} \to \mathbf{H}$ and $\widehat{H}(t) : \mathbb{R} \times \mathbf{H} \to \mathbf{H}$.

Maslov local quasiclassical asymptotic of the Green function of the Schrödinger equation (I.0.1.1) wih potential $V(x,t)$ is [26] :

$$G(t,x,y;\hbar) = (2\pi i\hbar)^{-n/2}\sqrt{\left|\frac{\partial^2 S}{\partial x \partial y}\right|}\exp\left[\frac{i}{\hbar}S(t,x,y)\right] + \hbar\zeta(t,x,y;\hbar),$$

$$\textbf{(I.0.1.3)}$$

$$S(t,x,y) = \int\limits_0^t \left[\frac{m}{2}\left(\frac{\partial X(\tau;t,x,y)}{\partial \tau}\right)^2 - V(X(\tau;t,x,y),\tau)\right]d\tau.$$

Here $t < T[V(x,t)], x,y \in \mathbb{R}^n, \eta(t,y) = \|\zeta(t,x,y;\hbar)\|_{\mathcal{L}_2(\mathbb{R}^n)} < C$ and

$$m\ddot{X} = -\frac{\partial V(X,t)}{\partial X},$$

$$\textbf{(I.0.1.4)}$$

$$X(0;t,x,y) = y, X(t;t,x,y) = x.$$

By using convolution one obtain [25]-[27]:

$$G(t_1 + t_2,x,y;\hbar) = \int\limits_{\mathbb{R}^n} G(t_1,\widetilde{y},y;\hbar)G(t_2,x,\widetilde{y};\hbar)d\widetilde{y}. \qquad \textbf{(I.0.1.5)}$$

We assume now that [25]-[27]:

(**i**) $V(x) \in C^\infty(\mathbb{R}^n)$ and

$$\sup_{x \in \mathbb{R}^n} |\mathbf{grad}_x V(x)| < \infty,$$

(0.1.6)

$$\inf_{x \in \mathbb{R}^n} V(x) \geq -\infty,$$

**Remark 0.1.1**. We remind that ander conditions (i) $H$ is essentially self-adjoint.

(**ii**) for a fixed $t$ and $y$ boundary problem is a system, (0.1.4) has a finite number of solutions $X_k(\tau; t, x, y), k = 1, \ldots, N,$

(**iii**) for a given $t, x, y$ the set $\{X_k(t; t, x, y) | k = 1, \ldots, N\}$ does not contain focal points. Then Maslov global quasiclassical asymptotic of the Green function obtained from Eq.(0.1.6) is [27]:

$$G(t, x, y; \hbar) =$$

$$(2\pi i\hbar)^{-n/2} \left[ \sum_{k=1}^{N} \sqrt{\left| \frac{\partial^2 S_k(t, x, y)}{\partial x \partial y} \right|} \exp\left[ \frac{i}{\hbar} S_k(t, x, y) - \frac{i\pi\gamma_k}{2} \right] + O(\hbar) \right],$$

(0.1.7)

$$S_k(t, x, y) = \int_0^t \left[ \frac{m}{2} \left( \frac{\partial X_k(\tau, t, x, y)}{\partial \tau} \right)^2 - V(X_k(\tau, t, x, y)) \right] d\tau.$$

Here the number $\gamma_k$ is the Morse index of the trajectory $X_k(\tau; t, x, y)$.

Unfortunately for real physical problems assuption (**i**) does not satisfies and consequently Maslov global quasiclassical asymptotic of the Green function (0.1.7) does not valid. We improve this difficulties using Colombeau approach [1]-[12],[42], [43].

Let us consider Colombeau-Schrödinger equation [42],[43]:

$$-i((\hbar_\varepsilon)_\varepsilon)\left(\frac{\partial\Psi_\varepsilon(t)}{\partial t}\right)_\varepsilon + \left(\widehat{H}_\varepsilon(t)\Psi_\varepsilon(t)\right)_\varepsilon = 0,$$

$$(0.1.8)$$

$$(\Psi_\varepsilon(0))_\varepsilon = (\Psi_{0,\varepsilon}(x))_\varepsilon, x \in \mathbb{R}^n,$$

where $\forall\varepsilon \in (0,1]$ operator $H_\varepsilon(t)$ given via formula

$$H_\varepsilon(t) = -\frac{\hbar_\varepsilon^2}{2m}\Delta + V_\varepsilon(x),$$

$$(V_\varepsilon(x))_\varepsilon \in G(\mathbb{R}^n),$$

$$(0.1.9)$$

$$\inf_{x\in\mathbb{R}^n} (V_\varepsilon(x))_\varepsilon \geq -\infty$$

and $V_0(x) = V(x)$.

**Remark 0.1.2**. We note that $(V_\varepsilon(x))_\varepsilon \in G(\mathbb{R}^n) \Rightarrow$

$$\forall\varepsilon(\varepsilon \in (0,1])\left[\sup_{x\in\mathbb{R}^n} |\mathbf{grad}_x V_\varepsilon(x)| < \infty\right].$$

$$(0.1.10)$$

Maslov local quasiclassical asymptotic of the Green function of the Colombeau-Schrödinger equation (0.1.8) wih potential $(V_\varepsilon(x))_\varepsilon \in G(\mathbb{R}^n)$ is :

$$(G_\varepsilon(t_\varepsilon, x, y; \hbar_\varepsilon))_\varepsilon =$$

$$(2\pi i((\hbar_\varepsilon)_\varepsilon))^{-n/2} \sqrt{\left| \left( \frac{\partial^2 S_\varepsilon(t_\varepsilon, x, y)}{\partial x \partial y} \right)_\varepsilon \right|} \exp\left[ \frac{i}{(\hbar_\varepsilon)_\varepsilon}(S_\varepsilon(t_\varepsilon, x, y))_\varepsilon \right] +$$

$$(0.1.11)$$

$$((\hbar_\varepsilon)_\varepsilon)((\zeta_\varepsilon(t_\varepsilon, x, y; \hbar))_\varepsilon),$$

$$(S_\varepsilon(t_\varepsilon, x, y))_\varepsilon = \left( \int_0^{t_\varepsilon} \left[ \frac{m}{2} \left( \frac{\partial X_\varepsilon(\tau; t_\varepsilon, x, y)}{\partial \tau} \right)^2 - V_\varepsilon(X_\varepsilon(\tau; t_\varepsilon, x, y), \tau) \right] d\tau \right)_\varepsilon.$$

Here $(t_\varepsilon)_\varepsilon < (T[V_\varepsilon(x, t_\varepsilon)])_\varepsilon, x, y \in \mathbb{R}^n, \left( \|\zeta(t_\varepsilon, x, y; \hbar)\|_{\mathcal{L}_2(\mathbb{R}^n)} \right)_\varepsilon < (C_\varepsilon)_\varepsilon \in \widetilde{\widetilde{\mathbb{R}}}$ and

$$m\left( (\ddot{X}_\varepsilon)_\varepsilon \right) = -\left( \frac{\partial V_\varepsilon(X_\varepsilon, t)}{\partial X_\varepsilon} \right)_\varepsilon,$$

$$(0.1.12)$$

$$(X_\varepsilon(0, t, x, y))_\varepsilon = y, (X_\varepsilon(t; t, x, y))_\varepsilon = x.$$

By using convolution we obtain

$$(G_\varepsilon(t_1 + t_2, x, y; \hbar_\varepsilon))_\varepsilon = \left( \int_{\mathbb{R}^n} G_\varepsilon(t_1, \widetilde{y}, y; \hbar) G_\varepsilon(t_2, x, \widetilde{y}; \hbar) d\widetilde{y} \right)_\varepsilon. \qquad (0.1.13)$$

We assume now that :
**(i)** for a fixed $\varepsilon \in (0, 1]$, $t$ and $y$ boundary problem is a system, (0.1.12) has a finite number of solutions $X_{k,\varepsilon}(\tau; t, x, y), k = 1, \ldots, N(\varepsilon)$,
**(ii)** for a given $t, x, y : \forall \varepsilon \in (0, 1]$ the set $\{X_{k,\varepsilon}(t; t, x, y) | k = 1, \ldots, N(\varepsilon)\}$ does not contain focal points.
Then Maslov global quasiclassical asymptotic of the Green function obtained from Eq.(0.1.13) is:

$$(G_\varepsilon(t,x,y;\hbar_\varepsilon))_\varepsilon =$$

$$(2\pi i(\hbar_\varepsilon)_\varepsilon)^{-n/2} \times$$

$$\left[ \left( \sum_{k=1}^{N(\varepsilon)} \sqrt{\left| \left( \frac{\partial^2 S_{k,\varepsilon}(t,x,y)}{\partial x \partial y} \right) \right|} \right)_\varepsilon \times \right. \tag{0.1.14}$$

$$\left. \exp\left[ \frac{i}{(\hbar_\varepsilon)_\varepsilon} (S_{k,\varepsilon}(t,x,y))_\varepsilon - \frac{i\pi(\gamma_{k,\varepsilon})_\varepsilon}{2} \right] \right] + (O(\hbar_\varepsilon))_\varepsilon,$$

$$(S_{k,\varepsilon}(t,x,y))_\varepsilon = \left( \int_0^t \left[ \frac{m}{2} \left( \frac{\partial X_k(\tau,t,x,y)}{\partial \tau} \right)^2 - V_\varepsilon(X_{k,\varepsilon}(\tau,t,x,y)) \right] d\tau \right)_\varepsilon .$$

Here the number $\gamma_{k,\varepsilon}$ is the Morse index of the trajectory $X_{k,\varepsilon}(\tau;t,x,y)$.

# I.0.2. The Main result .

The main purpose of this paper is:

**(I)** To calculate exact quasiclassical asimptotic of the complex $\widetilde{\mathbb{C}}$-valued generalized quantum averages

$$((\langle s,t,\mathbf{x}_0;\hbar_\varepsilon\rangle)_\varepsilon = \frac{\left( \int x_s \Psi_\varepsilon(\mathbf{x},t,\mathbf{x}_0;\hbar_\varepsilon) d^n x \right)_\varepsilon}{\left( \int \Psi_\varepsilon(\mathbf{x},t,\mathbf{x}_0;\hbar_\varepsilon) d^n x \right)_\varepsilon},$$

$$(\Psi_\varepsilon(\mathbf{x},0,\mathbf{x}_0;\hbar))_\varepsilon = \delta(\mathbf{x} - \mathbf{x}_0), \tag{0.2.1}$$

$$s = 1,\ldots,n$$

without any reference to the corresponding exact quasiclassical asimptotic of the Schrödinger generalized wave function $(\Psi_\varepsilon(\mathbf{x}, t, \mathbf{x}_0; \hbar))_\varepsilon$ given from Colombeau-Schrödinger Eq.(0.1.8)-(0.1.9)

(**II**) To calculate exact quasiclassical asimptotic of the real $\widetilde{\mathbb{R}}$-valued generalized quantum averages

$$\left(\left\langle s, t, \mathbf{x}_0; \hbar_\varepsilon \right\rangle_\varepsilon\right)_\varepsilon = \frac{\left(\displaystyle\int_{\mathbb{R}^n} x_s |\Psi_\varepsilon(\mathbf{x}, t, \mathbf{x}_0; \hbar_\varepsilon)|^2 d^n x\right)_\varepsilon}{\left(\displaystyle\int_{\mathbb{R}^n} |\Psi_\varepsilon(\mathbf{x}, t, \mathbf{x}_0; \hbar_\varepsilon)|^2 d^n x\right)_\varepsilon}, \tag{0.2.2}$$

$$s = 1, \ldots, n$$

without any reference to the corresponding exact quasiclassical asimptotic of the Schrödinger wave function $(\Psi_\varepsilon(x, t, x_0; \hbar))_\varepsilon$ given from Colombeau-Schrödinger Eq.(0.1.8)-(0.1.9).

**Remark**.Note that such quasiclassical asimptotic very important from point of view of the "quantum jumps" problem, well known in modern quantum mechanics.The existence of such jumps was required by Bohr in his theory of the atom. He assumed that an atom remained in an atomic eigenstate until it made an instantaneous jump to another state with the emission or absorption of a photon. Since these jumps do not appear to occur in solutions of the Schrodinger equation, something similar to Bohr's idea has been added as an extra postulate in modern quantum mechanics. The question arises whether an explanation of these jumps can be found to result from a solution $\Psi(\mathbf{x}, t, \mathbf{x}_0; \hbar)$ of the Schrödinger equation (0.1.1)-(0.1.2) alone without additional postulates.

We suggest a new asymptotic representation for the quantum averages (0.2.2) with position variable with well localized initial data,i.e.

$$\left(\left\langle i, 0, \mathbf{x}_0; \hbar_\varepsilon \right\rangle_\varepsilon\right)_\varepsilon \simeq x_{0i}, \tag{0.2.3}$$

$$i = 1, \ldots, n.$$

The canonical physical interpretation of these asymptotics shows that the answer is "yes."

As is well-known,when the potential $V$ is sufficiently regular any solution of the Schrödinger equation (0.1.1) can be represented by Feynman path integral of the form:

$$\Psi(x,t,x_0;\hbar) = \int_{\substack{q(0)=y \\ q(t)=x}} \Psi_0(q(0)) \exp\left\{ \frac{i}{\hbar} S[\dot{q}(\tau),q(\tau),t] \right\} [Dq(\tau)] dy,$$

$$S[\dot{q}(\tau),q(\tau),t] = \int_0^t \mathcal{L}(\dot{q}(\tau),q(\tau),\tau) d\tau,$$

$$\mathcal{L}(\dot{q}(\tau),q(\tau),\tau) = m\frac{\dot{q}^2(\tau)}{2} - V(q(\tau),t).$$

(0.2.4)

Maslov quasiclassical asymptotic of the wave function $\Psi(x,t;\hbar)$ is [22],[25]:

$$\Psi(x,t;\hbar) = (2\pi i\hbar)^{-\frac{n}{2}} \sum_{k=1}^m \int dy \Psi_0(y) |D(q^k(\tau;y,x,t))|^{\frac{1}{2}} \times$$

$$\times \exp\left\{ \frac{i}{\hbar}\left[ S[q^k(\tau;y,x,t)] - \frac{i\pi\gamma_k}{2} \right] \right\} + O(\hbar).$$

(0.2.5)

Here

$$S[q^k(\tau;y,x,t),t] = \int\limits_0^t \mathcal{L}(\dot{q}^k(\tau;y,x,t),q^k(\tau;y,x,t),\tau)d\tau,$$

$$\mathcal{L}(\dot{q}(\tau),q(\tau),\tau) = m\frac{\dot{q}^2}{2} - V(q(\tau),\tau),$$

$$\delta S[q^k(\tau;y,x,t),t] = 0, q^k(0;y,x,t) = y, q^k(t;y,x,t) = x,$$

$$D(q^k(\tau;y,x,t)) = \det\left\| \frac{\partial S[q^k(\tau;y,x,t),t]}{\partial x \partial y} \right\| \neq 0, x \neq y.$$

(0.2.6)

and the number $\gamma_k = \gamma_k[q^k(\tau;y,x,t)]$ is the Morse index of the trajectory $q^k(\tau;y,x,t)$. Assume that $\Psi_0(x) = \frac{1}{\varepsilon}\exp\left(-\frac{(x-x_0)^2}{\varepsilon^2}\right), \varepsilon \ll \hbar \ll 1$, i.e. $\Psi_0(x) \simeq \delta(x-x_0)$. From Eq.(0.2.5) we obtain

$$\Psi(x,t,x_0;\hbar) = (2\pi i\hbar)^{-\frac{n}{2}} \sum_{k=1}^m \int dy \Psi_0(y-x_0) |D(q^k(\tau;y,x,t))|^{\frac{1}{2}} \times$$

$$\times \exp\left\{ \frac{i}{\hbar}\left[ S[q^k(\tau;y,x,t)] - \frac{i\pi\gamma_k}{2} \right] \right\} + O(\hbar) =$$

(0.2.7)

$$\sum_{k=1}^m |D_1(q^k(\tau;x_0,x,t))|^{\frac{1}{2}} \exp\left\{ \frac{i}{\hbar}\left[ S[q^k(\tau;x_0,x,t),t] - \frac{i\pi\gamma_k}{2} \right] \right\} + O(\hbar).$$

Let us consider quantum harmonic oscillator. Corresponding Schrödinger equation is:

$$-ih\frac{\partial \Psi(x,t)}{\partial t} = \frac{\hbar^2}{2}\Delta\Psi(x,t) - \frac{\omega^2}{2}x^2\Psi(x,t),$$

$$\Psi(0) \simeq \delta(x-x_0).$$

(0.2.8)

By canonical calculation one obtain

$$S[q^k(\tau;y,x,t),t] = \frac{\omega}{2\sin\omega t}[(\cos\omega t)(x^2+y^2)-2xy],$$

$$D(q^k(\tau;y,x,t)) = \det\left\|\frac{\partial^2 S}{\partial x\partial y}\right\| = \frac{\omega^n}{(\sin\omega t)^n}.$$

(0.2.9)

Substitution Eqs.(0.2.9) into Eq.(0.2.5) gives [22]:

$$\Psi(x,t,x_0;h) = (2\pi ih)^{-\frac{n}{2}}\int dy\,\delta(x-x_0)|D(q^k(\tau;y,x,t))|^{\frac{1}{2}}\times$$

$$\times\exp\left\{\frac{i}{\hbar}\left[S[q^k(\tau;y,x,t)]-\frac{i\pi\gamma_k}{2}\right]\right\}+O(h) =$$

$$\omega^{n/2}(2\pi ih\sin\omega t)^{-\frac{n}{2}}\times$$

$$\exp\left\{\frac{i}{\hbar}\left[\frac{\omega}{2\sin\omega t}[(\cos\omega t)(x^2+x_0^2)-2xx_0]-\frac{i\pi\gamma_k}{2}\right]\right\}+O(h).$$

(0.2.10)

From Eq.(0.2.10) one obtain

$$\Psi(x, t, x_0; \hbar) =$$

$$\omega^{n/2}(2\pi i\hbar \sin \omega t)^{-\frac{n}{2}} \times$$

$$\exp\left\{\frac{i}{\hbar}\left[\frac{\omega}{2\sin\omega t}[(\cos\omega t)(x^2 + x_0^2) - 2xx_0] - \frac{i\pi\gamma_k}{2}\right]\right\}$$

(0.2.11)

and therefore

$$\int x_s \Psi(x, t, x_0; \hbar)d^{\,n}x =$$

$$\omega^{n/2}(2\pi i\hbar \sin \omega t)^{-\frac{n}{2}} \times$$

$$\int d^{\,n}x x_s \exp\left\{\frac{i}{\hbar}\left[\frac{\omega}{2\sin\omega t}[(\cos\omega t)(\mathbf{x}^2 + \mathbf{x}_0^2) - 2\mathbf{x}\mathbf{x}_0] - \frac{i\pi\gamma_k}{2}\right]\right\} =$$

(0.2.12)

$$\frac{x_{0,s}}{\cos\omega t}\exp\left\{\frac{i}{\hbar}\left[\omega x_0^2(\tan\omega t) - \frac{i\pi\gamma_k}{2}\right]\right\}.$$

Substitution Eq.(0.2.11) and Eq.(0.2.12) into Eq.(0.2.1) gives

$$\langle x_i, t, x_0; \hbar \rangle = \frac{x_{s0}}{\cos\omega t},$$

$$s = 1, \dots, n.$$

(0.2.13)

Formula (0.2.5) can be written by using Maslov canonical operator $K^{1/\hbar}_{\Lambda^n(t)}[\cdot]$.

Let $\Lambda^n(t)$ be an $n$-dimensional Lagrangian manifold of class $\mathbf{C}^\infty$ in the phase space $\mathbb{R}^{2n}_{x,p}$ where $x \in \mathbb{R}^n$ and let $d\sigma$ be the volume element on $\Lambda$. A canonical atlas is a locally finite countable covering of $\Lambda^n(t)$ by bounded simply-connected domains $\Omega_j$ (the charts) in each of which one can take as coordinates either the variables $x$ or $p$ a mixed collection

$$(p_\alpha, x_\beta),$$

$$\alpha = \{\alpha_1, \ldots \alpha_s\}, \qquad (0.2.14)$$

$$\beta = \{\beta_1, \ldots, \beta_{n-s}\},$$

not containing dual pairs $(p_j, x_j)$. The Maslov canonical operator $K_{\Lambda^n(t)}^{1/h}$ maps $\mathbf{C}_0^\infty(\Lambda^n(t))$ into $\mathbf{C}(\mathbb{R}_x^n)$. The canonical operators are introduced as follows:

(**1**) Let the chart $\Omega_j$ be non-degenerate, that is, $\Omega_j$ is given by an equation $p = p(x, t)$ and

$$\left(K_{\Omega_j}^{1/h}[\varphi]\right)(x) = \sqrt{\left|\frac{d\sigma}{dx}\right|} \exp\left[\frac{i}{h}\int_{r^0}^{r}(p, dx)\right]\varphi(r),$$

$$r = (x, p(x, t)), \qquad (0.2.15)$$

$$(p, dx) = \sum_{j=1}^{n} p_j dx_j.$$

Here $h \ll 1$ is a parameter, $r^0 \in \Omega_j$ is a fixed point and $\varphi \in \mathbf{C}_0^\infty(\Omega)$.

(**2**) Let the local coordinates in the chart $\Omega_j$ be $p$, that is $\Omega_j$, is given by an equation $x = x(p, t)$, and let

$$\left(K_{\Omega_j}^{1/h}(\Omega_j[\varphi])\right)(x) =$$

$$\mathcal{F}_{1/h, p \to x}^{-1}\left\{\sqrt{\left|\frac{d\sigma}{dx}\right|} \exp\left[\frac{i}{h}\left(\int_{r^0}^{r}(p, dx) - (x(p, t), p)\right)\right]\varphi(r)\right\}, \qquad (0.2.16)$$

$$r = (x(p, t), p).$$

Here

$$\mathcal{F}^{-1}_{1/h,p\to x}[\phi](x) = \frac{1}{(-2\pi i h)^{n/2}} \int_{\mathbb{R}^n} d^n p\, \phi(p) \exp\left[\frac{i}{h}(x,p)\right]. \qquad (0.2.17)$$

$K^{1/h}_{\Omega_j}$ is defined analogously in the case when the coordinates in are some collection $(p_\alpha, x_\beta)$. Let $\oint\limits_{l}(p, dx) = 0$ and let the Maslov index $\mathbf{ind}(l) = 0$ for any closed path lying on $\Lambda^n(t)$. One introduces a partition of unity of class $C^\infty$ on $\Lambda^n(t)$

$$\sum_{j=1}^{\infty} e_j(x) = 1, \mathbf{supp}(e_j) \subset \Omega_j \qquad (0.2.18)$$

and one fixes a point $r^0 = \Omega_{j0}$. The Maslov canonical operator is defined by

$$\left(K^{1/h}_{\Lambda^n(t)}[\varphi(r)]\right)(x) = \sum_j c_j \left(K^{1/h}_{\Omega_j}[e_j\varphi]\right)(x),$$

$$(0.2.19)$$

$$c_j = \exp\left(\frac{i\pi}{2}\gamma_j^M\right).$$

Here $\gamma_j^M$ is the Maslov index of a chain of charts joining the charts $\Omega_{j0}$ and $\Omega_j$.
A point $r \in \Lambda^n(t)$ is called non-singular if it has a neighbourhood in given by an equation $p = p(x)$. Let the intersection of the charts $\Omega_i$ and $\Omega_j$ be non-empty and connected, let $r \in \Omega_i \cap \Omega_j$ be a non-singular point and let $(p_\alpha, x_\beta), \left(p_{\tilde{\alpha}}, x_{\tilde{\beta}}\right)$ be the coordinates in these charts. The number

$$\gamma_{ij}^M = \mathbf{inerdex}\left(\frac{\partial x_\alpha(r)}{\partial p_\alpha}\right) - \mathbf{inerdex}\left(\frac{\partial x_{\tilde{\alpha}}(r)}{\partial p_{\tilde{\alpha}}}\right) \qquad (0.2.20)$$

is the Maslov index of the pair of charts $\Omega_i$ and $\Omega_j$. The Maslov index of a chain of charts is defined by additivity. The Maslov index of a path is defined analogously. The Maslov index of a path (mod 4) on a Lagrangian manifold is an integer homotopy invariant. The Maslov canonical operator is invariant under the choice of the canonical atlas, of local coordinates in the charts and the partition of unity in the following sense: If $K_{\Lambda^n(t)}^{1/h}, \widetilde{K}_{\Lambda^n(t)}^{1/h}$ are two Maslov canonical operators, then in $\mathcal{L}_2(\mathbb{R}^n)$,

$$\left( K_{\Lambda^n(t)}^{1/h}\varphi - \widetilde{K}_{\Lambda^n(t)}^{1/h}\varphi - \right)(x) = O(\hbar) \qquad (0.2.21)$$

for any function $\varphi \in \mathbf{C}_0^\infty(\Lambda^n(t))$. Let $L(x, hD)$ be a differential operator with real symbol $L(x, p)$ of class $\mathbf{C}^\infty$ and suppose that $L(x, p) = 0$ on $\Lambda^n(t)$. Suppose that $\Lambda^n(t)$ and the volume element $d\sigma$ are invariant under the Hamiltonian system

$$\frac{dx}{d\tau} = \frac{\partial L}{\partial p}, \frac{dp}{d\tau} = -\frac{\partial L}{\partial x}. \qquad (0.2.22)$$

Then the following commutation formula is true

$$L(x, hD)\left( K_{\Lambda^n(t)}^{1/h}\varphi \right)(x) = \frac{\hbar}{i} K_{\Lambda^n(t)}^{1/h}[\Re\varphi + O(\hbar)],$$

$$\Re\varphi = \left[ \frac{d}{d\tau} - \frac{1}{2}\sum_{j=1}^{n} \frac{\partial^2 L(x,p)}{\partial x_j \partial p_j} \right]\varphi \qquad (0.2.23)$$

where $\varphi \in \mathbf{C}_0^\infty(\Lambda^n(t))$ and $d/d\tau$ is the derivative along the integral curves of the flow of the Hamiltonian system. The equation $\Re\varphi = 0$ is called the transport equation. The commutation formula implies that if $\Re\varphi = 0$, then the function is a formal asymptotic solution of the equation $L(x, hD)u = 0$.

In general case, when the potential $(V_\varepsilon)_\varepsilon$ in $G(\mathbb{R}^n)$ any solution $(\Psi_\varepsilon(x, t, x_0; \hbar_\varepsilon))_\varepsilon$ of the Colombeau-Schrödinger equation (0.1.8) can be represented by Colombeau-Feynman path integral of the form:

$$\left(\Psi_\varepsilon(x,t,x_0;\hbar_\varepsilon)\right)_\varepsilon =$$

$$\left(\int\limits_{\substack{q(0)=y\\q(t)=x}} \left(\Psi_{0,\varepsilon}(q(0))\right)_\varepsilon \exp\left\{\frac{i}{\hbar_\varepsilon}\mathbf{S}_\varepsilon[\dot{q}(\tau),q(\tau),t]\right\}[Dq(\tau),\varepsilon]dy\right)_\varepsilon, \qquad (0.2.24)$$

$$\left(\mathbf{S}_\varepsilon[\dot{q}(\tau),q(\tau),t]\right)_\varepsilon = \left(\int\limits_0^t \mathscr{L}_\varepsilon(\dot{q}(\tau),q(\tau),\tau)d\tau\right)_\varepsilon,$$

$$\left(\mathscr{L}_\varepsilon(\dot{q}(\tau),q(\tau),\tau)\right)_\varepsilon = m\frac{\dot{q}^2(\tau)}{2} - \left(V_\varepsilon(q(\tau),t)\right)_\varepsilon.$$

Maslov quasiclassical asymptotic of the generalized wave function $\left(\Psi_\varepsilon(x,t;\hbar)\right)_\varepsilon$ is:

$$\left(\Psi_\varepsilon(x,t;\hbar_\varepsilon)\right)_\varepsilon = \left(2\pi i(\hbar_\varepsilon)_\varepsilon\right)^{-\frac{n}{2}}\left(\sum_{k=1}^{m(\varepsilon)}\int dy \Psi_{0,\varepsilon}(y)|D(q_\varepsilon^k(\tau;t,y,x))|^{\frac{1}{2}} \times\right.$$

$$\left.\exp\left\{\frac{i}{\hbar_\varepsilon}\left[\mathbf{S}_\varepsilon[q_\varepsilon^k(\tau;t,y,x)] - \frac{i\pi\gamma_{k,\varepsilon}}{2}\right]\right\}\right)_\varepsilon + \left(O(\hbar_\varepsilon)\right)_\varepsilon. \qquad (0.2.25)$$

Here

$$\left(\mathbf{S}_\varepsilon[q_\varepsilon^k(\tau; t, y, x), t]\right)_\varepsilon = \left(\int\limits_0^t \mathcal{L}_\varepsilon(\dot{q}_\varepsilon^k(\tau; t, y, x), q_\varepsilon^k(\tau; t, y, x), \tau) d\tau\right)_\varepsilon,$$

$$\left(\mathcal{L}_\varepsilon(\dot{q}_\varepsilon(\tau; t, y, x), q_\varepsilon(\tau; t, y, x), \tau)\right)_\varepsilon = \left(m \frac{\dot{q}_\varepsilon^2(\tau; t, y, x)}{2} - V_\varepsilon(q_\varepsilon(\tau; t, y, x), \tau)\right)_\varepsilon,$$

$$\left(\delta \mathbf{S}_\varepsilon[q_\varepsilon^k(\tau; t, y, x), t]\right)_\varepsilon = 0, (q_\varepsilon^k(0; y, x, t))_\varepsilon = y, (q_\varepsilon^k(t; t, y, x))_\varepsilon = x, \qquad (0.2.26)$$

$$\left(D(q_\varepsilon^k(\tau; t, y, x))\right)_\varepsilon = \left(\det\left\|\frac{\partial \mathbf{S}[q_\varepsilon^k(\tau; t, y, x), t]}{\partial x \partial y}\right\|\right)_\varepsilon \neq 0, x \neq y.$$

Formula (0.2.25) can be written by using Colombeau-Maslov canonical operator $\left(K_{\Lambda_\varepsilon^n(t)}^{1/\hbar_\varepsilon}[\bullet]\right)_\varepsilon$. Let $\Lambda_\varepsilon^n(t), \forall \varepsilon \in (0, 1]$ be an $n$-dimensional Colombeau-Lagrangian manifold of class $G(\mathbb{R}^{2n})$ in the phase space $\mathbb{R}_{x,p}^{2n}$ where $x \in \mathbb{R}^n$ and let $d\sigma_\varepsilon$ be the volume element on $\Lambda_\varepsilon^n(t)$. A canonical atlas is a locally finite countable covering of $\Lambda_\varepsilon^n(t)$ by bounded simply-connected domains $\Omega_{j,\varepsilon}$ (the $\varepsilon$-charts) in each of which one can take as coordinates either the variables $x$ or $p$ a mixed collection

$$(p_\alpha, x_\beta),$$

$$\alpha = \{\alpha_1, \dots \alpha_s\}, \qquad (0.2.27)$$

$$\beta = \{\beta_1, \dots, \beta_{n-s}\},$$

not containing dual pairs $(p_j, x_j)$. The Colombeau-Maslov canonical operator $\left(K_{\Lambda_\varepsilon^n(t)}^{1/\hbar}\right)_\varepsilon$ maps $(\mathbf{C}_0^\infty(\Lambda^n(t)))^I$ into $(\mathbf{C}(\mathbb{R}_x^n))^I, I = (0, 1]$. The generalized canonical operators are introduced as follows:
(**1**) Let the $\varepsilon$-chart $\Omega_{j,\varepsilon}$ be non-degenerate, that is, $\Omega_{j,\varepsilon}$ is given by an equation $p_\varepsilon = p_\varepsilon(x, t)$ and

$$\left( \left( K_{\Omega_{j,\varepsilon}}^{1/h_\varepsilon}[\varphi_\varepsilon] \right)(x) \right)_\varepsilon = \sqrt{\left| \left( \frac{d\sigma_\varepsilon}{dx} \right)_\varepsilon \right|} \left( \exp\left[ \frac{i}{h_\varepsilon} \int_{r_\varepsilon^0}^{r_\varepsilon} (p_\varepsilon, dx) \right] \varphi_\varepsilon(r) \right)_\varepsilon,$$

$$(0.2.28)$$

$$r_\varepsilon = (x, p_\varepsilon(x, t)),$$

$$(p_\varepsilon, dx) = \sum_{j=1}^{n} p_{j,\varepsilon} dx_j.$$

Here $(h_\varepsilon)_\varepsilon \in \widetilde{\mathbb{R}}$ is a parameter, $r_\varepsilon^0 \in \Omega_{j,\varepsilon}$ is a fixed point and $\varphi_\varepsilon \in G(\Omega)$.

**(2)** Let the local coordinates in the $\varepsilon$-chart $\Omega_{j,\varepsilon}$ be $p$, that is $\Omega_{j,\varepsilon}$, is given by an equation $x_\varepsilon = x_\varepsilon(p, t)$, and let

$$\left( \left( K_{\Omega_{j,\varepsilon}}^{1/h,\varepsilon}(\Omega_{j,\varepsilon}[\varphi_\varepsilon]) \right)(x) \right)_\varepsilon =$$

$$\left( \mathcal{F}_{1/h,\varepsilon, p \to x}^{-1} \left\{ \sqrt{\left| \frac{d\sigma_\varepsilon}{dx} \right|} \exp\left[ \frac{i}{h} \left( \int_{r_\varepsilon^0}^{r_\varepsilon} (p, dx) - (x_\varepsilon(p, t), p) \right) \right] \varphi_\varepsilon(r) \right\} \right)_\varepsilon,$$

$$(0.2.29)$$

$$r_\varepsilon = (x_\varepsilon(p, t), p).$$

Here

$$\left( \left( \mathcal{F}_{1/h_\varepsilon, p \to x}^{-1}[\phi_\varepsilon] \right)_\varepsilon \right)(x) = \left( \frac{1}{(-2\pi i \hbar_\varepsilon)^{n/2}} \int_{\mathbb{R}^n} d^n p \, \phi_\varepsilon(p) \exp\left[ \frac{i}{h}(x, p) \right] \right)_\varepsilon.$$

$$(0.2.30)$$

$\left( K_{\Omega_{j,\varepsilon}}^{1/h_\varepsilon} \right)_\varepsilon$ is defined analogously in the case when the coordinates in are some collection $(p_\alpha, x_\beta)$. Let $\left( \oint_l (p_\varepsilon, dx) \right)_\varepsilon = 0$ and let the Maslov index $(\mathbf{ind}(l_\varepsilon))_\varepsilon = 0$ for any closed path $(l_\varepsilon)_\varepsilon$ lying on $\Lambda_\varepsilon^n(t)$. One introduces a partition of unity of class $C^\infty$ on $\Lambda_\varepsilon^n(t)$

$$\sum_{j=1}^{\infty} e_{j,\varepsilon}(x) = 1, \mathbf{supp}(e_{j,\varepsilon}) \subset \Omega_{j,\varepsilon} \qquad (0.2.31)$$

and one fixes a point $r_\varepsilon^0 = \Omega_{j_0,\varepsilon}$. The Colombeau-Maslov canonical operator is defined by

$$\left(\left(K_{\Lambda_\varepsilon^n(t)}^{1/\hbar_\varepsilon}[\varphi_\varepsilon(r)]\right)_\varepsilon\right)(x) = \left(\left(\sum_j c_{j_\varepsilon}\left(K_{\Omega_{j,\varepsilon}}^{1/\hbar,\varepsilon}[e_{j,\varepsilon}\varphi_\varepsilon]\right)\right)_\varepsilon\right)(x),$$

$$(0.2.32)$$

$$c_{j,\varepsilon} = \exp\left(\frac{i\pi}{2}\gamma_{j,\varepsilon}^M\right).$$

Here $\gamma_{j,\varepsilon}^M$ is the Maslov index of a chain of charts joining the $\varepsilon$-charts $\Omega_{j_0,\varepsilon}$ and $\Omega_{j,\varepsilon}$. A point $r_\varepsilon \in \Lambda_\varepsilon^n(t)$ is called non-singular if it has a neighbourhood in given by an equation $p = p_\varepsilon(x)$. Let the intersection of the $\varepsilon$-charts $\Omega_{i,\varepsilon}$ and $\Omega_{j,\varepsilon}$ be non-empty and connected, let $\forall \varepsilon \in (0,1] : r_\varepsilon \in \Omega_{i,\varepsilon} \cap \Omega_{j,\varepsilon}$ be a non-singular point and let $(p_\alpha, x_\beta), (p_{\widetilde{\alpha}}, x_{\widetilde{\beta}})$ be the coordinates in these charts. The number

$$\left(\gamma_{ij,\varepsilon}^M\right)_\varepsilon = \left(\mathbf{inerdex}\left(\frac{\partial x_\alpha(r_\varepsilon)}{\partial p_\beta}\right)\right)_\varepsilon - \left(\mathbf{inerdex}\left(\frac{\partial x_{\widetilde{\alpha}}(r_\varepsilon)}{\partial p_{\widetilde{\alpha}}}\right)\right)_\varepsilon \qquad (0.2.33)$$

is the generalized Maslov index of the pair of $\varepsilon$-nets $\Omega_{i,\varepsilon}$ and $\Omega_{j,\varepsilon}, \varepsilon \in (0,1]$. The Maslov index of a chain of netss is defined by additivity. The generalized Maslov index of a path is defined analogously. The Maslov index of a path (mod 4) on a Colombeau-Lagrangian manifold is an integer homotopy invariant. The Maslov canonical operator is invariant under the choice of the canonical atlas, of local coordinates in the charts and the partition of unity in the following sense: if $\left(K_{\Lambda_\varepsilon^n(t)}^{1/\hbar_\varepsilon}\right)_\varepsilon, \left(\widetilde{K}_{\Lambda_\varepsilon^n(t)}^{1/\hbar_\varepsilon}\right)_\varepsilon$ are two Colombeau-Maslov canonical operators, then in $G_{\mathcal{L}_2(\mathbb{R}^n)}$

$$\left(\left(K_{\Lambda_\varepsilon^n(t)}^{1/\hbar_\varepsilon}\varphi_\varepsilon - \widetilde{K}_{\Lambda_\varepsilon^n(t)}^{1/\hbar_\varepsilon}\varphi_\varepsilon - \right)_\varepsilon\right)(x) = \left(O(\hbar_\varepsilon)\right)_\varepsilon \qquad (0.2.34)$$

for any generalized function $(\varphi_\varepsilon)_\varepsilon \in G_{\mathcal{L}_2}(\Lambda_\varepsilon^n(t))$. Let $(L_\varepsilon(x, \hbar_\varepsilon D))_\varepsilon$ be a differential operator with real symbol $(L_\varepsilon(x,p))_\varepsilon$ of class $G(\mathbb{R}_{x,p}^{2n})$ and suppose that $\forall \varepsilon \in (0,1] : L_\varepsilon(x,p) = 0$ on $\Lambda_\varepsilon^n(t)$. Suppose that $(\Lambda_\varepsilon^n(t))_\varepsilon$ and the volume element $(d\sigma_\varepsilon)_\varepsilon$ are invariant under the Hamiltonian system

$$\left(\frac{dx_\varepsilon(\tau)}{d\tau}\right)_\varepsilon = \left(\left[\frac{\partial L_\varepsilon(x_\varepsilon(\tau),p)}{\partial p}\right]_{p=p_\varepsilon(\tau)}\right)_\varepsilon,$$

$$\left(\frac{dp_\varepsilon(\tau)}{d\tau}\right)_\varepsilon = -\left(\left[\frac{\partial L_\varepsilon(x,p_\varepsilon(\tau))}{\partial x}\right]_{x=x_\varepsilon(\tau)}\right)_\varepsilon.$$

(0.2.35)

Then the following commutation formula is true

$$\left(\left(L_\varepsilon(x,\hbar_\varepsilon D)\left(K^{1/\hbar_\varepsilon}_{\Lambda^n_\varepsilon(t)}\varphi_\varepsilon\right)\right)_\varepsilon\right)(x) =$$

$$\frac{(\hbar_\varepsilon)_\varepsilon}{i}\left(\left(K^{1/\hbar_\varepsilon}_{\Lambda^n_\varepsilon(t)}\right)_\varepsilon\right)[(\Re_\varepsilon\varphi_\varepsilon)_\varepsilon + (O(\hbar_\varepsilon))_\varepsilon],$$

(0.2.36)

$$(\Re_\varepsilon\varphi_\varepsilon)_\varepsilon = \left(\left[\frac{d_\varepsilon}{d_\varepsilon\tau} - \frac{1}{2}\sum_{j=1}^n\frac{\partial^2 L_\varepsilon(x,p)}{\partial x_j\partial p_j}\right]\varphi_\varepsilon\right)_\varepsilon$$

where $\forall\varepsilon\in(0,1]: \varphi_\varepsilon\in \mathbf{C}^\infty_0(\Lambda^n_\varepsilon(t))$ and $d_\varepsilon/d_\varepsilon\tau$ is the derivative along the integral curves of the flow of the Hamiltonian system (0.2.35). The equation $(\Re_\varepsilon\varphi_\varepsilon)_\varepsilon = 0$ is called the generalized (Colombeau) transport equation. The commutation formula implies that if $(\Re_\varepsilon\varphi_\varepsilon)_\varepsilon = 0$, then the function $(u_\varepsilon)_\varepsilon = \left(K^{1/\hbar_\varepsilon}_{\Lambda^n_\varepsilon(t)}\varphi_\varepsilon\right)$ is a formal Colombeau asymptotic solution of the equation $(L_\varepsilon(x,\hbar_\varepsilon D)u_\varepsilon)_\varepsilon = 0$.

Let $\mathbf{H}$ be a complex infinite dimensional separable Hilbert space, with inner product $\langle\cdot,\cdot\rangle$ and norm $\|\cdot\|$. Let us consider now Colombeau-Couchy problem:

$$-i((\hbar_\varepsilon)_\varepsilon)\left(\frac{\partial\Psi_\varepsilon(t)}{\partial t}\right)_\varepsilon + \left(\hat{H}_\varepsilon\Psi_\varepsilon(t)\right)_\varepsilon + ((\hbar_\varepsilon)_\varepsilon)(f_\varepsilon(t))_\varepsilon = 0,$$

(0.2.37)

$$(\Psi_\varepsilon(0))_\varepsilon = (\Psi_{0,\varepsilon})_\varepsilon, 0 \le t \le T.$$

**Lemma 0.2.1.** Assume that:

**(i)** $\forall \varepsilon \in (0,1]$ : operator $H_\varepsilon$ is essentially self-adjoint,

**(ii)** $\forall \varepsilon \in (0,1]$ : operator $\widehat{H}_\varepsilon$ is the closure of $H_\varepsilon$, $\Psi_\varepsilon(t) : \mathbb{R} \to \mathbf{H}$ and

**(iii)** $\forall \varepsilon \in (0,1]$ and $\forall t (0 \leq t \leq T)$ : function $f_\varepsilon(t) : \mathbb{R}_+ \to \mathbf{H}$ is continuous.

**(iv)** $(\Psi_{0,\varepsilon})_\varepsilon \in D\left(\left(\widehat{H}_\varepsilon\right)_\varepsilon\right)$.

Then there exists unique Colombeau solution $(\Psi_\varepsilon(t))_\varepsilon, 0 \leq t \leq T$ of the Colombeau-Couchy problem (0.2.37) and

$$\left(\|\Psi_\varepsilon(t)\|\right)_\varepsilon \leq \left(\|\Psi_\varepsilon(0)\|\right)_\varepsilon + \left(\int_0^T \|f_\varepsilon(t)\| dt\right)_\varepsilon. \tag{0.2.38}$$

**Proof.** From **Eq.**(0.2.37) one obtain

$$i\left(\frac{\partial \Psi_\varepsilon(t)}{\partial t}\right)_\varepsilon = \frac{1}{(\hbar_\varepsilon)_\varepsilon}\left(\widehat{H}_\varepsilon \Psi_\varepsilon(t)\right)_\varepsilon + (f_\varepsilon(t))_\varepsilon,$$

$$\tag{0.2.39}$$

$$\left(\Psi_\varepsilon(0)\right)_\varepsilon = (\Psi_{0,\varepsilon})_\varepsilon, 0 \leq t \leq T.$$

Under assumptions **(i)**-**(iv)** from Eq.(0.2.39) one obtain $\forall t \in [0,T]$ :

$$(\Psi_\varepsilon(t))_\varepsilon = \exp\left(-it\left(\frac{\widehat{H}_\varepsilon}{\hbar_\varepsilon}\right)_\varepsilon\right) + \left(\int_0^t \exp\left(\frac{i(\tau - t)\widehat{H}_\varepsilon}{\hbar_\varepsilon}\right) f_\varepsilon(\tau) d\tau\right)_\varepsilon. \tag{0.2.40}$$

Thus estimate (0.2.38) follows immediately from Eq.(0.2.40) and unitarity $\forall \alpha \left(\alpha \in \widetilde{\mathbb{R}}\right)$ of the operator $\exp\left(i\alpha\left(\widehat{H}_\varepsilon\right)_\varepsilon\right)$.

The next simple lemma for the classical case, given in [22] (see [22] chapter II, section 10.3).

**Lemma 0.2.2.** Let $(\Psi_\varepsilon(t,x;\hbar_\varepsilon))_\varepsilon, 0 \leq t \leq T$ be the Colombeau solution of the Couchy problem for the Colombeau-Schrödinger equation (0.2.42) and let $(\Psi_\varepsilon^{FAS}(t,x;\hbar_\varepsilon))_\varepsilon$ be formal asymptotic Colombeau solution of the Colombeau-Schrödinger equation (0.2.42) given by using Colombeau-Maslov operator,i.e.

$$\left(\Psi_\varepsilon^{FAS}(t,x;\hbar_\varepsilon)\right)_\varepsilon = \left(\left(K_{\Lambda_\varepsilon^t(t)}^{1/\hbar_\varepsilon}\Psi_{0,\varepsilon}\right)(x)\right)_\varepsilon. \qquad (0.2.41)$$

Then an asimptotic expansion [55],[56] of the net $(\Psi_\varepsilon(t,x;\hbar_\varepsilon))_\varepsilon$ coincide with $\left(\Psi_\varepsilon^{FAS}(t,x;\hbar_\varepsilon)\right)_\varepsilon$.

**The main result of this paper is**:

**Theorem 0.2.1**.Let **H** be a complex infinite dimensional separable Hilbertspace, with inner product $\langle\cdot,\cdot\rangle$ and norm $|\cdot|$. Let us consider Colombeau-Schrödinger equation:

$$-i((\hbar_\varepsilon)_\varepsilon)\left(\frac{\partial\Psi_\varepsilon(t)}{\partial t}\right)_\varepsilon + \left(\widehat{H}_\varepsilon(t)\Psi_\varepsilon(t)\right)_\varepsilon = 0,$$

$$(0.2.42)$$

$$(\Psi_\varepsilon(0))_\varepsilon = (\Psi_{0,\varepsilon}(x))_\varepsilon, x \in \mathbb{R}^n,$$

where $\forall\varepsilon \in (0,1]$ operator $H_\varepsilon(t)$ given via formula

$$H_\varepsilon(t) = -\frac{\hbar_\varepsilon^2}{2m}\Delta + V_\varepsilon(x,t), \qquad (0.2.43)$$

where $\hbar_\varepsilon = O(\varepsilon^k), k \in \mathbb{N}$ if $\varepsilon \to 0$ and :

**(i)** **cl**$[(V_\varepsilon(x,\cdot))_\varepsilon] \in G(\mathbb{R}^n)$,

**(ii)** $\forall\varepsilon \in (0,1]$ : operator $H_\varepsilon$ is essentially self-adjoint,

**(iii)** $\forall\varepsilon \in (0,1]$ : operator $\widehat{H}_\varepsilon$ is the closure of $H_\varepsilon, \Psi_\varepsilon(t) : \mathbb{R} \to \mathbf{H}$ and

**(iv)** $\forall\varepsilon \in (0,1]$ : $\widehat{H}_\varepsilon(t) : \mathbb{R}_+ \times \mathbf{H} \to \mathbf{H}$.

Let us consider Colombeau-Schrödinger equation (0.2.42) with initial condition

$(\Psi_{0,\varepsilon}(\mathbf{x}))_\varepsilon = \left(\sqrt{\Psi_\varepsilon(\mathbf{x},\mathbf{x_0})}\right)_\varepsilon$, where

$$\Psi_\varepsilon(\mathbf{x}, \mathbf{x}_0) = \frac{\eta_\varepsilon^{n/2}}{(2\pi)^{n/2} \hbar_\varepsilon^{n/2}} \exp\left(-\frac{\eta_\varepsilon(\mathbf{x} - \mathbf{x}_0)^2}{\hbar_\varepsilon}\right),$$

(0.2.44)

$$0 < \hbar_\varepsilon \ll \eta_\varepsilon \ll 1 \text{ if } \varepsilon \to 0,$$

i.e. $(\Psi_{0,\varepsilon}^2(x))_\varepsilon = \delta(x - x_0)$. We assume that $V_\varepsilon(\cdot, t)|_{\varepsilon=0} = V(\cdot, t) : \mathbb{R}^n \to \mathbb{R}$ is an polynomial function on variable $\mathbf{x} = (x_1, x_2, \ldots, x_n)$, i.e.

$$V(\mathbf{x}, t) = \sum_{\|\alpha\| \leqslant n} g_\alpha(t) x^\alpha, \alpha = (i_1, \ldots, i_n), \|\alpha\| = \sum_r i_r, i = 1, \ldots, n,$$

(0.2.45)

$$x^\alpha = x_1^{i_1} x_2^{i_2} \ldots x_n^{i_n}$$

Let $\mathbf{u}(\tau, t, \lambda, \mathbf{x}) = (u_1(\tau, t, \lambda, \mathbf{x}), \ldots, u_n(\tau, t, \lambda, \mathbf{x}))$ be the solution of the linear differential boundary problem (*differential master equation*):

$$\frac{d^2 \mathbf{u}(\tau, t, \lambda, \mathbf{y}, \mathbf{x})}{d\tau^2} = \mathbf{Hess}[V(\lambda, \tau)] \mathbf{u}^{\mathbf{T}}(\tau, t, \lambda, \mathbf{y}, \mathbf{x}) + [\mathbf{V}'(\lambda, \tau)]^{\mathbf{T}},$$

(0.2.46)

$$\mathbf{u}(0, t, \lambda, \mathbf{y}, \mathbf{x}) = \mathbf{y}, \mathbf{u}(t, t, \lambda, \mathbf{y}, \mathbf{x}) = \mathbf{x}.$$

Here $\mathbf{y}, \mathbf{x} \in \mathbb{R}^n$, $\lambda = (\lambda_1, \ldots, \lambda_n)$, $\mathbf{u}^{\mathbf{T}}(\tau, t, \lambda, \mathbf{y}, \mathbf{x}) = (u_1(\tau, t, \lambda, \mathbf{y}, \mathbf{x}), \ldots, u_n(\tau, t, \lambda, \mathbf{y}, \mathbf{x}))^{\mathbf{T}}$,

$$\mathbf{V}'(\lambda, \tau) = \left(\frac{\partial V(\mathbf{x}, \tau)}{\partial x_1}, \ldots, \frac{\partial V(\mathbf{x}, \tau)}{\partial x_n}\right)_{\mathbf{x} = \lambda}$$

(0.2.47)

and here $\mathbf{Hess}[V(\lambda, \tau)]$ is the Hessian of the function $V(\cdot, \tau) : \mathbb{R}^n \to \mathbb{R}$ at point $\mathbf{x} = \lambda$, i.e.

$$\mathbf{Hess}[V(\boldsymbol{\lambda}, \tau)] =$$

$$\mathbf{Hess}_{ij}[V(\boldsymbol{\lambda}, \tau)] = \left[ \frac{\partial^2 V(\mathbf{x}, \tau)}{\partial x_i \partial x_j} \right]_{\mathbf{x}=\boldsymbol{\lambda}} \qquad (0.2.48)$$

and here $\mathbf{Hess}_{ij}[\mathbf{V}(\boldsymbol{\lambda}, t)]$ is the Hessian of the functiontion $\mathbf{V}(\cdot, t) : \mathbb{R}^n \to \mathbb{R}$ at point $\mathbf{x} = \boldsymbol{\lambda}$. Let $\mathbf{S}(t, \boldsymbol{\lambda}, \mathbf{y}, \mathbf{x})$ be a function (*master action* corresponding to Eq.(0.2.46))

$$\mathbf{S}(t, \boldsymbol{\lambda}, \mathbf{y}, \mathbf{x}) = \int_0^t \widetilde{\mathcal{L}}(\dot{\mathbf{u}}(\tau, t, \boldsymbol{\lambda}, \mathbf{y}, \mathbf{x}), \mathbf{u}(\tau, t, \boldsymbol{\lambda}, \mathbf{y}, \mathbf{x}), \tau) d\tau, \qquad (0.2.49)$$

where function $\widetilde{\mathcal{L}}(\dot{\mathbf{u}}(\tau), \mathbf{u}(\tau), \tau)$ (*master Lagrangian* corresponding to Eq.(0.2.46)) is

$$\widetilde{\mathcal{L}}(\dot{\mathbf{u}}(\tau), \mathbf{u}(\tau), \tau) = m\frac{\dot{\mathbf{u}}^2(\tau)}{2} - \widehat{V}(\mathbf{u}(\tau), \tau),$$

$$\widehat{V}(\mathbf{u}(\tau), \tau) = [V'(\boldsymbol{\lambda}, \tau)]\mathbf{u}^{\mathbf{T}}(\tau) + \frac{1}{2}\mathbf{u}(\tau, \boldsymbol{\lambda})\mathbf{Hess}[\mathbf{V}(\boldsymbol{\lambda}, \tau)]\mathbf{u}^{\mathbf{T}}(\tau) \qquad (0.2.50)$$

Let $\mathbf{y_{cr}}(t, \boldsymbol{\lambda}, \mathbf{x})$ be solution for a system of equations:

$$\frac{\partial \mathbf{S}(t, \boldsymbol{\lambda}, \mathbf{x}, \mathbf{y_{cr}}(t, \boldsymbol{\lambda}, \mathbf{y}))}{\partial y_i} = 0,$$

$$i = 1, \ldots, n. \qquad (0.2.51)$$

Let $\widehat{\mathbf{u}}(t, \boldsymbol{\lambda}, \mathbf{x}_0) = \mathbf{x_{cr}}(t, \boldsymbol{\lambda}, \mathbf{x}_0)$ be solution for a system of equations:

$$\mathbf{y_{cr}}(t, \boldsymbol{\lambda}, \mathbf{x_{cr}}) + \boldsymbol{\lambda} - \mathbf{x}_0 = 0. \qquad (0.2.52)$$

We assume now that: for a fixed $t, \mathbf{x}_0$ and $\boldsymbol{\lambda}$ point $x_{\mathbf{cr}}(t, \boldsymbol{\lambda}, \mathbf{x}_0)$ does not focal point on the trajectory, given via solution of the linear boundary problem (0.2.46).

Then for any Colombeau solution of the Colombeau-Schrödinger equation Eq.(0.2.42) with initial condition (0.2.43) the inequalities is satisfied:

$$\limsup_{\varepsilon \to 0} |\langle i, t, \mathbf{x}_0; \hbar_\varepsilon, \varepsilon \rangle - \lambda_i| \leq |\hat{u}_i(t, \boldsymbol{\lambda}, \mathbf{x}_0)|,$$

$$(0.2.53)$$

$$i = 1, \ldots, n.$$

Here the $\langle i, t, \mathbf{x}_0; \hbar_\varepsilon, \varepsilon \rangle$ quantum averages given by Eq.(0.2.2),i.e.

$$\langle i, t, \mathbf{x}_0; \hbar_\varepsilon, \varepsilon \rangle = \frac{\displaystyle\int_{\mathbb{R}^n} x_i |\Psi_\varepsilon(\mathbf{x}, t, \mathbf{x}_0; \hbar)|^2 d^n x}{\displaystyle\int_{\mathbb{R}^n} |\Psi_\varepsilon(\mathbf{x}, t, \mathbf{x}_0; \hbar)|^2 d^n x},$$

$$(0.2.54)$$

$$i = 1, \ldots, n.$$

**Example**.Let's calculate master equation for the one dimensional quantum anharmonic oscillator with a cubic potential,supplemented by an additive sinusoidal driving i.e.

$$V(x, \tau) = \frac{m\omega^2}{2} x^2 + ax^3 - bx - [A \sin(\Omega \tau)]x. \qquad (0.2.55)$$

The corresponding classical Lagrangian is

$$\mathcal{L}(\dot{x}(\tau), x(\tau), \tau) = \frac{m}{2}\dot{x}^2(\tau) - \frac{m\omega^2}{2}x^2(\tau) - ax^3(\tau) + bx(\tau) + [A\sin(\Omega\tau)]x(\tau). \quad (0.2.56)$$

From Eq.(0.2.48 ) and Eq.(0.2.50 ) using Eq.(0.2.55) we obtain

$$\mathbf{V}'(\lambda, \tau) = \left[\frac{\partial V(x, \tau)}{\partial x}\right]_{x=\lambda} = [m\omega^2 x + 3ax^2 - b - A\sin(\Omega\tau)]_{x=\lambda} =$$

$$m\omega^2\lambda + 3a\lambda^2 - b - A\sin(\Omega\tau),$$

$$\mathbf{Hess}[V(\lambda, \tau)] = \left[\frac{\partial^2 V(x, \tau)}{\partial x^2}\right]_{x=\lambda} = [m\omega^2 + 6ax]_{x=\lambda} = m\omega^2 + 6a\lambda, \quad (0.2.57)$$

$$\widehat{V}(u, \tau) = [V'(\lambda, \tau)]u + \frac{1}{2}\mathbf{Hess}[\mathbf{V}(\lambda, \tau)]u^2 =$$

$$\frac{1}{2}(m\omega^2 + 6a\lambda)u^2 + (m\omega^2\lambda + 3a\lambda^2 - b - A\sin(\Omega\tau))u.$$

Therefore differential master equation (0.2.46) is

$$\frac{d^2u(\tau, t, \lambda, y, x)}{d\tau^2} = (m\omega^2 + 6a\lambda)u(\tau, t, \lambda, y, x) + (m\omega^2\lambda + 3a\lambda^2 - b - A\sin(\Omega\tau)),$$
$$(0.2.58)$$

$$u(0, t, \lambda, y, x) = y, u(t, t, \lambda, y, x) = x.$$

The corresponding master Lagrangian $\widetilde{\mathcal{L}}(\dot{u}, u)$ is

$$\widetilde{\mathcal{L}}(\dot{u}, u) =$$

$$\widetilde{\mathcal{L}}(\dot{u}, u) = \frac{m}{2}\dot{u}^2 - \frac{1}{2}(m\omega^2 + 6a\lambda)u^2 - (m\omega^2\lambda + 3a\lambda^2 - b - A\sin(\Omega t))u = \qquad (0.2.59)$$

$$\frac{m}{2}\dot{u}^2 - \frac{m}{2}\left(\omega^2 + \frac{6a\lambda}{m}\right)u^2 - (m\omega^2\lambda + 3a\lambda^2 - b - A\sin(\Omega t))u.$$

We assume now that $\omega^2 + \frac{6a\lambda}{m} \geq 0$ and rewrite master Lagrangian $\widetilde{\mathcal{L}}(\dot{u}, u)$ given by Eq.(0.2.59) in the next form

$$\widetilde{\mathcal{L}}(\dot{u}, u) = \frac{m}{2}\dot{u}^2 - \frac{m\varpi^2(\lambda)}{2}u^2 + g(\lambda, \tau)u,$$

$$\varpi(\lambda) = \left| \sqrt{\omega^2 + \frac{6a\lambda}{m}} \right|, \qquad (0.2.60)$$

$$g(\lambda, \tau) = -(m\omega^2\lambda + 3a\lambda^2 - b - A\sin(\Omega\tau)).$$

The Euler-Lagrange equation of $\widetilde{\mathcal{L}}(\dot{u}, u)$ with corresponding boundary conditions is

$$\frac{d}{dt}\left( \frac{\partial\widetilde{\mathcal{L}}}{\partial\dot{u}} \right) - \frac{\partial\widetilde{\mathcal{L}}}{\partial u} =$$

$$m\ddot{u} + m\varpi^2(\lambda) - g(\lambda, \tau) = 0, \qquad (0.2.61)$$

$$u(0, \tau, \lambda, y, x) = y, u(t, t, \lambda, y, x) = x.$$

Corresponding master action $\mathbf{S}(t, \lambda, y, x)$ given by Eq.(0.2.49) is :

$$\mathbf{S}(t,\lambda,y,x) = \frac{m\varpi}{2\sin\varpi T}\Bigg[(\cos(\varpi t))(y^2 + x^2) - 2yx + \frac{2x}{m\varpi}\int_0^t g(\lambda,\tau)\sin(\varpi\tau)d\tau -$$

$$+ \frac{2y}{m\varpi}\int_0^t g(\lambda,\tau)\sin[\varpi(t-\tau)]d\tau - \qquad (0.2.62)$$

$$- \frac{2}{m^2\varpi^2}\int_0^t\int_0^\tau g(\lambda,\tau)g(\lambda,s)[\sin\varpi(t-\tau)]\sin(\varpi s)dsd\tau\Bigg].$$

System of the Eqs.(0.2.51) reduces to a single equation

$$\frac{\partial\mathbf{S}(t,\lambda,y,x)}{\partial x} = x_{\mathbf{cr}}(t,\lambda,y)\cos(\varpi t) - y + \frac{1}{m\varpi}\int_0^t g(\lambda,\tau)\sin(\varpi\tau)d\tau = 0. \qquad (0.2.63)$$

Therefore function $\hat{u}(t,\lambda,y) = x_{\mathbf{cr}}(t,\lambda,y)$ is

$$x_{\mathbf{cr}}(t,\lambda,y) = \frac{1}{\cos(\varpi t)}\left(y - \frac{1}{m\varpi}\int_0^t g(\lambda,\tau)\sin(\varpi\tau)d\tau\right) \qquad (0.2.64)$$

System of the Eqs.(0.2.52) reduces to a single equation

$$\frac{1}{\cos(\varpi t)}\left(y_{\mathbf{cr}}(t,\lambda,x_0) - \frac{1}{m\varpi}\int_0^t g(\lambda,\tau)\sin(\varpi\tau)d\tau\right) + \lambda - x_0 = 0. \qquad (0.2.65)$$

# I.0.3. Colombeau-Wiener measure $\left(\widetilde{D}_x^{\mathbb{R}}[\omega(\tau),\nu_\varepsilon]\right)_\varepsilon$ and regularized Colombeau-Feynman Path Integral

**based on measure** $\left( \widetilde{D}_x^{\mathbb{R}}[\omega(\tau), \nu_\varepsilon] \right)_\varepsilon$ .

Let $x = (x_1, x_2, \ldots, x_n,)$ denote a point in $\mathbb{R}^d$ and $|x|$ denote its Euclidean length from the origin. Let us consider the initial-value problem for the Schrödinger equation (0.1.1) with $V \equiv 0$, i.e.

$$\frac{\partial \psi(t,x)}{\partial t} = \frac{i}{2m} \Delta \psi(t,x), \psi(0) = \psi_0, x \in \mathbb{R}^n. \qquad (0.3.1)$$

On the domain $\mathbf{D}(\Delta)$ of all square-integrable $\psi$, such that also $-|p|^2 \mathcal{F} \psi \in \mathcal{L}_2(\mathbb{R}^d)$ one can write

$$\Delta \psi = \mathcal{F}^{-1}[(-|p|^2)[\mathcal{F}\psi](p)]. \qquad (0.3.2)$$

Here $p$ denotes the variable in momentum space and $|p|^2 = \sum_{i=1}^d p_i^2$. Therefore $\Delta$ is self-adjoint, and if we denote

$$\Xi_m^t = \exp\left[ \frac{it}{2m} \Delta \right], \qquad (0.3.3)$$

then $\psi(t,x) = \Xi_m^t \psi_0(x)$ is the solution of (0.3.1). Explicitly,

$$\Xi_m^t \psi_0(x) = \left( \frac{2\pi it}{m} \right)^{-d/2} \int_{\mathbb{R}^n} d^d y\, \psi_0(y) \exp\left[ i\frac{m}{2} \left( \frac{|x-y|^2}{t} \right) \right] \qquad (0.3.4)$$

The integral in (0.3.4) exists for all $\psi \in \mathcal{L}_1 \cap \mathcal{L}_2$ for a general $\psi_0 \in \mathcal{L}_2$ we define $\Xi_m^t \psi_0(x)$ as the limit $\Xi_{m,\Lambda}^t \psi_0(x)$ in $\mathcal{L}_2$ as $\Lambda \to \infty$ of the integral (0.3.4) restricted to $|y| \leq \Lambda$. Suppose now that $1/2m$ is replaced by $0$ in Eq.(0.1.1). The operator $V$ of multiplication by the function $V(t,x)$, on the domain $\mathbf{D}(V)$ of all $\psi_0$ in $\mathcal{L}_2$ such that $V\psi_0$ is also in $\mathcal{L}_2$, is self-adjoint, and if

$$M_V^t = \exp\{-itV\},$$ 

(0.3.5)

then $\psi(t, x) = M_V^t \psi_0(x)$ is the solution of Eq.(0.3.1) with $1/2m$ replaced by $0$.

Asuume that Kato conditions (see Apendix **V**) under which the operator $(1/2m)\Delta - V$ is self-adjoint is satisfied.It well known [41] that under Kato conditions if we let

$$U_{m,V}^t = \exp\left[it\left(\frac{1}{2m}\Delta - V\right)\right]$$

(0.3.6)

then Trotter's theorem asserts that for all $\psi_0 \in \mathscr{L}_2$ :

$$U_{m,V}^t \psi_0 = \lim_{n \to \infty} \left(\Xi_{\frac{t}{m}}^{\frac{t}{n}} M_V^{\frac{t}{n}}\right)^n \psi_0.$$

(0.3.7)

Using (0.3.4) and (0.3.5) one find that

$$U_{m,V}^t \psi_0 = \lim_{n \to \infty} \left[\left(\frac{2\pi it}{nm}\right)^{-\frac{1}{2}dn} \int \ldots \int \exp[iS(x_0, \ldots, x_n)]\psi_0(x_n)dx_1 \ldots dx_n\right].$$

(0.3.8)

Here we set $x_0 = x$ and

$$S(x_0, \ldots, x_n) = \left(\frac{t}{n}\right) \sum_{j=1}^n \left[\frac{m}{2} \frac{(x_i - x_{j-1})^2}{(t/n)^2} - V(x_i)\right].$$

(0.3.9)

The limit in RHS of the Eq.(0.3.8) is well-defined, and according to Eq.(0.3.7) converges in $\mathscr{L}_2$ to the solution of the Schrödinger equation (0.1.1).Therefore using the results of Kato and Trotter one obtain a precise meaning to the Feynman integral (0.2.4) as limit in RHS of the Eq.(0.3.8), when the potential $V$ is sufficiently regular.

Let $\omega(\tau)$ be a trajectory; that is, a function from $[0, t]$ to $\mathbb{R}^d$, with $\omega(0) = x$ and set $x_j = \omega(jt/n)$ for $j = 0, \ldots, n$. Then (0.3.9) is, formally, the Riemann sum for the classical action

$$S_t(\dot{\omega}(\tau), \omega(\tau)) = \int\limits_0^t \left[ \frac{m}{2} \dot{\omega}^2(\tau) - V(\omega(\tau)) \right] d\tau \qquad (0.3.10)$$

of a particle of mass $m$ with the trajectory $\omega(\tau)$ in the presence of the potential $V$. Taking the limit formally in (0.3.8) yields the Feynman path integral (0.2.4)

$$N \times \left( \int\limits_{\Omega_x} \exp\left[ \frac{i}{\hbar} S_t(\dot{\omega}(\tau), \omega(\tau)) \right] D[\omega(\tau)] \right),$$
$$\qquad (0.3.11)$$

$$N = \textbf{constant},$$

where $\Omega_x \subset \mathbf{C}[0, t]$ is the set of all continuous trajectories with $\omega(0) = x$. As is well-known, there are difficulties in interpreting formal expression (0.3.11) as well defined integral of a functional $\exp\left[ \frac{i}{\hbar} S_t(\dot{\omega}(\tau), \omega(\tau)) \right] \in (\mathbf{C}^2[0, t])^*$ over subset $\Omega_x$:

(**i**) the "constant" is infinite (of infinite order),

(**ii**) the action (0.3.10) exists only if $\omega(\tau)$ and $V$ are sufficiently smooth, and

(**iii**) $D[\omega(\tau)] = \prod\limits_{0 \leq \tau \leq t} dx_\tau$ has no meaning as measure on $\mathbf{C}[0, t]$.

We wish to improve this difficultness using Colombeau approach. As a framework we use Colombeau algebra of generalized functionals, see Definition 0.3.1. We use the notation

$$\left( p_\varepsilon^t(x, y) dy \right)_\varepsilon = (4\pi t((\nu_\varepsilon)_\varepsilon))^{-d/2} \exp\left( -\frac{|x - y|^2}{4t((\nu_\varepsilon)_\varepsilon)} \right) dy, \qquad (0.3.12)$$

where $x$ are in $d$-dimensional Euclidean space $\mathbb{R}^d$, $\|x\| = \sqrt{x^2}$

$$dx = dx_1 \ldots dx_d, x^2 = \sum_{i=1}^{d} x_i^2, \qquad (0.3.13)$$

and $v_\varepsilon$, is a strictly positive Colombeau constant such that $v_\varepsilon = O(\varepsilon^r)$, if

$\varepsilon \to 0, r \in \mathbb{Z}$.

# Colombeau-Wiener measure $\left( \widetilde{D}_x^{\mathbb{R}} [\omega(\tau), v_\varepsilon] \right)_\varepsilon$.

We wish to construct now real-valued Colombeau-Wiener measure $\left( \widetilde{D}_x^{\mathbb{R}} [\omega(\tau), v_\varepsilon] \right)_\varepsilon$ on the space $\Omega_t$ of all trajectories. It is convenient to introduce the one-point compactification $\dot{\mathbb{R}}^d$ of $\mathbb{R}^d$ by a point $\infty$, and let $\Omega_t$ be the Cartesian product space:

$$\Omega_t = \prod_{0 \leq \tau \leq t} \dot{\mathbb{R}}^d. \qquad (0.3.14)$$

Thus an element $\omega(\tau)$ of $\Omega_t$ is an entirely arbitrary function from time $[0,t]$ to space $\mathbb{R}^d$. Let us give $\Omega_t$ the product topology, so that by the Tychonoff theorem, $\Omega_t$ is a compact Hausdorff space.

**Definition 0.3.1.** Let $\mathbf{C}(\Omega_t)$ be the set of all continuous functions on $\Omega_t$. Set $E(\mathbf{C}(\Omega_t)) = (\mathbf{C}(\Omega_t))^I, I = (0,1]$

$$E_{M,n}(\mathbf{C}(\Omega_t)) =$$

$$\{(\varphi_\varepsilon[\omega(\tau)])_\varepsilon \in E(\mathbf{C}(\Omega_t)) | (\exists n \in \mathbb{N}) \exists (\tau_1 < \ldots < \tau_n)(\exists p \in \mathbb{N})[((\varphi_\varepsilon[\omega(\tau)])_\varepsilon = \qquad (0.3.15)$$

$$= (\mathcal{F}_\varepsilon[\omega(\tau_1), \ldots, \omega(\tau_n)])_\varepsilon) \wedge (\sup_{\omega \in \Omega_t} |\varphi_\varepsilon[\omega(\tau)]| = O(\varepsilon^{-p}))]\},$$

where $p$ does not depend on $n$,

$$N_n(\mathbf{C}(\Omega_t)) =$$

$$\{(\varphi_{n,\varepsilon}[\omega(\tau)])_\varepsilon \in E(\mathbf{C}(\Omega_t)) | (\exists n \in \mathbb{N})\exists(\tau_1 < \ldots < \tau_n)(\forall q \in \mathbb{N})[((\varphi_{n,\varepsilon}[\omega(\tau)])_\varepsilon =$$

$$= (\mathcal{F}_\varepsilon[\omega(\tau_1),\ldots,\omega(\tau_n)])_\varepsilon) \wedge (\sup_{\omega \in \Omega_t}|\varphi_{n,\varepsilon}[\omega(\tau)]| = O(\varepsilon^q))]\}, \tag{0.3.16}$$

$$G_{M,n}(\mathbf{C}(\Omega_t)) = E_{M,n}(\mathbf{C}(\Omega_t))/N_n(\mathbf{C}(\Omega_t)).$$

Elements of $E_{M,n}(\mathbf{C}(\Omega_t))$ and $N_n(\mathbf{C}(\Omega_t))$ are called moderate, resp. negligible simple functionals.Note that $E_{M,n}(\mathbf{C}(\Omega_t))$ is an algebra with pointwise operations and in which $N_n(\mathbf{C}(\Omega_t))$ is ideal. Thus, $G_{M,n}(\mathbf{C}(\Omega_t))$ is an associative, commutative algebra.If $(\varphi_{n,\varepsilon}[\omega(\tau)])_\varepsilon \in E(\mathbf{C}(\Omega_t))$ is a representative of $\varphi_n \in G_{M,n}(\mathbf{C}(\Omega_t))$,we write $\varphi_n = \mathbf{cl}[(\varphi_{n,\varepsilon}[\omega(\tau)])_\varepsilon]$ or simply $\varphi_n = [(\varphi_{n,\varepsilon}[\omega(\tau)])_\varepsilon]$.

**Definition 0.3.2**.We denote

$$E_M(\mathbf{C}(\Omega_t)) = \bigoplus_{n=1}^{\infty} E_{M,n}(\mathbf{C}(\Omega_t)),$$

$$N(\mathbf{C}(\Omega_t)) = \bigoplus_{n=1}^{\infty} N_n(\mathbf{C}(\Omega_t)), \tag{0.3.17}$$

$$G_M(\mathbf{C}(\Omega_t)) = E_M(\mathbf{C}(\Omega_t))/N(\mathbf{C}(\Omega_t)).$$

Elements of $E_M(\mathbf{C}(\Omega_t))$ and $N(\mathbf{C}(\Omega_t))$ are called moderate,resp.negligible functionals. Note that $E_M(\mathbf{C}(\Omega_t))$ is an algebra with pointwise operations and in which $N(\mathbf{C}(\Omega_t))$is ideal. Thus, $G_M(\mathbf{C}(\Omega_t))$ is an associative, commutative algebra. Now we can define the Colombeau integral of Colombeau generalized functional $(\varphi_\varepsilon[\omega(\tau)])_\varepsilon$ i.e.,Colombeau integral of Colombeau $\widetilde{\mathbb{R}}$-valued generalized function $(\varphi_\varepsilon[\omega(\tau)])_{\varepsilon \in I}$ defined on $\Omega_t, I = (0,1]$.

Let's consider first Colombeau generalized functional $(\varphi_\varepsilon[\omega(\tau)])_\varepsilon : \Omega_t \to \widetilde{\mathbb{R}}$ of the simple form,i.e.

$$(\varphi_\varepsilon[\omega(\tau)])_\varepsilon = (\mathcal{F}_\varepsilon[\omega(\tau_1),\ldots,\omega(\tau_n)])_\varepsilon,$$

(0.3.18)

$$0 = \tau_1 < \ldots < \tau_n = t.$$

Let $x$ be a point in $\dot{\mathbb{R}}^d$ and let

$$(\mathbf{J}_\varepsilon^{\nu_\varepsilon}(\varphi_\varepsilon[\omega(\tau)]))_\varepsilon = \left(\int\limits_\Omega \varphi_\varepsilon[\omega(\tau)]\widetilde{\overline{D}}_x^{\mathbb{R}}[\omega(\tau),\nu_\varepsilon]\right)_\varepsilon =$$

$$\left(\int\ldots\int p_{\nu_\varepsilon}^{\tau_1}(x,x_1)dx_1 p_{\nu_\varepsilon}^{\tau_2-\tau_1}(x_1,x_2)dx_2\ldots p_{\nu_\varepsilon}^{\tau_n-\tau_{n-1}}(x_{n-1},x_n)dx_n \mathcal{F}_\varepsilon[\omega(\tau_1),\ldots,\omega(\tau_n)]\right)_\varepsilon$$

(0.3.19)

$$\int\ldots\int((p_{\nu_\varepsilon}^{\tau_1}(x)p_{\nu_\varepsilon}^{\tau_2-\tau_1}(x_1)\ldots p_{\nu_\varepsilon}^{\tau_n-\tau_{n-1}}(x_{n-1,x_n})\mathcal{F}_\varepsilon[\omega(\tau_1),\ldots,\omega(\tau_n)])_\varepsilon)dx_1\ldots dx_n$$

provided the integrals exist. Let $\mathbf{C}(\Omega_t) = \mathbf{C}^*[0,t]$ be the set of all continuous functions on $\Omega_t$ and let $\mathbf{C}^\bullet(\Omega_t)$ the set of all those of the form (0.3.15) where $\forall\varepsilon(\varepsilon \in (0,1])\mathcal{F}_\varepsilon$, is continuous. Any function in $\mathbf{C}(\Omega_t)$ may be uniformly approximated by functions in $\mathbf{C}^\bullet(\Omega_t)$, by the Stone-Weierstrass theorem. The mapping $\varphi_\varepsilon \to \mathbf{J}_\varepsilon^x(\varphi_\varepsilon[\omega(\tau)]) = \int\limits_{\Omega_x} \varphi_\varepsilon[\omega(\tau)]\widetilde{\overline{D}}_x^{\mathbb{R}}[\omega(\tau),\varepsilon], \ \varepsilon \in (0,1]$ is a linear functional on $\mathbf{C}^\bullet(\Omega_t)$, which is such that
$\mathbf{J}_\varepsilon^x(1) = N_\varepsilon$ and $\mathbf{J}_\varepsilon^x(\varphi_\varepsilon[\omega(\tau)]) \geq 0$ whenever $\varphi_\varepsilon[\omega(\tau)] \geq 0$. Consequently,

$$\forall\varepsilon(\varepsilon \in (0,1])[|\mathbf{J}_\varepsilon^x(\varphi_\varepsilon[\omega(\tau)])| \leq c_\varepsilon \sup_{\omega(\tau)}|\varphi_\varepsilon[\omega(\tau)]|],$$

(0.3.20)

i.e.

$$(|\mathbf{J}_\varepsilon^x(\varphi_\varepsilon[\omega(\tau)])|)_\varepsilon \leq (c_\varepsilon)_\varepsilon(\sup_{\omega(\tau)}|\varphi_\varepsilon[\omega(\tau)]|)_\varepsilon,$$

(0.3.21)

where $\mathbf{cl}[(c_\varepsilon)_\varepsilon] \in \widetilde{\mathbb{R}}$ and so the mapping $\varphi_\varepsilon \to \mathbf{J}_\varepsilon^x(\varphi_\varepsilon[\omega(\tau)])$ has a unique extension to a positive linear functional defined for all $\varphi_\varepsilon$ in $\mathbf{C}(\Omega_t)$. The Riesz representation theorem asserts that $\forall \varepsilon (\varepsilon \in (0,1])$ there is exist a regular measure, denoted $\widetilde{D}_x^{\mathbb{R}}[\omega(\tau), \nu_\varepsilon]$ or for short $\widetilde{D}_x[\omega(\tau), \nu_\varepsilon]$, such that

$$\mathbf{J}_\varepsilon^x(\varphi_\varepsilon[\omega(\tau)]) = \int\limits_{\Omega_t} \varphi_\varepsilon[\omega(\tau)] \widetilde{D}_x^{\mathbb{R}}[\omega(\tau), \nu_\varepsilon] \tag{0.3.22}$$

Therefore there is exist a regular Colombeau measure, denoted $\left(\widetilde{D}_x^{\mathbb{R}}[\omega(\tau), \nu_\varepsilon]\right)_\varepsilon$, or for short $\left(\widetilde{D}_x[\omega(\tau), \nu_\varepsilon]\right)_\varepsilon$ such that

$$(\mathbf{J}_\varepsilon^x(\varphi_\varepsilon[\omega(\tau)]))_\varepsilon =$$

$$\left(\int\limits_{\Omega_t} \varphi_\varepsilon[\omega(\tau)] \widetilde{D}_x^{\mathbb{R}}[\omega(\tau), \nu_\varepsilon]\right)_\varepsilon = \tag{0.3.23}$$

$$\int\limits_{\Omega_t} ((\varphi_\varepsilon[\omega(\tau)])_\varepsilon) \left(\widetilde{D}_x^{\mathbb{R}}[\omega(\tau), \nu_\varepsilon]\right)_\varepsilon.$$

Note that the elements of $\Omega_t$ were allowed to be arbitrary trajectories, even taking the value $\infty$. However it easy to see that the Colombeau measure $\left(\widetilde{D}_x^{\mathbb{R}}[\omega(\tau), \nu_\varepsilon]\right)_\varepsilon$ is concentrated on the continuous trajectories $\omega(\tau) \in \mathbf{C}[0,t]$ taking values in $\mathbb{R}^d$. This is based on the fact (see [41]) that if $\epsilon > 0, \varepsilon \in (0,1]$ and

$$(\rho_\varepsilon(\epsilon, \delta))_\varepsilon = \left(\sup_{t \le \delta(\varepsilon)} \int\limits_{\|y-x\| \ge \epsilon} p_\varepsilon^t(x,y) dy\right)_\varepsilon \tag{0.3.24}$$

then by Eq.(0.3.12) $(\rho_\varepsilon(\epsilon, \delta))_\varepsilon = o(\delta(\varepsilon))$. We assume now that the function $V$ is regular or satisfy Nelson's assumption (see [41] ). Then for almost every $\omega \in \Omega_t, V(\omega(\tau))$ is a continuous function of $\tau$ for $0 \le \tau < \infty$ and so we may form the Colombeau-Wiener integral:

$$(\mathbf{J}_\varepsilon^{x,t}[V(\omega(\tau)),\psi_\varepsilon(\omega(\tau))])_\varepsilon =$$

$$\left(\int_{\Omega_t}\exp\left[-i\int_0^t V(\omega(\tau))d\tau\right]\psi_\varepsilon(\omega(t))\widetilde{D}_x^{\mathbb{R}}[\omega(\tau),\varepsilon]\right)_\varepsilon = \tag{0.3.25}$$

$$\int_{\Omega_t}\exp\left[-i\int_0^t V(\omega(\tau))d\tau\right]((\psi_\varepsilon(\omega(t)))_\varepsilon)\left(\widetilde{D}_x^{\mathbb{R}}[\omega(\tau),\varepsilon]\right)_\varepsilon$$

for any $(\psi_\varepsilon(\omega(\tau)))_\varepsilon \in G(\mathbb{R}^d)$ and any $t \geq 0$. In fact, the integrand is defined for almost every $\omega \in \Omega_t$, and it is bounded in absolute value by $(|\psi_\varepsilon(\omega(t))|)_\varepsilon$. But this is Colombeau-integrable, for by definition and (0.3.20) one obtain

$$\left|(\mathbf{J}_\varepsilon^{x,t}[V(\omega(\tau)),\psi_\varepsilon(\omega(\tau))])_\varepsilon\right| \leq (\mathbf{J}_\varepsilon^x[\psi_\varepsilon(\omega(\tau))])_\varepsilon \leq ((c_\varepsilon)_\varepsilon)(\sup_{\omega(t)}|\psi_\varepsilon(\omega(t))|)_\varepsilon. \tag{0.3.26}$$

Since $V(\omega(\tau))$ is continuous for almost every $\omega(\tau)$,

$$\int_0^t V(\omega(\tau))d\tau = \lim_{n\to\infty}\left(\frac{t}{n}\sum_{j=1}^n V\left(\omega\left(j\frac{t}{n}\right)\right)\right) \tag{0.3.27}$$

for almost every $\omega$. By the Lebesgue dominated convergence theorem, this implies that

$$(\mathbf{J}_\varepsilon^{x,t}[V(\omega(\tau)),\psi_\varepsilon(\omega(\tau))])_\varepsilon =$$

$$\left(\lim_{n\to\infty}\int_\Omega\exp\left[-i\sum_{j=1}^n\frac{t}{n}V\left(\omega\left(j\frac{t}{n}\right)\right)\right]\psi_\varepsilon(\omega(t))\widetilde{D}_x^{\mathbb{R}}[\omega(\tau),\varepsilon]\right)_\varepsilon \tag{0.3.28}$$

for all $x \in \mathbb{R}^d \backslash F_V$, where $F_V$ is a closed set in $\mathbb{R}^d$ which depend only on $V$. By definition of the Colombeau-Wiener integral (0.3.25), this means that for all $x \notin F_V$

$$\left( \left( U_{\nu_\varepsilon, V}^t \right)_\varepsilon \right)(\psi_\varepsilon(x))_\varepsilon = \left( \lim_{n \to \infty} \left( \Xi_{\nu_\varepsilon}^{\frac{t}{n}} M_V^{\frac{t}{n}} \right)^n \psi_\varepsilon(x) \right)_\varepsilon, \qquad (0.3.29)$$

where

$$\left( U_{\nu_\varepsilon, V}^t \right)_\varepsilon = \exp[t(\nu_\varepsilon \Delta - iV)], \qquad (0.3.30)$$

$$\Xi_{\nu_\varepsilon}^t = \exp(t\nu_\varepsilon \Delta).$$

Explicitly,

$$(\Xi_{\nu_\varepsilon}^t \varphi_\varepsilon(x))_\varepsilon = (2\pi t((\nu_\varepsilon)_\varepsilon))^{-d/2} \int_{\mathbb{R}^d} \exp\left( -\frac{|x-y|^2}{2t((\nu_\varepsilon)_\varepsilon)} \right) \varphi(y) dy, \qquad (0.3.31)$$

A set $F$ of capacity $0$ has Lebesgue measure $0$, and so if $\psi_{n,\varepsilon} = \left( \Xi_{\nu_\varepsilon}^{\frac{t}{n}} M_V^{\frac{t}{n}} \right)^n \psi_\varepsilon$ then $\forall \varepsilon \in (0,1], \psi_{n,\varepsilon}$ converges almost everywhere in $\mathbb{R}^d$ to $U_{\nu_\varepsilon, V}^t \psi_\varepsilon$. By (0.3.31)

$$\left( \left| \Xi_{\nu_\varepsilon}^{\frac{t}{n}} \varphi_\varepsilon(x) \right| \right)_\varepsilon \leq \left( \Xi_{\nu_\varepsilon}^{\frac{t}{n}} |\varphi_\varepsilon(x)| \right)_\varepsilon. \qquad (0.3.32)$$

By definition $M_V^{\frac{t}{n}}$ is multiplication by a function of modulus $1$. Consequently, $\forall \varepsilon \in (0,1] : \left| \psi_{n,\varepsilon} \right| \leq \Xi_{\nu_\varepsilon}^t |\psi_{n,\varepsilon}|$, and so $\left| \psi_{n,\varepsilon} - \psi_{m,\varepsilon} \right|^2 \leq 4\left( \Xi_{\nu_\varepsilon}^t \psi_\varepsilon \right)^2$, which is in $\mathcal{L}_2(\mathbb{R}^n)$ Therefore, by the Lebesgue dominated convergence theorem, $\forall \varepsilon \in (0,1]$, $\psi_{n,\varepsilon}$ is a Cauchy sequence in $\mathcal{L}_2(\mathbb{R}^n)$, and so

$$\left(\left(U^{t}_{\nu_{\varepsilon},V}\right)_{\varepsilon}\right)(\psi_{\varepsilon})_{\varepsilon} = \left(\mathcal{L}_{2}\text{-}\lim_{n\to\infty}\left(\Xi^{\frac{t}{n}}_{\nu_{\varepsilon}}M^{\frac{t}{n}}_{V}\right)^{n}\psi_{\varepsilon}\right)_{\varepsilon}, \qquad (0.3.33)$$

where the limit now being taken in $\mathcal{L}_{2}(\mathbb{R}^{n})$.

## I.0.4. Regularized Colombeau-Feynman Path Integlal $(\mathbf{J}^{x,t}_{\varepsilon,V}[\mathcal{L}_{\varepsilon}(\dot{\omega}_{\varepsilon}(\tau),\omega(\tau)),\psi_{0}(\omega(t))])_{\varepsilon}$ based on Colombeau-Wiener measure $\left(\widetilde{D}^{\mathbb{R}}_{x}[\omega(\tau),\nu_{\varepsilon}]\right)_{\varepsilon}$.

Let us consider now the regularised Schrödinger equation

$$\frac{\partial\psi_{\varepsilon}(t,x)}{\partial t} = \frac{1}{2}(i+\nu_{\varepsilon})\Delta\psi_{\varepsilon}(t,x) - iV(x)\psi_{\varepsilon}(t,x),$$

$$((\psi_{\varepsilon}(0,x)))_{\varepsilon} = \psi_{0}(x), \qquad (0.4.1)$$

where $\nu_{\varepsilon}\to 0$ if $\varepsilon\to 0$. If the operator $H = -\Delta/2 + V(x)$ is:
(**i**) self-adjoint and
(**ii**) bounded from below, by the spectral theorem one obtain

$$\exp(-itH) = \lim_{\varepsilon\to 0}\exp[-(i+\nu_{\varepsilon})tH] \qquad (0.4.2)$$

strongly for all positive $t$. Therefore

$$\left(U_{\nu_\varepsilon,V}^t\right)_\varepsilon = \exp\{t[(i+(\nu_\varepsilon)_\varepsilon)\Delta - iV]\},$$

$$\Xi_{\nu_\varepsilon}^t = \exp(t(i+(\nu_\varepsilon)_\varepsilon)\Delta).$$

$$(0.4.3)$$

Explicitly,

$$\left(\left(\Xi_{\nu_\varepsilon}^t \varphi_\varepsilon(x)\right)_\varepsilon\right)_\varepsilon = \left(2\pi t(i + (\nu_\varepsilon)_\varepsilon)\right)^{-d/2} \int_{\mathbb{R}^d} \exp\left(-\frac{|x-y|^2}{2t(i+(\nu_\varepsilon)_\varepsilon)}\right)\varphi(y)dy =$$

$$\left(2\pi t(i+(\nu_\varepsilon)_\varepsilon)\right)^{-d/2} \int_{\mathbb{R}^d} \exp\left(-\frac{|x-y|^2((\nu_\varepsilon)_\varepsilon - i)}{2t((\nu_\varepsilon)_\varepsilon + i)((\nu_\varepsilon)_\varepsilon - i)}\right)\varphi(y)dy =$$

$$\left(2\pi t(i+(\nu_\varepsilon)_\varepsilon)\right)^{-d/2} \int_{\mathbb{R}^d} \exp\left(-\frac{|x-y|^2((\nu_\varepsilon)_\varepsilon - i)}{2t((\nu_\varepsilon^2)_\varepsilon + 1)}\right)\varphi(y)dy =$$

$$\left(2\pi t(i+(\nu_\varepsilon)_\varepsilon)\right)^{-d/2} \times$$

$$\int_{\mathbb{R}^d} \exp\left(-\frac{((\nu_\varepsilon)_\varepsilon)|x-y|^2}{2t((\nu_\varepsilon^2)_\varepsilon + 1)}\right)\exp\left(\frac{i|x-y|^2}{2t((\nu_\varepsilon^2)_\varepsilon + 1)}\right)\varphi(y)dy =$$

$$\left(2\pi t(i+(\nu_\varepsilon)_\varepsilon)\right)^{-d/2} \times \left[\frac{((\nu_\varepsilon)_\varepsilon)}{(\nu_\varepsilon^2)_\varepsilon + 1}\right]^{-d/2} \times \left[\frac{((\nu_\varepsilon)_\varepsilon)}{(\nu_\varepsilon^2)_\varepsilon + 1}\right]^{d/2} \times \qquad (0.4.4)$$

$$\int_{\mathbb{R}^d} \exp\left(-\frac{((\nu_\varepsilon)_\varepsilon)|x-y|^2}{2t((\nu_\varepsilon^2)_\varepsilon + 1)}\right)\exp\left(\frac{i|x-y|^2}{2t((\nu_\varepsilon^2)_\varepsilon + 1)}\right)\varphi(y)dy =$$

$$\left(i+(\nu_\varepsilon)_\varepsilon\right)^{-d/2} \times \left[\frac{((\nu_\varepsilon)_\varepsilon)}{((\nu_\varepsilon^2)_\varepsilon + 1)}\right]^{-d/2} \times$$

$$\int_{\mathbb{R}^d} \exp\left(\frac{i|x-y|^2}{((\nu_\varepsilon^2)_\varepsilon + 1)}\right)\varphi(y)p(x,dy) =$$

$$\left[\frac{((\nu_\varepsilon)_\varepsilon)}{(i-(\nu_\varepsilon^2)_\varepsilon)}\right]^{-d/2} \int_{\mathbb{R}^d} \exp\left(\frac{i|x-y|^2}{((\nu_\varepsilon^2)_\varepsilon + 1)}\right)\varphi(y)p_{\nu_\varepsilon}(x,dy),$$

where

$$p_{\nu_\varepsilon}(x, dy) = \left[ \frac{((\nu_\varepsilon)_\varepsilon)}{2\pi t((\nu_\varepsilon^2)_\varepsilon + 1)} \right]^{d/2} \exp\left( -\frac{((\nu_\varepsilon)_\varepsilon)|x - y|^2}{2t((\nu_\varepsilon^2)_\varepsilon + 1)} \right). \qquad (0.4.5)$$

**Definition 0.4.1.** We define Colombeau-Feynman Path Integlal $(\mathbf{J}_{\nu_\varepsilon, V}^{x,t}[\cdot])_{\varepsilon,(0,1]}$ based on real-valued Colombeau-Wiener measure $\left( \widetilde{D}_x^{\mathbb{R}}[\omega(\tau), \nu_\varepsilon] \right)_\varepsilon$ by formula

$$(\mathbf{J}_{\nu_\varepsilon, V}^{x,t}[\mathcal{L}_\varepsilon(\dot\omega(\tau), \omega(\tau)), \psi_\varepsilon(\omega(\tau))])_{\varepsilon,(0,1]} =$$

$$\left( \lim_{n\to\infty} \left[ \frac{((\nu_\varepsilon)_\varepsilon)}{(i - (\nu_\varepsilon^2)_\varepsilon)} \right]^{-nd/2} \times \right.$$

$$\int_{\Omega_t} \exp\left\{ i\frac{t}{n} \sum_{j=1}^n \left[ \frac{m}{2}\left(\frac{t}{n}\right)^{-2}(\omega_{j,n}(t) - \omega_{j-1,n}(t))^2 \right] - i\int_0^t V(\omega(\tau))d\tau \right\} \times$$

$$\left. \psi_0(\omega(t)) \widetilde{D}_x^{\mathbb{R}}[\omega(\tau), \nu_\varepsilon] \right)_\varepsilon$$

$$(0.4.6)$$

or by equivalent complete Colombeau formula:

$$(\mathbf{J}_{\nu_\varepsilon, n_\varepsilon, V}^{x,t}[\mathcal{L}_\varepsilon(\dot\omega(\tau), \omega(\tau)), \psi_0(\omega(t))])_{\varepsilon,(0,1]} =$$

$$\left( \left[ \frac{((\nu_\varepsilon^{n_\varepsilon})_\varepsilon)}{((i - (\nu_\varepsilon^2))^{n_\varepsilon})_\varepsilon} \right]^{-d/2} \times \right.$$

$$\int_{\Omega_t} \exp\left\{ i\frac{t}{n_\varepsilon} \sum_{j=1}^{n_\varepsilon} \left[ \frac{m}{2}\left(\frac{t}{n_\varepsilon}\right)^{-2}(\omega_{j,n_\varepsilon}(t) - \omega_{j-1,n_\varepsilon}(t))^2 \right] - i\int_0^t V(\omega(\tau))d\tau \right\} \times$$

$$(0.4.7)$$

$$\left. \psi_\varepsilon(\omega(t)) \widetilde{D}_x^{\mathbb{R}}[\omega(\tau), \nu_\varepsilon] \right)_\varepsilon,$$

where $n_\varepsilon \to \infty$ if $\varepsilon \to 0$.

We wish to expand now Colombeau-Wiener integral (0.3.25) on functionals given by

Eq.(0.3.10) using smoothing approach based on the mollification of the continuous trajectories $\omega(\tau) \in \mathbf{C}[0,t]$. We use a net of mollifiers $\phi_{n,\varepsilon}(\tau) \in \mathbf{S}(\mathbb{R}_+), \varepsilon \in (0,1]$ of the form

$$\phi_{n,\varepsilon}(\tau) = \frac{1}{\varepsilon} \phi_n \left( \frac{\tau}{\varepsilon} \right), \qquad (0.4.7)$$

where $\phi_n$ has the properties $\forall n \in \mathbb{N}$ :

(i) $\int_0^t \phi_n(\tau) d\tau = 1$,

(ii) $\phi_n(\tau) \geq 0$ for all $\tau \in [0,t]$,

(iii) $\phi_n \in C^\infty$.

**Definition 0.4.1**.[46].For $\varepsilon > 0$ the $\varepsilon$-mollification $\omega_{n,\varepsilon}(\tau)$ of a locally integrable function $\omega : [0,t] \to \mathbb{R}^d, \omega(\tau) = (\omega_1(\tau),\ldots,\omega_d(\tau))$ is the convolution on $\mathbb{R}$ such that

$$\omega_{n,\varepsilon}(\tau) = \phi_{n,\varepsilon} * \omega(\tau), \qquad (0.4.8)$$

that is,

$$\omega_{n,\varepsilon}(\tau) = \int_{\mathbb{R}} \phi_{n,\varepsilon}(\tau - z) \omega(z) dz. \qquad (0.4.9)$$

We remind that [46]: (i) for all $\varepsilon > 0$, the $\varepsilon$-mollification $\omega_{n,\varepsilon}(\tau)$ in $C^\infty(\mathbb{R}_+, \mathbb{R}^d)$,

(ii) If $\omega$ is continuous on $\mathbb{R}_+$, then $\omega_{n,\varepsilon}$ converges to $\omega$ uniformly on any compact sets for $\varepsilon \to 0$.

**Definition 0.4.2**. We define smoothed classical action by formula

$$S_{t,n,\varepsilon}(\dot\omega_\varepsilon(\tau),\omega(\tau)) = \int\limits_0^t \mathcal{L}_{n,\varepsilon}(\dot\omega_{n,\varepsilon}(\tau),\omega(\tau))d\tau =$$

$$\int\limits_0^t \left[ \frac{m}{2}\dot\omega_{n,\varepsilon}^2(\tau) - V(\omega(\tau)) \right]d\tau,$$

where $\dot\omega_{n,\varepsilon}(\tau) = \dfrac{d\omega_{n,\varepsilon}(\tau)}{d\tau} = \omega(\tau) * \dfrac{d\phi_{n,\varepsilon}(\tau)}{d\tau}$, and so we may form the complex Colombeau-Wiener integral:

$$\left(\mathbf{J}_{v_\varepsilon,n_\varepsilon}^{x,t}V\big[\mathcal{L}_\varepsilon(\dot\omega(\tau),\omega(\tau)),\psi_0(\omega(t))\big]\right)_\varepsilon =$$

$$\left(\left[\frac{((v_\varepsilon^{n_\varepsilon})_\varepsilon)}{((i-(v_\varepsilon^2))_\varepsilon^{n_\varepsilon})}\right]^{-d/2}\times\right.$$

$$\left. \int\limits_{\Omega_t}\exp\left[i\int\limits_0^t\left[\frac{m}{2}\dot\omega_{n_\varepsilon,\varepsilon}^2(\tau)-V(\omega(\tau))\right]d\tau\right]\psi_0(\omega(t))\widetilde{D}_x^{\mathbb{R}}[\omega(\tau),v_\varepsilon]\right)_\varepsilon =$$

$$\left[\frac{((v_\varepsilon^{n_\varepsilon})_\varepsilon)}{((i-(v_\varepsilon^2))_\varepsilon^{n_\varepsilon})}\right]^{-d/2}\times\int\limits_{\Omega_t}\exp\left[i\left(\int\limits_0^t\left[\frac{m}{2}\dot\omega_{n_\varepsilon,\varepsilon}^2(\tau)-V(\omega(\tau))\right]d\tau\right)_\varepsilon\right]\times$$

$$((\psi_0(\omega(t)))_\varepsilon)\left(\widetilde{D}_x^{\mathbb{R}}[\omega(\tau),v_\varepsilon]\right)_\varepsilon$$

and in more genelar case

$$\left(\mathbf{J}^{x,t}_{v_\varepsilon,n_\delta,V}\left[\mathcal{L}_{\varepsilon,\delta}(\dot{\omega}_{n,\delta}(\tau),\omega(\tau)),\psi_0(\omega(t))\right]\right)_{\varepsilon\in(0,1],\delta\in(0,1]}=$$

$$\left(\left[\frac{\left((v_\varepsilon^{n_\delta})_\varepsilon\right)}{\left((i-(v_\varepsilon^2))^{n_\delta}_\varepsilon\right)}\right]^{-d/2}\right)_{\varepsilon(0,1],\delta\in(0,1]}\times$$

$$\left(\int\limits_{\mathbf{\Omega}_t}\exp\left[i\int\limits_0^t\left[\frac{m}{2}\dot{\omega}^2_{n_\delta,\delta}(\tau)-V(\omega(\tau))\right]d\tau\right]\psi_0(\omega(t))\widetilde{D}^{\mathbb{R}}_x[\omega(\tau),v_\varepsilon]\right)_{\varepsilon(0,1],\delta\in(0,1]}= \qquad (0.4.12)$$

$$\int\limits_{\mathbf{\Omega}_t}\exp\left[i\left(\int\limits_0^t\left[\frac{m}{2}\dot{\omega}^2_{n_\delta,\delta}(\tau)-V(\omega(\tau))\right]d\tau\right)_{\delta\in(0,1]}\right]\times$$

$$\left((\psi_\varepsilon(\omega(t)))_{\varepsilon\in(0,1]}\right)\left(\widetilde{D}^{\mathbb{R}}_x[\omega(\tau),v_\varepsilon]\right)_{\varepsilon\in(0,1]}$$

**Definition 0.4.3.** We define *smoothed* Colombeau-Feynman Path Integlal based on real-valued Colombeau-Wiener measure $\left(\widetilde{D}^{\mathbb{R}}_x[\omega(\tau),v_\varepsilon]\right)_\varepsilon$ by formula

$$\int\limits_{\mathbf{\Omega}_t}\exp\left[i\int\limits_0^t\left[\frac{m}{2}\dot{\omega}^2(\tau)-V(\omega(\tau))\right]d\tau\right]((\psi_0(\omega(t)))_\varepsilon)\left(\widetilde{D}^{\mathbb{R}}_x[\omega(\tau),v_\varepsilon]\right)_\varepsilon \overset{\triangle}{=} \qquad (0.4.13)$$

$$\left(\mathbf{J}^{x,t}_{v_\varepsilon,n_\varepsilon,V}\left[\mathcal{L}_\varepsilon(\dot{\omega}_{n_\varepsilon,\varepsilon}(\tau),\omega(\tau)),\psi_0(\omega(t))\right]\right)_\varepsilon$$

and in more genelar case

$$\int_{\Omega_t} \exp\left[ i \int_0^t \left[ \frac{m}{2} \dot{\omega}^2(\tau) - V(\omega(\tau)) \right] d\tau \right] \times$$

$$\left( (\psi_0(\omega(t)))_{\varepsilon \in (0,1]} \right) \left( \widetilde{D}_x^{\mathbb{R}} [\omega(\tau), \nu_\varepsilon] \right)_{\varepsilon \in (0,1], \delta \in (0,1]} \triangleq \tag{0.4.14}$$

$$(\mathbf{J}_{\nu_\varepsilon, n_\delta, V}^{x,t} [\mathcal{L}_\varepsilon(\dot{\omega}_{n_\delta, \delta}(\tau), \omega(\tau)), \psi_0(\omega(t))])_{\varepsilon \in (0,1], \delta \in (0,1]}.$$

Since $\forall \varepsilon \in (0,1]$

$$\mathcal{L}_\varepsilon(\dot{\omega}_{n_\varepsilon, \varepsilon}(\tau), \omega(\tau)) = \frac{m}{2} \dot{\omega}_{n_\varepsilon, \varepsilon}^2(\tau) - V(\omega(\tau)) \tag{0.4.15}$$

is continuous for almost every $\omega(\tau)$, one obtain

$$\int_0^t \mathcal{L}_\varepsilon(\dot{\omega}_{n_\varepsilon, \varepsilon}(\tau), \omega(\tau)) d\tau = \lim_{n_\varepsilon \to \infty} \left( \frac{t}{n_\varepsilon} \sum_{j=1}^{n_\varepsilon} \left[ \frac{m}{2} \dot{\omega}_{n_\varepsilon, \varepsilon}^2 \left( j \frac{t}{n_\varepsilon} \right) - V\left( \omega\left( j \frac{t}{n_\varepsilon} \right) \right) \right] \right) \tag{0.4.16}$$

for almost every $\omega \in \Omega_t$. By the Lebesgue dominated convergence theorem, this implies that

$$(\mathbf{J}_{\nu_\varepsilon, n_\varepsilon, V}^{x,t} [\mathcal{L}_\varepsilon(\dot{\omega}_{n_\varepsilon, \varepsilon}(\tau), \omega(\tau)), \psi_0(\omega(t))])_\varepsilon =$$

$$\left[ \frac{((\nu_\varepsilon^{n_\varepsilon})_\varepsilon)}{((i - (\nu_\varepsilon^2))^{n_\varepsilon}_\varepsilon)} \right]^{-d/2} \tag{0.4.17}$$

$$\left( \int_{\Omega_t} \exp\left[ i \frac{t}{n} \sum_{j=1}^n \left[ \frac{m}{2} \dot{\omega}_{n_\varepsilon, \varepsilon}^2 \left( j \frac{t}{n} \right) - V\left( \omega\left( j \frac{t}{n} \right) \right) \right] \right] \psi_\varepsilon(\omega(t)) \widetilde{D}_x^{\mathbb{R}} [\omega(\tau), \nu_\varepsilon] \right)_\varepsilon$$

for all $x \in \mathbb{R}^d \backslash F_V$, where $F_V$ is a closed set in $\mathbb{R}^d$ which depend only on $V$. By definition of the Colombeau-Wiener integral (0.3.25), this means that for all $x \notin F_V$

$$\left( \left( U^t_{v_\varepsilon, n_\varepsilon, V} \psi_0(x) \right)_\varepsilon \right) = \left( \left( \Xi^{\frac{t}{n_\varepsilon}}_{v_\varepsilon, n_\varepsilon} M_V^{\frac{t}{n_\varepsilon}} \right)^{n_\varepsilon} \psi_0(x) \right)_\varepsilon. \qquad (0.4.18)$$

Here

$$\left( U^t_{v_\varepsilon, n_\varepsilon, V} \right)_\varepsilon = \exp\{t[(i + (v_\varepsilon)_\varepsilon)(\Delta[\phi_{n_\varepsilon, \varepsilon}])_\varepsilon - iV]\},$$

$$\Xi^t_{v_\varepsilon, n_\varepsilon} = \exp(t(i + v_\varepsilon)\Delta[\phi_{n_\varepsilon, \varepsilon}]), \qquad (0.4.19)$$

where we denote by $\Delta[\phi_{n_\varepsilon, \varepsilon}]$ smoothed Laplase operator, which defined by obvious way.

## I.0.5. Colombeau-Feynman Path Integlal $(\mathbf{J}^{x,t}_{v_\varepsilon, n_\varepsilon, V_\varepsilon}[\mathcal{L}_\varepsilon(\dot{\omega}_\varepsilon(\tau), \omega(\tau)), \psi_\varepsilon(\omega(t))])_\varepsilon$ based on Colombeau-Wiener measure $\left( \widetilde{D}_x[\omega(\tau), v_\varepsilon] \right)_\varepsilon$.

In our subsequent consideration, the potential $V(t, x)$ does not satisfy Kato conditions. We remove it by using cut-off or by another type of the regularization. As a framework we use Colombeau's algebra of generalized functions, i.e. we change potential $V(t, x)$ by Colombeau generalized function $(V_\varepsilon(t, x))_{\varepsilon \in (0,1]}$ such that $V_0(t, x) = V(t, x)$. Therefore classical action is

$$(S_{t, \varepsilon}(\dot{\omega}(\tau), \omega(\tau)))_\varepsilon = \left( \int_0^t \left[ \frac{m}{2} \dot{\omega}^2(\tau) - V_\varepsilon(\omega(\tau)) \right] d\tau \right)_\varepsilon \qquad (0.5.1)$$

We assume that $\forall \varepsilon \in (0,1]$ the potential $V_\varepsilon(t,x)$ is satisfy Kato conditions. Therefore the operator $V_\varepsilon$ of multiplication by the function $V_\varepsilon(t,x)$, on the domain $\mathbf{D}(V_\varepsilon)$ of all $\psi_{0,\varepsilon}$ in $\mathcal{L}_2$ such that $\forall \varepsilon \in (0,1]: V_\varepsilon \psi_{0,\varepsilon}$ is also in $\mathcal{L}_2$, operator $V_\varepsilon$ is self-adjoint, and if

$$\left( M_{V_\varepsilon}^t \right)_\varepsilon = \left( \exp\{-itV_\varepsilon\} \right)_\varepsilon =$$

$$\exp\{-it(V_\varepsilon)_\varepsilon\}, \tag{0.5.2}$$

then $(\psi_\varepsilon(t,x))_\varepsilon = \left( M_{V_\varepsilon}^t \psi_{0,\varepsilon}(x) \right)_\varepsilon$ is the Colombeau solution of Eq.(0.2.21) with $1/2m$ replaced by $0$. If we let now

$$(U_{m,V_\varepsilon}^t)_\varepsilon = \exp\left[ it\left( \frac{1}{2m}\Delta - (V_\varepsilon)_\varepsilon \right) \right] \tag{0.5.3}$$

then Trotter's theorem asserts that for all $(\psi_{0,\varepsilon})_\varepsilon$ such that $\psi_{0,\varepsilon} \in \mathcal{L}_2$ :

$$((U_{m,V_\varepsilon}^t)_\varepsilon)(\psi_{0,\varepsilon})_\varepsilon = (U_{m,V_\varepsilon}^t \psi_{0,\varepsilon})_\varepsilon = \left( \lim_{n\to\infty} \left( \Xi^{\frac{t}{m}} M_{V_\varepsilon}^{\frac{t}{n}} \right)^n \psi_{0,\varepsilon} \right)_\varepsilon . \tag{0.5.4}$$

Using (0.5.2)-(0.5.4) we find that

$$((U_{m,V_\varepsilon}^t)_\varepsilon)(\psi_{0,\varepsilon})_\varepsilon =$$

$$\left( \lim_{n\to\infty} \left( \frac{2\pi it}{nm} \right)^{-\frac{1}{2}dn} \int \ldots \int \exp[iS_\varepsilon(x_0,,\ldots,x_n)]\psi_{0,\varepsilon}(x_n)dx_1\ldots dx_n \right)_\varepsilon . \tag{0.5.5}$$

Here we set $x_0 = x$ and

$$S_\varepsilon(x_0,\,,\ldots,x_n) = \left(\frac{t}{n}\right)\sum_{j=1}^{n}\left[\frac{m}{2}\frac{(x_i - x_{j-1})^2}{(t/n)^2} - V_\varepsilon(x_i)\right]. \qquad (0.5.6)$$

$x \in \bigcup\limits_{\varepsilon \in (0,1]} \mathbb{R}^d \backslash F_{V_\varepsilon}$, where $F_{V_\varepsilon}$ is a closed set in $\mathbb{R}^d$ which depend on $V_\varepsilon$.

We assume now that $\forall \varepsilon \in (0,1]$, the function $V_\varepsilon(x)$ is regular or satisfy Nelson's assumption (see [41]). Then for almost every $\omega \in \Omega_t$ and $\forall \varepsilon \in (0,1], V_\varepsilon(\omega(\tau))$ is a continuous function of $\tau$ for $0 \leq \tau < \infty$ and so (as above in section I.0.4) we may form the Colombeau-Wiener integral:

$$(\mathbf{J}_{V_\varepsilon}^{x,t}[V_\varepsilon(\omega(\tau)), \psi_\varepsilon(\omega(\tau))])_\varepsilon =$$

$$\left(\int\limits_{\Omega_t} \exp\left[-i\int\limits_0^t V_\varepsilon(\omega(\tau))d\tau\right]\psi_\varepsilon(\omega(t))\widetilde{D}_x^{\mathbb{R}}[\omega(\tau), \nu_\varepsilon]\right)_\varepsilon =$$

$$\int\limits_{\Omega_t} \exp\left[-i\left(\int\limits_0^t V_\varepsilon(\omega(\tau))d\tau\right)_\varepsilon\right]((\psi_\varepsilon(\omega(t)))_\varepsilon)\left(\widetilde{D}_x^{\mathbb{R}}[\omega(\tau), \nu_\varepsilon]\right)_\varepsilon$$

$$(0.5.7)$$

for any $(\psi_\varepsilon(\omega(\tau)))_\varepsilon \in G(\mathbb{R}^d)$ and any $t \geq 0$. In fact, the integrand is defined $\forall \varepsilon \in (0,1]$ for almost every $\omega \in \Omega_t$, and it is bounded in absolute value by $(|\psi_\varepsilon(\omega(t))|)_\varepsilon$. But this is Colombeau-integrable, for by definition and (0.5.7) one obtain

$$\left|(\mathbf{J}_{V_\varepsilon}^{x,t}[V_\varepsilon(\omega(\tau)), \psi_\varepsilon(\omega(\tau))])_\varepsilon\right| \leq (\mathbf{J}_{V_\varepsilon}^{x,t}[\psi_\varepsilon(\omega(\tau))])_\varepsilon \leq ((c_\varepsilon)_\varepsilon)(\sup\nolimits_{\omega(\tau)}|\psi_\varepsilon(\omega(t))|)_\varepsilon. \qquad (0.5.8)$$

Since $\forall \varepsilon \in (0,1], V_\varepsilon(\omega(\tau))$ is continuous for almost every $\omega(\tau)$,

$$\int\limits_{0}^{t} V_\varepsilon(\omega(\tau))d\tau \;=\; \lim_{n\to\infty}\left(\frac{t}{n}\sum_{j=1}^{n}V_\varepsilon\Big(\omega\Big(j\frac{t}{n}\Big)\Big)\right) \tag{0.5.9}$$

for almost every $\omega$. By the Lebesgue dominated convergence theorem, this implies that

$$(\mathbf{J}_\varepsilon^{x,t}[V(\omega(\tau)),\psi_\varepsilon(\omega(\tau))])_\varepsilon =$$

$$\left(\lim_{n\to\infty}\int_{\Omega}\exp\left[-i\sum_{j=1}^{n}\frac{t}{n}V_\varepsilon\Big(\omega\Big(j\frac{t}{n}\Big)\Big)\right]\psi_\varepsilon(\omega(t))\widetilde{D}_x^{\mathbb{R}}[\omega(\tau),\varepsilon]\right)_\varepsilon \tag{0.5.10}$$

for all $x \in \mathbb{R}^d\setminus\bigcup_{\varepsilon\in(0,1]}F_{V_\varepsilon}$, where $F_{V_\varepsilon}$ is a closed set in $\mathbb{R}^d$ which depend only on $V_\varepsilon$. By definition of the Colombeau-Wiener integral (0.3.25), this means that for all $x \notin \bigcup_{\varepsilon\in(0,1]}F_{V_\varepsilon}$

$$\left(\Big(U_{\nu_\varepsilon,V_\varepsilon}^{t}\Big)_\varepsilon\right)(\psi_\varepsilon(x))_\varepsilon = \left(\lim_{n\to\infty}\Big(\Xi_{\nu_\varepsilon}^{\frac{t}{n}}M_{V_\varepsilon}^{\frac{t}{n}}\Big)^n\psi_\varepsilon(x)\right)_\varepsilon, \tag{0.5.11}$$

where

$$\left(U_{\nu_\varepsilon,V_\varepsilon}^{t}\right)_\varepsilon = \exp[t(\nu_\varepsilon\Delta - iV_\varepsilon)], \tag{0.5.12}$$

$$\Xi_{\nu_\varepsilon}^{t} = \exp(t\nu_\varepsilon\Delta).$$

Explicitly,

$$(\Xi^t_{V_\varepsilon}\varphi_\varepsilon(x))_\varepsilon = (2\pi t((\nu_\varepsilon)_\varepsilon))^{-d/2}\int_{\mathbb{R}^d}\exp\left(-\frac{|x-y|^2}{2t((\nu_\varepsilon)_\varepsilon)}\right)\varphi(y)dy, \qquad (0.5.13)$$

A set $F_{V_\varepsilon}$ of capacity $0$ has Lebesgue measure $0$, and so if $\psi_{n,\varepsilon} = \left(\Xi^{\frac{t}{n}}_{V_\varepsilon}M^{\frac{t}{n}}_V\right)^n\psi_\varepsilon$ then $\forall\varepsilon\in(0,1],\psi_{n,\varepsilon}$ converges almost everywhere in $\mathbb{R}^d$ to $U^t_{\nu_\varepsilon,V_\varepsilon}\psi_\varepsilon$. By (0.5.13)

$$\left(\left|\Xi^{\frac{t}{n}}_{V_\varepsilon}\varphi_\varepsilon(x)\right|\right)_\varepsilon \leq \left(\Xi^{\frac{t}{n}}_{V_\varepsilon}|\varphi_\varepsilon(x)|\right)_\varepsilon. \qquad (0.5.14)$$

By definition $M^{\frac{t}{n}}_V$ is multiplication by a function of modulus $1$. Consequently, $\forall\varepsilon\in(0,1]:\left|\psi_{n,\varepsilon}\right|\leq\Xi^t_{V_\varepsilon}|\psi_{n,\varepsilon}|$, and so $\left|\psi_{n,\varepsilon}-\psi_{m,\varepsilon}\right|^2\leq 4\left(\Xi^t_{V_\varepsilon}\psi_\varepsilon\right)^2$, which is in $\mathcal{L}_2(\mathbb{R}^n)$ Therefore, by the Lebesgue dominated convergence theorem, $\forall\varepsilon\in(0,1]$, $\psi_{n,\varepsilon}$ is a Cauchy sequence in $\mathcal{L}_2(\mathbb{R}^n)$, and so

$$\left(\left(U^t_{\nu_\varepsilon,V_\varepsilon}\right)_\varepsilon\right)(\psi_\varepsilon)_\varepsilon = \left(\mathcal{L}_2\text{-}\lim_{n\to\infty}\left(\Xi^{\frac{t}{n}}_{V_\varepsilon}M^{\frac{t}{n}}_{V_\varepsilon}\right)^n\psi_\varepsilon\right)_\varepsilon, \qquad (0.5.15)$$

where the limit now being taken in $\mathcal{L}_2(\mathbb{R}^n)$.

Let us consider now the regularised Colombeau-Schrödinger equation

$$\frac{\partial\psi_\varepsilon(t,x)}{\partial t} = \frac{1}{2}(i+\nu_\varepsilon)\Delta\psi_\varepsilon(t,x) - iV_\varepsilon(x)\psi_\varepsilon(t,x),$$
$$(0.5.16)$$
$$((\psi_\varepsilon(0,x)))_\varepsilon = (\psi_{0,\varepsilon}(x))_\varepsilon,$$

where $\nu_\varepsilon \to 0$ if $\varepsilon \to 0$. If the operator $H_\varepsilon = -\Delta/2 + V_\varepsilon(x)$ is:

(**i**) $\forall \varepsilon \in (0,1]$, self-adjoint and

(**ii**) $\forall \varepsilon \in (0,1]$, bounded from below, by the spectral theorem one obtain

$$\exp(-itH) = \lim_{\varepsilon \to 0} \exp[-(i + \nu_\varepsilon)tH] \qquad (0.5.17)$$

strongly for all positive $t$. Therefore

$$\left(U^t_{\nu_\varepsilon, V_\varepsilon}\right)_\varepsilon = \exp\{t[(i + (\nu_\varepsilon)_\varepsilon)\Delta - iV_\varepsilon]\},$$

$$\Xi^t_{\nu_\varepsilon} = \exp(t(i + (\nu_\varepsilon)_\varepsilon)\Delta). \qquad (0.5.18)$$

Explicitly,

$$\left(\left(\Xi_{\nu_\varepsilon}^t \varphi_\varepsilon(x)\right)_\varepsilon\right)_\varepsilon = \left(2\pi t(i+(\nu_\varepsilon)_\varepsilon)\right)^{-d/2} \int_{\mathbb{R}^d} \exp\left(-\frac{|x-y|^2}{2t(i+(\nu_\varepsilon)_\varepsilon)}\right)\varphi(y)dy =$$

$$\left(2\pi t(i+(\nu_\varepsilon)_\varepsilon)\right)^{-d/2} \int_{\mathbb{R}^d} \exp\left(-\frac{|x-y|^2((\nu_\varepsilon)_\varepsilon - i)}{2t((\nu_\varepsilon)_\varepsilon + i)((\nu_\varepsilon)_\varepsilon - i)}\right)\varphi(y)dy =$$

$$\left(2\pi t(i+(\nu_\varepsilon)_\varepsilon)\right)^{-d/2} \int_{\mathbb{R}^d} \exp\left(-\frac{|x-y|^2((\nu_\varepsilon)_\varepsilon - i)}{2t((\nu_\varepsilon^2)_\varepsilon + 1)}\right)\varphi(y)dy =$$

$$\left(2\pi t(i+(\nu_\varepsilon)_\varepsilon)\right)^{-d/2} \times$$

$$\int_{\mathbb{R}^d} \exp\left(-\frac{((\nu_\varepsilon)_\varepsilon)|x-y|^2}{2t((\nu_\varepsilon^2)_\varepsilon + 1)}\right)\exp\left(\frac{i|x-y|^2}{2t((\nu_\varepsilon^2)_\varepsilon + 1)}\right)\varphi(y)dy =$$

$$\left(2\pi t(i+(\nu_\varepsilon)_\varepsilon)\right)^{-d/2} \times \left[\frac{((\nu_\varepsilon)_\varepsilon)}{(\nu_\varepsilon^2)_\varepsilon + 1}\right]^{-d/2} \times \left[\frac{((\nu_\varepsilon)_\varepsilon)}{(\nu_\varepsilon^2)_\varepsilon + 1}\right]^{d/2} \times \quad (0.5.19)$$

$$\int_{\mathbb{R}^d} \exp\left(-\frac{((\nu_\varepsilon)_\varepsilon)|x-y|^2}{2t((\nu_\varepsilon^2)_\varepsilon + 1)}\right)\exp\left(\frac{i|x-y|^2}{2t((\nu_\varepsilon^2)_\varepsilon + 1)}\right)\varphi(y)dy =$$

$$(i+(\nu_\varepsilon)_\varepsilon)^{-d/2} \times \left[\frac{((\nu_\varepsilon)_\varepsilon)}{((\nu_\varepsilon^2)_\varepsilon + 1)}\right]^{-d/2} \times$$

$$\int_{\mathbb{R}^d} \exp\left(\frac{i|x-y|^2}{((\nu_\varepsilon^2)_\varepsilon + 1)}\right)\varphi(y)p(x,dy) =$$

$$\left[\frac{((\nu_\varepsilon)_\varepsilon)}{(i-(\nu_\varepsilon^2)_\varepsilon)}\right]^{-d/2} \int_{\mathbb{R}^d} \exp\left(\frac{i|x-y|^2}{((\nu_\varepsilon^2)_\varepsilon + 1)}\right)\varphi(y)p_{\nu_\varepsilon}(x,dy),$$

where

$$p_{v_\varepsilon}(x, dy) = \left[ \frac{((v_\varepsilon)_\varepsilon)}{2\pi t((v_\varepsilon^2)_\varepsilon + 1)} \right]^{d/2} \exp\left( -\frac{((v_\varepsilon)_\varepsilon)|x - y|^2}{2t((v_\varepsilon^2)_\varepsilon + 1)} \right). \qquad (0.5.20)$$

**Definition 0.5.1.** We define Colombeau-Feynman Path Integlal $(\mathbf{J}_{v_\varepsilon, V_\varepsilon}^{x,t}[\bullet])_{\varepsilon,(0,1]}$ based on real-valued Colombeau-Wiener measure $\left( \widetilde{D}_x^{\mathbb{R}}[\omega(\tau), v_\varepsilon] \right)_\varepsilon$ by formula

$$(\mathbf{J}_{v_\varepsilon, V_\varepsilon}^{x,t}[\mathcal{L}_\varepsilon(\dot{\omega}(\tau), \omega(\tau)), \psi_\varepsilon(\omega(t))])_{\varepsilon,(0,1]} =$$

$$\left( \lim_{n\to\infty} \left[ \frac{((v_\varepsilon)_\varepsilon)}{(i - (v_\varepsilon^2)_\varepsilon)} \right]^{-nd/2} \times \right.$$

$$\int_{\Omega_t} \exp\left\{ i\frac{t}{n} \sum_{j=1}^n \left[ \frac{m}{2} \left( \frac{t}{n} \right)^{-2} (\omega_{j,n}(t) - \omega_{j-1,n}(t))^2 \right] - i\int_0^t V_\varepsilon(\omega(\tau))d\tau \right\} \times \qquad (0.5.21)$$

$$\left. \psi_{0,\varepsilon}(\omega(t)) \widetilde{D}_x^{\mathbb{R}}[\omega(\tau), v_\varepsilon] \right)_\varepsilon$$

or by equivalent complete Colombeau formula:

$$\left(\mathbf{J}^{x,t}_{V_\varepsilon,n_\delta,V_\varepsilon}[\mathcal{L}_\varepsilon(\dot\omega(\tau),\omega(\tau)),\psi_0(\omega(t))]\right)_{\varepsilon\in(0,1],\delta\in(0,1]} =$$

$$\left(\left[\frac{((v^{n_\delta}_\varepsilon)_\varepsilon)}{((i-(v^2_\varepsilon))^{n_\delta})_\varepsilon}\right]^{-d/2} \times \right.$$

$$(0.5.22)$$

$$\int_{\Omega_t}\exp\left\{i\frac{t}{n_\varepsilon}\sum_{j=1}^{n_\delta}\left[\frac{m}{2}\left(\frac{t}{n_\delta}\right)^{-2}(\omega_{j,n_\delta}(t)-\omega_{j-1,n_\delta}(t))^2\right]-i\int_0^t V_\varepsilon(\omega(\tau))d\tau\right\} \times$$

$$\left.\psi_\varepsilon(\omega(t))\widetilde{D}^{\mathbb{R}}_x[\omega(\tau),v_\varepsilon]\right)_{\varepsilon,\delta},$$

where $n_\delta \to \infty$ if $\delta \to 0$.

We wish to expand now Colombeau-Wiener integral (0.6.22) on functionals given by Eq.(0.3.10) using smoothing approach based on the mollification of the continuous trajectories $\omega(\tau) \in \mathbf{C}[0,t]$. We use a net of mollifiers $\phi_{n,\delta}(\tau) \in \mathbf{S}(\mathbb{R}_+), \delta \in (0,1]$ of the form

$$\phi_{n,\delta}(\tau) = \frac{1}{\delta}\phi_n\left(\frac{\tau}{\delta}\right), \qquad (0.5.23)$$

where $\phi_n$ has the properties $\forall n \in \mathbb{N}$ :

**(i)** $\int_0^t \phi_n(\tau)d\tau = 1,$
**(ii)** $\phi_n(\tau) \geq 0$ for all $\tau \in [0,t],$
**(iii)** $\phi_n \in C^\infty.$
**Definition 0.5.2.**[46].For $\delta > 0$ the $\delta$-mollification $\omega_{n,\delta}(\tau)$ of a locally integrable function $\omega : [0,t] \to \mathbb{R}^d, \omega(\tau) = (\omega_1(\tau),\dots,\omega_d(\tau))$ is the convolution on $\mathbb{R}$ such that

$$\omega_{n,\delta}(\tau) = \phi_{n,\delta} * \omega(\tau), \qquad (0.5.24)$$

that is,

$$\omega_{n,\delta}(\tau) = \int\limits_{\mathbb{R}} \phi_{n,\delta}(\tau - z)\omega(z)dz. \qquad (0.5.25)$$

We remind that [46]: (i) for all $\delta > 0$, the $\delta$-mollification $\omega_{n,\delta}(\tau)$ in $C^{\infty}(\mathbb{R}_+, \mathbb{R}^d)$,
(ii) If $\omega$ is continuous on $\mathbb{R}_+$, then $\omega_{n,\delta}$ converges to $\omega$ uniformly on any compact sets for $\delta \to 0$.

**Definition 0.5.3.** We define smoothed classical action by formula

$$S_{t,n,\varepsilon,\delta}(\dot{\omega}_{\delta}(\tau), \omega(\tau)) = \int\limits_0^t \mathscr{L}_{n,\varepsilon,\delta}(\dot{\omega}_{n,\delta}(\tau), \omega(\tau))d\tau =$$

$$\qquad (0.5.26)$$

$$\int\limits_0^t \left[ \frac{m}{2}\dot{\omega}_{n,\delta}^2(\tau) - V_{\varepsilon}(\omega(\tau)) \right] d\tau,$$

where $\dot{\omega}_{n,\delta}(\tau) = \dfrac{d\omega_{n,\delta}(\tau)}{d\tau} = \omega(\tau) * \dfrac{d\phi_{n,\delta}(\tau)}{d\tau}$, and so we may form the complex

Colombeau-Wiener path integral $(\mathbf{J}_{V_{\varepsilon}, n_{\varepsilon}V_{\varepsilon}}^{x,t}[\bullet])_{\varepsilon \in (0,1], \delta \in (0,1]}$ :

$$(\mathbf{J}^{x,t}_{v_\varepsilon,n_\varepsilon V_\varepsilon}[\mathcal{L}_{\varepsilon,\delta}(\dot{\omega}(\tau),\omega(\tau)),\psi_0(\omega(t))])_{\varepsilon\in(0,1],\delta\in(0,1]} =$$

$$\left(\left[\frac{v_\varepsilon^{n_\delta}}{(i-(v_\varepsilon^2))^{n_\delta}}\right]^{-d/2}\times\right.$$

$$\left.\int\limits_{\Omega_t}\exp\left[i\int\limits_0^t\left[\frac{m}{2}\dot{\omega}^2_{n_\delta,\delta}(\tau)-V(\omega(\tau))\right]d\tau\right]\psi_\varepsilon(\omega(t))\widetilde{D}^{\mathbb{R}}_x[\omega(\tau),v_\varepsilon]\right)_{\varepsilon,\delta} = \qquad (0.5.27)$$

$$\left[\frac{(v_\varepsilon^{n_\delta})_{\varepsilon,\delta}}{((i-(v_\varepsilon^2))^{n_\delta})_{\varepsilon,\delta}}\right]^{-d/2}\times\int\limits_{\Omega_t}\exp\left[i\left(\int\limits_0^t\left[\frac{m}{2}\dot{\omega}^2_{n_\delta,\varepsilon}(\tau)-V_\varepsilon(\omega(\tau))\right]d\tau\right)_{\varepsilon,\delta}\right]\times$$

$$((\psi_\varepsilon(\omega(t)))_\varepsilon)\left(\widetilde{D}^{\mathbb{R}}_x[\omega(\tau),v_\varepsilon]\right)_\varepsilon$$

**Definition 0.5.4.** We define *smoothed* Colombeau-Feynman Path Integlal based on real-valued Colombeau-Wiener measure $\left(\widetilde{D}^{\mathbb{R}}_x[\omega(\tau),v_\varepsilon]\right)_\varepsilon$ by formula

$$\int\limits_{\Omega_t}\exp\left[i\int\limits_0^t\left[\frac{m}{2}\dot{\omega}^2(\tau)-V_\varepsilon(\omega(\tau))\right]d\tau\right]((\psi_\varepsilon(\omega(t)))_\varepsilon)\left(\widetilde{D}^{\mathbb{R}}_x[\omega(\tau),v_\varepsilon]\right)_\varepsilon\triangleq \qquad (0.5.28)$$

$$(\mathbf{J}^{x,t}_{v_\varepsilon,n_\varepsilon,V_\varepsilon}[\mathcal{L}_\varepsilon(\dot{\omega}_{n_\varepsilon,\varepsilon}(\tau),\omega(\tau)),\psi_0(\omega(t))])_\varepsilon$$

and in more genelar case by formula

$$\int\limits_{\Omega_t} \exp\left[ i \int\limits_0^t \left[ \frac{m}{2} \dot{\omega}^2(\tau) - V_\varepsilon(\omega(\tau)) \right] d\tau \right] \times$$

$$\left( (\psi_\varepsilon(\omega(t)))_{\varepsilon \in (0,1]} \right) \left( \widetilde{D}_x^{\mathbb{R}} [\omega(\tau), v_\varepsilon] \right)_{\varepsilon \in (0,1], \delta \in (0,1]} \triangleq \qquad (0.5.29)$$

$$(\mathbf{J}_{\nu_\varepsilon, n_\delta, V_\varepsilon}^{x,t} [\mathcal{L}_{\varepsilon,\delta}(\dot{\omega}_{n_\delta,\delta}(\tau), \omega(\tau)), \psi_\varepsilon(\omega(t))])_{\varepsilon \in (0,1], \delta \in (0,1]}.$$

Since $\forall \varepsilon, \delta \in (0,1]$

$$\mathcal{L}_{\varepsilon,\delta}(\dot{\omega}_{n_\delta,\delta}(\tau), \omega(\tau)) = \frac{m}{2} \dot{\omega}_{n_\delta,\delta}^2(\tau) - V_\varepsilon(\omega(\tau)) \qquad (0.5.30)$$

is continuous for almost every $\omega(\tau)$, one obtain

$$\int\limits_0^t \mathcal{L}_\varepsilon(\dot{\omega}_{n_\delta,\delta}(\tau), \omega(\tau)) d\tau =$$

$$\lim_{n_\delta \to \infty, n_\varepsilon \to \infty} \left( \frac{t}{n_\delta} \sum_{j=1}^{n_\delta} \left[ \frac{m}{2} \dot{\omega}_{n_\delta,\delta}^2 \left( j \frac{t}{n_\delta} \right) - V_\varepsilon \left( \omega \left( j \frac{t}{n_\delta} \right) \right) \right] \right) \qquad (0.5.31)$$

for almost every $\omega \in \Omega_t$. By the Lebesgue dominated convergence theorem, this implies that

$$(\mathbf{J}_{v_\varepsilon,n_\delta,V_\varepsilon}^{x,t}[\mathcal{L}_\varepsilon(\dot{\omega}_{n_\delta,\varepsilon}(\tau),\omega(\tau)),\psi_0(\omega(t))])_{\varepsilon\in(0,1],\delta\in(0,1]} =$$

$$\left(\left[\frac{v_\varepsilon^{n_\delta}}{((i-(v_\varepsilon^2))^{n_\delta})}\right]^{-d/2}\right)_{\varepsilon\in(0,1],\delta\in(0,1]} \times$$

$$\left(\int\limits_{\Omega_t}\exp\left[i\frac{t}{n_\delta}\sum_{j=1}^{n_\delta}\left[\frac{m}{2}\dot{\omega}_{n_\varepsilon,\varepsilon}^2\left(j\frac{t}{n_\delta}\right)-V_\varepsilon\left(\omega\left(j\frac{t}{n_\delta}\right)\right)\right]\right] \times \right.$$

$$\left. \psi_\varepsilon(\omega(t))\widetilde{D}_x^{\mathbb{R}}[\omega(\tau),v_\varepsilon]\right)_{\varepsilon\in(0,1],\delta\in(0,1]}$$

(0.5.32)

for all $x \in \mathbb{R}^d \setminus \bigcup_{\varepsilon\in(0,1]} F_{V_\varepsilon}$, where $F_{V_\varepsilon}$ is a closed set in $\mathbb{R}^d$ which depend only on $V_\varepsilon$. By definition of the Colombeau-Wiener integral (0.3.25), this means that for all $x \notin \bigcup_{\varepsilon\in(0,1]} F_{V_\varepsilon}$

$$\left(\left(U_{v_\varepsilon,n_\delta,V_\varepsilon}^t\psi_\varepsilon(x)\right)_{\varepsilon\in(0,1],\delta\in(0,1]}\right) = \left(\left(\Xi_{v_\varepsilon,n_\delta}^{\frac{t}{n_\delta}}M_{V_\varepsilon}^{\frac{t}{n_\delta}}\right)^{n_\delta}\psi_\varepsilon(x)\right)_{\varepsilon\in(0,1],\delta\in(0,1]}.$$

(0.5.33)

Here

$$\left(U_{v_\varepsilon,n_\delta,V_\varepsilon}^t\right)_{\varepsilon\in(0,1],\delta\in(0,1]} = (\exp\{t[(i+v_\varepsilon)\Delta[\phi_{n_\delta,\delta}]-iV_\varepsilon]\})_{\varepsilon\in(0,1],\delta\in(0,1]},$$

(0.5.34)

$$\Xi_{v_\varepsilon,n_\delta}^t = \exp(t(i+v_\varepsilon)\Delta[\phi_{n_\delta,\delta}]),$$

where we denote by $\Delta[\phi_{n_\delta,\delta}]$ smoothed Laplase operator, which defined by obvious way.

## I.0.6. Colombeau-Feynman Path Integlal based on complex Colombeau-Wiener measure $\widetilde{D}^{\mathbb{C}}[x(\tau),v_\varepsilon]$.

Let $B(\Omega)$ denote the class of all Borel sets of a topological space $\Omega$, i.e. $B(\Omega)$ is the $\sigma$-algebra of sets generated by all open sets. If $\Omega$ is locally compact we denote by $C_0(\Omega)$ the space of all continuous complex-valued functions on $\Omega$-vanishing at infinity.
Equipped with the uniform norm $\|f\| = \sup_{x\in\Omega}|f(x)|$ this space is a Banach space.

**Theorem 0.6.1.** (**Riesz-Markov theorem**) Let $\Omega$ be a locally compact space, then the set $\mathbf{M}(\Omega)$ of all finite complex regular Borel measures on $\Omega$ equipped with the norm

$$\|\mu\| = \sup_{f\in C_0(\Omega),\|f\|\leq 1}\left|\int_\Omega f(x)d\mu\right| \qquad (0.6.1)$$

is a Banach space, which coincides with the set of all continuous linear functionals on $C_0(\Omega)$.

Any complex $\sigma$-additive measure $\mu$ on $\mathbb{R}^d$ has a representation of the form

$$\mu(dy) = f(y)M(dy) \qquad (0.6.2)$$

with a positive measure $M(dy)$ and a bounded complex-valued function $f(y)$. Moreover, the measure $M(dy)$ in (0.6.2) is uniquely defined under additional assumption that $|f(y)| = 1$ for all $y$. If this condition holds, the positive measure $M(dy)$ is called the total variation measure of the complex measure $\mu(dy)$ and is denoted by $|\mu|(dy)$. In general, if (0.6.2) holds, then

$$\|\mu\| = \int |f(y)|M(dy) \qquad (0.6.3)$$

**Definition 0.6.1.** We say that a map $\mu^{\mathbb{C}} : \mathbb{R}^d \times B(\mathbb{R}^d) \to \mathbb{C}$ is a complex transition kernel, if for every $x$, the map $A \mapsto \mu^{\mathbb{C}}(x,A)$ is a finite complex measure on $\mathbb{R}^d$, and for every $A \in \mathbf{B}(\mathbb{R}^d)$ the map $x \to \mu^{\mathbb{C}}(x,A)$ is $B$-measurable.

**Definition 0.6.2.** A time homogeneous complex transition function (abbreviated CTF) on $\mathbb{R}^d$ is a family $\mu_t^{\mathbb{C}}, t \geq 0$ of complex transition kernels such that $\mu_0^{\mathbb{C}}(x,dy) = \delta(x-y)$ for all $x$, where $\delta(x-y)$ is the Dirac measure in $x$, and such that for every non-negative $s,t$, the Chapman-Kolmogorov equation

$$\int \mu_s^C(x, dy) \mu_t^C(y, A) = \mu_{s+t}^C(x, A) \qquad (0.6.4)$$

is satisfied.

**Definition 0.6.3.** A CTF will be called regular, if there exists a positive constant $K$ such that for all $x$ and $t > 0$, the norm $\|\mu_t^C(x, \cdot)\|$ of the measure $A \mapsto \mu_t^C(x, A)$ does not exceed $\exp(Kt)$.

**Remark 0.6.1.** We remaind that if $T_t$ is a strongly continuous semigroup of bounded linear operators in $C_0(\mathbb{R}^d)$, then there exists a timehomogeneous CTF $\mu_t^C(x, dy)$ such that

$$T_t f(x) = \int \mu_t^C(x, dy) f(y). \qquad (0.6.5)$$

**Remark 0.6.2.** In fact, the existence of a measure $\mu_t^C(x, \cdot)$ such that (2.2) is satisfied follows from the Riesz-Markov theorem, and the semigroup identity $T_s T_t = T_{s+t}$ is equivalent to the Chapman-Kolmogorov equation. Since $\int \mu_t^C(x, dy) f(y)$ is continuous for all $f \in C_0(\mathbb{R}^d)$, it follows by the monotone convergence theorem and the fact that each complex measure is a linear combination of four positive measures that $\mu_t^C(x, A)$ is a Borel function of $x$.

**Definition 0.6.4.** We say that the semigroup $T_t$ is regular, if the corresponding CTF is regular.

**Remark 0.6.3.** Clearly, this Definition 0.6.4 is equivalent to the assumption that $\|T_t\| \leq \exp(Kt)$ for all $t > 0$ and some constant $K$.

**Definition 0.6.5.** Let $\Omega$ be a locally compact space. Set $E_0(\Omega) = (C_0(\Omega))^I, I = (0, 1]$

$$E_M^0(\Omega) = \left\{ (f_\varepsilon)_\varepsilon \in E_0(\Omega) | (\forall \alpha \in \mathbb{N}_0^n)(\exists p \in \mathbb{N}) \sup_{x \in \Omega} |f(x)| = O(\varepsilon^{-p}) \right\},$$

$$N_0(\Omega) = \left\{ (f_\varepsilon)_\varepsilon \in E_0(\Omega) | (\forall \alpha \in \mathbb{N}_0^n)(\forall q \in \mathbb{N}) \sup_{x \in \Omega} |f(x)| = O(\varepsilon^q) \right\}, \qquad (0.6.6)$$

$$G_0(\Omega) = E_M^0(\Omega) / N_0(\Omega).$$

Elements of $E_M^0(\Omega)$ and $N_0(\Omega)$ are called moderate, resp. negligible functions. Note that $E_M^0(\Omega)$ is an algebra with pointwise operations. It is the largest subalgebra of $E_0(\Omega)$ in which $N_0(\Omega)$ is ideal. Thus, $G_0(\Omega)$ is an associative, commutative algebra.

**Definition 0.6.6.**Let $\Omega$ be a locally compact space.Set

$$E(\mathbf{M}(\Omega)) = (\mathbf{M}(\Omega))^I, I = (0,1]$$

$$E_M(\mathbf{M}(\Omega)) = \left\{ (\mu_\varepsilon)_\varepsilon \in E(\mathbf{M}(\Omega)) | (\forall \alpha \in \mathbb{N}_0^n)(\exists p \in \mathbb{N}) \sup_{B \subseteq \Omega} |\mu_\varepsilon(B)| = O(\varepsilon^{-p}) \right\},$$

$$N_0(\mathbf{M}(\Omega)) = \left\{ (\mu_\varepsilon)_\varepsilon \in E(\mathbf{M}(\Omega)) | (\forall \alpha \in \mathbb{N}_0^n)(\forall q \in \mathbb{N}) \sup_{B \subseteq \Omega} |\mu_\varepsilon(B)| = O(\varepsilon^q) \right\}, \quad (0.6.7)$$

$$G(E_M(\mathbf{M}(\Omega))) = E_M(\mathbf{M}(\Omega))/N_0(\mathbf{M}(\Omega)).$$

Elements of $E_M^0(\Omega)$ and $N_0(\Omega)$ are called moderate,resp.negligible Colombeau-Borel measures on $\Omega$.

**Theorem 0.6.2.**(**Generalized Riesz-Markov theorem**) Let $\Omega$ be a locally compact space,then the set $E_M(\mathbf{M}(\Omega))$ of all finite complex regular Colombeau-Borel measures on $\Omega$, which coincides with the set of all continuous linear Colombeau generalized functionals on $E_M^0(\Omega)$.

Any complex $\sigma$-additive Colombeau measure $(\mu_\varepsilon)_\varepsilon$ on $\mathbb{R}^d$ has a representation of the form

$$(\mu_\varepsilon(dy))_\varepsilon = (f_\varepsilon(y)M_\varepsilon(dy))_\varepsilon \quad (0.6.8)$$

with a positive Colombeau measure $(M_\varepsilon(dy))_\varepsilon$ and a bounded for every $\varepsilon \in (0,1]$ complex-valued generalized function $(f_\varepsilon(y))_\varepsilon$.Moreover, the Colombeau measure $(M_\varepsilon(dy))_\varepsilon$ in (0.6.8) is uniquely defined under additional assumption that $(|f_\varepsilon(y)|)_\varepsilon = 1$ for all $y$. If this condition holds, the positive measure $(M_\varepsilon(dy))_\varepsilon$ is called the total variation measure of the complex measure $(\mu_\varepsilon^\mathbb{C}(dy))_\varepsilon$ and is denoted by $(|\mu_\varepsilon^\mathbb{C}|(dy))_\varepsilon$.In general, if (0.6.8) holds, then

$$(\|\mu_\varepsilon^\mathbb{C}\|)_\varepsilon = \left( \int |f_\varepsilon(y)| M_\varepsilon(dy) \right)_\varepsilon \quad (0.6.9)$$

**Definition 0.6.7.**We say that a map $(\mu_\varepsilon^\mathbb{C})_\varepsilon : \mathbb{R}^d \times B(\mathbb{R}^d) \to \widetilde{\mathbb{C}}$ is a Colombeau complex transition kernel, if for every $\varepsilon \in (0,1]$ and $x$, the map $A \mapsto \mu_\varepsilon^\mathbb{C}(x,A)$ is a finite complex Colombeau measure on $\mathbb{R}^d$, and for every $\varepsilon \in (0,1]$ and $A \in \mathbf{B}(\mathbb{R}^d)$ the map $x \to \mu_\varepsilon^C(x, A)$ is $B$-measurable.

**Definition 0.6.8.** A time homogeneous Colombeau complex transition function (abbreviated CCTF) on $\mathbb{R}^d$ is a family $\mu_{t,\varepsilon}^C, t \geq 0$ of Colombeau complex transition kernels such that for every $\varepsilon \in (0, 1]$, $\mu_{0,\varepsilon}^C(x, dy) = \delta(x - y)$ for all $x$, where $\delta(x - y)$ is the Dirac measure in $x$, and such that for every non-negative $s, t$, the Colombeau-Chapman-Kolmogorov equation

$$\left( \int \mu_{s,\varepsilon}^C(x, dy) \mu_{t,\varepsilon}^C(y, A) \right)_\varepsilon = \left( \mu_{s+t,\varepsilon}^C(x, A) \right)_\varepsilon \qquad (0.6.10)$$

is satisfied.

**Definition 0.6.9.** A CCTF will be called regular, if there exists a positive constant $K_\varepsilon = K(\varepsilon)$ such that $\forall \varepsilon \in (0, 1]$ for all $x$ and $t > 0$, the norm $\|\mu_{t,\varepsilon}^C(x, \cdot)\|$ of the measure $A \mapsto \mu_{t,\varepsilon}^C(x, A)$ does not exceed $\exp(K_\varepsilon t)$.

**Remark 0.6.4.** We note that if $\forall \varepsilon \in (0, 1]$, $T_{t,\varepsilon}$ is a strongly continuous semigroup of bounded linear operators in $E_M^0(\mathbb{R}^d)$, then there exists a timehomogeneous CCTF $\mu_{t,\varepsilon}^C(x, dy)$ such that

$$\left( T_{t,\varepsilon} f_\varepsilon(x) \right)_\varepsilon = \left( \int \mu_{t,\varepsilon}^C(x, dy) f_\varepsilon(y) \right)_\varepsilon. \qquad (0.6.11)$$

Now we construct Colombeau measure on the path space corresponding to each regular CCTF, introducing first some (rather standard) notations. Let $\dot{\mathbb{R}}^d$ denote the one point compactification of the Euclidean space $\mathbb{R}^d$ (i.e. $\dot{\mathbb{R}}^d = \mathbb{R}^d \cup \{\infty\}$ and is homeomorphic to the sphere $\mathbf{S}^d$). Let $\dot{\mathbb{R}}_{[s,t]}^d = \prod_{s \leq \tau \leq t} \dot{\mathbb{R}}_\tau^d$ denote the infinite product of $[s, t]$, i.e. it is the set of all functions from $[s, t]$ to $\dot{\mathbb{R}}^d$, the path space.

**Definition 0.6.10.** Let $\mathbf{B}\left( \left( \dot{\mathbb{R}}^d \right)^m \right)$ be the set of all bounded complex Borel functions on $\left( \dot{\mathbb{R}}^d \right)^m$. Set $E\left( \mathbf{B}\left( \left( \dot{\mathbb{R}}^d \right)^m \right) \right) = \left( \mathbf{B}\left( \left( \dot{\mathbb{R}}^d \right)^m \right) \right)^I, I = (0, 1]$

$$E_M\big(\mathbf{B}\big(\big(\dot{\mathbb{R}}^d\big)^m\big)\big) =$$

$$\Big\{(f_\varepsilon(x))_\varepsilon \in E\big(\mathbf{B}\big(\big(\dot{\mathbb{R}}^d\big)^m\big)\big)\big|(\exists p \in \mathbb{N}) = (\sup_{x\in\mathbb{R}^{d\times m}}|f_\varepsilon(x_1,\ldots,x_n)| = O(\varepsilon^{-p}))\Big\},$$

$$N\big(\mathbf{B}\big(\big(\dot{\mathbb{R}}^d\big)^m\big)\big) = \tag{0.6.12}$$

$$\Big\{(f_\varepsilon(x))_\varepsilon \in E\big(\mathbf{B}\big(\big(\dot{\mathbb{R}}^d\big)^m\big)\big)\big|(\forall p \in \mathbb{N}) = (\sup_{x\in\mathbb{R}^{d\times m}}|f_\varepsilon(x_1,\ldots,x_n)| = O(\varepsilon^{-p}))\Big\},$$

$$G\big(\mathbf{B}\big(\big(\dot{\mathbb{R}}^d\big)^m\big)\big) = E_M\big(\mathbf{B}\big(\big(\dot{\mathbb{R}}^d\big)^m\big)\big)/N\big(\mathbf{B}\big(\big(\dot{\mathbb{R}}^d\big)^m\big)\big).$$

**Definition 0.6.11.** Let $\mathfrak{R}^k_{y(\cdot)}\big(\dot{\mathbb{R}}^d_{[s,t]}\big)$ denote the set of a generalized functions on $\dot{\mathbb{R}}^d_{[s,t]}$ having the form

$$\Big(\Phi^f_{t_0,\ldots,t_{k+1,\varepsilon}}(y(\cdot))\Big)_\varepsilon = (f_\varepsilon(y(t_0),\ldots,y(t_{k+1})))_\varepsilon, \tag{0.6.13}$$

where $(f_\varepsilon)_\varepsilon \in G\big(\big(\dot{\mathbb{R}}^d\big)^{k+2}\big)$. The union $\mathfrak{R}^\infty_{y(\cdot)}\big(\dot{\mathbb{R}}^d_{[s,t]}\big) = \bigcup_{k\in\mathbb{N}} \mathfrak{R}^k_{y(\cdot)}\big(\dot{\mathbb{R}}^d_{[s,t]}\big)$ is called the set of cylindrical generalized functions (or generalized functionals) on $\dot{\mathbb{R}}^d_{[s,t]}$. Any **CCTF** $\mu^\mathbb{C}_{t,\varepsilon}(x,dy)$ defines a family of linear generalized functionals $(\mu^x_{s,t,\varepsilon})_\varepsilon, x \in \mathbb{R}^d$ on $\mathfrak{R}^\infty_{y(\cdot)}\big(\dot{\mathbb{R}}^d_{[s,t]}\big)$ by the formula

$$\Big(\mu^x_{s,t,\varepsilon}\Big[\Phi^f_{t_0,\ldots,t_{k+1,\varepsilon}}\Big]\Big)_\varepsilon =$$
$$\tag{0.6.14}$$
$$\Big(\int f_\varepsilon(x,y_1,\ldots y_{k+1})\mu^\mathbb{C}_{t_1-t_0,\varepsilon}(x,dy_1)\mu^\mathbb{C}_{t_2-t_1,\varepsilon}(y_1,dy_2)\ldots\mu^\mathbb{C}_{t_{k+1}-t_{k1},\varepsilon}(y_k,dy_{k+1})\Big)_\varepsilon.$$

**Theorem 0.6.3.** If the Colombeau semigroup $(T_{t,\varepsilon})_\varepsilon$ in $E^0_M(\mathbb{R}^d)$ is regular and $\mu^\mathbb{C}_{t,\varepsilon}(x,dy)$ is its corresponding **CCTF**, then $\forall \varepsilon \in (0,1]$ the functional $\mu^x_{s,t,\varepsilon}\Big[\Phi^f_{t_0,\ldots,t_{k+1,\varepsilon}}\Big]$ is bounded. Hence, it can be extended by continuity to a unique bounded linear functional $\tilde{\mu}^x_{s,t,\varepsilon}$ on $E^0_M\big(\dot{\mathbb{R}}^d_{[s,t]}\big)$, and consequently there exists a (regular) complex Colombeau-Borel measure $\mathbf{D}^x_{s,t,\varepsilon}$ on the path space $\dot{\mathbb{R}}^d_{[s,t]}$ such that

$$\left( \widetilde{\mu}_{s,t,\varepsilon}^{x}[F_{\varepsilon}] \right)_{\varepsilon} = \left( \int F_{\varepsilon}(y(\tau)) \mathbf{D}_{s,t,\varepsilon}^{x}[y(\tau)] \right)_{\varepsilon} \qquad (0.6.15)$$

for all $(F_{\varepsilon})_{\varepsilon} \in G_0\left( \dot{\mathbb{R}}_{[s,t]}^d \right)$. In particular,

$$((T_{t,\varepsilon} f_{\varepsilon})(x))_{\varepsilon} = \left( \int f_{\varepsilon}(y(t)) \mathbf{D}_{s,t,\varepsilon}^{x}[y(\tau)] \right)_{\varepsilon}. \qquad (0.6.16)$$

**Theorem 0.6.4.**Let $(B_{\varepsilon})_{\varepsilon \in (0,1]}$ and $(A_{\varepsilon})_{\varepsilon \in (0,1]}$ be the nets of linear operators in $E_M^0\left( \dot{\mathbb{R}}_{[s,t]}^d \right)$
such that $\forall \varepsilon \in (0,1]$, $A_{\varepsilon}$ is bounded and $(B_{\varepsilon})_{\varepsilon}$ is the generator of a strongly continuous regular Colombeau semigroup $(T_{t,\varepsilon})_{\varepsilon}$. Then $(A_{\varepsilon} + B_{\varepsilon})_{\varepsilon \in (0,1]}$ is also the generator of a regular semigroup, which we denote by $\left( \widetilde{T}_{t,\varepsilon} \right)_{\varepsilon}$.

**Definition 0.6.12.**Let $\mathscr{F}(\mathbb{R}^d)$ denote the Banach space of Fourier transforms of elements of $\mathbf{M}(\mathbb{R}^d)$, i.e.the space of continuous functions on $\mathbb{R}^d$ of the form

$$V_{\mu}(x) = \int_{\mathbb{R}^d} \exp(ip \cdot x) \mu(dp) \qquad (0.6.17)$$

for some $\mu \in \mathbf{M}(\mathbb{R}^d)$, with the induced norm $\|V_{\mu}\| = \|\mu\|$. Since $\mathbf{M}(\mathbb{R}^d)$ is a Banach algebra with convolution as the multiplication, it follows that $\mathscr{F}(\mathbb{R}^d)$ is also a Banach algebra with respect to the pointwise multiplication. We say that an element $f \in \mathscr{F}(\mathbb{R}^d)$
is infinitely divisible if there exists a family $f_t, t \geq 0$ of elements of $\mathscr{F}(\mathbb{R}^d)$ such that $f_0 = 1, f_1 = f,$ and $f_{t+s} = f_t \cdot f_s$ for all positive $s, t$.Clearly if $f$ is infinitely divisible, then it has no zeros and a continuous function $g = \ln f$ is well defined (and is unique up to an imaginary shift). Moreover, the family $f_t$ has the form $f_t = \exp(tg)$
**Definition 0.6.13.**We say that a continuous function $g$ on $\mathbb{R}^d$ is a complex characteristic exponent (abbreviated **CCE**), if $\exp(g)$ is an infinitely divisible element of
$\mathscr{F}(\mathbb{R}^d),$ or equivalently, if $\forall t, t \geq 0 : \exp(tg) \in \mathscr{F}(\mathbb{R}^d).$

**Remark 0.6.5.**It follows from the definitions that the set of spatially homogeneous

**CTF**'s $\mu_t^C(dx)$ is in one-to-one correspondence with **CCE** $g$, in such a way that for any positive $t$ the function $\exp(tg)$ is the Fourier transform of the transition measure $\mu_t^C(dx)$.

**Definition 0.6.14.**[54].(**Hardy spaces for the unit disk**).The Hardy space $H^p(\mathbf{D})$ for $0 < p < \infty$ is the class of holomorphic functions $f(z)$ on the open unit disk satisfying

$$\sup_{0 < r < 1} \left( \frac{1}{2\pi} \int_0^{2\pi} \left| f(re^{i\theta}) \right|^p d\theta \right)^{\frac{1}{p}} < \infty. \qquad (0.6.18)$$

This class $\mathbf{H}^p(\mathbf{D})$ is a vector space.The number on the left side of the above inequality is the Hardy space p-norm for f, denoted by $\|f\|_{\mathbf{H}^p}$ It is a norm when $p \geq 1$, but not when $0 < p < 1$.

**Definition 0.6.15.**[54].(**Hardy spaces for the upper half plane**).The Hardy space on the upper half-plane $H^p(\mathbf{H})$ is defined to be the space of holomorphic functions $f(z)$ on $\mathbf{H}$ with bounded (quasi-)norm,the norm being given by

$$\|f\|_{H^p} = \sup_{0 < y} \left( \int \left| f(x+iy) \right|^p dx \right)^{\frac{1}{p}} < \infty. \qquad (0.6.19)$$

Note that the unit disk $\mathbf{D}$ and the upper half-plane $\mathbf{H}$ can be mapped to one-another by means of Möbius transformations, they are not interchangeable as domains for Hardy spaces. Contributing to this difference is the fact that the unit circle has finite (one-dimensional) Lebesgue measure while the real line does not. However, for $H^2$, one may still state the following theorem: Given the Möbius transformation $\mathbf{m} : \mathbf{D} \to \mathbf{H}$
with $\mathbf{m}(z) = i(1+z)(1-z)^{-1}$ then there is an isometric isomorphism $M : H^p(\mathbf{D}) \to H^p(\mathbf{H})$
with $[Mf](x) = \sqrt{\pi}\,(1-z)^{-1}f(m(z))$.

**Definition 0.6.16.**A sequence of points $\{a_n\}_{n=1}^{\infty}$ inside the unit disk $\mathbf{D}$ is said to satisfy the Blaschke condition when $\sum_{n=1}^{\infty}(1 - |a_n|) < \infty$. Given a sequence obeying the Blaschke condition, the Blaschke product $B(z)$ is defined as

$$B(z) = \prod_n B_n(z), B_n(z) = -\frac{|a_n|}{a_n} \frac{z - a_n}{1 - \overline{a}_n \cdot z}. \qquad (0.6.20)$$

The Blaschke product $B(z)$ is analytic in the open unit disc, and is zero at the $a_n$

only.

**Definition 0.6.17.**[54].(**i**) One says that $G(z)$ is an **outer** (**exterior**) function if it takes the form

$$G(z) = c \exp\left[\frac{1}{2\pi}\int_{-\pi}^{\pi}\frac{e^{i\theta}+z}{e^{i\theta}-z}\log\varphi(e^{i\theta})d\theta\right] \qquad (0.6.21)$$

for some complex number $c$ with $|c|=1$, and some positive measurable function $\varphi$ on the unit circle such that $\log\varphi$ is integrable on the circle.

(**ii**) One says that $h(z)$ is an **inner** (**interior**) function if and only if $|h(z)|\leq 1$ on the unit disk and the limit $\lim_{r\to 1^-}h(re^{i\theta})$ exists for almost all $\theta$ and its modulus is equal to $1$.

Let $G$ be an outer function represented as above from a function $\varphi$ on the circle. Replacing $\varphi$ by $\varphi^{\alpha}, \alpha > 0$, a family $G_{\alpha}$ of outer functions is obtained, with the properties:

$$G_1 = G, G_{\alpha+\beta} = G_{\alpha}G_{\beta},$$

$$\qquad (0.6.22)$$

$$|G_{\alpha}| = |G|^{\alpha}.$$

**Remark 0.6.6.**For $0 < p \leq \infty$, every non-zero function $f(z)$ in $H^p$ can be written as the product $f(z) = G(z)h(z)$ where $G(z)$ is an outer function and $h(z)$ is an inner function.

**Remark 0.6.7.**It is easy to show that if $f_1 \in \mathcal{F}(\mathbb{R})$ is infinite divisible and such that the measures corresponding to all functions $f_t, t > 0$, are concentrated on the half line $\mathbb{R}_+$ (complex generalisation of subordinators) and have densities from $\mathcal{L}_2(\mathbb{R}_+)$, then $f_1$ belongs to the Hardy space $H^2$ of analytic functions on the upper half plane, which have no Blaschke product in its canonical decomposition.

**Theorem 0.6.5.**[49]. If $V$ is a **CCE**, then the solution to the Cauchy problem

$$\frac{\partial u}{\partial t} = V\left(\frac{1}{i}\frac{\partial}{\partial y}\right)u \qquad (0.6.23)$$

defines a strongly continuous and spatially homogeneous semigroup $T_t$ of bounded linear operators in $C_0(\mathbb{R}^d)$ (i.e. $(T_t u_0)(y)$ is the solution to equation (0.6.23) with the initial function $u_0$). Conversely, each such semigroup is the solution to the Cauchy problem of an equation of type (0.6.23) with some **CCE** $g$.

**Definition 0.6.18.**[49].We say that a **CCE** $V$ is regular, if equation (0.6.23) defines a regular semigroup.

**Theorem 0.6.6.**[49].Let $V \in \mathcal{F}(\mathbb{R}^d)$, i.e. it is given by (0.6.17) with $\mu^\mathbb{C} \in \mathbf{M}(\mathbb{R}^d)$. Then

$V$ is a regular **CCE**. Moreover, if the positive measure $M(dy)$ in the representation (0.6.2) for $\mu^\mathbb{C}$ has no atom at the origin, i.e. $M(\{0\}) = 0$, then the corresponding measure $\mathbf{D}_{0,t}^x$ on the path space is concentrated on the set of piecewise-constant paths in $\mathbb{R}_{[s,t]}^d$ with a finite number of jumps. In other words,$\mathbf{D}_{0,t}^x$ is the measure of a jump- process.

Let us consider the pseudo-differential equation of the Schrödinger type

$$\frac{\partial \tilde{u}}{\partial t} = -G[(-\Delta)^\alpha]\tilde{u} + \left(A, \frac{\partial}{\partial x}\right)\tilde{u} + V(x)\tilde{u}, \qquad (0.6.24)$$

where $G$ is a complex constant with a non-negative real part,$\alpha$ is any positive constant,

$A$ is a real-valued vector (if $\operatorname{Re} G > 0$, then $A$ can be also complex- valued), and $V(x)$ is a complex-valued function of form (0.6.17).The equation on the inverse Fourier transform

$$u(y) = \int\limits_{\mathbb{R}^d} dx \exp(-iyx)\tilde{u}(x) \qquad (0.6.25)$$

of $\tilde{u}(x)$ has the form

$$\frac{\partial u}{\partial t} = -G(y^2)^\alpha + i(A,y)u + V\left(\frac{1}{i}\frac{\partial}{\partial y}\right)u. \qquad (0.6.26)$$

**Theorem 0.6.7.**[49].The solution to the Cauchy problem of equation (0.6.26) can be written in the form of a complex Feynman-Kac formula

$$u(t, u) = \int \exp\left\{ -\int_0^t \left[ (G(q^2(\tau))^\alpha - (A, q(\tau))) \right] \right\} u_0(q(t)) \mathbf{D}_{0,t}^x[q(\cdot)], \qquad (0.6.27)$$

where $\mathbf{D}_{0,t}^x[q(\cdot)]$ is the measure of the jump process corresponding to equation (0.6.23).

As another example, let us consider the case of complex anharmonic oscillator, i.e. the equation

$$\frac{\partial \widetilde{\psi}}{\partial t} = \frac{1}{2}(G\Delta - x^2 - iV(x))\widetilde{\psi}, \qquad (0.6.28)$$

where $V = V_\mu \in \mathcal{F}(\mathbb{R}^d)$. The Fourier transform of this equation has the form

$$\frac{\partial \psi(t, p)}{\partial t} = \frac{1}{2}\left( \Delta - Gp^2 - iV\left( \frac{1}{i}\frac{\partial}{\partial p} \right) \right) \psi(t, p). \qquad (0.6.29)$$

**Theorem 0.6.8**.[49].If $\operatorname{Re} G \geq 0$, the Cauchy problem of equation $(0.6.29)$ defines a regular semigroup of operators in $C_0(\mathbb{R}^d)$, and thus can be presented as the path integral from (Proposition 2.1.[49]).

**Remark 0.6.8**.[49].The statement of the Proposition can be generalised easily to the following situation, which includes all Schrödinger equations, namely to the case of the equation

$$\frac{\partial \psi}{\partial t} = i(A - B)\psi. \qquad (0.6.30)$$

where $A$ is selfadjoint operator, for which therefore exists (according to spectral theory) a unitary transformation $U$ such that $UAU^{-1}$ is the multiplication operator in some
$\mathcal{L}_2(X, d\mu)$, where $X$ is locally compact, and $B$ is such that $UBU^{-1}$ is a bounded operator
in $C_0(X)$.
We shall consider now a formal Schrödinger operator with a magnetic field, namely the operator

$$H = \left(i\frac{\partial}{\partial x} + A(x)\right)^2 + V(x) \qquad (0.6.31)$$

where $x = (x_1,\ldots,x_d)$ and $\frac{\partial}{\partial x} = \left(\frac{\partial}{\partial x_1},\ldots,\frac{\partial}{\partial x_d}\right)$ is the gradient operator in $\mathbb{R}^d$, under the following conditions [49]:

(**C**1) the magnetic vector-potential $\mathbf{A} = (A^1,\ldots,A^d)$ is a bounded measurable $\mathbb{R}^d$-valued function on $\mathbb{R}^d$,

(**C**2) the potential $V$ and the divergence $div\mathbf{A} = \sum_{j=1}^{d}\frac{\partial A^j}{\partial x_j}$ of $\mathbf{A}$, defined in the sense of distributions, are both Borel measures,

(**C**3) if $d > 1$ there exist $\alpha > d - 2$ and $C > 0$ such that for all $x \in \mathbb{R}^d$ and $r \in (0,1]$

$$|\mathbf{div}\mathbf{A}|(B_r(x)) \le Cr^\alpha, |V|(B_r(x)) \le Cr^\alpha \qquad (0.6.32)$$

where $|\mathbf{div}\mathbf{A}|(\cdot)$ and $|V|(\cdot)$ denote the total variations of the (possibly non-positive) measures $V$ and $\mathbf{div}\mathbf{A}$ respectively, and $B_r(x)$ denotes the ball in $\mathbb{R}^d$ of the radius $r$ centered at $x$; if $d = 1,$ then the same holds for $\alpha = 0$.

Let us consider the following more general Cauchy problem

$$\frac{\partial u(t,x)}{\partial t} = -DHu(t,x), u(0,x) = u_0(x), \qquad (0.6.33)$$

where the (generalised) diffusion coefficient $D$ is an arbitrary complex number such that $\epsilon = \mathrm{Re}\, D \ge 0, |D| > 0$. In the interaction representation, equation (0.6.32) takes the form

$$u(t,x) =$$

$$\qquad (0.6.34)$$

$$[\exp(-Dt\Delta/2)]u_0(x) - D\int_0^t [\exp(-D(t-s)\Delta/2)](W(x) - 2i(\mathbf{A}(x),\nabla))u(s,x)ds.$$

Here $W(x) = V(x) + |\mathbf{A}(x)|^2 - i\mathbf{div}\mathbf{A}(x).$ Equation (0.6.34) can be formally solved by iterations. In the case of the Green function $G^D(t,x,y)$ of equation (0.6.34),i.e. its solution with the Dirac initial condition $u_0(x) = \delta(x - y)$, the iteration procedure leads to the following representation:

$$G^D(t,x,y) = \sum_{k=0}^{\infty} I_k^D(t,x,y),$$

$$I_k^D(t,x,y) =$$

$$-D \int_0^t \int_{\mathbb{R}^d} I_{k-1}^D(t-s,x,\xi)(W(d\xi) + 2iD^{-1}s^{-1}(\mathbf{A},(\xi-y))d\xi)G_{free}^D(s,\xi-y)ds,$$

$$(0.6.35)$$

$$G_{free}^D(s,\xi-y) = (2\pi tD)^{-d/2} \exp\left[-\frac{(\xi-y)^2}{2Dt}\right].$$

**Theorem 0.6.9.**[49].**(i)** If $\epsilon = \operatorname{Re} D > 0$, then all terms of series (0.6.35) are well defined
as absolutely convergent integrals, the series itself is absolutely convergent and its sum $G^D(t,x,y)$ is continuous in $x,y \in \mathbb{R}^d$, $t > 0$ (and $D$) and satisfies the following estimate

$$G^D(t,x,y) \le K G_{free}^{\frac{|D|^2}{\epsilon}}(t,x-y) \exp(B|x-y|) \qquad (0.6.36)$$

uniformly for $t \le t_0$ with any fixed $t_0$, where $B,K$ are constants.
**(ii)** The integral operators

$$[U^D(t)u_0](t,x) = \int u^D(t,x,y)u_0(y)dy \qquad (0.6.37)$$

defining the solutions to equation (0.6.34) for $t \in [0,t_0]$ form a uniformly bounded family of operators $\mathcal{L}_2(\mathbb{R}^d) \to \mathcal{L}_2(\mathbb{R}^d)$.
**(iii)** If $D$ is real, i.e. $D = \epsilon > 0$, then there exists a constant $\omega > 0$ such that the Green function $G^\epsilon$ has the asymptotic representation

$$G^{\epsilon}(t,x,y) = G^{\epsilon}_{free}(t,x,y)(1 + O(t^{\omega}) + O(|x-y|)) \qquad (0.6.38)$$

for small $t$ and $x - y$. In case of vanishing $\mathbf{A}$, the multiplier $\exp(B|x-y|)$ in (0.6.36) can be omitted, and the term $O(|x-y|)$ in (0.6.38) can be dropped. In this case, formula (0.6.38) gives global (uniform for all $x,y$) small time asymptotics for $G^{\epsilon}(t,x,y)$.
(**iv**) One can give rigorous meaning to formal expression (0.6.31) as a bounded below self adjoint operator in such a way that the family (0.6.37) of operators $U^D(t)$ (giving solutions to equation (0.6.34), which is formally equivalent to the evolutionary equation (0.6.33) with the formal generator $DH$) coincides with the semigroup $\exp(-tDH)$ defined by means of the functional calculus. Hence for the integral kernel of the operators $\exp(-tDH)$ the estimates (0.6.36) and (0.6.38) hold.

The rigorous construction a path integral representation for the solution $U^D(t)u_0$ of Eq.(0.6.33), by construction a measure on a path space that is supported on the set of continuous piecewise linear paths, was given in paper [49]. Denote this set by $\mathbf{CPL}^{\mathbb{R}^d}$. Let $\mathbf{CPL}^{x,y}(0,t)$ denote the class of paths $q : [0,t] \to \mathbb{R}^d$ from $\mathbf{CPL}^{\mathbb{R}^d}$ joining $x$ and $y$ in time $t$, i.e. such that $q(0) = x, q(t) = y$. By $\mathbf{CPL}_n^{x,y}(0,t)$ we denote its subclass consisting of all paths from $\mathbf{CPL}^{x,y}(0,t)$ that have exactly $n$ jumps of their derivative. Clearly, each $\mathbf{CPL}_n^{x,y}(0,t)$ is parametrised by the simplex

$$\mathbf{sim}_t^n = \{s_1, \ldots, s_n | 0 < s_1 < s_2 < \ldots < s_n < t\} \qquad (0.6.39)$$

of the times of jumps $s_1, \ldots, s_n$ of the derivatives of a path and by $n$ positions $q(s_j)$, $j = 1, \ldots, n$, of this path at these points. In other words, an arbitrary path in $\mathbf{CPL}_n^{x,y}(0,t)$
has the form

$$q_n(s) = q_{\eta_1,\ldots,\eta_n}^{s_1,\ldots,s_n} = \eta_j + (s - s_j)\frac{\eta_{j+1} - \eta_j}{s_{j+1} - s_j},$$

$$\qquad (0.6.40)$$

$$s_j \in [s_j, s_{j+1}],$$

where $s_0 = 0, s_{n+1} = t, \eta_0 = x, \eta_{n+1} = y$. Therefore $\mathbf{CPL}^{x,y}(0,t) = \cup_{n=0}^{\infty} \mathbf{CPL}_n^{x,y}(0,t)$. To any $\mathbb{R}^{d+1}$-valued Borel measure $\widetilde{M} = (\mu, \nu) = (\mu, \nu^1, \ldots, \nu^n)$ on $\mathbb{R}^{d+1}$ there corresponds a $\sigma$-finite complex measure $\mathbf{D}^{\mathbf{CPL}}[q(\tau)]$ on $\mathbf{CPL}^{x,y}(0,t)$, which is defined as the sum of the measures $\mathbf{M}_n^{\mathbf{CPL}}$ $n \in \mathbb{N}$ on the finite-dimensional spaces $\mathbf{CPL}_n^{x,y}(0,t)$ such that $\mathbf{M}_0^{\mathbf{CPL}}$ is
just the unit measure on the one-point set $\mathbf{CPL}_0^{x,y}(0,t)$ and each $\mathbf{M}_n^{\mathbf{CPL}}, n > 0,$ is defined
in the following way: if $\Phi(\cdot)$ is a functional on $\mathbf{CPL}^{x,y}(0,t)$, then:

$$\int_{\mathbf{CPL}_n^{x,y}(0,t)} \Phi(q(\tau)) \mathbf{D}_n^{\mathbf{CPL}}[q(\tau)] = \int_{\mathbf{sim}_t^n} ds_1 \ldots ds_n \times$$

$$(0.6.41)$$

$$\int_{\mathbb{R}^d} \ldots \int_{\mathbb{R}^d} \left( \mu + 2i \frac{(\eta_2 - \eta_1, \nu)}{s_2 - s_1} \right) (d\eta_1) \ldots \left( \mu + 2i \frac{(y - \eta_n, \nu)}{t - s_n} \right) (d\eta_n) \Phi(q(\tau)).$$

Now, we let $\widetilde{M} = (W, D^{-1}\mathbf{A}(x)dx)$.

**Theorem 0.6.10**. For any $D$ with $\epsilon = \operatorname{Re} D > 0$, the Green function $G^D(t,x,y)$ of equation
(0.6.33) has the following path integral representation:

$$G^D(t,x,y) = \int_{\mathbf{CPL}^{x,y}(0,t)} \Phi_D(q(\tau)) \exp\left\{ \frac{1}{2D} \int_0^t \dot{q}^2(s)ds \right\} \mathbf{D}^{\mathbf{CPL}}[q(\tau)] \qquad (0.6.42)$$

with $q(s)$ given by Eq.(0.6.40) and

$$\Phi_D(q(\tau)) = D^n \prod_{j=1}^{n+1} (2\pi D(s_j - s_{j-1}))^{-d/2}. \qquad (0.6.43)$$

For any $u_0 \in \mathcal{L}_2(\mathbb{R}^d)$ the solution $u(t,x)$ of the Cauchy problem (0.6.33) with $D = i$
has the form

$$u(t,x) = \lim_{\delta \to 0} \int_{\mathbf{CPL}^{x,y}(0,t)} \int_{\mathbb{R}^d} u_0(x) \Phi_{i+\delta}(q(\tau)) \exp\left\{ \frac{1}{2(i+\delta)} \int_0^t \dot{q}^2(s)\,ds \right\} \mathbf{D}^{\mathbf{CPL}}[q(\tau)]. \quad (0.6.44)$$

Here the limit is understood in $\mathcal{L}_2(\mathbb{R}^d)$ sence.

**Theorem 0.6.11.** Suppose the vector potential A vanishes and V satisfies assumption (**C**3). Suppose additionally that $V$ has no atom at the origin and is a finite positive measure so that $\mu_V = V(\mathbb{R}^d) > 0$. Let the paths of **CPL** are parametrised by Eq.(0.6.40). Let $\mathbf{E}_t$ denote the expectation with respect to the process of jumps $\eta_j$ which are identically independently distributed according to the probability measure $V/\mu_V$, and which occur at times $s_j$ from $[0,t]$ distributed according to Poisson process of intensity $\mu_V^{\mathbb{R}}$. Then the integral (0.6.42) can be written in the form

$$G^D(t,x,y) = e^{\mu_V} \mathbf{E}_t \left[ \Phi_D(q(\tau)) \exp\left\{ -\frac{1}{2D} \int_0^t \dot{q}^2(s)\,ds \right\} \mathbf{D}^{\mathbf{CPL}}[q(\tau)] \right]. \quad (0.6.45)$$

**Definition 0.6.19.** (**Colombeau-Hardy spaces for the unit disk**). Let $H^p(\mathbf{D})$ be the Hardy space for $0 < p < \infty$. Set $E(H^p(\mathbf{D})) = (H^p(\mathbf{D}))^I, I = (0,1]$,

$$E_M(H^p(\mathbf{D})) =$$

$$\left\{ (f_\varepsilon)_\varepsilon \in E(H^p(\mathbf{D})) | (\forall \alpha \in \mathbb{N})(\exists p \in \mathbb{N}) \left[ \sup_{z \in \mathbf{D}} \left| \frac{\partial^\alpha f_\varepsilon(z)}{\partial z^\alpha} \right| = O(\varepsilon^{-p}) \right] \right\},$$

$$N(H^p(\mathbf{D})) = \qquad\qquad (0.6.46)$$

$$\left\{ (f_\varepsilon)_\varepsilon \in E(H^p(\mathbf{D})) | (\forall \alpha \in \mathbb{N})(\forall q \in \mathbb{N}) \left[ \sup_{z \in \mathbf{D}} \left| \frac{\partial^\alpha f_\varepsilon(z)}{\partial z^\alpha} \right| = O(\varepsilon^q) \right] \right\},$$

$$G(H^p(\mathbf{D})) = E_M(H^p(\mathbf{D}))/N(H^p(\mathbf{D})).$$

**Definition 0.6.20.**(**Colombeau-Hardy spaces for the upper half plane**).Let $H^p(\mathbf{H})$ be the Hardy space for $0 < p < \infty$. Set $E(H^p(\mathbf{H})) = (H^p(\mathbf{H}))^I, I = (0,1],$

$$E_M(H^p(\mathbf{H})) =$$

$$\left\{ (f_\varepsilon)_\varepsilon \in E(H^p(\mathbf{H})) | (\forall \alpha \in \mathbb{N})(\exists p \in \mathbb{N}) \left[ \sup_{z \in \mathbf{D}} \left| \frac{\partial^\alpha f_\varepsilon(z)}{\partial z^\alpha} \right| = O(\varepsilon^{-p}) \right] \right\},$$

$$N(H^p(\mathbf{H})) = \tag{0.6.47}$$

$$\left\{ (f_\varepsilon)_\varepsilon \in E(H^p(\mathbf{H})) | (\forall \alpha \in \mathbb{N})(\forall q \in \mathbb{N}) \left[ \sup_{z \in \mathbf{D}} \left| \frac{\partial^\alpha f_\varepsilon(z)}{\partial z^\alpha} \right| = O(\varepsilon^q) \right] \right\},$$

$$G(H^p(\mathbf{H})) = E_M(H^p(\mathbf{H}))/N(H^p(\mathbf{H})).$$

**Definition 0.6.21.**A sequence of points $\left\{ (a_{n,\varepsilon})_\varepsilon \right\}_{n=1}^\infty$ inside the unit disk $\widetilde{\mathbf{D}} \subset \widetilde{\mathbb{C}}$ is said to satisfy the Blaschke condition when $\left( \sum_{n=1}^\infty (1 - |a_{n,\varepsilon}|) \right)_\varepsilon < \infty$. Given a sequence obeying the Blaschke condition, the Blaschke generalized product $(B_\varepsilon(z))_\varepsilon$ is defined as

$$B_\varepsilon(z) = \prod_n (B_{n,\varepsilon}(z))_\varepsilon,$$

$$\tag{0.6.48}$$

$$B_{n,\varepsilon}(z) = -\frac{|a_{n,\varepsilon}|}{a_{n,\varepsilon}} \frac{z - a_{n,\varepsilon}}{1 - \overline{a}_{n,\varepsilon} \cdot z}.$$

The Blaschke generalized product $(B_\varepsilon(z))_\varepsilon$ is analytic in the open unit disc $\widetilde{\mathbf{D}}$, and is zero at the $(a_{n,\varepsilon})_\varepsilon$ only.

**Definition 0.6.22.**(**i**) One says that $(G_\varepsilon(z))_\varepsilon$ is an **outer** (**exterior**) generalized function if it takes the form

$$(G_\varepsilon(z))_\varepsilon = ((c_\varepsilon)_\varepsilon \left( \exp\left[ \frac{1}{2\pi} \int\limits_{-\pi}^{\pi} \frac{e^{i\theta}+z}{e^{i\theta}-z} \log \varphi_\varepsilon(e^{i\theta}) d\theta \right] \right)_\varepsilon \qquad (0.6.49)$$

**Definition 0.6.23.** Let $\mathscr{F}^{\blacktriangle}(\mathbb{R}^d) = G(\mathscr{F}(\mathbb{R}^d))$ denote the Colombeau algebra of Colombeau-Fourier transforms of elements of $\mathbf{M}^{\blacktriangle}(\mathbb{R}^d) = G(\mathbf{M}(\mathbb{R}^d))$, i.e. the algebra of Colombeau generalized functions on $\mathbb{R}^d$ of the form

$$(V_{\mu,\varepsilon}(x))_\varepsilon = \left( \int\limits_{\mathbb{R}^d} \exp(ip \cdot x) \mu_\varepsilon(dp) \right)_\varepsilon \qquad (0.6.50)$$

for some $(\mu_\varepsilon)_\varepsilon \in \mathbf{M}^{\blacktriangle}(\mathbb{R}^d)$, with the induced $\widetilde{\mathbb{R}}$-valued norm $\|(V_{\mu,\varepsilon})_\varepsilon\| = (\|\mu_\varepsilon\|)_\varepsilon$. Since $\mathbf{M}^{\blacktriangle}(\mathbb{R}^d)$ is a Colombeau-Banach algebra with convolution as the multiplication, it follows that $\mathscr{F}^{\blacktriangle}(\mathbb{R}^d)$ is also a Colombeau-Banach algebra with respect to the pointwise multiplication. We say that an element $(f_\varepsilon)_\varepsilon \in \mathscr{F}^{\blacktriangle}(\mathbb{R}^d)$ is infinitely divisible if there exists a family $(f_{t,\varepsilon})_\varepsilon, t \geq 0$ of elements of $\mathscr{F}^{\blacktriangle}(\mathbb{R}^d)$ such that $(f_{0,\varepsilon})_\varepsilon = 1, (f_{1,\varepsilon})_\varepsilon = (f_\varepsilon)_\varepsilon$, and $(f_{t+s,\varepsilon})_\varepsilon = (f_{t,\varepsilon} \cdot f_{s,\varepsilon})_\varepsilon$ for all positive $s, t$. Clearly if $(f_\varepsilon)_\varepsilon$ is infinitely divisible, then it has no zeros and a continuous generalized function $(g_\varepsilon)_\varepsilon = (\ln f_\varepsilon)_\varepsilon$ is well defined (and is unique up to an imaginary shift). Moreover, the family $(f_{t,\varepsilon})_\varepsilon$ has the form $(f_{t,\varepsilon})_\varepsilon = (\exp(tg_\varepsilon))_\varepsilon$.

**Definition 0.6.24.** We say that a Colombeau generalized function $(g_\varepsilon)_\varepsilon \in G_0(\mathbb{R}^d)$ is a Colombeau complex characteristic exponent (abbreviated **CCCE**), if $\exp(g_\varepsilon)_\varepsilon$ is an infinitely divisible element of $\mathscr{F}^{\blacktriangle}(\mathbb{R}^d)$, or equivalently, if $\forall t, t \geq 0$ : $(\exp(tg_\varepsilon))_\varepsilon \in \mathscr{F}^{\blacktriangle}(\mathbb{R}^d)$.

**Remark 0.6.8.** It follows from the definitions that the set of spatially homogeneous **CCTF**'s $(\mu_{t,\varepsilon}(dx))_\varepsilon$ is in one-to-one correspondence with **CCCE** $(g_\varepsilon)_\varepsilon$, in such a way that for any positive $t$ the generalized function $(\exp(tg_\varepsilon))_\varepsilon$ is the Colombeau-Fourier transform of the transition Colombeau measure $(\mu_{t,\varepsilon}(dx))_\varepsilon$.

**Remark 0.6.9.** It is easy to show that if $(f_{1,\varepsilon})_\varepsilon \in \mathscr{F}^{\blacktriangle}(\mathbb{R})$ is infinite divisible and such that the measures corresponding to all functions $(f_{t,\varepsilon})_\varepsilon, t > 0$, are concentrated on the half line $\mathbb{R}_+$ (complex generalisation of subordinators) and have densities from Colombeau module $G_{\mathscr{L}_2(\mathbb{R}_+)}$, then $(f_{1,\varepsilon})_\varepsilon$ belongs to the Colombeau algebra $G(H^2)$ of Colombeau generalized analytic functions on the upper half plane, which have no Colombeau-Blaschke product in its canonical decomposition.

**Theorem 0.6.12.** If $(V_\varepsilon)_\varepsilon$ is a **CCCE**, then the Colombeau solution to the

Colombeau- Cauchy problem

$$\left(\frac{\partial u_\varepsilon}{\partial t}\right)_\varepsilon = \left(V_\varepsilon\left(\frac{1}{i}\frac{\partial}{\partial y}\right)u_\varepsilon\right)_\varepsilon,$$

$$(u_\varepsilon(0,y))_\varepsilon = (u_{0,\varepsilon}(y))_\varepsilon$$

(0.6.51)

defines a strongly continuous and spatially homogeneous Colombeau semigroup $(T_{t,\varepsilon})_\varepsilon$ of bounded linear operators in $C_0(\mathbb{R}^d)$ (i.e. $((T_{t,\varepsilon}u_{0,\varepsilon})(y))_\varepsilon$ is the Colombeau solution to equation (0.6.51) with the initial Colombeau generalized function $(u_{0,\varepsilon})_\varepsilon$). Conversely, each such semigroup is the solution to the Cauchy problem of an equation of type (0.6.51) with some **CCCE** $(g_\varepsilon)_\varepsilon$.

**Definition 0.6.25.**We say that a **CCCE** $(V_\varepsilon)_\varepsilon$ is regular, if equation (0.6.51) defines a regular Colombeau semigroup.

**Theorem 0.6.13.**Let $(V_\varepsilon)_\varepsilon \in \mathcal{F}^{\blacktriangle}(\mathbb{R}^d)$, i.e. it is given by (0.6.50) with $(\mu_\varepsilon^{\mathbb{C}})_\varepsilon \in \mathbf{M}^{\blacktriangle}(\mathbb{R}^d)$. Then $(V_\varepsilon)_\varepsilon$ is a regular **CCCE**. Moreover, if the positive Colombeau measure $(M_\varepsilon(dy))_\varepsilon$ in the representation (0.6.8) for $(\mu_\varepsilon^{\mathbb{C}})_\varepsilon$ has no atom at the origin, i.e. $(M_\varepsilon(\{0\}))_\varepsilon = 0$, then the corresponding Colombeau measure $(\mathbf{D}_{0,t,\varepsilon}^x[q(\cdot)])_\varepsilon$ on the path space is concentrated on the set of piecewise-constant paths in $\mathbb{R}_{[s,t]}^d$ with a finite number of jumps. In other words,$(\mathbf{D}_{0,t,\varepsilon}^x[q(\cdot)])_\varepsilon$ is the measure of a Colombeau generalizwd jump-process.

Let us consider now the pseudo-differential equation of the Colombeau-Schrödinger type

$$\left(\frac{\partial \widetilde{u}_\varepsilon}{\partial t}\right)_\varepsilon = -G[(-\Delta)^\alpha](\widetilde{u}_\varepsilon)_\varepsilon + \left(\left(A(x),\frac{\partial}{\partial x}\right)\widetilde{u}_\varepsilon\right)_\varepsilon + (V_\varepsilon(x)\widetilde{u}_\varepsilon)_\varepsilon,$$

(0.6.52)

where $G$ is a complex constant with a non-negative real part,$\alpha$ is any positive constant,

$A$ is a real-valued vector (if $\operatorname{Re} G > 0$, then $(A_\varepsilon(x))_\varepsilon$ can be also $\widetilde{\mathbb{C}}$-valued), and $(V_\varepsilon(x))_\varepsilon$ is a $\widetilde{\mathbb{C}}$-valued function of form (0.6.50).The equation on the inverse Colombeau-Fourier transform

$$(u_\varepsilon(y))_\varepsilon = \int\limits_{\mathbb{R}^d} dx\exp(-iyx)(\widetilde{u}_\varepsilon(x))_\varepsilon$$

(0.6.53)

of $(\widetilde{u}_\varepsilon(x))_\varepsilon$ has the form

$$\left(\frac{\partial u_\varepsilon}{\partial t}\right)_\varepsilon = -G(y^2)^\alpha + i((A_\varepsilon, y)u_\varepsilon)_\varepsilon + \left(V_\varepsilon\left(\frac{1}{i}\frac{\partial}{\partial y}\right)u_\varepsilon\right)_\varepsilon. \tag{0.6.54}$$

**Theorem 0.6.15.**The solution to the Colombeau-Cauchy problem of equation (0.6.54) can be written in the form of a complex Colombeau-Feynman-Kac formula

$$(u_\varepsilon(t,u))_\varepsilon = \left(\int \exp\left\{-\int_0^t \left[G(q^2(\tau))^\alpha - (A_\varepsilon, q(\tau))\right]\right\}u_{0,\varepsilon}(q(t))\mathbf{D}_{0,t,\varepsilon}^x[q(\tau)]\right)_\varepsilon, \tag{0.6.55}$$

where $(\mathbf{D}_{0,t,\varepsilon}^x[q(\tau)])_\varepsilon$ is the measure of the Colombeau generalized jump process corresponding to equation (0.6.51).
As another example, let us consider the case of generalized complex anharmonic oscillator,i.e. the equation

$$\left(\frac{\partial\widetilde{\psi}_\varepsilon}{\partial t}\right)_\varepsilon = \frac{1}{2}((G\Delta - x^2 - iV_\varepsilon(x))\widetilde{\psi}_\varepsilon)_\varepsilon, \tag{0.6.56}$$

where $(V_\varepsilon)_\varepsilon = (V_{\mu,\varepsilon})_\varepsilon \in \mathcal{F}^{\blacktriangle}(\mathbb{R}^d)$. The Colombeau-Fourier transform of this equation has the form

$$\left(\frac{\partial\psi_\varepsilon(t,p)}{\partial t}\right) = \frac{1}{2}\left(\left(\Delta - Gp^2 - iV_\varepsilon\left(\frac{1}{i}\frac{\partial}{\partial p}\right)\right)\psi_\varepsilon(t,p)\right)_\varepsilon. \tag{0.6.57}$$

**Theorem 0.6.16.**If $\operatorname{Re}G \geq 0$, the Cauchy problem of equation (0.6.57) defines a Colombeau regular semigroup of operators in $C_0(\mathbb{R}^d)$, and thus can be presented as the Colombeau path integral from .
**Remark 0.6.10.**The statement of the Proposition can be generalised easily to the following situation,which includes all Schrödinger equations, namely to the case of the equation

$$\left(\frac{\partial \psi_\varepsilon}{\partial t}\right)_\varepsilon = i((A_\varepsilon - B_\varepsilon)\psi_\varepsilon)_\varepsilon. \tag{0.6.58}$$

where $(A_\varepsilon)_\varepsilon$ is Colombeau selfadjoint operator (i.e. $\forall \varepsilon \in (0, 1]$, $A_\varepsilon$ is selfadjoint operator), for which therefore exists (according to spectral theory) a unitary transformations $U_\varepsilon$ such that $\forall \varepsilon \in (0, 1], U_\varepsilon A_\varepsilon U_\varepsilon^{-1}$ is the multiplication operator in some

$\mathscr{L}_2(X, d\mu)$, where $X$ is locally compact, and $(B_\varepsilon)_\varepsilon$ is such that $(U_\varepsilon B_\varepsilon U_\varepsilon^{-1})_\varepsilon$ is a bounded operator in $C_0(X)$.

We shall consider now Colombeau-Schrödinger type operator with a magnetic field, namely the operator

$$(H_\varepsilon)_\varepsilon = \left(\left(i\frac{\partial}{\partial x} + A_\varepsilon(x)\right)^2 + V_\varepsilon(x)\right)_\varepsilon \tag{0.6.59}$$

where $x = (x_1, \ldots, x_d)$ and $\frac{\partial}{\partial x} = \left(\frac{\partial}{\partial x_1}, \ldots, \frac{\partial}{\partial x_d}\right)$ is the gradient operator in $\mathbb{R}^d$, under the following conditions :

(**C**[*]1) the generalized magnetic vector-potential $\mathbf{A}_\varepsilon = (A_\varepsilon^1, \ldots, A_\varepsilon^d)$ is a bounded measurable $\mathbb{R}^d$-valued Colombeau generalized function on $\mathbb{R}^d$,

(**C**[*]2) the potential $(V_\varepsilon)_\varepsilon$ and the divergence $\mathbf{div}(\mathbf{A}_\varepsilon)_\varepsilon = \sum_{j=1}^{d}\left(\frac{\partial A_\varepsilon^j}{\partial x_j}\right)_\varepsilon$ of $(\mathbf{A}_\varepsilon)_\varepsilon$,

defined in the sense of distributions, are both Colombeau-Borel measures,

(**C**[*]3) if $d > 1$ there exist $\alpha > d - 2$ and $C > 0$ such that for all $x \in \mathbb{R}^d$ and $r \in (0, 1]$

$$|\mathbf{div}(\mathbf{A}_\varepsilon)_\varepsilon|(B_r(x)) \le C(r + \varepsilon)^{-\alpha}, |(V_r)_\varepsilon|(B_r(x)) \le C(r + \varepsilon)^{-\alpha} \tag{0.6.60}$$

where $|\mathbf{div}(\mathbf{A}_\varepsilon)_\varepsilon|(\cdot)$ and $|V|(\cdot)$ denote the total variations of the (possibly non-positive) Colombeau-Borel measures $(V_\varepsilon)_\varepsilon$ and $\mathbf{div}\mathbf{A}$ respectively, and $B_r(x)$ denotes the ball in $\mathbb{R}^d$ of the radius $r$ centered at $x$; if $d = 1$, then the same holds for $\alpha = 0$.

Let us consider the following more general Cauchy problem

$$\left(\frac{\partial u_\varepsilon(t, x)}{\partial t}\right)_\varepsilon = -D(H_\varepsilon u_\varepsilon(t, x))_\varepsilon, (u_\varepsilon(0, x))_\varepsilon = (u_{0,\varepsilon}(x))_\varepsilon, \tag{0.6.61}$$

where the (generalised) diffusion coefficient $D$ is an arbitrary complex number such that $\epsilon = \operatorname{Re} D \geq 0, |D| > 0$. In the interaction representation, equation (0.6.61) takes the form

$$(u_\varepsilon(t,x))_\varepsilon = [\exp(-Dt\Delta/2)]((u_{0,\varepsilon}(x))_\varepsilon -$$

$$D\left( \int_0^t [\exp(-D(t-s)\Delta/2)](W_\varepsilon(x) - 2i(\mathbf{A}_\varepsilon(x),\nabla))u_\varepsilon(s,x)ds \right)_\varepsilon . \qquad (0.6.62)$$

Here $(W_\varepsilon(x))_\varepsilon = (V_\varepsilon(x))_\varepsilon + \left( |\mathbf{A}_\varepsilon(x)|^2 \right)_\varepsilon - i(\mathbf{div}\mathbf{A}_\varepsilon(x))_\varepsilon$. Equation (0.6.62) can be formally solved by iterations. In the case of the Green function $(G_\varepsilon^D(t,x,y))_\varepsilon$ of equation (0.6.62), i.e. its solution with the Dirac initial condition $(u_{0,\varepsilon}(x))_\varepsilon = \delta(x-y)$, the iteration procedure leads to the following representation:

$$(G_\varepsilon^D(t,x,y))_\varepsilon = \left( \sum_{k=0}^\infty I_{k,\varepsilon}^D(t,x,y) \right)_\varepsilon ,$$

$$(I_{k,\varepsilon}^D(t,x,y))_\varepsilon =$$

$$-D\int_0^t \int_{\mathbb{R}^d} (I_{k-1,\varepsilon}^D(t-s,x,\xi)(W_\varepsilon(d\xi) + 2iD^{-1}s^{-1}(\mathbf{A}_\varepsilon(x),(\xi-y))d\xi))_\varepsilon G_{free}^D(s,\xi-y)ds, \qquad (0.6.63)$$

$$G_{free}^D(s,\xi-y) = (2\pi tD)^{-d/2} \exp\left[ -\frac{(\xi-y)^2}{2Dt} \right].$$

**Theorem 0.6.17.(i)** If $\epsilon = \operatorname{Re} D > 0$, then all terms of series (0.6.35) are well defined as absolutely convergent Colombeau integrals, the series itself is absolutely convergent and its sum $(G_\varepsilon^D(t,x,y))_\varepsilon$ is continuous in $x,y \in \mathbb{R}^d$, $t > 0$ (and $D$) and satisfies the following estimate

$$(G_\varepsilon^D(t,x,y))_\varepsilon \leq ((K_\varepsilon)_\varepsilon) G_{free}^{\frac{|D|^2}{\epsilon}}(t,x-y) \exp(((B_\varepsilon)_\varepsilon)|x-y|) \qquad (0.6.64)$$

uniformly for $t \leq t_0$ with any fixed $t_0$, where $\mathbf{cl}[(B_\varepsilon)_\varepsilon], \mathbf{cl}[(K_\varepsilon)_\varepsilon] \in \widetilde{\mathbb{R}}$ are Colombeau constants.

**(ii)** The integral operators

$$([U_\varepsilon^D(t)u_{0,\varepsilon}](t,x))_\varepsilon = \left(\int u_\varepsilon^D(t,x,y)u_{0,\varepsilon}(y)dy\right)_\varepsilon \qquad (0.6.65)$$

defining the solutions to equation (0.6.62) for $t \in [0,t_0]$ form a uniformly bounded family of Colombeau operators $\mathcal{L}_2(\mathbb{R}^d) \to \mathcal{L}_2(\mathbb{R}^d)$.

**(iii)** If $D$ is real, i.e. $D = \epsilon > 0$, then there exists a constant $\omega > 0$ such that the Green function $(G_\varepsilon^\epsilon(t,x,y))_\varepsilon$ has the asymptotic representation

$$(G_\varepsilon^\epsilon(t,x,y))_\varepsilon = G_{free}^\epsilon(t,x,y)(1 + (O_\varepsilon(t^\omega))_\varepsilon + (O_\varepsilon(|x-y|))_\varepsilon) \qquad (0.6.66)$$

for small $t$ and $x - y$. In case of vanishing $\mathbf{A}$, the multiplier $\exp(B|x-y|)$ in (0.6.64) can be omitted, and the term $(O_\varepsilon(|x-y|))_\varepsilon$ in (0.6.66) can be dropped. In this case, formula (0.6.66) gives global (uniform for all $x,y$) small time asymptotics for $(G_\varepsilon^\epsilon(t,x,y))_\varepsilon$.

**(iv)** One can give rigorous meaning to formal expression (0.6.63) as a bounded below Colombeau self adjoint operator in such a way that the family (0.6.65) of operators $U^D(t)$ (giving solutions to equation (0.6.62), which is formally equivalent to the evolutionary equation (0.6.61) with the formal generator $D(H_\varepsilon)_\varepsilon$) coincides with the semigroup $(\exp(-tDH_\varepsilon))_\varepsilon$ defined by means of the functional calculus. Hence for the integral kernel of the operators $(\exp(-tDH_\varepsilon))_\varepsilon$ the estimates (0.6.65) and (0.6.66) hold.

The rigorous construction a path integral representation for the Colombeau solution $(U_\varepsilon^D(t)u_{0,\varepsilon})_\varepsilon$ of Eq.(0.6.61), gives by construction an Colombeau-Borel measure on a path space that is supported on the set of continuous piecewise linear paths. Denote this set by $\mathbf{CPL}^{\mathbb{R}^d}$. Let $\mathbf{CPL}^{x,y}(0,t)$ denote the class of paths $q : [0,t] \to \mathbb{R}^d$ from $\mathbf{CPL}^{\mathbb{R}^d}$ joining $x$ and $y$ in time $t$, i.e. such that $q(0) = x, q(t) = y$. By $\mathbf{CPL}_n^{x,y}(0,t)$ we denote its subclass consisting of all paths from $\mathbf{CPL}^{x,y}(0,t)$ that have exactly $n$ jumps of their derivative. Clearly, each $\mathbf{CPL}_n^{x,y}(0,t)$ is parametrised by the simplex

$$\mathbf{sim}_t^n = \{s_1,\ldots,s_n | 0 < s_1 < s_2 < \ldots < s_n < t\} \qquad (0.6.67)$$

of the times of jumps $s_1, \ldots, s_n$ of the derivatives of a path and by $n$ positions $q(s_j)$, $j = 1, \ldots, n$, of this path at these points. In other words, an arbitrary path in $\mathbf{CPL}_n^{x,y}(0,t)$
has the form

$$q_n(s) = q_{\eta_1,\ldots,\eta_n}^{s_1,\ldots,s_n} = \eta_j + (s - s_j)\frac{\eta_{j+1} - \eta_j}{s_{j+1} - s_j},$$

$$(0.6.68)$$

$$s_j \in [s_j, s_{j+1}],$$

where $s_0 = 0, s_{n+1} = t, \eta_0 = x, \eta_{n+1} = y$. Therefore $\mathbf{CPL}^{x,y}(0,t) = \cup_{n=0}^{\infty} \mathbf{CPL}_n^{x,y}(0,t)$. To any $\mathbb{R}^{d+1}$-valued Colombeau-Borel measure $\left(\widetilde{M}_\varepsilon\right)_\varepsilon = \{(\mu_\varepsilon)_\varepsilon, (\nu_\varepsilon)_\varepsilon\} = \{(\mu_\varepsilon)_\varepsilon, (\nu_\varepsilon^1)_\varepsilon, \ldots, (\nu_\varepsilon^n)_\varepsilon\}$ on $\mathbb{R}^{d+1}$ there corresponds a $\sigma$-finite complex measure $\mathbf{D}^{\mathbf{CPL}}[q(\tau)]$ on $\mathbf{CPL}^{x,y}(0,t)$, which is defined as the sum of the Colombeau-Borel measures $(\mathbf{M}_{n,\varepsilon}^{\mathbf{CPL}})_\varepsilon$ $n \in \mathbb{N}$ on the finite-dimensional spaces $\mathbf{CPL}_n^{x,y}(0,t)$ such that $\mathbf{M}_0^{\mathbf{CPL}}$ is just the unit measure on the one-point set $\mathbf{CPL}_0^{x,y}(0,t)$ and each $\mathbf{M}_n^{\mathbf{CPL}}, n > 0$, is defined in the following way: if $\Phi(\cdot)$ is a functional on $\mathbf{CPL}^{x,y}(0,t)$, then:

$$\left(\int_{\mathbf{CPL}_n^{x,y}(0,t)} \Phi(q(\tau))\mathbf{D}_{n,\varepsilon}^{\mathbf{CPL}}[q(\tau)]\right)_\varepsilon = \left(\int_{\mathbf{sim}_t^n} ds_1 \ldots ds_n \times\right.$$

$$(0.6.69)$$

$$\left.\int_{\mathbb{R}^d} \ldots \int_{\mathbb{R}^d} \left(\mu_\varepsilon + 2i\frac{(\eta_2 - \eta_1, \nu_\varepsilon)}{s_2 - s_1}\right)(d\eta_1)\ldots\left(\mu + 2i\frac{(y - \eta_n, \nu_\varepsilon)}{t - s_n}\right)(d\eta_n)\Phi(q(\tau))\right)_\varepsilon.$$

Now, we let $\left(\widetilde{M}_\varepsilon\right)_\varepsilon = \{(W_\varepsilon(x))_\varepsilon, D^{-1}(\mathbf{A}_\varepsilon(x)dx)_\varepsilon\}$.

**Theorem 0.6.10.** For any $D$ with $\epsilon = \operatorname{Re} D > 0$, the Green function $G_\varepsilon^D(t,x,y)$ of equation
(0.6.62) has the following path integral representation:

$$G^D(t,x,y) = \int\limits_{\mathbf{CPL}^{x,y}(0,t)} \Phi_D(q(\tau)) \exp\left\{\frac{1}{2D}\int_0^t \dot{q}^2(s)ds\right\} \mathbf{D^{CPL}}[q(\tau)] \qquad (0.6.70)$$

with $q(s)$ given by Eq.(0.6.68) and

$$\Phi_D(q(\tau)) = D^n \prod_{j=1}^{n+1}(2\pi D(s_j - s_{j-1}))^{-d/2}. \qquad (0.6.71)$$

For any $(u_{0,\varepsilon})_\varepsilon \in G_{\mathcal{L}_2(\mathbb{R}^d)}$ the Colombeau solution $(u_\varepsilon(t,x))_\varepsilon$ of the Colombeau-Cauchy problem (0.6.61) with $D = i$ has the form

$$(u_\varepsilon(t,x))_\varepsilon =$$

$$\left(\lim_{\delta\to 0}\int\limits_{\mathbf{CPL}^{x,y}(0,t)}\int_{\mathbb{R}^d} u_{0,\varepsilon}(x)\Phi_{i+\delta}(q(\tau))\exp\left\{\frac{1}{2(i+\delta)}\int_0^t \dot{q}^2(s)ds\right\}\mathbf{D_\varepsilon^{CPL}}[q(\tau)]\right)_\varepsilon. \qquad (0.6.72)$$

Here the limit is understood in $\mathcal{L}_2(\mathbb{R}^d)$ sence.

**Theorem 0.6.18.** Suppose the vector potential A vanishes and V satisfies assumption (**C**3). Suppose additionally that $(V_\varepsilon)_\varepsilon$ has no atom at the origin and is a finite positive measure so that $(\mu_{V_\varepsilon})_\varepsilon = (V_\varepsilon(\mathbb{R}^d))_\varepsilon > 0$. Let the paths of **CPL** are parametrised by Eq.(0.6.40). Let $(\mathbf{E}_{t,\varepsilon}(\bullet))_\varepsilon$ denote the expectation with respect to the generalized process of jumps $\eta_j$ which are identically independently distributed according to the probability Colombeau measure $(V_\varepsilon/\mu_{V_\varepsilon})_\varepsilon$, and which occur at times $s_j$ from $[0,t]$ distributed according to Colombeau generalized Poisson process of intensity $(\mu_{V_\varepsilon})_\varepsilon$. Then the integral (0.6.70) can be written in the form

$$(G_\varepsilon^D(t,x,y))_\varepsilon = \left(e^{\mu_{V_\varepsilon}}\mathbf{E}_{t,\varepsilon}\left[\Phi_D(q(\tau))\exp\left\{-\frac{1}{2D}\int_0^t \dot{q}^2(s)ds\right\}\mathbf{D_\varepsilon^{CPL}}[q(\tau)]\right]\right)_\varepsilon. \qquad (0.6.73)$$

$$(0.6.)$$

$$(0.6.)$$

## I.0.7.Colombeau-Wiener Path Integral and generalized Haba theorem.

Let us consider the heat equation with $D \geq 0$ and complex term

$Z(x) = W(x) + iV(x).$

$$\frac{\partial u(t,x)}{\partial t} = D\Delta u(t,x) + Z(x)u(t,x), u(0) = u_0(x), x \in \mathbb{R}^n. \qquad (0.7.1)$$

Nelson [40] has shown how to solve such equations by means of Wiener integrals for the purely imaginary $Z(x) = iV(x)$ under following assumption (**N**) on the function $V$

(**N**): There is a closed set $F$ in $\mathbb{R}^d$ of capacity $0$ such that $V$ is continuous and real on the complement of $F$.

For example, if $d > 1$ and $V$ is a function which is continuous except at one point, say the origin (with an arbitrary singularity at the origin and arbitrary growth at infinity), then $V$ satisfies (**N**).

Let $x$ be any point in the complement of $F$ and $\mathbf{Pr}_x$ the corresponding Wiener measure, with diffusion constant $D > 0$, on the space $\Omega$ of all trajectories. Then for almost every $\omega \in \Omega$, $V(\omega(t))$ is a continuous function of $t$ for $0 \leq t < \infty$, and so one may form the Wiener integral

$$[U^t_{m,V}u_0](x) = \int_\Omega \exp\left[-i\int_0^t V(\omega(s))ds\right] u_0(\omega(t))\mathbf{Pr}_x(d\omega) \qquad (0.7.2)$$

for any $u_0 \in \mathcal{L}_2$ and any $t \geq 0$. In fact, the integrand is defined for almost every $\omega \in \Omega$, and it is bounded in absolute value by $|u_0(\omega(t))|$. By definition one obtain

$$\int_\Omega |u_0(\omega(t))|\mathbf{Pr}_x(d\omega) = \Xi^t_m|u_0|(x) < \infty. \qquad (0.7.3)$$

Since $V(\omega(t))$ is continuous for almost every $\omega \in \Omega$, one obtain

$$\int_0^t V(\omega(s))ds = \lim_{n \to \infty} \sum_{j=1}^n \frac{t}{n} V\left(\omega\left(j\frac{t}{n}\right)\right). \qquad (0.7.4)$$

Therefore by the Lebesgue dominated convergence theorem one obtain $\forall x \notin F$ :

$$[U^t_{m,V}u_0](x) = \lim_{n \to \infty} \int_\Omega \exp\left[-i\sum_{j=1}^n \frac{t}{n} V\left(\omega\left(j\frac{t}{n}\right)\right)\right] u_0(\omega(t))\mathbf{Pr}_x(d\omega) \qquad (0.7.5)$$

and by definition of the Wiener integral Eq.(0.7.5) gives

$$[U_{m,V}^t u_0](x) = \left( \Xi_m^{\frac{t}{n}} M_V^{\frac{t}{n}} \right)^n u_0(x).$$ (0.7.6)

**Theorem 0.7.1.**(Nelson [41]). Let $D > 0$, $Z$ be purely imaginary and let $V$ satisfy condition (**N**).Then for all $u_0(x) \in \mathcal{L}_2(\mathbb{R}^n)$ and $t \geq 0$, the limit (0.4.6) exists, where $\Xi_m^t$ and $M_V^t$ are defined by (0.3.4) and (0.3.5), respectively.The operators $U_{m,V}^t$ defined (0.7.6) have the semigroup property

$$U_{m,V}^t U_{m,V}^s = U_{m,V}^{t+s}, 0 \leq t, s < \infty$$ (0.7.7)

and for all $u_0(x) \in \mathcal{L}_2(\mathbb{R}^n), t \to U_{m,V}^t$ is continuous from $[0,\infty)$ to $\mathcal{L}_2(\mathbb{R}^n)$. Let us consider now the heat equation with $D \geq 0$ and generalized Colombeau complex term $(Z_\varepsilon(x))_\varepsilon = (W_\varepsilon(x))_\varepsilon + i(V_\varepsilon(x))_\varepsilon$ such that $\mathbf{cl}[(W_\varepsilon(x))_\varepsilon], \mathbf{cl}[(V_\varepsilon(x))_\varepsilon],$ $\mathbf{cl}[(u_{0,\varepsilon}(x))_\varepsilon] \in G(\mathbb{R}^n)$

$$\frac{\partial}{\partial t}((u_\varepsilon(t,x))_\varepsilon) = D\Delta((u_\varepsilon(t,x))_\varepsilon) + ((Z_\varepsilon(x))_\varepsilon)((u_\varepsilon(t,x))_\varepsilon),$$ (0.7.8)

$$(u_\varepsilon(0))_\varepsilon = (u_{0,\varepsilon}(x))_\varepsilon, x \in \mathbb{R}^n.$$

One can solve such equations by means of Colombeau-Wiener integrals for the $D > 0$ under following assumption ($\widetilde{\mathbf{N}}$) on the Colombeau generalized function $(Z_\varepsilon)_\varepsilon = (W_\varepsilon(x))_\varepsilon + i(V_\varepsilon(x))_\varepsilon$

($\widetilde{\mathbf{N}}$) : $\mathbf{cl}[(W_\varepsilon(x))_\varepsilon], \mathbf{cl}[(V_\varepsilon(x))_\varepsilon] \in G(\mathbb{R}^n)$. There is a closed set $F$ in $\mathbb{R}^d$ of capacity $0$ such that $\forall \varepsilon \in (0,1]$, $Z_\varepsilon(x)$ is continuous on the complement of $F$.
Let $x$ be any point in the complement of $F$ and $\mathbf{Pr}_x$ the corresponding Wiener measure, with diffusion constant $D > 0$, on the space $\Omega$ of all trajectories. Then for almost every $\omega \in \Omega$, $Z_\varepsilon(\omega(t)), \varepsilon \in (0,1]$ is a continuous function of $t$ for $0 \leq t < \infty$, and so one may form the Colombeau-Wiener integral:

$$([U^t_{m,Z_\varepsilon} u_0](x))_\varepsilon =$$

$$\left( \int\limits_\Omega \exp\left[ -i \int\limits_0^t Z_\varepsilon(\omega(s)) ds \right] u_{0,\varepsilon}(\omega(t)) \mathbf{Pr}_x(d\omega) \right)_\varepsilon = \tag{0.7.9}$$

$$\int\limits_\Omega \exp\left[ -i \int\limits_0^t ((Z_\varepsilon(\omega(s)))_\varepsilon) ds \right] ((u_{0,\varepsilon}(\omega(t)))_\varepsilon) \mathbf{Pr}_x(d\omega)$$

for any $(u_{0,\varepsilon})_\varepsilon$ such that $\forall \varepsilon \in (0,1]$, $u_{0,\varepsilon} \in \mathcal{L}_2$ and any $t \geq 0$. In fact, the integrand is defined for almost every $\omega \in \Omega$, and it is bounded in absolute value by $|u_{0,\varepsilon}(\omega(t))| \exp(\varepsilon^{-N})$ for some $N \in \mathbb{N}$. By definition one obtain $\forall \varepsilon \in (0,1]$ :

$$\int\limits_\Omega |u_{0,\varepsilon}(\omega(t))| \mathbf{Pr}_x(d\omega) = \Xi^t_m |u_{0,\varepsilon}|(x) < \infty. \tag{0.7.10}$$

Since $\forall \varepsilon \in (0,1]$, $Z_\varepsilon(\omega(t))$ is continuous for almost every $\omega \in \Omega$, one obtain

$$\left( \int\limits_0^t Z_\varepsilon(\omega(s)) ds \right)_\varepsilon = \left( \lim_{n\to\infty} \sum_{j=1}^n \frac{t}{n} Z_\varepsilon\left(\omega\left(j\frac{t}{n}\right)\right) \right)_\varepsilon. \tag{0.7.11}$$

Therefore by the Lebesgue dominated convergence theorem one obtain $\forall \varepsilon \in (0,1]$, $\forall x \notin F$ :

$$([U^t_{m,Z_\varepsilon} u_0](x))_\varepsilon =$$

$$\left( \lim_{n\to\infty} \int\limits_\Omega \exp\left[ -i \sum_{j=1}^n \frac{t}{n} Z_\varepsilon\left(\omega\left(j\frac{t}{n}\right)\right) \right] u_{0,\varepsilon}(\omega(t)) \mathbf{Pr}_x(d\omega) \right)_\varepsilon \tag{0.7.12}$$

and by definition of the Wiener integral Eq.(0.7.13) gives

$$\left( [U_{m,Z_\varepsilon}^t u_{0,\varepsilon}](x) \right)_\varepsilon = \left( \left( \Xi_m^{\frac{t}{n}} M_{Z_\varepsilon}^{\frac{t}{n}} \right)^n u_{0,\varepsilon}(x) \right)_\varepsilon . \tag{0.7.13}$$

Here

$$\left( [\Xi_m^t u_{0,\varepsilon}](x) \right)_\varepsilon = \left( \frac{2\pi i t}{m} \right)^{-d/2} \left( \int_{\mathbb{R}^n} d^d y\, u_{0,\varepsilon}(y) \exp\left[ i\frac{m}{2} \left( \frac{|x-y|^2}{t} \right) \right] \right)_\varepsilon , \tag{0.7.14}$$

$$\left( M_{Z_\varepsilon}^t \right)_\varepsilon = \left( \exp\{-t Z_\varepsilon\} \right)_\varepsilon .$$

Haba [44] has found conditions on $u_0(x)$ and $Z(x)$ under which Eq.(0.7.2) valid for complex $D \in \mathbb{C}$.

**Definition 0.7.1**. Let $\mathcal{F}_d$ denote the complex Banach algebra of functions which can be represented of the form

$$g(x) = \int_{\mathbb{R}^d} \exp(ipx) d\mu_g(p). \tag{0.7.15}$$

where $\mu_g(p)$ is a complex measure on $\mathbb{R}^d$ with a bounded variation $|\mu_g|(p)$.

**Definition 0.7.2**. Let $\mathcal{F}_d^\spadesuit$ denote the set of functions which can be represented of the form (0.7.15) with $\mu_g(p)$ such that

$$\forall \epsilon(\epsilon > 0) \exists K(\epsilon) \forall a(a \in \mathbb{R}) \left[ \int_{\mathbb{R}^d} \exp(a|p|) d|\mu_g|(p) \leq K(\epsilon) \exp(\epsilon a^2) \right]. \tag{0.7.16}$$

Note that $\mathcal{F}_d$ and $\mathcal{F}_d^\spadesuit$ are dense subsets in $\mathcal{L}_2(\mathbb{R}^d)$.

Let us consider now a complex extension of the Schrödinger equation (0.1.1) in the form

$$\frac{\partial \psi(t,x)}{\partial t} = \frac{1}{2}\sigma^2\lambda^2\Delta\psi(t,x) + g\hbar^{-1}\lambda^{-2}Z(x)\psi(t,x), \psi(0,x) = \psi_0(x), x \in \mathbb{R}^d. \qquad (0.7.17)$$

Here $\lambda \in \mathbb{C}, \sigma = \sqrt{\hbar m^{-1}}$. Schrödinger equation corresponds to $\lambda = \sqrt{i} = (1+i)/\sqrt{2}$.

Haba [44] has express the Feynman path integral corresponding to Eq.(0.7.17) using Wiener measure corresponding to Wiener process $\mathbf{b}(\tau) \triangleq \mathbf{b}_\tau \in \mathbb{R}^d$,(Brownian motion). Here $\mathbf{b}_\tau \in \mathbb{R}^d$ is the Gaussian process with the covariance

$$\mathbf{E}[b^k(\tau_1)b^r(\tau_2)] = \delta_{kr}\min(\tau_1,\tau_2). \qquad (0.7.18)$$
$$k,r = 1,\ldots,d.$$

We use the notation

$$\mathbf{E}[F] \triangleq \int_{\mathbf{\Omega}_t} d\mu_t^W(\omega)F(\omega), \qquad (0.7.19)$$

where $\mu_t^W(\omega)$ is the Wiener measure and the integral is over the Wiener space $\mathbf{\Omega}_t = \mathbf{C}([0,t])$ which consists of continuous functions defined on the interval $[0,t]$ with $\omega(0) = 0$.

**Lemma 0.7.1.**[44]. Let $\mu$ be a Gaussian measure on the Wiener space $\mathbf{\Omega}_t = \mathbf{C}([0,t])$.
Assume that $Q$ is a non-negative quadratic form on $\mathbf{\Omega}_t$ such that $Q = \lim_{m\to\infty} Q_m$ with probability 1, where $Q_m$ is a continuous non-negative quadratic form. Let $\mathcal{L}(\omega) = \sum_{k=1}^r a_k\omega(\tau_k)$, where $a_k \in \mathbb{R}$ $(k = 1,2,\ldots,r)$ and $0 < \tau_k < t$ are arbitrary. Then

$$\left|\int d\mu(\omega)\exp[-Q(\omega) - i\lambda\mathcal{L}(\omega)]\right| \leq 1 \qquad (0.7.20)$$

if $\mathrm{Re}(\lambda^2) \geq 0$.

Let us denote

$$\theta_R(x) = \exp\left(\frac{|x^2|}{2R}\right),$$

$$\Theta_t[g\theta_R(\omega(\tau)); Z(x + \lambda\sigma\omega(\tau))] = \qquad (0.7.21)$$

$$\exp\left\{\frac{g}{\hbar\lambda^2} \int_0^t d\tau[\theta_R(\omega(\tau))Z(x + \lambda\sigma\omega(\tau))]\right\}.$$

We shall also use a shorthand notation $\Theta_t[g\theta_R; Z]$ for the quantity $\Theta_t[g\theta_R(\omega(\tau)); Z(x + \lambda\sigma\omega(\tau))]$.

**Theorem 0.7.2.**(Haba [44]). Assume that $\psi_0, Z \in \mathcal{F}_d^{\blacktriangle}$. Then, $\Theta_t[g\theta_R; Z]$ is an analytic function of $g \in \mathbb{C}$ and $\lambda \neq 0$ as long as $\operatorname{Re}\lambda^2 \geq 0$. Define

$$\psi^R(t,x) = \mathbf{E}[\Theta_t[g\theta_R(\omega(\tau)); Z(x + \lambda\sigma\omega(\tau))]\psi_0(x + \lambda\sigma\omega(\tau))] \qquad (0.7.22)$$

Then, limit

$$\psi(t,x) = \lim_{R\to\infty} \psi^R(t,x) \qquad (0.7.23)$$

exists uniformly in $t$ and $x$. Function $\psi(t,x)$ solves the Schrödinger equation (0.7.17). The solution is an analytic function of $g$ and $\lambda$ in the same region where $\mathbf{E}[\Theta_t]$ is analytic.

Let $Z(t,x)$ be the function

$$Z(t,x) = \frac{1}{2}\varpi(t)x^2 + \vartheta(t)x + gZ(x), \qquad (0.7.24)$$

where $Z(x) \in \mathcal{F}_d^{\blacktriangle}$ and $\varpi(t), \vartheta(t)$ are continuous functions of time.

**Theorem 0.7.3.**(Haba [44]). Assume that $\varpi(t) \leq \varpi_0$ for each $t$. If $t < \pi/\varpi_0$ then theorem 0.7.2 holds true for potential given by Eq.(0.7.24).

Let us consider a decomposition of $\mathbb{R}^{2n}$ into $\mathbb{R}^n \bigoplus \mathbb{R}^n$ where

$z = (x_1, x_2), x_1, x_2 \in \mathbb{R}^n.$ Let
$\mathbf{G}_z$ be a subset $\mathbf{G}_z \subset \mathbb{R}^{2n}$

$$\mathbf{G}_z = \{w \in \mathbb{R}^{2n} | w = (x_1 + y, x_2 + y), y \in \mathbb{R}^n\} \qquad (0.7.25)$$

and let $Z(x_1, x_2)$ be a complex-valued potential defined on $\mathbb{R}^{2n}$ such that on $\mathbf{G}_z$

$$|\operatorname{Im} Z(x_1 + y, x_2 + y)| \leq A(x_1, x_2)(1 + |y|) \qquad (0.7.26)$$

with a certain continuous non-negative function $A.$ We define a semigroup $\mathbf{S}_t$ on the set of functions defined on $\mathbb{R}^{2n}$ such that on $\mathbf{G}_z$

$$|\Phi(x_1 + y, x_2 + y)| \leq M(x_1, x_2)\exp(|y|C(x_1, x_2)) \qquad (0.7.27)$$

with certain non-negative continuous functions $C(x_1, x_2)$ and $M(x_1, x_2)$. Let $\mathbf{S}_t,$ is defined by the Feynman-Kac formula:

$$(\mathbf{S}_t\Phi)(x_1,x_2) =$$

$$\mathbf{E}\Bigg[\exp\Bigg(-\frac{i}{\hbar}\int\limits_0^t d\tau Z\Big(x_1+\frac{\sigma\omega(\tau)}{\sqrt{2}},x_2+\frac{\sigma\omega(\tau)}{\sqrt{2}}\Big)\times$$

$$\Phi\Big(x_1+\frac{\sigma\omega(t)}{\sqrt{2}},x_2+\frac{\sigma\omega(t)}{\sqrt{2}}\Big)\Bigg)\Bigg] =$$

$$\int\limits_{\mathbb{R}^n}\int\limits_{\mathbb{R}^n}\mathbf{E}\Bigg[\delta\Big(y_1-x_1-\frac{\sigma\omega(\tau)}{\sqrt{2}}\Big)\delta\Big(y_2-x_2-\frac{\sigma\omega(\tau)}{\sqrt{2}}\Big)\times \qquad (0.7.28)$$

$$\exp\Bigg(-\frac{i}{\hbar}\int\limits_0^t d\tau \widetilde{Z}\Big(x_1+\frac{\sigma\omega(\tau)}{\sqrt{2}},x_2+\frac{\sigma\omega(\tau)}{\sqrt{2}}\Big)\Bigg)\Bigg]\times$$

$$\Phi(y_1,y_2)dy_1dy_2 =$$

$$\int\limits_{\mathbb{R}^n}\int\limits_{\mathbb{R}^n}S_t(x_1,x_2;y_1,y_2)\Phi(y_1,y_2)dy_1dy_2,$$

where $\sigma=\sqrt{\hbar/m}$. Using Eq.(0.7.28) one can express the kernel $S_t$ $(x_1,x_2;y_1,y_2)$ of $\mathbf{S}_t$, by the Brownian bridge $\boldsymbol{\alpha}(s)$. The Brownian bridge is the Gaussian process defined on the interval $[0,1]$ with the covariance

$$\mathbf{E}[\alpha_k(s)\alpha_r(s')] = \delta_{kr}(1-s') \text{ iff } s \leq s' \qquad (0.7.29)$$

and the boundary conditions $\alpha(0) = \alpha(1) = 0.$ Then the expectation value in the second line of Eq.(0.7.28) is equal to the expectation value

$$\mathbf{E}\left[\delta\left(y_1 - x_1 - \frac{\sigma\omega(\tau)}{\sqrt{2}}\right)\delta\left(y_2 - x_2 - \frac{\sigma\omega(\tau)}{\sqrt{2}}\right)\right] =$$

$$(4\pi t\sigma^2)^{-n/2}\exp\left(-\frac{(x_1 - x_2 - y_1 + y_2)^2}{4t\sigma^2}\right)\delta(x_1 - x_2 - y_1 + y_2)$$

(0.7.30)

times the expectation value of $\exp\left(-(i/h)\int_0^t Z\right),$ where the Brownian motion is constrained to end at $y$. A solution of this constraint gives a simple formula for the kernel

$$S_t(x_1, x_2; y_1, y_2) =$$

$$(4\pi t\sigma^2)^{-n/2}\exp\left(-\frac{(x_1 - x_2 - y_1 + y_2)^2}{4t\sigma^2}\right)\delta(x_1 - x_2 - y_1 + y_2) \times$$

$$\mathbf{E}\left[\exp\left(-\frac{i}{h}\int_0^t d\tau Z\left(x_1\left(1 - \frac{s}{t}\right) + y_1\left(\frac{s}{t}\right) + \sqrt{\frac{t}{2}}\,\sigma\boldsymbol{\alpha}\left(\frac{s}{t}\right)\right),\right.\right.$$

(0.7.31)

$$\left.\left.x_2\left(1 - \frac{s}{t}\right) + y_2\left(\frac{s}{t}\right) + \sqrt{\frac{t}{2}}\,\sigma\boldsymbol{\alpha}\left(\frac{s}{t}\right)\right)ds\right)\right]$$

Now we restrict oursevles to functions of the form

$$Z(x_1, x_2) = \widetilde{Z}(x_1 + ix_2), \Phi(x_1, x_2) = \varphi(x_1 + ix_2),$$

(0.7.32)

where $Z$ and $\varphi$ are holomorphic functions of their arguments. Subsituting equation (0.7.32) into equation (0.7.28) we obtain

$$(\mathbf{S}_t\varphi)(x_1 + ix_2) =$$

$$(4\pi t\sigma^2)^{-n/2} \int_{\mathbb{R}^n} dy\, E_t(y, x_1, x_2)\varphi\left(\frac{1}{2}(y + x_1 + x_2) + \frac{i}{2}(y - x_1 + x_2)\right) \times \qquad (0.7.33)$$

$$\exp\left(-\frac{(y - x_1 - x_2)^2}{4t\sigma^2}\right)$$

where by $E_t$ we have denoted $\mathbf{E}[\dots]$. We note that

$$(4\pi t\sigma^2)^{-n/2} \int_{\mathbb{R}^n} dy\, \chi\left(\frac{1}{2}(y + x) + \frac{i}{2}(y - x)\right)\exp\left(-\frac{(x - y)^2}{4t\sigma^2}\right) =$$

$$(4\pi t\sigma^2)^{-n/2} \int_{\mathbb{R}^n} dy\, \chi\exp\left(-\frac{(x - y)^2}{2it\sigma^2}\right)$$

(where $i^{n/2} = \exp(\frac{1}{2}in\pi)$) for any analytic function $\chi \in \mathscr{L}_1(\mathbb{R}^n)$ of the form $\chi(y) = \int dp\, \chi\tilde{\ }(p)\exp(ipy)$. Finally, we set $x_2 = 0$ in equation (0.7.33) and apply the identity (0.7.34). We obtain the final formula for the Feynman propagator:

$$K_t(x, y) =$$

$$(4\pi t\sigma^2)^{-n/2}\exp\left(-\frac{(x - y)^2}{2it\sigma^2}\right) \times$$

$$\mathbf{E}\left[\exp\left(-\frac{i}{\hbar}\int_0^t d\tau\, \widetilde{Z}\left(x\left(1 - \frac{s}{t}\right) + y\left(\frac{s}{t}\right) + \lambda\sqrt{t}\,\sigma\boldsymbol{\alpha}\left(\frac{s}{t}\right)\right)\right)\right] = \qquad (0.7.35)$$

$$K_t^0(x, y)\Re(t, x, y).$$

Here $\lambda = \sqrt{i}$ and we have denoted the free propagator by $K_t^0(x, y)$.

**Theorem 0.7.4.**(Haba [44]).Define

$$\Re_t^R(t,x,y) = \mathbf{E}\left[\exp\left\{\frac{g}{\hbar\lambda^2}\int\limits_0^t ds\theta_R\left(\alpha\left(\frac{s}{t}\right)\right)Z\left(x\left(1-\frac{s}{t}\right)+y\frac{s}{t}+\lambda\sigma\right)\right\}\right], \quad (0.7.36)$$

where $Z \in \mathcal{F}_n^\spadesuit$. The limit $R \to \infty$ exists uniformly in $t,x,y$. $\Re(t,x,y)$, is an analytic function
of $g$ and $\lambda$ if $\operatorname{Re}\lambda^2 \geq 0$. Formula $K_t(x,y) = K_t^0(x,y)\lim_{R\to\infty}\Re_t^R(t,x,y)$ gives the Feynman propagator (0.7.35) for any $Z \in \mathcal{F}_n^\spadesuit$.

**Definition 0.7.3.**Let $G_{\mathcal{F}_d}(\mathbb{R}^d)$ denote the Colombeau module over $\widetilde{\mathbb{C}}$ of Colombeau generalized functions which can be represented of the form

$$(g_\varepsilon(x))_\varepsilon = \left(\int\limits_{\mathbb{R}^d}\exp(ipx)d\mu_\varepsilon(p;\mathbf{cl}[(g_\varepsilon(x))_\varepsilon])\right)_\varepsilon =$$

$$\int\limits_{\mathbb{R}^d}\exp(ipx)(d\mu_\varepsilon(p;\mathbf{cl}[(g_\varepsilon(x))_\varepsilon]))_\varepsilon. \tag{0.7.37}$$

where $(\mu_\varepsilon(p;\mathbf{cl}[(g_\varepsilon(x))_\varepsilon]))_\varepsilon$ is a complex Colombeau measure on $\mathbb{R}^d$ with a variation
$(|\mu_\varepsilon|(p;\mathbf{cl}[(g_\varepsilon(x))_\varepsilon]))_\varepsilon$.

**Definition 0.7.4.** Let $G_{\mathcal{F}_d^\spadesuit}(\mathbb{R}^n)$ denote Colombeau module over $\widetilde{\mathbb{C}}$ of Colombeau generalized functions which can be represented of the form (0.7.37) with
$\mu_{g,\varepsilon}(p) = (\mu_\varepsilon(p;\mathbf{cl}[(g_\varepsilon(x))_\varepsilon]))_\varepsilon, g = \mathbf{cl}[(g_\varepsilon(x))_\varepsilon]$ such that

$$\forall\epsilon(\epsilon > 0)\exists K(\epsilon)\exists N(N \in \mathbb{N})\forall a(a \in \mathbb{R})$$

$$\forall\varepsilon(\varepsilon \in (0,1])\left[\int\limits_{\mathbb{R}^d}\exp(a|p|)d|\mu_{g,\varepsilon}|(p) \leq \varepsilon^{-N}K(\epsilon)\exp(\epsilon a^2)\right]. \tag{0.7.38}$$

Let us consider now a complex extension of the Colombeau-Schrödinger equation in the form

$$\left(\frac{\partial \psi_\varepsilon(t,x)}{\partial t}\right)_\varepsilon = \frac{1}{2}\sigma^2\lambda^2(\Delta\psi_\varepsilon(t,x))_\varepsilon + g\hbar^{-1}\lambda^{-2}((Z_\varepsilon(x))_\varepsilon)((\psi_\varepsilon(t,x))_\varepsilon),$$

(0.7.39)

$$(\psi_\varepsilon(0,x))_\varepsilon = (\psi_{0,\varepsilon}(x))_\varepsilon, x \in \mathbb{R}^d.$$

Here $\lambda \in \mathbb{C}, \sigma = \sqrt{\hbar m^{-1}}$. Schrödinger equation corresponds to $\lambda = \sqrt{i} = (1+i)/\sqrt{2}$.

We wish to express the Colombeau-Feynman path integral corresponding to Eq.(0.7.39) using Wiener measure corresponding to Wiener process $\mathbf{b}(\tau) \triangleq \mathbf{b}_\tau \in \mathbb{R}^d$, Here $\mathbf{b}_\tau \in \mathbb{R}^d$ is the Gaussian process with the covariance given by Eq.(0.7.18).We use the notation

$$(\mathbf{E}_\varepsilon[F_\varepsilon])_\varepsilon \triangleq \left(\int_{\Omega_t} d\mu_{t,\varepsilon}^W(\omega)F_\varepsilon(\omega)\right)_\varepsilon,$$

(0.7.40)

where $\mu_{t,\varepsilon}^W(\omega)$ is the Colombeau-Wiener measure and the Colombeau integral is over the Wiener space $\Omega_t = \mathbf{C}([0,t])$ which consists of continuous functions defined on the interval $[0,t]$ with $\omega(0) = 0$.

Let us denote

$$\theta_R(x) = \exp\left( \frac{|x^2|}{2R} \right),$$

$$\left( \Theta_{t,\varepsilon}[g\theta_R(\omega(\tau)); Z_\varepsilon(x + \lambda\sigma\omega(\tau))] \right)_\varepsilon = \tag{0.7.41}$$

$$\exp\left\{ \frac{g}{\hbar\lambda^2} \left( \int\limits_0^t d\tau [\theta_R(\omega(\tau)) Z_\varepsilon(x + \lambda\sigma\omega(\tau))] \right)_\varepsilon \right\}.$$

We shall also use a shorthand notation $\Theta_{t,\varepsilon}[g\theta_R; Z]$ for the Colombeau quantity $\left( \Theta_{t,\varepsilon}[g\theta_R(\omega(\tau)); Z_\varepsilon(x + \lambda\sigma\omega(\tau))] \right)_\varepsilon$.

**Theorem 0.7.5.** Assume that $\mathbf{cl}[(\psi_{0,\varepsilon})_\varepsilon], \mathbf{cl}[(Z_\varepsilon)_\varepsilon] \in G_{\mathscr{F}_d^\blacktriangle}(\mathbb{R}^d)$. Then, $\forall \varepsilon \in (0,1]$, $\Theta_{t,\varepsilon}[g\theta_R; Z_\varepsilon]$ is an analytic function of $g \in \mathbb{C}$ and $\lambda \neq 0$ as long as $\operatorname{Re}\lambda^2 \geq 0$. Define

$$\left( \psi_\varepsilon^R(t,x) \right)_\varepsilon = \left( \mathbf{E}_\varepsilon[\Theta_{t,\varepsilon}[g\theta_R(\omega(\tau)); Z_\varepsilon(x + \lambda\sigma\omega(\tau))]\psi_{0,\varepsilon}(x + \lambda\sigma\omega(\tau))] \right)_\varepsilon \tag{0.7.42}$$

Then, limit

$$\psi_\varepsilon(t,x) = \lim_{R\to\infty} \psi_\varepsilon^R(t,x) \tag{0.7.43}$$

exists uniformly in $t$ and $x$. Function $(\psi_\varepsilon(t,x))_\varepsilon$ solves the Colombeau-Schrödinger equation (0.7.39). The solution is an analytic function of $g$ and $\lambda$ in the same region where $\mathbf{E}[\Theta_t]$ is analytic. Let $(Z_\varepsilon(t,x))_\varepsilon$ be Colombeau generalized function

$$\left( Z_\varepsilon(t,x) \right)_\varepsilon = \frac{1}{2}\varpi_\varepsilon(t)x^2 + \vartheta_\varepsilon(t)x + g(Z_\varepsilon(x)), \tag{0.7.44}$$

where $\mathbf{cl}[(Z_\varepsilon(x))_\varepsilon] \in G_{\mathscr{F}_d^\blacktriangle}(\mathbb{R}^d)$ and $\forall \varepsilon \in (0,1]$, $\varpi_\varepsilon(t), \vartheta_\varepsilon(t)$ are continuous functions of time.

**Theorem 0.7.6.** Assume that $\varpi_\varepsilon(t) \leq \varpi_{0,\varepsilon}$ for each $t$. If $t < \pi/\varpi_{0,\varepsilon}$ then theorem 0.7.5 holds true for potential given by Eq.(0.7.44).

Let us consider a decomposition of $\mathbb{R}^{2n}$ into $\mathbb{R}^n \bigoplus \mathbb{R}^n$ where

$z = (x_1, x_2), x_1, x_2 \in \mathbb{R}^n.$ Let
$\mathbf{G}_z$ be a subset $\mathbf{G}_z \subset \mathbb{R}^{2n}$

$$\mathbf{G}_z = \{w \in \mathbb{R}^{2n} | w = (x_1 + y, x_2 + y), y \in \mathbb{R}^n\} \qquad (0.7.45)$$

and let $Z_\varepsilon(x_1, x_2)$ be a complex-valued potential defined on $\mathbb{R}^{2n}$ such that $\forall \varepsilon \in (0, 1]$
on $\mathbf{G}_z$

$$|\mathrm{Im} Z_\varepsilon(x_1 + y, x_2 + y)| \leq A_\varepsilon(x_1, x_2)(1 + |y|) \qquad (0.7.46)$$

with a certain continuous non-negative function $A_\varepsilon$. We define Colombeau
semigroup $\mathbf{S}_{t,\varepsilon}$ on the set of functions defined on $\mathbb{R}^{2n}$ such that on $\mathbf{G}_z$

$$|\Phi_\varepsilon(x_1 + y, x_2 + y)| \leq M_\varepsilon(x_1, x_2) \exp(|y| C_\varepsilon(x_1, x_2)) \qquad (0.7.47)$$

with certain non-negative continuous functions $C_\varepsilon(x_1, x_2)$ and $M_\varepsilon(x_1, x_2)$. Let
$(\mathbf{S}_{t,\varepsilon})_\varepsilon$, is defined by the Feynman-Kac formula:

$$((\mathbf{S}_{t,\varepsilon}\Phi_\varepsilon)(x_1,x_2))_\varepsilon =$$

$$\left(\mathbf{E}_\varepsilon\left[\exp\left(-\frac{i}{\hbar}\int_0^t d\tau Z_\varepsilon\left(x_1+\frac{\sigma\omega(\tau)}{\sqrt{2}},x_2+\frac{\sigma\omega(\tau)}{\sqrt{2}}\right)\times\right.\right.\right.$$

$$\left.\left.\left.\Phi_\varepsilon\left(x_1+\frac{\sigma\omega(t)}{\sqrt{2}},x_2+\frac{\sigma\omega(t)}{\sqrt{2}}\right)\right)\right]\right)_\varepsilon =$$

$$\left(\int_{\mathbb{R}^n}\int_{\mathbb{R}^n}\mathbf{E}_\varepsilon\left[\delta\left(y_1-x_1-\frac{\sigma\omega(\tau)}{\sqrt{2}}\right)\delta\left(y_2-x_2-\frac{\sigma\omega(\tau)}{\sqrt{2}}\right)\times\right.\right.$$

$$\left.\left.\exp\left(-\frac{i}{\hbar}\int_0^t d\tau\widetilde{Z}\left(x_1+\frac{\sigma\omega(\tau)}{\sqrt{2}},x_2+\frac{\sigma\omega(\tau)}{\sqrt{2}}\right)\right)\right]\Phi_\varepsilon(y_1,y_2)dy_1 dy_2\right)_\varepsilon =$$

$$\left(\int_{\mathbb{R}^n}\int_{\mathbb{R}^n}S_{t,\varepsilon}(x_1,x_2;y_1,y_2)\Phi_\varepsilon(y_1,y_2)dy_1 dy_2\right)_\varepsilon,$$

(0.7.48)

where $\sigma=\sqrt{\hbar/m}$. Using Eq.(0.7.48) one can express the kernel $(S_{t,\varepsilon}(x_1,x_2;y_1,y_2))_\varepsilon$ of $(\mathbf{S}_{t,\varepsilon})_\varepsilon$, by the Brownian bridge $\boldsymbol{\alpha}(s)$. The Brownian bridge is the Gaussian process defined on the interval $[0,1]$ with the covariance

$$\mathbf{E}[\alpha_k(s)\alpha_r(s')]=\delta_{kr}(1-s')\ \text{iff}\ s\leq s' \tag{0.7.49}$$

and the boundary conditions $\alpha(0)=\alpha(1)=0$. Then the expectation value in the second line of Eq.(0.7.48) is equal to the expectation value

$$\mathbf{E}_\varepsilon\left[\delta\left(y_1 - x_1 - \frac{\sigma\omega(\tau)}{\sqrt{2}}\right)\delta\left(y_2 - x_2 - \frac{\sigma\omega(\tau)}{\sqrt{2}}\right)\right] =$$

$$(4\pi t\sigma^2)^{-n/2}\exp\left(-\frac{(x_1 - x_2 - y_1 + y_2)^2}{4t\sigma^2}\right)\delta(x_1 - x_2 - y_1 + y_2)$$

$$(0.7.50)$$

times the expectation value of $\exp\left(-(i/h)\int_0^t d\tau(Z_\varepsilon(\tau))_\varepsilon\right)$, where the Brownian motion is constrained to end at $y$. A solution of this constraint gives a simple formula for the kernel

$$(S_{t,\varepsilon}(x_1, x_2; y_1, y_2))_\varepsilon =$$

$$(4\pi t\sigma^2)^{-n/2}\exp\left(-\frac{(x_1 - x_2 - y_1 + y_2)^2}{4t\sigma^2}\right)\delta(x_1 - x_2 - y_1 + y_2) \times$$

$$\left(\mathbf{E}\left[\exp\left(-\frac{i}{h}\int_0^t d\tau Z_\varepsilon\left(x_1\left(1 - \frac{s}{t}\right) + y_1\left(\frac{s}{t}\right) + \sqrt{\frac{t}{2}}\,\sigma\boldsymbol{\alpha}\left(\frac{s}{t}\right)\right),\right.\right.$$

$$\left.\left. x_2\left(1 - \frac{s}{t}\right) + y_2\left(\frac{s}{t}\right) + \sqrt{\frac{t}{2}}\,\sigma\boldsymbol{\alpha}\left(\frac{s}{t}\right)\right)ds\right]\right)_\varepsilon$$

$$(0.7.51)$$

Now we restrict ourselves to functions of the form

$$(Z_\varepsilon(x_1, x_2))_\varepsilon = \left(\widetilde{Z}_\varepsilon(x_1 + ix_2)\right)_\varepsilon, (\Phi_\varepsilon(x_1, x_2))_\varepsilon = (\varphi_\varepsilon(x_1 + ix_2))_\varepsilon, \qquad (0.7.52)$$

where $\forall\varepsilon \in (0, 1]$, $Z_\varepsilon$ and $\varphi_\varepsilon$ are holomorphic functions of their arguments. Substituting equation (0.7.52) into equation (0.7.48) we obtain

$$((\mathbf{S}_{t,\varepsilon}\varphi_\varepsilon)(x_1 + ix_2))_\varepsilon =$$

$$(4\pi t\sigma^2)^{-n/2}\Bigg(\int_{\mathbb{R}^n} dy\, E_{t,\varepsilon}(y, x_1, x_2)\varphi_\varepsilon\Big(\frac{1}{2}(y + x_1 + x_2) + \frac{i}{2}(y - x_1 + x_2)\Big) \times \tag{0.7.53}$$

$$\exp\Bigg(-\frac{(y - x_1 - x_2)^2}{4t\sigma^2}\Bigg)\Bigg)_\varepsilon$$

where by $E_{t,\varepsilon}$ we have denoted $\mathbf{E}[\dots]$. We note that

$$(4\pi t\sigma^2)^{-n/2}\int_{\mathbb{R}^n} dy\, \chi\Big(\frac{1}{2}(y + x) + \frac{i}{2}(y - x)\Big)\exp\Bigg(-\frac{(x - y)^2}{4t\sigma^2}\Bigg) =$$

$$(4\pi t\sigma^2)^{-n/2}\int_{\mathbb{R}^n} dy\, \chi \exp\Bigg(-\frac{(x - y)^2}{2it\sigma^2}\Bigg) \tag{0.7.54}$$

(where $i^{n/2} = \exp(\frac{1}{2}in\pi)$) for any analytic function $\chi_\varepsilon \in \mathcal{L}_1(\mathbb{R}^n)$ of the form

$\forall \varepsilon \in (0, 1], \chi_\varepsilon(y) = \int dp\, \tilde{\chi}_\varepsilon(p)\exp(ipy)$. Finally, we set $x_2 = 0$ in equation (0.7.53) and apply the identity (0.7.54). We obtain the final formula for the generalized Feynman propagator:

$$(K_{t,\varepsilon}(x,y))_\varepsilon =$$

$$(4\pi t\sigma^2)^{-n/2} \exp\left(-\frac{(x-y)^2}{2it\sigma^2}\right) \times$$

$$\left(\mathbf{E}\left[\exp\left(-\frac{i}{\hbar}\int_0^t d\tau \widetilde{Z}_\varepsilon\Big(x\Big(1-\frac{s}{t}\Big)+y\Big(\frac{s}{t}\Big)+\lambda\sqrt{t}\,\sigma\boldsymbol{\alpha}\Big(\frac{s}{t}\Big)\Big)\right)\right]\right)_\varepsilon =$$

$$K_t^0(x,y)(\mathfrak{R}_\varepsilon(t,x,y))_\varepsilon. \qquad (0.7.55)$$

Here $\lambda = \sqrt{i}$ and we have denoted the free propagator by $K_t^0(x,y)$.
**Theorem 0.7.7**. Define

$$(\mathfrak{R}_{t,\varepsilon}^R(t,x,y))_\varepsilon =$$

$$\mathbf{E}\left[\exp\left\{\frac{g}{\hbar\lambda^2}\int_0^t ds\,\theta_R\Big(\alpha\Big(\frac{s}{t}\Big)\Big)Z_\varepsilon\Big(x\Big(1-\frac{s}{t}\Big)+y\frac{s}{t}+\lambda\sigma\Big)\right\}\right], \qquad (0.7.56)$$

where $(Z_\varepsilon)_\varepsilon \in G_{\mathscr{F}_n^{\blacktriangle}}(\mathbb{R}^d)$. The limit $R \to \infty$ exists $\forall \varepsilon \in (0,1]$ uniformly in $t, x, y$.
$\forall \varepsilon \in (0,1]$, $\mathfrak{R}_\varepsilon(t,x,y)$, is an analytic function of $g$ and $\lambda$ if $\operatorname{Re}\lambda^2 \geq 0$. Formula

$$(K_{t,\varepsilon}(x,y))_\varepsilon = K_t^0(x,y)(\lim_{R\to\infty}\mathfrak{R}_{t,\varepsilon}^R(t,x,y))_\varepsilon \qquad (0.7.57)$$

gives the Feynman propagator (0.7.55) for any $(Z_\varepsilon)_\varepsilon \in G_{\mathscr{F}_{\blacktriangle}}(\mathbb{R}^d)$.

# I.0.8. Colombeau-Feynman Path Integral using Fujiwara approach.

We shall treat now the quantum dynamical system using Colombeau-Feynman Path Integlal. Let $\mathscr{L}(\dot{x}, x, t)$ be the Lagrangian of the form:

$$\mathscr{L}(\dot{x}, x, t) = \frac{1}{2}|\dot{x}|^2 - V(t, x). \tag{0.8.1}$$

Here potential $V(t, x)$ is assumed to satisfy the following assumptions:

**(i)** $V(t, x)$ is a real valued function of $(t, x) \in \mathbb{R}_+ \times \mathbb{R}^n$. For any fixed $t \in \mathbb{R}_+$, $V(t, x)$ is a function of $x$ of class $C^\infty$. $V(t, x)$ is a measurable function of $(t, x) \in \mathbb{R}_+ \times \mathbb{R}^n$.

**(ii)** For any multi-index $\alpha$, with its length $|\alpha| \geq 2$, the non-negative measurable function $\mu_\alpha(t)$ defined by

$$\mu_\alpha(t) = \sup_{x \in \mathbb{R}^n} \left| \left( \frac{\partial}{\partial x} \right)^\alpha V(t, x) \right| + \sup_{|x| \leq 1} |V(t, x)| \tag{0.8.2}$$

is essentially bounded on every bounded interval of $\mathbb{R}$.

**(iii)** We fix a sufficiently large $K$, and we let $T = \infty$ if

$$\mu^* = \sum_{2 \leq \alpha \leq K} ess. \sup_{t \in \mathbb{R}} \mu_\alpha(t) < \infty \tag{0.8.3}$$

Otherwise, we let $T$ be an arbitrarily fixed finite positive number. We shall discuss everything in the time interval $(-T, T)$.

Let $S(t, s, x, y)$ be the classical action along the classical orbit starting from $y$ at time $s$ and reaching $x$ at time $t$. It well known that there exists a positive constant $\delta(T)$ such that $S(t, s, x, y)$ is uniquely defined for any $x$ and $y$ in $\mathbb{R}^n$ if $|t - s| < \delta(T)$ [40]. Using this, we define the following oscillatory integral transformations used in [40]. For any $N = 0, 1, 2, \ldots$

$$E^{(N)}\left( \tilde{\lambda}, t, s \right) f(x) =$$

$$\left( \frac{-\tilde{\lambda}}{2\pi(t-s)} \right)^{n/2} \int_{\mathbb{R}^n} d^n y f(y) a^{(N)}\left( \tilde{\lambda}, t, s, x, y \right) \exp[\lambda S(t, s, x, y)]. \tag{0.8.4}$$

Here $\tilde{\lambda} = ih^{-1}, h$ being a small positive parameter. In Eq.(0.8.4), the amplitude function $a^{(N)}\big(\tilde{\lambda}, t, s, x, y\big)$ is defined by

$$a^{(N)}\big(\tilde{\lambda}, t, s, x, y\big) = 1,$$

$$a^{(N)}\big(\tilde{\lambda}, t, s, x, y\big) = \sum_{j=1}^{N} \tilde{\lambda}^{1-j} a_j(t, s, x, y),$$

(0.8.5)

where $a_j(t, s, x, y)$ will be defined by transport equation [40].Let $[s, t]$ be an arbitrary time interval in $(-T, T)$ and let $\Delta = \{s = t_0 < t_1 < \ldots t_L = t\}$ be an arbitrary subdivision of the interval $[s, t]$. We put $\delta(\Delta) = \max_{1 \leqslant j \leqslant L}(|t_j - t_{j-1}|)$. We introduce the iterated integral operators $E^{(N)}(\Delta|\tilde{\lambda}, t, s)$ and $E^{(N)}(\Delta|\tilde{\lambda}, s, t)$ associated with the subdivision $\Delta$ by the formulae

$$E^{(N)}(\Delta|\tilde{\lambda}, t, s) = E^{(N)}(\tilde{\lambda}, t, t_{L-1}) E^{(N)}(\tilde{\lambda}, t_{L-1}, t_{L-2}) \cdots E^{(N)}(\tilde{\lambda}, t_{L-1}, t),$$

$$E^{(N)}(\Delta|\tilde{\lambda}, s, t) = E^{(N)}(\tilde{\lambda}, s, t_1) E^{(N)}(\tilde{\lambda}, t_1, t_2) \cdots E^{(N)}(\tilde{\lambda}, t_L, t)$$

(0.8.6)

for $N = 0, 1, 2, \ldots$ . Let $\mathbf{l}^{(N)}(\Delta|\tilde{\lambda}, t, s, x, y)$ and $\mathbf{l}^{(N)}(\Delta|\tilde{\lambda}, s, t, x, y)$ be the kernel functions of $E^{(N)}(\Delta|\tilde{\lambda}, t, s)$ and $E^{(N)}(\Delta|\tilde{\lambda}, s, t)$, respectively. Then, we have

$$E^{(N)}(\Delta|\tilde{\lambda}, t, s) f(x) = \int_{\mathbb{R}^n} d^n y f(y) \mathbf{l}^{(N)}(\Delta|\tilde{\lambda}, t, s, x, y),$$

$$\mathbf{l}^{(N)}(\Delta|\tilde{\lambda}, t, s, x, y) =$$

$$\prod_{j=1}^{L} \left( \frac{-\tilde{\lambda}}{2\pi(t_j - t_{j-1})} \right)^{n/2} \int_{\mathbb{R}^n} \cdots \int_{\mathbb{R}^n} \prod_{j=1}^{L} a^{(N)}\big(\tilde{\lambda}, t_j, t_{j-1}, x^j, y^{j-1}\big) \times$$

(0.8.7)

$$\exp\left\{ \tilde{\lambda} \left[ \sum_{j=1}^{L} S(t_j, t_{j-1}, x^j, y^{j-1}) \right] \right\} \prod_{j=1}^{L-1} dx^j,$$

where $x^0 = y$ and $x^L = x$. Similar formulae are valid for $E^{(N)}(\Delta|\tilde{\lambda}, s, t)$.Feynman formally "stated" in [13] that limit

$$\lim_{\delta(\Delta)\to 0} \mathbf{I}^{(0)}(\Delta|\tilde{\lambda},t,s,x,y) \tag{0.8.8}$$

exists and equals the kernel function of the fundamental solution (Green's function) for the Schrödinger equation (0.1.1) or

$$\frac{\partial}{\tilde{\lambda}\partial t}u(t,x) + \frac{1}{2}\sum_{j=1}^{n}\left(\frac{\partial}{\tilde{\lambda}\partial x_j}\right)^2 u(t,x) + V(t,x)u(t,x) = 0, \tag{0.8.9}$$

$$u(0,x) = u_0(x).$$

Note that

$$\mathbf{I}^{(N)}(\Delta|\tilde{\lambda},t,s,x,y) = \left(\frac{-\tilde{\lambda}}{2\pi(t-s)}\right)^{n/2}a^{(N)}\big(\tilde{\lambda},t,s,x,y\big)\exp\big[\tilde{\lambda}S(t,s,x,y)\big] \tag{0.8.10}$$

and

$$\mathbf{I}^{(N)}(\Delta|\tilde{\lambda},s,t,x,y) = \left(\frac{-\tilde{\lambda}}{2\pi(s-t)}\right)^{n/2}a^{(N)}\big(\tilde{\lambda},s,t,x,y\big)\exp\big[\tilde{\lambda}S(s,t,x,y)\big] \tag{0.8.11}$$

if $|t-s| < \delta(T)$. These amplitude functions $a^{(N)}\big(\tilde{\lambda},t,s,x,y\big)$ and $a^{(N)}\big(\tilde{\lambda},s,t,x,y\big)$ belong to the function space $D(\mathbb{R}_x^n \times \mathbb{R}_y^n)$ of Schwartz [47] with parameters $\Delta,\tilde{\lambda},t,s$. As usual we shall employ the norm $\|f\|_m$ defined by

$$\|f\|_m = \sum_{|\alpha+\beta|\leq m}\sup_{(x,y)\in\mathbb{R}_x^n\times\mathbb{R}_y^n}\left|\left(\frac{\partial}{\partial x}\right)^\alpha\left(\frac{\partial}{\partial y}\right)^\beta f(x,y)\right| \tag{0.8.12}$$

for $m = 0,1,2\ldots$. We shall denote by $a^{(N)}\big(\tilde{\lambda},t,s\big)$ the function $a^{(N)}\big(\tilde{\lambda},t,s,x,y\big)$

considered as an element of the function space $D(\mathbb{R}^n_x \times \mathbb{R}^n_y)$ with parameters $\Delta, \tilde{\lambda}, t, s$.

Fujiwara [40] prove Feynman assertion (0.8.8) under the assumptions (**i**)-(**iii**) concerning the potential $V(t,x)$.

**Theorem**.**0**.**8**.**1**. (**i**) (Fujiwara [40]) Assume that (i) and (ii) hold. Assume further that $|t-s| < \delta(T)$ Then, for any $N = 0, 1, 2 \ldots$ both the limits

$$k\big(\tilde{\lambda}, t, s\big) = \lim_{\delta(\Delta) \to 0} a^{(N)}\big(\Delta|\tilde{\lambda}, t, s\big) \qquad (0.8.13)$$

and

$$k\big(\tilde{\lambda}, s, t\big) = \lim_{\delta(\Delta) \to 0} a^{(N)}\big(\Delta|\tilde{\lambda}, s, t\big) \qquad (0.8.14)$$

exist in the function space $D(\mathbb{R}^n_x \times \mathbb{R}^n_y)$. Moreover, for any integers $m > 0$ and $N > 0$, there exists a positive constant $\gamma = \gamma(m, N, T)$ such that

$$\|k\big(\tilde{\lambda}, t, s\big) - a^{(N)}\big(\Delta|\tilde{\lambda}, t, s\big)\|_m \leq \gamma |\lambda|^{-N} (\delta(\Delta))^{N+1} |t-s| \qquad (0.8.15)$$

and

$$\|k\big(\tilde{\lambda}, s, t\big) - a^{(N)}\big(\Delta|\tilde{\lambda}, s, t\big)\|_m \leq \gamma |\lambda|^{-N} (\delta(\Delta))^{N+1} |t-s| \qquad (0.8.16)$$

if $|t-s| < \delta(T)$. The constant $\gamma(m, N, T)$ is independent of any particular choice of $t, s, \Delta$ and $\tilde{\lambda}$ if $|\tilde{\lambda}|$ is bounded away from $0$. The function $k\big(\tilde{\lambda}, t, s, x, y\big)$ of is independent of $N$. For any $t \in [s - \delta(T), s + \delta(T)]$, we define

$$U\big(\tilde{\lambda}, t, s\big)\varphi(x) = \left(\frac{-\tilde{\lambda}}{2\pi(t-s)}\right)^{n/2} \int_{\mathbb{R}^n} dy \varphi(y) k\big(\tilde{\lambda}, t, s, x, y\big) \exp\big[\tilde{\lambda} S(t, s, x, y)\big] \qquad (0.8.17)$$

Then, $u(t,x) = U(\tilde{\lambda}, t, s)\varphi(x)$ satisfies the Schrödinger equation (0.8.9) and the initial condition $u(s,x) = \varphi(x)$.

**(ii)** Let $f(x) \in \mathcal{L}_2(\mathbb{R}^n)$ then

$$U(\tilde{\lambda}, t, s)f(x) = \int_{\mathbb{R}^n} \left\{ \lim_{\delta(\Delta) \to 0} \mathbf{I}^{(N)}(\Delta|\tilde{\lambda}, t, s, x, y) \right\} f(y)dy =$$

$$(0.8.18)$$

$$\lim_{\delta(\Delta) \to 0} E^{(N)}(\Delta|\tilde{\lambda}, t, s)f(x)$$

in the sense of convergence in $\mathcal{L}_2(\mathbb{R}^n)$, i.e. one can interchange the order of operations $\lim_{\delta(\Delta) \to 0}$ and $\int dy$ if $f(x) \in \mathcal{L}_2(\mathbb{R}^n)$.

**Theorem.0.8.2.**(Fujiwara [40]) If $|t - s| < \delta(T)$, then for any integers $N \geq 0$ and

$m \geq 0$, we have

$$\left\| k(\tilde{\lambda}, t, s) - a^{(N)}(\tilde{\lambda}, t, s) \right\|_m \leq \gamma |\lambda|^{-N}|t - s|^{N+1},$$

$$(0.8.19)$$

$$\left\| k(\tilde{\lambda}, s, t) - a^{(N)}(\tilde{\lambda}, s, t) \right\|_m \leq \gamma |\lambda|^{-N}|t - s|^{N+1}$$

with some constant $\gamma = \gamma(m, N, T)$. Here $a^{(N)}(\tilde{\lambda}, t, s)$ is the amplitude function of $E^{(N)}(\Delta|\tilde{\lambda}, t, s)$ of eq.(0.8.4).

In our subsequent consideration, the potential $V(t,x)$ does not satisfy Fujiwara condition **(i)-(iii)**. We remove it by cut-off or another type of the regularization. As a framework we use Colombeau's algebra of generalized functions, i.e. we change potential $V(t,x)$ by Colombeau generalized function $(V_\varepsilon(t,x))_{\varepsilon \in (0,1]}$ such that $V_0(t,x) = V(t,x)$ and

$(V_\varepsilon \mathbf{i})$ $V_\varepsilon(t,x)$ is a real valued function of $(t,x) \in \mathbb{R}_+ \times \mathbb{R}^n$. For any fixed $t \in \mathbb{R}_+$,

**cl**$\left[(V_\varepsilon(t,x))_\varepsilon\right] \in G(\mathbb{R}^n)$.

$\forall \varepsilon \in (0,1]\left[ V_\varepsilon(t,x) \text{ is a measurable function of } (t,x) \in \mathbb{R}_+ \times \mathbb{R}^n \right]$.

$(V_\varepsilon \mathbf{ii})$ For any multi-index $\alpha$, with its length $|\alpha| \geq 2$, the non-negative measurable function $\mu_{\alpha,\varepsilon}(t)$ defined by

$$\mu_{\alpha,\varepsilon}(t) = \sup_{x \in \mathbb{R}^n} \left| \left( \frac{\partial}{\partial x} \right)^{\alpha} V_{\varepsilon}(t,x) \right| + \sup_{|x| \leq 1} |V_{\varepsilon}(t,x)|, \tag{0.8.20}$$

is essentially bounded on every bounded interval of $\mathbb{R}$ and $\mu_{\alpha,\varepsilon}(t) = O(\varepsilon^{-M(\alpha)})$.

$(V_{\varepsilon}\textbf{iii})$ We fix a sufficiently large $K$, and we let $T = \infty$ if $\forall \varepsilon \in (0,1]$ :

$$\mu_{\varepsilon}^{*} = \sum_{2 \leq \alpha \leq K} ess. \sup_{t \in \mathbb{R}} \mu_{\alpha,\varepsilon}(t) < \infty \tag{0.8.21}$$

We shall discuss Feynman path integral representation of the Colombeau's solutions
for the Colombeau-Schrödinger type equation

$$\frac{\partial}{\tilde{\lambda}\partial t}(u_{\varepsilon}(t,x))_{\varepsilon} + \frac{1}{2}\sum_{j=1}^{n}\left(\frac{\partial}{\tilde{\lambda}\partial x_j}\right)^2 (u_{\varepsilon}(t,x))_{\varepsilon} + ((V_{\varepsilon}(t,x))_{\varepsilon})(u_{\varepsilon}(t,x))_{\varepsilon} = 0,$$
$$\tag{0.8.22}$$

$$(u_{\varepsilon}(0,x))_{\varepsilon} = (u_{0,\varepsilon}(x))_{\varepsilon}.$$

Initial data are strongly singular and regularized with delta sequence. It well known that there exists Colombeau's solutions of the Eq.(0.8.22),which are compatible with classical solutions in a limiting case when $\varepsilon \to 0$ [12],[42-43].
Let $S_{\varepsilon}(t,s,x,y)$ be the classical action along the classical orbit starting from $y$ at time $s$ and reaching $x$ at time $t$. It well known that there exists a positive constant $\delta(T)$ such that $S(t,s,x,y)$ is uniquely defined for any $x$ and $y$ in $\mathbb{R}^n$ if $|t-s| < \delta(T,\varepsilon)$ [40]. Using this, we define the following Colombeau oscillatory integral transformations. For any $N = 0, 1, 2,$

$$\left(E_\varepsilon^{(N)}\big(\tilde{\lambda},t,s\big)f_\varepsilon(x)\right)_\varepsilon =$$

$$\left(\left(\frac{-\tilde{\lambda}}{2\pi(t-s)}\right)^{n/2}\left(\int_{\mathbb{R}^n}d^n y f_\varepsilon(y)a_\varepsilon^{(N)}\big(\tilde{\lambda},t,s,x,y\big)\exp[\lambda S_\varepsilon(t,s,x,y)]\right)\right)_\varepsilon = \qquad (0.8.23)$$

$$\left(\frac{-\tilde{\lambda}}{2\pi(t-s)}\right)^{n/2}\int_{\mathbb{R}^n}d^n y((f_\varepsilon(y))_\varepsilon\Big(\big(a_\varepsilon^{(N)}\big(\tilde{\lambda},t,s,x,y\big)\big)_\varepsilon\Big)\exp[\lambda((S_\varepsilon(t,s,x,y))_\varepsilon)]$$

Here $\tilde{\lambda}=ih^{-1}$, $h$ being a small positive parameter. In Eq.(0.8.23), the generalized amplitude function is defined by

$$\left(a_\varepsilon^{(N)}\big(\tilde{\lambda},t,s,x,y\big)\right)_\varepsilon = 1,$$

$$\left(a_\varepsilon^{(N)}\big(\tilde{\lambda},t,s,x,y\big)\right)_\varepsilon = \sum_{j=1}^{N}\tilde{\lambda}^{1-j}(a_{j,\varepsilon}(t,s,x,y))_\varepsilon, \qquad (0.8.24)$$

where $(a_{j,\varepsilon}(t,s,x,y))_\varepsilon$ will be defined by transport equation.Let $[s,t]$ be an arbitrary time interval in $(-T,T)$ and let $\Delta(\varepsilon)=\{s=t_0<t_1<\ldots t_{L_\varepsilon}=t\}$ be an arbitrary subdivision of the interval $[s,t]$. We put $\delta(\Delta(\varepsilon))=\max_{1\le j\le L_\varepsilon}(|t_j-t_{j-1}|)$. We introduce the iterated generalized integral operators $\left[\left(E_\varepsilon^{(N)}(\Delta(\varepsilon)|\tilde{\lambda},t,s)\right)_\varepsilon\right]$ and $\left[\left(E_\varepsilon^{(N)}(\Delta(\varepsilon)|\tilde{\lambda},s,t)\right)_\varepsilon\right]$ associated with the subdivision $\Delta(\varepsilon)$ by the formulae:

$$\left[\left(E_\varepsilon^{(N)}(\Delta(\varepsilon)|\tilde{\lambda},t,s)\right)_\varepsilon\right] =$$

$$\left[\left(E_\varepsilon^{(N)}(\tilde{\lambda},t,t_{L_\varepsilon-1})\right)_\varepsilon\right]\left[\left(E_\varepsilon^{(N)}(\tilde{\lambda},t_{L_\varepsilon-1},t_{L_\varepsilon-2})\right)_\varepsilon\right]\cdots\left[\left(E_\varepsilon^{(N)}(\tilde{\lambda},t_{L_\varepsilon-1},t)\right)_\varepsilon\right],$$

$$\left[\left(E_\varepsilon^{(N)}(\Delta(\varepsilon)|\tilde{\lambda},s,t)\right)_\varepsilon\right] =$$

$$\left[\left(E_\varepsilon^{(N)}(\tilde{\lambda},s,t_1)\right)_\varepsilon\right]\left[\left(E_\varepsilon^{(N)}(\tilde{\lambda},t_1,t_2)\right)_\varepsilon\right]\cdots\left[\left(E_\varepsilon^{(N)}(\tilde{\lambda},t_{L_\varepsilon},t)\right)_\varepsilon\right]$$

(0.8.25)

for $N = 0,1,2,\ldots$ . Let $\left(\mathbf{I}_\varepsilon^{(N)}(\Delta(\varepsilon)|\tilde{\lambda},t,s,x,y)\right)_\varepsilon$ and $\left(\mathbf{I}_\varepsilon^{(N)}(\Delta(\varepsilon)|\tilde{\lambda},t,s,x,y)\right)_\varepsilon$ be the kernel functions of $\left(E_\varepsilon^{(N)}(A|\tilde{\lambda},t,s)\right)_\varepsilon$ and $\left(E_\varepsilon^{\ell(N)}(A|\tilde{\lambda},s,t)\right)_\varepsilon$, respectively. Then, we have

$$\left[\left(E_\varepsilon^{(N)}(\Delta(\varepsilon)|\tilde{\lambda},t,s)\right)_\varepsilon\right][(f_\varepsilon(x))_\varepsilon] = \int_{\mathbb{R}^n} d^{\,n}y\,[(f_\varepsilon(y))_\varepsilon]\left[\left(\mathbf{I}_\varepsilon^{(N)}(\Delta(\varepsilon)|\tilde{\lambda},t,s,x,y)\right)_\varepsilon\right],$$

$$\left(\mathbf{I}_\varepsilon^{(N)}(\Delta(\varepsilon)|\tilde{\lambda},t,s,x,y)\right)_\varepsilon =$$

$$\prod_{j=1}^{L_\varepsilon}\left(\frac{-\tilde{\lambda}}{2\pi(t_j-t_{j-1})}\right)^{n/2}\int_{\mathbb{R}^n}\cdots\int_{\mathbb{R}^n}\prod_{j=1}^{L_\varepsilon}\left[\left(a_\varepsilon^{(N)}(\tilde{\lambda},t_j,t_{j-1},x^j,y^{j-1})\right)_\varepsilon\right]\times$$

(0.8.26)

$$\exp\left\{\tilde{\lambda}\left[\sum_{j=1}^{L_\varepsilon}(S_\varepsilon(t_j,t_{j-1},x^j,y^{j-1}))_\varepsilon\right]\right\}\prod_{j=1}^{L_\varepsilon-1}dx^j,$$

where $x^0 = y$ and $x^{L_\varepsilon} = x$. Similar formulae are valid for $\left(E_\varepsilon^{(N)}(\Delta|\tilde{\lambda},s,t)\right)_\varepsilon$. We stated that $\forall\varepsilon \in (0,1]$

$$\lim_{\delta(\Delta(\varepsilon))\to 0}\mathbf{I}_\varepsilon^{(0)}(\Delta(\varepsilon)|\tilde{\lambda},t,s,x,y)$$

(0.8.27)

exists and equals the generalized kernel function of the fundamental solution (Green's function) for the Colombeau-Schrödinger equation (0.8.22).

We note that

$$\left(\mathbf{l}_{\varepsilon}^{(N)}(\Delta|\tilde{\lambda},t,s,x,y)\right)_{\varepsilon} = \left(\frac{-\tilde{\lambda}}{2\pi(t-s)}\right)^{n/2}\left(a_{\varepsilon}^{(N)}\left(\tilde{\lambda},t,s,x,y\right)\exp\left[\tilde{\lambda}S_{\varepsilon}(t,s,x,y)\right]\right)_{\varepsilon} \qquad (0.8.28)$$

and

$$\left(\mathbf{l}_{\varepsilon}^{(N)}(\Delta|\tilde{\lambda},s,t,x,y)\right)_{\varepsilon} = \left(\frac{-\tilde{\lambda}}{2\pi(s-t)}\right)^{n/2}\left(a_{\varepsilon}^{(N)}\left(\tilde{\lambda},s,t,x,y\right)\exp\left[\tilde{\lambda}S_{\varepsilon}(s,t,x,y)\right]\right)_{\varepsilon} \qquad (0.8.29)$$

if $|t-s| < \delta(T)$. These Colombeau generalized amplitude functions $a^{(N)}\left(\tilde{\lambda},t,s,x,y\right)$

and $\forall \varepsilon \in (0,1], a_{\varepsilon}^{(N)}\left(\tilde{\lambda},s,t,x,y\right)$ belong to the function space $D(\mathbb{R}_x^n \times \mathbb{R}_y^n)$ of Schwartz [47] with parameters $\varepsilon, \Delta(\varepsilon), \tilde{\lambda}, t, s$. As usual we shall employ the norm $\|f_{\varepsilon}\|_m$ defined by

$$\|f_{\varepsilon}\|_m = \sum_{|\alpha+\beta|\leq m} \sup_{(x,y)\in\mathbb{R}_x^n\times\mathbb{R}_y^n} \left|\left(\frac{\partial}{\partial x}\right)^{\alpha}\left(\frac{\partial}{\partial y}\right)^{\beta}f_{\varepsilon}(x,y)\right| \qquad (0.8.30)$$

for $m = 0,1,2\ldots$ . We shall denote by $a_{\varepsilon}^{(N)}\left(\tilde{\lambda},t,s\right)$ the function $a_{\varepsilon}^{(N)}\left(\tilde{\lambda},t,s,x,y\right)$ considered as an element of the function space $D(\mathbb{R}_x^n \times \mathbb{R}_y^n)$ with parameters $\varepsilon, \Delta(\varepsilon), \tilde{\lambda}, t, s$.

**Theorem**.**0.8.3**. (i) Assume that generalized Fujiwara conditions $(V_{\varepsilon}\mathbf{i})$ and $(V_{\varepsilon}\mathbf{ii})$ hold. Assume further that $|t-s| < \delta(T)$ Then, for any $N = 0,1,2\ldots$ both the limits

$$\left(k_{\varepsilon}\left(\tilde{\lambda},t,s\right)\right)_{\varepsilon} = \left(\lim_{\delta(\Delta(\varepsilon))\to 0} a_{\varepsilon}^{(N)}\left(\Delta(\varepsilon)|\tilde{\lambda},t,s\right)\right)_{\varepsilon} \qquad (0.8.31)$$

and

$$\left( k_\varepsilon\big(\tilde{\lambda}, s, t\big) \right)_\varepsilon = \left( \lim_{\delta(\Delta(\varepsilon)) \to 0} a^{(N)}\big(\Delta(\varepsilon)|\tilde{\lambda}, s, t\big) \right)_\varepsilon \qquad (0.8.32)$$

exist in the function space $D(\mathbb{R}^n_x \times \mathbb{R}^n_y)$. Moreover, for any integers $m > 0$ and $N > 0$, there exists a positive constant $\gamma_\varepsilon = \gamma(m, N, T, \varepsilon)$ such that

$$\left( \left\| k_\varepsilon\big(\tilde{\lambda}, t, s\big) - a^{(N)}_\varepsilon\big(\Delta(\varepsilon)|\tilde{\lambda}, t, s\big) \right\|_m \right)_\varepsilon \le \left( \big(\gamma_\varepsilon(\delta(\Delta(\varepsilon)))^{N+1}\big)_\varepsilon \right) |\lambda|^{-N} |t - s| \qquad (0.8.33)$$

and

$$\left( \left\| k_\varepsilon\big(\tilde{\lambda}, s, t\big) - a^{(N)}_\varepsilon\big(\Delta(\varepsilon)|\tilde{\lambda}, s, t\big) \right\|_m \right)_\varepsilon \le \left( \big(\gamma_\varepsilon(\delta(\Delta(\varepsilon)))^{N+1}\big)_\varepsilon \right) |\lambda|^{-N} |t - s| \qquad (0.8.34)$$

if $|t - s| < \delta(T)$. The constant $\gamma(m, N, T, \varepsilon)$ is independent of any particular choice of $\varepsilon, t, s, \Delta(\varepsilon)$ and $\tilde{\lambda}$ if $|\tilde{\lambda}|$ is bounded away from $0$. The function $\left( k_\varepsilon\big(\tilde{\lambda}, t, s, x, y\big) \right)_\varepsilon$ of is independent of $N$. For any $t \in [s - \delta(T), s + \delta(T)]$, we define

$$\left( U_\varepsilon\big(\tilde{\lambda}, t, s\big) \varphi_\varepsilon(x) \right)_\varepsilon =$$

$$\left( \frac{-\tilde{\lambda}}{2\pi(t - s)} \right)^{n/2} \left( \int_{\mathbb{R}^n} dy \varphi_\varepsilon(y) k_\varepsilon\big(\tilde{\lambda}, t, s, x, y\big) \exp\big[ \tilde{\lambda} S_\varepsilon(t, s, x, y) \big] \right)_\varepsilon \qquad (0.8.35)$$

Then generalized function $(u_\varepsilon(t, x))_\varepsilon = \left( U_\varepsilon\big(\tilde{\lambda}, t, s\big) \varphi_\varepsilon(x) \right)_\varepsilon$ satisfies the Colombeau-Schrödinger equation (0.8.22) and the initial condition $(u_\varepsilon(s, x))_\varepsilon = (\varphi_\varepsilon(x))_\varepsilon$.

**(ii)** Let $\forall \varepsilon \in (0, 1], f_\varepsilon(x) \in \mathcal{L}_2(\mathbb{R}^n)$ then

$$\left( U_\varepsilon\left(\tilde{\lambda},t,s\right)f_\varepsilon(x)\right)_\varepsilon = \left( \int_{\mathbb{R}^n}\left\{\lim_{\delta(\Delta)\to 0}\mathbf{l}_\varepsilon^{(N)}(\Delta(\varepsilon)|\tilde{\lambda},t,s,x,y)\right\}f_\varepsilon(y)dy\right)_\varepsilon =$$

$$\text{(0.8.36)}$$

$$\left( \lim_{\delta(\Delta)\to 0} E_\varepsilon^{(N)}(\Delta(\varepsilon)|\tilde{\lambda},t,s)f_\varepsilon(x)\right)_\varepsilon$$

in the sense of convergence in $\mathscr{L}_2(\mathbb{R}^n)$, i.e. one can interchange the order of operations $\lim\limits_{\delta(\Delta)\to 0}$ and $\int dy$ if $f_\varepsilon(x)\in \mathscr{L}_2(\mathbb{R}^n)$.

**Theorem**.**0.8.4**.If $|t-s|<\delta(T)$,then for any integers $N\geq 0$ and $m\geq 0$,we have

$$\left(\left\| k_\varepsilon\left(\tilde{\lambda},t,s\right)-a^{(N)}\left(\tilde{\lambda},t,s\right)\right\|_m\right)_\varepsilon \leq ((\gamma_\varepsilon)_\varepsilon)|\lambda|^{-N}|t-s|^{N+1},$$

$$\text{(0.8.37)}$$

$$\left(\left\| k_\varepsilon\left(\tilde{\lambda},s,t\right)-a_\varepsilon^{(N)}\left(\tilde{\lambda},s,t\right)\right\|_m\right)_\varepsilon \leq ((\gamma_\varepsilon)_\varepsilon)|\lambda|^{-N}|t-s|^{N+1}$$

with some constant $\gamma_\varepsilon = \gamma(m,N,T,\varepsilon)$. Here $\left( a_\varepsilon^{(N)}\left(\tilde{\lambda},t,s\right)\right)_\varepsilon$ is the generalized amplitude function of $\left( E_\varepsilon^{(N)}(\Delta(\varepsilon)|\tilde{\lambda},t,s)\right)_\varepsilon$ of Eq.(0.8.23).

# I.1.Colombeau generalized functions.

We use [1],[2],[3],[4],[5] as standard references for the foundations and various applications of standard Colombeau theory. We briefly recall the basic Colombeau construction. Throughout the paper $\Omega$ will denote an open subset of $\mathbb{R}^n$. Stanfard Colombeau generalized functions on $\Omega$ are defined as equivalence classes $u = [(u_\varepsilon)_{\varepsilon\in(0,1]}]$ of nets of smooth functions $u_\varepsilon\in C^\infty(\Omega)$ (regularizations) subjected to asymptotic norm conditions with respect to $\varepsilon\in(0,1]$ for their derivatives on compact sets.

The basic idea of *classical Colombeau's theory of nonlinear generalized functions* [1],[2] is regularization by sequences (nets) of smooth functions and the use of asymptotic estimates in terms of a regularization parameter $\varepsilon$. Let $(u_\varepsilon)_{\varepsilon\in(0,1]}$ with $u_\varepsilon \in C^\infty(M)$ for all $\varepsilon \in \mathbb{R}_+$, where $M$ a separable, smooth orientable Hausdorff manifold of dimension $n$.

**Definition 1.1.1.** The classical Colombeau's algebra of generalized functions on $M$ is defined as the quotient:

$$G(M) \triangleq E_M(M)/N(M) \tag{1.1.1}$$

of the space $E_M(M)$ of sequences of moderate growth modulo the space $N(M)$ of negligible sequences. More precisely the notions of moderateness resp. negligibility are defined by the following asymptotic estimates (where $X(M)$ denoting the space of smooth vector fields on $M$):

$$E_M(M) \triangleq \Big\{ (u_\varepsilon)_\varepsilon | \ \forall K(K \subsetneqq M) \ \forall k(k \in \mathbb{N}) \exists N(N \in \mathbb{N})$$

$$\forall \xi_1, \ldots, \xi_k (\xi_1, \ldots, \xi_k \in X(M)) \left[ \sup_{p \in K} |L_{\xi_1} \ldots L_{\xi_k} u_\varepsilon(p)| = O(\varepsilon^{-N}) \text{ as } \varepsilon \to 0 \right] \Big\}, \tag{1.1.2}$$

$$N(M) \triangleq \Big\{ (u_\varepsilon)_\varepsilon | \ \forall K(K \subsetneqq M), \ \forall k(k \in \mathbb{N}_0) \forall q(q \in N)$$

$$\forall \xi_1, \ldots, \xi_k (\xi_1, \ldots, \xi_k \in X(M)) \left[ \sup_{p \in K} |L_{\xi_1} \ldots L_{\xi_k} u_\varepsilon(p)| = O(\varepsilon^q) \text{ as } \varepsilon \to 0 \right] \Big\}. \tag{1.1.3}$$

**Remark 1.1.1.** In the definition the Landau symbol $a_\varepsilon = O(\psi(\varepsilon))$ appears, having the following meaning: $\exists C(C > 0) \exists \varepsilon_0 (\varepsilon_0 \in (0,1]) \forall \varepsilon(\varepsilon < \varepsilon_0)[a_\varepsilon \leq C\psi(\varepsilon)]$.

**Definition 1.1.2.** Elements of $G(M)$ are denoted by:

$$u = \mathbf{cl}[(u_\varepsilon)_\varepsilon] \triangleq (u_\varepsilon)_\varepsilon + N(M). \tag{1.1.4}$$

**Remark 1.1.2.** With componentwise operations $(\cdot, \pm)$ $G(M)$ is a fine sheaf of differential algebras with respect to the Lie derivative defined by $L_\xi u \triangleq \mathbf{cl}[(L_\xi u_\varepsilon)_\varepsilon]$. The spaces of moderate resp. negligible sequences and hence the algebra itself may be characterized locally, i.e., $u \in G(M)$ iff $u \circ \psi_\alpha \in G(\psi_\alpha(V_\alpha))$ for all charts $(V_\alpha, \psi_\alpha)$, where on the open set $\psi_\alpha(V_\alpha) \subset \mathbb{R}^n$ in the respective estimates Lie derivatives are replaced by partial derivatives.

The spaces of moderate resp. negligible sequences and hence the algebra itself may be characterized locally, i.e., $u \in G(M)$ iff $u \circ \psi_\alpha \in G(\psi_\alpha(V_\alpha))$ for all charts $(V_\alpha, \psi_\alpha)$, where on the open set $\psi_\alpha(V_\alpha) \subset \mathbb{R}^n$ in the respective estimates Lie derivatives are replaced by partial derivatives.

**Remark 1.1.3.** Smooth functions $f \in C^\infty(M)$ are embedded into $G(M)$ simply by the "constant" embedding $\sigma$, i.e., $\sigma(f) = \mathbf{cl}[(f)_\varepsilon]$, hence $C^\infty(M)$ is a faithful subalgebra of G($M$).

# I.1.1. Point Values of Generalized Functions on $M$. Generalized Numbers.

Within the classical distribution theory, distributions cannot be characterized by their point values in any way similar to classical functions. On the other hand, there is a very natural and direct way of obtaining the point values of the elements of Colombeau's algebra: points are simply inserted into representatives. The objects so obtained are sequences of numbers, and as such are not the elements in the field $\mathbb{R}$ or $\mathbb{C}$. Instead, they are the representatives of *Colombeau's generalized numbers* [6],[7],[8]. We give the exact definition of these "numbers".

**Definition 1.1.3.** Inserting $p \in M$ into $u \in$ G($M$) yields a well defined element of the ring of constants (also called generalized numbers) K (corresponding to **K** = R resp. C), defined as the set of moderate nets of numbers $((r_\varepsilon)_\varepsilon \in \mathbf{K}^{(0,1]}$ with $|r_\varepsilon| = O(\varepsilon^{-N})$ for some $N$) modulo negligible nets ($|r_\varepsilon| = O(\varepsilon^m)$ for each $m$); componentwise insertion of points of $M$ into elements of G($M$) yields well-defined generalized numbers, i.e., elements of the ring of constants:

$$\mathrm{K} = \mathrm{E_c}(M)/\mathrm{N_c}(M) \qquad (1.1.8)$$

with $\mathrm{K} = \widetilde{\mathbb{R}}$ or $\mathrm{K} = \widetilde{\mathbb{C}}$ for $\mathbf{K} = \mathbb{R}$ or $\mathbf{K} = \mathbb{C}$, where

$$\mathrm{E_c}(M) = \left\{ (r_\epsilon)_\epsilon \in \mathbf{K}^I | \exists n (n \in \mathbb{N}) \big[ |r_\epsilon| = O(\epsilon^{-n}) \text{ as } \varepsilon \to 0 \big] \right\}$$
$$\mathrm{N_c}(M) = \left\{ (r_\epsilon)_\epsilon \in \mathbf{K}^I | \forall m (m \in \mathbb{N}) \big[ |r_\epsilon| = O(\epsilon^m) \text{ as } \varepsilon \to 0 \big] \right\} \qquad (1.1.9)$$
$$I = (0,1].$$

Generalized functions on $M$ are characterized by their generalized point values, i.e., by their values on points in $\tilde{M}_c$, the space of equivalence classes of compactly supported nets $(p_\varepsilon)_\varepsilon \in M^{(0,1]}$ with respect to the relation $p_\varepsilon \sim p'_\varepsilon :\Leftrightarrow d_h(p_\varepsilon, p'_\varepsilon) = O(\varepsilon^m)$ for all $m$, where $d_h$ denotes the distance on $M$ induced by any Riemannian metric.

**Definition 1.1.3.** For $u \in$ G($M$) and $x_0 \in M$, the point value of $u$ at the point $x_0, u(x_0)$, is defined as the class of $(u_\varepsilon(x_0))_\varepsilon$ in K.

**Definition 1.1.4.** We say that an element $r \in$ K is *strictly nonzero* if there exists a representative $(r_\varepsilon)_\varepsilon$ and a $q \in \mathbb{N}$ such that $|r_\varepsilon| \geqslant \varepsilon^q$ for $\varepsilon$ sufficiently small. If $r$ is strictly nonzero, then it is also invertible with the inverse $[(1/r_\varepsilon)_\varepsilon]$. The converse is true as well.

Treating the elements of Colombeau algebras as a generalization of classical functions, the question arises whether the definition of point values can be extended in such a way that each element is characterized by its values. Such an extension is indeed possible.

**Definition 1.1.5**. Let $\Omega$ be an open subset of $\mathbb{R}^n$. On a set $\hat{\Omega}$ :

$$\hat{\Omega} = \left\{ (x_\varepsilon)_\varepsilon \in \Omega^I | \exists p(p > 0)[|x_\varepsilon| = O(\varepsilon^p)] \right\} =$$

$$\left\{ (x_\varepsilon)_\varepsilon \in \Omega^I | \exists p(p > 0) \exists \varepsilon_0(\varepsilon_0 > 0) \Big[ |x_\varepsilon| \leq \varepsilon^p, \text{ for } 0 < \varepsilon < \varepsilon_0 \Big] \right\}$$

(1.1.10)

we introduce an equivalence relation:

$$(x_\varepsilon)_\varepsilon \sim (y_\varepsilon)_\varepsilon \iff \forall q(q > 0) \forall \varepsilon(\varepsilon > 0) \Big[ |x_\varepsilon - y_\varepsilon| \leq \varepsilon^q, \text{ for } 0 < \varepsilon < \varepsilon_0 \Big] \quad (1.1.11)$$

and denote by $\widetilde{\Omega} = \hat{\Omega}/\sim$ the set of generalized points. The set of points with compact support is

$$\widetilde{\Omega}_c =$$

$$\left\{ \widetilde{x} = \mathbf{cl}[(x_\varepsilon)_\varepsilon] \in \widetilde{\Omega} | \exists K(K \subset \Omega) \exists \varepsilon_0(\varepsilon_0 > 0) \Big[ x_\varepsilon \in K \text{ for } 0 < \varepsilon < \varepsilon_0 \Big] \right\}$$

(1.1.12)

**Definition 1.1.6**. A generalized function $u \in G(M)$ is called associated to zero, $u \approx 0$ on $\Omega \subseteq M$ in *L. Schwartz's sense* if one (hence any) representative $(u_\epsilon)_\epsilon$ converges to zero weakly, i.e.

$$w\text{-}\lim_{\epsilon \to 0} u_\epsilon = 0 \qquad (1.1.13)$$

We shall often write:

$$u \underset{\text{Sch}}{\approx} 0.$$ 

<div align="right">(1.1.14)</div>

## I.1.2. Colombeau algebra of compactly supported generalized functions and Colombeau algebra of tempered generalized functions.

**Definition 1.1.7**. We define the spaces $E_{\mathbf{c},M}(\Omega)$ and $N_{\mathbf{c}}(\Omega)$ as the sets of all nets $(u_\varepsilon)_\varepsilon$ in $E_M(\Omega)$ and $N(\Omega)$ respectively such that

$$\exists K(K \Subset \Omega) \forall \varepsilon(\varepsilon \in (0,1]) \big[ \mathbf{supp}(u_\varepsilon) \subseteq \Omega \big]. \tag{1.1.15}$$

**Definition 1.1.8**. The Colombeau algebra $G_{\mathbf{c}}(\Omega)$ is isomorphic to the factor space $E_{\mathbf{c},M}(\Omega)/N_{\mathbf{c}}(\Omega)$ by canonical map

$$l : G_{\mathbf{c}}(\Omega) \to E_{\mathbf{c},M}(\Omega)/N_{\mathbf{c}}(\Omega) : (u_\varepsilon)_\varepsilon + N(\Omega) \to (\psi_u u_\varepsilon)_\varepsilon + N_{\mathbf{c}}(\Omega), \tag{1.1.16}$$

where $\psi_u = \psi u$, $\psi \in C_{\mathbf{c}}^\infty(\Omega)$ and $\psi$ is identically 1 in a neighborhood of $\mathbf{supp}(u)$.

**Definition 1.1.9**. We call the elements of the set

$$E_\tau(\mathbb{R}^n) \triangleq$$

$$\{(u_\varepsilon)_\varepsilon \in E_\tau[\mathbb{R}^n] | \forall \alpha \in \mathbb{N}^n \exists N \in \mathbb{N}$$

$$\left[\sup_{x \in \mathbb{R}^n} (1 + |x|)^{-N} |\partial^\alpha u_\varepsilon(x)| = O(\varepsilon^{-N}) \text{ as } \varepsilon \to 0\right]\}, \qquad (1.1.17)$$

$$E_\tau[\mathbb{R}^n] \triangleq O_M(\mathbb{R}^n)^{(0,1]}$$

$\tau$-moderate.
We call the elements of the set

$$N_\tau(\mathbb{R}^n) \triangleq$$

$$\{(u_\varepsilon)_\varepsilon \in E_\tau[\mathbb{R}^n] | \forall \alpha \in \mathbb{N}^n \exists N \in \mathbb{N} \forall q \in \mathbb{N}$$

$$\left[\sup_{x \in \mathbb{R}^n} (1 + |x|)^{-N} |\partial^\alpha u_\varepsilon(x)| = O(\varepsilon^{-N}) \text{ as } \varepsilon \to 0\right]\}, \qquad (1.1.18)$$

$\tau$-negligible.

**Definition 1.1.10.**[2].The Colombeau algebra of tempered generalized functions
is defined as the quotient

$$G_\tau(\mathbb{R}^n) \triangleq E_\tau(\mathbb{R}^n)/N_\tau(\mathbb{R}^n). \qquad (1.1.19)$$

**Theorem 1.1.1.** $(u_\varepsilon)_\varepsilon \in E_\tau(\mathbb{R}^n)$ is $\tau$-negligible if and only if

$$\exists N(N \in \mathbb{N}) \forall q(q \in \mathbb{N}) \left[\sup_{x \in \mathbb{R}^n} (1 + |x|)^{-N} |\partial^\alpha u_\varepsilon(x)| = O(\varepsilon^q) \text{ as } \varepsilon \to 0\right]. \qquad (1.1.20)$$

## I.1.3. Algebra of generalized functions based on $\mathscr{L}(\mathbb{R}^n)$. Colombeau algebra $G_{\tau,\mathscr{L}}(\mathbb{R}^n)$.

**Definition 1.1.11.** Let $\mathscr{L}[\mathbb{R}^n] \triangleq \mathscr{L}(\mathbb{R}^n)^{(0,1]}$. The elements of

$$E_{\mathscr{L}}[\mathbb{R}^n] \triangleq$$

$$\{(u_\varepsilon)_\varepsilon \in \mathscr{L}[\mathbb{R}^n] | \forall \alpha \forall \beta (\alpha, \beta \in \mathbb{N}^n)$$

$$\exists N(N \in \mathbb{N}) \left[ \sup_{x \in \mathbb{R}^n} |x^\alpha \partial^\beta u_\varepsilon(x)| O(\varepsilon^{-N}) \text{ as } \varepsilon \to 0 \right] \} \tag{1.1.21}$$

are called $\mathscr{L}$-moderate. The elements of

$$N_{\mathscr{L}}[\mathbb{R}^n] \triangleq$$

$$\{(u_\varepsilon)_\varepsilon \in \mathscr{L}[\mathbb{R}^n] | \forall \alpha \forall \beta (\alpha, \beta \in \mathbb{N}^n)$$

$$\forall q(q \in \mathbb{N}) \in \mathbb{N} \left[ \sup_{x \in \mathbb{R}^n} |x^\alpha \partial^\beta u_\varepsilon(x)| O(\varepsilon^q) \text{ as } \varepsilon \to 0 \right] \} \tag{1.1.22}$$

are called $\mathscr{L}$-negligible.

**Definition 1.1.12.** The algebra of generalized functions based on $\mathscr{L}(\mathbb{R}^n)$ is defined as the factor space

$$G_{\mathscr{L}}(\mathbb{R}^n) \triangleq E_{\mathscr{L}}[\mathbb{R}^n]/N_{\mathscr{L}}[\mathbb{R}^n]. \tag{1.1.23}$$

**Definition 1.1.13.**By factorizing $E_\tau(\mathbb{R}^n)$ with respect to the ideal $N_{\mathscr{L}}[\mathbb{R}^n]$ Thus we obtain the quotient algebra

$$G_{\tau,\mathscr{L}}(\mathbb{R}^n) \triangleq E_\tau(\mathbb{R}^n)/N_{\mathscr{L}}[\mathbb{R}^n]. \qquad (1.1.24)$$

It is clearly that $G_{\mathscr{L}}(\mathbb{R}^n)$ is contained in $G_{\tau,\mathscr{L}}(\mathbb{R}^n)$.

# I.1.4.Regular Colombeau algebras.

**Definition 1.1.14.**In the setting of Colombeau theory, the notion of regularity is based on the subalgebra $G^\infty(\Omega)$ of regular generalized functions in $G(\Omega)$. It is defined by those elements which have a representative $(u_\varepsilon)_\varepsilon$ satisfying

$$\forall K(K \Subset \Omega)\exists N(N \in \mathbb{N})\forall \alpha(\alpha \in \mathbb{N}^n)\left[\sup_{x \in K} |\partial^\alpha u_\varepsilon(x)| = O(\varepsilon^{-N}) \text{ as } \varepsilon \to 0\right]. \qquad (1.1.25)$$

**Remark 1.1.4.**Note that (1.1.25) is obtained by interchanging the quantifiers in the definition of $E_M(\Omega)$. It well known that

$$G^\infty(\Omega) \cap D'(\Omega) = C^\infty(\Omega). \qquad (1.1.26)$$

**Definition 1.1.15.**The Colombeau algebra $G_c^\infty(\Omega)$ of regular generalized functions with compact support is defined by intersection

$$G_c^\infty(\Omega) \triangleq G^\infty(\Omega) \cap G_c(\Omega). \qquad (1.1.27)$$

**Remark 1.1.5.** It well known that $G_c^\infty(\Omega) \cap E'(\Omega) = C_c^\infty(\Omega)$.

**Definition 1.1.16.** The Colombeau algebra $G_{\mathcal{L}}^\infty(\Omega)$ is given by the quotient

$$G_{\mathcal{L}}^\infty(\mathbb{R}^n) \triangleq E_{\mathcal{L}}^\infty(\mathbb{R}^n)/N_{\mathcal{L}}(\mathbb{R}^n), \qquad (1.1.28)$$

where

$$E_{\mathcal{L}}^\infty[\mathbb{R}^n] \triangleq$$

$$\{(u_\varepsilon)_\varepsilon \in \mathcal{L}[\mathbb{R}^n] | \exists N(N \in \mathbb{N}) \forall \alpha \forall \beta(\alpha, \beta \in \mathbb{N}^n)$$

$$\left[ \sup_{x \in \mathbb{R}^n} |x^\alpha \partial^\beta u_\varepsilon(x)| O(\varepsilon^{-N}) \text{ as } \varepsilon \to 0 \right] \}. \qquad (1.1.29)$$

# I.1.5. Integration in Colombeau theory

**Definition.1.1.17.** Let $K \subset \Omega$ be a compact set and $u \in \mathrm{G}(M)$. The linear map from $\mathrm{G}(M)$ into $\widetilde{\mathbb{C}}$ via formula

$$u \to \int_K u(x) dx \triangleq \left( \int_K u_\varepsilon(x) dx \right)_\varepsilon + N \qquad (1.1.30)$$

is called the Colombeau integral of $u$. Colombeau integral is $\widetilde{\mathbb{C}}$ linear by definition, i.e., $\forall \alpha_1, \alpha_2 \in \widetilde{\mathbb{C}}, \forall u_1, u_2 \in \mathrm{G}(M)$

$$\int_K [\alpha_1 u_1(x) + \alpha_2 u_2(x)] dx = \alpha_1 \int_K u_1(x) dx + \alpha_2 \int_K u_2(x) dx \qquad (1.1.31)$$

**Definition**.**1**.**1**.**18**. Let $u \in G(\Omega)$ with compact support $K \Subset \Omega$. We define

$$\int_\Omega u(x)dx \triangleq \int_{\widehat{K}} u(x)dx, \qquad (1.1.32)$$

where $\widehat{K}$ is any compact subset of $\Omega$ containing $K$ in its interior.

**Theorem**.**1**.**1**.**2**. (1) If $u, v \in G(\Omega)$ and either $u$ or $v$ has compact support then for any
$\alpha \in \mathbb{N}^n$

$$\int_\Omega u(x)\partial^\alpha v(x)dx = (-1)^{|\alpha|}\int_\Omega (\partial^\alpha u(x))v(x)dx \qquad (1.1.33)$$

(2) If $u \in G_{\mathscr{L}}(\mathbb{R}^n)$ then integral

$$\int_{\mathbb{R}^n} u(x)dx \triangleq \mathbf{cl}\left[\left(\int_{\mathbb{R}^n} u_\varepsilon(x)dx\right)_\varepsilon\right] \qquad (1.1.34)$$

is a well-defined element of $\widetilde{\mathbb{C}}$. Moreover,the following properties hold:

(i)

$$\int_{\mathbb{R}^n} (\iota g)(x)dx = \mathbf{cl}\left[\left(\int_{\mathbb{R}^n} g(x)dx\right)_\varepsilon\right] \qquad (1.1.34)$$

(ii)

$$\int_{\mathbb{R}^n} u(x)\partial^\alpha v(x)dx = (-1)^{|\alpha|}\int_{\mathbb{R}^n} (\partial^\alpha u(x))v(x)dx \qquad (1.1.35)$$

for all $u, v \in G_{\mathscr{L}}(\mathbb{R}^n)$ and $\alpha \in \mathbb{N}^n$.

**Definition.1.1.19**.Let $u \in G_{\tau,\mathcal{L}}(\mathbb{R}^n) \cup G_{\tau}(\mathbb{R}^n)$ and $\varphi$ be a mollifier in $\mathcal{L}(\mathbb{R}^n)$. We define the $\widehat{\varphi}$-integral of $u$ as the complex generalized number by

$$\int_{\mathbb{R}^n} u(x) d_{\widehat{\varphi}} x \triangleq \int_{\mathbb{R}^n} u_{\widehat{\varphi}}(x) dx = \mathbf{cl}\left[\left(\int_{\mathbb{R}^n} u_\varepsilon(x)\widehat{\varphi}_\varepsilon(x) dx\right)\right], \qquad (1.1.36)$$

where $u_{\widehat{\varphi}}(x) \triangleq \left(u_\varepsilon(x)\widehat{\varphi}_\varepsilon(x)\right)_\varepsilon + N_{\mathcal{L}}(\mathbb{R}^n)$.

**Theorem.1.1.3**. (1) If $u \in G_{\mathbf{c}}(\mathbb{R}^n) \subseteq G_{\tau,\mathcal{L}}(\mathbb{R}^n)$ then

$$\int_{\Omega} u(x) dx = \int_{\mathbb{R}^n} u(x) d_{\widehat{\varphi}} x. \qquad (1.1.37)$$

(2) If $u \in G_{\tau,\mathcal{L}}(\mathbb{R}^n)$ and $v \in G_{\mathcal{L}}^\infty(\mathbb{R}^n) \subseteq G_{r,\mathcal{L}}(\mathbb{R}^n)$ then

$$\int_{\mathbb{R}^n} u(x) v(x) dx = \int_{\mathbb{R}^n} u(x) v(x) d_{\widehat{\varphi}} x. \qquad (1.1.38)$$

(3) If $u \in G_{\tau,\mathcal{L}}(\mathbb{R}^n)$ and $v \in G_{\mathcal{L}}^\infty(\mathbb{R}^n)$ then for all $\alpha \in \mathbb{N}^\alpha$

$$\int_{\mathbb{R}^n} (\partial^\alpha u(x)) v(x) dx = (-1)^{|\alpha|} \int_{\mathbb{R}^n} u(x)(\partial^\alpha v(x)) d_{\widehat{\varphi}} x. \qquad (1.1.39)$$

**Definition.1.1.20**.Let $u,v$ be generalized functions in $G(\Omega)$. $u$ and $v$ are said to be equal in the generalized distributions sense ($u =_{\mathbf{g.d.}} v$ for short) iff

$$\int_{\mathbb{R}^n} [u(x) - v(x)][\iota(f)](x) dx = 0. \qquad (1.1.40)$$

for all $f \in C_c^\infty(\Omega)$.

**Definition**.**1**.**1**.**21**.Let $u, v$ be generalized functions in $G_{\tau,\mathcal{L}}(\Omega)$. $u$ and $v$ are said to be equal in the generalized distributions sense ($u =_{\mathbf{g.t.d.}} v$ for short)

$$\int_{\mathbb{R}^n} [u(x) - v(x)][\iota(f)](x)dx = 0. \qquad (1.1.41)$$

iff for all $f \in \mathcal{L}(\mathbb{R}^n)$.

**Theorem**.**1**.**1**.**4**. Let us consider the expression

$$\left(\widehat{k}u\right)(x) = \int_\Omega k(x,y)u(y)dy = 0. \qquad (1.1.42)$$

(1) If $k \in G(\Omega \times \Omega)$ then (1.1.42) defines a $\widetilde{\mathbb{C}}$-linear map $\widehat{k} : G_c(\Omega) \to G(\Omega)$ : $u \to \widehat{k}u$, where $\widehat{k}u$ is the generalized function with representative

$$\left(\left(\widehat{k}_\varepsilon u_\varepsilon\right)(x)\right)_\varepsilon = \left(\int_\Omega k_\varepsilon(x,y)u_\varepsilon(y)dy\right)_\varepsilon. \qquad (1.1.43)$$

(2) If $k \in G^\infty(\Omega \times \Omega)$ then $\widehat{k}$ maps $G_c(\Omega)$ into $G^\infty(\Omega)$.
(3) If $k \in G_c(\Omega \times \Omega)$ then $\widehat{k}$ maps $G(\Omega)$ into $G_c(\Omega)$.
(4) If $k \in G_c^\infty(\Omega \times \Omega)$ then $\widehat{k}$ maps $G(\Omega)$ into $G_c^\infty(\Omega)$.
(5) If $k \in G(\Omega \times \Omega)$ and $\pi_1, \pi_2 : \mathbf{supp}(k) \to \Omega$ are proper then $\mathbf{supp}\left(\widehat{k}u\right) \Subset \Omega$ for all $u \in G_c(\Omega)$ and $k$ can be uniquely extended to a $\widetilde{\mathbb{C}}$-linear map from $G(\Omega)$ into $G(\Omega)$ such that for all $u \in G(\Omega)$ and $v \in G_c(\Omega)$

$$\int_\Omega \left(\widehat{k}u\right)(x)v(x)dx = \int_\Omega v(y)\left({}^t\widehat{k}v\right)(y)dy,$$

(1.1.44)

$$\left({}^t\widehat{k}v\right)(y) = \int_\Omega \widehat{k}(x,y)v(x)dx.$$

(6) If $k \in G^\infty(\Omega \times \Omega)$ and $\pi_1, \pi_2 : \mathbf{supp}(k) \to \Omega$ are proper then the extension defined above maps $G(\Omega)$ into $G^\infty(\Omega)$.

(7) Let us now consider the expression

$$\left(\widehat{k}u\right)(y) = \int_{\mathbb{R}^n} \widehat{k}(x,y)u(y)dy.$$

(1.1.45)

If $k \in G_\tau(\mathbb{R}^{2n})$ has a representative $(k_\varepsilon)_\varepsilon$ such that $(k_\varepsilon)_\varepsilon$ satisfying the condition

$$\forall \alpha(\alpha \in \mathbb{N}^n)\forall s(s \in \mathbb{N})\exists N(N \in \mathbb{N})$$

(1.1.46)

$$\left[\sup_{x \in \mathbb{R}^n}(1+|x|)^{-N}\sup_{y \in \mathbb{R}^n, \beta \le s}(1+|y|)^{-s}\left|\partial_x^\alpha\partial_y^\beta k_\varepsilon(x,y)\right| = O(\varepsilon^{-N}) \text{ as } \varepsilon \to 0\right],$$

then $\widehat{k}$ maps $G_\tau(\mathbb{R}^n)$ into $G_\tau(\mathbb{R}^n)$.

(8) If $k \in G_{\mathcal{L}}(\mathbb{R}^{2n})$ then $\widehat{k}$ maps $G_\tau(\mathbb{R}^n)$ into $G_{\mathcal{L}}(\mathbb{R}^n)$.

(9) If $k \in G_{\mathcal{L}}^\infty(\mathbb{R}^{2n})$ then $\widehat{k}$ maps $G_\tau(\mathbb{R}^n)$ into $G_{\mathcal{L}}^\infty(\mathbb{R}^n)$.

# I.1.6. Fourier transform in Colombeau theory

**Definition**.1.1.22.Let $u \in G_{\mathcal{L}}(\mathbb{R}^n)$. The Fourier transform $\mathscr{F}u$ and the inverse Fourier transform $\mathscr{F}^*u$ of $u$ are defined by the corresponding transformations at the level of representatives.

Thus $\mathcal{F}u$ and $\mathcal{F}^*u$ map $G_{\mathcal{L}}(\mathbb{R}^n)$ into itself and extend the classical Fourier transform on $\mathcal{L}(\mathbb{R}^n)$.

**Theorem**.**1**.**1**.**5**. For all $u, v \in G_{\mathcal{L}}(\mathbb{R}^n)$ and $\alpha \in \mathbb{N}^n$

(1)

$$\int_{\mathbb{R}^n} (\mathcal{F}u)(x)v(v)dx = \int_{\mathbb{R}^n} u(x)(\mathcal{F}v)(v)dx,$$

$$(1.1.47)$$

$$\int_{\mathbb{R}^n} (\mathcal{F}^*u)(x)v(v)dx = \int_{\mathbb{R}^n} u(x)(\mathcal{F}^*v)(v)dx.$$

(2)

$$u = \mathcal{F}[\mathcal{F}^*u] = \mathcal{F}^*[\mathcal{F}u]. \qquad (1.1.48)$$

(3)

$$\mathcal{F}[\iota(y^\alpha)u] = i^{|\alpha|}\partial^\alpha\mathcal{F}[u],$$

$$(1.1.49)$$

$$\mathcal{F}^*[\iota(y^\alpha)u] = (-i)^{|\alpha|}\partial^\alpha\mathcal{F}^*[u].$$

(4)

$$(-i)^{|\alpha|}\mathcal{F}[\partial^\alpha u] = \iota(y^\alpha)\mathcal{F}u,$$

$$(1.1.50)$$

$$i^{|\alpha|}\mathcal{F}^*[\partial^\alpha u] = \iota(y^\alpha)\mathcal{F}^*u,$$

where $\iota(y^\alpha)$ is intended as the embedding of $y^\alpha$ in $G_{\tau,\mathcal{L}}(\mathbb{R}^n)$ and the product $\iota(y^\alpha)u \in G_{\mathcal{L}}(\mathbb{R}^n)$.

**Remark 1**.**1**.**6**.Note that the generalized function $\iota(y^\alpha)$ has $y^\alpha$ as a representative since $(\varphi_\varepsilon * y^\alpha)(x) = x^\alpha$.

**Definition**.**1**.**1**.**23**. Let $u \in G_{\tau,\mathcal{L}}(\mathbb{R}^n)$ and $\varphi$ be a mollifier in $\mathcal{L}(\mathbb{R}^n)$. We define the $\widehat{\varphi}$-Fourier transform of $u$ as the generalized function in $G_{\mathcal{L}}(\mathbb{R}^n)$ given by

$$\mathcal{F}_{\widehat{\varphi}}[u] \triangleq \mathcal{F}\left[u_{\widehat{\varphi}}\right] \tag{1.1.51}$$

and the inverse $\widehat{\varphi}$-Fourier transform of $u$ as

$$\mathcal{F}_{\widehat{\varphi}}^*[u] \triangleq \mathcal{F}^*\left[u_{\widehat{\varphi}}\right]. \tag{1.1.52}$$

**Remark 1**.**1**.**7**. Note that $\mathcal{F}_{\widehat{\varphi}}[u]$ and $\mathcal{F}_{\widehat{\varphi}}^*[u]$ have representatives of the form

$$\int_{\mathbb{R}^n} (\exp[-iy\xi])u_\varepsilon(y)\widehat{\varphi}_\varepsilon(y)dy,$$

$$\int_{\mathbb{R}^n} (\exp[iy\xi])u_\varepsilon(\xi)\widehat{\varphi}_\varepsilon(\xi)dy \tag{1.1.53}$$

respectively, with $d\xi = (2\pi)^{-n}d\xi$.

**Theorem**.**1**.**1**.**6**. (1) $\mathcal{F}_{\widehat{\varphi}}[u] = \mathcal{F}[u]$ on $G_{\mathcal{L}}^\infty(\mathbb{R}^n)$ and $G_{\mathbf{c}}(\Omega)$.
(2) If $u \in G_{\tau,\mathcal{L}}(\mathbb{R}^n)$ and $v \in G_{\mathcal{L}}^\infty(\mathbb{R}^n)$ then

$$\int_{\mathbb{R}^n} \mathcal{F}_{\widehat{\varphi}}[u](x)v(x)dx = \int_{\mathbb{R}^n} u(x)\mathcal{F}_{\widehat{\varphi}}[v](x)dx. \tag{1.1.54}$$

(3) If $u \in G_{\tau,\mathcal{L}}(\mathbb{R}^n)$ and $v \in G_{\mathcal{L}}^\infty(\mathbb{R}^n)$ then

$$\int_{\mathbb{R}^n} \mathcal{F}_{\widehat{\varphi}}^*[u](x)v(x)dx = \int_{\mathbb{R}^n} u(x)\mathcal{F}_{\widehat{\varphi}}^*[v](x)dx. \tag{1.1.55}$$

(4) If $u \in G_{\tau,\mathcal{L}}(\mathbb{R}^n)$ then

$$\mathscr{F}_{\widehat{\varphi}}[\iota(y^\alpha)u] = i^{|\alpha|}\partial^\alpha \mathscr{F}_{\widehat{\varphi}}[u],$$

(1.1.56)

$$\mathscr{F}^*_{\widehat{\varphi}}[\iota(y^\alpha)] = (-i)^{|\alpha|}\partial^\alpha \mathscr{F}^*_{\widehat{\varphi}}[u].$$

(5) If $u \in G_{\tau,\mathcal{L}}(\mathbb{R}^n)$ and $v \in G^\infty_{\mathcal{L}}(\mathbb{R}^n)$ then

$$(-i)^{|\alpha|} \int_{\mathbb{R}^n} \mathscr{F}_{\widehat{\varphi}}[\partial^\alpha u](x)v(x)dx = \int_{\mathbb{R}^n} \left(\iota(y^\alpha)\mathscr{F}_{\widehat{\varphi}}[u]\right)(x)v(x)dx,$$

(1.1.57)

$$i^{|\alpha|} \int_{\mathbb{R}^n} \mathscr{F}^*_{\widehat{\varphi}}[\partial^\alpha u](x)v(x)dx = \int_{\mathbb{R}^n} \left(\iota(y^\alpha)\mathscr{F}^*_{\widehat{\varphi}}[u]\right)(x)v(x)dx.$$

# I.1.7.Topological $\widetilde{\mathbb{C}}$-modules

**Definition**.**1.1.24**. An abelian group $(G,+)$ is a $\widetilde{\mathbb{C}}$-module if there is given a (product) map $\widetilde{\mathbb{C}} \times G \to G : (\lambda,u) \to \lambda u$ such that for all $\lambda,\mu \in \widetilde{\mathbb{C}}$ and $u,v \in G$ :

$$(\lambda + \mu)u = \lambda u + \mu u,$$

$$\mu(u+v) = \mu u + \mu v,$$

(1.1.58)

$$\lambda_1(\lambda_2)u = (\lambda_1\lambda_2)u,$$

$$1u = u.$$

**Definition**.**1**.**1**.**25**.Let us now introduce the function $\mathbf{v}[(u_\varepsilon)_\varepsilon]$ on $E(M)$:

$$\mathbf{v} : E(M) \to (-\infty, +\infty] : (u_\varepsilon)_\varepsilon \to \sup\left\{\gamma \in \mathbb{R} \| u_\varepsilon| = O(\varepsilon^\gamma) \text{ as } \varepsilon \to 0\right\} \qquad (1.1.59)$$

It satisfies the following properties:
(1) $\mathbf{v}[(u_\varepsilon)_\varepsilon] = -\infty$ iff $(u_\varepsilon)_\varepsilon \in N(M)$,
(2) $\mathbf{v}[(u_\varepsilon)_\varepsilon (v_\varepsilon)_\varepsilon] \geq \mathbf{v}[(u_\varepsilon)_\varepsilon] + \mathbf{v}[(v_\varepsilon)_\varepsilon]$,
(3) $\mathbf{v}[(u_\varepsilon)_\varepsilon + \mathbf{v}(v_\varepsilon)_\varepsilon] \geq \min\{\mathbf{v}[(u_\varepsilon)_\varepsilon], \mathbf{v}[(v_\varepsilon)_\varepsilon]\}$.

**Definition**.**1**.**1**.**26**.Let us define the valuation of the complex generalized numbers $u = \mathbf{cl}[(u_\varepsilon)_\varepsilon]$

$$\mathbf{v}_{\widetilde{\mathbb{C}}}(u) \triangleq \mathbf{v}[(u_\varepsilon)_\varepsilon] \qquad (1.1.60)$$

Note that all the previous properties hold for the elements of $\widetilde{\mathbb{C}}$.

**Definition**.**1**.**1**.**27**. We let

$$|\cdot|_{\mathbf{e}} \triangleq \widetilde{\mathbb{C}} \to [0, +\infty) : u \to |u|_{\mathbf{e}} = \exp\left[-\mathbf{v}_{\widetilde{\mathbb{C}}}(u)\right] \qquad (1.1.61)$$

**Remark 1**.**1**.**8**.Note that the properties of the valuation on $\widetilde{\mathbb{C}}$ makes the coarsest topology on $\widetilde{\mathbb{C}}$ such that the map $|\cdot|_{\mathbf{e}}$ is continuous compatible with the ring structure.

**Definition**.**1**.**1**.**28**.Let $G$ be a $\widetilde{\mathbb{C}}$-module. A valuation on $G$ is a function $\mathbf{v} : G \to (-\infty, +\infty]$ such that for all $\lambda \in \widetilde{\mathbb{C}}$ and $u, v \in G$ :

(1) $\mathbf{v}(0) = +\infty$,
(2) $\mathbf{v}(\lambda u) \geq \mathbf{v}_{\widetilde{\mathbb{C}}}(\lambda) + \mathbf{v}(u)$,
(2′) $\mathbf{v}(\lambda u) \geq \mathbf{v}_{\widetilde{\mathbb{C}}}(\lambda) + \mathbf{v}(u)$ for all $\lambda = \mathbf{cl}[(c\varepsilon^a)_\varepsilon], c \in \mathbb{C}, a \in \mathbb{R}$,

(3) $\mathbf{v}(u+v) \geq \min\{\mathbf{v}(u),\mathbf{v}(v)\}$.

**Definition**.**1**.**1**.**29**. An ultra-pseudo-seminorm on $G$ is a function $\wp : G \to [0,+\infty)$ such that for all $\lambda \in \widetilde{\mathbb{C}}$ and $u,v \in G$ :

(1) $\wp(0) = 0$,

(2) $\wp(\lambda u) \leq |\lambda|_\mathbf{e} \, \wp(u)$,

(2′) $\wp(\lambda u) = |\lambda|_\mathbf{e} \, \wp(u)$ for all $\lambda = \mathbf{cl}[(c\varepsilon^a)_\varepsilon], c \in \mathbb{C}, a \in \mathbb{R}$,

(3) $\wp(u+v) \leq \max\{\wp(u),\wp(v)\}$.

**Remark 1**.**1**.**9**. (**1**) We note that a typical example of an ultra-pseudo-seminorm obtained by means of a valuation $\mathbf{v}(\cdot)$ on $G$ is given via canonical formula

$$\wp(u) = \exp[-\mathbf{v}(u)]. \tag{1.1.62}$$

An ultra-pseudo-norm is an ultra-pseudo-seminorm $\wp(\cdot)$ such that $\wp(0) = 0$ implies $u = 0$.

(**2**) The function $|\cdot|_\mathbf{e}$ introduced in (1.1.61) is an ultra-pseudo-norm on $\widetilde{\mathbb{C}}$.

**Definition**.**1**.**1**.**30**.[9] (1) A subset $A$ of a $\widetilde{\mathbb{C}}$-module $G$ is $\widetilde{\mathbb{C}}$-absorbent if for all $u \in G$ there exists $a \in \mathbb{R}$ such that $u \in (\mathbf{cl}[(\varepsilon^b)_\varepsilon])A$ for all $b \leq a$.

(2) $A \subseteq G$ is $\widetilde{\mathbb{C}}$-balanced if $\lambda A \subseteq A$ for all $\lambda \in \widetilde{\mathbb{C}}$ with $|\lambda|_\mathbf{e} < 1$.

(3) $A \subseteq G$ is $\widetilde{\mathbb{C}}$-convex if $A + A \subseteq A$ and $(\mathbf{cl}[(\varepsilon^b)_\varepsilon])A \subseteq A$ all $b \geq 0$.

We now present another nontrivial example of a valuation on a $\widetilde{\mathbb{C}}$-module $G$.

**Definition**.**1**.**1**.**31**.[9]. A locally convex topological $\widetilde{\mathbb{C}}$-module is a topological $\widetilde{\mathbb{C}}$-module which has a base of $\widetilde{\mathbb{C}}$-convex neighborhoods of the origin.

**Theorem**.**1**.**1**.**7**.[9].Every locally convex topological $\widetilde{\mathbb{C}}$-module $G$ has a base of absolutely convex and absorbent neighborhoods of the origin.

**Theorem**.**1**.**1**.**8**.[9].Let $A$ be an absolutely convex and absorbent subset of a $\widetilde{\mathbb{C}}$-module $G$.Then

$$\mathbf{v}_A(u) = \sup\{\gamma \in \mathbb{R} | u \in (\mathbf{cl}[(\varepsilon^\gamma)])A\} \tag{1.1.63}$$

is a valuation on $G$. Moreover, for $\wp_A(u) = \exp[-\mathbf{v}_A(u)]$ and $\eta > 0$ the chain of

inclusion

$$\{u \in G | \wp_A(u) < \eta\} \subseteq (\mathbf{cl}[(\varepsilon^{-\ln \eta})_\varepsilon])A \subseteq \{u \in G | \wp_A(u) \leq \eta\} \qquad (1.1.64)$$

holds.

# I.1.9. Definition and basic properties of $G_E$ and $\widetilde{G}_E$

**Definition 1.1.33.** Let $E$ be a locally convex topological vector space topologized through the family of seminorms $\{p_i\}_{i \in I}$. The elements of

$$M_E \triangleq \left\{ (u_\varepsilon)_\varepsilon \in E^{(0,1]} | \forall i (i \in I) \exists N (N \in \mathbb{N}) \big[ p_i(u_\varepsilon) = O(\varepsilon^{-N}) \text{ as } \varepsilon \to 0 \big] \right\}, \qquad (1.1.65)$$

are called $E$-moderate. The elements of

$$N_E \triangleq \left\{ (u_\varepsilon)_\varepsilon \in E^{(0,1]} | \forall i (i \in I) \exists q (q \in \mathbb{N}) \big[ p_i(u_\varepsilon) = O(\varepsilon^q) \text{ as } \varepsilon \to 0 \big] \right\}, \qquad (1.1.66)$$

are called $E$-negligible. We define the space of generalized functions based on $E$ as the factor space

$$G_E \triangleq M_E / N_E. \qquad (1.1.67)$$

We note that the definition of $G_E$ does not depend on the family of seminorms which determines the locally convex topology of $E$. We adopt the notation $u = \mathbf{cl}[(u_\varepsilon)_\varepsilon]$ for the class $u$ of $(u_\varepsilon)_\varepsilon$ in $G_E$ and we embed $E$ into $G_E$ via the constant embedding $e \hookrightarrow \mathbf{cl}[(e_\varepsilon)_\varepsilon]$. By the properties of seminorms on E we may define the product between complex generalized numbers and elements of $G_E$ via the map

$$\widetilde{\mathbb{C}} \times G_E \to G_E : (\mathbf{cl}[(\lambda_\varepsilon)_\varepsilon], \mathbf{cl}[(u_\varepsilon)_\varepsilon]) \to \mathbf{cl}[(\lambda_\varepsilon u_\varepsilon)_\varepsilon], \qquad (1.1.68)$$

which equips $G_E$ with the structure of a $\widetilde{\mathbb{C}}$-module.

Since the growth in $\varepsilon$ of an $E$-moderate net is estimated in terms of any seminorm $p_i$ of $E$, it is natural to introduce the $p_i$-valuation of $(u_\varepsilon)_\varepsilon \in M_E$ as

$$\mathbf{v}_{p_i}[(u_\varepsilon)_\varepsilon] = \sup\big\{b \in \mathbb{R} \,|\, p_i(u_\varepsilon) = O(\varepsilon^b) \text{ as } \varepsilon \to 0\big\} \qquad (1.1.69)$$

Note that

$$\mathbf{v}_{p_i}[(u_\varepsilon)_\varepsilon] = \mathbf{v}_{\widetilde{\mathbb{C}}}[(p_i(u_\varepsilon))_\varepsilon] \qquad (1.1.70)$$

where the function $\mathbf{v}_{\widetilde{\mathbb{C}}}$ in (1.1.70) gives the valuation on $\widetilde{\mathbb{C}}$.

**Theorem.1.1.9**. The function $\mathbf{v}_{p_i}$ maps $M_E$ into $(-\infty, +\infty]$ and the following properties hold:

**(i)** $\mathbf{v}_{p_i}[(u_\varepsilon)_\varepsilon] = +\infty$ for all $i \in I$ iff $(u_\varepsilon)_\varepsilon \in N_E$,

**(ii)** $\mathbf{v}_{p_i}[(\lambda_\varepsilon u_\varepsilon)_\varepsilon] \geq \mathbf{v}[(\lambda_\varepsilon)_\varepsilon] + \mathbf{v}_{p_i}[(u_\varepsilon)_\varepsilon]$ for all $(\lambda_\varepsilon)_\varepsilon \in E_M, (u_\varepsilon)_\varepsilon \in M_E$,

**(ii)'** $\mathbf{v}_{p_i}[(\lambda_\varepsilon u_\varepsilon)_\varepsilon] = \mathbf{v}[(\lambda_\varepsilon)_\varepsilon] + \mathbf{v}_{p_i}[(u_\varepsilon)_\varepsilon]$ for all $(\lambda_\varepsilon)_\varepsilon = (c\varepsilon^b)_\varepsilon, c \in \mathbb{C}, (u_\varepsilon)_\varepsilon \in M_E$,

**(iii)** $\mathbf{v}_{p_i}[(w_\varepsilon + u_\varepsilon)_\varepsilon] \geq \inf\{\mathbf{v}_{p_i}[(w_\varepsilon)_\varepsilon], \mathbf{v}_{p_i}[(u_\varepsilon)_\varepsilon]\}$.

Property (i) combined with (iii) shows that $\mathbf{v}_{p_i}[(w_\varepsilon)_\varepsilon] = \mathbf{v}_{p_i}[(u_\varepsilon)_\varepsilon]$ iff $\mathbf{v}_{p_i}[(w_\varepsilon - u_\varepsilon)_\varepsilon]$ is $E$ − negligible.

Therefore we can use (1.1.70) for defining the $p_i$-valuation $\mathbf{v}_{p_i}[u] = \mathbf{v}_{p_i}[(u_\varepsilon)_\varepsilon]$ of a generalized function $u = \mathbf{cl}[(u_\varepsilon)_\varepsilon] \in G_E$. Thus

$$\wp_i[u] = \exp[-\mathbf{v}_{p_i}[u]] \qquad (1.1.71)$$

is an ultra-pseudo-seminorm on the $\widetilde{\mathbb{C}}$-module $G_E$.

**Theorem.1.1.10**.Let $E$ be a locally convex topological vector space.

(i) If $E$ is topologized through an increasing sequence $\{p_i\}_{i \in \mathbb{N}}$ of seminorms and

$$N_E = \left\{ (u_\varepsilon)_\varepsilon \in M_E \middle| \forall q \in \mathbb{N} \Big[ p_0(u_\varepsilon) = O(\varepsilon^q) \text{ as } \varepsilon \to 0 \Big] \right\} \tag{1.1.72}$$

then each $\wp_i$ is an ultra-pseudo-norm on $G_E$.

(ii) If $E$ has a countable base of neighborhoods of the origin then $G_E$ with the sharp topology is metrizable.

**Definition 1.1.34.** Let $E$ be an algebra on $\mathbb{C}$ with a family of seminorms $\{p_i\}_{i \in I}$ on $E$.

We recall that a $\widetilde{\mathbb{C}}$-module $G_E$ is a $\widetilde{\mathbb{C}}$-algebra $\widetilde{G}_E$ if there is given a multiplication $G_E \times G_E \to G_E : (u, w) \to uw$ such that $\forall u, w, v \in \widetilde{G}_E, \lambda \in \widetilde{\mathbb{C}}$ :

**(i)** $u(wv) = (uw)v$,

**(ii)** $u(w + v) = uw + uv$,

**(iii)** $(w + v)u = wu + vu$,

**(iv)** $\lambda(uv) = (\lambda u)v = u(\lambda v)$,

**(v)** $\forall i [i \in I]$ there exist finite subsets $I_0 \subset I, I'_0 \subset I$ and a constant $c_i > 0$ such that for all $u, v \in E$

$$p_i(uv) \leq c_i \left( \max_{j \in I_0} p_j(u) \right) \left( \max_{j \in I'_0} p_j(v) \right) \tag{1.1.73}$$

Then $\widetilde{G}_E$ with the sharp topology determined by the ultra-pseudo-seminorms $\{\wp_i\}_{i \in I}$

is a locally convex topological $\widetilde{\mathbb{C}}$-module and a topological $\widetilde{\mathbb{C}}$-algebra since from (1.1.73) it follows that $\forall i \in I, \forall u, v \in \widetilde{G}_E$ :

$$\wp_i(uv) \leq c_i \left( \max_{j \in I_0} \wp_j(u) \right) \left( \max_{j \in I'_0} \wp_j(v) \right) \tag{1.1.74}$$

Let us consider examples of the Colombeau algebras $\widetilde{G}_E$ obtained as $\widetilde{\mathbb{C}}$-modules $G_E$.

Particular choices of $E$ in Definition 1.1.33 give us well known algebras of generalized functions and the corresponding sharp topologies. This is of course the case for $E = \mathbb{C}, G_E = \widetilde{\mathbb{C}}$ and $\widetilde{G}_E = \widetilde{\mathbb{C}}$ which is an ultra-pseudo-normed $\widetilde{\mathbb{C}}$-module and more precisely a topological $\widetilde{\mathbb{C}}$-algebra.

Consider now an open subset $\Omega \subset \mathbb{R}^n$. $E = E(\Omega) = C^\infty(\Omega)$ topologized through the family of seminorms $p_{K_i,j}(f) = \sup_{x \in K_i, |\alpha| \leq j} |\partial^\alpha f(x)|$, where $K_0 K_1 \ldots K_i \ldots$ is a countable and exhausting sequence of compact subsets of $\Omega$, provides $\widetilde{G}_E = G_E(\Omega)$.

$G_E(\Omega)$ endowed with the sharp topology determined by $\left\{ \wp_{K_i,j}(\cdot) \right\}_{i,j \in \mathbb{N}}$ is a Fréchet $\widetilde{\mathbb{C}}$-module and a topological $\widetilde{\mathbb{C}}$-algebra. More examples of Fréchet $\widetilde{\mathbb{C}}$-modules which are also topological $\widetilde{\mathbb{C}}$-algebras are given by $G_E$ when E is $L(\mathbb{R}^n)$ or $W^{\infty,p}(\mathbb{R}^n), p \in [1, +\infty]$

In this way we construct the algebras $G_L(\mathbb{R}^n)$ and $G_{p,p}(\mathbb{R}^n)$ respectively, whose sharp topologies are obtained from $p_k(f) = \sup_{x \in \mathbb{R}^n, |\alpha| \leq k} |(1 + x^k)| \partial^\alpha f(x)||$ and $p_k(f) = \sup_{|\alpha| \leq k} \| \partial^\alpha f(x) \|_p$.

# I.1.10. Definition and basic properties of $G_B$ and $\widetilde{G}_B$.

**Definition 1.1.35.** Let $B$ be an Banach space topologized through the norm $\| \cdot \|$. The elements of

$$M_E \triangleq \left\{ (u_\varepsilon)_\varepsilon \in B^{(0,1]} | \forall i (i \in I) \exists N (N \in \mathbb{N}) \big[ \| u_\varepsilon \| = O(\varepsilon^{-N}) \text{ as } \varepsilon \to 0 \big] \right\}, \qquad (1.1.75)$$

are called $B$-moderate. The elements of

The elements of

$$N_B \triangleq \left\{ (u_\varepsilon)_\varepsilon \in B^{(0,1]} | \forall i (i \in I) \exists q (q \in \mathbb{N}) \big[ p_i(u_\varepsilon) = O(\varepsilon^q) \text{ as } \varepsilon \to 0 \big] \right\}, \qquad (1.1.76)$$

are called $B$-negligible. We define the space of generalized functions based on $B$ as the factor space

$$G_B \triangleq M_B/N_B. \qquad (1.1.77)$$

We adopt the notation $u = \mathbf{cl}[(u_\varepsilon)_\varepsilon]$ for the class $u$ of $(u_\varepsilon)_\varepsilon$ in $G_B$ and we embed $B$ into $G_B$ via the constant embedding $e \hookrightarrow \mathbf{cl}[(e_\varepsilon)_\varepsilon], e_\varepsilon = e$. By the properties of norm on $B$ we may define the product between complex generalized numbers and elements of $G_B$ via the map

$$\widetilde{\mathbb{C}} \times G_B \to G_B : (\mathbf{cl}[(\lambda_\varepsilon)_\varepsilon], \mathbf{cl}[(u_\varepsilon)_\varepsilon]) \to \mathbf{cl}[(\lambda_\varepsilon u_\varepsilon)_\varepsilon], \qquad (1.1.78)$$

which equips $G_B$ with the structure of a $\widetilde{\mathbb{C}}$-module.
Since the growth in $\varepsilon$ of an $B$-moderate net is estimated in terms of norm $\|\cdot\|$ of $B$, it is natural to introduce the $\|\cdot\|$-valuation of $(u_\varepsilon)_\varepsilon \in M_B$ as

$$\mathbf{v}_B[(u_\varepsilon)_\varepsilon] = \sup\left\{b \in \mathbb{R} | \|u_\varepsilon\| = O(\varepsilon^b) \text{ as } \varepsilon \to 0\right\} \qquad (1.1.79)$$

Note that

$$\mathbf{v}_B[(u_\varepsilon)_\varepsilon] = \mathbf{v}_{\widetilde{\mathbb{C}}}[(\|u_\varepsilon\|)_\varepsilon] \qquad (1.1.80)$$

where the function $\mathbf{v}_{\widetilde{\mathbb{C}}}$ in (1.1.80) gives the valuation on $\widetilde{\mathbb{C}}$.

**Theorem.1.1.11**. The function $\mathbf{v}_B$ maps $M_B$ into $(-\infty, +\infty]$ and the following properties hold:

**(i)** $\mathbf{v}_B[(u_\varepsilon)_\varepsilon] = +\infty$ iff $(u_\varepsilon)_\varepsilon \in N_B$,

**(ii)** $\mathbf{v}_B[(\lambda_\varepsilon u_\varepsilon)_\varepsilon] \geq \mathbf{v}_{\widetilde{\mathbb{C}}}[(\lambda_\varepsilon)_\varepsilon] + \mathbf{v}_B[(u_\varepsilon)_\varepsilon]$ for all $(\lambda_\varepsilon)_\varepsilon \in B_M, (u_\varepsilon)_\varepsilon \in M_B,$

**(ii$'$)** $\mathbf{v}_B[(\lambda_\varepsilon u_\varepsilon)_\varepsilon] = \mathbf{v}_{\widetilde{\mathbb{C}}}[(\lambda_\varepsilon)_\varepsilon] + \mathbf{v}_B[(u_\varepsilon)_\varepsilon]$ for all $(\lambda_\varepsilon)_\varepsilon = (c\varepsilon^b)_\varepsilon, c \in \mathbb{C}, (u_\varepsilon)_\varepsilon \in M_B,$

**(iii)** $\mathbf{v}_B[(w_\varepsilon + u_\varepsilon)_\varepsilon] \geq \inf\{\mathbf{v}_B[(w_\varepsilon)_\varepsilon], \mathbf{v}_B[(u_\varepsilon)_\varepsilon]\}.$

Property (i) combined with (iii) shows that $\mathbf{v}_B[(w_\varepsilon)_\varepsilon] = \mathbf{v}_B[(u_\varepsilon)_\varepsilon]$ iff $\mathbf{v}_B[(w_\varepsilon - u_\varepsilon)_\varepsilon]$ is
$B$-negligible.

Therefore we can use (1.1.70) for defining the $\|\cdot\|$-valuation $\mathbf{v}_B[u] = \mathbf{v}_B[\mathbf{cl}[(u_\varepsilon)_\varepsilon]]$ of
a generalized function $u = \mathbf{cl}[(u_\varepsilon)_\varepsilon] \in G_B.$ Thus

$$\wp_B[u] = \exp[-\mathbf{v}_B[u]] \tag{1.1.81}$$

is an ultra-pseudo-seminorm on the $\widetilde{\mathbb{C}}$-module $G_B.$

**Theorem**.**1**.**1**.**12**.Let $B$ be an Banach space.
If $B$ is topologized through of norm $\|\cdot\|$ and

$$N_B = \left\{(u_\varepsilon)_\varepsilon \in M_B | \forall q \in \mathbb{N}\big[\, \|u_\varepsilon\| = O(\varepsilon^q) \text{ as } \varepsilon \to 0 \,\big]\right\} \tag{1.1.82}$$

then $\wp_B$ is an ultra-pseudo-norm on $G_B.$

**Definition 1.1.36**.Let $B$ be an Banach algebra on $\mathbb{C}$ with norm $\|\cdot\|$ on $B.$
We recall that a $\widetilde{\mathbb{C}}$-module $G_B$ is a $\widetilde{\mathbb{C}}$-algebra $\widetilde{G}_E$ if there is given a multiplication
$G_B \times G_B \to G_B : (u,w) \to uw$ such that $\forall u,w,v \in \widetilde{G}_B, \lambda \in \widetilde{\mathbb{C}} :$

**(i)** $u(wv) = (uw)v,$

**(ii)** $u(w + v) = uw + uv,$

**(iii)** $(w + v)u = wu + vu,$

**(iv)** $\lambda(uv) = (\lambda u)v = u(\lambda v),$

Note that by definition of Banach algebra $B$ for all $u,v \in B$

$$\|uv\| \leq \|u\|\|v\|. \tag{1.1.83}$$

Therefore $\widetilde{G}_B$ with the sharp topology determined by the ultra-pseudo-seminorm
$\wp_B$
is a locally convex topological $\widetilde{\mathbb{C}}$-module and a topological $\widetilde{\mathbb{C}}$-algebra since from
(1.1.83) it follows that $\forall u,v \in \widetilde{G}_B :$

$$\wp_B(uv) \leq \wp_B(u)\,\wp_B(v). \tag{1.1.84}$$

Let us consider examples of the Colombeau algebras $\widetilde{G}_B$ obtained as $\widetilde{\mathbb{C}}$-modules $G_E$.

Particular choices of $B$ in Definition 1.1.35 give us well known algebras of generalized functions and the corresponding sharp topologies. This is of course the case for $B = \mathbb{C}, G_B = \widetilde{\mathbb{C}}$ and $\widetilde{G}_B = \widetilde{\mathbb{C}}$ which is an ultra-pseudo-normed $\widetilde{\mathbb{C}}$-module and more precisely a topological $\widetilde{\mathbb{C}}$-algebra.

**Definition 1.1.35.**(I) Consider now an open subset $\Omega \subseteq \mathbb{R}^n$. Suppose that $u \in L_p(\Omega)$ and there exist weak derivatives $\partial^\alpha u(x)$ for any $\alpha$ with $|\alpha| \leq l$ (all derivatives up to order $l$), such that $\partial^\alpha u \in L_p(\Omega)$. Then we say that $u \in W_p^l(\Omega)$. We introduce the (standard) norm in $W_p^l(\Omega)$ :

$$\|u\|_{W_p^l(\Omega)} = \left( \int_\Omega \sum_{|\alpha| \leq l} |\partial^\alpha u(x)| d^n x \right)^{1/p}. \tag{1.1.85}$$

We recall that (1) the norm $\sum_{|\alpha| \leq l} \|\partial^\alpha u(x)\|_{p,\Omega}$ is equivalent to the standard norm.

(2) $W_p^l(\Omega)$ is complete,i.e. $W_p^l(\Omega)$ is a Banach space.

(3) If $p = 2$, the space $W_2^l(\Omega)$ is a Hilbert space with the inner product

$$(u,w)_{W_p^l(\Omega)} = \left( \int_\Omega \sum_{|\alpha| \leq l} \partial^\alpha u(x) \overline{\partial^\alpha w(x)} d^n x \right). \tag{1.1.86}$$

Let $B = W_p^l(\Omega)$ topologized through by norm $\|u\|_{W_p^l(\Omega)}$ provides $G_{W_p^l(\Omega)}$. $G_{W_p^l(\Omega)}$ endowed with the sharp topology determined by $\wp_{W_p^l(\Omega)}(\bullet)$ is a $\widetilde{\mathbb{C}}$-module.

(II) Let $\omega(x)$ be a weight on $\mathbb{R}^n$ i.e.,a locally integrable function on $\mathbb{R}^n$ such that $\omega(x) > 0$ for a.e. $x \in \mathbb{R}^n$. Let $\Omega \subseteq \overline{\mathbb{R}}^n$ be open, $1 \leq p < 1$, and $l$ a nonnegative integer. The weighted Sobolev spaces $W_p^l(\Omega, \omega)$ consists of all functions $u$ with weak derivatives $\partial^\alpha u(x)$ for any $\alpha$ with $|\alpha| \leq l$ satisfying

$$\|u\|_{W_p^l(\Omega,\omega)} = \left( \int_\Omega \sum_{|\alpha| \leq l} |\partial^\alpha u(x)| \omega(x) d^n x \right)^{1/p} < \infty. \tag{1.1.86'}$$

Let $B = W_p^l(\Omega,\omega)$ topologized through by norm $\|u\|_{W_p^l(\Omega,\omega)}$ provides $G_{W_p^l(\Omega,\omega)}$. $G_{W_p^l(\Omega,\omega)}$ endowed with the sharp topology determined by $\wp_{W_p^l(\Omega,\omega)}(\cdot)$ is a $\widetilde{\mathbb{C}}$-module.

# I.1.11. Colombeau Sobolev spaces. Definition and basic properties of $\left( W_p^{l,k,\varepsilon}(\Omega) \right)_{\varepsilon \in (0,1]}$.

**Definition 1.1.36.** $(C^{l,k,\varepsilon}(\overline{\mathbb{R}}^n))_{\varepsilon \in (0,1]} = (C^{l,k,\varepsilon}(\overline{\mathbb{R}}^n))_\varepsilon$ is the linear space of all functions $u(x)$ in $\overline{\mathbb{R}}^n$ such that:
(1) $u(x)$ and $\partial^\alpha u(x)$ $|\alpha| \leq l$ are continuous in $\mathbb{R}^n$,
(2) $\sum_{|\alpha| \leq l} \sup_{x \in \overline{\mathbb{R}}^n} |\partial^\alpha u(x)| \leq \infty$ and the $\varepsilon$-norm $\|u\|_{C^{l,k,\varepsilon}}$ :

$$\|u\|_{C^{l,k,\varepsilon}} = \sum_{|\alpha| \leq l} \sup_{x \in \overline{\mathbb{R}}^n} \left| \partial^\alpha \left[ \frac{u(x)}{1 + \varepsilon^k \|x\|^{2k}} \right] \right| d^n x \tag{1.1.87}$$

is finite for any $\varepsilon \in (0,1]$. If $l = 0$, we denote $(C^{0,k,\varepsilon}(\overline{\mathbb{R}}^n))_\varepsilon = (C^{k,\varepsilon}(\overline{\mathbb{R}}^n))_\varepsilon$.

**Definition 1.1.37.** We can use (1.1.87) for defining the $\widetilde{\mathbb{R}}$-valued norm $(\|u\|_{C^{l,k,\varepsilon}})_\varepsilon$ :

$$\|u\|_{l,k,\widetilde{\mathbb{R}}} = \mathbf{cl}\left[ (\|u\|_{C^{l,k,\varepsilon}})_\varepsilon \right] \tag{1.1.88}$$

of a function $u \in (C^{l,k,\varepsilon}(\overline{\mathbb{R}}^n))_{\varepsilon \in (0,1]}$.

**Definition 1.1.38.** We can use (1.1.88) for defining the valuation $\mathbf{v}_{(C^{l,k,\varepsilon}(\overline{\mathbb{R}}^n))_\varepsilon}[u]$ :

$$\mathbf{v}_{(C^{l,k,\varepsilon}(\overline{\mathbb{R}}^n))_\varepsilon}[u] = \mathbf{v}_{\widetilde{\mathbb{R}}}\left[ \mathbf{cl}\left[ (\|u\|_{C^{l,k,\varepsilon}})_\varepsilon \right] \right] \tag{1.1.89}$$

of a function $u \in \left(C^{l,k,\varepsilon}(\overline{\mathbb{R}}^n)\right)_{\varepsilon \in (0,1]}$.

**Definition 1.1.39**.We can use (1.1.89) for defining the ultra pseudo-norm $\wp_{\left(C^{l,k,\varepsilon}(\overline{\mathbb{R}}^n)\right)_\varepsilon}[u]$ :

$$\wp_{\left(C^{l,k,\varepsilon}(\overline{\mathbb{R}}^n)\right)_\varepsilon}[u] = \exp\left[-\mathbf{v}_{\left(C^{l,k,\varepsilon}(\overline{\mathbb{R}}^n)\right)_\varepsilon}[u]\right] \qquad (1.1.90)$$

of a function $u \in \left(C^{l,k,\varepsilon}(\overline{\mathbb{R}}^n)\right)_{\varepsilon \in (0,1]}$.

**Definition 1.1.40.** Consider now an open subset $\Omega \subseteq \mathbb{R}^n$. $\left(W_p^{l,k,\varepsilon}(\Omega)\right)_{\varepsilon \in (0,1]} = \left(W_p^{l,k,\varepsilon}(\Omega)\right)_\varepsilon \left(\left(W_p^{l,k,\varepsilon}(\overline{\mathbb{R}}^n)\right)_{\varepsilon \in (0,1]} = \left(W_p^{l,k,\varepsilon}(\overline{\mathbb{R}}^n)\right)_\varepsilon\right)$ is the linear space of all functions $u(x)$ in $\Omega$ (in $\overline{\mathbb{R}}^n$) with weak derivatives $\partial^\alpha u(x)$ $|\alpha| \leq l$ satisfying $\sum_{|\alpha| \leq l} \sup_{x \in \overline{\mathbb{R}}^n} |\partial^\alpha u(x)| \leq \infty$ and the $\varepsilon$-norm $\|u\|_{W_p^{l,k,\varepsilon}(\Omega)}$ $\left(\|u\|_{W_p^{l,k,\varepsilon}(\overline{\mathbb{R}}^n)}\right)$ :

$$\|u\|_{W_p^{l,k,\varepsilon}(\Omega)} = \left(\int_\Omega \sum_{|\alpha| \leq l} \left|\partial^\alpha \left[\frac{u(x)}{1 + \varepsilon^k \|x\|^{2k}}\right]\right|^p d^n x\right)^{1/p},$$

$$\|u\|_{W_p^{l,k,\varepsilon}(\mathbb{R}^n)} = \left(\int_{\mathbb{R}^n} \sum_{|\alpha| \leq l} \left|\partial^\alpha \left[\frac{u(x)}{1 + \varepsilon^k \|x\|^{2k}}\right]\right|^p d^n x\right)^{1/p} \qquad (1.1.91)$$

is finite for any $\varepsilon \in (0,1]$.

**Definition 1.1.41**.We can use (1.1.91) for defining the $\widetilde{\mathbb{R}}$-valued norm $\left(\|u\|_{W_p^{l,k,\varepsilon}(\Omega)}\right)_\varepsilon$ :

$$\|u\|_{l,k,p,(\Omega)\widetilde{\mathbb{R}}} = \mathbf{cl}\left[\left(\|u\|_{W_p^{l,k,\varepsilon}(\Omega)}\right)_\varepsilon\right],$$

$$\|u\|_{l,k,p,(\Omega),\mathbb{R},\widetilde{\mathbb{R}}} = \mathbf{cl}\left[\left(\|u\|_{W_p^{l,k,\varepsilon}(\mathbb{R}^n)}\right)_\varepsilon\right] \qquad (1.1.92)$$

of a function $u \in \left(W_p^{l,k,\varepsilon}(\Omega)\right)_{\varepsilon \in (0,1]}$ $\left(u \in \left(W_p^{l,k,\varepsilon}(\overline{\mathbb{R}}^n)\right)_{\varepsilon \in (0,1]}\right)$.

**Definition 1.1.42**.We can use (1.1.92) for defining the valuation $\mathbf{v}_{\left(W_p^{l,k,\varepsilon}(\Omega)\right)_\varepsilon}[u]$

$$\left(\mathbf{v}_{\left(W_p^{l,k,\varepsilon}(\mathbb{R}^n)\right)_\varepsilon}[u]\right):$$

$$\mathbf{v}_{\left(W_p^{l,k,\varepsilon}(\Omega)\right)_\varepsilon}[u] = \mathbf{v}_{\widetilde{\mathbb{R}}}\left[\left[\left(\|u\|_{W_p^{l,k,\varepsilon}(\Omega)}\right)_\varepsilon\right]\right],$$

$$\mathbf{v}_{\left(C^{l,k,\varepsilon}(\mathbb{R}^n)\right)_\varepsilon}[u] = \mathbf{v}_{\widetilde{\mathbb{R}}}\left[\mathbf{cl}\left[\left(\|u\|_{W_p^{l,k,\varepsilon}(\mathbb{R}^n)}\right)_\varepsilon\right]\right]$$

(1.1.93)

of a function $u \in \left(W_p^{l,k,\varepsilon}(\Omega)\right)_{\varepsilon\in(0,1]}$ $\left(u \in \left(W_p^{l,k,\varepsilon}(\overline{\mathbb{R}}^n)\right)_{\varepsilon\in(0,1]}\right).$

**Definition 1.1.43.** We can use (1.1.93) for defining the pseudo-norm

$$\wp_{\left(W_p^{l,k,\varepsilon}(\Omega)\right)_\varepsilon}[u]\left(\wp_{\left(W_p^{l,k,\varepsilon}(\mathbb{R}^n)\right)_\varepsilon}[u]\right):$$

$$\wp_{\left(W_p^{l,k,\varepsilon}(\Omega)\right)_\varepsilon}[u] = \exp\left[-\mathbf{v}_{\left(W_p^{l,k,\varepsilon}(\Omega)\right)_\varepsilon}[u]\right],$$

$$\wp_{\left(W_p^{l,k,\varepsilon}(\mathbb{R}^n)\right)_\varepsilon}[u] = \exp\left[-\mathbf{v}_{\left(W_p^{l,k,\varepsilon}(\mathbb{R}^n)\right)_\varepsilon}[u]\right]$$

(1.1.94)

of a function $u \in \left(W_p^{l,k,\varepsilon}(\Omega)\right)_{\varepsilon\in(0,1]}$ $\left(u \in \left(W_p^{l,k,\varepsilon}(\overline{\mathbb{R}}^n)\right)_{\varepsilon\in(0,1]}\right).$

**Definition 1.1.44.** $\left(W_p^{l,\varepsilon}(\overline{\mathbb{R}}^n, \omega(x))\right)_{\varepsilon\in(0,1]}, \omega(x) = \exp(-\varepsilon x^2)$ is the linear space of all functions $u(x)$ in $\overline{\mathbb{R}}^n$ with weak derivatives $\partial^\alpha u(x)$ $|\alpha| \leq l$ $\sum_{|\alpha|\leq l}\sup_{\overline{\mathbb{R}}^n}|\partial^\alpha u(x)| \leq \infty$ and the $\varepsilon$-norm

$$\|u\|_{W_p^{l,\varepsilon}(\mathbb{R}^n,\omega)} = \left(\int_{\mathbb{R}^n}|\partial^\alpha u(x)|e^{-\varepsilon x^2}d^n x\right)^{1/p},$$

(1.1.95)

is finite for any $\varepsilon \in (0,1]$.

**Definition 1.1.**

## I.2. Pseudodifferential operators with coefficients are distributions.

Let us consider partiel differential equations of the form

$$\frac{\partial u_\varepsilon(t,x)}{\partial t} + A_\varepsilon(t,x,D)u_\varepsilon(t,x) = f(t,x),$$

$$u(0,x) = u_0(x),$$

$$x \in \Omega \subseteq \mathbb{R}^n.$$

(1.2.1)

We assume that

$$A(t,x,D) = \sum_{|\alpha| \leq m} a_\alpha(t,x)D^\alpha$$

(1.2.2)

with coefficients $a_\alpha(t,x)$ are distributions. Thus we may form the regularized operator

$$A_\varepsilon(t,x,D) = \sum_{|\alpha| \leq m} a_\alpha(t,x) * \varphi_\varepsilon(x)D^\alpha$$

(1.2.3)

where $\varphi_\varepsilon(x) \in C_c^\infty(\mathbb{R}^n)$ is a mollifier of the form

$$\varphi_\varepsilon(x) = \frac{1}{\varepsilon^n}\varphi(\frac{x}{\varepsilon}),$$

$$\int_{\mathbb{R}^n} \varphi(x)d^n x = 1.$$

(1.2.4)

Thus Eq.(1.2.1) becomes

$$\frac{\partial u_\varepsilon}{\partial t} + A_\varepsilon(t,x,D)u_\varepsilon = f_\varepsilon(t,x),$$

$$u_\varepsilon(0,x) = u_0(x) * \varphi_\varepsilon(x), \tag{1.2.5}$$

$$x \in \Omega \subseteq \mathbb{R}^n.$$

**Definition 1.2.1.** (1) We say that a family of smooth functions $(f_\varepsilon)_{\varepsilon \in (0,1]}$ satisfying an estimate of the type (**a**) iff

$$\forall \alpha(|\alpha| \in \mathbb{N}) \exists N(N \in \mathbb{N}) \forall \varepsilon(\varepsilon \in (0,1]) \left[ \sup_{x \in \mathbb{R}^n} |\partial^\alpha f(x)| = O(\varepsilon^{-N}) \right]. \tag{1.2.6}$$

(2) We say that a family of smooth functions $(f_\varepsilon)_{\varepsilon \in (0,1]}$ satisfying an estimate of the type (**b**) iff

$$\exists N(N \in \mathbb{N}) \forall \alpha(|\alpha| \in \mathbb{N}) \left[ \sup_{x \in \mathbb{R}^n} |\partial^\alpha f(x)| = O(\varepsilon^{-N}) \right]. \tag{1.2.7}$$

(3) We say that a family of smooth functions $(f_\varepsilon)_{\varepsilon \in (0,1]}$ satisfying an estimate of the type
(**c**) iff

$$\forall \alpha(|\alpha| \in \mathbb{N}) \forall q(q \in \mathbb{N}) \left[ \sup_{x \in \mathbb{R}^n} |\partial^\alpha f(x)| = O(\varepsilon^q) \right]. \tag{1.2.8}$$

**Remark 1.2.1.** The condition (2) signifies a regularity property of the family

$(f_\varepsilon)_{\varepsilon\in(0,1]}.$

It is well known that if the regularizations $f_\varepsilon = f * \varphi_\varepsilon(x)$ of a distribution $f \in E'(\mathbb{R}^n)$ satisfy (1.2.7) then $f$ actually is an infinitely differentiable function.

**Remark 1.2.2.** A family of smooth functions $(f_\varepsilon)_{\varepsilon\in(0,1]}$ satisfying an estimate of the type (**c**) will be considered as asymptotically negligible.

**Remark 1.2.3.** The regularized coefficients of the operator $A_\varepsilon(t,x,D)$ given by Eq.(1.2.3) will satisfy an estimate of type (**a**) or type (**b**),at least locally on compact sets, and the action of

$A_\varepsilon(t,x,D)$ on nets $(u_\varepsilon)_{\varepsilon\in(0,1]}$ preserves the asymptotic properties (1.2.6) and (1.2.7).
Thus if $(u_\varepsilon)_{\varepsilon\in(0,1]}$ enjoys either of these properties, so does $(A_\varepsilon(t,x,D)u_\varepsilon)_{\varepsilon\in(0,1]}.$

# I.3. The generalized oscillatory integral.

In this section we argue the meaning and the most important properties of oscillatory integrals of the type

$$\int_{K\times\mathbb{R}^p} \exp[i\varphi_\varepsilon(y,\xi)]a_\varepsilon(y,\xi)dyd\xi, \qquad (1.3.1)$$

where $K \Subset \Omega, \Omega$ an open subset of $\mathbb{R}^n$. The function $\varphi_\varepsilon(y,\xi)$ for any value of the parametr $\varepsilon \in (0,1]$ is assumed to be a phase function,i.e.,it is smooth on $\Omega \times \mathbb{R}^p \backslash \{0\}$
real valued, positively homogeneous of degree $1$ in $\xi$ and $\nabla\varphi_\varepsilon(y,\xi) \neq 0$ for all $y \in \Omega,$
$\xi \in \mathbb{R}^p \backslash \{0\}, \varepsilon \in (0,1].$ In the sequel of this paper we shall use the canonical square bracket notation

$$\mathbf{S}^m_{\rho,\delta}[\Omega \times \mathbb{R}^p] \qquad (1.3.3)$$

for the space of nets $\mathbf{S}^m_{\rho,\delta}(\Omega \times \mathbb{R}^p)^{(0,1]}$ where $\mathbf{S}^m_{\rho,\delta}(\Omega \times \mathbb{R}^p), m \in \mathbb{N}, \rho, \delta \in [0,1]$ is the usual space of Hörmander symbols.

**Definition1.3.1.** An element of $\mathbf{S}^m_{\rho,\delta,M}(\Omega \times \mathbb{R}^p)$ is a net $(u_\varepsilon)_\varepsilon \in \mathbf{S}^m_{\rho,\delta}[\Omega \times \mathbb{R}^p]$ such that

$$\forall \alpha(\alpha \in \mathbb{N}^p) \forall \beta(\beta \in \mathbb{N}^n) \forall K(K \Subset \Omega) \exists N(N \in \mathbb{N}) \exists \eta(\eta \in (0,1]) \exists c(c > 0)$$

$$\forall y(y \in K) \forall \xi(\xi \in \mathbb{R}^p) \forall \varepsilon(\varepsilon \in (0,1]) \Big[ \partial^\alpha_\xi \partial^\beta_y u_\varepsilon(y,\xi) \le c \langle \xi \rangle^{m-\rho|\alpha|+\delta|\beta|} \varepsilon^{-N} \Big], \tag{1.3.4}$$

$$\langle \xi \rangle = \sqrt{1 + |\xi|^2} \,.$$

**Remark 1.3.1.** The subscript $M$ underlines the moderate growth property, i. e., the bound of type $\varepsilon^{-N}$ as $\varepsilon \to 0$.

**Definition1.3.2.** An element of $\mathbf{N}^m_{\rho,\delta}(\Omega \times \mathbb{R}^p)$ is a net $(u_\varepsilon)_\varepsilon \in \mathbf{S}^m_{\rho,\delta}[\Omega \times \mathbb{R}^p]$ such that

$$\forall \alpha(\alpha \in \mathbb{N}^p) \forall \beta(\beta \in \mathbb{N}^n) \forall K(K \Subset \Omega) \exists q(q \in \mathbb{N}) \exists \eta(\eta \in (0,1]) \exists c(c > 0)$$

$$\forall y(y \in K) \forall \xi(\xi \in \mathbb{R}^p) \forall \varepsilon(\varepsilon \in (0,1]) \Big[ \partial^\alpha_\xi \partial^\beta_y u_\varepsilon(y,\xi) \le c \langle \xi \rangle^{m-\rho|\alpha|+\delta|\beta|} \varepsilon^q \Big]. \tag{1.3.5}$$

**Remark 1.3.2.** Nets of this type are called negligible.

**Definition1.3.3.** A generalized symbol of order $m$ and type $(\rho, \delta)$ is an element of the factor space

$$\tilde{\mathbf{S}}^m_{\rho,\delta}(\Omega \times \mathbb{R}^p) \triangleq \mathbf{S}^m_{\rho,\delta}(\Omega \times \mathbb{R}^p)/\mathbf{N}^m_{\rho,\delta}(\Omega \times \mathbb{R}^p) \tag{1.3.6}$$

**Remark 1.3.3.** We denote an arbitrary representative of $v \in \tilde{\mathbf{S}}^m_{\rho,\delta}[\Omega \times \mathbb{R}^p]$ by $(v_\varepsilon)_{\varepsilon \in (0,1]}$ or simply by $(v_\varepsilon)_\varepsilon$.

**Definition1.3.4.** An element $v \in \tilde{\mathbf{S}}^m_{\rho,\delta}(\Omega \times \mathbb{R}^p)$ is called regular if it has a representative $(v_\varepsilon)_\varepsilon$ with the following property:

$$\forall K(K \Subset \Omega) \exists N(N \in \mathbb{N}) \forall \alpha(\alpha \in \mathbb{N}^p) \forall \beta(\beta \in \mathbb{N}^n) \exists \eta(\eta \in (0,1]) \exists c(c > 0)$$

$$\text{(1.3.7)}$$

$$\forall y(y \in K) \forall \xi(\xi \in \mathbb{R}^p) \forall \varepsilon(\varepsilon \in (0,1]) \left[ \partial_\xi^\alpha \partial_y^\beta u_\varepsilon(y,\xi) \le c \langle \xi \rangle^{m-\rho|\alpha|+\delta|\beta|} \varepsilon^{-N} \right].$$

We denote by

$$\tilde{\mathbf{S}}_{\rho,\delta,\mathbf{rg}}^m(\Omega \times \mathbb{R}^p) \tag{1.3.8}$$

the subspace of regular elements of the space $\tilde{\mathbf{S}}_{\rho,\delta}^m(\Omega \times \mathbb{R}^p)$.

If the property (1.3.7) holds for one representative of $v$, it holds for every representa- tive.Consequently, if $\tilde{\mathbf{S}}_{\rho,\delta,\mathbf{rg}}^m(\Omega \times \mathbb{R}^p)$ is the space defined by (1.3.7) we can introduce $\tilde{\mathbf{S}}_{\rho,\delta,\mathbf{rg}}^m(\Omega \times \mathbb{R}^p)$ as the factor space:

$$\tilde{\mathbf{S}}_{\rho,\delta,\mathbf{rg}}^m(\Omega \times \mathbb{R}^p) \triangleq \mathbf{S}_{\rho,\delta,\mathbf{rg}}^m(\Omega \times \mathbb{R}^p)/\mathbf{N}_{\rho,\delta}^m(\Omega \times \mathbb{R}^p) \tag{1.3.9}$$

**Remark 1.3.4.** Note that (1) $(v_\varepsilon)_\varepsilon \in \mathbf{S}_{\rho,\delta,\mathbf{M}}^m(\Omega \times \mathbb{R}^p)$ implies $\left( \partial_\xi^\alpha \partial_x^\beta v_\varepsilon \right)_\varepsilon \in \mathbf{S}_{\rho,\delta,\mathbf{M}}^{m-\rho|\alpha|-\delta|\beta|}(\Omega \times \mathbb{R}^p)$,(2) $(v_\varepsilon)_\varepsilon \in \mathbf{S}_{\rho,\delta,\mathbf{M}}^{m_1}(\Omega \times \mathbb{R}^p)$ and $(w_\varepsilon)_\varepsilon \in \mathbf{S}_{\rho,\delta,\mathbf{M}}^{m_2}(\Omega \times \mathbb{R}^p)$ implies $(v_\varepsilon + w_\varepsilon)_\varepsilon \in \mathbf{S}_{\rho,\delta,\mathbf{M}}^{\max\{m_1,m_2\}}(\Omega \times \mathbb{R}^p)$ and $(v_\varepsilon w_\varepsilon)_\varepsilon \in \mathbf{S}_{\rho,\delta,\mathbf{M}}^{m_1+m_2}(\Omega \times \mathbb{R}^p)$. Since the results concerning derivatives and sums hold with $\mathbf{S}_{\rho,\delta,\mathbf{rg}}^m(\Omega \times \mathbb{R}^p)$ and $N_{\rho,\delta}$ in place of $\mathbf{S}_{\rho,\delta,\mathbf{M}}^m(\Omega \times \mathbb{R}^p)$, we can define derivatives and sums on the corresponding factor spaces.Moreover,$(v_\varepsilon)_\varepsilon \in \mathbf{S}_{\rho,\delta,\mathbf{M}}^{m_1}(\Omega \times \mathbb{R}^p)$ and $(w_\varepsilon)_\varepsilon \in \mathbf{N}_{\rho,\delta}^{m_2}(\Omega \times \mathbb{R}^p)$ imply $(v_\varepsilon w_\varepsilon)_\varepsilon \in \mathbf{S}_{\rho,\delta}^{m_1+m_2}(\Omega \times \mathbb{R}^p)$, thus we obtain that the product is a well-defined map from the space $\tilde{\mathbf{S}}_{\rho,\delta}^{m_1}(\Omega \times \mathbb{R}^p) \times \tilde{\mathbf{S}}_{\rho,\delta}^{m_2}(\Omega \times \mathbb{R}^p)$ into $\tilde{\mathbf{S}}_{\rho,\delta}^{m_1+m_2}(\Omega \times \mathbb{R}^p)$.Similarly, it is well-defined as a map from $\tilde{\mathbf{S}}_{\rho,\delta,\mathbf{rg}}^{m_1}(\Omega \times \mathbb{R}^p) \times \tilde{\mathbf{S}}_{\rho,\delta,\mathbf{rg}}^{m_2}(\Omega \times \mathbb{R}^p)$ into $\tilde{\mathbf{S}}_{\rho,\delta,\mathbf{rg}}^{m_1+m_2}(\Omega \times \mathbb{R}^p)$.

**Remark 1.3.5.** The classical space $\mathbf{S}_{\rho,\delta}^{m_1}(\Omega \times \mathbb{R}^p)$ is contained in $\tilde{\mathbf{S}}_{\rho,\delta,\mathbf{rg}}^m(\Omega \times \mathbb{R}^p)$. We assume from now on that $\rho > 0$ and $\delta < 1$ and return to the meaning of the oscillatory integral

$$\int_{K \times \mathbb{R}^p} \exp[i\varphi_\varepsilon(y,\xi)]a_\varepsilon(y,\xi)dyd\xi. \tag{1.3.10}$$

Obviously, if $(\varphi_\varepsilon)_\varepsilon, (a_\varepsilon)_\varepsilon \in \mathbf{S}_{\rho,\delta,\mathbf{M}}^{m_1}(\Omega \times \mathbb{R}^p)$ then (1.3.10) makes sense as an oscillatory integral for every $\varepsilon \in (0,1]$. We recall now that given the phase function $\varphi_\varepsilon, \varepsilon \in (0,1]$, there exists an operator

$$L_\varepsilon = \sum_{i=1}^{p} a_{i,\varepsilon}(y,\xi)\frac{\partial}{\partial \xi_i} + \sum_{k=1}^{n} b_{k,\varepsilon}(y,\xi)\frac{\partial}{\partial y_k} + c_\varepsilon(y,\xi) \tag{1.3.11}$$

such that $a_{i,\varepsilon}(y,\xi) \in \mathbf{S}^0(\Omega \times \mathbb{R}^p), b_{k,\varepsilon}(y,\xi) \in \mathbf{S}^{-1}(\Omega \times \mathbb{R}^p), c_\varepsilon(y,\xi) \in \mathbf{S}^{-1}(\Omega \times \mathbb{R}^p)$, and such that ${}^tL\exp[i\varphi_\varepsilon] = \exp[i\varphi_\varepsilon]$, where ${}^tL$ is the formal adjoint. We note that for $m - js < -p$ and $\chi \in C_c^\infty(\mathbb{R}^p)$ identically equal to $1$ in a neighborhood of the origin, the oscillatory integral $I_K(\varphi_\varepsilon, a_\varepsilon)$ at fixed $\varepsilon$, can be defined by either of the two expressions on the right hand-side of Eq.(1.3.12)

$$I_K(\varphi_\varepsilon, a_\varepsilon) = \int_{K \times \mathbb{R}^p} \exp[i\varphi_\varepsilon(y,\xi)]a_\varepsilon(y,\xi)dyd\xi =$$

$$\int_{K \times \mathbb{R}^p} \exp[i\varphi_\varepsilon(y,\xi)]L^j a_\varepsilon(y,\xi)dyd\xi = \tag{1.3.12}$$

$$\lim_{h \to 0+} \int_{K \times \mathbb{R}^p} \exp[i\varphi_\varepsilon(y,\xi)]a_\varepsilon(y,\xi)\chi(h\xi)dyd\xi,$$

where the equalities hold for all $\varepsilon \in (0,1]$.

**Theorem 1.3.1**. Let $s = \min\{\rho, 1 - \delta\}$ and $j \in \mathbb{N}$. Then the following statements hold:
(1) if $(a_\varepsilon)_\varepsilon \in \mathbf{S}_{\rho,\delta,\mathbf{M}}^m(\Omega \times \mathbb{R}^p)$ then $(L^j a_\varepsilon)_\varepsilon \in \mathbf{S}_{\rho,\delta,\mathbf{M}}^{m-js}(\Omega \times \mathbb{R}^p)$,
(2) statement (1) is valid with $\mathbf{S}_{\rho,\delta,\mathbf{rg}}^m(\Omega \times \mathbb{R}^p)$ in place of $\mathbf{S}_{\rho,\delta,\mathbf{M}}^m(\Omega \times \mathbb{R}^p)$.

**Theorem 1.3.2**. Let $K$ be a compact set contained in $\Omega$ and let $\varphi_\varepsilon, \varepsilon \in (0,1]$ be a phase function on $\Omega \times \mathbb{R}^p$ and $\varphi = \mathbf{cl}[(\varphi_\varepsilon)_\varepsilon]$, $a = \mathbf{cl}[(a_\varepsilon)_\varepsilon]$ an element of $\tilde{\mathbf{S}}_{\rho,\delta}^m(\Omega \times \mathbb{R}^p)$. Then the oscillatory integral

$$I_K(\varphi, a) = \int_{K \times \mathbb{R}^p} \exp[i\varphi_\varepsilon(y, \xi)] a_\varepsilon(y, \xi) dy d\xi =$$

$$(1.3.13)$$

$$(I_K(\varphi_\varepsilon, a_\varepsilon))_\varepsilon + \mathbf{N}$$

is a well-defined element of $\widetilde{\mathbb{C}}$.

**Definition 1.3.5.** Let $\varphi, a \in \widetilde{\mathbf{S}}_{\rho, \delta}^m(\Omega \times \mathbb{R}^p)$ with $\mathrm{supp}_y(a)$. We define the generalized oscillatory integral

$$I(\varphi, a) \triangleq \int_{\Omega \times \mathbb{R}^p} \exp[i\varphi(y, \xi)] a(y, \xi) dy d\xi \triangleq$$

$$(1.3.14)$$

$$\int_{K \times \mathbb{R}^p} \exp[i\varphi(y, \xi)] a(y, \xi) dy d\xi,$$

where $K$ is any compact subset of $\Omega$ containing $\mathrm{supp}_y(a)$ in its interior.

**Remark 1.3.6.** We note that definition 1.3.5 does not depend on the choice of $K$.

**Definition 1.3.6.** Let us consider now phase functions and symbols depending on an additional parameter.

$$I(\varphi, a)(x) = \int_{K \times \mathbb{R}^p} \exp[i\varphi(x, y, \xi)] a(x, y, \xi) dy d\xi, \qquad (1.3.15)$$

where $x \in \Omega'$, an open subset of $\mathbb{R}^{n'}$. If for any fixed $\tilde{x} \in \Omega', \varphi(\tilde{x}, y, \xi) \in \widetilde{\mathbf{S}}_{\rho, \delta}^{m_1}(\Omega \times \mathbb{R}^p)$ is a phase function with respect to the variables $(y, \xi)$ and $a(\tilde{x}, y, \xi) \in \widetilde{\mathbf{S}}_{\rho, \delta}^{m_2}(\Omega \times \mathbb{R}^p)$ the oscillatory integral $I(\varphi, a)(x)$ defines a map from $\Omega'$ to $\widetilde{\mathbb{C}}$.

**Remark 1.3.7.** Let for any $\varepsilon \in (0, 1]$, $\varphi_\varepsilon(x, y, \xi)$ be a real valued continuous function on $\Omega' \times \Omega \times \mathbb{R}^p$, smooth on $\Omega' \times \Omega \times \mathbb{R}^p \backslash \{0\}$ such that for all $x \in \Omega'$, such that for any $\varepsilon \in (0, 1], \varphi(\tilde{x}, y, \xi)$ is a phase function with respect to the variables $(y, \xi)$. Obviously we have that for all $j \in \mathbb{N}, (\varphi_\varepsilon(x, y, \xi))_\varepsilon \in \mathbf{S}_{\rho, \delta}^{m_1}(\Omega' \times \Omega \times \mathbb{R}^p)$ and $(a_\varepsilon(x, y, \xi))_\varepsilon \in \mathbf{S}_{\rho, \delta}^{m_2}(\Omega' \times \Omega \times \mathbb{R}^p)$ implies $(L_x^j a_\varepsilon(x, y, \xi))_\varepsilon \in \mathbf{S}_{\rho, \delta, \mathbf{M}}^{m_2 - j s}(\Omega' \times \Omega \times \mathbb{R}^p)$. The same result holds with $\mathbf{S}_{\rho, \delta, \mathbf{rg}}^m(\Omega' \times \Omega \times \mathbb{R}^p)$ in place of $\mathbf{S}_{\rho, \delta, \mathbf{M}}^m(\Omega' \times \Omega \times \mathbb{R}^p)$.

**Theorem 1.3.3.** Let $\varphi(x, y, \xi) = \mathbf{cl}[(\varphi_\varepsilon(x, y, \xi))_\varepsilon]$ be as in Remark 1.3.7.

(1) If $a(x,y,\xi) \in \tilde{\mathbf{S}}_{\rho,\delta}^{m_2}(\Omega' \times \Omega \times \mathbb{R}^p)$ then for all $K \Subset \Omega$

$$W_K(x) = \int_{K \times \mathbb{R}^p} \exp[i\varphi(x,y,\xi)]a(x,y,\xi)dyd\xi \qquad (1.3.16)$$

then $W_K(x) \in G(\Omega')$.
(2) If $\varphi(x,y,\xi) \in \tilde{\mathbf{S}}_{\rho,\delta,\mathbf{rg}}^{m_1}(\Omega' \times \Omega \times \mathbb{R}^p)$ and $a(x,y,\xi) \in \tilde{\mathbf{S}}_{\rho,\delta,\mathbf{rg}}^{m_2}(\Omega' \times \Omega \times \mathbb{R}^p)$ then $W_K(x) \in G^\infty(\Omega')$ for all $K \Subset \Omega$.

**Proof**. An arbitrary representative of $W_K(x)$ is given by the oscillatory integral

$$(W_{K,\varepsilon}(x))_\varepsilon = \left( \int_{K \times \mathbb{R}^p} \exp[i\varphi_\varepsilon(x,y,\xi)]a_\varepsilon(x,y,\xi)dyd\xi \right)_\varepsilon . \qquad (1.3.17)$$

From Remark 1.3.7. it follows that $(W_{K,\varepsilon}(x))_\varepsilon \in E(\Omega')$. At this point computing the $x$-derivatives of the expression $\exp[i\varphi_\varepsilon(x,y,\xi)]L_x^j a_\varepsilon(x,y,\xi)$ for $\xi \neq 0$, we conclude that

$$\forall\alpha\left(\alpha \in \mathbb{N}^{n'}\right)\forall K'\left(K' \Subset \Omega'\right)\exists N(N \in \mathbb{N})\exists\eta(\eta \in (0,1])\forall x(x \in K')$$

$$\forall y(y \in K)\forall\xi(\xi \in \mathbb{R}^p\backslash\{0\})\exists\varepsilon(\varepsilon \in (0,1]) \qquad (1.3.18)$$

$$\left[\left|\partial_x^\alpha\left(\exp[i\varphi_\varepsilon(x,y,\xi)]L_x^j a_\varepsilon(x,y,\xi)\right)\right| \leq \langle\xi\rangle^{m_1+m_2-js+|\alpha|}\varepsilon^{-N}\right].$$

Now if $m_1 + m_2 - js + |\alpha| < -p$ then we obtain for any $x \in K'$ and $\varepsilon \in (0,1]$ the inequality

$$|\partial_x^\alpha W_{K,\varepsilon}(x)| \leq \varepsilon^{-N}. \qquad (1.3.19)$$

Hence $(W_{K,\varepsilon}(x))_\varepsilon \in E_{\mathbf{M}}(\Omega')$. Obviously if $(a_\varepsilon)_\varepsilon \in \mathbf{N}_{\rho,\delta}^m(\Omega' \times \Omega \times \mathbb{R}^p)$ then $(W_{K,\varepsilon}(x))_\varepsilon \in \mathbf{N}(\Omega')$. If $\varphi \in \tilde{\mathbf{S}}^{m_1}_{\rho,\delta,\mathbf{rg}}(\Omega' \times \Omega \times \mathbb{R}^p)$ and $a \in \tilde{\mathbf{S}}^{m_2}_{\rho,\delta,\mathbf{rg}}(\Omega' \times \Omega \times \mathbb{R}^p)$ the exponent $N$ in (1.3.19) does not depend on the derivatives and then $W_K(x) \in G^\infty(\Omega')$. This result completes the proof.

**Remark 1.3.8.** If $\varphi = \mathbf{cl}[(\varphi_\varepsilon)_\varepsilon]$ is a generalized phase function with respect to $(y,\xi)$ and there exists a compact subset $K$ of $\Omega$ such that $\forall x (x \in \Omega') \big[ \mathrm{supp}_y(x,y,\xi) \subseteq K \big]$ then the oscillatory integral

$$(W_{K,\varepsilon}(x))_\varepsilon = \left( \int_{\Omega \times \mathbb{R}^p} \exp[i\varphi_\varepsilon(x,y,\xi)] a_\varepsilon(x,y,\xi) dy d\xi \right)_\varepsilon \tag{1.3.20}$$

defines a generalized function belonging to $G(\Omega')$.

**Remark 1.3.9.** We recall that for each generalized phase function $(\varphi_\varepsilon)_\varepsilon$

$$C_{\varphi_\varepsilon} = \{(x,\xi) \in \Omega \times \mathbb{R}^p \backslash\{0\} | \nabla_\xi \varphi_\varepsilon(x,\xi) = 0, \varepsilon \neq 0\} \tag{1.3.21}$$

is a cone-shaped subset of $\Omega \times \mathbb{R}^p \backslash\{0\}$. Let $\pi : \Omega \times \mathbb{R}^p \backslash\{0\} \to \Omega$ be the projection onto $\Omega$ and put $S_{\varphi_\varepsilon} = \pi(C_{\varphi_\varepsilon}), R_{\varphi_\varepsilon} = \Omega \backslash S_{\varphi_\varepsilon}$. Interpreting $x \in \Omega$ as a parameter we have from Theorem 1.3.3 that

$$W(x) =$$

$$\int_{\mathbb{R}^p} \exp[i\varphi(x,\xi)] a(x,\xi) d\xi =$$

$$\left( \int_{\mathbb{R}^p} \exp[i\varphi_\varepsilon(x,\xi)] a_\varepsilon(x,\xi) d\xi \right)_\varepsilon + \mathbf{N} = \tag{1.3.22}$$

$$(W_\varepsilon(x))_\varepsilon + \mathbf{N}.$$

makes sense as an generalized oscillatory integral for $x \in R_{\varphi_\varepsilon}$.

(1) If $\varphi \in \tilde{\mathbf{S}}^{m_1}_{\rho,\delta}(\Omega' \times \Omega \times \mathbb{R}^p)$ and $a \in \tilde{\mathbf{S}}^{m_2}_{\rho,\delta}(\Omega' \times \Omega \times \mathbb{R}^p)$ then $W(x) \in G(R_\varphi)$, where $R_\varphi = \bigcap_{\varepsilon \in (0,1]} R_{\varphi_\varepsilon}$.

(2) If $\varphi \in \tilde{\mathbf{S}}^{m_1}_{\rho,\delta}(\Omega' \times \Omega \times \mathbb{R}^p)$ and $a \in \tilde{\mathbf{S}}^{m_2}_{\rho,\delta}(\Omega' \times \Omega \times \mathbb{R}^p)$ and $\mathrm{supp}_x a_\varepsilon \subseteq \Omega \backslash \Omega'_\varepsilon$ where $\Omega'_\varepsilon$ is an open neighborhood of $R_{\varphi_\varepsilon}$ then $W(x)$ can be extended to a generalized function on $\Omega$ with support contained in $\Omega \backslash \cup_{\varepsilon \in (0,1]} \Omega'_\varepsilon$.

We just give some details concerning the proof of the assertion (2). Let

$\{\Omega'_{j,\varepsilon}\}_{j\in\mathbb{N}\backslash\{0\}}$ be an open covering of $\Omega'_\varepsilon$ such that $\Omega'_{j,\varepsilon}$ is relatively compact and $\Omega'_{j,\varepsilon} \subseteq \overline{\Omega}'_{j,\varepsilon} \subseteq \Omega'_{j+1,\varepsilon}$ for all $j$. Choosing now cut-off functions $\{\psi_{j,\varepsilon}\}_{j\in\mathbb{N}\backslash\{0\}}$ such that $\psi_{j,\varepsilon} \in C^\infty_c(\Omega'_\varepsilon)$ and $\psi_{j,\varepsilon} \equiv 1$ in a neighborhood of $\overline{\Omega}'_{j,\varepsilon}$, we observe that $((1-\psi_{j,\varepsilon}(x))a_\varepsilon(x,\xi))_\varepsilon$ is a representative of a identically equal to $0$ on $\Omega'_{j,\varepsilon}$ for all $\varepsilon \in (0,1]$. Hence we see that

$$W_0(x) = \left(\int_{\mathbb{R}^p}\exp[i\varphi_\varepsilon(x,\xi)]a_\varepsilon(x,\xi)d\xi\right)_\varepsilon + \mathbf{N}(R_\varphi),$$

$$W_j(x) = \left(\int_{\mathbb{R}^p}(1-\psi_{j,\varepsilon}(x))\exp[i\varphi_\varepsilon(x,\xi)]a_\varepsilon(x,\xi)d\xi\Big|_{\Omega'_{j,\varepsilon}}\right)_\varepsilon + \mathbf{N}(\Omega'_{j,\varepsilon}),$$

$$j \geq 1$$

(1.3.23)

is a coherent family of generalized functions which defines $W(x) \in G(\Omega)$ such that $W_0(x)$ and $W_j(x)$ are its restrictions to $R_\varphi$ and $\cup_{\varepsilon\in(0,1]}\Omega'_{j,\varepsilon}$ respectively.

# I.4. Pseudodifferential operators with generalized symbol.

**Definition 1.4.1.** Let $m,\mu,\delta$, be real numbers, $\rho,\delta \in [0,1]$. Let $\omega(\tau)$ be a real valued function $\omega : (0,1] \to \mathbb{R}, \omega(\tau) > 0$ on the interval $(0,1]$, such that for some $r$ in $\mathbb{R}$, for some $C > 0$ and for all $\varepsilon \in (0,1], \omega(\varepsilon) > C\varepsilon^r$. We denote by $\mathbf{S}^{m,\mu}_{\rho,\delta,\omega}(\Omega' \times \mathbb{R}^p)$, where $\Omega'$ an open subset of $\mathbb{R}^{n'}$, the set of all $(a_\varepsilon)_\varepsilon \in \mathbf{S}^m_{\rho,\delta}(\Omega' \times \mathbb{R}^p)$ such that the following statement holds:

$$\forall K(K \Subset \Omega')\exists N(N \in \mathbb{N})\forall\alpha(\alpha \in \mathbb{N}^p)\forall\beta\left(\beta \in \mathbb{N}^{n'}\right)\exists\eta(\eta \in (0.1])$$

$$\exists c(c > 0)\forall x(x \in K)\forall\xi(\xi \in \mathbb{R}^p)\forall\varepsilon(\varepsilon \in (0,1])$$

(1.4.1)

$$\left|\partial^\alpha_\xi\partial^\beta_x a_\varepsilon(x,\xi)\right| \leq c\langle\xi\rangle^{m-\rho|\alpha|+\delta|\beta|}\varepsilon^{-N}(\omega(\varepsilon))^{-(|\beta|-\mu)_+},$$

where the exponent $-(|\beta|-\mu)_{+} = -\max\{0,|\beta|-\mu\}$ reflects differentiability up to order $\mu$ in the case when $(a_\varepsilon)_\varepsilon$ is obtained from a non-smooth, classical symbol by convolution with a mollifier with scale $\omega(\varepsilon)$.

**Definition 1.4.2.** An element of $\mathbf{N}^{m,\mu}_{\rho,\delta,\omega}(\Omega' \times \mathbb{R}^p)$ is a net in $\mathbf{S}^m_{\rho,\delta}[\Omega' \times \mathbb{R}^p]$ such that the following statement holds:

$$\forall K(K \Subset \Omega')\forall\alpha(\alpha \in \mathbb{N}^p)\forall\beta\left(\beta \in \mathbb{N}^{n'}\right)\forall q(q \in \mathbb{N})\exists\eta(\eta \in (0,1])$$

$$\exists c(c > 0)\forall x(x \in K)\forall\xi(\xi \in \mathbb{R}^p)\forall\varepsilon(\varepsilon \in (0.1]) \tag{1.4.2}$$

$$\left|\partial_\xi^\alpha \partial_x^\beta a_\varepsilon(x,\xi)\right| \leq c\langle\xi\rangle^{m-\rho|\alpha|+\delta|\beta|}\varepsilon^q(\omega(\varepsilon))^{-(|\beta|-\mu)_+}.$$

Nets with this property are called negligible.

**Definition 1.4.3.** The factor space $\mathbf{S}^{m,\mu}_{\rho,\delta,\omega}(\Omega' \times \mathbb{R}^p)/\mathbf{N}^{m,\mu}_{\rho,\delta,\omega}(\Omega' \times \mathbb{R}^p)$ will be denoted by

$$\tilde{\mathbf{S}}^{m,\mu}_{\rho,\delta,\omega}(\Omega' \times \mathbb{R}^p). \tag{1.4.3}$$

**Remark 1.4.1.** The following mapping properties hold

$$\partial_\xi^\alpha \partial_x^\beta : \tilde{\mathbf{S}}^{m,\mu}_{\rho,\delta,\omega}(\Omega' \times \mathbb{R}^p) \to \tilde{\mathbf{S}}^{m-\rho|\alpha|+\delta|\beta|,\mu-|\beta|}_{\rho,\delta,\omega}(\Omega' \times \mathbb{R}^p),$$

$$(\cdot + \cdot) : \tilde{\mathbf{S}}^{m_1,\mu}_{\rho,\delta,\omega}(\Omega' \times \mathbb{R}^p) \times \tilde{\mathbf{S}}^{m_2,\mu}_{\rho,\delta,\omega}(\Omega' \times \mathbb{R}^p) \to \tilde{\mathbf{S}}^{\max(m_1,m_2),\mu}_{\rho,\delta,\omega}(\Omega' \times \mathbb{R}^p). \tag{1.4.4}$$

**Definition 1.4.4.** We denote by $\mathbf{S}^{-\infty}_{\mathbf{rg}}(\Omega' \times \mathbb{R}^p)$ the set of all $(a_\varepsilon)_\varepsilon \in \mathbf{S}^{-\infty}[\Omega' \times \mathbb{R}^p]$ such that

$$\forall K(K \Subset \Omega')\exists N(N \in \mathbb{N})\forall m(m \in \mathbb{R})\forall \alpha(\alpha \in \mathbb{N}^p)\forall \beta\big(\beta \in \mathbb{N}^{n'}\big)\exists \eta(\eta \in (0.1])$$

$$\exists c(c > 0)\forall x(x \in K)\forall \xi(\xi \in \mathbb{R}^p)\forall \varepsilon(\varepsilon \in (0,1]) \qquad (1.4.5)$$

$$\left|\partial_\xi^\alpha \partial_x^\beta a_\varepsilon(x,\xi)\right| \le c\langle \xi\rangle^{m-|\alpha|}\varepsilon^{-N}.$$

We denote by $\mathbf{N}^{-\infty}(\Omega' \times \mathbb{R}^p)$ the set of all $(a_\varepsilon)_\varepsilon \in \mathbf{S}^{-\infty}[\Omega' \times \mathbb{R}^p]$ such that

$$\forall K(K \Subset \Omega')\exists N(N \in \mathbb{N})\forall q(q \in \mathbb{N})\forall \alpha(\alpha \in \mathbb{N}^p)\forall \beta\big(\beta \in \mathbb{N}^{n'}\big)\exists \eta(\eta \in (0.1])$$

$$\exists c(c > 0)\forall x(x \in K)\forall \xi(\xi \in \mathbb{R}^p)\forall \varepsilon(\varepsilon \in (0,1]) \qquad (1.4.6)$$

$$\left|\partial_\xi^\alpha \partial_x^\beta a_\varepsilon(x,\xi)\right| \le c\langle \xi\rangle^{m-|\alpha|}\varepsilon^{-q}.$$

The factor space $\mathbf{S}_{\mathbf{rg}}^{-\infty}(\Omega' \times \mathbb{R}^p)/\mathbf{N}^{-\infty}(\Omega' \times \mathbb{R}^p)$ we denoted by $\tilde{\mathbf{S}}_{\mathbf{rg}}^{-\infty}(\Omega' \times \mathbb{R}^p)$.

**Remark 1**.**4**.**2**. If in addition to the usual assumptions on $\omega(\varepsilon)$ we assume that $(\omega^{-1}(\varepsilon))_\varepsilon$ is a slow scale net, we can uniformly bound the contributions of the derivatives in (1.4.1) by a single negative power of $\varepsilon$,obtaining for a suitable constant $c$, for $x \in K$ and for all $\varepsilon$ small enough,hence

$$\left|\partial_\xi^\alpha \partial_x^\beta a_\varepsilon(x,\xi)\right| \le c\langle \xi\rangle^{m-\rho|\alpha|+\delta|\beta|}\varepsilon^{-N-1}.$$

**Definition 1**.**4**.**5**.Let $\Omega$ be an open subset of $\mathbb{R}^n$. The elements of the factor spaces $\tilde{\mathbf{S}}_{\rho,\delta,\omega}^{m,\mu}(\Omega \times \mathbb{R}^p)$ and $\tilde{\mathbf{S}}_{\rho,\delta,\omega}^{m,\mu}(\Omega \times \Omega \times \mathbb{R}^p)$ will be called symbols and amplitudes of order $m$ and type $(\rho,\delta,\mu,\omega)$ respectively. The elements of the factor spaces $\tilde{\mathbf{S}}_{\mathbf{rg}}^{-\infty}(\Omega \times \mathbb{R}^p)$ and $\tilde{\mathbf{S}}_{\mathbf{rg}}^{-\infty}(\Omega \times \Omega \times \mathbb{R}^p)$ will be called smoothing symbols and smoothing amplitudes, respectively.

**Theorem 1**.**4**.**1**. Let $a \in \tilde{\mathbf{S}}_{\rho,\delta,\omega}^{m,\mu}(\Omega \times \Omega \times \mathbb{R}^n)$. The oscillatory integral

$$Au = \int_{\Omega \times \mathbb{R}^n} \exp[i(x-y)\xi]a(x,y,\xi)u(y)dyd\xi =$$

$$(A_\varepsilon u_\varepsilon(x))_\varepsilon + \mathbf{N}(\Omega), \qquad (1.4.7)$$

$$A_\varepsilon u_\varepsilon(x) = \int_{\Omega \times \mathbb{R}^n} \exp[i(x-y)\xi]a_\varepsilon(x,y,\xi)u_\varepsilon(y)dyd\xi$$

defines:

(**1**) a $\widetilde{\mathbb{C}}$-linear continuous map from $G_c(\Omega)$ into $G(\Omega)$,

(**2**) a $\widetilde{\mathbb{C}}$-linear continuous map from $G_c^\infty(\Omega)$ into $G^\infty(\Omega)$ if $(\omega^{-1}(\varepsilon))_\varepsilon$ is a slow scale net,

(**3**) a $\widetilde{\mathbb{C}}$-linear continuous map from $G_c(\Omega)$ into $G^\infty(\Omega)$ if $a$ belongs to

$\mathbf{S}_{\mathbf{rg}}^{-\infty}(\Omega \times \Omega \times \mathbb{R}^n)$.

**Definition 1.4.6.** (1) Let $a \in \widetilde{\mathbf{S}}_{\rho,\delta,\omega}^{m,\mu}(\Omega \times \Omega \times \mathbb{R}^n)$. The linear map $A : G_c(\Omega) \to G(\Omega)$ defined via formula

$$u \to (Au)(x) = \int_{\Omega \times \mathbb{R}^n} \exp[i(x-y)\xi]a(x,y,\xi)u(y)dyd\xi \qquad (1.4.8)$$

will be called a (generalized) pseudodifferential operator with amplitude $a$.
(2) The formal transpose of the operator $A$ is the pseudodifferential operator $^tA : G_c(\Omega) \to G(\Omega)$ defined via formula

$$u \to ({}^tAu)(y) = \int_{\Omega \times \mathbb{R}^n} \exp[i(x-y)\xi]a(x,y,\xi)u(y)dxd\xi \qquad (1.4.9)$$

**Definition 1.4.7.** Let $A$ be a pseudodifferential operator. The map $k_A \in L\left(G_c(\Omega \times \Omega), \widetilde{\mathbb{C}}\right)$ defined by

$$k_A(u) = \int (Au(x, \cdot))(x)dx \qquad (1.4.10)$$

is called the kernel of $A$.

**Theorem 1.4.2**. Let $a \in \tilde{\mathbf{S}}_{\rho,\delta,\omega}^{m,\mu}(\Omega \times \Omega \times \mathbb{R}^n)$ and $A$ be the corresponding pseudo-differential operator.

(1) For all $u_1 \in G_{\mathbf{c}}(\Omega)$ and $u_2 \in G_{\mathbf{c}}(\Omega)$

$$k_A(u_2 \otimes u_1) = \int_\Omega (Au_1)(x)u_2(x)dx = \int_\Omega u_1(x)(^tAu_2)(x)dx \qquad (1.4.11)$$

where $u_2 \otimes u_1 = (u_{2,\varepsilon}(x)u_{1,\varepsilon}(y))_\varepsilon + N_{\mathbf{c}}(\Omega \times \Omega)$.

(2) $k_A \in G(\Omega \times \Omega \backslash \Delta)$ where is the diagonal of $\Omega \times \Omega$. Moreover, for open subsets $\mathbf{W}$ and $\mathbf{W}'$ of with $\mathbf{W} \times \mathbf{W}' \subseteq \Omega \times \Omega \backslash \Delta$, and for all $u \in G_{\mathbf{c}}(\mathbf{W}')$

$$A[u(x)|_{\mathbf{W}}] = \int_\Omega k_A(x,y)u(y)dy. \qquad (1.4.12)$$

(3) If $\mathbf{supp}_{x,y}(a) \in \Omega \times \Omega \backslash \Omega'$, where $\Omega'$ is an open neighborhood of $\Delta$, then $k_A \in G(\Omega \times \Omega)$.

(4) If $(\omega^{-1}(\varepsilon))_\varepsilon$ is a slow scale net then (2) and (3) are valid with $G^\infty(\Omega \times \Omega \backslash \Delta)$ and $G^\infty(\Omega \times \Omega)$ in place of $G(\Omega \times \Omega \backslash \Delta)$ and $G(\Omega \times \Omega)$ respectively.

(5) If $a \in \tilde{\mathbf{S}}_{\mathbf{rg}}^{-\infty}(\Omega \times \Omega \times \mathbb{R}^p)$ then $k_A \in G^\infty(\Omega \times \Omega)$.

**Theorem 1.4.3**. Let $a \in \tilde{\mathbf{S}}_{\rho,\delta,\omega}^{m,\mu}(\Omega \times \Omega \times \mathbb{R}^n)$ and $(\omega^{-1}(\varepsilon))_\varepsilon$ a slow scale net. Let $A$ be the corresponding pseudodifferential operator. Then for all $u \in G_{\mathbf{c}}(\Omega)$

$$\mathrm{singsupp}_{\mathbf{g}}(Au) \subseteq \mathrm{singsupp}_{\mathbf{g}}(u). \qquad (1.4.13)$$

**Theorem 1.4.4**. Let us consider a linear operator $A : G_{\mathbf{c}}(\Omega) \to G(\Omega)$ of the form

$$A[u(x)] = \int_\Omega k_A(x,y)u(y)dy. \tag{1.4.14}$$

where $k_A \in G^\infty(\Omega \times \Omega).A$ is an operator with regular generalized kernel if and only if it is a pseudodifferential operator with smoothing amplitude in $\tilde{\mathbf{S}}_{\mathbf{rg}}^{-\infty}(\Omega \times \Omega \times \mathbb{R}^p)$.

**Definition 1.4.8.** A pseudodifferential operator $A$ is properly supported if and only if $\mathrm{supp}_{x,y}(k_A)$ is a proper set. An amplitude $a \in \tilde{\mathbf{S}}_{\rho,\delta,\omega}^{m,\mu}(\Omega \times \Omega \times \mathbb{R}^n)$ is called properly supported if and only if $\mathrm{supp}_{x,y}(a)$ is a proper set.

Note that if $A$ is properly supported then ${}^t\!A$ is properly supported.

**Theorem 1.4.4.** Let $A$ be a pseudodifferential operator with amplitude $a \in \tilde{\mathbf{S}}_{\rho,\delta,\omega}^{m,\mu}(\Omega \times \Omega \times \mathbb{R}^n).(1)$ If $(\omega(\varepsilon))_\varepsilon$ is bounded then $A$ is properly supported if and only if it can be written with a properly supported amplitude in $\tilde{\mathbf{S}}_{\rho,\delta,\omega}^{m,\mu}(\Omega \times \Omega \times \mathbb{R}^n)$. (2) If $(\omega^{-1}(\varepsilon))_\varepsilon$ is a slow scale net then A is properly supported if and only if it can be written with a properly supported amplitude in $\tilde{\mathbf{S}}_{\rho,\mathbf{rg}}^{m}(\Omega \times \Omega \times \mathbb{R}^n)$.

# I.6. Oscillatory integrals with generalized phase function on $\mathbb{R}^n$.

In this section we argue the most important properties of the integral, depending on a real parameter $\varepsilon$, of the type

$$\left( \int_{\mathbb{R}^n} a_\varepsilon(y) \exp[i\varphi_\varepsilon(y)] dy \right)_\varepsilon. \tag{1.6.1}$$

We assume that $\varphi_\varepsilon \in C^\infty(\mathbb{R}^n \backslash \{0\}), \varepsilon \in (0,1]$ is a phase function of order $k > 0$, $\forall \varepsilon(\varepsilon \in (0,1])[\varphi_\varepsilon \in \Phi^k(\mathbb{R}^n)]$ for short, if it is real valued, positively homogeneous of order $k$, i.e. $\varphi_\varepsilon(tx) = t^k \varphi_\varepsilon(x)$ for $t > 0, \varepsilon \in (0,1]$ and $\nabla\varphi_\varepsilon \neq 0$ for $x \neq 0$.

**Definition 1.6.1.** A function $a \in C^\infty(\mathbb{R}^n)$ is an element of $\mathbf{A}_l^m(\mathbb{R}^n), l, m \in \mathbb{R}$, iff

$$\sup_{y \in \mathbb{R}^n} \langle y \rangle^{l|\alpha|-m} |\partial^\alpha a(y)| < +\infty. \tag{1.6.2}$$

We denote the set of nets $(a_\varepsilon)_\varepsilon \in \mathbf{A}_l^m(\mathbb{R}^n)^{(0,1]}$ by $\mathbf{A}_l^m[\mathbb{R}^n]$. Let $\mathbf{A}_{l,\mathbf{M}}^m(\mathbb{R}^n)$ be the set of all nets $(a_\varepsilon)_\varepsilon \in \mathbf{A}_l^m[\mathbb{R}^n]$ such that

$$\forall \alpha (\alpha \in \mathbb{N}^n) \exists N(N \in \mathbb{N}) \left[ \sup_{y \in \mathbb{R}^n} \langle y \rangle^{l|\alpha|-m} |\partial^\alpha a_\varepsilon(y)| = O(\varepsilon^{-N}) \text{ as } \varepsilon \to 0 \right] \tag{1.6.3}$$

and $\mathbf{N}_l^m(\mathbb{R}^n)$ be the set of all nets $(a_\varepsilon)_\varepsilon \in \mathbf{A}_l^m[\mathbb{R}^n]$ such that

$$\forall \alpha (\alpha \in \mathbb{N}^n) \forall q(q \in \mathbb{N}) \left[ \sup_{y \in \mathbb{R}^n} \langle y \rangle^{l|\alpha|-m} |\partial^\alpha a_\varepsilon(y)| = O(\varepsilon^q) \text{ as } \varepsilon \to 0 \right]. \tag{1.6.4}$$

The classes $a = \mathbf{cl}[(a_\varepsilon)_\varepsilon]$ of the factor space $\tilde{\mathbf{A}}_l^m(\mathbb{R}^n) = \mathbf{A}_{l,\mathbf{M}}^m(\mathbb{R}^n)/\mathbf{N}_l^m(\mathbb{R}^n)$ are called generalized amplitudes on $\mathbb{R}^n$.

**Definition 1.6.2.** A function $\varphi \in C^\infty(\mathbb{R}^n \backslash \{0\})$ is an element of $\mathbf{A}_{l,\mathbf{ph}}^m(\mathbb{R}^n \backslash \{0\}), l, m \in \mathbb{R}$, iff

$$\nabla \varphi \neq 0 \text{ for } x \neq 0,$$

$$\sup_{y \in \mathbb{R}^n} \langle y \rangle^{l|\alpha|-m} |\partial^\alpha \varphi(y)| < +\infty. \tag{1.6.5}$$

We denote the set of nets $(\varphi_\varepsilon)_\varepsilon \in \mathbf{A}_{l,\mathbf{ph}}^m(\mathbb{R}^n \backslash \{0\})^{(0,1]}$ by $\mathbf{A}_{l,\mathbf{ph}}^m[\mathbb{R}^n]$. Let $\mathbf{A}_{l,\mathbf{ph},\mathbf{M}}^m(\mathbb{R}^n)$ be the set of all nets $(\varphi_\varepsilon)_\varepsilon \in \mathbf{A}_{l,\mathbf{ph}}^m[\mathbb{R}^n]$ such that

$$\forall \alpha (\alpha \in \mathbb{N}^n) \exists N (N \in \mathbb{N}) \left[ \sup_{y \in \mathbb{R}^n} \langle y \rangle^{l|\alpha|-m} |\partial^\alpha \varphi_\varepsilon(y)| = O(\varepsilon^{-N}) \text{ as } \varepsilon \to 0 \right]. \qquad (1.6.6)$$

The classes $\varphi = \mathbf{cl}[(\varphi_\varepsilon)_\varepsilon]$ of the factor space $\tilde{\mathbf{A}}^m_{l,\mathbf{ph}}(\mathbb{R}^n \backslash \{0\}) = \mathbf{A}^m_{l,\mathbf{ph},\mathbf{M}}(\mathbb{R}^n \backslash \{0\})/\mathbf{N}^m_l(\mathbb{R}^n \backslash \{0\})$ are called generalized phase function on $\mathbb{R}^n$.

**Definition 1.6.3**. An element $a \in \tilde{\mathbf{A}}^m_l(\mathbb{R}^n)$ is called regular if it has a representative $(a_\varepsilon)_\varepsilon$ fulfilling the following condition:

$$\exists N (N \in \mathbb{N}) \forall \alpha (\alpha \in \mathbb{N}^n) \left[ \sup_{y \in \mathbb{R}^n} \langle y \rangle^{l|\alpha|-m} |\partial^\alpha a_\varepsilon(y)| = O(\varepsilon^{-N}) \text{ as } \varepsilon \to 0 \right]. \qquad (1.6.7)$$

We denote the set of all regular generalized amplitudes on $\mathbb{R}^n$ by $\tilde{\mathbf{A}}^m_{l,\mathbf{rg}}(\mathbb{R}^n)$

**Definition 1.6.4**. An element $\varphi \in \tilde{\mathbf{A}}^m_{l,\mathbf{ph}}(\mathbb{R}^n \backslash \{0\})$ is called regular if it has a representative $(\varphi_\varepsilon)_\varepsilon$ fulfilling the following condition:

$$\exists N (N \in \mathbb{N}) \forall \alpha (\alpha \in \mathbb{N}^n) \left[ \sup_{y \in \mathbb{R}^n} \langle y \rangle^{l|\alpha|-m} |\partial^\alpha \varphi_\varepsilon(y)| = O(\varepsilon^{-N}) \text{ as } \varepsilon \to 0 \right]. \qquad (1.6.8)$$

We denote the set of all regular generalized phase functions on $\mathbb{R}^n \backslash \{0\}$ by $\tilde{\mathbf{A}}^m_{l,\mathbf{ph},\mathbf{rg}}(\mathbb{R}^n \backslash \{0\})$.

**Remark**.**1.6.1**.Let $\mathbf{A}^m_{l,N}(\mathbb{R}^n)$ be the set of nets in $\mathbf{A}^m_l[\mathbb{R}^n]$ defined by (1.6.7) by fixing $N$, and $\mathbf{A}^m_{l,\mathbf{rg}}(\mathbb{R}^n) = \cup_{N \in \mathbb{N}} \mathbf{A}^m_{l,N}(\mathbb{R}^n)$. Since (1.6.7) holds for all representatives of $a$ once it holds for one, we can introduce $\tilde{\mathbf{A}}^m_{l,\mathbf{rg}}(\mathbb{R}^n)$ as the quotient $\tilde{\mathbf{A}}^m_{l,\mathbf{rg}}(\mathbb{R}^n) = \mathbf{A}^m_{l,\mathbf{rg}}(\mathbb{R}^n)/\mathbf{N}^m_l(\mathbb{R}^n)$.

**Theorem 1.6.1**.Let $(a_\varepsilon)_\varepsilon \in \mathbf{A}^m_l[\mathbb{R}^n]$ and $\varphi_\varepsilon \in \Phi^k(\mathbb{R}^n)$ with $\varepsilon \in (0,1]$ $1-k < l \le 1$. Let $\psi$ be an arbitrary function of $\mathcal{L}(\mathbb{R}^n)$ with $\psi(0) = 1$ and $\phi$ any function in $G^\infty_c(\mathbb{R}^n)$

**Definition 1.6.5**. We define the generalized oscillatory integral of the generalized phase function $\varphi = \mathbf{cl}[(\varphi_\varepsilon)_\varepsilon]$ and amplitude $a = \mathbf{cl}[(a_\varepsilon)_\varepsilon]$ as the class in $\widetilde{\mathbb{C}}$ :

$$I_\varphi[a] = \int_{\mathbb{R}^n} a(y) \exp[i\varphi(y)] dy = \mathbf{cl}\left[\left(\int_{\mathbb{R}^n} a_\varepsilon(y) \exp[i\varphi_\varepsilon(y)] dy\right)_\varepsilon\right]. \qquad (1.6.9)$$

# I.7. Generalized symbols and amplitudes of weighted type.

**Definition 1.7.1.** In the following for real functions $f$ and $g$ on $\mathbb{R}^n$, we write $f(z) \preceq g(z)$ on $A \subseteq \mathbb{R}^n$, if there exists a positive constant $c > 0$ such that $f(z) \leq cg(z)$, for all $z \in A$.
If $f(z) \preceq g(z)$ and $g(z) \preceq f(z)$, we write $f(z) \sim g(z)$.

**Definition 1.7.2.** A continuous real function $\Lambda(z)$ on $\mathbb{R}^{2n}$ is a weight function if and only if:
(1) there exists $\lambda > 0$ such that $\langle z \rangle^\lambda \preceq \Lambda(z) \preceq \langle z \rangle$ on $\mathbb{R}^{2n}$.
(2) $\Lambda(z) \sim \Lambda(\zeta)$ on $A_{z,\zeta} = \{(z,\zeta) \| z - \zeta \| \leq \lambda \Lambda(z)\}$.
(3) $\forall t(t \to \mathbb{R}^{2n}) : \Lambda(t \bullet z) \preceq \Lambda(z)$ on $\mathbb{R}^{2n}$, where $t \bullet z = (t_1 z_1, \ldots, t_{2n} z_{2n})$.

**Definition 1.7.3.** Let $m \in \mathbb{R}, \rho \in (0,1), \Lambda(z)$ be a weight function and $z = (x,\xi) \in \mathbb{R}^{2n}$. We denote by $S^m_{\Lambda,\rho,\mathbf{M}}(\mathbb{R}^{2n})$ or $S^m_{\Lambda,\rho,M}$, the set of all nets $(a_\varepsilon)_\varepsilon \in E[\mathbb{R}^{2n}]$ with the property that

$$\forall \alpha \in \mathbb{N}^{2n} \exists N \in \mathbb{N} \exists c(c > 0) \forall \varepsilon(\varepsilon \in (0,1]) \forall z \in \mathbb{R}^{2n} \big[ |\partial^\alpha a_\varepsilon(z)| \leq c(\Lambda(z))^{m-\rho|\alpha|} \varepsilon^{-N} \big] \qquad (1.7.1)$$

We denote by $N^m_{\Lambda,\rho}(\mathbb{R}^{2n})$ or $N^m_{\Lambda,\rho}$, the set of all nets nets $(a_\varepsilon)_\varepsilon \in E[\mathbb{R}^{2n}]$ with the property

$$\forall \alpha \in \mathbb{N}^{2n} \forall q \in \mathbb{N} \exists c(c > 0) \forall \varepsilon(\varepsilon \in (0,1]) \forall z \in \mathbb{R}^{2n} \big[ |\partial^\alpha a_\varepsilon(z)| \leq c(\Lambda(z))^{m-\rho|\alpha|} \varepsilon^q \big] \qquad (1.7.2)$$

**Definition 1.7.4.** The elements of the factor space

$$\widetilde{S}^m_{\Lambda,\rho} = S^m_{\Lambda,\rho,\mathbf{M}}/N^m_{\Lambda,\rho} \tag{1.7.3}$$

are called generalized symbols of order $m$ and type $(\Lambda, \rho)$.

We say that a generalized symbol $a \in \widetilde{S}^m_{\Lambda,\rho}$ is regular if it has a representative $(a_\varepsilon)_\varepsilon$ satisfying the following condition

$$\exists N \in \mathbb{N} \forall \alpha \in \mathbb{N}^{2n} \exists c(c > 0) \forall \varepsilon(\varepsilon \in (0,1]) \forall z \in \mathbb{R}^{2n} \big[ |\partial^\alpha a_\varepsilon(z)| \leq c(\Lambda(z))^{m-\rho|\alpha|} \varepsilon^{-N} \big]. \tag{1.7.4}$$

The set of all regular symbols in $\widetilde{S}^m_{\Lambda,\rho}$ will be denoted by $\widetilde{S}^m_{\Lambda,\rho,\mathbf{reg}}$.

# I.8. Pseudodifferential operators acting on $G_{\tau,\mathcal{L}}(\mathbb{R}^n)$.

In the classical theory an integral

$$A[u(x)] = \int_{\mathbb{R}^{2n}} (\exp[i(x-y)]) a(x,y,\xi) u(y) d\xi, \tag{1.8.1}$$

where $a \in \bar{S}^m_{\Lambda,\rho}$ and $u \in \mathcal{L}(\mathbb{R}^n)$, defines a continuous map $\mathcal{L}(\mathbb{R}^n) \to \mathcal{L}(\mathbb{R}^n)$, which can be extended as a continuous map $\mathcal{L}'(\mathbb{R}^n) \to \mathcal{L}'(\mathbb{R}^n)$.

**Theorem 1.8.1.** Let $a \in \widetilde{\bar{S}}^m_{\Lambda,\rho}, u \in G_{\tau,\mathcal{L}}(\mathbb{R}^n)$ and $\varphi$ be a mollifier in $\mathcal{L}(\mathbb{R}^n)$. Then for every $x \in \mathbb{R}^n$

$$(\exp(ix\xi)) a(x,y,\xi) u_{\widehat{\varphi}}(y) \triangleq \mathbf{cl}\Big[ \big( (\exp(ix\xi)) a_\varepsilon(x,y,\xi) u_\varepsilon(y) \widehat{\varphi}_\varepsilon \big)_\varepsilon \Big] \tag{1.8.2}$$

is a well-defined element of $\bar{\mathbf{A}}^\nu_0(\mathbb{R}^n_y \times \mathbb{R}^n_\xi)$ with $\nu = m_+ + m'_+$.

Therefore one can define the oscillatory integral

$$\int_{\mathbb{R}^{2n}} (\exp[i(x-y)]) a(x,y,\xi) u_{\widehat{\varphi}}(y) dy d\xi =$$

$$(1.8.3)$$

$$\int_{\mathbb{R}^{2n}} \exp(-iy\xi) \Big( a(x,y,\xi) u_{\widehat{\varphi}}(y) \exp(ix\xi) \Big) dy d\xi.$$

where $-y\xi$ is the phase function of order 2 and $a(x,y,\xi) u_{\widehat{\varphi}}(y) \exp(ix\xi)$ the generalized amplitude in $\widetilde{A}_0^y(\mathbb{R}_y^n \times \mathbb{R}_\xi^n)$.

# II.1.1.Feynman-Colombeau path integral.

Let us consider $L_2$ transition probability amplitude $\mathbf{K}(\vec{x}, t|\vec{x}_0, 0)$

$$\mathbf{K}(\vec{x}, t|\vec{x}_0, 0) = \left\langle \psi_1^*, \exp\left(-\frac{it\widehat{H}}{\epsilon}\right)\psi_2 \right\rangle =$$

$$(2.1.1.1)$$

$$\int_{\mathbb{R}^n} \left( \psi_1^*(\vec{x}) \left[ \exp\left(-\frac{it\widehat{H}}{\epsilon}\right)\psi_2(\vec{x}) \right] \right) d\vec{x},$$

where $\psi_1, \psi_2 \in L_2$,

$$H = -\frac{\epsilon^2}{2m}\Delta + V(\vec{x}, t) \qquad (2.1.1.2)$$

and $H$ is essentially self-adjoint, $\widehat{H}$ is the closure of $H$. Propagator (2.1.1.1) is formally given via formula:

$$\mathbf{K}(\vec{x}, t | \vec{x}_0, 0) = \lim_{k \to \infty} \Omega_k(n) \int_{\mathbb{R}^{(k-1)n}} \exp\left[ \frac{i\Delta\tau}{\epsilon} S_k(\vec{x}, \dots \vec{x}_0) \right] d\vec{x}_1 \dots d\vec{x}_{k-1}, \qquad (2.1.1.3)$$

where

$$\Omega_{n,k} = \left( \frac{m}{2i\pi\epsilon\Delta\tau} \right)^{\frac{nk}{2}}, \Delta\tau = \frac{t}{k},$$

$$\qquad (2.1.1.4)$$

$$S_k(\vec{x}, \dots \vec{x}_0) = \sum_{j=1}^{k} \left[ \frac{m}{2} \left( \frac{\vec{x}_j - \vec{x}_{j-1}}{\Delta\tau} \right) - V(\vec{x}_j) \right], \vec{x} = \vec{x}_k$$

Let $\mathbf{H}$ be a complex separable infinite dimensional Hilbert space, with inner product $\langle , \rangle$ and norm $|\cdot|$. Let us consider Schrödinger equation:

$$-i\frac{\partial \Psi(t)}{\partial t} + \widehat{H}(t)\Psi(t) = 0,$$

$$\qquad (2.1.1.5)$$

$$\Psi(0) = \Psi_0\left(\vec{x}\right),$$

where operator $\widehat{H}(t)$ given via formula (2.1.2), $\Psi(t) : \mathbb{R} \to \mathbf{H}$ and $\widehat{H}(t) : \mathbb{R} \times \mathbf{H} \to \mathbf{H}$.

**Theorem 2.1.1**. [20].We assume that:
(1) $\Psi_0\left(\vec{x}\right) \in L_2(\mathbb{R}^n)$,
(2) $V\left(\vec{x}\right)$ is continuous and $\sup_{\vec{x} \in \mathbb{R}^n}\left(\left|V\left(\vec{x}\right)\right|\right) < +\infty$.
Then corresponding solution of the Schrödinger equation (2.1.5) exist and can be represented via formulae (2.1.1.3)-(2.1.1.4).

Let $\mathbf{H}$ be a complex separable infinite dimensional Hilbert space, with inner product

$\langle , \rangle$ and norm $|\cdot|$. Let us consider Colombeau-Schrödinger type equation

$$-i\left(\frac{\partial \Psi_\varepsilon(t)}{\partial t}\right)_{\varepsilon \in (0,1]} + \left(\widehat{H}_\varepsilon(t)\Psi_\varepsilon(t)\right)_{\varepsilon \in (0,1]} = 0, \qquad (2.1.1.6)$$

where: (1) $\forall \varepsilon(\varepsilon \neq 0)[\Psi_\varepsilon(t) : \mathbb{R} \to \mathbf{H}]$, (2) $\widehat{H}_0(t) = \widehat{H}(t)$ and operator $\widehat{H}(t)$ is given via Eq.(2.1.1.2) and $\sup_{\vec{x} \in \mathbb{R}^n}\left(\left|V(\vec{x})\right|\right) = \infty$,

(3) $\forall \varepsilon(\varepsilon \neq 0)\left[\widehat{H}_\varepsilon(t) : \mathbb{R} \times \mathbf{H} \to \mathbf{H}\right]$,(4) $\mathbf{cl}\left(\left(\widehat{H}_\varepsilon(t)\right)_\varepsilon\right) : \widetilde{\mathbb{R}} \times G_\mathbf{H} \to G_\mathbf{H}$.

From Eq.(2.1.1.6) one obtain formal "infinitesimal" solution

$$\left(\Psi_\varepsilon(t + dt)\right)_{\varepsilon \in (0,1]} = \left(\exp\left\{-idt\left[\widehat{H}_\varepsilon(t)\right]\right\}\Psi_\varepsilon(t)\right)_{\varepsilon \in (0,1]} \qquad (2.1.1.7)$$

Using simple formal calculation one obtain formal solution of the Eq.(2.1.1.6) via formula :

$$\left(\Psi_\varepsilon(t)\right)_{\varepsilon \in (0,1]} =$$

$$\left(\lim_{\delta_N \to 0}\left[\exp\left\{-i\Delta t_{N-1}\widehat{H}_\varepsilon(t_{N-1})\right\}\ldots\exp\left\{-i\Delta t_0\widehat{H}_\varepsilon(t_0)\right\}\Psi_\varepsilon(t_0)\right]\right)_{\varepsilon \in (0,1]}, \qquad (2.1.1.8)$$

where $t_0 < t_1 < \ldots < t_N = t, \Delta t_j = t_{j+1} - t_j, \delta_N = \max_{0 \leqslant j \leqslant N-1} \Delta t_j$
If $\forall \varepsilon[\varepsilon \in (0,1]]$ limit $\Psi_\varepsilon(t) = \lim_{\delta_N \to 0}[\circ]$ in the formula (2.1.1.8) exists we obtain Colombeau generalized function $\left(\Psi_\varepsilon(t)\right)_{\varepsilon \in (0,1]}$ We denote RH of Eq.(2.1.1.8) as

$$\left(\Psi_\varepsilon(t)\right)_{\varepsilon \in (0,1]} = \left(\left[\prod_{\tau=t_0}^{t}\exp\left\{-id\tau\widehat{H}_\varepsilon(\tau)\right\}\right]\Psi_\varepsilon(t_0)\right)_{\varepsilon \in (0,1]} \qquad (2.1.1.9)$$

and call it *Maslov-Colombeau chronological exponential* or *Maslov-Colombeau* $T$-exponent.
Let $\mathbf{B}$ be a complex separable infinite dimensional Banach space with the norm

$\|\circ\|$. The Banach space $\mathbf{C}_{\mathbf{B}}^{t} = \mathbf{C}([0,t], \mathbf{B})$ consist of continuous functions $\psi(t), \psi : [0,t] \to \mathbf{B}$, with the norm

$$\|\psi\|_{\mathbf{C}_{\mathbf{B}}^{t}} = \sup_{0 \le \tau \le t} \|\psi(\tau)\|_{\mathbf{B}}. \qquad (2.1.1.10)$$

If $t' < t$ we denote $\psi'(t) \triangleq \psi(t) \upharpoonright [0,t']$, thus $\psi'(t) \in \mathbf{C}([0,t'], \mathbf{B})$. Suppose that:

(1) $\varepsilon_0 > 0, T_0 > 0$,
(2) $\forall t, \widetilde{\varepsilon}$ such that $0 \le t \le T, (-t) \le \widetilde{\varepsilon} \le \varepsilon_0$ there exists a function

$$S_{t,\widetilde{\varepsilon}}[\circ] : \mathbf{C}([0,t'], \mathbf{B}) \times [-t, \varepsilon_0] \to \mathbf{B} \qquad (2.1.1.11)$$

such that $\forall v \in \mathbf{C}([0,t'], \mathbf{B}), \forall t \in [0, T_0]$ function $q(t, \widetilde{\varepsilon}) = \|S_{t,\widetilde{\varepsilon}}[v]\|_{\mathbf{B}}$ is continuous on variable $\widetilde{\varepsilon}$,
(3) if $\widetilde{\varepsilon} \le 0$, then $S_{t,\widetilde{\varepsilon}}[v] = S_t[v] \triangleq v(t)$.
Let's consider now finite splitting $\Delta_N$ of the closed interval $[0, T_0]$ :

$$\Delta_N = \{0 = t_0 < t_1 < \ldots < t_N = T_0\}. \qquad (2.1.1.12)$$

Note that for every $v_0 \in \mathbf{B}$ there exists an sequence $v_j(t), j = 0, 1, \ldots, N-1$ such that $v_{j+1}(t) = S_{t_j, t-t_j}[v_j], 0 \le t \le t_{j+1}$.
**Definition 2.1.1.1**. [20].

$$v_N(t) = S_{t_N, \Delta t_{N-1}} \circ \ldots \circ S_{t_j, \Delta t_j} \circ \ldots \circ S_{t_0, \Delta t_0}[v_0],$$
$$\qquad (2.1.1.13)$$
$$\Delta t_j = t_{j+1} - t_j, j = 0, 1, \ldots, N-1.$$

**Definition 2.1.1.2**. [20] Suppose that:
(1) $v_0 \in D \subsetneqq \mathbf{B}$ and
(2) $s\text{-}\lim_{\delta_N \to 0} v_N(t)$ there exists.

We denote this limit by

$$v(t) = \left( \prod_{\tau=0}^{t} S_{\tau, d\tau} \right) [v_0] \qquad (2.1.1.14)$$

and call it *chronological exponential with a generator* $S_{t, \widetilde{z}}[v]$.

Let $\mathbf{B}_1 \subsetneqq \mathbf{B}_2$ be a complex separable infinite dimensional Banach space with norm $\| \circ \|_1 \geq \| \circ \|_2$ and $\mathbf{C}_{\mathbf{B}_j}^t = \mathbf{C}([0, t], \mathbf{B}_j), j = 1, 2.$ Thus $\mathbf{C}_{\mathbf{B}_1}^t \subsetneqq \mathbf{C}_{\mathbf{B}_2}^t.$
Suppose that for every $t \in [0, T_0]$ there exists continuous operator

$$\Re_t[\circ] : \mathbf{C}_{\mathbf{B}_1}^t \to \mathbf{B}_2 \qquad (2.1.1.15)$$

such that for every $v \in \mathbf{C}_{\mathbf{B}_1}^t$ function $\| \Re_t[v] \|_{\mathbf{B}_2}$ is continuous on a variable $t \in [0, T_0].$
Let's consider the Cauchy problem with initial condition $\Phi_0 \in \mathbf{B}_1$ :

$$\frac{d}{dt} \Phi(t) = \Re_t[\Phi(t)], \Phi(0) = \Phi_0, \qquad (2.1.1.16)$$

where

$$\frac{d}{dt} \Phi(t) \triangleq \lim_{\Delta t \to 0} \left\| \frac{\Phi(t + \Delta t) - \Phi(t)}{\Delta t} \right\|_2, \qquad (2.1.1.17)$$

i.e. $\forall t$ and $\forall v_t \in \mathbf{B}_2$

$$\frac{d}{dt} \Phi(t) = v_t \Leftrightarrow \lim_{\Delta t \to 0} \left\| v_t - \frac{\Phi(t + \Delta t) - \Phi(t)}{\Delta t} \right\|_2 = 0. \qquad (2.1.1.18)$$

**Definition 2.1.1.3.**[20]. Let be $S_{t,\widetilde{\varepsilon}}[\circ] : \mathbf{C}([0,t], \mathbf{B}) \times [-t, \varepsilon_0] \to \mathbf{B}$ the infinitesimal generator of the corresponding $T$-exponent. Suppose that for every $t \in [0, T_0]$

the next condition is satisfied

$$\lim_{\widetilde{\varepsilon} \to 0} \left\| \frac{S_{t,\widetilde{\varepsilon}}[v] - v(t)}{\widetilde{\varepsilon}} - \Re_t[v] \right\|_2 = 0. \tag{2.1.19}$$

locally uniformly on variable $v(t) \in \mathbf{C}_{\mathbf{B}_1}^t$. We call operator function $S_{t,\widetilde{\varepsilon}}[\circ]$ as master operator for Cauchy problem (2.1.16).

**Theorem 2.1.1.2.**[20]. Let be $S_{t,\widetilde{\varepsilon}}[\circ]$ the master operator for Cauchy problem (2.1.1.16)

and for some $\psi_0 \in \mathbf{B}_1$ there exists $\psi(t) = \left( \prod_{\tau=0}^t S_{\tau, d\tau} \right) [\psi_0], \psi(t) \in \mathbf{C}_1(T_0)$.

Then $\psi(t)$ is a solution of the Eq.(2.1.1.16).

**Proof**. Let's consider finite splitting $\Delta_N = \{0 = t_0 < t_1 < \ldots < t_N = T_0\}$ of the closed interval $[0, T_0]$ and $v_N(t) = S_{t_N, \Delta t_{N-1}} \circ \ldots \circ S_{t_j, \Delta t_j} \circ \ldots \circ S_{t_0, \Delta t_0}[\psi_0]$. Then for every

$t_{j+1} \geq t \geq t_j$, one obtain $\psi_N(t) = S_{t_j, t-t_j}[\psi_N]$. Thus there exists $N_0$ such that

$$\lim_{\Delta t \to +0} \left\| \frac{\psi_N(t_j + \Delta t) - \psi_N(t_j)}{\Delta t} - \Re_{t_j}[\psi_N] \right\|_2 = 0 \tag{2.1.1.20}$$

uniformly on $N \geq N_0$.

Let $t$ be a fixed point $t \in [0, T_0)$ and let $\{\Delta_N\}_{N=1}^\infty$ be the sequence of the splittings such that for any $N : t \in \Delta_N$.

Now take into account that operator $\Re_t$ is continuous finally we obtain

$$\lim_{\Delta t \to +0} \left\| \frac{\psi(t + \Delta t) - \psi(t)}{\Delta t} - \Re_t[\psi] \right\|_2 = 0. \tag{2.1.1.21}$$

Let $\Re_t[\psi]$ be operator such that

$$\Re_t[\psi] = A_t[\psi]\psi(t). \tag{2.1.1.22}$$

where $\forall \psi(\psi \in \mathbf{C}_{\mathbf{B}_1}^t)$ operator $A_t[\psi] : \mathbf{B}_1 \to \mathbf{B}_2$ is continuous. Thus $A_t[\psi]$ is a infinitesimal generator of the one-parameter semigroup

$$\exp[\widetilde{\varepsilon} A_t[\psi]].\qquad(2.1.1.23)$$

Hence master operator of the Cauchy problem (2.1.1.16) is

$$S_{t,\widetilde{\varepsilon}}[v] = \exp[\widetilde{\varepsilon} A_t[\psi]]v(t).\qquad(2.1.1.24)$$

Corresponding $T$-exp is

$$\psi(t) = \prod_{\tau}^{t} \exp[A_\tau[\psi]]d\tau.\qquad(2.1.1.25)$$

**Definition 2.1.1.4.** Let **B** be a complex separable infinite dimensional Banach space with the norm $\|\circ\|$. The $\widetilde{\mathbb{C}}$-module $\mathbf{C}_{G_\mathbf{B}}^t = \mathbf{C}([0,t], G_\mathbf{B})$ consist of continuous functions $(\psi_\varepsilon(t))_\varepsilon, (\psi_\varepsilon)_\varepsilon : [0,t] \to G_\mathbf{B}$, with the ultra-pseudo-norm $\|\circ\|_{\mathbf{C}_{G_\mathbf{B}}^t}$ :

$$\|(\psi_\varepsilon)_\varepsilon\|_{\mathbf{C}_{G_\mathbf{B}}^t} = \sup_{0\leq\tau\leq t} \|(\psi_\varepsilon(\tau))_\varepsilon\|_{G_\mathbf{B}}.\qquad(2.1.1.26)$$

If $t' < t$ we denote $(\psi'_\varepsilon(t))_\varepsilon \triangleq (\psi_\varepsilon(t))_\varepsilon \upharpoonright [0,t']$, thus $(\psi'_\varepsilon(t))_\varepsilon \in \mathbf{C}([0,t'], G_\mathbf{B})$.

Suppose that:

**(1)** $\varepsilon_0 > 0, T_0 > 0$,

**(2)** $\forall t, \widetilde{\varepsilon}$ such that $0 \leq t \leq T, (-t) \leq \widetilde{\varepsilon} \leq \varepsilon_0$ there exists Colombeau generalized function

$$(S_{t,\widetilde{\varepsilon},\varepsilon}[\circ])_{\varepsilon\in(0,1]} : \mathbf{C}([0,t'], G_\mathbf{B}) \times [-t, \varepsilon_0] \to G_\mathbf{B}\qquad(2.1.1.27)$$

such that $\forall (v_\varepsilon)_\varepsilon \in \mathbf{C}([0, t'], G_\mathbf{B}), \forall t \in [0, T_0]$ function $(q_\varepsilon(t, \widetilde{\varepsilon}))_\varepsilon = \|(S_{t,\widetilde{\varepsilon},\varepsilon}[v_\varepsilon])_\varepsilon\|_{G_\mathbf{B}}$ is continuous on variable $\widetilde{\varepsilon}$,

(3) if $\widetilde{\varepsilon} \leq 0$, then $(S_{t,\widetilde{\varepsilon},\varepsilon}[v_\varepsilon])_\varepsilon = (S_{t,\varepsilon}[v_\varepsilon])_\varepsilon \triangleq (v_\varepsilon(t))_\varepsilon$.

Let's consider now finite splitting $\Delta_N$ of the closed interval $[0, T_0]$ :

$$\Delta_N = \{0 = t_0 < t_1 < \ldots < t_N = T_0\}. \tag{2.1.1.28}$$

Note that for every $(v_{0,\varepsilon})_\varepsilon \in G_\mathbf{B}$ there exists an sequence $(v_{j,\varepsilon}(t))_\varepsilon, j = 0, 1, \ldots, N-1$ such that $(v_{j+1,\varepsilon}(t))_\varepsilon = (S_{t_j, t-t_j, \varepsilon}[v_{j,\varepsilon}])_\varepsilon, \ 0 \leq t \leq t_{j+1}$.

**Definition 2.1.1.5**.

$$(v_{N,\varepsilon}(t))_\varepsilon = \left( S_{t_N, \Delta t_{N-1}, \varepsilon} \circ \ldots \circ S_{t_j, \Delta t_j, \varepsilon} \circ \ldots \circ S_{t_0, \Delta t_0, \varepsilon}[v_{0,\varepsilon}] \right)_\varepsilon,$$

$$\tag{2.1.1.29}$$

$$\Delta t_j = t_{j+1} - t_j, \ j = 0, 1, \ldots, N-1.$$

**Definition 2.1.1.6**. Suppose that:

(1) $(v_{0,\varepsilon})_\varepsilon \in D \subsetneqq G_\mathbf{B}$ and

(2) $s\text{-}\lim_{\delta_N \to 0} (v_{N,\varepsilon}(t))_\varepsilon$ there exists.

We denote this limit by

$$(v_\varepsilon(t))_\varepsilon = \left( \left( \prod_{\tau=0}^{t} S_{\tau, d\tau, \varepsilon} \right) [v_{0,\varepsilon}] \right)_\varepsilon \tag{2.1.1.30}$$

and call it *chronological exponential with a generator* $(S_{t,\widetilde{\varepsilon},\varepsilon}[v_\varepsilon])_\varepsilon$.

Let $\mathbf{B}_1 \subsetneqq \mathbf{B}_2$ be a complex separable infinite dimensional Banach spaces with norm $\|\circ\|_1 \geq \|\circ\|_2$ and let $\mathbf{C}^t_{G_{\mathbf{B}_j}}$ be a $\widetilde{\mathbb{C}}$-modules $\mathbf{C}^t_{G_{\mathbf{B}_j}} = \mathbf{C}([0, t], G_{\mathbf{B}_j}), j = 1, 2$. Thus $\mathbf{C}^t_{G_{\mathbf{B}_1}} \subsetneqq \mathbf{C}^t_{G_{\mathbf{B}_2}}$. Suppose that for every $t \in [0, T_0]$ there exists continuous operator

$$(\Re_{t,\varepsilon}[\circ])_\varepsilon : \mathbf{C}^t_{G_{\mathbf{B}_1}} \to G_{\mathbf{B}_2} \tag{2.1.1.31}$$

such that for every $(v_\varepsilon)_\varepsilon \in \mathbf{C}^t_{G_{\mathbf{B}_1}}$ function $\|(\Re_{t,\varepsilon}[v_\varepsilon])_\varepsilon\|_{G_{\mathbf{B}_2}}$ is continuous on a

variable $t \in [0, T_0]$.

Let's consider the Cauchy problem with initial condition $(\Phi_{0,\varepsilon})_\varepsilon \in G_{\mathbf{B}_1}$ :

$$\frac{d}{dt}(\Phi_\varepsilon(t))_\varepsilon = (\Re_{t,\varepsilon}[\Phi_\varepsilon(t)])_\varepsilon, (\Phi_\varepsilon(0))_\varepsilon = (\Phi_{0,\varepsilon})_\varepsilon, \qquad (2.1.1.32)$$

where

$$\frac{d}{dt}(\Phi_\varepsilon(t))_\varepsilon \triangleq \lim_{\Delta t \to 0} \left\| \left( \frac{\Phi_\varepsilon(t + \Delta t) - \Phi_\varepsilon(t)}{\Delta t} \right)_\varepsilon \right\|_{G_{\mathbf{B}_2}}, \qquad (2.1.1.33)$$

i.e. $\forall t$ and $\forall (v_{t,\varepsilon})_\varepsilon \in G_{\mathbf{B}_2}$

$$\frac{d}{dt}(\Phi_\varepsilon(t))_\varepsilon = (v_{t,\varepsilon})_\varepsilon \Longleftrightarrow \lim_{\Delta t \to 0} \left\| \left( v_{t,\varepsilon} - \frac{\Phi_\varepsilon(t + \Delta t) - \Phi_\varepsilon(t)}{\Delta t} \right)_\varepsilon \right\|_{G_{\mathbf{B}_2}} = 0. \qquad (2.1.1.34)$$

**Definition 2.1.1.7.** Let be $(S_{t,\widetilde{\varepsilon},\varepsilon}[\circ])_\varepsilon : \mathbf{C}([0,t], G_{\mathbf{B}}) \times [-t, \varepsilon_0] \to G_{\mathbf{B}}$ the infinitesimalgenerator of the corresponding $T$-exponent. Suppose that for every $t \in [0, T_0]$ the next condition is satisfied

$$\lim_{\widetilde{\varepsilon} \to 0} \left\| \left( \frac{S_{t,\widetilde{\varepsilon},\varepsilon}[v] - v_\varepsilon(t)}{\widetilde{\varepsilon}} - \Re_{t,\varepsilon}[v_\varepsilon] \right)_\varepsilon \right\|_{G_{\mathbf{B}_2}} = 0. \qquad (2.1.1.35)$$

locally uniformly on variable $(v_\varepsilon(t))_\varepsilon \in G_{\mathbf{C}^t_{\mathbf{B}_1}}$. We call operator function $(S_{t,\widetilde{\varepsilon},\varepsilon}[\circ])_\varepsilon$ as *master operator* for Cauchy problem (2.1.1.32).

**Theorem 2.1.1.3.** Let be $(S_{t,\widetilde{\varepsilon},\varepsilon}[\circ])_\varepsilon$ the master operator for Cauchy problem (2.1.1.32) and for some $\psi_{0,\varepsilon} \in G_{\mathbf{B}_1}$ there exists $(\psi_\varepsilon(t))_\varepsilon = \left( \left( \prod_{\tau=0}^t S_{\tau, d\tau, \varepsilon} \right) [\psi_{0,\varepsilon}] \right)_\varepsilon$, $(\psi_\varepsilon(t))_\varepsilon \in G_{\mathbf{C}_1(T_0)}.$ Then $(\psi_\varepsilon(t))_\varepsilon$ is a solution of the Eq.(2.1.1.32).

**Proof**. Let's consider finite splitting $\Delta_N = \{0 = t_0 < t_1 < \ldots < t_N = T_0\}$ of the closed interval $[0, T_0]$ and $(v_{N,\varepsilon}(t))_\varepsilon = (S_{t_N,\Delta t_{N-1},\varepsilon} \circ \ldots \circ S_{t_j,\Delta t_j,\varepsilon} \circ \ldots \circ S_{t_0,\Delta t_0,\varepsilon}[\psi_{0,\varepsilon}])_\varepsilon.$ Then for every

$t_{j+1} \geq t \geq t_j,$ one obtain $(\psi_{N,\varepsilon}(t))_\varepsilon = (S_{t_j,t-t_j,\varepsilon}[\psi_{N,\varepsilon}])_\varepsilon.$ Thus there exists $N_0$ such that

$$\lim_{\Delta t \to +0} \left\| \left( \frac{\psi_{N,\varepsilon}(t_j + \Delta t) - \psi_{N,\varepsilon}(t_j)}{\Delta t} - \mathfrak{R}_{t_j,\varepsilon}[\psi_{N,\varepsilon}] \right)_\varepsilon \right\|_{G_{\mathbf{B}_2}} = 0 \qquad (2.1.1.36)$$

uniformly on $N \geq N_0$.

Let $t$ be a fixed point $t \in [0, T_0)$ and let $\{\Delta_N\}_{N=1}^\infty$ be the sequence of the splittings such that for any $N : t \in \Delta_N$.
Now take into account that operator $(\mathfrak{R}_{t,\varepsilon})_\varepsilon$ is continuous finally we obtain

$$\lim_{\Delta t \to +0} \left\| \left( \frac{\psi_\varepsilon(t + \Delta t) - \psi_\varepsilon(t)}{\Delta t} - \mathfrak{R}_{t,\varepsilon}[\psi_\varepsilon] \right)_\varepsilon \right\|_{G_{\mathbf{B}_2}} = 0. \qquad (2.1.1.37)$$

Let $\mathfrak{R}_t[\psi]$ be operator such that

$$(\mathfrak{R}_{t,\varepsilon}[\psi])_\varepsilon = (A_{t,\varepsilon}[\psi_\varepsilon]\psi_\varepsilon(t))_\varepsilon. \qquad (2.1.1.38)$$

where $\forall (\psi_\varepsilon)_\varepsilon \left( (\psi_\varepsilon)_\varepsilon \in G_{\mathbf{C}_{\mathbf{B}_1}^t} \right)$ operator $(A_{t,\varepsilon}[\psi_\varepsilon])_\varepsilon : G_{\mathbf{B}_1} \to G_{\mathbf{B}_2}$ is continuous. Thus $(A_{t,}[\psi_\varepsilon])_\varepsilon$ is a infinitesimal generator of the one-parameter semigroup

$$\exp([\widetilde{\varepsilon} A_{t,\varepsilon}[\psi_\varepsilon]])_\varepsilon. \qquad (2.1.1.39)$$

Hence master operator of the Cauchy problem (2.1.1.32) is

$$(S_{t,\tilde{\varepsilon},\varepsilon}[v])_\varepsilon = \exp([\tilde{\varepsilon}A_{t,\varepsilon}[\psi_\varepsilon]]v_\varepsilon(t))_\varepsilon. \tag{2.1.1.40}$$

Corresponding $T$-exp is

$$(\psi_\varepsilon(t))_\varepsilon = \left(\prod_{\tau}^{t} \exp[A_{\tau,\varepsilon}[\psi_\varepsilon]]d\tau\right)_\varepsilon. \tag{2.1.1.41}$$

## II.1.2.Generalized Maslov theorem.

Let us consider Cauchy problem:

$$-i\hbar\frac{d}{dt}\Psi_k(t,x) - \hbar^2\Delta\Psi_k(t,x) + \sum_{l=1}^{m} V_{k,l}(x)\Psi_l(t,x) +$$

$$\sum_{l,j=1}^{m} \Psi_l(t,x)\int \rho_{k,l,j}(x,y)|\Psi_j(t,y)|^2 dy = 0, \tag{2.1.2.1}$$

$$\Psi_k(0,x) = \Psi_{0,k}(x),$$

$$k = 1,\ldots,m,$$

where

$$\Psi_0 = (\Psi_{0,1}, \ldots, \Psi_{0,m}) \in \mathcal{L}_{2,m}(\mathbb{R}^n),$$

$$V_{k,l}(x) = \overline{V_{l,k}(x)}, \rho_{k,l,j}(x,y) = \overline{\rho_{l,k,j}(x,y)}, \tag{2.1.2.2}$$

$$V_{k,l}(x) \in C(\mathbb{R}^n), \rho_{k,l,j}(x,y) \in C(\mathbb{R}^n \times \mathbb{R}^n).$$

We remind that $\mathcal{L}_{2,m}(\mathbb{R}^n) = \{\Psi(x) | \|\Psi\| < \infty\}$, where

$$\|\Psi\|^2 = \sum_{j=1}^{m} \int |\Psi_j(x)|^2 dx < \infty \tag{2.1.2.3}$$

and $W_2^{k,m}(\mathbb{R}^n) = \left\{\Psi(x) | \|\Psi\|_{W_2^{k,m}} < \infty\right\}$, where

$$\|\Psi\|_{W_2^{k,m}} = \sum_{j=1}^{m} \|\Psi_j\|_{W_2^k}. \tag{2.1.2.4}$$

Let $\Psi(t) = (\Psi_1(t), \ldots, \Psi_m(t))$ be a function $\Psi : \mathbb{R} \to \mathcal{L}_{2,m}(\mathbb{R}^n)$. We abbreviate $\Psi(t) \in \mathcal{L}_{2,m}(\mathbb{R}^n)$ iff $\forall t[\Psi(t) \in \mathcal{L}_{2,m}(\mathbb{R}^n)]$.
Let $v \in \mathcal{L}_{2,m}(\mathbb{R}^n)$ and let $A_{ij}[x,v]$ be the continuous matrix-function

$$A_{ij}[x,v] = V_{ij}(x) + \sum_{j=1}^{m} \int \rho_{k,l,j}(x,y)|v_l(y)|^2 dy. \tag{2.1.2.5}$$

Let $\widehat{A}[v]$ be a linear bounded operator on $\mathcal{L}_{2,m}(\mathbb{R}^n)$, given via formula

$$\left(\widehat{A}[v]u\right)_i = \sum_{j=1}^{m} A_{ij}[x,v]u_j. \tag{2.1.2.6}$$

Let's rewrite now Cauchy problem (2.1.2.1) in the form

$$-i\hbar \frac{d}{dt}\Psi(t) - \hbar^2 \Delta \Psi(t) + \widehat{A}[\Psi(t)]\Psi(t) = 0,$$

(2.1.2.7)

$$\Psi(0) = \Psi_0.$$

**Theorem 2.1.2.1.Maslov**.[20]. Cauchy problem (2.1.2.1) has a solution $\Psi(t)$ such that $\Psi(t) \in \mathcal{L}_{2,m}(\mathbb{R}^n)$ and the function $\Psi(t)$ is given via formula:

$$\Psi(t) = \left(\prod_{\tau=0}^{t} G_{d\tau}[\Psi(\tau)]\right)[\Psi_0],$$

(2.1.2.8)

where

$$G_{\overline{\varepsilon}}[v] = \exp\{i\overline{\varepsilon}\Delta\}\exp\left\{i\overline{\varepsilon}\widehat{A}[v]\right\},$$

(2.1.2.9)

$$v \in \mathcal{L}_{2,m}(\mathbb{R}^n).$$

We let: $\delta = \frac{t}{N}$, $N = 1,2,\ldots$ and

$$\Psi_{N,0} = \Psi_0,$$

$$\Psi_{N,k+1} = G_\delta \Psi_{N,k},$$

(2.1.2.10)

$$k = 0,1,2,\ldots,N-1,$$

where

$$G_\delta u(x) = \exp(i\delta\hbar\Delta)\exp(i\delta\hbar^{-1}A[u])u(x),$$

(2.1.2.11)

$$u(x) \in L_2(\mathbb{R}^n).$$

and we set below

$$\Psi_N \triangleq \Psi_{N,N}.$$

(2.1.2.12)

Suppose that $\Psi_0(x) \in W_2^{2,m}(\mathbb{R}^n)$. Then every solution $\Psi(t,x)$ of the Eq.(2.1.2.7) such that $\Psi(0,x) = \Psi_0$, has the representation:

$$\Psi(t,x) = s\text{-}\lim_{N\to\infty}\Psi_N,$$

(2.1.2.13)

where the strong limit in (2.1.2.13) exist uniformly on $[0,T]$ for any $T < \infty$ in the sence of the convergence in $L_{2,m}(\mathbb{R}^n)$ norm.
**Proof**.

**1**. We note that

$$\|\exp(i\tilde{\varepsilon}\hbar^{-1}A[v])\| = 1,$$

(2.1.2.14)

$$\|\exp(i\tilde{\varepsilon}\hbar\Delta)\| = 1,$$

where we denote

$$\|\cdot\| \triangleq \|\cdot\|_{L_{2,m}}. \tag{2.1.2.15}$$

Hence one obtain

$$\|\Psi_{N,k}\| = \|\Psi_0\|. \tag{2.1.2.16}$$

**2**. Let's now to prove the formula

$$\Psi_{N,k} = \exp(ik\delta h\Delta)\Psi_0 +$$

$$i\delta h^{-1} \sum_{j=0}^{k-1} \exp[i(k-j-1)\delta h\Delta] A[\Psi_{N,j}]\Psi_{N,j} + O\left(\frac{1}{hN}\right) \tag{2.1.2.17}$$

Note that $\|v\| \leqslant r \Rightarrow \|A[v]\| \leqslant C_r.$ Hence

$$\exp(i\delta h^{-1}A[v]) = 1 + i\delta h^{-1}A[\Psi_{N,j}] + O\left(\frac{1}{h^2N^2}\right) \tag{2.1.2.18}$$

From Eq.(2.1.2.18) by simple calculation one obtain

$$\Psi_{N,k} = \exp(ik\delta h\Delta)\Psi_{N,k-1} +$$

$$i\delta h^{-1} \exp[i\delta h\Delta] A[\Psi_{N,k-1}]\Psi_{N,k-1} + O\left(\frac{1}{h^2N^2}\right), \tag{2.1.2.19}$$

where

$$\left\| O\left(\frac{1}{\hbar^2 N^2}\right) \right\|_2 \leqslant \frac{C}{\hbar^2 N^2} \qquad (2.1.2.20)$$

and where constant $C$ does not depend on $k = 1, 2, \ldots, N.$ From Eq.(2.1.2.19) by simple calculation one obtain Eq.(2.1.2.17).

**3**. Let's now to prove the inequality:

$$\|\Psi_{N,l+k} - \Psi_{N,l}\| \leqslant$$

$$C_0[\exp(C_1 t)](k\hbar^{-1}\delta + \|\exp(ik\delta\hbar\Delta)\Psi_0 - \Psi_0\|), \qquad (2.1.2.21)$$

$$l + k \leqslant N,$$

where $C_0 = const$ and $C_1 = const$.

By using formula (2.1.2.17) one obtain

$$\Psi_{N,l+k} - \Psi_{N,l} = \exp(ik\delta\hbar\Delta)(\exp(ik\delta\hbar\Delta)\Psi_0 - \Psi_0) +$$

$$+ i\delta\hbar^{-1}\sum_{j=0}^{k-1}\exp[i\delta\hbar(k+l-j-1)\Delta]A\Psi_{N,j} + \qquad (2.1.2.22)$$

$$+ i\delta\hbar^{-1}\sum_{j=0}^{l-1}\exp[i\delta\hbar(l-j-1)\Delta](A\Psi_{N,j+k} - A\Psi_{N,j}) + O\left(\frac{1}{\hbar N}\right)$$

From inequality: $\|A[v]\| \leqslant C_r$ and Eqs.(2.1.2.22) we obtain

$$a_0 \leqslant \epsilon,$$

$$a_1 \leqslant C\left[\epsilon + k\delta h^{-1} + h^{-1}\sum_{j=0}^{l-1}\delta a_j\right],$$

$$l = 1, 2, \ldots, N, \qquad\qquad (2.1.2.23)$$

$$a_j = \|\Psi_{N,j+k} - \Psi_{N,j}\|,$$

$$\epsilon = \|\exp(ik\delta h\Delta)\Psi_0 - \Psi_0\|$$

**4**. Suppose that $M \in \mathbb{N}, M \geq 1$ and $M = const$ then for every $j \in \mathbb{N}$ one obtain

$$\frac{j}{M} - \left[\frac{j}{M}\right] \leq 1. \qquad\qquad (2.1.2.24)$$

Hence from (2.1.2.21) we obtain

$$\Psi_{MN,j} = \Psi_{M,\left[\frac{j}{M}\right]} + O\left(\frac{1}{hN}\right), \qquad\qquad (2.1.2.25)$$

where

$$\left\|O\left(\frac{1}{hN}\right)\right\| \leq$$

$$\left(\delta h^{-1} + \left\|\Psi_0 - \exp\left[i\delta h\left(\frac{j}{M} - \left[\frac{j}{M}\right]\right)\Delta\right]\Psi_0\right\|\right) \leq \qquad\qquad (2.1.2.26)$$

$$\delta h^{-1}\|\Psi_0\|_{W_2^{2m}(\mathbb{R}^n)}.$$

We denote

$$\gamma_N(\hbar^{-1}) \triangleq \delta\hbar^{-1} \|\Psi_0\|_{W_2^2(\mathbb{R}^3)},$$

(2.1.2.27)

$$\delta = \delta(N).$$

Thus we have two finite splittings $\Delta_N^{(1)} = \left\{0 = t_0 < t_1^{(1)} < \ldots < t_N^{(1)} = t\right\}$ and $\Delta_{NM}^{(2)} = \left\{0 = t_0 < t_1^{(2)} < \ldots < t_{NM}^{(2)} = t\right\}$ of the closed interval $[0, t]$ such that $\left\{t_i^{(1)}\right\}_{i=0}^{N} \subsetneqq \left\{t_j^{(2)}\right\}_{j=0}^{NM}$.

Let's rewrite formula (2.1.2.17) by using replacement $N \to NM$ for $k_0 = M \times l$. Thus we obtain

$$\Psi_{MN,Ml} = \exp(i\delta\hbar\Delta)\Psi_0 + i\delta\hbar^{-1}\sum_{j=0}^{l-1}\exp[i\delta\hbar(l-j-1)\Delta] \times$$

(2.1.2.28)

$$A\Psi_{MN,Mj} + O(\gamma_N(\hbar^{-1})).$$

Let's denote

$$b_j = \|\Psi_{N,j} - \Psi_{MN,Mj}\|$$

(2.1.2.29)

By using subtraction Eq.(2.1.2.17) from Eq.(2.1.2.29) we obtain

$$b_0 \leq \alpha_N$$

$$b_l \leq \alpha_N + c\delta\hbar^{-1}\sum_{j=0}^{l-1}b_j,$$

(2.1.2.30)

$$l = 1, \ldots, N,$$

where $\alpha_N \leq C_1 \gamma_N(\hbar^{-1})$. From (2.1.2.30) we obtain

$$b_l \leq \alpha_N \exp[C_2 t] \qquad (2.1.2.31)$$

In particulars for $l = N$ we obtain:

$$\| \Psi_N(t) - \Psi_{NM}(t) \| \leq \frac{t}{N} C_1 \exp(C_2 t) \| \Psi_0 \|_{W_2^{2,m}} \qquad (2.1.2.32)$$

and therefore sequence $\Psi_N(t)$ uniformly convergence on interval $[0, T]$ for any $T < \infty$ in the sence of the convergence in $L_{2,m}(\mathbb{R}^n)$ norm.

**5**. Any function $\Psi_N(t)$ is continuous on $[0, T]$. Therefore $\Psi(t) = s\text{-lim}_{N \to \infty} \Psi_N(t)$ is continuous on $[0, T]$. Hence the Riemann sum

$$S_N(t) = \sum_{k=0}^{N-1} \frac{t}{N} \exp\left\{ i \frac{t}{N}(N - k - 1)\Delta \right\} \widehat{A}[\Psi_{N,k}] \Psi_{N,k} \qquad (2.1.2.33)$$

is convergence in the norm of the $L_{2,m}(\mathbb{R}^n)$ to integral

$$\int_0^t d\tau \exp\{i(t - \tau)\} \widehat{A}[\Psi(\tau)] \Psi(\tau). \qquad (2.1.2.34)$$

if $N \to \infty$. Finally from Eq.(2.1.2.17) in the limit $N \to \infty$ we obtain integral equation

$$\Psi(t) = \exp\{it\Delta\} \Psi_0 + \int_0^t d\tau \exp\{i(t - \tau)\} \widehat{A}[\Psi(\tau)] \Psi(\tau). \qquad (2.1.2.35)$$

Eq.(2.1.2.35) completed the proof.

Let us consider now Cauchy problem:

$$-i\hbar\frac{d}{dt}(\Psi_{k,\varepsilon}(t,x))_\varepsilon - \hbar^2\Delta(\Psi_{k,\varepsilon}(t,x))_\varepsilon + \sum_{l=1}^{m}(V_{k,l,\varepsilon}(x)\Psi_{l,\varepsilon}(t,x))_\varepsilon +$$

$$\sum_{l,j=1}^{m}(\Psi_{l,\varepsilon}(t,x))_\varepsilon \int \Big[\Big(\rho_{k,l,j,\varepsilon}(x,y)|\Psi_{j,\varepsilon}(t,y)|^2\Big)_\varepsilon\Big]dy = 0, \qquad (2.1.2.36)$$

$$(\Psi_{k,\varepsilon}(0,x))_\varepsilon = (\Psi_{0,k,\varepsilon}(x))_\varepsilon,$$

$$k = 1,\ldots,m,$$

where

$$(\Psi_{0,\varepsilon})_\varepsilon = ((\Psi_{0,1,\varepsilon})_\varepsilon,\ldots,(\Psi_{0,m,\varepsilon})_\varepsilon) \in G_{\mathcal{L}_{2,m}(\mathbb{R}^n)},$$

$$(V_{k,l,\varepsilon}(x))_\varepsilon = \left(\overline{V_{l,k,\varepsilon}(x)}\right)_\varepsilon,$$

$$(\rho_{k,l,j,\varepsilon}(x,y))_\varepsilon = \left(\overline{\rho_{l,k,j,\varepsilon}(x,y)}\right)_\varepsilon, \qquad (2.1.2.37)$$

$$(V_{k,l,\varepsilon}(x))_\varepsilon \in G_{C(\mathbb{R}^n)},$$

$$(\rho_{k,l,j,\varepsilon}(x,y))_\varepsilon \in G_{C(\mathbb{R}^n\times\mathbb{R}^n)}.$$

Let $(\Psi_\varepsilon(t))_\varepsilon = ((\Psi_{1,\varepsilon}(t))_\varepsilon,\ldots,(\Psi_{m,\varepsilon}(t))_\varepsilon)$ be a function $(\Psi_\varepsilon)_\varepsilon : \mathbb{R} \to G_{\mathcal{L}_{2,m}(\mathbb{R}^n)}$. We abbreviate $(\Psi_\varepsilon(t))_\varepsilon \in G_{\mathcal{L}_{2,m}(\mathbb{R}^n)}$ iff $\forall t[\mathbf{cl}[(\Psi_\varepsilon(t))_\varepsilon] \in G_{\mathcal{L}_{2,m}(\mathbb{R}^n)}]$.
Let $(v_\varepsilon)_\varepsilon \in G_{\mathcal{L}_{2,m}(\mathbb{R}^n)}$ and let $(A_{ij,\varepsilon}[x,v_\varepsilon])_\varepsilon$ be the continuous matrix-function

$$(A_{ij,\varepsilon}[x,v_\varepsilon])_\varepsilon = (V_{ij,\varepsilon}(x))_\varepsilon + \sum_{j=1}^{m} \int \Big[ \big( \rho_{k,l,j,\varepsilon}(x,y)|v_{l,\varepsilon}(y)|^2 \big)_\varepsilon \Big] dy. \qquad (2.1.2.38)$$

Let $\left( \widehat{A}_\varepsilon[(v_\varepsilon)] \right)_\varepsilon$ be a linear bounded operator on $G_{\mathcal{L}_{2,m}(\mathbb{R}^n)}$, given via formula

$$\left( \left( \widehat{A}_\varepsilon[v_\varepsilon] u_\varepsilon \right)_i \right)_\varepsilon = \sum_{j=1}^{m} (A_{ij,\varepsilon}[x,v_\varepsilon] u_{j,\varepsilon})_\varepsilon. \qquad (2.1.2.39)$$

Let's rewrite now Cauchy problem (2.1.2.36) in the form:

$$-ih\frac{d}{dt}(\Psi_\varepsilon(t))_\varepsilon - h^2\Delta(\Psi_\varepsilon(t))_\varepsilon + \Big[ \left( \widehat{A}_\varepsilon[\Psi_\varepsilon(t)] \right)_\varepsilon \Big](\Psi_\varepsilon(t))_\varepsilon = 0,$$
$$(2.1.2.40)$$

$$(\Psi_\varepsilon(0))_\varepsilon = (\Psi_{0,\varepsilon})_\varepsilon.$$

**Theorem 2**.**1**.**2**.**2**.(**Generalized Maslov Theorem**).Generalized Cauchy problem

(2.1.2.40) has a solution $(\Psi_\varepsilon(t))_\varepsilon$ such that $(\Psi_\varepsilon(t))_\varepsilon \in G_{\mathcal{L}_{2,m}(\mathbb{R}^n)}$ and the generalized function $(\Psi_\varepsilon(t))_\varepsilon$ is given via formula:

$$(\Psi_\varepsilon(t))_\varepsilon = \left( \left( \prod_{\tau=0}^{t} G_{\tau,d\tau,\varepsilon}[\Psi_\varepsilon(\tau)] \right)[\Psi_{0,\varepsilon}] \right)_\varepsilon, \qquad (2.1.2.41)$$

where

$$G_{\widetilde{\varepsilon},\varepsilon}[v_\varepsilon] = \exp\{i\widetilde{\varepsilon}\Delta\}\exp\left\{ \left( i\widetilde{\varepsilon}\widehat{A}_\varepsilon[v_\varepsilon] \right)_\varepsilon \right\},$$
$$(2.1.2.42)$$

$$(v_\varepsilon)_\varepsilon \in G_{\mathcal{L}_{2,m}(\mathbb{R}^n)}.$$

We let: $\delta = \frac{t}{N}$, $N = 1, 2, \ldots$ and

$$(\Psi_{N,0,\varepsilon})_\varepsilon = (\Psi_{0,\varepsilon})_\varepsilon,$$

$$(\Psi_{N,k+1,\varepsilon})_\varepsilon = (G_{\delta,\varepsilon}\Psi_{N,k,\varepsilon})_\varepsilon, \qquad (2.1.2.43)$$

$$k = 0, 1, 2, \ldots, N - 1,$$

where

$$(G_{\delta,\varepsilon}u_\varepsilon(x))_\varepsilon = \exp(i\delta\hbar\Delta)\exp((i\delta\hbar^{-1}A_\varepsilon[u_\varepsilon])u_\varepsilon(x))_\varepsilon,$$

$$\qquad (2.1.2.44)$$

$$(u_\varepsilon(x))_\varepsilon \in G_{L_2(\mathbb{R}^n)}.$$

and we set below

$$(\Psi_{N,\varepsilon})_\varepsilon \triangleq (\Psi_{N,N,\varepsilon})_\varepsilon. \qquad (2.1.2.45)$$

Suppose that $\Psi_{0,\varepsilon}(x) \in G_{W_2^{2m}(\mathbb{R}^n)}$. Then every solution $(\Psi_\varepsilon(t,x))_\varepsilon$ of the Eq.(2.1.2.40) such that $(\Psi_\varepsilon(0,x))_\varepsilon = (\Psi_{0,\varepsilon}(x))_\varepsilon$ has the representation:

$$(\Psi_\varepsilon(t,x))_\varepsilon = s\text{-}\lim_{N \to \infty}\Psi_N, \qquad (2.1.2.46)$$

where the strong limit in (2.1.2.46) exist uniformly on $[0, T]$ for any $T < \infty$ in the sence of the convergence in ultra-pseudo-norm norm on $G_{L_{2,m}(\mathbb{R}^n)}$.
**Proof**.

**1**. We note that

$$\left\| \exp\left(\left(i\widetilde{\varepsilon}\hbar^{-1}A_\varepsilon[v_\varepsilon]\right)\right)_\varepsilon \right\|_\# = 1,$$

$$(2.1.2.47)$$

$$\left\| \exp(i\widetilde{\varepsilon}\hbar\Delta) \right\|_\# = 1,$$

where we denote

$$\|\cdot\|_\# \triangleq \|\cdot\|_{G_{L_{2,m}}}^{\mathbf{op}}.$$

$$(2.1.2.48)$$

Hence one obtain

$$\left\| (\Psi_{N,k,\varepsilon})_\varepsilon \right\|_{G_{L_{2,m}}} = \left\| (\Psi_{0,\varepsilon})_\varepsilon \right\|_{G_{L_{2,m}}}.$$

$$(2.1.2.49)$$

**2**. Let's now to prove the formula

$$(\Psi_{N,k,\varepsilon})_\varepsilon = \exp(ik\delta\hbar\Delta)(\Psi_{0,\varepsilon})_\varepsilon +$$

$$(2.1.2.50)$$

$$i\delta\hbar^{-1}\sum_{j=0}^{k-1}\exp[i(k-j-1)\delta\hbar\Delta](A_\varepsilon[\Psi_{N,j,\varepsilon}])_\varepsilon(\Psi_{N,j,\varepsilon})_\varepsilon + O\left(\frac{1}{\hbar N}\right)$$

Note that $\|v\|_{G_{L_{2,m}}} \leqslant r \Rightarrow \left\| (A_\varepsilon[v_\varepsilon])_\varepsilon \right\|_{G_{L_{2,m}}} \leqslant C_r.$ Hence

$$\exp(i\delta\hbar^{-1}(A_\varepsilon[v_\varepsilon])_\varepsilon) = 1 + i\delta\hbar^{-1}(A_\varepsilon[\Psi_{N,j,\varepsilon}])_\varepsilon + O\left(\frac{1}{\hbar^2 N^2}\right) \qquad (2.1.2.51)$$

From Eq.(2.1.2.51) by simple calculation one obtain

$$(\Psi_{N,k,\varepsilon})_\varepsilon = \exp(ik\delta\hbar\Delta)(\Psi_{N,k-1,\varepsilon})_\varepsilon +$$

$$i\delta\hbar^{-1}\exp[i\delta\hbar\Delta](A_\varepsilon[\Psi_{N,k-1,\varepsilon}])_\varepsilon(\Psi_{N,k-1,\varepsilon})_\varepsilon + O\left(\frac{1}{\hbar^2 N^2}\right), \quad (2.1.2.52)$$

where

$$\left\| O\left(\frac{1}{\hbar^2 N^2}\right) \right\|_{G_{L_{2,m}}} \leqslant \frac{C}{\hbar^2 N^2} \qquad (2.1.2.53)$$

and where constant $C$ does not depend on $k = 1, 2, \ldots, N.$ From Eq.(2.1.2.52) by simple calculation one obtain Eq.(2.1.2.50).

**3**. Let's now to prove the inequality:

$$\left\| (\Psi_{N,l+k,\varepsilon})_\epsilon - (\Psi_{N,l,\varepsilon})_\varepsilon \right\|_{G_{L_{2,m}}} \leqslant$$

$$C_0[\exp(C_1 t)]\left( kh^{-1}\delta + \left\| \exp(ik\delta\hbar\Delta)(\Psi_{0,\varepsilon})_\varepsilon - (\Psi_{0,\varepsilon})_\varepsilon \right\|_{G_{L_{2,m}}} \right), \qquad (2.1.2.54)$$

$$l + k \leqslant N,$$

where $C_0 = const$ and $C_1 = const.$

By using formula (2.1.2.50) one obtain

$$(\Psi_{N,l+k,\varepsilon})_\varepsilon - (\Psi_{N,l,\varepsilon})_\varepsilon = \exp(ik\delta\hbar\Delta)(\exp(ik\delta\hbar\Delta)(\Psi_{0,\varepsilon})_\varepsilon - (\Psi_{0,\varepsilon})_\varepsilon) +$$

$$+i\delta\hbar^{-1}\sum_{j=0}^{k-1}\exp[i\delta\hbar(k+l-j-1)\Delta](A_\varepsilon\Psi_{N,j,\varepsilon})_\varepsilon + \qquad (2.1.2.55)$$

$$+i\delta\hbar^{-1}\sum_{j=0}^{l-1}\exp[i\delta\hbar(l-j-1)\Delta]((A_\varepsilon\Psi_{N,j+k,\varepsilon})_\varepsilon - (A_\varepsilon\Psi_{N,j,\varepsilon})_\varepsilon) + O\left(\frac{1}{\hbar N}\right)$$

From inequality: $\|(A_\varepsilon[v_\varepsilon])_\varepsilon\|_\# \leqslant C_r$ and Eqs.(2.1.2.55) we obtain

$$a_0 \leqslant \epsilon,$$

$$a_1 \leqslant C\Big[\epsilon + k\delta\hbar^{-1} + \hbar^{-1}\sum_{j=0}^{l-1}\delta a_j\Big],$$

$$l = 1, 2, \ldots, N, \qquad (2.1.2.56)$$

$$a_j = \|(\Psi_{N,j+k,\varepsilon})_\varepsilon - (\Psi_{N,j,\varepsilon})_\varepsilon\|_{G_{L_{2,m}}},$$

$$\epsilon = \|\exp(ik\delta\hbar\Delta)(\Psi_{0,\varepsilon})_\varepsilon - (\Psi_{0,\varepsilon})_\varepsilon\|_{G_{L_{2,m}}}$$

**4**. Suppose that $M \in \mathbb{N}, M \geq 1$ and $M = const$ then for every $j \in \mathbb{N}$ one obtain

$$\frac{j}{M} - \left[\frac{j}{M}\right] \leq 1. \qquad (2.1.2.57)$$

Hence from (2.1.2.21) we obtain

$$\left(\Psi_{MN,j,\varepsilon}\right)_\varepsilon = \left(\Psi_{M,\left[\frac{j}{M}\right],\varepsilon}\right)_\varepsilon + O\left(\frac{1}{\hbar N}\right), \qquad (2.1.2.58)$$

where

$$\left\|O\left(\frac{1}{\hbar N}\right)\right\|_{G_{L_{2,m}}} \leq$$

$$\left(\delta\hbar^{-1} + \left\|\left(\Psi_{0,\varepsilon}\right)_\varepsilon - \exp\left[i\delta\hbar\left(\frac{j}{M} - \left[\frac{j}{M}\right]\right)\Delta\right]\left(\Psi_{0,\varepsilon}\right)_\varepsilon\right\|_{G_{L_{2,m}}}\right) \leq \qquad (2.1.2.59)$$

$$\delta\hbar^{-1}\left\|\left(\Psi_{0,\varepsilon}\right)_\varepsilon\right\|_{G_{W_2^{2,m}(\mathbb{R}^n)}}.$$

We denote

$$\gamma_N(\hbar^{-1}) \triangleq \delta\hbar^{-1}\left\|\left(\Psi_{0,\varepsilon}\right)_\varepsilon\right\|_{G_{W_2^{2,m}(\mathbb{R}^n)}},$$

$$\qquad (2.1.2.60)$$

$$\delta = \delta(N).$$

Thus we have two finite splittings $\Delta_N^{(1)} = \left\{0 = t_0 < t_1^{(1)} < \ldots < t_N^{(1)} = t\right\}$ and $\Delta_{NM}^{(2)} = \left\{0 = t_1^{(2)} < t_1^{(2)} < \ldots < t_{NM}^{(2)} = t\right\}$ of the closed interval $[0, t]$ such that $\left\{t_i^{(1)}\right\}_{i=0}^{N} \subsetneqq \left\{t_j^{(2)}\right\}_{j=0}^{NM}$.

Let's rewrite now formula (2.1.2.50) by using replacement $N \to NM$ for $k_0 = M \times l$. Thus we obtain

$$(\Psi_{MN,Ml,\varepsilon})_\varepsilon = \exp(i\delta h\Delta)(\Psi_{0,\varepsilon})_\varepsilon + i\delta h^{-1}\sum_{j=0}^{l-1}\exp[i\delta h(l-j-1)\Delta]\times$$

(2.1.2.61)

$$(A_\varepsilon\Psi_{MN,Mj,\varepsilon})_\varepsilon + O(\gamma_N(h^{-1})).$$

Let's denote

$$b_j = \|(\Psi_{N,j,\varepsilon})_\varepsilon - (\Psi_{MN,Mj,\varepsilon})_\varepsilon\|_{G_{L_{2,m}}}$$

(2.1.2.62)

By using subtraction Eq.(2.1.2.50) from Eq.(2.1.2.61) we obtain

$$b_0 \le \alpha_N$$

$$b_l \le \alpha_N + c\delta h^{-1}\sum_{j=0}^{l-1}b_j,$$

(2.1.2.63)

$$l = 1,\ldots,N,$$

where $\alpha_N \le C_1\gamma_N(h^{-1})$. From (2.1.2.63) we obtain

$$b_l \le \alpha_N\exp[C_2 t].$$

(2.1.2.64)

In particulars for $l = N$ we obtain:

$$\|(\Psi_{N,\varepsilon}(t))_\varepsilon - (\Psi_{NM,\varepsilon}(t))_\varepsilon\|_{G_{L_{2,m}}} \le \frac{t}{N}C_1\exp(C_2 t)\|(\Psi_{0,\varepsilon})_\varepsilon\|_{G_{W_2^{2,m}}}$$

(2.1.2.65)

and therefore sequence $(\Psi_{N,\varepsilon}(t))_\varepsilon$ uniformly convergence on interval $[0,T]$ for any $T < \infty$ in the sence of the convergence in ultra-pseudo-seminorm on $G_{L_{2,m}(\mathbb{R}^n)}$.

**5**. Any function $\Psi_N(t)$ is continuous on $[0,T]$. Therefore $\Psi(t) = s\text{-}\lim_{N\to\infty}(\Psi_{N,\varepsilon}(t))_\varepsilon$ is continuous on $[0,T]$. Hence the Riemann sum

$$(S_N,\varepsilon(t))_\varepsilon = \sum_{k=0}^{N-1} \frac{t}{N} \exp\left\{ i\frac{t}{N}(N-k-1)\Delta \right\} \left( \widehat{A}_\varepsilon[\Psi_{N,k,\varepsilon}] \right)_\varepsilon (\Psi_{N,k,\varepsilon})_\varepsilon \qquad (2.1.2.66)$$

is convergence in the ultra-pseudo-seminorm of the $G_{L_{2,m}(\mathbb{R}^n)}$ to integral

$$\int_0^t d\tau \exp\{i(t-\tau)\} \left( \widehat{A}_\varepsilon[\Psi_\varepsilon(\tau)] \right)_\varepsilon (\Psi_\varepsilon(\tau))_\varepsilon. \qquad (2.1.2.67)$$

Finally from Eq.(**II**.1.2.50) in the limit $N \to \infty$ we obtain integral equation

$$(\Psi_\varepsilon(t))_\varepsilon = \exp\{it\Delta\}(\Psi_{0,\varepsilon})_\varepsilon + \int_0^t d\tau \exp\{i(t-\tau)\} \left[ \left( \widehat{A}_\varepsilon[\Psi_\varepsilon(\tau)] \right)_\varepsilon \right] (\Psi_\varepsilon(\tau))_\varepsilon. \quad (\mathbf{II}.1.2.68)$$

Eq.(**II**.1.2.68) completed the proof.

# II.1.3. Representation of the Maslov-Colombeau T-exponent by using Feynman-Colombeau path Integral.

**Definition II**.1.3.1. We define formal pseudodifferential operator with symbol $A(x,p)$ via formula

$$A\left( \overset{2}{x}, -ih\, \overset{1}{\frac{\partial}{\partial x}} \right) v(x) = \frac{1}{2\pi h} \int dp A(x,p) \int dy v(y) \exp\left( -\frac{i}{h} yp \right). \qquad (\mathbf{II}.1.3.1)$$

Let's consider formal $T$-exponent with a generator

$$G_{t,\overline{z}} = G_{t,\overline{z}}\left(\overset{2}{x}, -ih\ \overset{1}{\frac{\partial}{\partial x}}, t, [v]\right)v(x,t). \tag{2.1.3.2}$$

Let $\Delta_N = \{0 = t_0 < t_1 < \ldots < t_{N+1} = t\}, \Delta t_j = t_{j+1} - t_j, \delta_N = \max_j \Delta t_j, x \in \mathbb{R}^n$.

$T$-exponent

$$\Psi(x,t) = \left(\prod_{\tau=0}^{t} G_{\tau,d\tau}\right)[\Psi_0] \tag{2.1.3.3}$$

is defined via formulae

$$\Psi(x,t) = s\text{-}\lim_{\delta_N \to 0} \Psi_{N+1}(x,t),$$

$$\Psi_0(x,t) = \Psi_0,$$

$$\Psi_{k+1}(x,t) = \Psi_k(x,t) \text{ iff } 0 \le t \le t_k, \tag{2.1.3.4}$$

$$\Psi_{k+1}(x,t) = G_{t_k, t-t_k}\left(\overset{2}{x}, -ih\ \overset{1}{\frac{\partial}{\partial x}}, t, [\Psi_k]\right)\Psi_k(x,t) \text{ iff } t_k \le t \le t_{k+1}$$

$$k = 0, 1, \ldots, N$$

**Theorem 2.1.3.1.**[20]. Suppose that $T$-exponent given via Eq.(2.1.3.3) exist and $\Psi(x,t) \in W_2^k(\mathbb{R}^n)$. Then function $\Psi(x,t)$ has the representation

$$\Psi(x,t) = s\text{-} \lim_{\delta_N \to 0} \frac{1}{(2\pi h)^{(N+1)n}} \int \ldots \int \exp\left[\frac{i}{h}(x_{k+1} - x_k)\xi_k\right] \times$$

$$\times \prod_{k=0}^{N} G_{t_k, \Delta t_k}(x_{k+1}, \xi_k; [\Psi])\Psi_0(x_0)d^{N+1}xd^{N+1}\xi, \tag{2.1.3.5}$$

where $x_{N+1} = x$ and $d^{N+1} = dx_0 \ldots dx_N, d^{N+1}\xi = d\xi_0 \ldots d\xi_N, x_j \in \mathbb{R}^n, \xi_j \in \mathbb{R}^n$.
We abbreviate Eq.(2.1.3.5) in canonical symbolic form

$$\Psi(x,t) = \int_{q(t)=x} D[q(\tau)] \int D[p(\tau)]\Psi_0(q(0))\exp\left\{\frac{i}{h}\int_0^t p(\tau)\dot{q}(\tau)d\tau\right\} \times$$

$$\times \prod_{\tau=0}^{t} G_{\tau, d\tau}(q(\tau), p(\tau); [\Psi]). \tag{2.1.3.6}$$

From Theorem 2.1.3.1 and Theorem 2.1.2.1 one obtain directly:
**Theorem 2.1.3.2.**[20]. Suppose that $V(x) \in C(\mathbb{R}^n), a(x,y) \in C(\mathbb{R}^n \times \mathbb{R}^n)$ and $\Psi_0(x) \in W_2^2(\mathbb{R}^n)$. Let $\Psi(x,t)$ be the solution of the Cauchy problem (2.1.2.7). Then $\Psi(x,t)$ can be expressed via formula:

$$\Psi(x,t) = \lim_{N \to \infty} (4\pi i h)^{-\frac{(N+1)}{2}} \int \ldots \int \Psi_0(x_0) \frac{dx_0 \ldots dx_N}{\sqrt{\Delta t_0 \ldots \Delta t_N}} \times$$

$$\exp\left\{\frac{i}{h}\sum_{k=0}^{N}\left[\frac{1}{4}\left(\frac{x_{k+1} - x_k}{\Delta t_k}\right)^2 - V(x_{k+1}) - \int a(x_k, y)\Psi(y, t_k)dy\right]\Delta t_k\right\}, \tag{2.1.3.6$'$}$$

$$x_{N+1} = x.$$

Asuume that $a(x,y) = 0$. Then for the linear Schrödinger equation (0.1) we obtain

$$\Psi(x,t) = \lim_{N\to\infty} (4\pi i\hbar)^{-\frac{(N+1)}{2}} \int \ldots \int \Psi_0(x_0) \frac{dx_0 \ldots dx_N}{\sqrt{\Delta t_0 \ldots \Delta t_N}} \times$$

$$\exp\left\{ \frac{i}{\hbar} \sum_{k=0}^{N} \left[ \frac{1}{4}\left( \frac{x_{k+1}-x_k}{\Delta t_k} \right)^2 - V(x_{k+1}) - \int a(x_k,y)\Psi(y,t_k)dy \right] \Delta t_k \right\}, \qquad (2.1.3.6^{\#})$$

$$x_{N+1} = x.$$

We abbreviate Eq.(2.1.3.6$^{\#}$) in canonical symbolic form

$$\Psi(x,t) = \int_{q(t)=x} D[q(t)]\Psi_0(q(0)) \exp\left\{ \frac{i}{\hbar} \int_0^t \left[ \frac{1}{4}\dot{q}^2(\tau) - V(q(\tau)) \right] d\tau \right\} =$$

$$\int_{\substack{q(t)=x \\ q(0)=y}} D[q(t)]dy\Psi_0(y) \exp\left\{ \frac{i}{\hbar} \int_0^t \left[ \frac{1}{4}\dot{q}^2(\tau) - V(q(\tau)) \right] d\tau \right\}. \qquad (2.1.3.6^{\#})$$

**Definition 2.1.3.2**. Assume that $v(x) \in G(\mathbb{R}^n), v(x) = \mathbf{cl}[(v_\varepsilon(x))_\varepsilon], A(x,p) \in G(\mathbb{R}^n \times \mathbb{R}^n)$,

$A(x,p) = \mathbf{cl}[(A_\varepsilon(x,p))_\varepsilon]$. We define formal Colombeau pseudodifferential operator with symbol $A(x,p)$ via formula:

$$\left( A_\varepsilon\left( \overset{2}{x}, -ih\ \overset{1}{\frac{\partial}{\partial x}} \right) v_\varepsilon(x) \right)_{\varepsilon \in (0,1]} =$$

$$\frac{1}{2\pi h} \left( \int dp A_\varepsilon(x,p) \int dy v_\varepsilon(y) \exp\left( -\frac{i}{h}yp \right) \right)_{\varepsilon \in (0,1]}. \qquad (2.1.3.7)$$

Let's consider formal Maslov-Colombeau $T$-exponent with a generator

$$\left( G_{t,\overline{\varepsilon},\varepsilon} \right)_{\varepsilon \in (0,1]} = \left( G_{t,\overline{\varepsilon},\varepsilon} \left( \overset{2}{x}, -ih \; \overset{1}{\frac{\partial}{\partial x}}, t, [v_\varepsilon] \right) v_\varepsilon(x,t) \right)_{\varepsilon \in (0,1]}. \tag{2.1.3.8}$$

Let $\Delta_N = \{ 0 = t_0 < t_1 < \ldots < t_{N+1} = t \}, \Delta t_j = t_{j+1} - t_j, \delta_N = \max_j \Delta t_j, x \in \mathbb{R}^n.$
Maslov-Colombeau $T$-exponent

$$\left( \Psi_\varepsilon(x,t) \right)_{\varepsilon \in (0,1]} = \left( \left( \prod_{\tau=0}^{t} G_{\tau, d\tau, \varepsilon} \right) [\Psi_{0,\varepsilon}] \right)_{\varepsilon \in (0,1]} \tag{2.1.3.9}$$

is defined via formulae

$$\left( \Psi_\varepsilon(x,t) \right)_\varepsilon = s\text{-} \lim_{\delta_N \to 0} \left( \Psi_{N+1,\varepsilon}(x,t) \right)_\varepsilon,$$

$$\left( \Psi_{0,\varepsilon}(x,t) \right)_\varepsilon = \left( \Psi_{0,\varepsilon}(x) \right)_\varepsilon,$$

$$\left( \Psi_{k+1,\varepsilon}(x,t) \right)_\varepsilon = \left( \Psi_{k,\varepsilon}(x,t) \right)_\varepsilon \text{ iff } 0 \leq t \leq t_k, \tag{2.1.3.10}$$

$$\left( \Psi_{k+1,\varepsilon}(x,t) \right)_\varepsilon = \left( G_{t_k, t-t_k, \varepsilon} \left( \overset{2}{x}, -ih \; \overset{1}{\frac{\partial}{\partial x}}, t, [\Psi_{k,\varepsilon}] \right) \Psi_{k,\varepsilon}(x,t) \right)_\varepsilon \text{ iff } t_k \leq t \leq t_{k+1}$$

$$k = 0, 1, \ldots, N$$

**Theorem 2.1.3.3**. Suppose that $T$-exponent given via Eq.(2.1.3.10) exist and

$\mathbf{cl}[(\Psi_\varepsilon(x,t))_\varepsilon] \in G_{W_2^k(\mathbb{R}^n)}$. Then Colombeau generalized function $(\Psi_\varepsilon(x,t))_\varepsilon$ has the representation

$$(\Psi_\varepsilon(x,t))_\varepsilon = s\text{-}\lim_{\delta_N \to 0}\left(\frac{1}{(2\pi h)^{(N+1)n}}\int\ldots\int\exp\left[\frac{i}{h}(x_{k+1}-x_k)\xi_k\right]\times\right.$$

$$\left.\times\prod_{k=0}^{N}G_{t_k,\Delta t_k,\varepsilon}(x_{k+1},\xi_k;[\Psi_\varepsilon])\Psi_{0,\varepsilon}(x_0)d^{N+1}xd^{N+1}\xi\right)_\varepsilon,$$

(2.1.3.11)

We rewrite Eq.(2.1.3.11) in symbolic form

$$(\Psi_\varepsilon(x,t))_\varepsilon = \left(\int_{q(t)=x}D[q(\tau)]\int D[p(\tau)]\Psi_0(q(0))\exp\left\{\frac{i}{h}\int_0^t p(\tau)\dot{q}(\tau)d\tau\right\}\times\right.$$

$$\left.\times\prod_{\tau=0}^{t}G_{\tau,d\tau}(q(\tau),p(\tau);[\Psi],\varepsilon)\right)_\varepsilon.$$

(2.1.3.12)

From Theorem 2.1.2.2 and Theorem 2.1.3.3 one obtain directly:

**Theorem 2.1.3.4.** Suppose that $V(x)\in G_{C(\mathbb{R}^n)}, a_\varepsilon(x,y)\in G_{C(\mathbb{R}^n\times\mathbb{R}^n)}$ and $\Psi_0(x)\in G_{W_2^2(\mathbb{R}^n)}$. Let $\Psi(x,t)$ be the solution of the Cauchy problem (2.1.2.7). Then $(\Psi_\varepsilon(x,t))_\varepsilon$ can be expressed via formula:

$$(\Psi_\varepsilon(x,t))_\varepsilon = \lim_{N\to\infty}\left((4\pi ih)^{-\frac{(N+1)}{2}}\int\ldots\int\Psi_0(x_0)\frac{dx_0\ldots dx_N}{\sqrt{\Delta t_0\ldots\Delta t_N}}\times\right.$$

(2.1.3.13)

$$\left.\exp\left\{\frac{i}{h}\sum_{k=0}^{N}\left[\frac{1}{4}\left(\frac{x_{k+1}-x_k}{\Delta t_k}\right)^2-V_\varepsilon(x_{k+1})-\int a_\varepsilon(x_{k,\varepsilon},y)\Psi(y,t_k)dy\right]\Delta t_k\right\}\right)_\varepsilon,$$

# II.2.Master Equation for the Feynman-Colombeau path integral.

## II.2.1. Master Equation for the Feynman-Colombeau path integral corresponding with a Schrödinger equation with non-Hermitian Hamiltonian.

The evolution of a closed quantum system during a time interval $T$ is described by the evolution operator $(\mathbf{U}_\varepsilon(T))_\varepsilon$. The matrix element of the operator $(\mathbf{U}_\varepsilon(T))_\varepsilon$ between the states with definite positions is called the propagator and may be expressed in the form of the corresponding Feynman-Colombeau path integral (see

Apendix I.3):

$$(\mathbf{U}_\varepsilon(x, T|x_0, 0))_{\varepsilon \in (0,1]} =$$

$$\left( \int_{q(T)=x} D_\mathbb{C}[q(t), \varepsilon] \Psi_{0,\varepsilon}(q(0)) \times \right.$$

$$\left. \exp\left[ \frac{i}{\hbar_\varepsilon} \left( \int_0^T dt \left[ \frac{m}{2} \left( \frac{d_{\widetilde{\mathbb{R}},\varepsilon} q(t)}{d_{\widetilde{\mathbb{R}},\varepsilon} t} \right)^2 - V_\varepsilon(q(t)) \right] \right) \right] \right)_{\varepsilon \in (0,1]},$$

$$(\Psi_{0,\varepsilon}(x))_\varepsilon = \delta(x - x_0).$$

(2.1.1)

If the system with the same Hamiltonian undergoes a continuous measurement and therefore is considered as being open, interacting with a measuring device or environment, its evolution may be described by the set of partial evolution operators $(\mathbf{U}_{T,\varepsilon}^\alpha)_\varepsilon$ depending on the output of the measurement:

$$(|\Psi_\varepsilon(T,\alpha)\rangle\rangle)_\varepsilon = (\mathbf{U}_{T,\varepsilon}^\alpha|\Psi_0\rangle)_\varepsilon \qquad (2.1.2)$$

The partial propagators are expressed by restricted path integrals. This means [32]-[35] that the path integral for $(\mathbf{U}_\varepsilon(T))_\varepsilon$ must be given via formula (2.1.1), but restricted according to the information given by the measurement readout $\alpha$. The information given by may be described by a weight functional $(w_{\alpha,\varepsilon}[q])_\varepsilon$ (positive, with values between 0 and 1) so that the partial propagator has to be written as a weighted Feynman-Colombeau path integral:

$$(\mathbf{U}_{\alpha,\varepsilon}(x,T|x_0,0))_{\varepsilon\in(0,1]} =$$

$$\left( \int_{q(T)=x} D_\mathbb{C}[q(t),\varepsilon]w_{\alpha,\varepsilon}[q]\Psi_{0,\varepsilon}(q(0)) \times \right.$$

$$\left. \exp\left[ \frac{i}{\hbar_\varepsilon}\left( \int_0^T dt\left[ \frac{m}{2}\left( \frac{d_{\widetilde{\mathbb{R}},\varepsilon}q(t)}{d_{\widetilde{\mathbb{R}},\varepsilon}t} \right)^2 - V_\varepsilon(q(t)) \right] \right) \right] \right)_{\varepsilon\in(0,1]}, \qquad (2.1.3)$$

$$V_\varepsilon(q(t)) = V(q_\varepsilon(t)),$$

$$q_\varepsilon(t) = \frac{q(t)}{1+\varepsilon^k q^m(t)}, k \geq 1, m \geq 2,$$

$$(\Psi_{0,\varepsilon}(x))_\varepsilon = \delta(x-x_0).$$

If the readout is unknown (nonselective description), the evolution of the measured system may be presented by the generalized density matrix:

$$(\rho_{T,\varepsilon})_\varepsilon = \int d\mu(\alpha)(\rho_{T,\varepsilon}^\alpha)_\varepsilon = \int d\mu(\alpha)[(\mathbf{U}_{\alpha,\varepsilon}(T))_\varepsilon][(\rho_{0,\varepsilon})_\varepsilon]\left[ \left( (\rho_{T,\varepsilon}^\alpha)^\dagger \right)_\varepsilon \right] \qquad (2.1.4)$$

and the generalized unitarity condition

$$\int d\mu(\alpha) \Big[ \big( (\mathbf{U}_{\alpha,\varepsilon}(T))^{\dagger} \big)_{\varepsilon} \Big] \big[ (\mathbf{U}_{\alpha,\varepsilon}(T))_{\varepsilon} \big] = 1 \qquad (2.1.5)$$

provides conservation of probabilities.

Let us consider the special case, when the weight functional $(w_{\alpha,\beta,\varepsilon}[q(t),q(0),q(T),\lambda])_{\varepsilon} = (w_{\alpha,\beta,\varepsilon}[q(t),q_0,q_T,\lambda])_{\varepsilon}$ describing an measurement may be taken of the symple Gaussian form such that:

$$(w_{\alpha,\beta,\varepsilon}[q])_{\varepsilon} =$$

$$\left( \exp\left[ -\nu \int_0^T \big[ q(t) - \alpha(t) \big]^2 dt - \frac{\upsilon}{\hbar_{\varepsilon}} \int_0^T \big[ q(t) - \beta(t,q_0,q_T,\lambda) \big]^2 dt - \right.\right.$$

$$\left.\left. \frac{1}{\hbar_{\varepsilon}} \gamma(q_0,q_T,\lambda) \right] \right)_{\varepsilon} =$$

$$\left( \exp\left[ -\frac{1}{\hbar_{\varepsilon}} \int_0^T \Big\{ \hbar_{\varepsilon} \nu [q(t) - \alpha(t)]^2 + \upsilon [q(t) - \beta(t,q_0,q_T,\lambda)]^2 \Big\} dt - \right.\right. \qquad (2.1.6)$$

$$\left.\left. \frac{1}{\hbar_{\varepsilon}} \gamma(q_0,q_T,\lambda) \right] \right)_{\varepsilon} =$$

$$\left( \exp\left[ -\frac{1}{\hbar_{\varepsilon}} \int_0^T R_{\alpha,\beta,\varepsilon}(q(t),\nu,\upsilon,t,q_0,q_T,\lambda) dt \right] \right)_{\varepsilon},$$

$$\gamma(q_0,q_T,\lambda) = \upsilon[q_0 - c - \lambda]^2 + \upsilon[q_T - c - \lambda]^2,$$

$$\nu,\upsilon \in \mathbb{R}_+, \lambda, c \in \mathbb{R}.$$

The corresponding Feynman-Colombeau path integral is

$$(\mathbf{U}_{\alpha,\varepsilon}(x,T|x_0,0))_{\varepsilon\in(0,1]} =$$

$$\left( \int\limits_{q(T)=x} D_{\mathbb{C}}[q(t),\varepsilon]\Psi_{0,\varepsilon}(q(0)) \times \right.$$

$$\exp\left[ \frac{i}{\hbar_\varepsilon}\left( \int\limits_0^T dt\left[ \frac{m}{2}\dot{q}^2_{\widetilde{\mathbb{R}},\varepsilon}(t) - V(q_\varepsilon(t)) \right] \right) - \right.$$

$$\left. \left. -\frac{1}{\hbar_\varepsilon}\int\limits_0^T R_{\alpha,\beta,\varepsilon}(q(t),v,v,t,q_0,q_T,\lambda)dt \right] \right)_\varepsilon =$$

(2.1.7)

$$\left( \int\limits_{q(T)=x} D_{\mathbb{C}}[q(t),\varepsilon]\Psi_{0,\varepsilon}(q(0)) \times \right.$$

$$\left. \exp\left[ \frac{i}{\hbar_\varepsilon}\left( \int\limits_0^T dt\left[ \frac{m}{2}\dot{q}^2_{\widetilde{\mathbb{R}},\varepsilon}(t) - V(q_\varepsilon(t)) + iR_{\alpha,\beta,\varepsilon}(q(t),v,v,t,q_0,q_T,\lambda) \right] \right) \right] \right)_\varepsilon$$

$$q_\varepsilon(t) = \frac{q(t)}{1+\varepsilon^k q^m(t)}, k\geq 1, m\geq 2,$$

$$(\Psi^2_{0,\varepsilon}(x-x_0))_\varepsilon = \delta(x-x_0).$$

Path integral (2.1.7) has the form of a canonical non-weighted path integral but with a non-Hermitian Hamiltonian

$$(\mathbf{H}_{\alpha,\beta,\varepsilon}(p,q,t))_\varepsilon = (\mathbf{H}_\varepsilon(p,q,t))_\varepsilon - \frac{1}{2i\hbar_\varepsilon m_2}p^2 - iR_{\alpha,\beta,\varepsilon}(q,v,v,t,q_0,q_T,\lambda) \qquad (2.1.8)$$

instead of the original Hamiltonian $(\mathbf{H}_\varepsilon)_\varepsilon$. Therefore, instead of calculating a restricted path integral one may solve the Colombeau-Schrödinger equation with a

non- Hermitian Hamiltonian:

$$\frac{\partial}{\partial t}\left(\Psi_{\alpha,\beta,\varepsilon}(t)\right)_{\varepsilon} = \left(-\frac{i}{(\hbar_{\varepsilon})_{\varepsilon}}(\mathbf{H}_{\varepsilon})_{\varepsilon} - R_{\alpha,\beta,\varepsilon}(q,v,\upsilon,t,q_0,q_T,\lambda,c)\right)\left(\Psi_{\alpha,\beta,\varepsilon}(t)\right)_{\varepsilon},$$

$$\left(\Psi_{\alpha,\varepsilon}(0)\right)_{\varepsilon} = \left(\Psi_{\varepsilon}(x-x_0)\right)_{\varepsilon}.$$

$$(2.1.9)$$

Let us consider now quantum average $(\langle T, x_0, \alpha, \beta; \hbar_{\varepsilon}, \varepsilon\rangle)_{\varepsilon}$ :

$$(\langle T, x_0, \alpha; \hbar, \varepsilon\rangle)_{\varepsilon} = \int dx x \left(\left(\mathbf{U}_{\alpha,,\beta,\varepsilon}(x,T|x_0,0)\mathbf{U}^{\dagger}_{\alpha,\varepsilon}(x,T|x_0,0)\right)\right)_{\varepsilon} =$$

$$Z^{-1}\int dx x \left(|\mathbf{U}_{\alpha,,\beta,\varepsilon}(x,T|x_0,0)|^2\right)_{\varepsilon}$$

$$(2.1.10)$$

for the case: (i) $v = 0$. Thus we set $v = 0$ in Eq.(2.1.7). In the sequel we denote for short

$$R_{\alpha,\beta,\varepsilon}(q,0,\upsilon,t,q_0,q_T,\lambda) \triangleq Q_{\beta,\varepsilon}(q,\upsilon,t,q_0,q_T,\lambda) \triangleq Q_{\beta,\varepsilon}(q(t),\upsilon) \qquad (2.1.11)$$

From Eq.(2.1.7) and Eq.(2.1.11) $\forall \varepsilon \in (0,1]$ we obtain

$$\langle T, x_0, \alpha, \beta; \hbar_\varepsilon, \varepsilon \rangle =$$

$$\int dx \left| \int_{q(T)=x} |D_\mathbb{C}|[q(t), \varepsilon] \Psi_{0,\varepsilon}(q(0)) \times \right.$$

$$\left. \exp\left[ \frac{i}{\hbar_\varepsilon} \left( \int_0^T dt \left[ \frac{m}{2} \dot{q}^2_{\widetilde{\mathbb{R}}, \varepsilon}(t) - V(q_\varepsilon(t)) \right] \right) \right] \right|^2 =$$

$$\int dx \left\{ \left[ \int_{q(T)=x} D^+[q(t), \varepsilon] \Psi_0(q(0)) \exp\left( -\frac{1}{\hbar_\varepsilon} \int_0^T Q_{\beta,\varepsilon}(q(t), v) dt \right) \times \right. \right.$$

$$\left. \left. \cos\left[ \frac{1}{\hbar_\varepsilon} \left( \int_0^T dt \left[ \frac{m}{2} \dot{q}^2_{\widetilde{\mathbb{R}}, \varepsilon}(t) - V(q_\varepsilon(t)) \right] \right) \right] \right]^2 + \right. \tag{2.1.12}$$

$$\left[ \int_{q(T)=x} D^+[q(t), \varepsilon] \Psi_{0,\varepsilon}(q(0)) \exp\left( -\frac{1}{\hbar_\varepsilon} \int_0^T Q_{\beta,\varepsilon}(q, v) dt \right) \times \right.$$

$$\left. \left. \sin\left[ \frac{1}{\hbar_\varepsilon} \left( \int_0^T dt \left[ \frac{m_1}{2} \dot{q}^2_{\widetilde{\mathbb{R}}, \varepsilon}(t) - V(q_\varepsilon(t)) \right] \right) \right] \right]^2 \right\} =$$

$$\eta_\varepsilon^{(1)}(T, \hbar_\varepsilon, \alpha, \beta) + \eta_\varepsilon^{(2)}(T, \hbar_\varepsilon, \alpha, \beta).$$

Here

$$\eta_\varepsilon^{(1)}(T, \hbar_\varepsilon, \alpha, \beta) =$$

$$\int dx x \left[ \int\limits_{q(T)=x} D^+[q(t), \varepsilon] \Psi_{0,\varepsilon}(q(0)) \exp\left( -\frac{1}{\hbar_\varepsilon} \int\limits_0^T Q_{\beta,\varepsilon}(q, v) dt \right) \times \right.$$

$$\left. \cos\left[ \frac{1}{\hbar_\varepsilon} \left( \int\limits_0^T dt \left[ \frac{m}{2} \dot{q}_{\tilde{\mathbb{R}},\varepsilon}^2(t) - V(q_\varepsilon(t)) \right] \right) \right] \right]^2 =$$

$$\int dx x \left[ \hat{\eta}_\varepsilon^{(1)}(T, x, \hbar_\varepsilon, \alpha, \beta) \right]^2,$$

$$(2.1.12)$$

$$\eta_\varepsilon^{(2)}(T, \hbar_\varepsilon, \alpha, \beta) =$$

$$\int dx x \left[ \int\limits_{q(T)=x} D^+[q(t), \varepsilon] \Psi_{0,\varepsilon}(q(0)) \exp\left( -\frac{1}{\hbar_\varepsilon} \int\limits_0^T Q_{\beta,\varepsilon}(q, v) dt \right) \times \right.$$

$$\left. \sin\left[ \frac{1}{\hbar_\varepsilon} \left( \int\limits_0^T dt \left[ \frac{m}{2} \dot{q}_{\tilde{\mathbb{R}},\varepsilon}^2(t) - V(q_\varepsilon(t)) \right] \right) \right] \right]^2 =$$

$$\int dx x \left[ \hat{\eta}_\varepsilon^{(2)}(T, x, \hbar_\varepsilon, \alpha, \beta) \right]^2.$$

Here

$$\hat{\eta}_\varepsilon^{(1)}(T, x, x_0, \hbar_\varepsilon, \alpha, \beta) =$$

$$\int\limits_{q(T)=x} D^+[q(t), \varepsilon] \Psi_{0,\varepsilon}(q(0)) \exp\left(-\frac{1}{\hbar_\varepsilon} \int\limits_0^T Q_{\beta,\varepsilon}(q, v) dt\right) \times \qquad (2.1.13)$$

$$\cos\left[\frac{1}{\hbar_\varepsilon}\left(\int\limits_0^T dt\left[\frac{m}{2}\left(\dot{q}_{\widetilde{\mathbb{R}},\varepsilon}^2(t)\right)^2 - V(q_\varepsilon(t))\right]\right)\right]$$

and

$$\hat{\eta}_\varepsilon^{(2)}(T, x, x_0, \hbar_\varepsilon, \alpha, \beta) =$$

$$\int\limits_{q(T)=x} D^+[q(t), \varepsilon] \Psi_{0,\varepsilon}(q(0)) \exp\left(-\frac{1}{\hbar_\varepsilon} \int\limits_0^T Q_{\beta,\varepsilon}(q, v) dt\right) \times \qquad (2.1.14)$$

$$\sin\left[\frac{1}{\hbar_\varepsilon}\left(\int\limits_0^T dt\left[\frac{m}{2}\dot{q}_{\widetilde{\mathbb{R}},\varepsilon}^2(t) - V(q_\varepsilon(t))\right]\right)\right].$$

Let us now to evaluate the quantities $\hat{\eta}_\varepsilon^{(1)}(T, x, x_0, \hbar_\varepsilon, \alpha)$ and $\hat{\eta}_\varepsilon^{(2)}(T, x, x_0, \hbar_\varepsilon, \alpha)$. We assume now that: $V(q_\varepsilon(t)) = V_0(q_\varepsilon(t)) + V_1(q_\varepsilon(t))$. Here

$$V_0(q_\varepsilon(t), t) = aq_\varepsilon^2(t) + b(t)(q_\varepsilon(t)),$$

$$V_1(q_\varepsilon(t)) = a_3 q_\varepsilon^3(t) + \ldots + a_r q_\varepsilon^r(t),$$

$$a > 0, \tag{2.1.15}$$

$$q_\varepsilon(t) = \frac{q(t)}{1 + \varepsilon^k q^m(t)},$$

$$k \geq 1, m \geq 2$$

Substitution Eq.(2.1.15) into Eq.(2.1.13) gives

$$\widehat{\eta}_{\varepsilon}^{(1)}(T,x,x_0,\hbar_{\varepsilon},\alpha,\beta) =$$

$$\int\limits_{q(T)=x} D^{+}[q(t),\varepsilon]\Psi_{0,\varepsilon}(q(0))\exp\left(-\frac{1}{\hbar_{\varepsilon}}\int\limits_{0}^{T}Q_{\beta,\varepsilon}(q,\upsilon)dt\right)\times$$

$$\cos\left[\frac{1}{\hbar_{\varepsilon}}\left(\int\limits_{0}^{T}dt\Big[\frac{m}{2}\dot{q}_{\widetilde{\mathbb{R}},\varepsilon}^{2}(t)-V(q_{\varepsilon}(t))\Big]\right)\right] =$$

$$\int\limits_{q(T)=x} D^{+}[q(t),\varepsilon]\Psi_{0,\varepsilon}(q(0))\exp\left(-\frac{1}{\hbar_{\varepsilon}}\int\limits_{0}^{T}Q_{\beta,\varepsilon}(q,\upsilon)dt\right)\times$$

$$\cos\left[\frac{1}{\hbar_{\varepsilon}}\left(\int\limits_{0}^{T}dt\Big[\frac{m}{2}\dot{q}_{\mathbb{R},\varepsilon}^{2}(t)-V_{0}(q_{\varepsilon}(t))-V_{1}(q_{\varepsilon}(t))\Big]\right)\right] =$$

$$(2.1.16)$$

$$\int\limits_{q(T)=x} D^{+}[q(t),\varepsilon]\Psi_{0,\varepsilon}(q(0))\exp\left(-\frac{1}{\hbar_{\varepsilon}}\int\limits_{0}^{T}Q_{\beta,\varepsilon}(q,\upsilon)dt\right)\times$$

$$\cos\left[\frac{1}{\hbar_{\varepsilon}}\left(\int\limits_{0}^{T}dt\Big[\frac{m}{2}\dot{q}_{\widetilde{\mathbb{R}},\varepsilon}^{2}(t)-V_{0}(q_{\varepsilon}(t))\Big]\right)\right]\cos\left(\frac{1}{\hbar_{\varepsilon}}\int\limits_{0}^{T}V_{1}(q_{\varepsilon}(t))dt\right)+$$

$$\int\limits_{q(T)=x} D^{+}[q(t),\varepsilon]\Psi_{0,\varepsilon}(q(0))\exp\left(-\frac{1}{\hbar_{\varepsilon}}\int\limits_{0}^{T}Q_{\beta,\varepsilon}(q,\upsilon)dt\right)\times$$

$$\sin\left[\frac{1}{\hbar_{\varepsilon}}\left(\int\limits_{0}^{T}dt\Big[\frac{m}{2}\dot{q}_{\mathbb{R},\varepsilon}^{2}(t)-V_{0}(q_{\varepsilon}(t))\Big]\right)\right]\sin\left(\frac{1}{\hbar_{\varepsilon}}\int\limits_{0}^{T}V_{1}(q_{\varepsilon}(t))dt\right) =$$

$$\chi_{\varepsilon}^{(1)}(T,x,x_0,\hbar_{\varepsilon})+\chi_{\varepsilon}^{(2)}(T,x,x_0,\hbar_{\varepsilon}).$$

Here

$$\chi_{\varepsilon}^{(1)}(T,x,x_0,\hbar_{\varepsilon}) =$$

$$\int\limits_{q(T)=x} D^+[q(t),\varepsilon]\Psi_{0,\varepsilon}(q(0))\exp\left(-\frac{1}{\hbar_{\varepsilon}}\int\limits_0^T Q_{\beta,\varepsilon}(q,\upsilon)dt\right)\times \qquad (2.1.17)$$

$$\cos\left[\frac{1}{\hbar_{\varepsilon}}\left(\int\limits_0^T dt\left[\frac{m}{2}\dot{q}_{\widetilde{\mathbb{R}},\varepsilon}^2(t)-V_0(q_{\varepsilon}(t))\right]\right)\right]\cos\left(\frac{1}{\hbar_{\varepsilon}}\int\limits_0^T V_1(q_{\varepsilon}(t))dt\right)$$

and

$$\chi_{\varepsilon}^{(2)}(T,x,x_0,\hbar_{\varepsilon}) =$$

$$\int\limits_{q(T)=x} D^+[q(t),\varepsilon]\Psi_{0,\varepsilon}(q(0))\exp\left(-\frac{1}{\hbar_{\varepsilon}}\int\limits_0^T Q_{\beta,\varepsilon}(q,\upsilon)dt\right)\times \qquad (2.1.18)$$

$$\sin\left[\frac{1}{\hbar_{\varepsilon}}\left(\int\limits_0^T dt\left[\frac{m}{2}\dot{q}_{\widetilde{\mathbb{R}},\varepsilon}^2(t)-V_0(q_{\varepsilon}(t))\right]\right)\right]\sin\left(\frac{1}{\hbar_{\varepsilon}}\int\limits_0^T V_1(q_{\varepsilon}(t))dt\right)$$

Let us now evaluate the quantities $\chi_{\varepsilon}^{(1)}(T,x,x_0,\hbar_{\varepsilon})$ and $\chi_{\varepsilon}^{(2)}(T,x,x_0,\hbar_{\varepsilon})$.
Using replacement $q(t) \rightarrow \sqrt{\hbar_{\varepsilon}}\,q(t)$, from Eq.(2.1.17) and Eq.(2.1.18) we obtain:

$$\chi_\varepsilon^{(1)}(T,x,x_0,\hbar_\varepsilon) =$$

$$\int\limits_{q(T)=\frac{x}{\sqrt{\hbar_\varepsilon}}} D^+[q(t),\varepsilon]\Psi_{0,\varepsilon}\left(\sqrt{\hbar_\varepsilon}\,q(0)\right)\exp\left(-\frac{1}{\hbar_\varepsilon}\int\limits_0^T Q_{\beta,\varepsilon}\left(\sqrt{\hbar_\varepsilon}\,q,\upsilon\right)dt\right)\times$$

$$\cos\left[\frac{1}{\hbar_\varepsilon}\left(\int\limits_0^T dt\left[\frac{m}{2}\hbar_\varepsilon\dot{q}^2_{\tilde{\mathbb{R}},\varepsilon}(t)-V_0\left(\sqrt{\hbar_\varepsilon}\,q_\varepsilon(t)\right)\right]\right)\right]\times$$

$$\cos\left[\frac{1}{\hbar_\varepsilon}\left(\int\limits_0^T V_1\left(\sqrt{\hbar_\varepsilon}\,q_\varepsilon(t)\right)dt\right)\right]=$$

$$\int\limits_{q(T)=\frac{x}{\sqrt{\hbar_\varepsilon}}} D^+[q(t),\varepsilon]\Psi_{0,\varepsilon}\left(\sqrt{\hbar_\varepsilon}\,q(0)\right)\exp\left(-\frac{1}{\hbar_\varepsilon}\int\limits_0^T Q_{\beta,\varepsilon}\left(\sqrt{\hbar_\varepsilon}\,q,\upsilon\right)dt\right)\times$$

$$\cos\left(\int\limits_0^T dt\left[\frac{m}{2}\dot{q}^2_{\tilde{\mathbb{R}},\varepsilon}(t)-aq^2_{\varepsilon,\hbar_\varepsilon}(t)-b(t)\frac{q_{\varepsilon,\hbar_\varepsilon}(t)}{\sqrt{\hbar_\varepsilon}}\right]\right)\times$$

$$\cos\left[\left(\int\limits_0^T V_{1,\hbar_\varepsilon}(q_\varepsilon(t))dt\right)\right]$$

(2.1.19)

and

$$\chi_\varepsilon^{(2)}(T,x,x_0,\hbar_\varepsilon) =$$

$$\int_{q(T)=\frac{x}{\sqrt{\hbar_\varepsilon}}} D^+[q(t),\varepsilon]\,\Psi_{0,\varepsilon}\left(\sqrt{\hbar_\varepsilon}\,q(0)\right)\exp\left(-\frac{1}{\hbar_\varepsilon}\int_0^T Q_{\beta,\varepsilon}\left(\sqrt{\hbar_\varepsilon}\,q,\upsilon\right)dt\right)\times$$

$$\sin\left[\frac{1}{\hbar_\varepsilon}\left(\int_0^T dt\left[\frac{m}{2}\hbar_\varepsilon\dot{q}_{\mathbb{R},\varepsilon}^2(t)-V_0\left(\sqrt{\hbar_\varepsilon}\,q_\varepsilon(t)\right)\right]\right)\right]\times$$

$$\sin\left[\frac{1}{\hbar_\varepsilon}\left(\int_0^T V_1\left(\sqrt{\hbar_\varepsilon}\,q_\varepsilon(t)\right)dt\right)\right] = \qquad (2.1.20)$$

$$\int_{q(T)=\frac{x}{\sqrt{\hbar_\varepsilon}}} D^+[q(t),\varepsilon]\,\Psi_{0,\varepsilon}\left(\sqrt{\hbar_\varepsilon}\,q(0)\right)\exp\left(-\frac{1}{\hbar_\varepsilon}\int_0^T Q_{\beta,\varepsilon}\left(\sqrt{\hbar_\varepsilon}\,q,\upsilon\right)dt\right)\times$$

$$\sin\left(\int_0^T dt\left[\frac{m}{2}\dot{q}_{\widetilde{\mathbb{R}},\varepsilon}^2(t)-aq_{\varepsilon,\hbar_\varepsilon}^2(t)-b(t)\frac{q_{\varepsilon,\hbar_\varepsilon}(t)}{\sqrt{\hbar_\varepsilon}}\right]\right)\times$$

$$\sin\left[\left(\int_0^T V_{1,\hbar_\varepsilon}(q_\varepsilon(t))dt\right)\right]$$

Here

$$V_{1,\hbar_\varepsilon}(q_\varepsilon(t)) = \frac{1}{\hbar_\varepsilon}V_1\left(\sqrt{\hbar_\varepsilon}\,q_{\varepsilon,\epsilon}(t)\right) = a_3\hbar_\varepsilon^{1/2}q_{\varepsilon,\hbar_\varepsilon}^3(t)+\ldots+a_k\hbar_\varepsilon^{(r/2)-1}q_{\varepsilon,\hbar_\varepsilon}^r(t)$$

$$(2.1.21)$$

$$q_{\varepsilon,\hbar_\varepsilon}(t) = \frac{q(t)}{1+\varepsilon^k\hbar_\varepsilon^m q^m(t)}.$$

Let us rewrite the equations Eq.(2.1.19) and Eq.(2.1.20) in the following equivalent

form

$$\chi_\varepsilon^{(1)}(T,x,x_0,\hbar_\varepsilon) =$$

$$\int\limits_{q(T)=\frac{x}{\sqrt{\hbar_\varepsilon}}} D^+[q(t),\varepsilon]\Psi_{0,\varepsilon}\left(\sqrt{\hbar_\varepsilon}\,q(0)\right)\exp\left(-\frac{1}{\hbar_\varepsilon}\int\limits_0^T Q_{\beta,\varepsilon}\left(\sqrt{\hbar_\varepsilon}\,q,v\right)dt\right)\times$$

$$[\mathbf{sign}^+\Re_1(q_{\varepsilon,h_\varepsilon}(t)) - \mathbf{sign}^-\Re_1(q_{\varepsilon,h_\varepsilon}(t))]\times$$

$$\left|\cos\left(\int\limits_0^T dt\left[\frac{m}{2}\dot{q}_{\widetilde{\mathbb{R}},\varepsilon}^2(t) - aq_{\varepsilon,h_\varepsilon}^2(t) - b(t)\frac{q_{\varepsilon,h_\varepsilon}(t)}{\sqrt{\hbar_\varepsilon}}\right]\right)\right|\times \tag{2.1.22}$$

$$\cos\left[\left(\int\limits_0^T V_{1,h_\varepsilon}(q_\varepsilon(t))dt\right)\right],$$

$$\Re_1(q_{\varepsilon,h_\varepsilon}(t)) = \cos\left(\int\limits_0^T dt\left[\frac{m}{2}\dot{q}_{\widetilde{\mathbb{R}},\varepsilon}^2(t) - aq_{\varepsilon,h_\varepsilon}^2(t) - b(t)\frac{q_{\varepsilon,h_\varepsilon}(t)}{\sqrt{\hbar_\varepsilon}}\right]\right)$$

and

$$\chi_\varepsilon^{(2)}(T, x, x_0, \hbar_\varepsilon) =$$

$$\int_{q(T)=\frac{x}{\sqrt{\hbar_\varepsilon}}} D^+[q(t), \varepsilon] \Psi_{0,\varepsilon}\left(\sqrt{\hbar_\varepsilon}\, q(0)\right) \exp\left(-\frac{1}{\hbar_\varepsilon}\int_0^T Q_{\beta,\varepsilon}\left(\sqrt{\hbar_\varepsilon}\, q, \upsilon\right) dt\right) \times$$

$$[\mathbf{sign}^+ \Re_2(q_{\varepsilon,h_\varepsilon}(t)) - \mathbf{sign}^- \Re_2(q_{\varepsilon,h_\varepsilon}(t))] \times$$

$$\left| \sin\left(\int_0^T dt \left[\ \frac{m}{2}\dot{q}^2_{\widetilde{\mathbb{R}},\varepsilon}(t) - a q^2_{\varepsilon,h_\varepsilon}(t) - b(t)\frac{q_{\varepsilon,h_\varepsilon}(t)}{\sqrt{\hbar_\varepsilon}}\ \right]\right) \right| \times \qquad (2.1.23)$$

$$\sin\left[\ \left(\int_0^T V_{1,h_\varepsilon}(q_\varepsilon(t)) dt\right)\ \right],$$

$$\Re_2(q_{\varepsilon,h_\varepsilon}(t)) = \sin\left(\int_0^T dt \left[\ \frac{m}{2}\dot{q}^2_{\widetilde{\mathbb{R}},\varepsilon}(t) - a q^2_{\varepsilon,h_\varepsilon}(t) - b(t)\frac{q_{\varepsilon,h_\varepsilon}(t)}{\sqrt{\hbar_\varepsilon}}\ \right]\right)$$

From Eq.(2.1.22) we obtain:

$$\chi_\varepsilon^{(1)}(T,x,x_0,\hbar_\varepsilon) =$$

$$\int\limits_{q(T)=\frac{x}{\sqrt{\hbar_\varepsilon}}} D^+[q(t),\varepsilon]\Psi_{0,\varepsilon}\left(\sqrt{\hbar_\varepsilon}\,q(0)\right)\exp\left(-\frac{1}{\hbar_\varepsilon}\int\limits_0^T Q_{\beta,\varepsilon}\left(\sqrt{\hbar_\varepsilon}\,q,\nu\right)dt\right)\times$$

$$\{[\mathbf{sign}^+\Re_1(q_{\varepsilon,h_\varepsilon}(t))-\mathbf{sign}^-\Re_1(q_{\varepsilon,h_\varepsilon}(t))+1]-1\}\times$$

$$\left|\cos\left(\int\limits_0^T dt\left[\frac{m}{2}\dot{q}_{\overline{\mathbb{R}},\varepsilon}^2(t)-aq_{\varepsilon,h_\varepsilon}^2(t)-b(t)\frac{q_{\varepsilon,h_\varepsilon}(t)}{\sqrt{\hbar_\varepsilon}}\right]\right)\right|\times$$

$$\cos\left[\left(\int\limits_0^T V_{1,h_\varepsilon}(q_\varepsilon(t))dt\right)\right]=$$

$$\int\limits_{q(T)=\frac{x}{\sqrt{\hbar_\varepsilon}}} \widetilde{D}_1^+[q(t),\varepsilon]\Psi_{0,\varepsilon}\left(\sqrt{\hbar_\varepsilon}\,q(0)\right)\times$$

$$\left|\cos\left(\int\limits_0^T dt\left[\frac{m}{2}\dot{q}_{\overline{\mathbb{R}},\varepsilon}^2(t)-aq_{\varepsilon,h_\varepsilon}^2(t)-b(t)\frac{q_{\varepsilon,h_\varepsilon}(t)}{\sqrt{\hbar_\varepsilon}}\right]\right)\right|\times \qquad (2.1.24)$$

$$\cos\left[\left(\int\limits_0^T V_{1,h_\varepsilon}(q_\varepsilon(t))dt\right)\right]-$$

$$\int\limits_{q(T)=\frac{x}{\sqrt{\hbar_\varepsilon}}} D^+[q(t),\varepsilon]\Psi_{0,\varepsilon}\left(\sqrt{\hbar_\varepsilon}\,q(0)\right)\exp\left(-\frac{1}{\hbar_\varepsilon}\int\limits_0^T Q_{\beta,\varepsilon}\left(\sqrt{\hbar_\varepsilon}\,q,\nu\right)dt\right)\times$$

$$\left|\cos\left(\int\limits_0^T dt\left[\frac{m}{2}\dot{q}_{\overline{\mathbb{R}},\varepsilon}^2(t)-aq_{\varepsilon,h_\varepsilon}^2(t)-b(t)\frac{q_{\varepsilon,h_\varepsilon}(t)}{\sqrt{\hbar_\varepsilon}}\right]\right)\right|\times$$

$$\cos\left[\left(\int\limits_0^T V_{1,h_\varepsilon}(q_\varepsilon(t))dt\right)\right]=$$

where

$$\int\limits_{q(T)=\frac{x}{\sqrt{\hbar_\varepsilon}}} \widetilde{D}_1^+[q(t),\varepsilon](\cdot) =$$

$$\int\limits_{q(t)=\frac{x}{\sqrt{\hbar_\varepsilon}}} D^+[q(t),\varepsilon][\mathbf{sign}^+\Re_1(q_{\varepsilon,\hbar_\varepsilon}(t)) - \mathbf{sign}^-\Re_1(q_{\varepsilon,\hbar_\varepsilon}(t)) + 1](\cdot),$$

$$\widetilde{\chi}_\varepsilon^{(1)}(T,x,x_0,\hbar_\varepsilon) =$$

$$\int\limits_{q(T)=\frac{x}{\sqrt{\hbar_\varepsilon}}} \widetilde{D}^+[q(t),\varepsilon]\Psi_{0,\varepsilon}\left(\sqrt{\hbar_\varepsilon}\,q(0)\right)\exp\left(-\frac{1}{\hbar_\varepsilon}\int\limits_0^T Q_{\beta,\varepsilon}\left(\sqrt{\hbar_\varepsilon}\,q,\upsilon\right)dt\right) \qquad (2.1.25)$$

$$\left|\cos\left(\int\limits_0^T dt\left[\frac{m}{2}\dot{q}_{\widetilde{\mathbb{R}},\varepsilon}^2(t) - aq_{\varepsilon,\hbar_\varepsilon}^2(t) - b(t)\frac{q_{\varepsilon,\hbar_\varepsilon}(t)}{\sqrt{\hbar_\varepsilon}}\right]\right)\right| \times$$

$$\cos\left[\left(\int\limits_0^T V_{1,\hbar_\varepsilon}(q_\varepsilon(t))dt\right)\right]$$

and

$$v_\varepsilon^{(1)}(T, x, x_0, \hbar_\varepsilon) =$$

$$\int\limits_{q(T) = \frac{x}{\sqrt{\hbar_\varepsilon}}} D^+[q(t), \varepsilon] \Psi_{0,\varepsilon}\left(\sqrt{\hbar_\varepsilon}\, q(0)\right) \exp\left(-\frac{1}{\hbar_\varepsilon} \int\limits_0^T Q_{\beta,\varepsilon}\left(\sqrt{\hbar_\varepsilon}\, q, \upsilon\right) dt\right) \times$$

$$\left| \cos\left(\int\limits_0^T dt \left[\, \frac{m}{2} \dot{q}_{\mathbb{R},\varepsilon}^2(t) - a q_{\varepsilon,\hbar_\varepsilon}^2(t) - b(t) \frac{q_{\varepsilon,\hbar_\varepsilon}(t)}{\sqrt{\hbar_\varepsilon}} \,\right]\right) \right| \times$$

$$\cos\left[\left(\int\limits_0^T V_{1,\hbar_\varepsilon}(q_\varepsilon(t)) dt\right)\right].$$

$(2.1.26)$

From Eq.(2.1.23) we obtain:

$$\int\limits_{q(T)=\frac{x}{\sqrt{\hbar_\varepsilon}}} \widetilde{D}_1^+[q(t),\varepsilon](\cdot) =$$

$$\int\limits_{q(T)=\frac{x}{\sqrt{\hbar_\varepsilon}}} D^+[q(t),\varepsilon][\mathbf{sign}^+\Re_2(q_{\varepsilon,h_\varepsilon}(t)) - \mathbf{sign}^-\Re_2(q_{\varepsilon,h_\varepsilon}(t)) + 1](\cdot),$$

$$\chi_\varepsilon^{(2)}(T,x,x_0,\hbar_\varepsilon) =$$

$$\int\limits_{q(T)=\frac{x}{\sqrt{\hbar_\varepsilon}}} D^+[q(t),\varepsilon]\Psi_{0,\varepsilon}\left(\sqrt{\hbar_\varepsilon}\,q(0)\right)\exp\left(-\frac{1}{\hbar_\varepsilon}\int\limits_0^T Q_{\beta,\varepsilon}\left(\sqrt{\hbar_\varepsilon}\,q,\upsilon\right)dt\right)\times$$

$$\{[\mathbf{sign}^+\Re_2(q_{\varepsilon,h_\varepsilon}(t)) - \mathbf{sign}^-\Re_2(q_{\varepsilon,h_\varepsilon}(t)) + 1] - 1\}\times$$

$$\left|\sin\left(\int\limits_0^T dt\left[\frac{m}{2}\dot{q}_{\widetilde{\mathbb{R}},\varepsilon}^2(t) - aq_{\varepsilon,h_\varepsilon}^2(t) - b(t)\frac{q_{\varepsilon,h_\varepsilon}(t)}{\sqrt{\hbar_\varepsilon}}\right]\right)\right|\times$$

$$\sin\left[\left(\int\limits_0^T V_{1,h_\varepsilon}(q_\varepsilon(t))dt\right)\right] =$$

$$\int\limits_{q(T)=\frac{x}{\sqrt{\hbar_\varepsilon}}} \widetilde{D}_2^+[q(t),\varepsilon]\Psi_{0,\varepsilon}\left(\sqrt{\hbar_\varepsilon}\,q(0)\right)\exp\left(-\frac{1}{\hbar_\varepsilon}\int\limits_0^T Q_{\beta,\varepsilon}\left(\sqrt{\hbar_\varepsilon}\,q,\upsilon\right)dt\right)\times$$

(2.1.27)

$$\left|\sin\left(\int\limits_0^T dt\left[\frac{m}{2}\dot{q}_{\widetilde{\mathbb{R}},\varepsilon}^2(t) - aq_{\varepsilon,h_\varepsilon}^2(t) - b(t)\frac{q_{\varepsilon,h_\varepsilon}(t)}{\sqrt{\hbar_\varepsilon}}\right]\right)\right|\times$$

$$\sin\left[\left(\int\limits_0^T V_{1,h_\varepsilon}(q_\varepsilon(t))dt\right)\right] -$$

$$\int\limits_{q(T)=\frac{x}{\sqrt{\hbar_\varepsilon}}} D^+[q(t),\varepsilon]\Psi_{0,\varepsilon}\left(\sqrt{\hbar_\varepsilon}\,q(0)\right)\exp\left(-\frac{1}{\hbar_\varepsilon}\int\limits_0^T Q_{\beta,\varepsilon}\left(\sqrt{\hbar_\varepsilon}\,q,\upsilon\right)dt\right)\times$$

where

$$\widetilde{\chi}_{\varepsilon}^{(2)}(T,x,x_0,\hbar_{\varepsilon}) =$$

$$\int\limits_{q(T)=\frac{x}{\sqrt{\hbar_{\varepsilon}}}} \widetilde{D}_2^+[q(t),\varepsilon]\Psi_{0,\varepsilon}\left(\sqrt{\hbar_{\varepsilon}}\,q(0)\right)\exp\left(-\frac{1}{\hbar_{\varepsilon}}\int\limits_0^T Q_{\beta,\varepsilon}\left(\sqrt{\hbar_{\varepsilon}}\,q,v\right)dt\right)\times$$

$$\left|\sin\left(\int\limits_0^T dt\left[\frac{m}{2}\dot{q}_{\widetilde{\mathbb{R}},\varepsilon}^2(t) - aq_{\varepsilon,\hbar_{\varepsilon}}^2(t) - b(t)\frac{q_{\varepsilon,\hbar_{\varepsilon}}(t)}{\sqrt{\hbar_{\varepsilon}}}\right]\right)\right|\times$$

$$\sin\left[\left(\int\limits_0^T V_{1,\hbar_{\varepsilon}}(q_{\varepsilon}(t))dt\right)\right]$$

(2.1.28)

and

$$v_{\varepsilon}^{(2)}(T,x,x_0,\hbar_{\varepsilon}) =$$

$$\int\limits_{q(T)=\frac{x}{\sqrt{\hbar_{\varepsilon}}}} D^+[q(t),\varepsilon]\Psi_{0,\varepsilon}\left(\sqrt{\hbar_{\varepsilon}}\,q(0)\right)\exp\left(-\frac{1}{\hbar_{\varepsilon}}\int\limits_0^T Q_{\beta,\varepsilon}\left(\sqrt{\hbar_{\varepsilon}}\,q,v\right)dt\right)\times$$

$$\left|\sin\left(\int\limits_0^T dt\left[\frac{m}{2}\dot{q}_{\widetilde{\mathbb{R}},\varepsilon}^2(t) - aq_{\varepsilon,\hbar_{\varepsilon}}^2(t) - b(t)\frac{q_{\varepsilon,\hbar_{\varepsilon}}(t)}{\sqrt{\hbar_{\varepsilon}}}\right]\right)\right|\times$$

(2.1.29)

$$\sin\left[\left(\int\limits_0^T V_{1,\hbar_{\varepsilon}}(q_{\varepsilon}(t))dt\right)\right].$$

Let us estimate now the quantities $\widetilde{\chi}_{\varepsilon}^{(1)}(T,x,x_0,\hbar_{\varepsilon})$ and $\widetilde{\chi}_{\varepsilon}^{(1)}(T,x,x_0,\hbar_{\varepsilon})$.

From Eq.(2.1.25) one obtain directly the inequality:

$$\left| \widetilde{\chi}_{\varepsilon}^{(1)}(T,x,x_0,\hbar_{\varepsilon}) \right| \leq$$

$$\int\limits_{q(T)=\frac{x}{\sqrt{\hbar_{\varepsilon}}}} \widetilde{D}^{+}[q(t),\varepsilon] \Psi_{0,\varepsilon}\left(\sqrt{\hbar_{\varepsilon}}\,q(0)\right) \exp\left(-\frac{1}{\hbar_{\varepsilon}}\int\limits_{0}^{T} Q_{\beta,\varepsilon}\left(\sqrt{\hbar_{\varepsilon}}\,q,\upsilon\right)dt\right) \times$$

$$\left| \cos\left(\int\limits_{0}^{T} dt \left[ \frac{m}{2}\dot{q}_{\mathbb{R},\varepsilon}^{2}(t) - aq_{\varepsilon,h_{\varepsilon}}^{2}(t) - b(t)\frac{q_{\varepsilon,h_{\varepsilon}}(t)}{\sqrt{\hbar_{\varepsilon}}} \right]\right) \right| \times$$

$$\left| \cos\left[ \left(\int\limits_{0}^{T} V_{1,h_{\varepsilon}}(q_{\varepsilon}(t))dt\right) \right] \right|$$

(2.1.30)

From Eq.(2.1.30) and Hölder's inequality () for $p_2 = \dfrac{1}{\sqrt{\hbar_{\varepsilon}}}$, $p_1 = \dfrac{1}{1-\sqrt{\hbar_{\varepsilon}}} = 1 + \sqrt{\hbar_{\varepsilon}} + o(\hbar_{\varepsilon}) \simeq 1$ we obtain:

$$\left| \widetilde{\chi}_\varepsilon^{(1)}(T, x, x_0, h_\varepsilon) \right| \leq$$

$$\left[ \int\limits_{q(T) = \frac{x}{\sqrt{h_\varepsilon}}} \widetilde{D}^+[q(t), \varepsilon] \Psi_{0,\varepsilon}\left(\sqrt{h_\varepsilon}\, q(0)\right) \exp\left( -\frac{1}{h_\varepsilon} \int\limits_0^T Q_{\beta,\varepsilon}\left(\sqrt{h_\varepsilon}\, q, \upsilon\right) dt \right) \times \right.$$

$$\left. \left| \cos\left( \int\limits_0^T dt \left[ \frac{m}{2} \dot{q}_{\widetilde{\mathbb{R}},\varepsilon}^2(t) - a q_{\varepsilon,h_\varepsilon}^2(t) - b(t) \frac{q_{\varepsilon,h_\varepsilon}(t)}{\sqrt{h_\varepsilon}} \right] \right) \right|^{p_1} \right]^{\frac{1}{p_1}} \times \tag{2.1.31}$$

$$\left[ \int\limits_{q(T) = \frac{x}{\sqrt{h_\varepsilon}}} \widetilde{D}^+[q(t), \varepsilon] \Psi_{0,\varepsilon}\left(\sqrt{h_\varepsilon}\, q(0)\right) \left| \cos\left[ \left( \int\limits_0^T V_{1,h_\varepsilon}(q_\varepsilon(t)) dt \right) \right] \right|^{p_2} \right]^{\frac{1}{p_2}} \leq$$

Let us

$$\Theta_\varepsilon(T, x, x_0, \hbar_\varepsilon) =$$

$$\int\limits_{q(T)=\frac{x}{\sqrt{\hbar_\varepsilon}}} D^+[q(t), \varepsilon]\, \Psi_{0,\varepsilon}\left(\sqrt{\hbar_\varepsilon}\, q(0)\right) \exp\left(-\frac{1}{\hbar_\varepsilon} \int\limits_0^T Q_{\beta,\varepsilon}\left(\sqrt{\hbar_\varepsilon}\, q, \upsilon\right) dt\right) \times \qquad (2.1.)$$

$$\cos\left(\int\limits_0^T dt\left[\frac{m}{2}\dot{q}^2_{\mathbb{R},\varepsilon}(t) - a q^2_{\varepsilon,\hbar_\varepsilon}(t) - b(t)\frac{q_{\varepsilon,\hbar_\varepsilon}(t)}{\sqrt{\hbar_\varepsilon}}\right]\right)$$

$$\widehat{\Theta}_{\beta,\varepsilon}(T,x,x_0,\hbar_\varepsilon) =$$

$$\int |x|dx \int\limits_{q(T)=x} D^+[q(t),\varepsilon]\Psi_{0,\varepsilon}(q(0)) \exp\left(-\frac{1}{\hbar_\varepsilon}\int\limits_0^T Q_{\beta,\varepsilon}(q,\upsilon)dt\right) \times$$

$$\exp\left(\frac{i}{\hbar_\varepsilon}\int\limits_0^T dt\Big[\frac{m}{2}\dot{q}^2_{\widetilde{\mathbb{K}},\varepsilon}(t) - aq^2_{\varepsilon,h_\varepsilon}(t) - b(t)q_{\varepsilon,h_\varepsilon}(t)\Big]\right) =$$

$$\int |x|dx \int\limits_{q(T)=x} D^+[q(t),\varepsilon]\Psi_{0,\varepsilon}(q(0)) \times$$

$$\exp\left(-\frac{\upsilon}{\hbar_\varepsilon}\int\limits_0^T\Big[(q(t)-\beta(t,q(0),q(T)))^2\Big]dt - \frac{\upsilon}{\hbar_\varepsilon}q^2(T) - \frac{\upsilon}{\hbar_\varepsilon}q^2(0)\right)$$

$$\exp\left(\frac{i}{\hbar_\varepsilon}\int\limits_0^T \mathcal{L}\Big(\dot{q}_{\widetilde{\mathbb{K}},\varepsilon}, q_{\varepsilon,h_\varepsilon}\Big)dt\right) \simeq$$

$$\text{(2.1.)}$$

$$\int |x|dx \int\limits_{\substack{q(T)=x \\ q(0)=y}} \int dyD^+[q(t),\varepsilon]\Psi_{0,\varepsilon}(y) \exp\left(\frac{i}{\hbar_\varepsilon}\int\limits_0^T \mathcal{L}(\dot{q},q)dt\right) \times$$

$$\exp\left(-\frac{\upsilon}{\hbar_\varepsilon}\int\limits_0^T[(q(t)-\beta(t,q(0),q(T)))]^2dt - \frac{\upsilon}{\hbar_\varepsilon}x^2 - \frac{\upsilon}{\hbar_\varepsilon}y^2\right),$$

$$\mathcal{L}\Big(\dot{q}_{\widetilde{\mathbb{K}},\varepsilon}, q_{\varepsilon,h_\varepsilon}\Big) = \frac{m}{2}\dot{q}^2_{\widetilde{\mathbb{K}},\varepsilon}(t) - aq^2_{\varepsilon,h_\varepsilon}(t) - b(t)q_{\varepsilon,h_\varepsilon}(t) \simeq$$

$$\mathcal{L}\Big(\dot{q}_{\widetilde{\mathbb{K}},\varepsilon}, q\Big) = \frac{m}{2}\dot{q}^2_{\widetilde{\mathbb{K}},\varepsilon}(t) - aq^2(t) - b(t)q(t) \text{ if } \varepsilon \to 0.$$

Let us estimate now path integral

$$\widehat{\Theta}_{\beta,\varepsilon}(T,x,x_0,\hbar_\varepsilon) =$$

$$\int\limits_{\substack{q(T)=x \\ q(0)=y}} D^+[q(t),\varepsilon]\Psi_{0,\varepsilon}(q(0))q(T)\exp\left(-\frac{1}{\hbar_\varepsilon}\int\limits_0^T Q_{\beta,\varepsilon}(q,\upsilon)dt\right)\times$$

$$\exp\left(\frac{i}{\hbar_\varepsilon}\int\limits_0^T dt\left[\frac{m}{2}\dot{q}_{\widetilde{\mathbb{R}},\varepsilon}^2(t) - aq_{\varepsilon,h_\varepsilon}^2(t) - b(t)q_{\varepsilon,h_\varepsilon}(t)\right]\right) =$$

$$\int\limits_{\substack{q(T)=x \\ q(0)=y}} D^+[q(t),\varepsilon]\Psi_{0,\varepsilon}(q(0))q(T)\times$$

$$\exp\left(\frac{i}{\hbar_\varepsilon}\int\limits_0^T \mathcal{L}\left(\dot{q}_{\widetilde{\mathbb{R}},\varepsilon},q_{\varepsilon,h_\varepsilon}\right)dt\right)\times$$

$$\exp\left(-\frac{\upsilon}{\hbar_\varepsilon}\int\limits_0^T[(q(t)-c-\beta(t,q(0)-c,q(T)-c))]^2 dt \right. \qquad (2.1.)$$

$$\left. -\frac{\upsilon}{\hbar_\varepsilon}q^2(T) - \frac{\upsilon}{\hbar_\varepsilon}q^2(0)\right) \simeq$$

$$\int\limits_{\substack{q(T)=x \\ q(0)=y}}\int dy D^+[q(t),\varepsilon]\Psi_{0,\varepsilon}(y)q(T)\exp\left(\frac{i}{\hbar_\varepsilon}\int\limits_0^T\mathcal{L}\left(\dot{q}_{\widetilde{\mathbb{R}},\varepsilon},q\right)dt\right)\times$$

$$\exp\left(-\frac{\upsilon}{\hbar_\varepsilon}\int\limits_0^T[(q(t)-c-\beta(t,q(0)-c,q(T)-c))]^2 dt - \right.$$

$$\left. \frac{\upsilon}{\hbar_\varepsilon}[q(T)-c]^2 - \frac{\upsilon}{\hbar_\varepsilon}[q(0)-c]\right),$$

$$\mathcal{L}\left(\dot{q}_{\widetilde{\mathbb{R}},\varepsilon},q_{\varepsilon,h_\varepsilon}\right) = \frac{m}{2}\dot{q}_{\widetilde{\mathbb{R}},\varepsilon}^2(t) - aq_{\varepsilon,h_\varepsilon}^2(t) - b(t)q_{\varepsilon,h_\varepsilon}(t) \simeq$$

We set now for short $b(t) = b = const,$ without loss of generality. Substitution

$$q(t) = q_1(t) + c,$$

$$c = -\frac{b}{2a}$$

(2.1.)

into (2.1.) gives

$$\mathcal{L}\left(\dot{q}_{\widetilde{\mathbb{R}},\varepsilon}, q\right) = \frac{m}{2}\dot{q}^2_{1\widetilde{\mathbb{R}},\varepsilon} - a(q_1 + c)^2 - b(q_1 + c) =$$

$$\frac{m}{2}\dot{q}^2_{1\widetilde{\mathbb{R}},\varepsilon} - aq_1^2 - 2acq_1 - bq_1 - ac^2 - bc =$$

(2.1.)

$$\mathcal{L}_1\left(\dot{q}_{1\widetilde{\mathbb{R}},\varepsilon}, q_1\right) = \frac{m}{2}\dot{q}^2_{1\widetilde{\mathbb{R}},\varepsilon} - aq_1^2 + \frac{b^2}{2a}.$$

and

$$\Theta_{\beta,\varepsilon}(T,x,x_0,\hbar_\varepsilon) =$$

$$\exp\left(\frac{i}{\hbar_\varepsilon}\,\frac{b^2 T}{2a}\right)\int\limits_{\substack{q_1(T)=x-c\\q_1(0)=y-c}} D^+[q(t),\varepsilon]\Psi_{0,\varepsilon}(q_1(0)+c)[q_1(T)+c]\times$$

$$\exp\left(\frac{i}{\hbar_\varepsilon}\int\limits_0^T \mathcal{L}_1\!\left(\dot{q}_{1\widetilde{\mathbb{R}},\varepsilon},q_1\right)dt\right)\times \qquad\qquad (2.1.)$$

$$\exp\left(-\frac{\upsilon}{\hbar_\varepsilon}\int\limits_0^T [(q_1(t)-\beta(t,q_1(0),q(T)))]^2 dt -\right.$$

$$\left.\frac{\upsilon}{\hbar_\varepsilon}q_1^2(T)-\frac{\upsilon}{\hbar_\varepsilon}q_1^2(0)\right)$$

Let us estimate now path integral

$$\Theta'_{\beta,\varepsilon}(T, x, x_0, \hbar_\varepsilon) =$$

$$\int\limits_{\substack{q_1(T)=x-c \\ q_1(0)=y-c}} D^+[q(t), \varepsilon] \Psi_{0,\varepsilon}(q_1(0) + c)[q_1(T) + c] \times$$

$$\exp\left( \frac{i}{\hbar_\varepsilon} \int\limits_0^T \mathcal{L}_1\left(\dot{q}_{1\widetilde{\mathbb{R}},\varepsilon}, q_1\right) dt \right) \times \qquad (2.1.)$$

$$\exp\left( -\frac{\upsilon}{\hbar_\varepsilon} \int\limits_0^T [(q_1(t) - \beta(t, q_1(0), q(T)))]^2 dt - \right.$$

$$\left. \frac{\upsilon}{\hbar_\varepsilon} q_1^2(T) - \frac{\upsilon}{\hbar_\varepsilon} q_1^2(0) \right)$$

using stationary phase method.

The Euler-Lagrange equation of $\mathcal{L}\left(\dot{q}_{\widetilde{\mathbb{R}},\varepsilon}, q\right)$ with corresponding boundary conditions is

$$\frac{d}{dt}\left( \frac{\partial \mathcal{L}}{\partial \dot{q}_{\widetilde{\mathbb{R}},\varepsilon}} \right) - \frac{\partial \mathcal{L}}{\partial q} = 0,$$

$$q(0) = y - c, q(T) = x - c. \qquad (2.1.)$$

Thus if $\varepsilon \to 0$, Euler-Lagrange equation of $\mathcal{L}\left(\dot{q}_{\widetilde{\mathbb{R}},\varepsilon}, q\right)$ is

$$\frac{d}{dt}\left(\frac{\partial \mathcal{L}}{\partial \dot{q}}\right) - \frac{\partial \mathcal{L}}{\partial q} = 0,$$

$$(2.1.)$$

$$q(0) = y - c, q(T) = x - c.$$

Let $q_{\mathbf{cr}}(t, y, x)$ be the solution of the Euler-Lagrange equation (2.1. ). Substitution

$$q(t) = \delta(t) + q_{\mathbf{cr}}(t, y, x)$$

$$(2.1.)$$

$$\delta(0) = 0, \delta(T) = 0$$

into (2.1. ) gives

$$\Theta'_{\beta, \varepsilon}(T, x_0, \hbar_\varepsilon) =$$

$$\int\limits_{\substack{\delta(T)=0 \\ \delta(0)=0}} D^+[\delta(t), \varepsilon] x \Psi^2_{0, \varepsilon}(y) \exp\left(\frac{i}{\hbar_\varepsilon} \widetilde{S}_{\mathbf{cl}}(T, y - c, x - c)\right) \times$$

$$\exp\left(\frac{i}{\hbar_\varepsilon} \int\limits_0^T \mathcal{L}\left(\dot{\delta}_{\widetilde{\mathbb{R}}, \varepsilon}, \delta\right) dt\right) \times \qquad (2.1.)$$

$$\exp\left(-\frac{\upsilon}{\hbar_\varepsilon} \int\limits_0^T \Big[ (\delta(t) + q_{\mathbf{cr}}(t, y - c, x - c) - \beta(t, y - c, x - c))^2 \Big] dt$$

$$-\frac{\upsilon}{\hbar_\varepsilon}(x - c)^2 - \frac{\upsilon}{\hbar_\varepsilon}(y - c)^2\Big).$$

Where

$$\widetilde{S}_{\mathbf{el}}(T, y - c, x - c) = \tag{2.1.}$$

We set now

$$\beta(t, y - c, x - c) = q_{\mathbf{cr}}(t, y - c, x - c), \tag{2.1.}$$

From Eq.(2.1 ) and Eq.(2.1 ) we obtain

$$\widehat{\Theta}_{\beta,\varepsilon}(T, x_0, \hbar_\varepsilon) =$$

$$\int |x| dx \int dy \exp\left( \frac{i}{\hbar_\varepsilon} \Big[ \widetilde{S}_{\mathbf{el}}(T, y, x) + iv x^2 + iv y^2 \Big] \right)$$

$$\int_{\substack{\delta(T)=0 \\ \delta(0)=0}} D^+[\delta(t), \varepsilon] \exp\left( \frac{i}{\hbar_\varepsilon} \left[ \int_0^T \mathcal{L}\left( \dot{\delta}_{\widetilde{\mathbb{R}},\varepsilon}, \delta \right) dt + iv \int_0^T \delta^2(t) dt \right] \right), \tag{2.1.}$$

From Eq.(2.2.1.24) and Hölder's inequality (2.2.2.17) for $p_2 = \dfrac{1}{\sqrt{\hbar_\varepsilon}}$,

$p_1 = \dfrac{1}{1 - \sqrt{\hbar_\varepsilon}} = 1 + \sqrt{\hbar_\varepsilon} + o(\hbar_\varepsilon) \simeq 1$ we obtain:

$$|\widetilde{v}_{2,\varepsilon}(T,x,x_0,\alpha,\hbar)| \leq$$

$$\left| \int\limits_{q(T)=\frac{x}{\sqrt{\varepsilon}}} D[q(t)]\Psi_0(\sqrt{\varepsilon}\,q(0))\rho_{q,\varepsilon}(\bullet)\exp\left(-\frac{1}{\hbar}\int\limits_0^T\left(\sqrt{\hbar}\,q(t)-\alpha(t)\right)^2 dt\right) \times \right.$$

$$\sin\left(\int\limits_0^T dt\left[\frac{m}{2}\left(\frac{dq(t)}{dt}\right)^2 - aq_\varepsilon^2(t)\right]\right) \times$$

$$\left. \sin\left(\int\limits_0^T dt\left[b(t)\frac{q_\varepsilon(t)}{\sqrt{\varepsilon}}\right]\right)\cos\left(\int\limits_0^T V_{\varepsilon,\epsilon}(q_\varepsilon(t))dt\right)\right| \leq$$

$$\int\limits_{q(T)=\frac{x}{\sqrt{\varepsilon}}} D[q(t)]\Psi_0(\sqrt{\varepsilon}\,q(0))|\rho_{q,\varepsilon}(\bullet)|\exp\left(-\frac{1}{\hbar}\int\limits_0^T\left(\sqrt{\hbar}\,q(t)-\alpha(t)\right)^2 dt\right) \times$$

$$\left|\sin\left(\int\limits_0^T dt\left[\frac{m}{2}\left(\frac{dq(t)}{dt}\right)^2 - aq_\varepsilon^2(t)\right]\right)\right| \times$$

$$\left|\sin\left(\int\limits_0^T dt\left[b(t)\frac{q_\varepsilon(t)}{\sqrt{\varepsilon}}\right]\right)\right| \times \left|\cos\left(\int\limits_0^T V_{\varepsilon,\epsilon}(q_\varepsilon(t))dt\right)\right| \leq$$

$$\tag{2.2.1.26}$$

$$\left[\int\limits_{q(T)=\frac{x}{\sqrt{\varepsilon}}} D[q(t)]|\Psi_0(\sqrt{\varepsilon}\,q(0))|^{1+\sqrt{\epsilon}}|\rho_{q,\varepsilon}(\bullet)|^{1+\sqrt{\epsilon}} \times \right.$$

$$\exp\left(-\frac{1}{\hbar}\int\limits_0^T\left(\sqrt{\hbar}\,q(t)-\alpha(t)\right)^2 dt\right) \times$$

$$\left.\left|\sin\left(\int\limits_0^T dt\left[b(t)\frac{q_\varepsilon(t)}{\sqrt{\varepsilon}}\right]\right)\sin\left(\int\limits_0^T dt\left[\frac{m}{2}\left(\frac{dq(t)}{dt}\right)^2 - aq_\varepsilon^2(t)\right]\right)\right|^{1+\sqrt{\epsilon}}\right]^{1-\sqrt{\epsilon}} \times$$

Note that $|\rho_q(\cdot)|^{1+\sqrt{\epsilon}} = 1$. Let us calculate now the oscillatory path integral

$$\theta_\varepsilon(T, x, x_0, \alpha, \hbar) =$$

$$\int\limits_{q(T)=\frac{x}{\sqrt{\epsilon}}} D[q(t)] |\Psi_0(\sqrt{\epsilon}\, q(0))| \exp\left( -\frac{1}{\hbar} \int\limits_0^T \left( \sqrt{\hbar}\, q(t) - \alpha(t) \right)^2 dt \right) \times$$

$$\left| \cos\left( \int\limits_0^T dt \left[ \frac{m}{2} \left( \frac{dq(t)}{dt} \right)^2 - a q_\varepsilon^2(t) \right] \right) \right| \times$$

$$\left| \cos\left( \int\limits_0^T dt \left[ b(t) \frac{q_\varepsilon(t)}{\sqrt{\hbar}} \right] \right) \right| = \tag{2.2.1.27}$$

$$\int\limits_{q(T)=\frac{x}{\sqrt{\epsilon}}} D[q(t)] |\Psi_0(\sqrt{\epsilon}\, q(0))| \exp\left( -\frac{1}{\hbar} \int\limits_0^T \left( \sqrt{\hbar}\, q(t) - \alpha(t) \right)^2 dt \right) \times$$

$$\mathbf{sign}\left[ S_0^{(1)}(T, x, x_0, q_\varepsilon(t), \epsilon) \right] \times \mathbf{sign}\left[ S_0^{(2)}(T, x, x_0, q_\varepsilon(t), \epsilon) \right] \times$$

$$\cos\left( \int\limits_0^T dt \left[ b(t) \frac{q_\varepsilon(t)}{\sqrt{\sqrt{\hbar}}} \right] \right) \cos\left( \int\limits_0^T dt \left[ \frac{m}{2} \left( \frac{dq(t)}{dt} \right)^2 - a q_\varepsilon^2(t) \right] \right).$$

Here

$$S_0^{(1)}(T, x, x_0, q(t), \hbar) = \cos\left( \int\limits_0^T dt \left[ \frac{m}{2} \left( \frac{dq(t)}{dt} \right)^2 - a q_{\varepsilon,\epsilon}^2(t) \right] \right),$$

$$\tag{2.2.1.28}$$

$$S_0^{(2)}(T, x, x_0, q(t), \hbar) = \cos\left( \int\limits_0^T dt \left[ b(t) \frac{q_{\varepsilon,\epsilon}(t)}{\sqrt{\hbar}} \right] \right).$$

From Eq.(2.2.1.28) we obtain

$$\theta_\varepsilon(T,x,x_0,\alpha,\hbar) =$$

$$\int\limits_{q(T)=\frac{x}{\sqrt{\epsilon}}} D[q(t)]|\Psi_0(\sqrt{\epsilon}\,q(0))|\widetilde{\rho}_{q,\varepsilon}^{(1)}(\cdot)\exp\left(-\frac{1}{\hbar}\int\limits_0^T\left(\sqrt{\hbar}\,q(t)-\alpha(t)\right)^2 dt\right)\times$$

$$\left\{\left[\exp\left[i\widetilde{S}_0^{(1)}(T,x,x_0,q_{\varepsilon,\epsilon}(t),\epsilon)\right]+\exp\left[-i\widetilde{S}_0^{(1)}(T,x,x_0,q_{\varepsilon,\epsilon}(t),\epsilon)\right]\right]\times\right.$$

$$\left.\left[\exp\left[i\widetilde{S}_0^{(2)}(T,x,x_0,q_{\varepsilon,\epsilon}(t),\epsilon)\right]+\exp\left[-i\widetilde{S}_0^{(2)}(T,x,x_0,q_{\varepsilon,\epsilon}(t),\epsilon)\right]\right]\right\} = \quad (2.2.1.29)$$

$$\int\limits_{q(T)=\frac{x}{\sqrt{\epsilon}}} D[q(t)]|\Psi_0(\sqrt{\epsilon}\,q(0))|\widetilde{\rho}_{q,\varepsilon,\epsilon}^{(1)}(\cdot)\exp\left(-\frac{1}{\hbar}\int\limits_0^T\left(\sqrt{\hbar}\,q(t)-\alpha(t)\right)^2 dt\right)\times$$

$$\left\{\exp\left[i\left(\widetilde{S}_0^{(1)}(\cdot)+\widetilde{S}_0^{(2)}(\cdot)\right)\right]+\exp\left[i\left(-\widetilde{S}_0^{(1)}(\cdot)+S_0^{(2)}(\cdot)\right)\right]+\right.$$

$$\left.\exp\left[i\left(\widetilde{S}_0^{(1)}(\cdot)-\widetilde{S}_0^{(2)}(\cdot)\right)\right]+\exp\left[i\left(-\widetilde{S}_0^{(1)}(\cdot)-\widetilde{S}_0^{(2)}(\cdot)\right)\right]\right\}.$$

Here

$$\widetilde{\rho}_{q,\varepsilon,\epsilon}^{(1)}(\cdot) = \mathbf{sign}\Big[ S_0^{(1)}(T,x,x_0,q_{\varepsilon,\epsilon}(t),\hbar) \Big]\mathbf{sign}\Big[ S_0^{(2)}(T,x,x_0,q_\varepsilon(t),\hbar) \Big],$$

$$\widetilde{S}_0^{(1)}(T,x,x_0,q_\varepsilon(t),\epsilon) = \int\limits_0^T dt\Bigg[ \frac{m}{2}\left( \frac{dq(t)}{dt} \right)^2 - aq_{\varepsilon,\epsilon}^2(t) \Bigg], \qquad (2.2.1.30)$$

$$\widetilde{S}_0^{(2)}(T,x,x_0,q_\varepsilon(t),\epsilon) = \int\limits_0^T dt\Bigg[ b(t)\frac{q_\varepsilon(t)}{\sqrt{\hbar}} \Bigg].$$

From Eq.(2.2.1.29) we obtain

$$\theta_\varepsilon(T,x,x_0,\hbar) =$$

$$\theta_\varepsilon^{(1)}(T,x,x_0,\alpha,\hbar) + \theta_\varepsilon^{(2)}(T,x,x_0,\alpha,\hbar) + \qquad (2.2.1.30)$$

$$+\theta_\varepsilon^{(3)}(T,x,x_0,\alpha,\hbar) + \theta_\varepsilon^{(4)}(T,x,x_0,\alpha,\hbar).$$

Here

$$\theta_\varepsilon^{(1)}(T,x,x_0,\alpha,\hbar) =$$

$$\int\limits_{q(T)=\frac{x}{\sqrt{\hbar}}} D[q(t)] \left|\Psi_0\left(\sqrt{\hbar}\,q(0)\right)\right| \widetilde{\rho}_{q,\varepsilon}^{(1)}(\cdot) \exp\left(-\frac{1}{\hbar}\int\limits_0^T \left(\sqrt{\hbar}\,q(t)-\alpha(t)\right)^2 dt\right) \times$$

$$\exp\left[i\left(\widetilde{S}_0^{(1)}(\cdot) + \widetilde{S}_0^{(2)}(\cdot)\right)\right],$$

$$\theta_\varepsilon^{(2)}(T,x,x_0,\alpha,\hbar) =$$

$$\int\limits_{q(T)=\frac{x}{\sqrt{\hbar}}} D[q(t)] \left|\Psi_0\left(\sqrt{\hbar}\,q(0)\right)\right| \widetilde{\rho}_{q,\varepsilon}^{(1)}(\cdot) \exp\left(-\frac{1}{\hbar}\int\limits_0^T \left(\sqrt{\hbar}\,q(t)-\alpha(t)\right)^2 dt\right) \times$$

$$\exp\left[-i\left(\widetilde{S}_0^{(1)}(\cdot) - \widetilde{S}_0^{(2)}(\cdot)\right)\right],$$

$$\theta_\varepsilon^{(3)}(T,x,x_0,\epsilon) = \tag{2.2.1.31}$$

$$\int\limits_{q(T)=\frac{x}{\sqrt{\hbar}}} D[q(t)] \left|\Psi_0\left(\sqrt{\hbar}\,q(0)\right)\right| \widetilde{\rho}_{q,\varepsilon}^{(1)}(\cdot) \exp\left(-\frac{1}{\hbar}\int\limits_0^T \left(\sqrt{\hbar}\,q(t)-\alpha(t)\right)^2 dt\right) \times$$

$$\exp\left[i\left(\widetilde{S}_0^{(1)}(\cdot) - \widetilde{S}_0^{(2)}(\cdot)\right)\right],$$

$$\theta_\varepsilon^{(4)}(T,x,x_0,\alpha,\hbar) =$$

$$\int\limits_{q(T)=\frac{x}{\sqrt{\epsilon}}} D[q(t)] \left|\Psi_0\left(\sqrt{\hbar}\,q(0)\right)\right| \widetilde{\rho}_{q,\varepsilon}^{(1)}(\cdot) \exp\left(-\frac{1}{\hbar}\int\limits_0^T \left(\sqrt{\hbar}\,q(t)-\alpha(t)\right)^2 dt\right) \times$$

$$\exp\left[-i\left(\widetilde{S}_0^{(1)}(\cdot) + \widetilde{S}_0^{(2)}(\cdot)\right)\right].$$

Let's calculate now path integral $\theta_\varepsilon^{(1)}(T, x, x_0, \alpha, \hbar)$ :

$$\theta_\varepsilon^{(1)}(T, x, x_0, \epsilon) =$$

$$\int_{q(T)=\frac{x}{\sqrt{\epsilon}}} D[q(t)] \left| \Psi_0\left(\sqrt{\hbar}\,q(0)\right) \right| \widetilde{\rho}_{q,\varepsilon}^{(1)}(\cdot) \exp\left( -\frac{1}{\hbar} \int_0^T \left(\sqrt{\hbar}\,q(t) - \alpha(t)\right)^2 dt \right) \times$$

$$\exp\left[ i\left( \widetilde{S}_0^{(1)}(\cdot) + \widetilde{S}_0^{(2)}(\cdot) \right) \right] =$$

$$\int_{q(T)=\frac{x}{\sqrt{\epsilon}}} D[q(t)] \left| \Psi_0\left(\sqrt{\hbar}\,q(0)\right) \right| \widetilde{\rho}_{q,\varepsilon}^{(1)}(\cdot) \times \qquad (2.1.32)$$

$$\exp\left\{ i \int_0^T dt \left[ \frac{m}{2}\left(\frac{dq(t)}{dt}\right)^2 - a q_\varepsilon^2(t) + b(t)\frac{q_\varepsilon(t)}{\sqrt{\hbar}} \right] \right\} \times$$

$$\exp\left( -\frac{1}{\hbar} \int_0^T \left(\sqrt{\hbar}\,q(t) - \alpha(t)\right)^2 dt \right)$$

by using saddle point method.Using replacement $q_\varepsilon(t) := \dfrac{q_{\varepsilon, \hbar_\varepsilon}(t)}{\sqrt{\hbar_\varepsilon}}$ we rewrite

(2.2.1.32) in the next equivalent form

$$\theta_\varepsilon^{(1)}(T, x, x_0, \hbar_\varepsilon) =$$

$$\int\limits_{q(T)=x} D^+[q(t), \varepsilon] |\Psi_{0,\varepsilon}(q(0))| \times$$

$$\exp\left\{ \frac{i}{\hbar} \int\limits_0^T dt \left[ \frac{m}{2} \dot{q}_{\mathbb{R},\varepsilon}^2(t) - aq_\varepsilon^2(t) + b(t)q_\varepsilon(t) \right] \right\} \times$$

$$\exp\left( -\nu\hbar_\varepsilon \int\limits_0^T (q_\varepsilon(t) - \alpha(t))^2 dt \right) \exp\left[ -\frac{1}{\hbar_\varepsilon}(\gamma_T^2 q(0) + \theta_T^2 q(T)) \right] =$$

$$\int\limits_{q(T)=x} D^+[q(t), \varepsilon] |\Psi_{0,\varepsilon}(q(0))| \times \tag{2.1.33}$$

$$\exp\left\{ \frac{i}{\hbar} \int\limits_0^T dt \left[ \frac{m}{2} \dot{q}_{\mathbb{R},\varepsilon}^2(t) - aq_\varepsilon^2(t) + b(t)q_\varepsilon(t) + i\hbar_\varepsilon(q_\varepsilon(t) - \alpha(t))^2 \right] \right\} =$$

$$\int\limits_{q(T)=x} D^+[q(t), \varepsilon] |\Psi_{0,\varepsilon}(q(0))| \exp\left[ -\frac{1}{\hbar_\varepsilon}(\gamma_T^2 q(0) + \theta_T^2 q(T)) \right]$$

$$\exp\left\{ \frac{i}{\hbar} \int\limits_0^T dt \left[ \frac{m}{2} \dot{q}_{\mathbb{R},\varepsilon}^2(t) - (a - \hbar_\varepsilon i)q_\varepsilon^2(t) + (b(t) - 2i\hbar_\varepsilon\alpha(t))q_\varepsilon(t) + i\hbar_\varepsilon\alpha^2(t) \right] \right\}.$$

We assume now that

$$\Psi(x, 0) = \frac{1}{\sqrt[4]{2\pi\hbar_\varepsilon}} \exp\left( -\frac{x^2}{2\hbar_\varepsilon} \right), \tag{2.1.34}$$

Substitution Eq.(2.1.34) into Eq.(2.1.33) gives

$$\theta_\varepsilon^{(1)}(T, x, \alpha, \hbar_\varepsilon) =$$

$$\frac{1}{\sqrt[4]{2\pi\hbar_\varepsilon}} \int\limits_{\substack{q(T)=x \\ q(0)=y}} dy D^+[q(t), \varepsilon] \exp\left(-\frac{(y-x_0)^2}{2\hbar_\varepsilon}\right) \exp\left[-\frac{1}{\hbar_\varepsilon}(\gamma_T^2 q^2(0) + \theta_T^2 q^2(T))\right] \times$$

$$\exp\left\{\frac{i}{\hbar} \int\limits_0^T dt\left[\frac{m}{2} q_{\widetilde{\mathbb{R}},\varepsilon}^2(t) - \hat{a} q_\varepsilon^2(t) + \hat{b}(t) q_\varepsilon(t)\right]\right\} =$$

$$\frac{1}{\sqrt[4]{2\pi\hbar_\varepsilon}} \int\limits_{\substack{q(T)=x \\ q(0)=y}} dy D^+[q(t), \varepsilon] \exp\left(-\frac{(y-x_0)^2}{2\hbar_\varepsilon}\right) \exp\left[-\frac{1}{\hbar_\varepsilon}(\gamma_T^2 y^2 + \theta_T^2 x^2)\right] \times \qquad (2.1.35)$$

$$\exp\left\{\frac{i}{\hbar_\varepsilon} \int\limits_0^T \mathcal{L}\left(\dot{q}_{\widetilde{\mathbb{R}},\varepsilon}(t), q_\varepsilon(t)\right) dt\right\},$$

$$\mathcal{L}\left(\dot{q}_{\widetilde{\mathbb{R}},\varepsilon}(t), q_\varepsilon(t)\right) = \frac{m}{2} q_{\widetilde{\mathbb{R}},\varepsilon}^2(t) - \hat{a} q_\varepsilon^2(t) + \hat{b}(t) q_\varepsilon(t) + i\hbar_\varepsilon \alpha^2(t),$$

$$\hat{a} = a - \hbar_\varepsilon i,$$

$$\hat{b}(t) = b(t) - 2i\hbar_\varepsilon \alpha(t).$$

We rewrite $\mathcal{L}\left(\dot{q}_{\widetilde{\mathbb{R}},\varepsilon}(t), q_\varepsilon(t)\right)$ in the next equivalent form:

$$\mathcal{L}\left(\dot{q}_{\widetilde{\mathbb{R}},\varepsilon}(t), q_\varepsilon(t)\right) = \frac{m}{2} \dot{q}_{\widetilde{\mathbb{R}},\varepsilon}^2(t) - \frac{m\varpi^2}{2} q_\varepsilon^2(t) + \hat{b}(t) q_\varepsilon(t) + i\hbar_\varepsilon \alpha^2(t),$$

$$(2.1.36)$$

$$\frac{m\varpi^2}{2} = \hat{a}.$$

The Euler-Lagrange equation of $\mathcal{L}_\varepsilon(\dot{q}, q_\varepsilon)$ with corresponding boundary conditions is

$$\frac{d}{dt}\left(\frac{\partial \mathcal{L}_\varepsilon}{\partial \dot{q}_{\widetilde{\mathbb{R}},\varepsilon}}\right) - \frac{\partial \mathcal{L}_\varepsilon}{\partial q} = 0,$$

(2.1.37)

$$q(0) = y, q(T) = x.$$

Corresponding complex action $S(y, x, T)$ is:

$$S(y, x, T) = \frac{m\varpi}{2\sin\varpi T}\Bigg[ (\cos\varpi T)(y^2 + x^2) - 2xy + \frac{2x}{m\varpi}\int_0^T \widehat{b}(t)\sin(\varpi t)dt +$$

$$+ \frac{2y}{m\varpi}\int_0^T \widehat{b}(t)\sin\varpi(T-t)dt -$$

$$- \frac{2}{m^2\varpi^2}\int_0^T\int_0^t \widehat{b}(t)b(s)\sin\varpi(T-t)\sin(\varpi s)dsdt + i\hbar_\varepsilon\int_0^T \alpha^2(t)dt\Bigg] =$$

$$\frac{m\varpi}{2\sin\varpi T}[(\cos\varpi T)(y^2 + x^2) - 2xy + 2w_2x + 2w_1y + \Theta_T]$$

(2.1.38)

$$\Theta_T = -\frac{2}{m^2\varpi^2}\int_0^T\int_0^t \widehat{b}(t)b(s)\sin\varpi(T-t)\sin(\varpi s)dsdt + i\hbar_\varepsilon\int_0^T \alpha^2(t)dt,$$

$$w_1 = \frac{1}{m\varpi}\int_0^T \widehat{b}(t)\sin\varpi(T-t)dt,$$

$$w_2 = \frac{1}{m\varpi}\int_0^T \widehat{b}(t)\sin(\varpi t)dt$$

$$y = \frac{1}{\cos\varpi T}\left(x - \frac{1}{m\varpi}\int_0^T \widehat{b}(t)\sin\varpi(T-t)dt\right), w_T = \frac{1}{m\varpi}\int_0^T \widehat{b}(t)\sin\varpi(T-t)dt$$

$$y = \frac{1}{\cos\varpi T}(x - w_T), \frac{m\varpi}{2\sin\varpi T} = \alpha_T$$

$$S'(y,x,T) = i\alpha_T[(\cos\varpi T)(y^2+x^2) - 2xy + 2xw_T + 2yw_T] - \gamma_T^2 y^2$$

$$\frac{\partial S'}{\partial y} = i\alpha_T[2y\cos\varpi T - 2x + 2w_T] - 2y\gamma_T^2 = 0$$

$$2y\cos\varpi T - 2x + 2w_T - 2y\frac{\gamma_T^2}{i\alpha_T} = 0, 2y\cos\varpi T - 2x + 2w_T + 2y\frac{i\gamma_T^2}{\alpha_T} = 0$$

$$y\left[\cos\varpi T + \frac{i\gamma_T^2}{\alpha_T}\right] = x - w_T,$$

$$y = \frac{x}{\beta_T} - \frac{w_T}{\beta_T}, \quad \beta_T = \cos\varpi T + \frac{i\gamma_T^2}{\alpha_T},$$

$$S'(y,x,T) = i\alpha_T[(\cos\varpi T)(y^2+x^2) - 2xy + 2xw_T + 2yw_T] - \gamma_T^2 y^2 =$$

$$i\alpha_T\left[(\cos\varpi T)\left[\left(\frac{x}{\beta_T} - \frac{w_T}{\beta_T}\right)^2 + x^2\right] - 2x\left(\frac{x}{\beta_T} - \frac{w_T}{\beta_T}\right) + 2xw_T + 2w_T\left(\frac{x}{\beta_T} - \frac{w_T}{\beta_T}\right)\right] -$$
$$-\gamma_T^2\left(\frac{x}{\beta_T} - \frac{w_T}{\beta_T}\right)^2 =$$

$$i\alpha_T\cos\varpi T\left[\left(\frac{x}{\beta_T} - \frac{w_T}{\beta_T}\right)^2 + x^2\right] - 2xi\alpha_T\left(\frac{x}{\beta_T} - \frac{w_T}{\beta_T}\right) + 2xi\alpha_T w_T + 2i\alpha_T w_T\left(\frac{x}{\beta_T} - \frac{w_T}{\beta_T}\right)$$
$$-\gamma_T^2\left(\frac{x}{\beta_T} - \frac{w_T}{\beta_T}\right)^2 =$$

$$i\alpha_T\cos\varpi T\left(\frac{x}{\beta_T} - \frac{w_T}{\beta_T}\right)^2 + x^2 i\alpha_T\cos\varpi T - 2xi\alpha_T\left(\frac{x}{\beta_T} - \frac{w_T}{\beta_T}\right) + 2xi\alpha_T w_T + 2i\alpha_T w_T\left(\frac{x}{\beta_T} - \right.$$
$$-\gamma_T^2\left(\frac{x}{\beta_T} - \frac{w_T}{\beta_T}\right)^2 =$$

$$\left(\frac{x}{\beta_T} - \frac{w_T}{\beta_T}\right)^2(i\alpha_T\cos\varpi T - \gamma_T^2) + x^2 i\alpha_T\cos\varpi T - 2i\alpha_T w_T\left(\frac{x}{\beta_T} - \frac{w_T}{\beta_T}\right)(x - w_T) + 2xi\alpha_T w_T$$
$$\beta_T^{-2}(i\alpha_T\cos\varpi T - \gamma_T^2)(x - w_T)^2 + x^2 i\alpha_T\cos\varpi T - 2i\alpha_T w_T\beta_T^{-1}(x - w_T)^2 + 2xi\alpha_T w_T =$$
$$(x - w_T)^2\underbrace{[\beta_T^{-2}(i\alpha_T\cos\varpi T - \gamma_T^2) - 2i\alpha_T w_T\beta_T^{-1}]}_{\delta_T} + x^2 i\alpha_T\cos\varpi T + 2xi\alpha_T w_T =$$

$$(x - w_T)^2\delta_T + x^2 i\alpha_T\cos\varpi T + 2xi\alpha_T w_T = \delta_T(x^2 - 2xw_T + w_T^2) + x^2 i\alpha_T\cos\varpi T + 2xi\alpha_T w_T =$$

$$\delta_T x^2 - 2x\delta_T w_T + \delta_T w_T^2 + x^2 i\alpha_T\cos\varpi T + 2xi\alpha_T w_T =$$
$$S(x) = x^2(\delta_T + i\alpha_T\cos\varpi T) - 2xw_T(\delta_T - i\alpha_T) + \delta_T w_T^2 - \theta_T^2 x^2 =$$
$$x^2[(\delta_T + i\alpha_T\cos\varpi T) - \theta_T^2] - 2xw_T(\delta_T - i\alpha_T) + \delta_T w_T^2$$
$$\frac{\partial S(x)}{\partial x} = 2x(\delta_T + i\alpha_T\cos\varpi T) - 2w_T(\delta_T - i\alpha_T) - 2x\theta_T^2 = 0$$
$$x[(\delta_T + i\alpha_T\cos\varpi T) - \theta_T^2] = (\delta_T w_T - i\alpha_T w_T), x_{cr} = \frac{w_T(\delta_T - i\alpha_T)}{[(\delta_T + i\alpha_T\cos\varpi T) - \theta_T^2]}$$

$$S(x_{cr}) = \frac{w_T^2(\delta_T - i\alpha_T)^2}{\left[(\delta_T + i\alpha_T \cos \varpi T) - \theta_T^2\right]} - \frac{2w_T^2(\delta_T - i\alpha_T)^2}{\left[(\delta_T + i\alpha_T \cos \varpi T) - \theta_T^2\right]} = \frac{w_T^2(\delta_T - i\alpha_T)^2}{\left[(\delta_T + i\alpha_T \cos \varpi T) - \theta_T^2\right]}$$

Substitution Eq.(2.1.38) into Eq.(2.1.35) gives

$$\theta_\varepsilon^{(1)}(T, x, x_0, \alpha, \hbar_\varepsilon) =$$

$$\frac{1}{\sqrt[4]{2\pi\hbar_\varepsilon}} \int\limits_{\substack{q(T)=x \\ q(0)=y}} dy D^+[q(t), \varepsilon] \exp\left(-\frac{y^2}{2\hbar_\varepsilon}\right) \exp\left[-\frac{1}{\hbar_\varepsilon}(\gamma_T^2 q^2(0) + \theta_T^2 q^2(T))\right] \times$$

$$\exp\left\{\frac{i}{\hbar_\varepsilon} \int\limits_0^T \mathcal{L}\left(\dot{q}_{\widetilde{\mathbb{R}},\varepsilon}(t), q_\varepsilon(t)\right) dt\right\} =$$

$$= \widetilde{\theta}_\varepsilon^{(1)}(T, x, x_0, \hbar_\varepsilon).$$

(2.1.39)

Here

$$\widetilde{\theta}_{\varepsilon}^{(1)}(T,x,\alpha,\hbar_{\varepsilon}) =$$

$$\frac{1}{\sqrt[4]{2\pi\hbar_{\varepsilon}}}\sqrt{\frac{m\varpi}{2\pi i\hbar_{\varepsilon}\sin\varpi T}}\int\limits_{-\infty}^{\infty}dy\exp\left(-\frac{y^2}{2\hbar_{\varepsilon}}\right)\times$$

$$\exp\left\{\frac{1}{\hbar_{\varepsilon}}\left[i\frac{m\varpi}{2\sin\varpi T}S(y,x,T)-(\gamma_T^2 y^2+\theta_T^2 x^2)\right]\right\} =$$

$$\frac{1}{\sqrt[4]{2\pi\hbar_{\varepsilon}}}\sqrt{\frac{\alpha_T}{2\pi i\hbar_{\varepsilon}}}\int\limits_{-\infty}^{\infty}dy\exp\left\{\frac{1}{\hbar_{\varepsilon}}[i\alpha_T S(y,x,T)-(\overline{\gamma}_T^2 y^2+\theta_T^2 x^2)]\right\} =$$

$$\frac{1}{\sqrt[4]{2\pi\hbar_{\varepsilon}}}\sqrt{\frac{\alpha_T}{2\pi i\hbar_{\varepsilon}}}\int\limits_{-\infty}^{\infty}dy\exp\left\{\frac{1}{\hbar_{\varepsilon}}\widetilde{S}(y,x,T)\right\}, \qquad (2.1.40)$$

$$\widetilde{S}(y,x,T) = i\alpha_T S(y,x,T)-(\overline{\gamma}_T^2 y^2+\theta_T^2 x^2),$$

$$S(y,x,T) = (\cos\varpi T)(y^2+x^2)-2xy+2w_2 x+2w_1 y+\Theta_T,$$

$$\overline{\gamma}_T^2 = 2\gamma_T^2,$$

$$\alpha_T = \frac{m\varpi}{2\sin\varpi T}$$

$$\sin\varpi T \neq 0.$$

Let's calculate now integral $\widetilde{\theta}_{\varepsilon}^{(1)}(T,x,\alpha,\hbar_{\varepsilon})$ by using saddle point method.

$$\cos\varpi T+\frac{i\gamma_T^2}{\alpha_T}$$

$$\frac{\partial \widetilde{S}}{\partial y} = i\alpha_T[2y\cos\varpi T - 2x + 2w_1] - 2\overline{\gamma}_T^2 y = 0,$$

$$y(i\alpha_T\cos\varpi T - \overline{\gamma}_T^2) = i\alpha_T(x - w_1),$$

$$y\beta_T = i\alpha_T(x - w_1), \tag{2.1.41}$$

$$\beta_T = i\alpha_T\cos\varpi T - \overline{\gamma}_T^2.$$

Therefore saddle point $y_{cr}(x, T)$ is

$$y_{cr}(x, T) = \frac{i\alpha_T(x - w_1)}{i\alpha_T\cos\varpi T - \overline{\gamma}_T^2} = i\alpha_T\beta_T^{-1}(x - w_1),$$

$$\beta_T = i\alpha_T\cos\varpi T - \overline{\gamma}_T^2, \tag{2.1.42}$$

$$y_{cr}(x, T) = i\alpha_T\beta_T^{-1}\left(x - \frac{1}{m\varpi}\int_0^T \widehat{b}(t)\sin\varpi(T - t)dt\right),$$

$$\cos\varpi T \neq 0.$$

Hence

$$y_{cr}(x, T) = i\alpha_T\beta_T^{-1}(x - w_1)$$
$$y_{cr}^2(x, T) = -\alpha_T^2\beta_T^{-2}(x - w_1)^2$$

$$\widetilde{S}(y_{\mathbf{cr}}, x, T) = i\alpha_T S(y_{\mathbf{cr}}, x, T) - (\overline{\gamma}_T^2 y_{\mathbf{cr}}^2 + \theta_T^2 x^2) =$$

$$i\alpha_T[(\cos\varpi T)(y_{\mathbf{cr}}^2 + x^2) - 2xy_{\mathbf{cr}} + 2w_2 x + 2w_1 y_{\mathbf{cr}} + \Theta_T] - (\overline{\gamma}_T^2 y_{\mathbf{cr}}^2 + \theta_T^2 x^2) =$$

$$i\alpha_T(\cos\varpi T)y_{\mathbf{cr}}^2 + i\alpha_T(\cos\varpi T)x^2 - 2i\alpha_T xy_{\mathbf{cr}} + 2i\alpha_T w_2 x + 2i\alpha_T w_1 y_{\mathbf{cr}} + i\alpha_T \Theta_T -$$

$$-\overline{\gamma}_T^2 y_{\mathbf{cr}}^2 - \theta_T^2 x^2 =$$

$$(i\alpha_T \cos\varpi T - \overline{\gamma}_T^2)y_{\mathbf{cr}}^2 + (i\alpha_T \cos\varpi T - \theta_T^2)x^2 - 2i\alpha_T y_{\mathbf{cr}}(x - w_1) + 2i\alpha_T w_2 x + i\alpha_T \Theta_T =$$

$$-\alpha_T^2 \beta_T^{-2} \delta_T^{(1)}(x - w_1)^2 + 2\alpha_T^2 \beta_T^{-1}(x - w_1)^2 + \delta_T^{(2)} x^2 + 2i\alpha_T w_2 x + i\alpha_T \Theta_T = \qquad (2.1.)$$

$$(x - w_1)^2 \left(2\alpha_T^2 \beta_T^{-1} - \alpha_T^2 \beta_T^{-2} \delta_T^{(1)}\right) + \delta_T^{(2)} x^2 + 2i\alpha_T w_2 x + i\alpha_T \Theta_T =$$

$$(x - w_1)^2 \eta_T + \delta_T^{(2)} x^2 + 2i\alpha_T w_2 x + i\alpha_T \Theta_T,$$

$$\eta_T = 2\alpha_T^2 \beta_T^{-1} - \alpha_T^2 \beta_T^{-2} \delta_T^{(1)} = \alpha_T^2 \beta_T^{-1}\left(2 - \beta_T^{-1} \delta_T^{(1)}\right) = \alpha_T^2 \beta_T^{-1},$$

$$\delta_T^{(1)} = (i\alpha_T \cos\varpi T - \overline{\gamma}_T^2) = \beta_T,$$

$$\delta_T^{(2)} = i\alpha_T \cos\varpi T - \theta_T^2$$

Thus
$$\delta_T^{(2)} = i\alpha_T \cos\varpi T - \theta_T^2$$

$$\widetilde{S}(y_{\mathbf{cr}}, x, T) =$$

$$\eta_T(x - w_1)^2 + \delta_T^{(2)} x^2 + 2i\alpha_T w_2 x + i\alpha_T \Theta_T =$$

$$\eta_T x^2 - 2\eta_T w_1 x + \eta_T w_1^2 + \delta_T^{(2)} x^2 + 2i\alpha_T w_2 x + i\alpha_T \Theta_T =$$

$$\left(\eta_T + \delta_T^{(2)}\right) x^2 + 2(i\alpha_T w_2 - \eta_T w_1) x + \eta_T w_1^2 + i\alpha_T \Theta_T =$$

$$\vartheta_T x^2 + \sigma_T x + \eta_T w_1^2 + \Theta_T^{(1)},$$

$$\vartheta_T = \eta_T + \delta_T^{(2)} = \eta_T + i\alpha_T \cos \varpi T - \theta_T^2,$$

$$\sigma_T = 2(i\alpha_T w_2 - \eta_T w_1),$$

$$\Theta_T^{(1)} = i\alpha_T \Theta_T,$$

$$\eta_T = \alpha_T^2 \beta_T^{-1} = \frac{\alpha_T^2}{i\alpha_T \cos \varpi T - \overline{\gamma}_T^2}.$$

$$\eta_T = 2\alpha_T^2 \beta_T^{-1} - \alpha_T^2 \beta_T^{-2} \delta^{(1)} = \alpha_T^2 \beta_T^{-1} \left(2 - \beta_T^{-1} \delta^{(1)}\right) =$$

$$\alpha_T^2 \beta_T^{-1} (2 - \beta_T^{-1} \beta_T) = \alpha_T^2 \beta_T^{-1}.$$

(2.1. )

From Eq.(2.1.42) and Eq.(2.1.40) we obtain

$$\widetilde{\theta}_\varepsilon^{(1)}(T, x, \alpha, \hbar_\varepsilon) =$$

$$\frac{1}{\sqrt[4]{2\pi\hbar_\varepsilon}} \sqrt{\frac{\alpha_T}{2\pi i \widetilde{S}_{yy}(x)}} \exp\left\{\frac{1}{\hbar_\varepsilon} \widetilde{S}(y_{\mathbf{cr}}, x, T)\right\}.$$

(2.1. )

Let's consider integral

$$U(T,\alpha,\hbar_\varepsilon) = \int dx |x| \left| \widetilde{\vartheta}_\varepsilon^{(1)}(T,x,\alpha,\hbar_\varepsilon) \right|^2 =$$

$$\frac{1}{\sqrt{2\pi\hbar_\varepsilon}} \sqrt{\frac{|\alpha_T|}{2\pi|\widetilde{S}_{yy}(x)|}} \times$$

$$\int dx |x| \exp\left[ \frac{1}{\hbar_\varepsilon} \operatorname{Re}\left[ \widetilde{S}(y_{\mathbf{cr}}(x,T),x) \right] \right] = \tag{2.1.44}$$

$$\frac{1}{\sqrt{2\pi\hbar_\varepsilon}} \sqrt{\frac{|\alpha_T|}{2\pi|\widetilde{S}_{yy}(x)|}} \int dx |x| \exp\left[ \frac{1}{\hbar_\varepsilon}\left[ \vartheta'_T x^2 + \sigma'_T x + \eta'_T w_1'^2 + \Theta'_T \right] \right],$$

$$\vartheta'_T = \operatorname{Re}\vartheta_T, \sigma'_T = \operatorname{Re}\sigma_T, \Theta'_T = \operatorname{Re}\Theta_T^{(2)}.$$

$$\frac{\partial \widetilde{S}(y_{\mathbf{cr}}(x,T),x)}{\partial x} = 2\vartheta'_T x + \sigma'_T = 0,$$

$$\tag{2.1.44}$$

$$x_{\mathbf{cr}} = -\frac{\sigma'_T}{2\vartheta'_T}.$$

$$\widetilde{S}(y_{\mathbf{cr}}(x_{\mathbf{cr}},T),x_{\mathbf{cr}}) = \vartheta'_T x_{\mathbf{cr}}^2 + \sigma'_T x_{\mathbf{cr}} + \eta'_T w_1'^2 + \Theta'_T =$$

$$\frac{\sigma_T'^2}{4\vartheta'_T} - \frac{\sigma_T'^2}{2\vartheta'_T} + \eta'_T w_T'^2 + \Theta'_T = -\frac{\sigma_T'^2}{4\vartheta'_T} + \eta'_T w_1'^2 + \Theta'_T. \tag{2.1.}$$

Let us calculate Gaussian integral (2.1.44) using Laplase method.
From Eq.(2.2.1.25) and Eq.(2.2.1.43) we obtain

From Eq.(2.2.1.45) we obtain master inequelity for the Feynman-Colombeau path integral.

**Theorem 2.2.1.1**. (**Master Inequelity for the Feynman-Colombeau Path Integral**).

$$\int_{\mathbb{R}} |x| |\widetilde{\nu}_{1,\varepsilon}(T,x,x_0,\alpha,\hbar)|^2 \leq 4U(T,x_0).$$

(2.1.46)

# II.2.2. Master Equation for the Feynman-Colombeau path integral corresponding with a Schrödinger equation with Hermitian Hamiltonian.

Let us consider regularized propagator

$$(\mathbf{K}_{\varepsilon,\eta,g}(x,T|x_0,0))_\varepsilon =$$

$$\left( \overline{\int\limits_{\substack{q(T)=x \\ q(0)=y}} dy D^+[q(t),\varepsilon]\Psi_0(q(0))} \times \exp\left[ \frac{i}{\epsilon}\left( \int\limits_0^T dt\left[ \frac{m}{2}\left( \frac{d\widetilde{q}_\varepsilon(t)}{dt} \right)^2 - V(q_\varepsilon(t)) \right] \right) \right. \right.$$

$$\left. \left. -\frac{1}{\epsilon}\int\limits_0^T dt\left[ \{q(t)\}_\eta + \eta g(t,q(T),q(0),\lambda) - \eta\lambda \right]^2 \right] \right)_\varepsilon,$$

$$\widetilde{q}_\varepsilon(t) = \frac{q(t)}{1 + \varepsilon\left\| \dfrac{\dot{q}(t)}{q(t)} \right\|_{2,T}^2 + \varepsilon^2\left\| \dfrac{\ddot{q}(t)}{q(t)} \right\|_{2,T}^2}, \tag{2.2.2.1}$$

$$q_\varepsilon(t) = \frac{q(t)}{1 + \varepsilon^k q^m(t)}, k \geq 1, m \geq 2,$$

$$\Psi_0(x) \simeq \delta(x - x_0),$$

$$\epsilon(\varepsilon) << \eta << 1, \eta << m,$$

$$\epsilon(\varepsilon) << \varepsilon.$$

Here: (1) $g(0,q(T),q(0),\lambda) = q(0)$, and (2)

$$\{q(t)\}_\eta = \{\delta(t)\}_\eta + \eta u(t,T,x,y)$$
$$\textbf{iff}$$
$$q(t) = \delta(t) + u(0,T,x,y), \tag{2.2.2.2}$$

$$\delta(t) = \sum_{n=1}^\infty a_n \sin\left( \frac{n\pi t}{T} \right).$$

Here $u(t,x,y,\lambda)$ is any continuous function such that $u(0,x,y,\lambda) = y, u(T,x,y,\lambda) = x$ and

$$\langle \delta(t) \rangle_\eta = \sum_{n=1}^{N(\varepsilon)} \sqrt{\eta}\, a_n \sin\left(\frac{n\pi t}{T}\right) + \sum_{N(\varepsilon)}^{\infty} a_n \sin\left(\frac{n\pi t}{T}\right) \qquad (2.2.2.3)$$

Let us now consider quantum average $\langle T, x_0; \epsilon, \eta, \varepsilon, g \rangle$ :

$$\langle T, x_0; \epsilon, \eta, \varepsilon, g \rangle = \int dx x (\mathbf{K}_{\varepsilon,\eta,g}(x, T|x_0, 0) \mathbf{K}_{\varepsilon,\eta,g}^*(x, T|x_0, 0)) =$$

$$\int dx x |\mathbf{K}_{\varepsilon,\eta,g}(x, T|x_0, 0)|^2 \qquad (2.2.2.4)$$

From Eq.(2.2.2.1) we obtain

$$|\langle T, x_0; \epsilon, \eta, \varepsilon, g\rangle| \leq$$

$$\overline{\int} dx |x| \times$$

$$\left| \overline{\int_{\substack{q(T)=x \\ q(0)=y}}} dy D^+[q(t),\varepsilon]\Psi_0(q(0)) \times \right.$$

$$\exp\left[ \frac{i}{\epsilon}\left( \int\limits_0^T dt\left[ \frac{m}{2}\left(\frac{d\widetilde{q}_\varepsilon(t)}{dt}\right)^2 - V(q_\varepsilon(t)) \right] \right) \right.$$

$$\left. \left. -\frac{1}{\epsilon}\int\limits_0^T dt\Big[ \{q(t)\}_\eta + \eta g(t,q(T),q(0),\lambda) - \eta\lambda \Big]^2 \right] \right|^2 =$$

$$\overline{\int} dx |x| \times$$

$$\left\{ \left[ \overline{\int_{\substack{q(T)=x \\ q(0)=y}}} dy \widetilde{D}^+_\eta[q(t),\varepsilon]\Psi_0(q(0)) \cos\left[ \frac{1}{\epsilon}\left( \int\limits_0^T dt\left[ \frac{m}{2}\left(\frac{d\widetilde{q}_\varepsilon(t)}{dt}\right)^2 - V(q_\varepsilon(t)) \right] \right) \right] \right]^2 + \right. \qquad (2.2.2.5)$$

$$\left. \left[ \overline{\int_{\substack{q(T)=x \\ q(0)=y}}} dy \widetilde{D}^+_\eta[q(t),\varepsilon]\Psi_0(q(0)) \sin\left[ \frac{1}{\epsilon}\left( \int\limits_0^T dt\left[ \frac{m}{2}\left(\frac{d\widetilde{q}_\varepsilon(t)}{dt}\right)^2 - V(q_\varepsilon(t)) \right] \right) \right] \right]^2 \right\} =$$

$$\overline{\int} dx |x| \left[ \overline{\int_{\substack{q(T)=x \\ q(0)=y}}} dy \widetilde{D}^+_\eta[q(t),\varepsilon]\Psi_0(q(0)) \cos\left[ \frac{1}{\epsilon}\left( \int\limits_0^T dt\left[ \frac{m}{2}\left(\frac{dq_\varepsilon(t)}{dt}\right)^2 - V_\varepsilon(q(t)) \right] \right) \right] \right]^2 +$$

Here

$$\overline{\int\limits_{\substack{q(T)=x\\q(0)=y}}} \widetilde{D}_{\eta}^{+}[q(t),\varepsilon]\Psi_0(q(0))(\cdot) =$$

$$\overline{\int\limits_{\substack{q(T)=x\\q(0)=y}}} D^{+}[q(t),\varepsilon]\Psi_0(q(0)) \times \qquad (2.2.2.6)$$

$$\exp\left\{-\frac{1}{\epsilon}\int\limits_0^T dt\Big[\left\{q(t)\right\}_\eta + \eta g(t,q(T),q(0),\lambda) - \eta\lambda\Big]^2\right\} \times (\cdot)$$

and

$$\chi_1(T, \epsilon, \varepsilon, \eta) = \overline{\int} dx |x| \times$$

$$\left[ \overline{\int_{q(T)=x}} \widetilde{D}_\eta^+[q(t), \varepsilon] \Psi_0(q(0)) \times \right.$$

$$\left. \cos\left[ \frac{1}{\epsilon} \left( \int_0^T dt \left[ \frac{m}{2} \left( \frac{d\widetilde{q}_\varepsilon(t)}{dt} \right)^2 - V(q_\varepsilon(t)) \right] \right) \right] \right]^2 =$$

$$\overline{\int} dx \left\{ \left[ \overline{\int_{q(T)=x}} \widetilde{D}_\eta^+[q(t), \varepsilon] \Psi_0(q(0)) \sqrt{|q(T)|} \times \right. \right.$$

$$\left. \left. \cos\left[ \frac{1}{\epsilon} \left( \int_0^T dt \left[ \frac{m}{2} \left( \frac{d\widetilde{q}_\varepsilon(t)}{dt} \right)^2 - V(q_\varepsilon(t)) \right] \right) \right] \right]^2 \right\},$$

$$(2.2.2.7)$$

$$\chi_2(T, \epsilon, \varepsilon, \eta) = \overline{\int} dx |x| \times$$

$$\left[ \overline{\int_{q(T)=x}} \widetilde{D}_\eta^+[q(t), \varepsilon] \Psi_0(q(0)) \times \right.$$

$$\left. \sin\left[ \frac{1}{\epsilon} \left( \int_0^T dt \left[ \frac{m}{2} \left( \frac{d\widetilde{q}_\varepsilon(t)}{dt} \right)^2 - V(q_\varepsilon(t)) \right] \right) \right] \right]^2 =$$

$$\overline{\int} dx \left\{ \left[ \overline{\int_{q(T)=x}} \widetilde{D}_\eta^+[q(t), \varepsilon] \Psi_0(q(0)) \sqrt{|q(T)|} \times \right. \right.$$

$$\left. \left. \sin\left[ \frac{1}{\epsilon} \left( \int_0^T dt \left[ \frac{m}{2} \left( \frac{d\widetilde{q}_\varepsilon(t)}{dt} \right)^2 - V(q_\varepsilon(t)) \right] \right) \right] \right]^2 \right\}.$$

and here

$$\overline{\int_{q(T)=x} \widetilde{D}_\eta^+[q(t),\varepsilon]\Psi_0(q(0))(\cdot)} =$$

$$\int_{\substack{q(T)=x \\ q(0)=y}} dy D^+[q(t),\varepsilon]\Psi_0(q(0)) \times$$

$$\exp\left\{-\frac{1}{\epsilon}\int_0^T dt\Big[\{q(t)\}_\eta + \eta g(t,q(T),q(0),\lambda) - \eta\lambda\Big]^2\right\} \times (\cdot) \quad (2.2.2.8)$$

Let us now evaluate a quantity $(\chi_1(T,\epsilon,\varepsilon,\eta))_\varepsilon$ :

$$(\chi_1(T,\epsilon,\varepsilon,\eta))_\varepsilon =$$

$$\left(\bar{\int} dx|x| \times\right.$$

$$\left[\overline{\int_{q(T)=x} \widetilde{D}_\eta^+[q(t),\varepsilon]\Psi_0(q(0))} \times \right. \quad (2.2.2.9)$$

$$\left.\cos\left[\frac{1}{\epsilon}\left(\int_0^T dt\Big[\frac{m}{2}\left(\frac{d\widetilde{q}_\varepsilon(t)}{dt}\right)^2 - V(q_\varepsilon(t))\Big]\right)\right]\right]^2\right)_\varepsilon =$$

$$\left(\bar{\int} dx|x|\Pi_{1,\varepsilon}^2(T,x,x_0,\epsilon,\eta)\right)_\varepsilon.$$

Here

$$(\Pi_{1,\varepsilon}(T,x,x_0,\epsilon,\eta))_\varepsilon =$$

$$\left( \overline{\int\limits_{q(T)=x} \widetilde{D}_\eta[q(t),\varepsilon]} \Psi_0(q(0)) \cos\left[ \frac{1}{\epsilon}\left( \int\limits_0^T dt\left[ \frac{m}{2}\left( \frac{d\widetilde{q}_\varepsilon(t)}{dt} \right)^2 - V(q_\varepsilon(t)) \right] \right) \right] \right)_\varepsilon \qquad (2.2.2.10)$$

We assume now that: $(V(q_\varepsilon(t),t))_\varepsilon = (V_0(q_\varepsilon(t),t))_\varepsilon + (V_1(q_\varepsilon(t),t))_\varepsilon$. Here

$$(V_0(q_\varepsilon(t),t))_\varepsilon = (a_\varepsilon(t))_\varepsilon q_\varepsilon^2(t) + (b_\varepsilon(t))_\varepsilon (q_\varepsilon(t))_\varepsilon,$$

$$V_{1,\varepsilon}(q_\varepsilon(t),t) = (a_{3,\varepsilon}(t))_\varepsilon (q_\varepsilon^3(t))_\varepsilon + \ldots + (a_{k,\varepsilon}(t))_\varepsilon (q_\varepsilon^k(t))_\varepsilon,$$

$$(a_\varepsilon)_\varepsilon > 0, \qquad\qquad (2.2.2.11)$$

$$(q_\varepsilon(t))_\varepsilon = \left( \frac{q(t)}{1+\varepsilon^k q^m(t)} \right)_\varepsilon,$$

$$k \geq 1, m \geq 2$$

Substitution Eq.(2.2.2.11) into Eq.(2.2.2.9) gives

$$(\Pi_{1,\varepsilon}(T,x,x_0,\epsilon,\eta))_\varepsilon =$$

$$\overline{\int\limits_{q(T)=x}} \widetilde{D}_{\mathfrak{h}}^+[q(t),\varepsilon]\Psi_0(q(0)) \times$$

$$\cos\left[\frac{1}{\epsilon}\left(\int\limits_0^T dt\left[\frac{m}{2}\left(\frac{d\widetilde{q}_\varepsilon(t)}{dt}\right)^2 - V(q_\varepsilon(t),t)\right]\right)\right] =$$

$$\overline{\int\limits_{q(T)=x}} \widetilde{D}_{\mathfrak{h}}^+[q(t),\varepsilon]\Psi_0(q(0)) \times$$

$$\cos\left[\frac{1}{\epsilon}\left(\int\limits_0^T dt\left[\frac{m}{2}\left(\frac{d\widetilde{q}_\varepsilon(t)}{dt}\right)^2 - V_0(q_\varepsilon(t),t) - V_1(q_\varepsilon(t),t)\right]\right)\right] =$$

$$\overline{\int\limits_{q(T)=x}} \widetilde{D}_{\mathfrak{h}}^+[q(t),\varepsilon]\Psi_0(q(0))\rho_{q,\varepsilon}(\cdot) \times$$

$$\cos\left[\frac{1}{\epsilon}\left(\int\limits_0^T dt\left[\frac{m}{2}\left(\frac{d\widetilde{q}_\varepsilon(t)}{dt}\right)^2 - V_0(q_\varepsilon(t),t)\right]\right)\right] \times \qquad (2.2.2.12)$$

$$\cos\left[\frac{1}{\epsilon}\left(\int\limits_0^T V_1(q_\varepsilon(t),t)dt\right)\right] +$$

$$\overline{\int\limits_{q(T)=x}} \widetilde{D}_{\mathfrak{h}}^+[q(t),\varepsilon]\Psi_0(q(0)) \times$$

$$\sin\left[\frac{1}{\epsilon}\left(\int\limits_0^T dt\left[\frac{m}{2}\left(\frac{d\widetilde{q}_\varepsilon(t)}{dt}\right)^2 - V_0(q_\varepsilon(t),t)\right]\right)\right] \times$$

$$\sin\left[-\frac{1}{\epsilon}\left(\int\limits_0^T V_1(q_\varepsilon(t),t)dt\right)\right] =$$

Here

$$v_{1,\varepsilon}(T,x,x_0,\epsilon,\eta) =$$

$$\overline{\int\limits_{q(T)=x} \widetilde{D}_\eta^+[q(t),\varepsilon]} \Psi_0(q(0)) \times$$

$$\cos\left[\frac{1}{\epsilon}\left(\int\limits_0^T dt\left[\frac{m}{2}\left(\frac{d\widetilde{q}_\varepsilon(t)}{dt}\right)^2 - V_0(q_\varepsilon(t),t)\right]\right)\right] \times$$

$$\cos\left[\frac{1}{\epsilon}\left(\int\limits_0^T V_1(q_\varepsilon(t),t)dt\right)\right],$$

$$v_{2,\varepsilon}(T,x,x_0,\epsilon,\eta) =$$

$$\overline{\int\limits_{q(T)=x} \widetilde{D}_\eta^+[q(t),\varepsilon]} \Psi_0(q(0)) \times$$

$$\sin\left[\frac{1}{\epsilon}\left(\int\limits_0^T dt\left[\frac{m}{2}\left(\frac{d\widetilde{q}_\varepsilon(t)}{dt}\right)^2 - V_0(q_\varepsilon(t),t)\right]\right)\right] \times$$

$$\sin\left[-\frac{1}{\epsilon}\left(\int\limits_0^T V_1(q_\varepsilon(t),t)dt\right)\right].$$

$$(2.2.2.13)$$

Let us now evaluate a quantity $(\chi_1(T,\epsilon,\varepsilon,\eta))_\varepsilon$ (see Eq.(2.2.9)).By using replacement $q(t) \rightarrow \sqrt{\epsilon}\,q(t),$ from Eq.(2.2.9) we obtain

$$(\chi_1(T,\epsilon,\varepsilon,\eta))_\varepsilon =$$

$$\left(\overline{\int} dx |x| \times \right.$$

$$\left[ \overline{\int_{q(T)=\frac{x}{\sqrt{\epsilon}}} \widetilde{D}_\eta[q(t),\varepsilon\epsilon]\Psi_0(\sqrt{\epsilon}\,q(0))} \times \right. \tag{2.2.14}$$

$$\left. \left. \cos\left[ \frac{1}{\epsilon}\left( \int_0^T dt\left[ \frac{m}{2}\left( \frac{\sqrt{\epsilon}\,d\widetilde{q}_{\varepsilon\epsilon}(t)}{dt} \right)^2 - V(\sqrt{\epsilon}\,q_{v(\varepsilon,\epsilon)}(t)) \right] \right) \right]^2 \right)_\varepsilon = \right.$$

$$\left( \overline{\int} dx |x| \Pi_{1,\varepsilon}^2\left( T, \frac{x}{\sqrt{\epsilon}}, x_0, \epsilon, \eta \right) \right)_\varepsilon.$$

By using replacement $q(t) \to \sqrt{\epsilon}\,q(t),$ from Eq.(2.2.2.13) we obtain

$$\nu_{1,\varepsilon}\left(T, \frac{x}{\sqrt{\epsilon}}, x_0, \epsilon, \eta\right) =$$

$$\overline{\int\limits_{q(T)=\frac{x}{\sqrt{\epsilon}}} \widetilde{D}_\eta^+[q(t), \varepsilon, \epsilon]} \Psi_0(\sqrt{\epsilon}\, q(0)) \times$$

$$\cos\left[\frac{1}{\epsilon}\left(\int\limits_0^T dt\left[\frac{m}{2}\left(\frac{\sqrt{\epsilon}\, d\widetilde{q}_{\varepsilon\epsilon}(t)}{dt}\right)^2 - V_0(\sqrt{\epsilon}\, q_{\nu(\varepsilon,\epsilon)}(t), t)\right]\right)\right] \times$$

$$\cos\left[\frac{1}{\epsilon}\left(\int\limits_0^T V_1(\sqrt{\epsilon}\, q_{\nu(\varepsilon,\epsilon)}(t), t) dt\right)\right] = \qquad (2.2.2.15.a)$$

$$\overline{\int\limits_{q(T)=\frac{x}{\sqrt{\epsilon}}} \widetilde{D}_\eta^+[q(t), \varepsilon, \epsilon]} \Psi_0(\sqrt{\epsilon}\, q(0)) \times$$

$$\cos\left(\int\limits_0^T dt\left[\frac{m}{2}\left(\frac{d\widetilde{q}_{\varepsilon\epsilon}(t)}{dt}\right)^2 - a_\varepsilon(t) q_{\nu(\varepsilon,\epsilon)}^2(t) + b_\varepsilon(t)\frac{q_{\nu(\varepsilon,\epsilon)}(t)}{\sqrt{\epsilon}}\right]\right) \times$$

$$\cos\left(\int\limits_0^T V_1(\sqrt{\epsilon}\, q_{\nu(\varepsilon,\epsilon)}(t)) dt\right)$$

and

$$v_{2,\varepsilon}\left(T, \frac{x}{\sqrt{\epsilon}}, x_0, \epsilon, \eta\right) =$$

$$\overline{\int\limits_{q(T)=\frac{x}{\sqrt{\epsilon}}} \widetilde{D}_\eta^+[q(t), \varepsilon, \epsilon] \Psi_0(\sqrt{\epsilon}\, q(0)) \times}$$

$$(2.2.2.15.b)$$

$$\sin\left[\frac{1}{\epsilon}\left(\int\limits_0^T dt\left[\frac{m}{2}\left(\frac{\sqrt{\epsilon}\, d\widetilde{q}_{\varepsilon\epsilon}(t)}{dt}\right)^2 - V_0(\sqrt{\epsilon}\, q_{\nu(\varepsilon,\epsilon)}(t), t)\right]\right)\right] \times$$

$$\sin\left[-\frac{1}{\epsilon}\left(\int\limits_0^T V_1(\sqrt{\epsilon}\, q_{\nu(\varepsilon,\epsilon)}(t), t) dt\right)\right].$$

Here

$$\nu(\varepsilon, \epsilon) = \varepsilon^k \epsilon^{\frac{m}{2}},$$

$$\frac{1}{\epsilon} V_1(\sqrt{\epsilon}\, q_{\nu(\varepsilon,\epsilon)}(t)) = \frac{1}{\epsilon} V_1(\sqrt{\epsilon}\, q_{\nu(\varepsilon,\epsilon)}(t)) = a_3(t)\epsilon^{1/2} q_{\nu(\varepsilon,\epsilon)}^3(t) + \ldots$$

$$+ a_k(t)\epsilon^{(k/2)-1} q_{\nu(\varepsilon,\epsilon)}^k(t) = \widehat{V}_1(\sqrt{\epsilon}\, q_{\nu(\varepsilon,\epsilon)}(t)). \qquad (2.2.2.16)$$

$$\widehat{V}_1(\sqrt{\epsilon}\, q_{\nu(\varepsilon,\epsilon)}) = \sqrt{\epsilon}\, \widetilde{V}_1(q_{\nu(\varepsilon,\epsilon)}(t)),$$

$$\widetilde{V}_1(q_{\nu(\varepsilon,\epsilon)}(t)) = a_3(t) q_{\nu(\varepsilon,\epsilon)}^3(t) + \ldots + a_k(t)\epsilon^{(k/2)-3/2} q_{\nu(\varepsilon,\epsilon)}^k(t).$$

Let us rewrite Eqs.(2.2.15.a,b) in the following equivalent form

$$v_{1,\varepsilon}\left(T, \frac{x}{\sqrt{\epsilon}}, x_0, \epsilon, \eta\right) =$$

$$\overline{\int\limits_{q(T)=\frac{x}{\sqrt{\epsilon}}} \widetilde{D}_{\mathbf{sm}}^{(1)}[q(t),\varepsilon,\epsilon]\Psi_0(\sqrt{\epsilon}\,q(0)) \times}$$

$$\left|\cos\left(\int\limits_0^T dt\left[\frac{m}{2}\left(\frac{d\widetilde{q}_{\varepsilon\epsilon}(t)}{dt}\right)^2 - a_\varepsilon(t)q_{\nu(\varepsilon,\epsilon)}^2(t) + b_\varepsilon(t)\frac{q_{\nu(\varepsilon,\epsilon)}(t)}{\sqrt{\epsilon}}\right]\right)\right| \times$$

$$\cos\left(\int\limits_0^T \widehat{V}_1(q_{\nu(\varepsilon,\epsilon)}(t))dt\right),$$

$$(2.2.2.17)$$

$$v_{2,\varepsilon}\left(T, \frac{x}{\sqrt{\epsilon}}, x_0, \epsilon, \eta\right) =$$

$$\overline{\int\limits_{q(T)=\frac{x}{\sqrt{\epsilon}}} \widetilde{D}_{\mathbf{sm}}^{(2)}[q(t),\varepsilon,\epsilon]\Psi_0(\sqrt{\epsilon}\,q(0)) \times}$$

$$\left|\sin\left(\int\limits_0^T dt\left[\frac{m}{2}\left(\frac{d\widetilde{q}_{\varepsilon\epsilon}(t)}{dt}\right)^2 - a_\varepsilon(t)q_{\nu(\varepsilon,\epsilon)}^2(t) + b_\varepsilon(t)\frac{q_{\nu(\varepsilon,\epsilon)}(t)}{\sqrt{\epsilon}}\right]\right)\right| \times$$

$$\sin\left(-\int\limits_0^T \widehat{V}_1(q_{\nu(\varepsilon,\epsilon)}(t))dt\right).$$

Here $\widetilde{D}_{\mathbf{sm}}^{(1)}[q(t),\varepsilon,\epsilon]$ and $\widetilde{D}_{\mathbf{sm}}^{(2)}[q(t),\varepsilon,\epsilon]$ is the corresponding signed Feynman's "measures", i.e.

$$\overline{\int_{q(T)=\frac{x}{\sqrt{\epsilon}}} \widetilde{D}_{\mathbf{sm}}^{(1)}[q(t),\varepsilon,\epsilon]} \Psi_0(\sqrt{\epsilon}\,q(0)) \times (\cdot) =$$

$$\overline{\int_{q(T)=\frac{x}{\sqrt{\epsilon}}} \widetilde{D}_{\eta}^{+}[q(t),\varepsilon,\epsilon]} \Psi_0(\sqrt{\epsilon}\,q(0)) \times$$

$$[\cos^+(S_\varepsilon(T,x,q(t),\epsilon)) + \cos^-(S_\varepsilon(T,x,q(t),\epsilon))] \times (\cdot) =$$

$$\overline{\int_{\substack{q(T)=\frac{x}{\sqrt{\epsilon}} \\ q(0)=\frac{y}{\sqrt{\epsilon}}}} dy\, D^+[q(t),\varepsilon]} \Psi_0(\sqrt{\epsilon}\,q(0)) \times$$

$$[\cos^+(S_\varepsilon(T,x,q(t),\epsilon)) + \cos^-(S_\varepsilon(T,x,q(t),\epsilon))] \times$$

$$\exp\left\{-\frac{1}{\epsilon}\int_0^T dt\Big[\left\{\sqrt{\epsilon}\,q(t)\right\}_\eta + \eta g(t,\sqrt{\epsilon}\,q(T),\sqrt{\epsilon}\,q(0),\lambda) - \eta\lambda\Big]^2\right\} \times (\cdot) = \tag{2.2.2.18}$$

$$\overline{\int_{q(T)=\frac{x}{\sqrt{\epsilon}}} D^+[q(t),\varepsilon]} \Psi_0(q(0)) \times$$

$$[\cos^+(S_\varepsilon(T,x,q(t),\epsilon)) + \cos^-(S_{1,\varepsilon}(T,x,q(t),\epsilon))] \times$$

$$\exp\left\{-\frac{1}{\epsilon}\int_0^T dt\Big[\left\{\sqrt{\epsilon}\,q(t)\right\}_\eta - \eta g(t,\sqrt{\epsilon}\,q(T),\sqrt{\epsilon}\,q(0)) - \eta\lambda\Big]^2\right\} \times (\cdot),$$

$$S_\varepsilon(T,x,q(t),\epsilon) = \int_0^T dt\left[\frac{m}{2}\left(\frac{d\widetilde{q}_{\varepsilon\epsilon}(t)}{dt}\right)^2 - a_\varepsilon(t)q_\varepsilon^2(t) + b_\varepsilon(t)\frac{q_{\nu(\varepsilon,\epsilon)}(t)}{\sqrt{\epsilon}}\right]$$

and

$$\overline{\int\limits_{q(T)=\frac{x}{\sqrt{\epsilon}}}} \widetilde{D}_{sm}^{(2)}[q(t),\varepsilon,\epsilon]\Psi_0(\sqrt{\epsilon}\,q(0)) \times (\cdot) =$$

$$\overline{\int\limits_{q(T)=\frac{x}{\sqrt{\epsilon}}}} \widetilde{D}_{\eta}^{+}[q(t),\varepsilon]\Psi_0(q(0)) \times$$

$$[\sin^+(S_\varepsilon(T,x,q(t),\epsilon)) + \sin^-(S_\varepsilon(T,x,q(t),\epsilon))] \times (\cdot) =$$

$$\overline{\int\limits_{\substack{q(T)=\frac{x}{\sqrt{\epsilon}}\\q(0)=\frac{y}{\sqrt{\epsilon}}}}} dy D^+[q(t),\varepsilon]\Psi_0(\sqrt{\epsilon}\,q(0)) \times$$

$$[\sin^+(S_\varepsilon(T,x,q(t),\epsilon)) + \sin^-(S_\varepsilon(T,x,q(t),\epsilon)) \times$$

$$\exp\left\{-\frac{1}{\epsilon}\int\limits_0^T dt\Big[\,\{q(t)\}_\eta - \eta g(t,\sqrt{\epsilon}\,q(T),\sqrt{\epsilon}\,q(0)) - \eta\lambda\,\Big]^2\right\} \times (\cdot) =$$

$$\overline{\int\limits_{q(T)=\frac{x}{\sqrt{\epsilon}}}} D^+[q(t),\varepsilon]\Psi_0(\sqrt{\epsilon}\,q(0)) \times$$

$$[\sin^+(S_{1,\varepsilon}(T,x,q(t),\epsilon)) + \sin^-(S_\varepsilon(T,x,q(t),\epsilon))] \times$$

$$\exp\left\{-\frac{1}{\epsilon}\int\limits_0^T dt\Big[\,\{\sqrt{\epsilon}\,q(t)\}_\eta + \eta g(t,\sqrt{\epsilon}\,q(T),\sqrt{\epsilon}\,q(0),\lambda) - \eta\lambda\,\Big]^2\right\} \times (\cdot),$$

(2.2.2.19)

$$S_\varepsilon(T,x,q(t),\epsilon) = \int\limits_0^T dt\left[\,\frac{m}{2}\left(\frac{d\widetilde{q}_\varepsilon(t)}{dt}\right)^2 - a_\varepsilon(t)q_\varepsilon^2(t) + b_\varepsilon(t)\frac{q_{v(\varepsilon,\epsilon)}(t)}{\sqrt{\epsilon}}\,\right].$$

Let us rewrite Eq.(2.2.18)-Eq.(2.2.19) in the following equivalent form

$$\overline{\int\limits_{q(T)=\frac{x}{\sqrt{\epsilon}}}} \widetilde{D}_{\mathbf{sm}}^{(1)}[q(t),\varepsilon,\epsilon]\Psi_0(\sqrt{\epsilon}\,q(0)) \times (\cdot) =$$

$$\overline{\int\limits_{q(T)=\frac{x}{\sqrt{\epsilon}}}} \widetilde{D}_{\eta}^{+}[q(t),\varepsilon,\epsilon]\Psi_0(\sqrt{\epsilon}\,q(0)) \times$$

$$[\cos^+(S_\varepsilon(T,x,q(t),\epsilon)) + \cos^-(S_\varepsilon(T,x,q(t),\epsilon)) + 1 - 1] \times (\cdot) =$$

$$\overline{\int\limits_{q(T)=\frac{x}{\sqrt{\epsilon}}}} \widetilde{D}_{\mathbf{m}}^{(1)}[q(t),\varepsilon]\Psi_0(\sqrt{\epsilon}\,q(0)) \times (\cdot) - \overline{\int\limits_{q(T)=\frac{x}{\sqrt{\epsilon}}}} \widetilde{D}_{\eta}^{+}[q(t),\varepsilon]\Psi_0(\sqrt{\epsilon}\,q(0)) \times (\cdot) = \quad (2.2.20)$$

$$\overline{\int\limits_{q(T)=\frac{x}{\sqrt{\epsilon}}}} \widetilde{D}_{\mathbf{m}}^{(1)}[q(t),\varepsilon,\epsilon]\Psi_0(q(0)) \times (\cdot) - \overline{\int\limits_{q(T)=\frac{x}{\sqrt{\epsilon}}}} D^{+}[q(t),\varepsilon]\Psi_0(q(0)) \times$$

$$\exp\left\{-\int\limits_0^T dt\Big[\left\{\sqrt{\epsilon}\,q(t)\right\}_\eta + \eta g(t,q(T),q(0),\lambda) - \eta\lambda\Big]^2\right\} \times (\cdot)$$

and

$$\overline{\int\limits_{q(T)=\frac{x}{\sqrt{\epsilon}}}}\widetilde{D}_{\mathbf{sm}}^{(2)}[q(t),\varepsilon,\epsilon]\Psi_0(\sqrt{\epsilon}\,q(0))\times(\cdot)=$$

$$\overline{\int\limits_{q(T)=\frac{x}{\sqrt{\epsilon}}}}\widetilde{D}_{\eta}^{+}[q(t),\varepsilon,\epsilon]\Psi_0(\sqrt{\epsilon}\,q(0))\times$$

$$[\sin^{+}(S_{\varepsilon}(T,x,q(t),\epsilon))+\sin^{-}(S_{\varepsilon}(T,x,q(t),\epsilon))+1-1]\times(\cdot)=$$

$$\overline{\int\limits_{q(T)=\frac{x}{\sqrt{\epsilon}}}}\widetilde{D}_{\mathbf{m}}^{(2)}[q(t),\varepsilon,\epsilon]\Psi_0(\sqrt{\epsilon}\,q(0))\times(\cdot)$$

$$\tag{2.2.21}$$

$$-\overline{\int\limits_{q(T)=\frac{x}{\sqrt{\epsilon}}}}\widetilde{D}_{\eta}^{+}[q(t),\varepsilon,\epsilon]\Psi_0(\sqrt{\epsilon}\,q(0))\times(\cdot)=$$

$$\overline{\int\limits_{q(T)=\frac{x}{\sqrt{\epsilon}}}}\widetilde{D}_{\mathbf{m}}^{(2)}[q(t),\varepsilon]\Psi_0(\sqrt{\epsilon}\,q(0))\times(\cdot)-\overline{\int\limits_{q(T)=\frac{x}{\sqrt{\epsilon}}}}D^{+}[q(t),\varepsilon]\Psi_0(\sqrt{\epsilon}\,q(0))\times$$

$$\exp\left\{-\frac{1}{\epsilon}\int\limits_{0}^{T}dt\Big[\left\{\sqrt{\epsilon}\,q(t)\right\}_{\eta}+\eta g(t,\sqrt{\epsilon}\,q(T),\sqrt{\epsilon}\,q(0),\lambda)-\eta\lambda\Big]^2\right\}\times(\cdot).$$

Here $\widetilde{D}_{\mathbf{m}}^{(1)}$ and $\widetilde{D}_{\mathbf{m}}^{(2)}$ is the corresponding positive Feynman submeasures, i.e.

$$\overline{\int\limits_{q(T)=\frac{x}{\sqrt{\epsilon}}} \widetilde{D}_{\mathbf{m}}^{(1)}[q(t),\varepsilon,\epsilon]} \Psi_0(\sqrt{\epsilon}\,q(0)) \times (\cdot) =$$

$$\overline{\int\limits_{q(T)=\frac{x}{\sqrt{\epsilon}}} D^+[q(t),\varepsilon,\epsilon]} \Psi_0(\sqrt{\epsilon}\,q(0)) \times$$

(2.2.22)

$$[\cos^+(S_\varepsilon(T,x,q(t),\epsilon)) + \cos^-(S_\varepsilon(T,x,q(t),\epsilon)) + 1] \times$$

$$\exp\left\{-\frac{1}{\epsilon}\int\limits_0^T dt\Big[\{\sqrt{\epsilon}\,q(t)\}_\eta + \eta g(t,\sqrt{\epsilon}\,q(T),\sqrt{\epsilon}\,q(0),\lambda) - \eta\lambda\Big]^2\right\} \times (\cdot).$$

and

$$\overline{\int\limits_{q(T)=\frac{x}{\sqrt{\epsilon}}} \widetilde{D}_{\mathbf{m}}^{(2)}[q(t),\varepsilon,\epsilon]} \Psi_0(\sqrt{\epsilon}\,q(0)) \times (\cdot) =$$

$$\overline{\int\limits_{q(T)=\frac{x}{\sqrt{\epsilon}}} D^+[q(t),\varepsilon]} \Psi_0(\sqrt{\epsilon}\,q(0)) \times$$

(2.2.2.23)

$$[\sin^+(S_\varepsilon(T,x,q(t),\epsilon)) + \sin^-(S_{1,\varepsilon}(T,x,q(t),\epsilon)) + 1] \times$$

$$\exp\left\{-\frac{1}{\epsilon}\int\limits_0^T dt\Big[\{\sqrt{\epsilon}\,q(t)\}_\eta - \eta g(t,\sqrt{\epsilon}\,q(T),\sqrt{\epsilon}\,q(0),\lambda) - \eta\lambda\Big]^2\right\} \times (\cdot).$$

From Eq.(2.2.2.18)-Eq.(2.2.2.23) we obtain

$$v_{1,\varepsilon}\left(T, \frac{x}{\sqrt{\epsilon}}, x_0, \epsilon, \eta\right) =$$

$$\overline{\int\limits_{q(T)=\frac{x}{\sqrt{\epsilon}}} \widetilde{D}_{\mathbf{sm}}^{(1)}[q(t), \varepsilon, \epsilon] \Psi_0(\sqrt{\epsilon}\, q(0))} \times$$

$$\left| \cos\left(\int\limits_0^T dt \left[ \frac{m}{2} \left(\frac{d\widetilde{q}_{\varepsilon\varepsilon}(t)}{dt}\right)^2 - a_\varepsilon(t) q_{v(\varepsilon,\epsilon)}^2(t) + b_\varepsilon(t) \frac{q_{v(\varepsilon,\epsilon)}(t)}{\sqrt{\epsilon}} \right] \right) \right| \times$$

$$\cos\left(\int\limits_0^T \widehat{V}_1(\sqrt{\epsilon}\, q_{v(\varepsilon,\epsilon)}(t)) dt\right) =$$

$$\overline{\int\limits_{q(T)=\frac{x}{\sqrt{\epsilon}}} \widetilde{D}_{\mathbf{m}}^{(1)}[q(t), \varepsilon, \epsilon] \Psi_0(\sqrt{\epsilon}\, q(0))} \times$$

$$\left| \cos\left(\int\limits_0^T dt \left[ \frac{m}{2} \left(\frac{d\widetilde{q}_{\varepsilon\varepsilon}(t)}{dt}\right)^2 - a_\varepsilon(t) q_{v(\varepsilon,\epsilon)}^2(t) + b_\varepsilon(t) \frac{q_{v(\varepsilon,\epsilon)}(t)}{\sqrt{\epsilon}} \right] \right) \right| \times$$

$$\cos\left(\int\limits_0^T \widehat{V}_1(\sqrt{\epsilon}\, q_{v(\varepsilon,\epsilon)}(t)) dt\right) -$$

$$- \overline{\int\limits_{q(T)=\frac{x}{\sqrt{\epsilon}}} \widetilde{D}_\eta^+[q(t), \varepsilon\epsilon] \Psi_0(\sqrt{\epsilon}\, q(0))} \times$$

$$\left| \cos\left(\int\limits_0^T dt \left[ \frac{m}{2} \left(\frac{d\widetilde{q}_{\varepsilon\varepsilon}(t)}{dt}\right)^2 - a_\varepsilon(t) q_{v(\varepsilon,\epsilon)}^2(t) + b_\varepsilon(t) \frac{q_{v(\varepsilon,\epsilon)}(t)}{\sqrt{\epsilon}} \right] \right) \right| \times$$

$$\cos\left(\int\limits_0^T \widehat{V}_1(\sqrt{\epsilon}\, q_{v(\varepsilon,\epsilon)}(t)) dt\right)$$

(2.2.2.24)

and

$$v_{2,\varepsilon}\left(T, \frac{x}{\sqrt{\epsilon}}, x_0, \epsilon, \eta\right) =$$

$$\overline{\int\limits_{q(T)=\frac{x}{\sqrt{\epsilon}}} \widetilde{D}_{\mathbf{sm}}^{(2)}[q(t), \varepsilon, \epsilon] \Psi_0(\sqrt{\epsilon}\, q(0)) \times}$$

$$\left| \sin\left( \int\limits_0^T dt \left[ \frac{m}{2}\left( \frac{d\widetilde{q}_{\varepsilon\epsilon}(t)}{dt} \right)^2 - a_\varepsilon(t) q_{\nu(\varepsilon,\epsilon)}^2(t) + b_\varepsilon(t) \frac{q_{\nu(\varepsilon,\epsilon)}(t)}{\sqrt{\epsilon}} \right] \right) \right| \times$$

$$\cos\left( \int\limits_0^T \widehat{V}_1(\sqrt{\epsilon}\, q_{\nu(\varepsilon,\epsilon)}(t)) dt \right) =$$

$$\overline{\int\limits_{q(T)=\frac{x}{\sqrt{\epsilon}}} \widetilde{D}_{\mathbf{m}}^{(2)}[q(t), \varepsilon, \epsilon] \Psi_0(\sqrt{\epsilon}\, q(0)) \times}$$

(2.2.2.25)

$$\left| \sin\left( \int\limits_0^T dt \left[ \frac{m}{2}\left( \frac{d\widetilde{q}_{\varepsilon\epsilon}(t)}{dt} \right)^2 - a_\varepsilon(t) q_{\nu(\varepsilon,\epsilon)}^2(t) + b_\varepsilon(t) \frac{q_{\nu(\varepsilon,\epsilon)}(t)}{\sqrt{\epsilon}} \right] \right) \right| \times$$

$$\sin\left( -\int\limits_0^T \widehat{V}_1(\sqrt{\epsilon}\, q_{\nu(\varepsilon,\epsilon)}(t)) dt \right) -$$

$$-\overline{\int\limits_{q(T)=\frac{x}{\sqrt{\epsilon}}} \widetilde{D}_\eta^+[q(t), \varepsilon\epsilon] \Psi_0(\sqrt{\epsilon}\, q(0)) \times}$$

$$\left| \sin\left( \int\limits_0^T dt \left[ \frac{m}{2}\left( \frac{d\widetilde{q}_{\varepsilon\epsilon}(t)}{dt} \right)^2 - a_\varepsilon(t) q_{\nu(\varepsilon,\epsilon)}^2(t) + b_\varepsilon(t) \frac{q_{\nu(\varepsilon,\epsilon)}(t)}{\sqrt{\epsilon}} \right] \right) \right| \times$$

$$\sin\left( \int\limits_0^T \widehat{V}_1(\sqrt{\epsilon}\, q_{\nu(\varepsilon,\epsilon)}(t)) dt \right).$$

From Eq.(2.2.2.25) and Eq.(2.2.2.26) we obtain

$$v_{1,\varepsilon}\left(T,\frac{x}{\sqrt{\epsilon}},x_0,\epsilon,\eta\right) = \widetilde{v}_{1,\varepsilon}\left(T,\frac{x}{\sqrt{\epsilon}},x_0,\epsilon,\eta\right) - v^*_{1,\varepsilon}\left(T,\frac{x}{\sqrt{\epsilon}},x_0,\epsilon,\eta\right),$$

$$\widetilde{v}_{1,\varepsilon}\left(T,\frac{x}{\sqrt{\epsilon}},x_0,\epsilon,\eta\right) = \overline{\int_{q(T)=\frac{x}{\sqrt{\epsilon}}} \widetilde{D}^{(1)}_{\mathbf{m}}[q(t),\varepsilon,\epsilon]\Psi_0(\sqrt{\epsilon}\,q(0)) \times}$$

$$\left|\cos\left(\int_0^T dt\left[\frac{m}{2}\left(\frac{d\widetilde{q}_{\varepsilon\epsilon}(t)}{dt}\right)^2 - a_\varepsilon(t)q^2_{v(\varepsilon,\epsilon)}(t) + b_\varepsilon(t)\frac{q_{v(\varepsilon,\epsilon)}(t)}{\sqrt{\epsilon}}\right]\right)\right| \times$$

$$\cos\left(\int_0^T \widehat{V}_1(\sqrt{\epsilon}\,q_{v(\varepsilon,\epsilon)}(t))dt\right), \tag{2.2.2.26}$$

$$v^*_{1,\varepsilon}\left(T,\frac{x}{\sqrt{\epsilon}},x_0,\epsilon,\eta\right) = \overline{\int_{q(T)=\frac{x}{\sqrt{\epsilon}}} \widetilde{D}^+_{\eta}[q(t),\varepsilon,\epsilon]\Psi_0(\sqrt{\epsilon}\,q(0)) \times}$$

$$\left|\cos\left(\int_0^T dt\left[\frac{m}{2}\left(\frac{d\widetilde{q}_{\varepsilon\epsilon}(t)}{dt}\right)^2 - a_\varepsilon(t)q^2_{v(\varepsilon,\epsilon)}(t) + b_\varepsilon(t)\frac{q_{v(\varepsilon,\epsilon)}(t)}{\sqrt{\epsilon}}\right]\right)\right| \times$$

$$\cos\left(\int_0^T \widehat{V}_1(\sqrt{\epsilon}\,q_{v(\varepsilon,\epsilon)}(t))dt\right)$$

and

$$v_{2,\varepsilon}\left(T, \frac{x}{\sqrt{\epsilon}}, x_0, \epsilon, \eta\right) = \widetilde{v}_{2,\varepsilon}\left(T, \frac{x}{\sqrt{\epsilon}}, x_0, \epsilon, \eta\right) - v_{2,\varepsilon}^*\left(T, \frac{x}{\sqrt{\epsilon}}, x_0, \epsilon, \eta\right),$$

$$\widetilde{v}_{2,\varepsilon}\left(T, \frac{x}{\sqrt{\epsilon}}, x_0, \epsilon, \eta\right) = \overline{\int\limits_{q(T)=\frac{x}{\sqrt{\epsilon}}} \widetilde{D}_{\mathbf{m}}^{(2)}[q(t), \varepsilon\epsilon]} \Psi_0(\sqrt{\epsilon}\, q(0)) \times$$

$$\left| \sin\left( \int\limits_0^T dt \left[ \frac{m}{2}\left( \frac{d\widetilde{q}_{\varepsilon\epsilon}(t)}{dt} \right)^2 - a_\varepsilon(t) q_{v(\varepsilon,\epsilon)}^2(t) + b_\varepsilon(t) \frac{q_{v(\varepsilon,\epsilon)}(t)}{\sqrt{\epsilon}} \right] \right) \right| \times$$

$$\sin\left( -\int\limits_0^T V_{1,\epsilon}(\sqrt{\epsilon}\, q_{v(\varepsilon,\epsilon)}(t)) dt \right),$$

(2.2.2.27)

$$v_{2,\varepsilon}^*\left(T, \frac{x}{\sqrt{\epsilon}}, x_0, \epsilon, \eta\right) = \overline{\int\limits_{q(T)=\frac{x}{\sqrt{\epsilon}}} \widetilde{D}_\eta^+[q(t), \varepsilon, \epsilon]} \Psi_0(\sqrt{\epsilon}\, q(0)) \times$$

$$\left| \sin\left( \int\limits_0^T dt \left[ \frac{m}{2}\left( \frac{d\widetilde{q}_{\varepsilon\epsilon}(t)}{dt} \right)^2 - a_\varepsilon(t) q_{v(\varepsilon,\epsilon)}^2(t) + b_\varepsilon(t) \frac{q_{v(\varepsilon,\epsilon)}(t)}{\sqrt{\epsilon}} \right] \right) \right| \times$$

$$\sin\left( -\int\limits_0^T \widehat{V}_1(\sqrt{\epsilon}\, q_{v(\varepsilon,\epsilon)}(t)) dt \right).$$

Let us avaluate now path integral $v_{1,\varepsilon}^*\left(T, \frac{x}{\sqrt{\epsilon}}, x_0, \epsilon, \eta\right)$ for the case

$g(t, q(T), q(0), \lambda) \equiv 0, \ \lambda = 0$ :

$$v_{1,\varepsilon}^*\left(T, \frac{x}{\sqrt{\epsilon}}, x_0, \epsilon, \eta\right) = \overline{\int\limits_{\substack{q(T)=\frac{x}{\sqrt{\epsilon}} \\ q(0)=\frac{y}{\sqrt{\epsilon}}}} dy D^+[q(t), \varepsilon, \epsilon] \Psi_0(\sqrt{\epsilon}\, q(0)) \Xi_1(q(t), \varepsilon, \epsilon) \times}$$

$$\cos\left(\int\limits_0^T \widehat{V}_1(\sqrt{\epsilon}\, q_{\nu(\varepsilon,\epsilon)}(t)) dt\right) \exp\left\{-\frac{1}{\epsilon}\int\limits_0^T dt\Big[\big\{\sqrt{\epsilon}\, q(t)\big\}_\eta\Big]^2\right\} \qquad (2.2.2.28)$$

$$\Xi_1(q(t), \varepsilon, \epsilon) =$$

$$\left|\cos\left(\int\limits_0^T dt\left[\frac{m}{2}\left(\frac{d\widetilde{q}_{\varepsilon\epsilon}(t)}{dt}\right)^2 - a_\varepsilon(t) q_{\nu(\varepsilon,\epsilon)}^2(t) + b_\varepsilon(t)\frac{q_{\nu(\varepsilon,\epsilon)}(t)}{\sqrt{\epsilon}}\right]\right)\right|.$$

From Eqs.(2.2.2.28) we obtain

$$v_{1,\varepsilon}^*\left(T, \frac{x}{\sqrt{\epsilon}}, x_0, \epsilon, \eta\right) = \overline{\int\limits_{\substack{q(T)=\frac{x}{\sqrt{\epsilon}} \\ q(0)=\frac{y}{\sqrt{\epsilon}}}} dy D^+[q(t),\varepsilon,\epsilon] \Psi_0(\sqrt{\epsilon}\, q(0)) \Xi_1(q(t),\varepsilon,\epsilon) \times}$$

$$\cos\left(\int\limits_0^T \sqrt{\epsilon}\, \widetilde{V}_1(q_{v(\varepsilon,\epsilon)}(t)) dt\right) \exp\left\{-\int\limits_0^T dt \left[\left\{q(t)\right\}_\eta\right]^2\right\} =$$

$$\overline{\int\limits_{\substack{q(T)=\frac{x}{\sqrt{\epsilon}} \\ q(0)=\frac{y}{\sqrt{\epsilon}}}} dy D^+[q(t),\varepsilon,\epsilon] \Psi_0(\sqrt{\epsilon}\, q(0)) \Xi_1(q(t),\varepsilon,\epsilon) \times}$$

$$\exp\left\{-\int\limits_0^T dt \left[\left\{q(t)\right\}_\eta\right]^2\right\} \left[\sum_{n=0}^M \frac{(-1)^n \epsilon^n}{(2n)!}\left(\int\limits_0^T \widetilde{V}_1(q_{v(\varepsilon,\epsilon)}(t)) dt\right)^{2n} + \right.$$

$$\left. \epsilon^{n+1} O\left(\left(\int\limits_0^T \widetilde{V}_1(q_{v(\varepsilon,\epsilon)}(t)) dt\right)^{2n+1}\right)\right] = \qquad (2.2.2.29)$$

$$\overline{\int\limits_{\substack{q(T)=\frac{x}{\sqrt{\epsilon}} \\ q(0)=\frac{y}{\sqrt{\epsilon}}}} dy D^+[q(t),\varepsilon,\epsilon] \Psi_0(y) \Xi_1(q(t),\varepsilon,\epsilon) \exp\left\{-\int\limits_0^T dt \left[\left\{q(t)\right\}_\eta\right]^2\right\}} +$$

$$\sum_{n=1}^M \frac{(-1)^n \epsilon^n}{(2n)!}\left[\overline{\int\limits_{\substack{q(T)=\frac{x}{\sqrt{\epsilon}} \\ q(0)=\frac{y}{\sqrt{\epsilon}}}} dy D^+[q(t),\varepsilon,\epsilon] \Psi_0(y) \Xi_1(q(t),\varepsilon,\epsilon) \times}\right.$$

$$\left.\exp\left\{-\int\limits_0^T dt \left[\left\{q(t)\right\}_\eta\right]^2\right\}\left(\int\limits_0^T \widetilde{V}_1(q_{v(\varepsilon,\epsilon)}(t)) dt\right)^{2n}\right]$$

Let us avaluate now path integral $v_{2,\varepsilon}^*\left(T, \dfrac{x}{\sqrt{\epsilon}}, x_0, \epsilon, \eta\right)$ for the case:

$g(t, q(T), q(0), \lambda) \equiv 0, \ \lambda = 0 :$

$$\widetilde{v}_{2,\varepsilon}^*\left(T, \frac{x}{\sqrt{\epsilon}}, x_0, \epsilon, \eta\right) = \overline{\int\limits_{\substack{q(T)=\frac{x}{\sqrt{\epsilon}} \\ q(0)=\frac{y}{\sqrt{\epsilon}}}} dy D^+[q(t), \varepsilon, \epsilon] \Psi_0(\sqrt{\epsilon}\, q(0)) \Xi_2(q(t), \varepsilon, \epsilon) \times}$$

$$\sin\left(-\int\limits_0^T \widehat{V}_1(\sqrt{\epsilon}\, q_{v(\varepsilon,\epsilon)}(t)) dt\right) \exp\left\{-\int\limits_0^T dt \left[\left\{q(t)\right\}_\eta\right]^2\right\} \qquad (2.2.2.30)$$

$$\Xi_2(q(t), \varepsilon, \epsilon) =$$

$$\left| \sin\left(\int\limits_0^T dt \left[\frac{m}{2}\left(\frac{d\widetilde{q}_{\varepsilon\epsilon}(t)}{dt}\right)^2 - a_\varepsilon(t) q_{v(\varepsilon,\epsilon)}^2(t) + b_\varepsilon(t) \frac{q_{v(\varepsilon,\epsilon)}(t)}{\sqrt{\epsilon}}\right]\right)\right| \times$$

From Eq.(2.2.2.30) we obtain

$$v_{2,\varepsilon}^{*}\left(T,\frac{x}{\sqrt{\epsilon}},x_0,\epsilon,\eta\right)=\overline{\int\limits_{\substack{q(T)=\frac{x}{\sqrt{\epsilon}}\\q(0)=\frac{y}{\sqrt{\epsilon}}}}dyD^{+}[q(t),\varepsilon,\epsilon]\Psi_0(\sqrt{\epsilon}\,q(0))\Xi_2(q(t),\varepsilon,\epsilon)}\times$$

$$\sin\left(-\int\limits_0^T\sqrt{\epsilon}\,\widetilde{V}_1(q_{v(\varepsilon,\epsilon)}(t))dt\right)\exp\left\{-\int\limits_0^T dt\Big[\,\{q(t)\}_\eta\,\Big]^2\right\}=$$

$$\overline{\int\limits_{\substack{q(T)=\frac{x}{\sqrt{\epsilon}}\\q(0)=\frac{y}{\sqrt{\epsilon}}}}dyD^{+}[q(t),\varepsilon,\epsilon]\Psi_0(\sqrt{\epsilon}\,q(0))\Xi_2(q(t),\varepsilon,\epsilon)}\times$$

$$\exp\left\{-\int\limits_0^T dt\Big[\,\{q(t)\}_\eta\,\Big]^2\right\}\left[\sum_{n=0}^M\frac{(-1)^n\epsilon^{n+1/2}}{(2n+1)!}\left(\int\limits_0^T\widetilde{V}_1(q_{v(\varepsilon,\epsilon)}(t))dt\right)^{2n+1}+\right.$$

$$\left.\epsilon^{n+1}O\left(\left(\int\limits_0^T\widetilde{V}_1(q_{v(\varepsilon,\epsilon)}(t))dt\right)^{2n+2}\right)\right]= \qquad (2.2.2.31)$$

$$\overline{\int\limits_{\substack{q(T)=\frac{x}{\sqrt{\epsilon}}\\q(0)=\frac{y}{\sqrt{\epsilon}}}}dyD^{+}[q(t),\varepsilon,\epsilon]\Psi_0(y)\Xi_2(q(t),\varepsilon,\epsilon)}\exp\left\{-\int\limits_0^T dt\Big[\,\{q(t)\}_\eta\,\Big]^2\right\}+$$

$$\sum_{n=1}^M\frac{(-1)^n\epsilon^n}{(2n)!}\left[\overline{\int\limits_{\substack{q(T)=\frac{x}{\sqrt{\epsilon}}\\q(0)=\frac{y}{\sqrt{\epsilon}}}}dyD^{+}[q(t),\varepsilon,\epsilon]\Psi_0(y)\Xi_2(q(t),\varepsilon,\epsilon)}\times\right.$$

$$\left.\exp\left\{-\int\limits_0^T dt\Big[\,\{q(t)\}_\eta\,\Big]^2\right\}\left(\int\limits_0^T\widetilde{V}_1(q_{v(\varepsilon,\epsilon)}(t))dt\right)^{2n}\right].$$

**Lemma**.2.2.2.1.(**Generalized Hölder's inequality for the Feynman-Colombeau path integral**) Let $\frac{1}{p_1} + \frac{1}{p_2} = 1$ with $p_1, p_2 > 1$. Assume that

$$\left( F_{3,\varepsilon}\left[ \frac{d\widetilde{q}_\varepsilon(t)}{dt}, q_\varepsilon(t), t_i, t_f \right] \right)_\varepsilon > 0,$$

$$\left( \overline{\int_{\substack{q(t_f)=x \\ q(t_i)=y}}} \left| F_{1,\varepsilon}\left[ \frac{d\widetilde{q}_\varepsilon(t)}{dt}, q_\varepsilon(t), t_i, t_f \right] \right|^{p_1} \times \right.$$

$$\left. \left\{ F_{3,\varepsilon}\left[ \frac{d\widetilde{q}_\varepsilon(t)}{dt}, q_\varepsilon(t), t_i, t_f \right] \right\} dy D^+[q(t)] \right)_{\varepsilon \in (0,1]} < \infty, \qquad (2.2.2.32)$$

$$\left( \overline{\int_{\substack{q(t_f)=x \\ q(t_i)=y}}} \left| F_{2,\varepsilon}\left[ \frac{d\widetilde{q}_\varepsilon(t)}{dt}, q_\varepsilon(t), t_i, t_f \right] \right|^{p_2} \times \right.$$

$$\left. \left\{ F_{3,\varepsilon}\left[ \frac{d\widetilde{q}_\varepsilon(t)}{dt}, q_\varepsilon(t), t_i, t_f \right] \right\} dy D^+[q(t)] \right)_{\varepsilon \in (0,1]} < \infty.$$

Then Hölder's inequality for Feinman-Colombeau path integral states that:

$$\left( \overline{\int_{\substack{q(t_f)=x \\ q(t_i)=y}}} \left| \left( F_{1,\varepsilon}\left[ \frac{d\widetilde{q}_\varepsilon(t)}{dt}, q_\varepsilon(t), t_i, t_f \right] F_{2,\varepsilon}\left[ \frac{d\widetilde{q}_\varepsilon(t)}{dt}, q_\varepsilon(t), t_i, t_f \right] \right) \right| \times \right.$$

$$\left. \left\{ \left( F_{3,\varepsilon}\left[ \frac{d\widetilde{q}_\varepsilon(t)}{dt}, q_\varepsilon(t), t_i, t_f \right] \right) \right\} dy D^+[q(t)] \right)_{\varepsilon\in(0,1]} =$$

$$\int_{\substack{q(t_f)=x \\ q(t_i)=y}} \left| \left( F_{1,\varepsilon}\left[ \frac{d\widetilde{q}_\varepsilon(t)}{dt}, q_\varepsilon(t), t_i, t_f \right] F_{2,\varepsilon}\left[ \frac{d\widetilde{q}_\varepsilon(t)}{dt}, q_\varepsilon(t), t_i, t_f \right] \right)_{\varepsilon\in(0,1]} \right| \times$$

$$\left\{ \left( F_{3,\varepsilon}\left[ \frac{d\widetilde{q}_\varepsilon(t)}{dt}, q_\varepsilon(t), t_i, t_f \right] \right)_{\varepsilon\in(0,1]} \right\} dy D^+[q(t)] \le \tag{2.2.2.33}$$

$$\left( \left[ \overline{\int_{\substack{q(t_f)=x \\ q(t_i)=y}}} \left| F_{1,\varepsilon}\left[ \frac{d\widetilde{q}_\varepsilon(t)}{dt}, q_\varepsilon(t) \right] \right|^{p_1} \times \right.\right.$$

$$\left. \left[ \left\{ F_{3,\varepsilon}\left[ \frac{d\widetilde{q}_\varepsilon(t)}{dt}, q_\varepsilon(t), t_i, t_f \right] \right\} dy D^+[q(t)] \right]^{\frac{1}{p_1}} \right)_{\varepsilon\in(0,1]} \times$$

$$\left( \left[ \overline{\int_{\substack{q(t_f)=x \\ q(t_i)=y}}} \left| F_2\left[ \frac{d\widetilde{q}_\varepsilon(t)}{dt}, q_\varepsilon(t), t_i, t_f \right] \right|^{p_2} \times \right.\right.$$

$$\left. \left\{ F_3\left[ \frac{d\widetilde{q}_\varepsilon(t)}{dt}, q_\varepsilon(t), t_i, t_f \right] \right\} dy D^+[q(t)] \right]^{\frac{1}{p_2}} \right)_{\varepsilon\in(0,1]},$$

Let us avaluate path integral $\tilde{v}_{1,\varepsilon}\left(T, \frac{x}{\sqrt{\epsilon}}, x_0, \epsilon, \eta\right)$ and $\tilde{v}_{2,\varepsilon}\left(T, \frac{x}{\sqrt{\epsilon}}, x_0, \epsilon, \eta\right)$. Assume now that $p_2 >> \epsilon^{-1}$, i.e. $p_1 \approx 1$. Then from Eqs.(2.2.2.26)-(2.2.2.27) and Hölder's inequality (2.2.2.33) we obtain:

$$|\widetilde{v}_{1,\varepsilon}(T,x,x_0,\epsilon,\eta)| \leq$$

$$\overline{\int\limits_{q(T)=\frac{x}{\sqrt{\epsilon}}} \widetilde{D}_{\mathbf{m}}^{(1)}[q(t),\varepsilon,\epsilon]\Psi_0(\sqrt{\epsilon}\,q(0))\times}$$

$$\left|\cos\left(\int\limits_0^T dt\left[\frac{m}{2}\left(\frac{d\widetilde{q}_{\varepsilon\epsilon}(t)}{dt}\right)^2 - a_\varepsilon(t)q_{v(\varepsilon,\epsilon)}^2(t) + b_\varepsilon(t)\frac{q_{v(\varepsilon,\epsilon)}(t)}{\sqrt{\epsilon}}\right]\right)\right| \times$$

$$\left|\cos\left(\int\limits_0^T \widehat{V}_1(\sqrt{\epsilon}\,q_{v(\varepsilon,\epsilon)}(t))dt\right)\right| \leq$$

$$\left[\overline{\int\limits_{q(T)=\frac{x}{\sqrt{\epsilon}}} \widetilde{D}_{\mathbf{m}}^{(1)}[q(t),\varepsilon,\epsilon]\Psi_0(\sqrt{\epsilon}\,q(0))\times}\right.$$

$$\left.\left|\cos\left(\int\limits_0^T dt\left[\frac{m}{2}\left(\frac{d\widetilde{q}_{\varepsilon\epsilon}(t)}{dt}\right)^2 - a_\varepsilon(t)q_{v(\varepsilon,\epsilon)}^2(t) + b_\varepsilon(t)\frac{q_{v(\varepsilon,\epsilon)}(t)}{\sqrt{\epsilon}}\right]\right)\right|\right]^{1/p_1} \times \tag{2.2.2.34}$$

$$\left[\overline{\int\limits_{q(T)=\frac{x}{\sqrt{\epsilon}}} \widetilde{D}_{\mathbf{m}}^{(1)}[q(t),\varepsilon,\epsilon]\Psi_0(\sqrt{\epsilon}\,q(0))\times}\right.$$

$$\left|\cos\left(\int\limits_0^T dt\left[\frac{m}{2}\left(\frac{d\widetilde{q}_{\varepsilon\epsilon}(t)}{dt}\right)^2 - a_\varepsilon(t)q_{v(\varepsilon,\epsilon)}^2(t) + b_\varepsilon(t)\frac{q_{v(\varepsilon,\epsilon)}(t)}{\sqrt{\epsilon}}\right]\right)\right| \times$$

$$\left.\left|\cos\left(\int\limits_0^T \widehat{V}_1(\sqrt{\epsilon}\,q_{v(\varepsilon,\epsilon)}(t))dt\right)\right|^{p_2}\right]^{1/p_2}.$$

and

$$|\widetilde{v}_{2,\varepsilon}(T,x,x_0,\epsilon,\eta)| \leq$$

$$\overline{\int\limits_{q(T)=\frac{x}{\sqrt{\epsilon}}} \widetilde{D}_{\mathbf{m}}^{(2)}[q(t),\varepsilon,\epsilon]} \Psi_0(\sqrt{\epsilon}\,q(0)) \times$$

$$\left|\sin\left(\int_0^T dt\left[\frac{m}{2}\left(\frac{d\widetilde{q}_{\varepsilon\epsilon}(t)}{dt}\right)^2 - a_\varepsilon(t)q_{\nu(\varepsilon,\epsilon)}^2(t) + b_\varepsilon(t)\frac{q_{\nu(\varepsilon,\epsilon)}(t)}{\sqrt{\epsilon}}\right]\right)\right| \times$$

$$\left|\sin\left(-\int_0^T \widehat{V}_1(\sqrt{\epsilon}\,q_{\nu(\varepsilon,\epsilon)}(t))dt\right)\right| \leq$$

$$\left[\overline{\int\limits_{q(T)=\frac{x}{\sqrt{\epsilon}}} \widetilde{D}_{\mathbf{m}}^{(2)}[q(t),\varepsilon,\epsilon]}\,\Psi_0(\sqrt{\epsilon}\,q(0)) \times\right.$$

$$\left|\sin\left(\int_0^T dt\left[\frac{m}{2}\left(\frac{d\widetilde{q}_{\varepsilon\epsilon}(t)}{dt}\right)^2 - a_\varepsilon(t)q_{\nu(\varepsilon,\epsilon)}^2(t) + b_\varepsilon(t)\frac{q_{\nu(\varepsilon,\epsilon)}(t)}{\sqrt{\epsilon}}\right]\right)\right|\right]^{1/p_1} \times$$

$$\left[\overline{\int\limits_{q(T)=\frac{x}{\sqrt{\epsilon}}} \widetilde{D}_{\mathbf{m}}^{(2)}[q(t),\varepsilon,\epsilon]}\,\Psi_0(\sqrt{\epsilon}\,q(0)) \times\right.$$

$$\left|\sin\left(\int_0^T dt\left[\frac{m}{2}\left(\frac{d\widetilde{q}_{\varepsilon\epsilon}(t)}{dt}\right)^2 - a_\varepsilon(t)q_{\nu(\varepsilon,\epsilon)}^2(t) + b_\varepsilon(t)\frac{q_{\nu(\varepsilon,\epsilon)}(t)}{\sqrt{\epsilon}}\right]\right)\right| \times$$

$$\left.\left|\sin\left(-\int_0^T \widehat{V}_1(\sqrt{\epsilon}\,q_{\nu(\varepsilon,\epsilon)}(t))dt\right)\right|^{p_2}\right]^{1/p_2}.$$

$$(2.2.2.35)$$

Therefore we obtain

$$|\widetilde{v}_{1,\varepsilon}(T,x,x_0,\epsilon,\eta)| \leq$$

$$\left[ \overline{\int\limits_{q(T)=\frac{x}{\sqrt{\epsilon}}} \widetilde{D}_{\mathbf{m}}^{(1)}[q(t),\varepsilon,\epsilon]\Psi_0(\sqrt{\epsilon}\,q(0)) \times} \right.$$

$$\left. \left| \cos\left( \int\limits_0^T dt \left[ \frac{m}{2}\left( \frac{d\widetilde{q}_{\varepsilon\epsilon}(t)}{dt} \right)^2 - a_\varepsilon(t)q_{v(\varepsilon,\epsilon)}^2(t) + b_\varepsilon(t)\frac{q_{v(\varepsilon,\epsilon)}(t)}{\sqrt{\epsilon}} \right] \right) \right| \right] \times \qquad (2.2.2.36.a)$$

$$\left[ \overline{\int\limits_{q(T)=\frac{x}{\sqrt{\epsilon}}} \widetilde{D}_{\mathbf{m}}^{(1)}[q(t),\varepsilon,\epsilon]\Psi_0(\sqrt{\epsilon}\,q(0)) \times} \right.$$

$$\left. \left| \cos\left( \int\limits_0^T dt \left[ \frac{m}{2}\left( \frac{d\widetilde{q}_{\varepsilon\epsilon}(t)}{dt} \right)^2 - a_\varepsilon(t)q_{v(\varepsilon,\epsilon)}^2(t) + b_\varepsilon(t)\frac{q_{v(\varepsilon,\epsilon)}(t)}{\sqrt{\epsilon}} \right] \right) \right| \right]^{1/p_2}.$$

and

$$\left| \widetilde{v}_{2,\varepsilon}(T,x,x_0,\epsilon,\eta) \right| \leq$$

$$\left[ \overline{\int\limits_{q(T)=\frac{x}{\sqrt{\epsilon}}} \widetilde{D}_{\mathbf{m}}^{(2)}[q(t),\varepsilon,\epsilon]\Psi_0(\sqrt{\epsilon}\,q(0)) \times} \right.$$

$$\left| \sin\left( \int\limits_0^T dt \left[ \frac{m}{2}\left( \frac{d\widetilde{q}_{\varepsilon\epsilon}(t)}{dt} \right)^2 - a_\varepsilon(t)q_{v(\varepsilon,\epsilon)}^2(t) + b_\varepsilon(t)\frac{q_{v(\varepsilon,\epsilon)}(t)}{\sqrt{\epsilon}} \right] \right) \right| \times \qquad (2.2.2.36.b)$$

$$\left[ \overline{\int\limits_{q(T)=\frac{x}{\sqrt{\epsilon}}} \widetilde{D}_{\mathbf{m}}^{(2)}[q(t),\varepsilon,\epsilon]\Psi_0(\sqrt{\epsilon}\,q(0)) \times} \right.$$

$$\left. \left| \sin\left( \int\limits_0^T dt \left[ \frac{m}{2}\left( \frac{d\widetilde{q}_{\varepsilon\epsilon}(t)}{dt} \right)^2 - a_\varepsilon(t)q_{v(\varepsilon,\epsilon)}^2(t) + b_\varepsilon(t)\frac{q_{v(\varepsilon,\epsilon)}(t)}{\sqrt{\epsilon}} \right] \right) \right| \right]^{1/p_2}.$$

From Eqs.(2.2.2.36.a,b) we obtain

$$\left| \widetilde{v}_{1,\varepsilon}\left( T,\frac{x}{\sqrt{\epsilon}},x_0,\epsilon,\eta \right) \right| \leq \mathbf{J}_{1,\varepsilon}(T,x,x_0,\epsilon,\eta) \times \Gamma_{1,\varepsilon}(T,x,x_0,\epsilon,\eta,p_2),$$

$$\qquad (2.2.2.37)$$

$$\left| \widetilde{v}_{2,\varepsilon}\left( T,\frac{x}{\sqrt{\epsilon}},x_0,\epsilon,\eta \right) \right| \leq \mathbf{J}_{2,\varepsilon}(T,x,x_0,\epsilon,\eta) \times \Gamma_{2,\varepsilon}(T,x,x_0,\epsilon,\eta,p_2).$$

Here

$$\mathbf{J}_{1,\varepsilon}(T,x,x_0,\epsilon,\eta) =$$

$$\overline{\int\limits_{q(T)=\frac{x}{\sqrt{\epsilon}}} \widetilde{D}_{\mathbf{m}}^{(1)}[q(t),\varepsilon,\epsilon]\Psi_0(\sqrt{\epsilon}\,q(0)) \times}$$

$$\left|\cos\left(\int\limits_0^T dt\left[\frac{m}{2}\left(\frac{d\widetilde{q}_{\varepsilon\epsilon}(t)}{dt}\right)^2 - a_\varepsilon(t)q_{\nu(\varepsilon,\epsilon)}^2(t) + b_\varepsilon(t)\frac{q_{\nu(\varepsilon,\epsilon)}(t)}{\sqrt{\epsilon}}\ \right]\right)\right|,$$

$$(2.2.2.38.a)$$

$$\Gamma_{1,\varepsilon}(T,x,x_0,\epsilon,\eta,p_2) =$$

$$\left[\overline{\int\limits_{q(T)=\frac{x}{\sqrt{\epsilon}}} \widetilde{D}_{\mathbf{m}}^{(1)}[q(t),\varepsilon\epsilon]\Psi_0(\sqrt{\epsilon}\,q(0)) \times}\right.$$

$$\left.\left|\cos\left(\int\limits_0^T dt\left[\frac{m}{2}\left(\frac{d\widetilde{q}_{\varepsilon\epsilon}(t)}{dt}\right)^2 - a_\varepsilon(t)q_{\nu(\varepsilon,\epsilon)}^2(t) + b_\varepsilon(t)\frac{q_{\nu(\varepsilon,\epsilon)}(t)}{\sqrt{\epsilon}}\ \right]\right)\right|\ \right]^{1/p_2}.$$

and

$$\mathbf{J}_{2,\varepsilon}(T,x,x_0,\epsilon,\eta) =$$

$$\overline{\int\limits_{q(T)=\frac{x}{\sqrt{\epsilon}}} \widetilde{D}_{\mathbf{m}}^{(2)}[q(t),\varepsilon\epsilon]\Psi_0(\sqrt{\epsilon}\,q(0)) \times}$$

$$\left|\sin\left(\int\limits_0^T dt\left[\frac{m}{2}\left(\frac{d\widetilde{q}_{\varepsilon\epsilon}(t)}{dt}\right)^2 - a_\varepsilon(t)q_{\nu(\varepsilon,\epsilon)}^2(t) + b_\varepsilon(t)\frac{q_{\nu(\varepsilon,\epsilon)}(t)}{\sqrt{\epsilon}}\right]\right)\right|,$$

$$(2.2.2.38.b)$$

$$\Gamma_{2,\varepsilon}(T,x,x_0,\epsilon,\eta,p_2) =$$

$$\left[\overline{\int\limits_{q(T)=\frac{x}{\sqrt{\epsilon}}} \widetilde{D}_{\mathbf{m}}^{(2)}[q(t),\varepsilon,\epsilon]\Psi_0(\sqrt{\epsilon}\,q(0)) \times}\right.$$

$$\left.\left|\sin\left(\int\limits_0^T dt\left[\frac{m}{2}\left(\frac{d\widetilde{q}_{\varepsilon\epsilon}(t)}{dt}\right)^2 - a_\varepsilon(t)q_{\nu(\varepsilon,\epsilon)}^2(t) + b_\varepsilon(t)\frac{q_{\nu(\varepsilon,\epsilon)}(t)}{\sqrt{\epsilon}}\right]\right)\right|\right]^{1/p_2}.$$

From Eq.(2.2.2.38.a,b) by using definition (see Eq.(2.2.2.23)) of the submeasures $\widetilde{D}_{\mathbf{m}}^{(1)}[q(t),\varepsilon,\epsilon], \widetilde{D}_{\mathbf{m}}^{(2)}[q(t),\varepsilon,\epsilon]$ we obtain

$$\mathbf{J}_{1,\varepsilon}(T,x,x_0,\epsilon,\eta) \leq$$

$$\overline{\int\limits_{q(T)=\frac{x}{\sqrt{\epsilon}}} \widetilde{D}_{\mathbf{sm}}^{(1)}[q(t),\varepsilon,\epsilon] \Psi_0(\sqrt{\epsilon}\,q(0)) \times}$$

$$\left| \cos\left( \int\limits_0^T dt \left[ \frac{m}{2}\left( \frac{d\widetilde{q}_{\varepsilon\epsilon}(t)}{dt} \right)^2 - a_\varepsilon(t)q_{\nu(\varepsilon,\epsilon)}^2(t) + b_\varepsilon(t)\frac{q_{\nu(\varepsilon,\epsilon)}(t)}{\sqrt{\epsilon}} \right] \right) \right| +$$

$$(2.2.2.39.a)$$

$$+ \overline{\int\limits_{q(T)=\frac{x}{\sqrt{\epsilon}}} \widetilde{D}_{\mathbf{\eta}}^{+}[q(t),\varepsilon,\epsilon] \Psi_0(\sqrt{\epsilon}\,q(0)) \times}$$

$$\left| \cos\left( \int\limits_0^T dt \left[ \frac{m}{2}\left( \frac{d\widetilde{q}_{\varepsilon\epsilon}(t)}{dt} \right)^2 - a_\varepsilon(t)q_{\nu(\varepsilon,\epsilon)}^2(t) + b_\varepsilon(t)\frac{q_{\nu(\varepsilon,\epsilon)}(t)}{\sqrt{\epsilon}} \right] \right) \right| =$$

$$\Theta_{1,\varepsilon}(T,x,x_0,\epsilon,\eta) + \Lambda_{1,\varepsilon}(T,x,x_0,\epsilon,\eta)$$

and

$$\mathbf{J}_{2,\varepsilon}(T,x,x_0,\epsilon,\eta) \leq$$

$$\overline{\int\limits_{q(T)=\frac{x}{\sqrt{\epsilon}}} \widetilde{D}_{\mathbf{sm}}^{(2)}[q(t),\varepsilon,\epsilon]\Psi_0(\sqrt{\epsilon}\,q(0))\times}$$

$$\left| \sin\left( \int\limits_0^T dt \left[ \frac{m}{2}\left( \frac{d\widetilde{q}_{\varepsilon\epsilon}(t)}{dt} \right)^2 - a_\varepsilon(t)q_{\nu(\varepsilon,\epsilon)}^2(t) + b_\varepsilon(t)\frac{q_{\nu(\varepsilon,\epsilon)}(t)}{\sqrt{\epsilon}} \right] \right) \right| +$$

$$(2.2.2.39.b)$$

$$+ \overline{\int\limits_{q(T)=\frac{x}{\sqrt{\epsilon}}} \widetilde{D}_{\mathbf{\eta}}^+[q(t),\varepsilon,\epsilon]\Psi_0(\sqrt{\epsilon}\,q(0))\times}$$

$$\left| \sin\left( \int\limits_0^T dt \left[ \frac{m}{2}\left( \frac{d\widetilde{q}_{\varepsilon\epsilon}(t)}{dt} \right)^2 - a_\varepsilon(t)q_{\nu(\varepsilon,\epsilon)}^2(t) + b_\varepsilon(t)\frac{q_{\nu(\varepsilon,\epsilon)}(t)}{\sqrt{\epsilon}} \right] \right) \right| \times$$

$$\Theta_{2,\varepsilon}(T,x,x_0,\epsilon,\eta) + \Lambda_{2,\varepsilon}(T,x,x_0,\epsilon,\eta).$$

Here

$$\Theta_{1,\varepsilon}(T,x,x_0,\epsilon,\eta) = \overline{\int\limits_{q(T)=\frac{x}{\sqrt{\epsilon}}} \widetilde{D}_{\mathbf{sm}}^{(1)}[q(t),\varepsilon,\epsilon]} \Psi_0(\sqrt{\epsilon}\,q(0)) \times$$

$$\left| \cos\left( \int_0^T dt \left[ \frac{m}{2}\left(\frac{d\widetilde{q}_{\varepsilon\epsilon}(t)}{dt}\right)^2 - a_\varepsilon(t)q_{\nu(\varepsilon,\epsilon)}^2(t) + b_\varepsilon(t)\frac{q_{\nu(\varepsilon,\epsilon)}(t)}{\sqrt{\epsilon}} \right] \right) \right| =$$

$$\overline{\int\limits_{q(T)=\frac{x}{\sqrt{\epsilon}}} \widetilde{D}_{\eta}^{+}[q(t),\varepsilon,\epsilon]} \Psi_0(\sqrt{\epsilon}\,q(0)) \times$$

$$\cos\left( \int_0^T dt \left[ \frac{m}{2}\left(\frac{d\widetilde{q}_{\varepsilon\epsilon}(t)}{dt}\right)^2 - a_\varepsilon(t)q_{\nu(\varepsilon,\epsilon)}^2(t) + b_\varepsilon(t)\frac{q_{\nu(\varepsilon,\epsilon)}(t)}{\sqrt{\epsilon}} \right] \right),$$

$$(2.2.2.40.a)$$

$$\Theta_{2,\varepsilon}(T,x,x_0,\epsilon,\eta) = \overline{\int\limits_{q(T)=\frac{x}{\sqrt{\epsilon}}} \widetilde{D}_{\mathbf{sm}}^{(2)}[q(t),\varepsilon,\epsilon]} \Psi_0(\sqrt{\epsilon}\,q(0)) \times$$

$$\left| \sin\left( \int_0^T dt \left[ \frac{m}{2}\left(\frac{d\widetilde{q}_{\varepsilon\epsilon}(t)}{dt}\right)^2 - a_\varepsilon(t)q_{\nu(\varepsilon,\epsilon)}^2(t) + b_\varepsilon(t)\frac{q_{\nu(\varepsilon,\epsilon)}(t)}{\sqrt{\epsilon}} \right] \right) \right| =$$

$$\overline{\int\limits_{q(T)=\frac{x}{\sqrt{\epsilon}}} \widetilde{D}_{\eta}^{+}[q(t),\varepsilon,\epsilon]} \Psi_0(\sqrt{\epsilon}\,q(0)) \times$$

$$\sin\left( \int_0^T dt \left[ \frac{m}{2}\left(\frac{d\widetilde{q}_{\varepsilon\epsilon}(t)}{dt}\right)^2 - a_\varepsilon(t)q_{\nu(\varepsilon,\epsilon)}^2(t) + b_\varepsilon(t)\frac{q_{\nu(\varepsilon,\epsilon)}(t)}{\sqrt{\epsilon}} \right] \right),$$

$$(2.2.2.40.b)$$

and

$$\Lambda_{1,\varepsilon}(T,x,x_0,\epsilon,\eta) = \overline{\int\limits_{q(T)=\frac{x}{\sqrt{\epsilon}}} \widetilde{D}_{\eta}^{+}[q(t),\varepsilon\epsilon]\Psi_0(\sqrt{\epsilon}\,q(0)) \times}$$

$$\left| \cos\left( \int\limits_0^T dt \left[ \frac{m}{2}\left( \frac{d\widetilde{q}_{\varepsilon\epsilon}(t)}{dt} \right)^2 - a_\varepsilon(t)q_{\nu(\varepsilon,\epsilon)}^2(t) + b_\varepsilon(t)\frac{q_{\nu(\varepsilon,\epsilon)}(t)}{\sqrt{\epsilon}} \right] \right) \right|,$$

$$(2.2.2.41.a)$$

$$\Lambda_{2,\varepsilon}(T,x,x_0,\epsilon,\eta) = \overline{\int\limits_{q(T)=\frac{x}{\sqrt{\epsilon}}} \widetilde{D}_{\eta}^{+}[q(t),\varepsilon,\epsilon]\Psi_0(\sqrt{\epsilon}\,q(0)) \times}$$

$$\left| \sin\left( \int\limits_0^T dt \left[ \frac{m}{2}\left( \frac{d\widetilde{q}_{\varepsilon\epsilon}(t)}{dt} \right)^2 - a_\varepsilon(t)q_{\nu(\varepsilon,\epsilon)}^2(t) + b_\varepsilon(t)\frac{q_{\nu(\varepsilon,\epsilon)}(t)}{\sqrt{\epsilon}} \right] \right) \right|.$$

$$(2.2.2.41.b)$$

From (2.2.2.37.a,b) and (2.2.2.39.a,b) we obtain the inequalities:

$$\left| \widetilde{v}_{1,\varepsilon}\left( T,\frac{x}{\sqrt{\epsilon}},x_0,\epsilon,\eta \right) \right| \le \mathbf{J}_{1,\varepsilon}(T,x,x_0,\epsilon,\eta) \le$$

$$\Theta_{1,\varepsilon}(T,x,x_0,\epsilon,\eta) + \Lambda_{1,\varepsilon}(T,x,x_0,\epsilon,\eta),$$

$$\left| \widetilde{v}_{2,\varepsilon}\left( T,\frac{x}{\sqrt{\epsilon}},x_0,\epsilon,\eta \right) \right| \le \mathbf{J}_{2,\varepsilon}(T,x,x_0,\epsilon,\eta) \le$$

$$(2.2.2.42)$$

$$\Theta_{2,\varepsilon}(T,x,x_0,\epsilon,\eta) + \Lambda_{2,\varepsilon}(T,x,x_0,\epsilon,\eta).$$

Hence

$$v_{1,\varepsilon}\left(T,\frac{x}{\sqrt{\epsilon}},x_0,\epsilon,\eta\right) = \widetilde{v}_{1,\varepsilon}\left(T,\frac{x}{\sqrt{\epsilon}},x_0,\epsilon,\eta\right) - v_{1,\varepsilon}^*\left(T,\frac{x}{\sqrt{\epsilon}},x_0,\epsilon,\eta\right) \leq$$

$$\mathbf{J}_{1,\varepsilon}(T,x,x_0,\epsilon,\eta) - v_{1,\varepsilon}^*\left(T,\frac{x}{\sqrt{\epsilon}},x_0,\epsilon,\eta\right) \leq \qquad (2.2.2.43.a)$$

$$\Theta_{1,\varepsilon}(T,x,x_0,\epsilon,\eta) + \Lambda_{1,\varepsilon}(T,x,x_0,\epsilon,\eta) - v_{1,\varepsilon}^*\left(T,\frac{x}{\sqrt{\epsilon}},x_0,\epsilon,\eta\right)$$

and

$$v_{2,\varepsilon}\left(T,\frac{x}{\sqrt{\epsilon}},x_0,\epsilon,\eta\right) = \widetilde{v}_{2,\varepsilon}\left(T,\frac{x}{\sqrt{\epsilon}},x_0,\epsilon,\eta\right) - v_{2,\varepsilon}^*\left(T,\frac{x}{\sqrt{\epsilon}},x_0,\epsilon,\eta\right) \leq$$

$$\mathbf{J}_{2,\varepsilon}(T,x,x_0,\epsilon,\eta) - v_{2,\varepsilon}^*\left(T,\frac{x}{\sqrt{\epsilon}},x_0,\epsilon,\eta\right) \leq \qquad (2.2.2.43.b)$$

$$\Theta_{2,\varepsilon}(T,x,x_0,\epsilon,\eta) + \Lambda_{2,\varepsilon}(T,x,x_0,\epsilon,\eta) - v_{2,\varepsilon}^*\left(T,\frac{x}{\sqrt{\epsilon}},x_0,\epsilon,\eta\right).$$

Substitution Eqs.(2.2.2.30,31) and Eqs.(2.2.2.41.a,b) into Eqs.(2.2.2.43.a,b) gives

$$\nu_{1,\varepsilon}\left(T, \frac{x}{\sqrt{\epsilon}}, x_0, \epsilon, \eta\right) \leq$$

$$\Theta_{1,\varepsilon}(T, x, x_0, \epsilon, \eta) + \Lambda_{1,\varepsilon}(T, x, x_0, \epsilon, \eta) -$$

$$\left\{ \overline{\int\limits_{\substack{q(T)=\frac{x}{\sqrt{\epsilon}} \\ q(0)=\frac{y}{\sqrt{\epsilon}}}} dy D^{+}[q(t), \varepsilon, \epsilon] \Psi_0(y) \Xi_1(q(t), \varepsilon, \epsilon) \exp\left\{ -\int\limits_0^T dt \left[ \{q(t)\}_\eta \right]^2 \right\} + \right.$$

$$\sum_{n=1}^{M} \frac{(-1)^n \epsilon^n}{(2n)!} \left[ \overline{\int\limits_{\substack{q(T)=\frac{x}{\sqrt{\epsilon}} \\ q(0)=\frac{y}{\sqrt{\epsilon}}}} dy D^{+}[q(t), \varepsilon, \epsilon] \Psi_0(y) \Xi_1(q(t), \varepsilon, \epsilon) \times \right.$$

$$\left. \left. \exp\left\{ -\int\limits_0^T dt \left[ \{q(t)\}_\eta \right]^2 \right\} \left( \int\limits_0^T \widetilde{V}_1(q_{\nu(\varepsilon,\epsilon)}(t)) dt \right)^{2n} \right] \right\} = \quad (2.2.2.44.a)$$

$$\Theta_{1,\varepsilon}(T, x, x_0, \epsilon, \eta) + \sum_{n=1}^{M} \frac{(-1)^n \epsilon^n}{(2n)!} \left[ \overline{\int\limits_{\substack{q(T)=\frac{x}{\sqrt{\epsilon}} \\ q(0)=\frac{y}{\sqrt{\epsilon}}}} dy D^{+}[q(t), \varepsilon, \epsilon] \Psi_0(y) \Xi_1(q(t), \varepsilon, \epsilon) \times \right.$$

$$\left. \exp\left\{ -\int\limits_0^T dt \left[ \{q(t)\}_\eta \right]^2 \right\} \left( \int\limits_0^T \widetilde{V}_1(q_{\nu(\varepsilon,\epsilon)}(t)) dt \right)^{2n} \right]$$

and

$$v_{2,\varepsilon}\left(T, \frac{x}{\sqrt{\epsilon}}, x_0, \epsilon, \eta\right) \leq$$

$$\Theta_{2,\varepsilon}(T, x, x_0, \epsilon, \eta) + \Lambda_{2,\varepsilon}(T, x, x_0, \epsilon, \eta) -$$

$$\left\{ \overline{\int\limits_{\substack{q(T)=\frac{x}{\sqrt{\epsilon}} \\ q(0)=\frac{y}{\sqrt{\epsilon}}}} dy D^+[q(t), \varepsilon, \epsilon] \Psi_0(y) \Xi_2(q(t), \varepsilon, \epsilon) \exp\left\{-\int\limits_0^T dt \left[\{q(t)\}_\eta\right]^2\right\}} + \right.$$

$$\sum_{n=1}^M \frac{(-1)^n \epsilon^{n+1/2}}{(2n)!} \left[ \overline{\int\limits_{\substack{q(T)=\frac{x}{\sqrt{\epsilon}} \\ q(0)=\frac{y}{\sqrt{\epsilon}}}} dy D^+[q(t), \varepsilon, \epsilon] \Psi_0(y) \Xi_2(q(t), \varepsilon, \epsilon) \times \right. \tag{2.2.2.44.b}$$

$$\left.\left.\exp\left\{-\int\limits_0^T dt \left[\{q(t)\}_\eta\right]^2\right\} \left(\int\limits_0^T \widetilde{V}_1(q_{v(\varepsilon,\epsilon)}(t)) dt\right)^{2n}\right]\right\} =$$

$$\Theta_{2,\varepsilon}(T, x, x_0, \epsilon, \eta) + \sum_{n=1}^M \frac{(-1)^n \epsilon^n}{(2n+1)!} \left[ \overline{\int\limits_{\substack{q(T)=\frac{x}{\sqrt{\epsilon}} \\ q(0)=\frac{y}{\sqrt{\epsilon}}}} dy D^+[q(t), \varepsilon, \epsilon] \Psi_0(y) \Xi_2(q(t), \varepsilon, \epsilon) \times \right.$$

$$\left.\exp\left\{-\int\limits_0^T dt \left[\{q(t)\}_\eta\right]^2\right\} \left(-\int\limits_0^T \widetilde{V}_1(q_{v(\varepsilon,\epsilon)}(t)) dt\right)^{2n+1}\right].$$

Therefore

$$v_{1,\epsilon}\left(T, \frac{x}{\sqrt{\epsilon}}, x_0, \epsilon, \eta\right) \leq \Theta_{1,\epsilon}(T, x, x_0, \epsilon, \eta) +$$

$$\sum_{n=1}^{M} \frac{(-1)^n \epsilon^n}{(2n)!} \left[ \overline{\int\limits_{\substack{q(T)=\frac{x}{\sqrt{\epsilon}} \\ q(0)=\frac{y}{\sqrt{\epsilon}}}} dy D^+[q(t), \varepsilon, \epsilon] \Psi_0(y) \Xi_1(q(t), \varepsilon, \epsilon) \times \times} \right. \qquad (2.2.2.45.a)$$

$$\left. \exp\left\{ -\int\limits_0^T dt \Big[ \{q(t)\}_\eta \Big]^2 \right\} \left( \int\limits_0^T \widetilde{V}_1(q_{\nu(\varepsilon,\epsilon)}(t)) dt \right)^{2n} \right].$$

and

$$v_{2,\epsilon}\left(T, \frac{x}{\sqrt{\epsilon}}, x_0, \epsilon, \eta\right) \leq \Theta_{2,\epsilon}(T, x, x_0, \epsilon, \eta) +$$

$$\sum_{n=1}^{M} \frac{(-1)^n \epsilon^{n+1/2}}{(2n+1)!} \left[ \overline{\int\limits_{\substack{q(T)=\frac{x}{\sqrt{\epsilon}} \\ q(0)=\frac{y}{\sqrt{\epsilon}}}} dy D^+[q(t), \varepsilon, \epsilon] \Psi_0(y) \Xi_2(q(t), \varepsilon, \epsilon) \times \times} \right. \qquad (2.2.2.45.b)$$

$$\left. \exp\left\{ -\int\limits_0^T dt \Big[ \{q(t)\}_\eta \Big]^2 \right\} \left( -\int\limits_0^T \widetilde{V}_1(q_{\nu(\varepsilon,\epsilon)}(t)) dt \right)^{2n+1} \right].$$

Let us calculate now oscillatory path integral $\Theta_{1,\epsilon}(T, x, x_0, \epsilon, \eta)$

$$\Theta_{1,\varepsilon}(T,x,x_0,\epsilon,\eta) =$$

$$\overline{\int\limits_{q(T)=\frac{x}{\sqrt{\epsilon}}} \widetilde{D}_{\eta}^{+}[q(t),\varepsilon\epsilon]\Psi_0(\sqrt{\epsilon}\,q(0))} \times \qquad (2.2.2.46)$$

$$\cos\left(\int\limits_{0}^{T} dt\left[\frac{m}{2}\left(\frac{d\widetilde{q}_{\varepsilon\epsilon}(t)}{dt}\right)^2 - a_{\varepsilon}(t)q_{\nu(\varepsilon,\epsilon)}^2(t) + b_{\varepsilon}(t)\frac{q_{\nu(\varepsilon,\epsilon)}(t)}{\sqrt{\epsilon}}\right]\right).$$

and oscillatory path integral $\Theta_{2,\varepsilon}(T,x,x_0,\epsilon,\eta)$

$$\Theta_{2,\varepsilon}(T,x,x_0,\epsilon,\eta) =$$

$$\overline{\int\limits_{q(T)=\frac{x}{\sqrt{\epsilon}}} \widetilde{D}_{\eta}^{+}[q(t),\varepsilon\epsilon]\Psi_0(\sqrt{\epsilon}\,q(0))} \times \qquad (2.2.2.47)$$

$$\sin\left[\left(\int\limits_{0}^{T} dt\left[\frac{m}{2}\left(\frac{d\widetilde{q}_{\varepsilon\epsilon}(t)}{dt}\right)^2 - a_{\varepsilon}(t)q_{\nu(\varepsilon,\epsilon)}^2(t) + b_{\varepsilon}(t)\frac{q_{\nu(\varepsilon,\epsilon)}(t)}{\sqrt{\epsilon}}\right]\right)\right].$$

From Eqs.(2.2.2.46,47) we obtain

$$\Theta_{1,\varepsilon}(T,x,x_0,\epsilon,\eta) =$$

$$\frac{1}{2}\overline{\int\limits_{q(T)=\frac{x}{\sqrt{\epsilon}}}\widetilde{D}_\eta^+[q(t),\varepsilon,\epsilon]\Psi_0(\sqrt{\epsilon}\,q(0))} \times$$

$$\left\{\left[\exp\left[i\widetilde{S}_0(T,x,x_0,q(t),\varepsilon,\epsilon)\right] + \exp\left[-i\widetilde{S}_0(T,x,x_0,q(t),\varepsilon,\epsilon)\right]\right]\right\} =$$

$$\frac{1}{2}\widetilde{\Theta}_{1,\varepsilon}(T,x,x_0,\epsilon,\eta) + \frac{1}{2}\widetilde{\Theta}_{2,\varepsilon}(T,x,x_0,\epsilon,\eta)$$

$(2.2.2.48.a)$

and

$$\Theta_{2,\varepsilon}(T,x,x_0,\epsilon,\eta) =$$

$$\frac{1}{2i}\overline{\int\limits_{q(T)=\frac{x}{\sqrt{\epsilon}}}\widetilde{D}_\eta^+[q(t),\varepsilon,\epsilon]\Psi_0(\sqrt{\epsilon}\,q(0))} \times$$

$$\left\{\left[\exp\left[i\widetilde{S}_0(T,x,x_0,q(t),\varepsilon,\epsilon)\right] - \exp\left[-i\widetilde{S}_0(T,x,x_0,q(t),\varepsilon,\epsilon)\right]\right]\right\} =$$

$$\frac{1}{2i}\Theta_{1,\varepsilon}(T,x,x_0,\epsilon,\eta) - \frac{1}{2i}\widetilde{\Theta}_{2,\varepsilon}(T,x,x_0,\epsilon,\eta).$$

$(2.2.2.48.b)$

Here

$$\widetilde{S}_0(T,x,x_0,q_\varepsilon(t),\varepsilon,\epsilon) = \int\limits_0^T dt\left[\frac{m}{2}\left(\frac{d\widetilde{q}_{\varepsilon\epsilon}(t)}{d\tau}\right)^2 - a_\varepsilon(t)q_\varepsilon^2(t) + b_\varepsilon(t)\frac{q_\varepsilon(t)}{\sqrt{\epsilon}}\right] \qquad (2.2.2.49)$$

and

$$\widetilde{\Theta}_{1,\varepsilon}(T,x,x_0,\epsilon,\eta) =$$

$$\overline{\int_{q(T)=\frac{x}{\sqrt{\epsilon}}} \widetilde{D}_\eta^+[q(t),\varepsilon\epsilon]\Psi_0(\sqrt{\epsilon}\,q(0))\exp\left[i\widetilde{S}_0(T,x,x_0,q(t),\varepsilon,\epsilon)\right],}$$

(2.2.2.50)

$$\widetilde{\Theta}_{2,\varepsilon}(T,x,x_0,\epsilon,\eta) =$$

$$\overline{\int_{q(T)=\frac{x}{\sqrt{\epsilon}}} \widetilde{D}_\eta^+[q(t),\varepsilon\epsilon]\Psi_0(\sqrt{\epsilon}\,q(0))\exp\left[-i\widetilde{S}_0(T,x,x_0,q(t),\varepsilon,\epsilon)\right].}$$

Let's calculate now path integral $\widetilde{\Theta}_{1,\varepsilon}(T,x,x_0,\epsilon)$ by using stationary phase method.

By using replacement $q(t) := \dfrac{q_\varepsilon(t)}{\sqrt{\epsilon}}$ we rewrite (2.2.2.50) in the next equivalent form

$$\widetilde{\Theta}_{1,\varepsilon}(T,x,x_0,\epsilon,\eta) = \overline{\int_{q(T)=\frac{x}{\sqrt{\epsilon}}} \widetilde{D}_\eta^+[q(t),\varepsilon]\Psi_0(\sqrt{\epsilon}\,q(0))} \times$$

(2.2.2.51)

$$\frac{1}{\epsilon}\int_0^T dt\left[\frac{m}{2}\left(\frac{d\widetilde{q}_\varepsilon(t)}{d\tau}\right)^2 - a_\varepsilon(t)q_\varepsilon^2(t) + b_\varepsilon(t)q_\varepsilon(t)\right].$$

We assume now that

$$\Psi(x,0,x_0) = \frac{1}{\sqrt[4]{2\pi\epsilon/\delta}}\exp\left(-\frac{x^2}{\epsilon/\delta}\right),$$

(2.2.2.52)

$$0 < \epsilon \ll \delta \ll 1.$$

Substitution Eq.(2.2.2.52) into Eq.(2.2.2.51) gives

$$\widetilde{\Theta}_{1,\varepsilon}(T,x,x_0,\epsilon,\eta) =$$

$$\frac{1}{\sqrt[4]{2\pi\epsilon/\delta}} \int\limits_{\substack{q(T)=x \\ q(0)=y}} dy \widetilde{D}_\eta^+[q(t),\varepsilon] \exp\left(-\frac{y^2}{\epsilon/\delta}\right) \times$$

$$\exp\left\{\frac{i}{\epsilon}\int\limits_0^T dt\left[\frac{m}{2}\left(\frac{d\widetilde{q}_\varepsilon(t)}{dt}\right)^2 - a_\varepsilon(t)q_\varepsilon^2(t) + b_\varepsilon(t)q_\varepsilon(t)\right]\right\} =$$

$$\frac{1}{\sqrt[4]{2\pi\epsilon/\delta}} \int\limits_{\substack{q(T)=x \\ q(0)=y}} dy \widetilde{D}_\eta^+[q(t),\varepsilon] \exp\left(-\frac{y^2}{\epsilon/\delta}\right) \times$$

$$\exp\left\{\frac{i}{\epsilon}\int\limits_0^T \mathcal{L}(\dot{q}(t),q_\varepsilon(t))dt\right\},$$

$$\mathcal{L}\left(\frac{d\widetilde{q}_\varepsilon(t)}{dt},q_\varepsilon(t)\right) = \frac{m}{2}\left(\frac{d\widetilde{q}_\varepsilon(t)}{dt}\right)^2 - a_\varepsilon(t)q_\varepsilon^2(t) + b_\varepsilon(t)q_\varepsilon(t).$$

We assume that $a_\varepsilon(t) > 0$ and rewrite $\mathcal{L}(\dot{q}(t),q_\varepsilon(t))$ in the next form

$$\mathcal{L}\left(\frac{d\widetilde{q}_\varepsilon(t)}{dt},q_\varepsilon(t)\right) = \frac{m}{2}\left(\frac{d\widetilde{q}_\varepsilon(t)}{dt}\right)^2 - \frac{m\varpi_\varepsilon^2}{2}q_\varepsilon^2(t) + b_\varepsilon(t)q_\varepsilon(t),$$

$$\frac{m\varpi_\varepsilon^2}{2} = a_\varepsilon(t).$$

From Eq.(2.2.2.1) and Eq.(2.2.2.54) for $\varepsilon \approx 0$ we obtain

$$\mathcal{L}\left(\frac{d\tilde{q}_\varepsilon(t)}{dt}, q_\varepsilon(t)\right) = \frac{m}{2}\left(\frac{dq(t)}{dt}\right)^2 - \frac{m\varpi_0^2}{2}q^2(t) + b_0(t)q(t) + o(\varepsilon) = \tag{2.2.2.55}$$

$$\mathcal{L}(\dot{q}, q) + o(\varepsilon).$$

We assume now that $\varpi_0(t) = \varpi = const.$ The Euler-Lagrange equation of $\mathcal{L}(\dot{q}, q)$ with corresponding boundary conditions is

$$\frac{d}{dt}\left(\frac{\partial \mathcal{L}}{\partial \dot{q}}\right) - \frac{\partial \mathcal{L}}{\partial q} = 0,$$

$$q(0) = y, q(T) = x. \tag{2.2.2.56}$$

Corresponding classic action $S_{\mathrm{cl}}$ is [13]:

$$S_{\mathrm{cl}}(y, x, T) = \frac{m\varpi}{2\sin \varpi T}\Bigg[(\cos \varpi T)(y^2 + x^2) - 2xy + \frac{2x}{m\varpi}\int_0^T b(t)\sin(\varpi t)dt +$$

$$+ \frac{2y}{m\varpi}\int_0^T b(t)\sin \varpi(T - t)dt - \tag{2.2.2.57}$$

$$- \frac{2}{m^2\varpi^2}\int_0^T \int_0^t b(t)b(s)\sin \varpi(T - t)\sin(\varpi s)dsdt\Bigg].$$

Substitution Eq.(2.2.2.57) into Eq.(2.2.2.53) gives

$$\widetilde{\Theta}_{1,\varepsilon}(T, x, 0, \epsilon, \eta) =$$

$$\frac{1}{\sqrt[4]{2\pi\epsilon/\delta}} \int\limits_{\substack{q(T)=x \\ q(0)=y}} dy \, \widetilde{D}_\eta[q(t), \varepsilon] \exp\left(-\frac{y^2}{\epsilon/\delta}\right) \times$$

$$\exp\left\{ \frac{i}{\epsilon} \int\limits_0^T \mathcal{L}(\dot{q}(t), q_\varepsilon(t)) dt \right\} =$$

$$= \widehat{\Theta}_{1,\varepsilon}(T, x, x_0, \epsilon, \eta).$$

(2.2.2.58)

Here

$$\widehat{\Theta}_{1,\varepsilon}(T, x, 0, \epsilon) =$$

$$\frac{1}{\sqrt[4]{2\pi\epsilon/\delta}} \sqrt{\frac{m\varpi}{2\pi i(\epsilon/\delta)\sin\varpi T}} \int dy \exp\left(-\frac{y^2}{\epsilon/\delta}\right) \times$$

(2.2.2.59)

$$\exp\left[ \frac{i}{\epsilon} S_{\mathbf{cl}}(y, x, T) + \frac{i\pi\gamma}{2} \right],$$

$$\sin\varpi T \neq 0.$$

Let's calculate now oscillatory integral $\widehat{\Theta}_{1,\varepsilon}(T, x, 0, \epsilon)$ by using stationary phase method.

$$\frac{\partial S_{\mathbf{cl}}}{\partial y} = 2y\cos\varpi T - 2x + \frac{2}{m\varpi} \int_0^T b(t)\sin\varpi(T - t)dt = 0.$$

(2.2.2.60)

Therefore critical point $y_{\mathbf{cr}}(x, T)$ is

$$y_{\mathrm{cr}}(x,T) = \frac{x}{\cos \varpi T} - \frac{1}{m \varpi \cos \varpi T} \int_0^T b(t) \sin \varpi (T-t) dt =$$

$$\frac{1}{\cos \varpi T} \left( x - \frac{1}{m \varpi} \int_0^T b(t) \sin \varpi (T-t) dt \right), \qquad (2.2.2.61)$$

$$\cos \varpi T \neq 0.$$

From Eq.(2.2.2.59) and Eq.(2.2.2.61) we obtain

$$\widehat{\Theta}_{1,\varepsilon}(T,x,0,\epsilon) = \frac{1}{\sqrt[4]{2\pi \epsilon / \delta} \, |\cos \varpi T|^{1/2}} \exp\left[ -\frac{y_{\mathrm{cr}}^2(x,T)}{\epsilon / \delta} \right] \times$$

$$\exp\left[ \frac{i}{\epsilon} \widetilde{S}_{\mathrm{cl}}(x,y_{\mathrm{cr}},T) + \frac{i\pi\gamma}{2} + \frac{i\pi\gamma'}{4} \right] + O(\epsilon/m) =$$

$$\frac{1}{\sqrt[4]{2\pi \epsilon / \delta} \, |\cos \varpi T|^{1/2}} \times \qquad (2.2.2.62)$$

$$\times \exp\left[ -\frac{1}{\epsilon / \delta} \left( \frac{x}{\cos \varpi T} - \frac{1}{m \varpi \cos \varpi T} \int_0^T b(t) \sin \varpi (T-t) dt \right)^2 \right] \times$$

$$\exp\left[ \frac{i}{\epsilon} \widetilde{S}_{\mathrm{cl}}(x,y_{\mathrm{cr}},\lambda,T) + \frac{i\pi\gamma}{2} + \frac{i\pi\gamma'}{4} \right] + O(\epsilon).$$

Let's calculate now Laplase integral

$$U(T) = \int dx |x| \left| \widehat{\Theta}_{1,\varepsilon}(T, x, 0, \epsilon) \right|^2 =$$

$$\int dx \left\{ \frac{|x|}{\sqrt{2\pi\epsilon/\delta} \, |\cos \varpi T|} \times \right. $$

(2.2.2.63)

$$\left. \exp\left[ -\frac{2}{(\epsilon/\delta)\cos \varpi T} \left( x - \frac{1}{m\varpi} \int_0^T b(t) \sin \varpi(T - t) dt \right)^2 \right] \right\}$$

using Laplase method. The maximum point is

$$x_{\max} = \frac{1}{m\varpi} \int_0^T b(t) \sin \varpi(T - t) dt.$$

(2.2.2.64)

Finally we obtain

$$U(T) \simeq \left| \frac{1}{m\varpi} \int_0^T b(t) \sin \varpi(T - t) dt \right|.$$

(2.2.2.65)

From Eq.(2.2.2.45.a,b) we obtain

$$v_{1,\varepsilon}\left( T, \frac{x}{\sqrt{\epsilon}}, x_0, \epsilon, \eta \right) \leq \Theta_{1,\varepsilon}(T, x, x_0, \epsilon, \eta) + \Re_{1,\varepsilon}(T, x, x_0, \epsilon, \eta),$$

(2.2.2.66)

$$v_{2,\varepsilon}\left( T, \frac{x}{\sqrt{\epsilon}}, x_0, \epsilon, \eta \right) \leq \Theta_{2,\varepsilon}(T, x, x_0, \epsilon, \eta) + \Re_{2,\varepsilon}(T, x, x_0, \epsilon, \eta),$$

Here

$$\Re_{1,\varepsilon}(T,x,x_0,\epsilon,\eta) =$$

$$\sum_{n=1}^{M} \frac{(-1)^n \epsilon^n}{(2n)!} \left[ \overline{\int\limits_{\substack{q(T)=\frac{x}{\sqrt{\epsilon}} \\ q(0)=\frac{y}{\sqrt{\epsilon}}}} dy D^+[q(t),\varepsilon,\epsilon] \Psi_0(y) \Xi_1(q(t),\varepsilon,\epsilon) \times \times \right. \qquad (2.2.2.67.a)$$

$$\left. \exp\left\{ -\int_0^T dt \left[ \{q(t)\}_\eta \right]^2 \right\} \left( \int_0^T \widetilde{V}_1(q_{v(\varepsilon,\epsilon)}(t)) dt \right)^{2n} \right].$$

and

$$\Re_{2,\varepsilon}(T,x,x_0,\epsilon,\eta) =$$

$$\sum_{n=1}^{M} \frac{(-1)^n \epsilon^{n+1/2}}{(2n+1)!} \left[ \overline{\int\limits_{\substack{q(T)=\frac{x}{\sqrt{\epsilon}} \\ q(0)=\frac{y}{\sqrt{\epsilon}}}} dy D^+[q(t),\varepsilon,\epsilon] \Psi_0(y) \Xi_2(q(t),\varepsilon,\epsilon) \times \times \right. \qquad (2.2.2.67.b)$$

$$\left. \exp\left\{ -\int_0^T dt \left[ \{q(t)\}_\eta \right]^2 \right\} \left( -\int_0^T \widetilde{V}_1(q_{v(\varepsilon,\epsilon)}(t)) dt \right)^{2n+1} \right].$$

From (2.2.2.66) we obtain

$$v_{1,\varepsilon}^2\left(T, \frac{x}{\sqrt{\epsilon}}, x_0, \epsilon, \eta\right) \leq$$

$$[\Theta_{1,\varepsilon}(T, x, x_0, \epsilon, \eta) + \Re_{1,\varepsilon}(T, x, x_0, \epsilon, \eta)]^2 \leq \qquad (2.2.2.68.a)$$

$$2[\Theta_{1,\varepsilon}^2(T, x, x_0, \epsilon, \eta) + \Re_{1,\varepsilon}^2(T, x, x_0, \epsilon, \eta)]$$

and

$$v_{2,\varepsilon}^2\left(T, \frac{x}{\sqrt{\epsilon}}, x_0, \epsilon, \eta\right) \leq$$

$$[\Theta_{2,\varepsilon}(T, x, x_0, \epsilon, \eta) + \Re_{2,\varepsilon}(T, x, x_0, \epsilon, \eta)]^2 \leq \qquad (2.2.2.68.b)$$

$$2[\Theta_{2,\varepsilon}^2(T, x, x_0, \epsilon, \eta) + \Re_{2,\varepsilon}^2(T, x, x_0, \epsilon, \eta)].$$

From Eq.(2.2.2.12) and (2.2.2.68.a,b) we obtain

$$\chi_1(T,\epsilon,\varepsilon,\eta) = \int dx |x| \left[ v_{1,\varepsilon}\left(T,\frac{x}{\sqrt{\epsilon}},x_0,\epsilon,\eta\right) + v_{2,\varepsilon}\left(T,\frac{x}{\sqrt{\epsilon}},x_0,\epsilon,\eta\right) \right]^2 \leq$$

$$2\int dx |x| \left[ v_{1,\varepsilon}^2\left(T,\frac{x}{\sqrt{\epsilon}},x_0,\epsilon,\eta\right) + v_{2,\varepsilon}^2\left(T,\frac{x}{\sqrt{\epsilon}},x_0,\epsilon,\eta\right) \right] =$$

$$2\left[ \int dx |x| v_{1,\varepsilon}^2\left(T,\frac{x}{\sqrt{\epsilon}},x_0,\epsilon,\eta\right) + \int dx |x| v_{2,\varepsilon}^2\left(T,\frac{x}{\sqrt{\epsilon}},x_0,\epsilon,\eta\right) \right] \leq$$

$$4\left[ \int dx |x| \Theta_{1,\varepsilon}^2(T,x,x_0,\epsilon,\eta) + \int dx |x| \Theta_{2,\varepsilon}^2(T,x,x_0,\epsilon,\eta) + \right. \tag{2.2.2.69}$$

$$\left. \int dx |x| \Re_{1,\varepsilon}^2(T,x,x_0,\epsilon,\eta) + \int dx |x| \Re_{2,\varepsilon}^2(T,x,x_0,\epsilon,\eta) \right] =$$

$$4U(T) + 4\mathscr{F}_{1,\varepsilon}(T,x_0,\epsilon,\eta) + 4\mathscr{F}_{2,\varepsilon}(T,x_0,\epsilon,\eta).$$

$$\mathscr{F}_{1,\varepsilon}(T,x_0,\epsilon,\eta) = \int dx |x| \Re_{1,\varepsilon}^2(T,x,x_0,\epsilon,\eta),$$

$$\mathscr{F}_{2,\varepsilon}(T,x_0,\epsilon,\eta) = \int dx |x| \Re_{2,\varepsilon}^2(T,x,x_0,\epsilon,\eta).$$

Let us evaluate now path integral $\mathscr{F}_{1,\varepsilon}(T,x_0,\epsilon,\eta)$ and $\mathscr{F}_{2,\varepsilon}(T,x_0,\epsilon,\eta)$. From Eq.(2.2.2.67.a,b) we obtain

$$\mathscr{F}_{1,\varepsilon}(T, x_0, \epsilon, \eta) = \int dx |x| \mathfrak{R}_{1,\varepsilon}^2(T, x, x_0, \epsilon, \eta) \leq$$

$$2 \sum_{n=1}^{M} \frac{\epsilon^{2n}}{[(2n)!]^2} \times$$

$$\int dx |x| \left[ \overline{\int\limits_{\substack{q(T)=\frac{x}{\sqrt{\epsilon}} \\ q(0)=\frac{y}{\sqrt{\epsilon}}}} dy D^+[q(t), \varepsilon, \epsilon] \Psi_0(y) \exp\left\{ -\int\limits_0^T dt \left[ \{q(t)\}_\eta \right]^2 \right\} } \right.$$

$$\left. \left( \int\limits_0^T \widetilde{V}_1(q_{v(\varepsilon,\epsilon)}(t)) dt \right)^n \right]^2 =$$

$$(2.2.2.70.a)$$

$$2 \sum_{n=1}^{M} \frac{\epsilon^{2n}}{[(2n)!]^2} \int dx |x| \mathfrak{I}_{1,n}^2(T, x, x_0, \epsilon, \eta) =$$

$$2 \sum_{n=1}^{M} \frac{\epsilon^{2n}}{[(2n)!]^2} \mathscr{L}_{1,n}(T, x_0, \epsilon, \eta).$$

$$\mathfrak{I}_{1,n}(T, x, x_0, \epsilon, \eta) =$$

$$\overline{\int\limits_{\substack{q(T)=\frac{x}{\sqrt{\epsilon}} \\ q(0)=\frac{y}{\sqrt{\epsilon}}}} dy D^+[q(t), \varepsilon, \epsilon] \Psi_0(y) \exp\left\{ -\int\limits_0^T dt \left[ \{q(t)\}_\eta \right]^2 \right\} \left( \int\limits_0^T \widetilde{V}_1(q_{v(\varepsilon,\epsilon)}(t)) dt \right)^n},$$

$$\mathscr{L}_{1,n}(T, x_0, \epsilon, \eta) = \int dx |x| \mathfrak{I}_{1,n}^2(T, x, x_0, \epsilon, \eta).$$

and

$$\mathscr{F}_{2,\varepsilon}(T,x_0,\epsilon,\eta) = \int dx |x| \mathfrak{R}_{2,\varepsilon}^2(T,x,x_0,\epsilon,\eta) \leq$$

$$2\sum_{n=1}^{M} \frac{\epsilon^{2n+1}}{\left[(2n+1)!\right]^2} \times$$

$$\int dx |x| \left[ \overline{\int\limits_{\substack{q(T)=\frac{x}{\sqrt{\epsilon}} \\ q(0)=\frac{y}{\sqrt{\epsilon}}}} dy D^+[q(t),\varepsilon,\epsilon]\Psi_0(y) \exp\left\{-\int_0^T dt\left[\left\{q(t)\right\}_\eta\right]^2\right\}} \right.$$

$$\left. \left(-\int_0^T \widetilde{V}_1(q_{\nu(\varepsilon,\epsilon)}(t))dt\right)^{2n+1} \right]^2 =$$

$$2\sum_{n=1}^{M} \frac{\epsilon^{2n+1}}{\left[(2n+1)!\right]^2} \int dx |x| \mathfrak{I}_{2,n}^2(T,x,x_0,\epsilon,\eta) = \qquad (2.2.2.70.b)$$

$$2\sum_{n=1}^{M} \frac{\epsilon^{2n+1}}{\left[(2n+1)!\right]^2} \mathscr{L}_{2,n}(T,x_0,\epsilon,\eta).$$

$$\mathfrak{I}_{2,n}(T,x,x_0,\epsilon,\eta) =$$

$$\overline{\int\limits_{\substack{q(T)=\frac{x}{\sqrt{\epsilon}} \\ q(0)=\frac{y}{\sqrt{\epsilon}}}} dy D^+[q(t),\varepsilon,\epsilon]\Psi_0(y)\Xi_2(q(t),\varepsilon,\epsilon) \exp\left\{-\int_0^T dt\left[\left\{q(t)\right\}_\eta\right]^2\right\}} \times$$

$$\left(-\int_0^T \widetilde{V}_1(q_{\nu(\varepsilon,\epsilon)}(t))dt\right)^{2n+1}.$$

$$\mathscr{L}_{2,n}(T,x_0,\epsilon,\eta) = \int dx |x| \mathfrak{I}_{2,n}^2(T,x,x_0,\epsilon,\eta)$$

Let us evaluate now path integral $\mathscr{L}_{1,n}(T, x_0, \epsilon, \eta)$ and $\mathscr{L}_{2,n}(T, x_0, \epsilon, \eta)$. From Eqs.(2.2.2.70.a,b) we obtain

$$\mathscr{L}_{1,n}(T, x_0, \epsilon, \eta) =$$

$$\int dx |x| \left[ \overline{\int\limits_{\substack{q(T)=\frac{x}{\sqrt{\epsilon}} \\ q(0)=\frac{y}{\sqrt{\epsilon}}}} dy D^+[q(t), \varepsilon, \epsilon] \Psi_0(y) \exp \times \right. \tag{2.2.2.71.a}$$

$$\left. \left\{ -\int\limits_0^T dt \left[ \{q(t)\}_\eta \right]^2 \right\} \left( \int\limits_0^T \widetilde{V}_1(q_{v(\varepsilon, \epsilon)}(t)) dt \right)^n \right]^2$$

and

$$\mathscr{L}_{2,n}(T, x_0, \epsilon, \eta) =$$

$$\int dx |x| \left[ \overline{\int\limits_{\substack{q(T)=\frac{x}{\sqrt{\epsilon}} \\ q(0)=\frac{y}{\sqrt{\epsilon}}}} dy D^+[q(t), \varepsilon, \epsilon] \Psi_0(y) \exp \times \right. \tag{2.2.2.71.b}$$

$$\left. \left\{ -\int\limits_0^T dt \left[ \{q(t)\}_\eta \right]^2 \right\} \left( \int\limits_0^T \widetilde{V}_1(q_{v(\varepsilon, \epsilon)}(t)) dt \right)^{n+1} \right]^2.$$

Substitution into Eq.(2.2.2.71.a,b) $q(t) := \delta(t) + z(t),$ here

$$z(t) = \frac{x-y}{T\sqrt{\epsilon}}t + \frac{y}{\sqrt{\epsilon}},$$

$$\delta(t) = \sum_{n=1}^{\infty} a_n \sin\left(\frac{n\pi t}{T}\right)$$

(2.2.2.72)

gives

$$\mathcal{L}_{1,n}(T, x_0, \epsilon, \eta) =$$

$$\frac{1}{\sqrt{2\pi\epsilon/\delta}} \int dx |x| \left[ \overline{\int\limits_{\substack{q(T)=\frac{x}{\sqrt{\epsilon}} \\ q(0)=\frac{y}{\sqrt{\epsilon}}}} dy D^+[q(t), \varepsilon, \epsilon]} \exp\left(-\frac{y^2}{\epsilon/\delta}\right) \exp \times \right.$$

$$\left. \left\{ -\int\limits_0^T dt \left[ \{q(t)\}_\eta \right]^2 \right\} \left( \int\limits_0^T \widetilde{V}_1(q_{\nu(\varepsilon,\epsilon)}(t)) dt \right)^n \right]^2 =$$

(2.2.2.73.a)

$$\frac{1}{\sqrt{2\pi\epsilon/\delta}} \int dx |x| \left[ \overline{\int\limits_{\substack{q(T)=0 \\ q(0)=0}} dy D^+[\delta(t), \varepsilon, \epsilon]} \exp\left(-\frac{y^2}{\epsilon/\delta}\right) \exp \times \right.$$

$$\left\{ -\int\limits_0^T dt \left[ \{\delta(t)\}_\eta + \eta\left(\frac{x-y}{T\sqrt{\epsilon}}t + \frac{y}{\sqrt{\epsilon}}\right) \right]^2 \right\} \times$$

$$\left. \left( \int\limits_0^T \widetilde{V}_1\left(\delta(t) + \left(\frac{x-y}{T\sqrt{\epsilon}}t + \frac{y}{\sqrt{\epsilon}}\right)\right) dt + o(\varepsilon\epsilon) \right)^n \right]^2$$

and

$$\mathcal{L}_{2,n}(T, x_0, \epsilon, \eta) =$$

$$\frac{1}{\sqrt{2\pi\epsilon/\delta}} \int dx |x| \left[ \overline{\int_{\substack{q(T)=\frac{x}{\sqrt{\epsilon}} \\ q(0)=\frac{y}{\sqrt{\epsilon}}}} dy D^+[q(t), \varepsilon, \epsilon] \exp\left(-\frac{y^2}{\epsilon/\delta}\right)} \exp \times \right.$$

$$\left. \left\{ -\int_0^T dt \left[ \{q(t)\}_\eta \right]^2 \right\} \left( \int_0^T \widetilde{V}_1(q_{v(\varepsilon,\epsilon)}(t)) dt \right)^{n+1} \right]^2 =$$

$$(2.2.2.73.b)$$

$$\frac{1}{\sqrt{2\pi\epsilon/\delta}} \int dx |x| \left[ \overline{\int_{\substack{q(T)=0 \\ q(0)=0}} dy D^+[\delta(t), \varepsilon, \epsilon] \exp\left(-\frac{y^2}{\epsilon/\delta}\right)} \exp \times \right.$$

$$\left\{ -\int_0^T dt \left[ \{\delta(t)\}_\eta + \eta\left(\frac{x-y}{T\sqrt{\epsilon}}t + \frac{y}{\sqrt{\epsilon}}\right) \right]^2 \right\} \times$$

$$\left. \left( \int_0^T \widetilde{V}_1\left(\delta(t) + \left(\frac{x-y}{T\sqrt{\epsilon}}t + \frac{y}{\sqrt{\epsilon}}\right)\right) dt + o(\varepsilon\epsilon) \right)^{n+1} \right]^2.$$

Substitution $x := x\sqrt{\epsilon}, y := y\sqrt{\epsilon}$, into Eq.(2.2.2.73.a,b) gives

$$\mathscr{L}_{1,n}(T, x_0, \epsilon, \eta) =$$

$$\frac{\epsilon(\sqrt{\epsilon})}{\sqrt{2\pi\epsilon/\delta}} \int dx |x| \left[ \overline{\int\limits_{\substack{q(T)=0 \\ q(0)=0}} dy \, D^+[\delta(t), \varepsilon, \epsilon]} \exp(-\delta y^2) \exp \times \right.$$

$$\left\{ -\int_0^T dt \left[ \{\delta(t)\}_\eta + \eta\left(\frac{x-y}{T}t + y\right) \right]^2 \right\} \times \qquad (2.2.2.74.a)$$

$$\left. \left( \int_0^T \widetilde{V}_1\left(\delta(t) + \left(\frac{x-y}{T}t + y\right)\right) dt + o(\varepsilon\epsilon) \right)^n \right]^2 =$$

$$\left( \sqrt{\delta} \, \epsilon \right) \eta^{-N(\varepsilon)} O(\gamma_{1,n}(T, \delta)).$$

and

$$\mathcal{L}_{2,n}(T, x_0, \epsilon, \eta) =$$

$$\frac{\epsilon(\sqrt{\epsilon})}{\sqrt{2\pi\epsilon/\delta}} \int dx |x| \left[ \overline{\int\limits_{\substack{q(T)=0 \\ q(0)=0}} dy D^+[\delta(t), \epsilon, \epsilon]} \exp(-\delta y^2) \exp \times \right.$$

$$\left\{ -\int_0^T dt \left[ \left\{ \delta(t) \right\}_\eta + \eta\left( \frac{x-y}{T} t + y \right) \right]^2 \right\} \times \qquad (2.2.2.74.b)$$

$$\left. \left( \int_0^T \widetilde{V}_1\left( \delta(t) + \left( \frac{x-y}{T} t + y \right) \right) dt + o(\varepsilon\epsilon) \right)^{n+1} \right]^2 =$$

$$\left( \sqrt{\delta}\,\epsilon \right) \eta^{-N(\varepsilon)} O(\gamma_{2,n}(T, \delta)).$$

We note that $\gamma_{1,n}(T, \delta) \simeq \sqrt{\delta}$, $\gamma_{2,n}(T, \delta) \simeq \sqrt{\delta}$ if $\delta \to 0$.

Substitution Eq.(2.2.2.74.a,b) into Eq.(2.2.2.70.a,b) gives

$$\mathcal{F}_{1,\varepsilon}(T, x_0, \epsilon, \eta) = \int dx |x| \Re_{1,\varepsilon}^2(T, x, x_0, \epsilon, \eta) \leq$$

$$2 \sum_{n=1}^M \frac{\epsilon^{2n}}{\left[ (2n)! \right]^2} \mathcal{L}_{1,n}(T, x_0, \epsilon, \eta) =$$

$$2 \sum_{n=1}^M \frac{\epsilon^{2n}}{\left[ (2n)! \right]^2} \epsilon \eta^{-N(\varepsilon)} O(\gamma_{1,n}(T, \delta)) = \qquad (2.2.2.75.a)$$

$$2\eta^{-N(\varepsilon)} \sum_{n=1}^M \frac{\epsilon^{2n+1}}{\left[ (2n)! \right]^2} O\left( \widetilde{\gamma}_{1,n}(T) \right)$$

and

$$\mathcal{F}_{2,\varepsilon}(T,x_0,\epsilon,\eta) = \int dx |x| \mathfrak{R}_{2,\varepsilon}^2(T,x,x_0,\epsilon,\eta) \leq$$

$$2\sum_{n=1}^{M} \frac{\epsilon^{2n+1}}{\left[(2n+1)!\right]^2} \mathcal{L}_{2,n}(T,x_0,\epsilon,\eta) =$$

$$2\sum_{n=1}^{M} \frac{\epsilon^{2n+1}}{\left[(2n+1)!\right]^2} \epsilon \eta^{-N(\varepsilon)} O(\gamma_n(T,\delta)) =$$

(2.2.2.75.b)

$$2\eta^{-N(\varepsilon)} \sum_{n=1}^{M} \frac{\epsilon^{2n+2}}{\left[(2n+1)!\right]^2} O(\widetilde{\gamma}_n(T)).$$

Finally from (2.2.2.69) and Eqs.(2.2.2.75.a,b) we obtain

$$\chi_1(T,\epsilon,\varepsilon,\eta) \leq$$

$$4U(T) + 4\mathcal{F}_{1,\varepsilon}(T,x_0,\epsilon,\eta) + 4\mathcal{F}_{2,\varepsilon}(T,x_0,\epsilon,\eta) =$$

$$4U(T) +$$

(2.2.2.76)

$$8\eta^{-N(\varepsilon)} \sum_{n=1}^{M} \frac{\epsilon^{2n+1}}{\left[(2n)!\right]^2} O\left(\widetilde{\gamma}_{1,n}(T)\right) + 8\eta^{-N(\varepsilon)} \sum_{n=1}^{M} \frac{\epsilon^{2n+2}}{\left[(2n+1)!\right]^2} O(\widetilde{\gamma}_n(T)) =$$

$$4U(T) + 8\eta^{-N(\varepsilon)} \sum_{n=1}^{M} \left[ \frac{\epsilon^{2n+1}}{\left[(2n)!\right]^2} O\left(\widetilde{\gamma}_{1,n}(T)\right) + \sum_{n=1}^{M} \frac{\epsilon^{2n+2}}{\left[(2n+1)!\right]^2} O(\widetilde{\gamma}_n(T)) \right].$$

Finally we obtain the master inequality:

**Theorem 2.2.2.1**.

$$|\langle x, t, 0; \epsilon, \delta \rangle| \leq U(T)(C + O(\delta)).$$  (2.2.2.77)

Here

$$U(T) = \left| \frac{1}{m\varpi} \int_0^T b(t) \sin \varpi(T - t) dt \right|.$$  (2.2.2.78)

Let us considered the quantity $\langle x_f, t_f, x_0; \epsilon \rangle - \lambda, \lambda \in \mathbb{R}$. We assume that $\int dx_f |\mathbf{K}(x_f, t_f | x_0, 0)|^2 = 1$. Then from Eq.(2.2.2.2) we obtain

$$\langle T, x_0, \lambda; \epsilon \rangle = \langle T, x_0; \epsilon \rangle - \lambda =$$

$$= \int dx_f x_f |\mathbf{K}(x_f, T | x_0, 0)|^2 - \lambda \int dx_f |\mathbf{K}(x_f, T | x_0, 0)|^2$$  (2.2.2.79)

$$= \int dx_f (x_f - \lambda) |\mathbf{K}(x_f, T | x_0, 0)|^2.$$

$$T = t_f.$$

From Eq.(2.2.2.79) we obtain

$$|\langle T, x_0, \lambda; \varepsilon, \epsilon \rangle| \leq \int dx |x - \lambda| \times$$

$$\left| \overline{\int_{q(T)=x} \widetilde{D}_\eta^+[q(t), \varepsilon] \Psi_0(q(0))} \exp\left[ \frac{i}{\epsilon} \left( \int_0^T dt \left[ \frac{m}{2} \left( \frac{d\widetilde{q}_\varepsilon(t)}{dt} \right)^2 - V(q_\varepsilon(t)) \right] \right) \right] \right|^2 =$$

$$\int dx \left\{ \left| \overline{\int_{q(T)=x} \widetilde{D}_\eta^+[q(t), \varepsilon] \Psi_0(q(0))} \sqrt{|x - \lambda|} \right. \right.$$

$$\left. \exp\left[ \frac{i}{\epsilon} \left( \int_0^T dt \left[ \frac{m}{2} \left( \frac{d\widetilde{q}_\varepsilon(t)}{dt} \right)^2 - V(q_\varepsilon(t)) \right] \right) \right] \right|^2 \right\} = \qquad (2.2.2.80)$$

$$\int dx \left\{ \left| \overline{\int_{q(T)=x} D^+[q(t), \varepsilon] \Psi_0(q(0))} \sqrt{|q(T) - \lambda|} \times \right. \right.$$

$$\exp\left\{ -\frac{m}{\epsilon} \int_0^T dt \left[ \{q(t)\}_\eta + \eta g(t, q(T), q(0), \lambda) \right]^2 \right\} \times$$

$$\left. \exp\left[ \frac{i}{\epsilon} \left( \int_0^T dt \left[ \frac{m}{2} \left( \frac{d\widetilde{q}_\varepsilon(t)}{dt} \right)^2 - V(q_\varepsilon(t)) \right] \right) \right] \right|^2 \right\}.$$

Here the function $V(q_\varepsilon(t), t)$ is given via formula

$$V(q_\varepsilon(t), t) = a_{1,\varepsilon}(t) q_\varepsilon(t) + a_{2,\varepsilon}(t) q_\varepsilon^2(t) + a_{3,\varepsilon} q_\varepsilon^3(t) + \ldots + a_{k,\varepsilon} q_\varepsilon^k(t),$$

$$q_\varepsilon(t) = \frac{q(t)}{1 + \varepsilon^k q^m(t)}, \qquad (2.2.2.80')$$

$$k \geq 1, m \geq 2$$

From Eq.(2.2.2.80) we obtain

$$|\langle T, x_0, \lambda; \varepsilon, \epsilon \rangle| \leq$$

$$\int dx \Bigg\{ \Bigg[ \int\limits_{q(T)=x} \widetilde{D}^+_\eta[q(t),\varepsilon]\Psi_0(q(0))\sqrt{|q(T)-\lambda|} \times$$

$$\cos\Bigg[\frac{1}{\epsilon}\Bigg(\int\limits_0^T dt\Bigg[\frac{m}{2}\Bigg(\frac{d\widetilde{q}_\varepsilon(t)}{dt}\Bigg)^2 - V(q_\varepsilon(t))\Bigg]\Bigg)\Bigg]\Bigg]^2 +$$

$$\Bigg[ \int\limits_{q(T)=x} \widetilde{D}^+_\eta[q(t),\varepsilon]\Psi_0(q(0))\sqrt{|q(T)-\lambda|} \times$$

$$\sin\Bigg[\frac{1}{\epsilon}\Bigg(\int\limits_0^T dt\Bigg[\frac{m}{2}\Bigg(\frac{d\widetilde{q}_\varepsilon(t)}{dt}\Bigg)^2 - V(q_\varepsilon(t))\Bigg]\Bigg)\Bigg]\Bigg]^2 \Bigg\} =$$

$$(2.2.2.81)$$

$$\int dx \Bigg\{ \Bigg[ \int\limits_{q(T)=x} \widetilde{D}^+_\eta[q(t),\varepsilon]\Psi_0(q(0))\sqrt{|q(T)-\lambda|} \times$$

$$\cos\Bigg[\frac{1}{\epsilon}\Bigg(\int\limits_0^T dt\Bigg[\frac{m}{2}\Bigg(\frac{d\widetilde{q}_\varepsilon(t)}{dt}\Bigg)^2 - V(q_\varepsilon(t))\Bigg]\Bigg)\Bigg]\Bigg]^2 \Bigg\} +$$

$$\int dx \Bigg\{ \Bigg[ \int\limits_{q(T)=x} \widetilde{D}^+_\eta[q(t),\varepsilon]\Psi_0(q(0))\sqrt{|q(T)-\lambda|} \times$$

$$\sin\Bigg[\frac{1}{\epsilon}\Bigg(\int\limits_0^T dt\Bigg[\frac{m}{2}\Bigg(\frac{d\widetilde{q}_\varepsilon(t)}{dt}\Bigg)^2 - V(q_\varepsilon(t))\Bigg]\Bigg)\Bigg]\Bigg]^2 \Bigg\} =$$

$$\chi_1(T,x_0,\lambda,\varepsilon,\epsilon,\eta) + \chi_{2,\varepsilon}(T,x_0,\lambda,\varepsilon,\epsilon,\eta).$$

Here

$$\chi_1(T, x_0, \lambda, \varepsilon, \epsilon, \eta) =$$

$$\int dx \left\{ \left[ \int_{q(T)=x} \widetilde{D}_\eta^+[q(t), \varepsilon] \Psi_0(q(0)) \sqrt{|q(T) - \lambda|} \times \right. \right.$$

$$\exp\left\{ -\frac{m}{\epsilon} \int_0^T dt \left[ \{q(t)\}_\eta + \eta g(t, q(T), q(0), \lambda) \right]^2 \right\} \times \qquad (2.2.2.82.a)$$

$$\left. \left. \cos\left[ \frac{1}{\epsilon} \left( \int_0^T dt \left[ \frac{m}{2} \left( \frac{d\widetilde{q}_\varepsilon(t)}{dt} \right)^2 - V(q_\varepsilon(t)) \right] \right) \right] \right]^2 \right\},$$

and

$$\chi_2(T, x_0, \lambda, \varepsilon, \epsilon, \eta) =$$

$$\int dx \left\{ \left[ \int_{q(T)=x} \widetilde{D}_\eta^+[q(t), \varepsilon] \Psi_0(q(0)) \sqrt{|q(T) - \lambda|} \times \right. \right.$$

$$\exp\left\{ -\frac{m}{\epsilon} \int_0^T dt \left[ \{q(t)\}_\eta + \eta g(t, q(T), q(0), \lambda) \right]^2 \right\} \times \qquad (2.2.2.82.b)$$

$$\left. \left. \sin\left[ \frac{1}{\epsilon} \left( \int_0^T dt \left[ \frac{m}{2} \left( \frac{d\widetilde{q}_\varepsilon(t)}{dt} \right)^2 - V(q_\varepsilon(t)) \right] \right) \right] \right]^2 \right\}.$$

By using replacement $q(t) - \lambda = u(t)$, from Eq.(2.2.2.80) we obtain

$$\chi_1(T, x_0, \lambda, \varepsilon, \epsilon, \eta) =$$

$$\int dx \left\{ \left[ \int_{u(T)=x-\lambda} D^+[u(t), \varepsilon] \Psi_0(u(0) + \lambda) \sqrt{|u(T)|} \times \right.\right.$$

$$\exp \left\{ -\frac{m}{\epsilon} \int_0^T dt \Big[ \{u(t) + \lambda\}_\eta + \eta g(t, u(T) + \lambda, u(0) + \lambda, \lambda) \Big]^2 \right\} \times$$

$$\cos \left[ \frac{1}{\epsilon} \left( \int_0^T dt \Big[ \frac{m}{2} \Big( \frac{d\widetilde{u}_\varepsilon(t)}{dt} \Big)^2 - V(u_\varepsilon(t) + \lambda) \Big] \right) \right] \Bigg] \Bigg\}^2 ,$$

$$(2.2.2.83.a)$$

and

$$\chi_2(T, x_0, \lambda, \varepsilon, \epsilon, \eta) =$$

$$\int dx \left\{ \left[ \int_{u(T)=x-\lambda} D^+[u(t), \varepsilon] \Psi_0(u(0) + \lambda) \sqrt{|u(T)|} \times \right.\right.$$

$$\exp \left\{ -\frac{m}{\epsilon} \int_0^T dt \Big[ \{u(t) + \lambda\}_\eta + \eta g(t, u(T) + \lambda, u(0) + \lambda, \lambda) \Big]^2 \right\} \times$$

$$\sin \left[ \frac{1}{\epsilon} \left( \int_0^T dt \Big[ \frac{m}{2} \Big( \frac{d\widetilde{u}_\varepsilon(t)}{dt} \Big)^2 - V(u_\varepsilon(t) + \lambda) \Big] \right) \right] \Bigg] \Bigg\}^2 ,$$

$$(2.2.2.83.b)$$

Let us rewrite a function $V(u_\varepsilon(t) + \lambda, t)$ in equivalent form:

$$V(u_\varepsilon(t) + \lambda, t) = V_0(u_\varepsilon(t), t, \lambda) + V_1(u_\varepsilon(t), t, \lambda)$$

$$V_0(u_\varepsilon(t), t, \lambda) = a_\varepsilon(\lambda)u_\varepsilon^2(t) + b_\varepsilon(t, \lambda)(q_\varepsilon(t)),$$

$$V_1(u_\varepsilon(t), t, \lambda) = a_{3,\varepsilon}(\lambda, t)u_\varepsilon^3(t) + \ldots + a_{k,\varepsilon}(\lambda, t)u_\varepsilon^k(t) + c(t, \lambda),$$

$$a_\varepsilon(\lambda, t) > 0,$$

$$u_\varepsilon(t) = \frac{u(t)}{1 + \varepsilon^k u^m(t)},$$

$$k \geq 1, m \geq 2$$

$$(2.2.2.84)$$

Let us now evaluate a quantities $(\chi_1(T, x_0, \lambda, \epsilon, \eta))_\varepsilon$:

$$\chi_1(T, x_0, \lambda, \varepsilon, \epsilon, \eta) =$$

$$\int dx \Bigg\{ \Bigg[ \int\limits_{u(T)=x-\lambda} D^+[u(t), \varepsilon] \Psi_0(u(0) + \lambda) \sqrt{|u(T)|} \ \times$$

$$\exp\Bigg\{ -\frac{m}{\epsilon} \int\limits_0^T dt \Big[ \{u(t) + \lambda\}_\eta + \eta g(t, u(T) + \lambda, u(0) + \lambda, \lambda) \Big]^2 \Bigg\} \times \qquad (2.2.2.85.a)$$

$$\cos\Bigg[ \frac{1}{\epsilon} \Bigg( \int\limits_0^T dt \Big[ \frac{m}{2}\Big( \frac{d\widetilde{u}_\varepsilon(t)}{dt} \Big)^2 - V(u_\varepsilon(t) + \lambda) \Big] \Bigg) \Bigg] \Bigg]^2 =$$

$$\int dx \Pi_1^2(T, x_0, \lambda, \epsilon, \varepsilon, \eta)$$

and $\chi_2(T, x_0, \lambda, \varepsilon, \epsilon, \eta)$ :

$$\chi_2(T, x_0, \lambda, \varepsilon, \epsilon, \eta) =$$

$$\int dx \Bigg\{ \Bigg[ \int\limits_{u(T)=x-\lambda} D^+[u(t), \varepsilon] \Psi_0(u(0) + \lambda) \sqrt{|u(T)|} \ \times$$

$$\exp\Bigg\{ -\frac{m}{\epsilon} \int\limits_0^T dt \Big[ \{u(t) + \lambda\}_\eta + \eta g(t, u(T) + \lambda, u(0) + \lambda, \lambda) \Big]^2 \Bigg\} \times \qquad (2.2.2.85.b)$$

$$\sin\Bigg[ \frac{1}{\epsilon} \Bigg( \int\limits_0^T dt \Big[ \frac{m}{2}\Big( \frac{d\widetilde{u}_\varepsilon(t)}{dt} \Big)^2 - V(u_\varepsilon(t) + \lambda) \Big] \Bigg) \Bigg] \Bigg]^2 =$$

$$\int dx \Pi_2^2(T, x_0, \lambda, \epsilon, \varepsilon, \eta).$$

Here

$$\Pi_1(T, x_0, \lambda, \epsilon, \varepsilon, \eta) =$$

$$\overline{\int\limits_{u(T)=x-\lambda} D^+[u(t), \varepsilon] \Psi_0(u(0) + \lambda) \sqrt{|u(T)|}} \times$$

$$\exp\left\{ -\frac{m}{\epsilon} \int\limits_0^T dt \Big[ \{u(t) + \lambda\}_\eta + \eta g(t, u(T) + \lambda, u(0) + \lambda, \lambda) \Big]^2 \right\} \times \qquad (2.2.2.86.a)$$

$$\cos\left[ \frac{1}{\epsilon} \left( \int\limits_0^T dt \Big[ \frac{m}{2} \Big( \frac{d\widetilde{u}_\varepsilon(t)}{dt} \Big)^2 - V(u_\varepsilon(t) + \lambda) \Big] \right) \right].$$

and

$$\Pi_2(T, x_0, \lambda, \epsilon, \varepsilon, \eta) =$$

$$\overline{\int\limits_{u(T)=x-\lambda} D^+[u(t), \varepsilon] \Psi_0(u(0) + \lambda) \sqrt{|u(T)|}} \times$$

$$\exp\left\{ -\frac{m}{\epsilon} \int\limits_0^T dt \Big[ \{u(t) + \lambda\}_\eta + \eta g(t, u(T) + \lambda, u(0) + \lambda, \lambda) \Big]^2 \right\} \times \qquad (2.2.2.86.b)$$

$$\sin\left[ \frac{1}{\epsilon} \left( \int\limits_0^T dt \Big[ \frac{m}{2} \Big( \frac{d\widetilde{u}_\varepsilon(t)}{dt} \Big)^2 - V(u_\varepsilon(t) + \lambda) \Big] \right) \right].$$

Substitution Eq.(2.2.2.84) into Eq.(2.2.2.86.a,b) gives

$$(\Pi_1(T, x, x_0, \lambda, \epsilon, \varepsilon, \eta))_\varepsilon =$$

$$\overline{\int\limits_{u(T)=x-\lambda} D^+[u(t), \varepsilon] \Psi_0(u(0) + \lambda)} \sqrt{|u(T)|} \times$$

$$\exp\left\{-\frac{m}{\epsilon} \int\limits_0^T dt \Big[ \{u(t) + \lambda\}_\eta + \eta g(t, u(T) + \lambda, u(0) + \lambda, \lambda) \Big]^2 \right\} \times$$

$$\cos\left[ \frac{1}{\epsilon}\left( \int\limits_0^T dt \Big[ \frac{m}{2}\left( \frac{d\widetilde{u}_\varepsilon(t)}{dt} \right)^2 - V_0(u_\varepsilon(t), t, \lambda) - V_1(u_\varepsilon(t), t, \lambda) \Big] \right) \right] =$$

$$\overline{\int\limits_{u(T)=x-\lambda} \widetilde{D}_\eta^+[u(t), \varepsilon] \Psi_0(u(0) + \lambda)} \sqrt{|u(T)|} \times$$

$$\cos\left[ \frac{1}{\epsilon}\left( \int\limits_0^T dt \Big[ \frac{m}{2}\left( \frac{d\widetilde{u}_\varepsilon(t)}{dt} \right)^2 - V_0(u_\varepsilon(t), t, \lambda) \Big] \right) \right] \times \qquad (2.2.2.87.a)$$

$$\cos\left[ \frac{1}{\epsilon}\left( \int\limits_0^T V_1(u_\varepsilon(t), t, \lambda) dt \right) \right] +$$

$$\overline{\int\limits_{u(T)=x-\lambda} \widetilde{D}_\eta^+[u(t), \varepsilon] \Psi_0(q(0) + \lambda)} \sqrt{|u(T)|} \times$$

$$\sin\left[ \frac{1}{\epsilon}\left( \int\limits_0^T dt \Big[ \frac{m}{2}\left( \frac{d\widetilde{u}_\varepsilon(t)}{dt} \right)^2 - V_0(u_\varepsilon(t), t, \lambda) \Big] \right) \right] \times$$

$$\sin\left[ -\frac{1}{\epsilon}\left( \int\limits_0^T V_1(u_\varepsilon(t), t, \lambda) dt \right) \right] =$$

$$v_{1,\varepsilon}(T, x, x_0, \lambda, \epsilon, \eta) + v_{1,\varepsilon}(T, x, x_0, \lambda, \epsilon, \eta).$$

Here

$$\nu_{1,\varepsilon}(T,x,x_0,\lambda,\epsilon,\eta) =$$

$$\overline{\int\limits_{u(T)=x-\lambda} \widetilde{D}_{\eta}^{+}[u(t),\varepsilon]\Psi_0(u(0)+\lambda)\,\sqrt{|u(T)|}} \times$$

$$\cos\left[\frac{1}{\epsilon}\left(\int\limits_0^T dt\left[\frac{m}{2}\left(\frac{d\widetilde{u}_\varepsilon(t)}{dt}\right)^2 - V_0(u_\varepsilon(t),t,\lambda)\right]\right)\right] \times$$

$$\cos\left[\frac{1}{\epsilon}\left(\int\limits_0^T V_1(u_\varepsilon(t),t,\lambda)dt\right)\right] =$$

$$\overline{\int} dy\left\{\overline{\int\limits_{u(T)=x-\lambda,u(0)=y-\lambda} \widetilde{D}_{\eta}^{+}[u(t),\varepsilon]\Psi_0(u(0)+\lambda)\,\sqrt{|u(T)|}} \times\right.$$

$$\sin\left[\frac{1}{\epsilon}\left(\int\limits_0^T dt\left[\frac{m}{2}\left(\frac{d\widetilde{u}_\varepsilon(t)}{dt}\right)^2 - V_0(u_\varepsilon(t),t,\lambda)\right]\right)\right] \times$$

$$\left.\cos\left[\frac{1}{\epsilon}\left(\int\limits_0^T V_1(u_\varepsilon(t),t,\lambda)dt\right)\right]\right\} = \qquad (2.2.2.88.a)$$

$$\overline{\int} dy\left\{\overline{\int\limits_{u(T)=x-\lambda,u(0)=y-\lambda} D^{+}[u(t),\varepsilon]\Psi_0(y)\,\sqrt{|u(T)|}} \times\right.$$

$$\exp\left\{-\frac{m}{\epsilon}\int\limits_0^T dt\left[\left\{u(t)+\lambda\right\}_\eta + \eta g(t,u(T)+\lambda,u(0)+\lambda,\lambda)\right]^2\right\} \times$$

$$\cos\left[\frac{1}{\epsilon}\left(\int\limits_0^T dt\left[\frac{m}{2}\left(\frac{d\widetilde{u}_\varepsilon(t)}{dt}\right)^2 - V_0(u_\varepsilon(t),t,\lambda)\right]\right)\right] \times$$

$$\left.\cos\left[\frac{1}{\epsilon}\left(\int\limits_0^T V_1(u_\varepsilon(t),t,\lambda)dt\right)\right]\right\}$$

and

$$v_{2,\varepsilon}(T, x, x_0, \lambda, \epsilon, \eta) =$$

$$\overline{\int\limits_{u(T)=x-\lambda} \widetilde{D}_\eta^+[u(t),\varepsilon]\Psi_0(u(0)+\lambda)\sqrt{|u(T)|}} \times$$

$$\sin\left[\frac{1}{\epsilon}\left(\int\limits_0^T dt\left[\frac{m}{2}\left(\frac{d\widetilde{u}_\varepsilon(t)}{dt}\right)^2 - V_0(u_\varepsilon(t),t,\lambda)\right]\right)\right] \times$$

$$\sin\left[-\frac{1}{\epsilon}\left(\int\limits_0^T V_1(u_\varepsilon(t),t,\lambda)dt\right)\right] =$$

$$\overline{\int} dy \left\{\overline{\int\limits_{u(T)=x-\lambda,u(0)=y-\lambda} \widetilde{D}_\eta^+[u(t),\varepsilon]\Psi_0(u(0)+\lambda)\sqrt{|u(T)|}} \times\right.$$

$$\sin\left[\frac{1}{\epsilon}\left(\int\limits_0^T dt\left[\frac{m}{2}\left(\frac{d\widetilde{u}_\varepsilon(t)}{dt}\right)^2 - V_0(u_\varepsilon(t),t,\lambda)\right]\right)\right] \times \qquad (2.2.2.88.b)$$

$$\left.\sin\left[-\frac{1}{\epsilon}\left(\int\limits_0^T V_1(u_\varepsilon(t),t,\lambda)dt\right)\right]\right\} =$$

$$\overline{\int} dy \left\{\overline{\int\limits_{u(T)=x-\lambda,u(0)=y-\lambda} D^+[u(t),\varepsilon]\Psi_0(y)\sqrt{|u(T)|}} \times\right.$$

$$\exp\left\{-\frac{m}{\epsilon}\int\limits_0^T dt\left[\{u(t)+\lambda\}_\eta + \eta g(t,u(T)+\lambda,u(0)+\lambda,\lambda)\right]^2\right\} \times$$

$$\sin\left[\frac{1}{\epsilon}\left(\int\limits_0^T dt\left[\frac{m}{2}\left(\frac{d\widetilde{u}_\varepsilon(t)}{dt}\right)^2 - V_0(u_\varepsilon(t),t,\lambda)\right]\right)\right] \times$$

$$\left.\sin\left[-\frac{1}{\epsilon}\left(\int\limits_0^T V_1(u_\varepsilon(t),t,\lambda)dt\right)\right]\right\}.$$

Let us now evaluate a quantities $v_{1,\varepsilon}(T,x_0,\lambda,\varepsilon,\epsilon,\eta)$ and $v_{2,\varepsilon}(T,x_0,\lambda,\varepsilon,\epsilon,\eta)$. By using replacement $u(t) \to \sqrt{\epsilon}\, u(t),$ from Eq.(2.2.2.88.a,b) we obtain

$$\nu_{1,\varepsilon}(T,x,x_0,\lambda,\epsilon,\eta) =$$

$$\bar{\int} dy \left\{ \overline{\int\limits_{u(T)=\frac{x-\lambda}{\sqrt{\epsilon}},\,u(0)=\frac{y-\lambda}{\sqrt{\epsilon}}} \widetilde{D}^+_{\eta}[u(t),\varepsilon]\Psi_0(y)\sqrt{\sqrt{\epsilon}\,|u(T)|}} \times \right.$$

$$\cos\left[ \frac{1}{\epsilon}\left( \int\limits_0^T dt \left[ \frac{m}{2}\left( \frac{\sqrt{\epsilon}\,d\tilde{u}_\varepsilon(t)}{dt} \right)^2 - V_0(\sqrt{\epsilon}\,u_\varepsilon(t),t,\lambda) \right] \right) \right] \times$$

$$\left. \cos\left[ \frac{1}{\epsilon}\left( \int\limits_0^T \widehat{V}_1(\sqrt{\epsilon}\,u_\varepsilon(t),t,\lambda)dt \right) \right] \right\} =$$

$$\bar{\int} dy \left\{ \overline{\int\limits_{u(T)=\frac{x-\lambda}{\sqrt{\epsilon}},\,u(0)=\frac{y-\lambda}{\sqrt{\epsilon}}} \widetilde{D}^+_{\eta}[u(t),\varepsilon]\Psi_0(y)\sqrt{\sqrt{\epsilon}\,|u(T)|}} \times \right.$$

$$\cos\left[ \left( \int\limits_0^T dt \left[ \frac{m}{2}\left( \frac{d\tilde{u}_{\varepsilon\epsilon}(t)}{dt} \right)^2 - \frac{1}{\epsilon}V_0(\sqrt{\epsilon}\,u_{\nu(\varepsilon,\epsilon)}(t),t,\lambda) \right] \right) \right] \times \qquad (2.2.2.89.a)$$

$$\left. \cos\left[ \frac{1}{\epsilon}\left( \int\limits_0^T \widehat{V}_1(\sqrt{\epsilon}\,u_{\nu(\varepsilon,\epsilon)}(t),t,\lambda)dt \right) \right] \right\} =$$

$$\bar{\int} dy \left\{ \overline{\int\limits_{u(T)=\frac{x-\lambda}{\sqrt{\epsilon}},\,u(0)=\frac{y-\lambda}{\sqrt{\epsilon}}} D^+[u(t),\varepsilon]\Psi_0(y)\sqrt{\sqrt{\epsilon}\,|u(T)|}} \times \right.$$

$$\exp\left\{ -\frac{m}{\epsilon}\int\limits_0^T dt \left[ \{u(t)+\lambda\}_\eta + \eta g(t,\sqrt{\epsilon}\,u(T)+\lambda,\sqrt{\epsilon}\,u(0)+\lambda,\lambda) \right]^2 \right\} \times$$

$$\cos\left[ \left( \int\limits_0^T dt \left[ \frac{m}{2}\left( \frac{d\tilde{u}_{\varepsilon\epsilon}(t)}{dt} \right)^2 - \frac{1}{\epsilon}V_0(\sqrt{\epsilon}\,u_{\nu(\varepsilon,\epsilon)}(t),t,\lambda) \right] \right) \right] \times$$

and

$$\nu_{2,\varepsilon}(T,x,x_0,\lambda,\epsilon,\eta) =$$

$$\overline{\int} dy \left\{ \overline{\int\limits_{u(T)=\frac{x-\lambda}{\sqrt{\epsilon}},u(0)=\frac{y-\lambda}{\sqrt{\epsilon}}}} \widetilde{D}^{+}_{\eta}[u(t),\varepsilon]\Psi_0(y)\sqrt{\sqrt{\epsilon}\,|u(T)|}\,\times \right.$$

$$\sin\left[\frac{1}{\epsilon}\left(\int\limits_0^T dt\left[\frac{m}{2}\left(\frac{\sqrt{\epsilon}\,d\widetilde{u}_\varepsilon(t)}{dt}\right)^2 - V_0(\sqrt{\epsilon}\,u_\varepsilon(t),t,\lambda)\right]\right)\right]\times$$

$$\left.\sin\left[-\frac{1}{\epsilon}\left(\int\limits_0^T \widehat{V}_1(\sqrt{\epsilon}\,u_\varepsilon(t),t,\lambda)dt\right)\right]\right\} =$$

$$\overline{\int} dy \left\{ \overline{\int\limits_{u(T)=\frac{x-\lambda}{\sqrt{\epsilon}},u(0)=\frac{y-\lambda}{\sqrt{\epsilon}}}} \widetilde{D}^{+}_{\eta}[u(t),\varepsilon]\Psi_0(y)\sqrt{\sqrt{\epsilon}\,|u(T)|}\,\times \right.$$

$$\sin\left[\left(\int\limits_0^T dt\left[\frac{m}{2}\left(\frac{d\widetilde{u}_{\varepsilon\epsilon}(t)}{dt}\right)^2 - \frac{1}{\epsilon}V_0(\sqrt{\epsilon}\,u_{\nu(\varepsilon,\epsilon)}(t),t,\lambda)\right]\right)\right]\times \qquad (2.2.2.89.b)$$

$$\left.\sin\left[\left(-\int\limits_0^T \widehat{V}_1(\sqrt{\epsilon}\,u_{\nu(\varepsilon,\epsilon)}(t),t,\lambda)dt\right)\right]\right\} =$$

$$\overline{\int} dy \left\{ \overline{\int\limits_{u(T)=\frac{x-\lambda}{\sqrt{\epsilon}},u(0)=\frac{y-\lambda}{\sqrt{\epsilon}}}} D^{+}[u(t),\varepsilon]\Psi_0(y)\sqrt{\sqrt{\epsilon}\,|u(T)|}\,\times \right.$$

$$\exp\left\{-\frac{m}{\epsilon}\int\limits_0^T dt\left[\{u(t)+\lambda\}_\eta + \eta g(t,\sqrt{\epsilon}\,u(T)+\lambda,\sqrt{\epsilon}\,u(0)+\lambda,\lambda)\right]^2\right\}\times$$

$$\sin\left[\left(\int\limits_0^T dt\left[\frac{m}{2}\left(\frac{d\widetilde{u}_{\varepsilon\epsilon}(t)}{dt}\right)^2 - \frac{1}{\epsilon}V_0(\sqrt{\epsilon}\,u_{\nu(\varepsilon,\epsilon)}(t),t,\lambda)\right]\right)\right]\times$$

From Eqs.(2.2.2.84) and Eqs.(2.2.2.89.a,b) we obtain

$$v_{1,\varepsilon}(T,x,x_0,\lambda,\epsilon,\eta) =$$

$$\overline{\int} dy \left\{ \overline{\int_{u(T)=\frac{x-\lambda}{\sqrt{\epsilon}},\, u(0)=\frac{y-\lambda}{\sqrt{\epsilon}}}} \widetilde{D}_{\eta}^{+}[u(t),\varepsilon]\Psi_0(y)\sqrt{\sqrt{\epsilon}\,|u(T)|} \times \right.$$

$$\cos\left[\left(\int_0^T dt\left[\frac{m}{2}\left(\frac{d\widetilde{u}_{\varepsilon\epsilon}(t)}{dt}\right)^2 - a_\varepsilon(t,\lambda)u_{v(\varepsilon,\epsilon)}^2(t) + b_\varepsilon(t,\lambda)\frac{u_{v(\varepsilon,\epsilon)}(t)}{\sqrt{\epsilon}}\right]\right)\right] \times$$

$$\left. \cos\left[\left(\sqrt{\epsilon}\int_0^T \widetilde{V}_1(u_{v(\varepsilon,\epsilon)}(t),t,\lambda)dt\right)\right]\right\} =$$

$$\overline{\int} dy \left\{ \overline{\int_{u(T)=\frac{x-\lambda}{\sqrt{\epsilon}},\, u(0)=\frac{y-\lambda}{\sqrt{\epsilon}}}} D^{+}[u(t),\varepsilon]\Psi_0(y)\sqrt{\sqrt{\epsilon}\,|u(T)|} \times \right.$$ 

(2.2.2.90.a)

$$\exp\left\{-\frac{m}{\epsilon}\int_0^T dt\left[\{u(t)+\lambda\}_\eta + \eta g(t,\sqrt{\epsilon}\,u(T)+\lambda,\sqrt{\epsilon}\,u(0)+\lambda,\lambda)\right]^2\right\} \times$$

$$\cos\left[\left(\int_0^T dt\left[\frac{m}{2}\left(\frac{d\widetilde{u}_{\varepsilon\epsilon}(t)}{dt}\right)^2 - a_\varepsilon(t,\lambda)u_{v(\varepsilon,\epsilon)}^2(t) + b_\varepsilon(t,\lambda)\frac{u_{v(\varepsilon,\epsilon)}(t)}{\sqrt{\epsilon}}\right]\right)\right] \times$$

$$\left. \cos\left[\left(\sqrt{\epsilon}\int_0^T \widetilde{V}_1(u_{v(\varepsilon,\epsilon)}(t),t,\lambda)dt\right)\right]\right\}$$

and

$$v_{2,\varepsilon}(T,x,x_0,\lambda,\epsilon,\eta) =$$

$$\overline{\int} dy \left\{ \overline{\int_{u(T)=\frac{x-\lambda}{\sqrt{\epsilon}},\,u(0)=\frac{y-\lambda}{\sqrt{\epsilon}}}} \widetilde{D}_\eta^+[u(t),\varepsilon]\Psi_0(y)\sqrt{\sqrt{\epsilon}\,|u(T)|} \times \right.$$

$$\sin\left[\left(\int_0^T dt\left[\frac{m}{2}\left(\frac{d\widetilde{u}_{\varepsilon\epsilon}(t)}{dt}\right)^2 - a_\varepsilon(t,\lambda)u_{\nu(\varepsilon,\epsilon)}^2(t) + b_\varepsilon(t,\lambda)\frac{u_{\nu(\varepsilon,\epsilon)}(t)}{\sqrt{\epsilon}}\right]\right)\right] \times$$

$$\left.\sin\left[\left(-\sqrt{\epsilon}\int_0^T \widetilde{V}_1(u_{\nu(\varepsilon,\epsilon)}(t),t,\lambda)dt\right)\right]\right\} =$$

$$(2.2.90.b)$$

$$\overline{\int} dy \left\{ \overline{\int_{u(T)=\frac{x-\lambda}{\sqrt{\epsilon}},\,u(0)=\frac{y-\lambda}{\sqrt{\epsilon}}}} D^+[u(t),\varepsilon]\Psi_0(y)\sqrt{\sqrt{\epsilon}\,|u(T)|} \times \right.$$

$$\exp\left\{-\frac{m}{\epsilon}\int_0^T dt\left[\{u(t)+\lambda\}_\eta + \eta g(t,\sqrt{\epsilon}\,u(T)+\lambda,\sqrt{\epsilon}\,u(0)+\lambda,\lambda)\right]^2\right\} \times$$

$$\sin\left[\left(\int_0^T dt\left[\frac{m}{2}\left(\frac{d\widetilde{u}_{\varepsilon\epsilon}(t)}{dt}\right)^2 - a_\varepsilon(t,\lambda)u_{\nu(\varepsilon,\epsilon)}^2(t) + b_\varepsilon(t,\lambda)\frac{u_{\nu(\varepsilon,\epsilon)}(t)}{\sqrt{\epsilon}}\right]\right)\right] \times$$

$$\left.\sin\left[\left(-\sqrt{\epsilon}\int_0^T \widetilde{V}_1(u_{\nu(\varepsilon,\epsilon)}(t),t,\lambda)dt\right)\right]\right\}.$$

Here

$$\nu(\varepsilon, \epsilon) = \varepsilon^k \epsilon^{\frac{m}{2}},$$

$$\frac{1}{\epsilon} V_1(\sqrt{\epsilon}\, u_{\nu(\varepsilon,\epsilon)}(t), t, \lambda) = \frac{1}{\epsilon} V_1(\sqrt{\epsilon}\, u_{\nu(\varepsilon,\epsilon)}(t)t, \lambda) = a_3(t, \lambda) \epsilon^{1/2} u^3_{\nu(\varepsilon,\epsilon)}(t) + \ldots$$

$$+ a_k(t, \lambda) \epsilon^{(k/2)-1} u^k_{\nu(\varepsilon,\epsilon)}(t) = \widehat{V}_1(\sqrt{\epsilon}\, u_{\nu(\varepsilon,\epsilon)}(t)t, \lambda) \qquad (2.2.91)$$

$$= \widehat{V}_1(\sqrt{\epsilon}\, u_{\nu(\varepsilon,\epsilon)}(t), t, \lambda) = \sqrt{\epsilon}\, \widetilde{V}_1(u_{\nu(\varepsilon,\epsilon)}(t), t, \lambda).$$

$$\widetilde{V}_1(u_{\nu(\varepsilon,\epsilon)}(t)) = a_3(t, \lambda) u^3_{\nu(\varepsilon,\epsilon)}(t) + \ldots + a_k(t, \lambda) \epsilon^{(k/2)-3/2} u^k_{\nu(\varepsilon,\epsilon)}(t).$$

Let us rewrite Eqs.(2.2.90.a,b) in the following equivalent form

$$\nu_{1,\varepsilon}(T, x, x_0, \lambda, \epsilon, \eta) =$$

$$\overline{\int} dy \left\{ \overline{\int_{u(T)=\frac{x-\lambda}{\sqrt{\epsilon}}, u(0)=\frac{y-\lambda}{\sqrt{\epsilon}}} \widetilde{D}^{(1)}_{\mathbf{sm}}[u(t), \varepsilon] \Psi_0(y) \sqrt{\sqrt{\epsilon}\, |u(T)|} \times} \right.$$

$$\left| \cos\left[ \left( \int_0^T dt \left[ \frac{m}{2} \left( \frac{d\widetilde{u}_{\varepsilon\epsilon}(t)}{dt} \right)^2 - a_\varepsilon(t, \lambda) u^2_{\nu(\varepsilon,\epsilon)}(t) + b_\varepsilon(t, \lambda) \frac{u_{\nu(\varepsilon,\epsilon)}(t)}{\sqrt{\epsilon}} \right] \right) \right] \right| \times \qquad (2.2.2.92.a)$$

$$\left. \cos\left[ \left( \sqrt{\epsilon} \int_0^T \widetilde{V}_1(u_{\nu(\varepsilon,\epsilon)}(t), t, \lambda) dt \right) \right] \right\}$$

and

$$\nu_{2,\varepsilon}(T,x,x_0,\lambda,\epsilon,\eta) =$$

$$\overline{\int} dy \left\{ \overline{\int_{u(T)=\frac{x-\lambda}{\sqrt{\epsilon}},u(0)=\frac{y-\lambda}{\sqrt{\epsilon}}}} \widetilde{D}_{\mathbf{sm}}^{(2)}[u(t),\varepsilon]\Psi_0(y)\sqrt{\sqrt{\epsilon}\,|u(T)|}\, \times \right.$$

$$\left| \sin\left[ \left( \int_0^T dt \left[ \frac{m}{2}\left( \frac{d\widetilde{u}_{\varepsilon\epsilon}(t)}{dt} \right)^2 - a_\varepsilon(t,\lambda)u_{\nu(\varepsilon,\epsilon)}^2(t) + b_\varepsilon(t,\lambda)\frac{u_{\nu(\varepsilon,\epsilon)}(t)}{\sqrt{\epsilon}} \right] \right) \right] \right| \times \qquad (2.2.2.92.b)$$

$$\left. \sin\left[ \left( -\sqrt{\epsilon}\int_0^T \widetilde{V}_1(u_{\nu(\varepsilon,\epsilon)}(t),t,\lambda)dt \right) \right] \right\}$$

Here $\widetilde{D}_{\mathbf{sm}}^{(1)}[u(t),\varepsilon,\epsilon]$ and $\widetilde{D}_{\mathbf{sm}}^{(2)}[u(t),\varepsilon,\epsilon]$ is the corresponding signed Feynman"measures", i.e.

$$\overline{\underset{u(T)=\frac{x-\lambda}{\sqrt{\epsilon}},\,u(0)=\frac{y-\lambda}{\sqrt{\epsilon}}}{\int}}\widetilde{\mathcal{D}}^{(1)}_{\mathbf{sm}}[u(t),\varepsilon,\epsilon]\Psi_0(y)\sqrt{\sqrt{\epsilon}\,|u(T)|}\times(\bullet)=$$

$$\overline{\underset{\substack{q(T)=\frac{x-\lambda}{\sqrt{\epsilon}}\\u(0)=\frac{y-\lambda}{\sqrt{\epsilon}}}}{\int}}\widetilde{\mathcal{D}}^{+}_{\eta}[u(t),\varepsilon,\epsilon]\Psi_0(y)\sqrt{\sqrt{\epsilon}\,|u(T)|}\times$$

$$[\cos^{+}(S_{\varepsilon}(T,x,\lambda,u(t),\epsilon))+\cos^{-}(S_{\varepsilon}(T,x,\lambda,u(t),\epsilon))]\times(\bullet)=$$

$$\overline{\underset{u(T)=\frac{x-\lambda}{\sqrt{\epsilon}},\,u(0)=\frac{y-\lambda}{\sqrt{\epsilon}}}{\int}}D^{+}[q(t),\varepsilon]\Psi_0(y)\sqrt{\sqrt{\epsilon}\,|u(T)|}\times \qquad (2.2.2.93.a)$$

$$[\cos^{+}(S_{\varepsilon}(T,x,\lambda,u(t),\epsilon))+\cos^{-}(S_{\varepsilon}(T,x,\lambda,u(t),\epsilon))]\times$$

$$\exp\left\{-\frac{1}{\epsilon}\int_0^T dt\Big[\left\{\sqrt{\epsilon}\,u(t)+\lambda\right\}_{\eta}-\eta g(t,\sqrt{\epsilon}\,u(T)+\lambda,\sqrt{\epsilon}\,u(0)+\lambda)\Big]^2\right\}\times(\bullet),$$

$$S_{\varepsilon}(T,x,\lambda,q(t),\epsilon)=\int_0^T dt\left[\frac{m}{2}\left(\frac{d\widetilde{q}_{\varepsilon\epsilon}(t)}{dt}\right)^2-a_{\varepsilon}(t,\lambda)q_{\varepsilon}^2(t)+b_{\varepsilon}(t,\lambda)\frac{q_{\nu(\varepsilon,\epsilon)}(t)}{\sqrt{\epsilon}}\right].$$

and

$$\overline{\int\limits_{u(T)=\frac{x-\lambda}{\sqrt{\epsilon}},u(0)=\frac{y-\lambda}{\sqrt{\epsilon}}}} \widetilde{D}_{\mathrm{sm}}^{(2)}[u(t),\varepsilon,\epsilon]\Psi_0(y)\sqrt{\sqrt{\epsilon}\,|u(T)|}\times(\bullet)=$$

$$\overline{\int\limits_{q(T)=\frac{x-\lambda}{\sqrt{\epsilon}},u(0)=\frac{y-\lambda}{\sqrt{\epsilon}}}} \widetilde{D}_{\eta}^{+}[u(t),\varepsilon,\epsilon]\Psi_0(y)\sqrt{\sqrt{\epsilon}\,|u(T)|}\times$$

$$[\sin^{+}(S_{\varepsilon}(T,x,\lambda,u(t),\epsilon))+\sin^{-}(S_{\varepsilon}(T,x,\lambda,u(t),\epsilon))]\times(\bullet)=$$

$$\overline{\int\limits_{u(T)=\frac{x-\lambda}{\sqrt{\epsilon}},u(0)=\frac{y-\lambda}{\sqrt{\epsilon}}}} D^{+}[q(t),\varepsilon]\Psi_0(y)\sqrt{\sqrt{\epsilon}\,|u(T)|}\times \qquad (2.2.2.93.b)$$

$$[\sin^{+}(S_{\varepsilon}(T,x,\lambda,u(t),\epsilon))+\sin^{-}(S_{\varepsilon}(T,x,\lambda,u(t),\epsilon))]\times$$

$$\exp\left\{-\frac{1}{\epsilon}\int\limits_0^T dt\Big[\left\{\sqrt{\epsilon}\,u(t)+\lambda\right\}_{\eta}-\eta g(t,\sqrt{\epsilon}\,u(T)+\lambda,\sqrt{\epsilon}\,u(0)+\lambda)\Big]^2\right\}\times(\bullet),$$

$$S_{\varepsilon}(T,x,\lambda,q(t),\epsilon)=\int\limits_0^T dt\left[\frac{m}{2}\left(\frac{d\widetilde{q}_{\varepsilon\epsilon}(t)}{dt}\right)^2-a_{\varepsilon}(t,\lambda)q_{\varepsilon}^2(t)+b_{\varepsilon}(t,\lambda)\frac{q_{\nu(\varepsilon,\epsilon)}(t)}{\sqrt{\epsilon}}\right].$$

Let us rewrite Eqs.(2.2.2.93.a,b) in the following equivalent form

$$\overline{\int_{u(T)=\frac{x-\lambda}{\sqrt{\epsilon}},\,u(0)=\frac{y-\lambda}{\sqrt{\epsilon}}}} \widetilde{D}^{(1)}_{\mathbf{sm}}[u(t),\varepsilon,\epsilon]\Psi_0(y)\sqrt{\sqrt{\epsilon}\,|u(T)|}\;\times(\cdot)=$$

$$\overline{\int_{q(T)=\frac{x-\lambda}{\sqrt{\epsilon}},\,u(0)=\frac{y-\lambda}{\sqrt{\epsilon}}}} \widetilde{D}^{+}_{\eta}[u(t),\varepsilon,\epsilon]\Psi_0(y)\sqrt{\sqrt{\epsilon}\,|u(T)|}\;\times$$

$$\left[\cos^{+}(S_{\varepsilon}(T,x,\lambda,u(t),\epsilon))+\cos^{-}(S_{\varepsilon}(T,x,\lambda,u(t),\epsilon))+1-1\right]\times(\cdot)=$$

$$\overline{\int_{\substack{u(T)=\frac{x-\lambda}{\sqrt{\epsilon}}\\ u(0)=\frac{y-\lambda}{\sqrt{\epsilon}}}}} \widetilde{D}^{(1)}_{\mathbf{m}}[u(t),\varepsilon,\epsilon]\Psi_0(y)\sqrt{\sqrt{\epsilon}\,|u(T)|}\;\times(\cdot)-$$

$$(2.2.2.94.a)$$

$$\overline{\int_{q(T)=\frac{x-\lambda}{\sqrt{\epsilon}},\,u(0)=\frac{y-\lambda}{\sqrt{\epsilon}}}} \widetilde{D}^{+}_{\eta}[u(t),\varepsilon,\epsilon]\Psi_0(y)\sqrt{\sqrt{\epsilon}\,|u(T)|}\;\times(\cdot)=$$

$$\overline{\int_{u(T)=\frac{x-\lambda}{\sqrt{\epsilon}},\,u(0)=\frac{y-\lambda}{\sqrt{\epsilon}}}} \widetilde{D}^{(1)}_{\mathbf{m}}[u(t),\varepsilon,\epsilon]\Psi_0(y)\sqrt{\sqrt{\epsilon}\,|u(T)|}\;\times(\cdot)$$

$$-\overline{\int_{u(T)=\frac{x-\lambda}{\sqrt{\epsilon}},\,u(0)=\frac{y-\lambda}{\sqrt{\epsilon}}}} D^{+}[q(t),\varepsilon]\Psi_0(y)\sqrt{\sqrt{\epsilon}\,|u(T)|}\;\times$$

$$\exp\left\{-\frac{1}{\epsilon}\int_0^T dt\left[\left\{\sqrt{\epsilon}\,u(t)+\lambda\right\}_{\eta}-\eta g(t,\sqrt{\epsilon}\,u(T)+\lambda,\sqrt{\epsilon}\,u(0)+\lambda)\right]^2\right\}\times(\cdot)$$

and

$$\overline{\int\limits_{u(T)=\frac{x-\lambda}{\sqrt{\epsilon}},\,u(0)=\frac{y-\lambda}{\sqrt{\epsilon}}}} \widetilde{D}_{\mathbf{sm}}^{(2)}[u(t),\varepsilon,\epsilon]\Psi_0(y)\sqrt{\sqrt{\epsilon}\,|u(T)|}\;\times(\cdot)=$$

$$\overline{\int\limits_{q(T)=\frac{x-\lambda}{\sqrt{\epsilon}},\,u(0)=\frac{y-\lambda}{\sqrt{\epsilon}}}} \widetilde{D}_{\eta}^{+}[u(t),\varepsilon,\epsilon]\Psi_0(y)\sqrt{\sqrt{\epsilon}\,|u(T)|}\;\times$$

$$[\sin^{+}(S_{\varepsilon}(T,x,\lambda,u(t),\epsilon))+\sin^{-}(S_{\varepsilon}(T,x,\lambda,u(t),\epsilon))+1-1]\times(\cdot)=$$

$$\overline{\int\limits_{u(T)=\frac{x-\lambda}{\sqrt{\epsilon}},\,u(0)=\frac{y-\lambda}{\sqrt{\epsilon}}}} \widetilde{D}_{\mathbf{m}}^{(2)}[u(t),\varepsilon,\epsilon]\Psi_0(y)\sqrt{\sqrt{\epsilon}\,|u(T)|}\;\times(\cdot)-$$

$$\overline{\int\limits_{q(T)=\frac{x-\lambda}{\sqrt{\epsilon}},\,u(0)=\frac{y-\lambda}{\sqrt{\epsilon}}}} \widetilde{D}_{\eta}^{+}[u(t),\varepsilon,\epsilon]\Psi_0(y)\sqrt{\sqrt{\epsilon}\,|u(T)|}\;\times(\cdot)=$$

$$\qquad\qquad\qquad\qquad (2.2.2.94.b)$$

$$\overline{\int\limits_{u(T)=\frac{x-\lambda}{\sqrt{\epsilon}},\,u(0)=\frac{y-\lambda}{\sqrt{\epsilon}}}} \widetilde{D}_{\mathbf{m}}^{(2)}[u(t),\varepsilon,\epsilon]\Psi_0(y)\sqrt{\sqrt{\epsilon}\,|u(T)|}\;\times(\cdot)$$

$$-\overline{\int\limits_{u(T)=\frac{x-\lambda}{\sqrt{\epsilon}},\,u(0)=\frac{y-\lambda}{\sqrt{\epsilon}}}} D^{+}[q(t),\varepsilon]\Psi_0(y)\sqrt{\sqrt{\epsilon}\,|u(T)|}\;\times$$

$$\exp\left\{-\frac{1}{\epsilon}\int\limits_{0}^{T}dt\Big[\big\{\sqrt{\epsilon}\,u(t)+\lambda\big\}_{\eta}-\eta g(t,\sqrt{\epsilon}\,u(T)+\lambda,\sqrt{\epsilon}\,u(0)+\lambda)\Big]^2\right\}\times(\cdot).$$

Here $\widetilde{D}_{\mathbf{m}}^{(1)}$ and $\widetilde{D}_{\mathbf{m}}^{(2)}$ is the corresponding positive Feynman-Colombeau submeasures,i.e.

$$\overline{\int\limits_{u(T)=\frac{x-\lambda}{\sqrt{\epsilon}},\,u(0)=\frac{y-\lambda}{\sqrt{\epsilon}}} \widetilde{D}_{\mathbf{m}}^{(1)}[u(t),\varepsilon,\epsilon]\Psi_0(y)\,\sqrt{\sqrt{\epsilon}\,|u(T)|}} \times (\bullet) =$$

$$\overline{\int\limits_{q(T)=\frac{x-\lambda}{\sqrt{\epsilon}},\,u(0)=\frac{y-\lambda}{\sqrt{\epsilon}}} \widetilde{D}_{\eta}^{+}[u(t),\varepsilon,\epsilon]\Psi_0(y)\,\sqrt{\sqrt{\epsilon}\,|u(T)|}} \times$$

$$\left[\cos^{+}(S_\varepsilon(T,x,\lambda,u(t),\epsilon)) + \cos^{-}(S_\varepsilon(T,x,\lambda,u(t),\epsilon)) + 1\right] \times (\bullet) =$$

$$(2.2.2.95.a)$$

$$\overline{\int\limits_{u(T)=\frac{x-\lambda}{\sqrt{\epsilon}},\,u(0)=\frac{y-\lambda}{\sqrt{\epsilon}}} D^{+}[q(t),\varepsilon]\Psi_0(y)\,\sqrt{\sqrt{\epsilon}\,|u(T)|}} \times$$

$$\left[\cos^{+}(S_\varepsilon(T,x,\lambda,u(t),\epsilon)) + \cos^{-}(S_\varepsilon(T,x,\lambda,u(t),\epsilon)) + 1\right] \times$$

$$\exp\left\{-\frac{1}{\epsilon}\int\limits_0^T dt\Big[\left\{\sqrt{\epsilon}\,u(t)+\lambda\right\}_\eta - \eta g(t,\sqrt{\epsilon}\,u(T)+\lambda,\sqrt{\epsilon}\,u(0)+\lambda)\Big]^2\right\} \times (\bullet)$$

and

$$\overline{\int\limits_{u(T)=\frac{x-\lambda}{\sqrt{\epsilon}},u(0)=\frac{y-\lambda}{\sqrt{\epsilon}}} \widetilde{D}_{\mathbf{m}}^{(2)}[u(t),\varepsilon,\epsilon]\Psi_0(y)\sqrt{\sqrt{\epsilon}\,|u(T)|}}\times(\bullet)=$$

$$\overline{\int\limits_{q(T)=\frac{x-\lambda}{\sqrt{\epsilon}},u(0)=\frac{y-\lambda}{\sqrt{\epsilon}}} \widetilde{D}_{\eta}^{+}[u(t),\varepsilon,\epsilon]\Psi_0(y)\sqrt{\sqrt{\epsilon}\,|u(T)|}}\times$$

$$[\sin^+(S_\varepsilon(T,x,\lambda,u(t),\epsilon))+\sin^-(S_\varepsilon(T,x,\lambda,u(t),\epsilon))+1]\times(\bullet)=$$

$$\qquad\qquad\qquad\qquad\qquad\qquad\qquad\qquad\qquad\qquad (2.2.2.95.b)$$

$$\overline{\int\limits_{u(T)=\frac{x-\lambda}{\sqrt{\epsilon}},u(0)=\frac{y-\lambda}{\sqrt{\epsilon}}} D^{+}[q(t),\varepsilon]\Psi_0(y)\sqrt{\sqrt{\epsilon}\,|u(T)|}}\times$$

$$[\sin^+(S_\varepsilon(T,x,\lambda,u(t),\epsilon))+\sin^-(S_\varepsilon(T,x,\lambda,u(t),\epsilon))+1]\times$$

$$\exp\left\{-\frac{1}{\epsilon}\int\limits_0^T dt\Big[\left\{\sqrt{\epsilon}\,u(t)+\lambda\right\}_\eta - \eta g(t,\sqrt{\epsilon}\,u(T)+\lambda,\sqrt{\epsilon}\,u(0)+\lambda)\Big]^2\right\}\times(\bullet)$$

From Eq.(2.2.2.92.a,b)-Eq.(2.2.2.95.a,b) we obtain

$$v_{1,\varepsilon}(T,x,x_0,\lambda,\epsilon,\eta) =$$

$$\bar{\int} dy \left\{ \overline{\int_{u(T)=\frac{x-\lambda}{\sqrt{\epsilon}},\, u(0)=\frac{y-\lambda}{\sqrt{\epsilon}}} \widetilde{D}_{\mathbf{sm}}^{(1)}[u(t),\varepsilon]\Psi_0(y)\sqrt{\sqrt{\epsilon}\,|u(T)|}} \times \right.$$

$$\left| \cos\left[ \left( \int_0^T dt\left[ \frac{m}{2}\left( \frac{d\widetilde{u}_{\varepsilon\epsilon}(t)}{dt} \right)^2 - a_\varepsilon(t,\lambda)u_{\nu(\varepsilon,\epsilon)}^2(t) + b_\varepsilon(t,\lambda)\frac{u_{\nu(\varepsilon,\epsilon)}(t)}{\sqrt{\epsilon}} \right] \right) \right] \right| \times$$

$$\left. \cos\left[ \left( \sqrt{\epsilon}\int_0^T \widetilde{V}_1(u_{\nu(\varepsilon,\epsilon)}(t),t,\lambda)dt \right) \right] \right\} =$$

$$\bar{\int} dy \left\{ \overline{\int_{\substack{u(T)=\frac{x-\lambda}{\sqrt{\epsilon}}\\ u(0)=\frac{y-\lambda}{\sqrt{\epsilon}}}} \widetilde{D}_{\mathbf{m}}^{(1)}[u(t),\varepsilon]\Psi_0(y)\sqrt{\sqrt{\epsilon}\,|u(T)|}} \times \right.$$

$$\text{(2.2.2.96.a)}$$

$$\left| \cos\left[ \left( \int_0^T dt\left[ \frac{m}{2}\left( \frac{d\widetilde{u}_{\varepsilon\epsilon}(t)}{dt} \right)^2 - a_\varepsilon(t,\lambda)u_{\nu(\varepsilon,\epsilon)}^2(t) + b_\varepsilon(t,\lambda)\frac{u_{\nu(\varepsilon,\epsilon)}(t)}{\sqrt{\epsilon}} \right] \right) \right] \right| \times$$

$$\left. \cos\left[ \left( \sqrt{\epsilon}\int_0^T \widetilde{V}_1(u_{\nu(\varepsilon,\epsilon)}(t),t,\lambda)dt \right) \right] \right\} -$$

$$\bar{\int} dy \left\{ \overline{\int_{u(T)=\frac{x-\lambda}{\sqrt{\epsilon}},\, u(0)=\frac{y-\lambda}{\sqrt{\epsilon}}} \widetilde{D}_{\eta}^{+}[u(t),\varepsilon]\Psi_0(y)\sqrt{\sqrt{\epsilon}\,|u(T)|}} \times \right.$$

$$\left| \cos\left[ \left( \int_0^T dt\left[ \frac{m}{2}\left( \frac{d\widetilde{u}_{\varepsilon\epsilon}(t)}{dt} \right)^2 - a_\varepsilon(t,\lambda)u_{\nu(\varepsilon,\epsilon)}^2(t) + b_\varepsilon(t,\lambda)\frac{u_{\nu(\varepsilon,\epsilon)}(t)}{\sqrt{\epsilon}} \right] \right) \right] \right| \times$$

and

$$\nu_{2,\varepsilon}(T,x,x_0,\lambda,\epsilon,\eta) =$$

$$\overline{\int} dy \left\{ \overline{\int_{u(T)=\frac{x-\lambda}{\sqrt{\epsilon}},u(0)=\frac{y-\lambda}{\sqrt{\epsilon}}}} \widetilde{D}_{\mathbf{sm}}^{(2)}[u(t),\varepsilon]\Psi_0(y)\sqrt{\sqrt{\epsilon}\,|u(T)|} \times \right.$$

$$\left| \sin\left[ \left( \int_0^T dt\left[ \frac{m}{2}\left(\frac{d\widetilde{u}_{\varepsilon\epsilon}(t)}{dt}\right)^2 - a_\varepsilon(t,\lambda)u_{\nu(\varepsilon,\epsilon)}^2(t) + b_\varepsilon(t,\lambda)\frac{u_{\nu(\varepsilon,\epsilon)}(t)}{\sqrt{\epsilon}} \right] \right) \right] \right| \times$$

$$\sin\left[ \left( -\sqrt{\epsilon}\int_0^T \widetilde{V}_1(u_{\nu(\varepsilon,\epsilon)}(t),t,\lambda)dt \right) \right] \Bigg\} =$$

$$\overline{\int} dy \left\{ \overline{\int_{u(T)=\frac{x-\lambda}{\sqrt{\epsilon}},u(0)=\frac{y-\lambda}{\sqrt{\epsilon}}}} \widetilde{D}_{\mathbf{m}}^{(1)}[u(t),\varepsilon]\Psi_0(y)\sqrt{\sqrt{\epsilon}\,|u(T)|} \times \right.$$

$$\left| \sin\left[ \left( \int_0^T dt\left[ \frac{m}{2}\left(\frac{d\widetilde{u}_{\varepsilon\epsilon}(t)}{dt}\right)^2 - a_\varepsilon(t,\lambda)u_{\nu(\varepsilon,\epsilon)}^2(t) + b_\varepsilon(t,\lambda)\frac{u_{\nu(\varepsilon,\epsilon)}(t)}{\sqrt{\epsilon}} \right] \right) \right] \right| \times$$

$$(2.2.2.96.b)$$

$$\sin\left[ \left( -\sqrt{\epsilon}\int_0^T \widetilde{V}_1(u_{\nu(\varepsilon,\epsilon)}(t),t,\lambda)dt \right) \right] \Bigg\} -$$

$$\overline{\int} dy \left\{ \overline{\int_{u(T)=\frac{x-\lambda}{\sqrt{\epsilon}},u(0)=\frac{y-\lambda}{\sqrt{\epsilon}}}} \widetilde{D}_{\eta}^{+}[u(t),\varepsilon]\Psi_0(y)\sqrt{\sqrt{\epsilon}\,|u(T)|} \times \right.$$

$$\left| \sin\left[ \left( \int_0^T dt\left[ \frac{m}{2}\left(\frac{d\widetilde{u}_{\varepsilon\epsilon}(t)}{dt}\right)^2 - a_\varepsilon(t,\lambda)u_{\nu(\varepsilon,\epsilon)}^2(t) + b_\varepsilon(t,\lambda)\frac{u_{\nu(\varepsilon,\epsilon)}(t)}{\sqrt{\epsilon}} \right] \right) \right] \right| \times$$

$$\sin\left[ \left( -\sqrt{\epsilon}\int_0^T \widetilde{V}_1(u_{\nu(\varepsilon,\epsilon)}(t),t,\lambda)dt \right) \right] \Bigg\} .$$

From Eq.(2.2.2.96.a,b) we obtain

$$\nu_{1,\varepsilon}(T,x,x_0,\lambda,\epsilon,\eta) = \widetilde{\nu}_{1,\varepsilon}(T,x,x_0,\lambda,\epsilon,\eta) - \nu^*_{1,\varepsilon}(T,x,x_0,\lambda,\epsilon,\eta).$$

$$\widetilde{\nu}_{1,\varepsilon}(T,x,x_0,\lambda,\epsilon,\eta) = \overline{\int} dy \Bigg\{ \overline{\int_{u(T)=\frac{x-\lambda}{\sqrt{\epsilon}},\ u(0)=\frac{y-\lambda}{\sqrt{\epsilon}}} \widetilde{D}^{(1)}_{\mathbf{m}}[u(t),\varepsilon]\Psi_0(y)\sqrt{\sqrt{\epsilon}\,|u(T)|}} \times$$

$$\left| \cos\left[ \left( \int_0^T dt \left[ \frac{m}{2}\left(\frac{d\widetilde{u}_{\varepsilon\epsilon}(t)}{dt}\right)^2 - a_\varepsilon(t,\lambda)u^2_{\nu(\varepsilon,\epsilon)}(t) + b_\varepsilon(t,\lambda)\frac{u_{\nu(\varepsilon,\epsilon)}(t)}{\sqrt{\epsilon}} \right] \right) \right] \right| \times$$

$$\cos\left[ \left( \sqrt{\epsilon}\int_0^T \widetilde{V}_1(u_{\nu(\varepsilon,\epsilon)}(t),t,\lambda)dt \right) \right] \Bigg\},$$

$$(2.2.2.97.a)$$

$$\nu^*_{1,\varepsilon}(T,x,x_0,\lambda,\epsilon,\eta) = \overline{\int} dy \Bigg\{ \overline{\int_{u(T)=\frac{x-\lambda}{\sqrt{\epsilon}},\ u(0)=\frac{y-\lambda}{\sqrt{\epsilon}}} \widetilde{D}^+_{\eta}[u(t),\varepsilon]\Psi_0(y)\sqrt{\sqrt{\epsilon}\,|u(T)|}} \times$$

$$\left| \cos\left[ \left( \int_0^T dt \left[ \frac{m}{2}\left(\frac{d\widetilde{u}_{\varepsilon\epsilon}(t)}{dt}\right)^2 - a_\varepsilon(t,\lambda)u^2_{\nu(\varepsilon,\epsilon)}(t) + b_\varepsilon(t,\lambda)\frac{u_{\nu(\varepsilon,\epsilon)}(t)}{\sqrt{\epsilon}} \right] \right) \right] \right| \times$$

$$\cos\left[ \left( \sqrt{\epsilon}\int_0^T \widetilde{V}_1(u_{\nu(\varepsilon,\epsilon)}(t),t,\lambda)dt \right) \right] \Bigg\}$$

and

$$v_{2,\varepsilon}(T,x,x_0,\lambda,\epsilon,\eta) = \widetilde{v}_{2,\varepsilon}(T,x,x_0,\lambda,\epsilon,\eta) - v_{2,\varepsilon}^{*}(T,x,x_0,\lambda,\epsilon,\eta).$$

$$\widetilde{v}_{2,\varepsilon}(T,x,x_0,\lambda,\epsilon,\eta) = \overline{\int} dy \left\{ \overline{\int_{u(T)=\frac{x-\lambda}{\sqrt{\epsilon}},\ u(0)=\frac{y-\lambda}{\sqrt{\epsilon}}} \widetilde{D}_{\mathbf{m}}^{(2)}[u(t),\varepsilon]\Psi_0(y)\sqrt{\sqrt{\epsilon}\,|u(T)|}} \times \right.$$

$$\left| \sin\left[ \left( \int_0^T dt \left[ \frac{m}{2}\left(\frac{d\widetilde{u}_{\varepsilon\epsilon}(t)}{dt}\right)^2 - a_\varepsilon(t,\lambda)u_{\nu(\varepsilon,\epsilon)}^2(t) + b_\varepsilon(t,\lambda)\frac{u_{\nu(\varepsilon,\epsilon)}(t)}{\sqrt{\epsilon}} \right] \right) \right] \right| \times$$

$$\left. \sin\left[ \left( -\sqrt{\epsilon}\int_0^T \widetilde{V}_1(u_{\nu(\varepsilon,\epsilon)}(t),t,\lambda)dt \right) \right] \right\}, \tag{2.2.2.97.$b$}$$

$$v_{2,\varepsilon}^{*}(T,x,x_0,\lambda,\epsilon,\eta) = \overline{\int} dy \left\{ \overline{\int_{u(T)=\frac{x-\lambda}{\sqrt{\epsilon}},\ u(0)=\frac{y-\lambda}{\sqrt{\epsilon}}} \widetilde{D}_{\eta}^{+}[u(t),\varepsilon]\Psi_0(y)\sqrt{\sqrt{\epsilon}\,|u(T)|}} \times \right.$$

$$\left| \sin\left[ \left( \int_0^T dt \left[ \frac{m}{2}\left(\frac{d\widetilde{u}_{\varepsilon\epsilon}(t)}{dt}\right)^2 - a_\varepsilon(t,\lambda)u_{\nu(\varepsilon,\epsilon)}^2(t) + b_\varepsilon(t,\lambda)\frac{u_{\nu(\varepsilon,\epsilon)}(t)}{\sqrt{\epsilon}} \right] \right) \right] \right| \times$$

$$\left. \sin\left[ \left( -\sqrt{\epsilon}\int_0^T \widetilde{V}_1(u_{\nu(\varepsilon,\epsilon)}(t),t,\lambda)dt \right) \right] \right\}$$

Let us avaluate now path integral $v_{1,\varepsilon}^{*}(T,x,x_0,\lambda,\epsilon,\eta)$

$$v_{1,\varepsilon}^*(T, x, x_0, \lambda, \epsilon, \eta) =$$

$$\overline{\int} dy \left\{ \overline{\underset{u(T)=\frac{x-\lambda}{\sqrt{\epsilon}},\ u(0)=\frac{y-\lambda}{\sqrt{\epsilon}}}{\int}} \widetilde{D}_\eta^+[u(t), \varepsilon] \Psi_0(y) \sqrt{\sqrt{\epsilon}\,|u(T)|}\, \Xi_1(q(t), \lambda, \varepsilon, \epsilon) \times \right.$$

$$\cos\left[ \left( \sqrt{\epsilon} \int\limits_0^T \widetilde{V}_1(u_{v(\varepsilon,\epsilon)}(t), t, \lambda) dt \right) \right], \tag{2.2.2.98.a}$$

$$\Xi_1(q(t), \lambda, \varepsilon, \epsilon) =$$

$$\left| \cos\left[ \left( \int\limits_0^T dt \left[ \frac{m}{2} \left( \frac{d\widetilde{u}_{\varepsilon\epsilon}(t)}{dt} \right)^2 - a_\varepsilon(t, \lambda) u_{v(\varepsilon,\epsilon)}^2(t) + b_\varepsilon(t, \lambda) \frac{u_{v(\varepsilon,\epsilon)}(t)}{\sqrt{\epsilon}} \right] \right) \right] \right|$$

and path integral $v_{2,\varepsilon}^*(T, x, x_0, \lambda, \epsilon, \eta)$

$$v_{2,\varepsilon}^*(T, x, x_0, \lambda, \epsilon, \eta) =$$

$$\overline{\int} dy \left\{ \overline{\underset{u(T)=\frac{x-\lambda}{\sqrt{\epsilon}},\ u(0)=\frac{y-\lambda}{\sqrt{\epsilon}}}{\int}} \widetilde{D}_\eta^+[u(t), \varepsilon] \Psi_0(y) \sqrt{\sqrt{\epsilon}\,|u(T)|}\, \Xi_2(q(t), \lambda, \varepsilon, \epsilon) \times \right.$$

$$\sin\left[ \left( \sqrt{\epsilon} \int\limits_0^T \widetilde{V}_1(u_{v(\varepsilon,\epsilon)}(t), t, \lambda) dt \right) \right], \tag{2.2.2.98.b}$$

$$\Xi_2(q(t), \lambda, \varepsilon, \epsilon) =$$

$$\left| \sin\left[ \left( \int\limits_0^T dt \left[ \frac{m}{2} \left( \frac{d\widetilde{u}_{\varepsilon\epsilon}(t)}{dt} \right)^2 - a_\varepsilon(t, \lambda) u_{v(\varepsilon,\epsilon)}^2(t) + b_\varepsilon(t, \lambda) \frac{u_{v(\varepsilon,\epsilon)}(t)}{\sqrt{\epsilon}} \right] \right) \right] \right|$$

From Eq.(2.2.2.98.a,b) we obtain

$$v_{1,\varepsilon}^{*}(T, x, x_0, \lambda, \epsilon, \eta) =$$

$$\overline{\int\limits_{u(T)=\frac{x-\lambda}{\sqrt{\epsilon}},\; u(0)=\frac{y-\lambda}{\sqrt{\epsilon}}}} dy\, D^{+}[u(t), \varepsilon]\, \Psi_0(y)\, \sqrt{\sqrt{\epsilon}\, |u(T)|}\; \Xi_1(q(t), \lambda, \varepsilon, \epsilon) \times$$

$$\cos\!\left( \int\limits_0^T \sqrt{\epsilon}\, \widetilde{V}_1(q_{v(\varepsilon,\epsilon)}(t))\, dt \right) \times$$

$$\exp\!\left\{ -\frac{1}{\epsilon} \int\limits_0^T dt \Big[ \big\{ \sqrt{\epsilon}\, u(t) + \lambda \big\}_\eta - \eta g(t, \sqrt{\epsilon}\, u(T) + \lambda,\, \sqrt{\epsilon}\, u(0) + \lambda) \Big]^2 \right\} =$$

$$\overline{\int\limits_{u(T)=\frac{x-\lambda}{\sqrt{\epsilon}},\; u(0)=\frac{y-\lambda}{\sqrt{\epsilon}}}} dy\, D^{+}[u(t), \varepsilon]\, \Psi_0(y)\, \sqrt{\sqrt{\epsilon}\, |u(T)|}\; \Xi_1(q(t), \lambda, \varepsilon, \epsilon) \times$$

$$\exp\!\left\{ -\frac{1}{\epsilon} \int\limits_0^T dt \Big[ \big\{ \sqrt{\epsilon}\, u(t) + \lambda \big\}_\eta - \eta g(t, \sqrt{\epsilon}\, u(T) + \lambda,\, \sqrt{\epsilon}\, u(0) + \lambda) \Big]^2 \right\} \times$$

$$\left[ \sum_{n=0}^{M} \frac{(-1)^n \epsilon^n}{(2n)!} \left( \int\limits_0^T \widetilde{V}_1(q_{v(\varepsilon,\epsilon)}(t))\, dt \right)^{2n} \right] =$$

$$\overline{\int\limits_{u(T)=\frac{x-\lambda}{\sqrt{\epsilon}},\; u(0)=\frac{y-\lambda}{\sqrt{\epsilon}}}} dy\, D^{+}[u(t), \varepsilon]\, \Psi_0(y)\, \sqrt{\sqrt{\epsilon}\, |u(T)|}\; \Xi_1(q(t), \lambda, \varepsilon, \epsilon) \times$$

$$\exp\!\left\{ -\frac{1}{\epsilon} \int\limits_0^T dt \Big[ \big\{ \sqrt{\epsilon}\, u(t) + \lambda \big\}_\eta - \eta g(t, \sqrt{\epsilon}\, u(T) + \lambda,\, \sqrt{\epsilon}\, u(0) + \lambda) \Big]^2 \right\} +$$

$$\sum_{n=1}^{M} \frac{(-1)^n \epsilon^n}{(2n)!} \times$$

$$\left[ \overline{\int} dy\, D^{+}[u(t), \varepsilon]\, \Psi_0(y)\, \sqrt{\sqrt{\epsilon}\, |u(T)|}\; \Xi_1(q(t), \lambda, \varepsilon, \epsilon) \times \right.$$

<div style="text-align:right">(2.2.2.99.a)</div>

and

$$v_{2,\varepsilon}^*(T, x, x_0, \lambda, \epsilon, \eta) =$$

$$\overline{\int_{u(T)=\frac{x-\lambda}{\sqrt{\epsilon}},\, u(0)=\frac{y-\lambda}{\sqrt{\epsilon}}}} dy\, D^+[u(t), \varepsilon]\, \Psi_0(y)\, \sqrt{\sqrt{\epsilon}\,|u(T)|}\; \Xi_2(q(t), \lambda, \varepsilon, \epsilon) \times$$

$$\sin\left( \int_0^T \sqrt{\epsilon}\, \widetilde{V}_1(q_{v(\varepsilon,\epsilon)}(t))\, dt \right) \times$$

$$\exp\left\{ -\frac{1}{\epsilon} \int_0^T dt \Big[ \{\sqrt{\epsilon}\, u(t) + \lambda\}_\eta - \eta g(t, \sqrt{\epsilon}\, u(T) + \lambda, \sqrt{\epsilon}\, u(0) + \lambda) \Big]^2 \right\} =$$

$$\overline{\int_{u(T)=\frac{x-\lambda}{\sqrt{\epsilon}},\, u(0)=\frac{y-\lambda}{\sqrt{\epsilon}}}} dy\, D^+[u(t), \varepsilon]\, \Psi_0(y)\, \sqrt{\sqrt{\epsilon}\,|u(T)|}\; \Xi_2(q(t), \lambda, \varepsilon, \epsilon) \times$$

$$\exp\left\{ -\frac{1}{\epsilon} \int_0^T dt \Big[ \{\sqrt{\epsilon}\, u(t) + \lambda\}_\eta - \eta g(t, \sqrt{\epsilon}\, u(T) + \lambda, \sqrt{\epsilon}\, u(0) + \lambda) \Big]^2 \right\} \times$$

$$\left[ \sum_{n=0}^{M} \frac{(-1)^n \epsilon^{n+1/2}}{(2n+1)!} \left( \int_0^T \widetilde{V}_1(q_{v(\varepsilon,\epsilon)}(t))\, dt \right)^{2n+1} + \right.$$

$$\left. \epsilon^{n+1} O\left( \left( \int_0^T \widetilde{V}_1(q_{v(\varepsilon,\epsilon)}(t))\, dt \right)^{2n+2} \right) \right] = \qquad (2.2.2.99.b)$$

$$\overline{\int_{u(T)=\frac{x-\lambda}{\sqrt{\epsilon}},\, u(0)=\frac{y-\lambda}{\sqrt{\epsilon}}}} dy\, D^+[u(t), \varepsilon]\, \Psi_0(y)\, \sqrt{\sqrt{\epsilon}\,|u(T)|}\; \Xi_2(q(t), \lambda, \varepsilon, \epsilon) \times$$

$$\exp\left\{ -\frac{1}{\epsilon} \int_0^T dt \Big[ \{\sqrt{\epsilon}\, u(t) + \lambda\}_\eta - \eta g(t, \sqrt{\epsilon}\, u(T) + \lambda, \sqrt{\epsilon}\, u(0) + \lambda) \Big]^2 \right\} +$$

$$+ \sum_{n=0}^{M} \frac{(-1)^n \epsilon^{n+1/2}}{(2n+1)!} \times$$

Let us avaluate path integral $\tilde{\nu}_{1,\varepsilon}(T,x,x_0,\lambda,\epsilon,\eta)$ and $\tilde{\nu}_{2,\varepsilon}(T,x,x_0,\lambda,\epsilon,\eta)$. Assume now that $p_2 \gg \epsilon^{-1}$, i.e. $p_1 \approx 1$. Then from Eqs.(2.2.2.97.a,b) and Hölder's inequality (2.2.2.33) we obtain:

$$|\widetilde{v}_{1,\varepsilon}(T,x,x_0,\lambda,\epsilon,\eta)| \leq$$

$$\overline{\int} dy \left\{ \overline{\int_{u(T)=\frac{x-\lambda}{\sqrt{\epsilon}},\ u(0)=\frac{y-\lambda}{\sqrt{\epsilon}}} \widetilde{D}_{\mathbf{m}}^{(1)}[u(t),\varepsilon]\Psi_0(y)\ \sqrt{\sqrt{\epsilon}\,|u(T)|}} \times \right.$$

$$\left|\cos\left[\left(\int_0^T dt\left[\frac{m}{2}\left(\frac{d\widetilde{u}_{\varepsilon\epsilon}(t)}{dt}\right)^2 - a_\varepsilon(t,\lambda)u_{\nu(\varepsilon,\epsilon)}^2(t) + \right.\right.\right.\right.$$

$$\left.\left.\left.\left. +b_\varepsilon(t,\lambda)\frac{u_{\nu(\varepsilon,\epsilon)}(t)}{\sqrt{\epsilon}}\ \right]\right)\right]\right| \times$$

$$\left.\left|\cos\left[\left(\sqrt{\epsilon}\int_0^T \widetilde{V}_1(u_{\nu(\varepsilon,\epsilon)}(t),t,\lambda)dt\right)\right]\right|\right\} \leq$$

$$\left[\ \overline{\int} dy \left\{ \overline{\int_{u(T)=\frac{x-\lambda}{\sqrt{\epsilon}},\ u(0)=\frac{y-\lambda}{\sqrt{\epsilon}}} \widetilde{D}_{\mathbf{m}}^{(1)}[u(t),\varepsilon]\Psi_0(y)\ \sqrt{\sqrt{\epsilon}\,|u(T)|}} \times \right.\right.$$

$$\left|\cos\left[\left(\int_0^T dt\left[\frac{m}{2}\left(\frac{d\widetilde{u}_{\varepsilon\epsilon}(t)}{dt}\right)^2 - a_\varepsilon(t,\lambda)u_{\nu(\varepsilon,\epsilon)}^2(t) + \right.\right.\right.\right.$$

$$(2.2.2.100.a)$$

$$\left.\left.\left.\left. +b_\varepsilon(t,\lambda)\frac{u_{\nu(\varepsilon,\epsilon)}(t)}{\sqrt{\epsilon}}\ \right]\right)\right]\right|\right]^{1/p_1} \times$$

$$\left[\ \overline{\int} dy \left\{ \overline{\int_{u(T)=\frac{x-\lambda}{\sqrt{\epsilon}},\ u(0)=\frac{y-\lambda}{\sqrt{\epsilon}}} \widetilde{D}_{\mathbf{m}}^{(1)}[u(t),\varepsilon]\Psi_0(y)\ \sqrt{\sqrt{\epsilon}\,|u(T)|}} \times \right.\right.$$

$$\left|\cos\left[\left(\int_0^T dt\left[\frac{m}{2}\left(\frac{d\widetilde{u}_{\varepsilon\epsilon}(t)}{dt}\right)^2 - a_\varepsilon(t,\lambda)u_{\nu(\varepsilon,\epsilon)}^2(t) + \right.\right.\right.\right.$$

and

$$|\widetilde{v}_{2,\varepsilon}(T,x,x_0,\lambda,\epsilon,\eta)| \leq$$

$$\overline{\int} dy \left\{ \overline{\int\limits_{u(T)=\frac{x-\lambda}{\sqrt{\epsilon}},\ u(0)=\frac{y-\lambda}{\sqrt{\epsilon}}} \widetilde{D}_{\mathbf{m}}^{(2)}[u(t),\varepsilon]\Psi_0(y)\sqrt{\sqrt{\epsilon}\,|u(T)|}} \times \right.$$

$$\left| \sin\left[ \left( \int\limits_0^T dt \left[ \frac{m}{2}\left(\frac{d\widetilde{u}_{\varepsilon\epsilon}(t)}{dt}\right)^2 - a_\varepsilon(t,\lambda)u_{\nu(\varepsilon,\epsilon)}^2(t) + \right. \right. \right. \right.$$

$$\left. \left. \left. \left. + b_\varepsilon(t,\lambda)\frac{u_{\nu(\varepsilon,\epsilon)}(t)}{\sqrt{\epsilon}} \right] \right) \right] \right| \times$$

$$\left. \left| \sin\left[ \left( -\sqrt{\epsilon}\int\limits_0^T \widetilde{V}_1(u_{\nu(\varepsilon,\epsilon)}(t),t,\lambda)dt \right) \right] \right| \right\} \leq$$

$$\left[ \overline{\int} dy \left\{ \overline{\int\limits_{u(T)=\frac{x-\lambda}{\sqrt{\epsilon}},\ u(0)=\frac{y-\lambda}{\sqrt{\epsilon}}} \widetilde{D}_{\mathbf{m}}^{(2)}[u(t),\varepsilon]\Psi_0(y)\sqrt{\sqrt{\epsilon}\,|u(T)|}} \times \right. \right.$$

$$\hspace{6cm} (2.2.2.100.b)$$

$$\left. \left| \sin\left[ \left( \int\limits_0^T dt \left[ \frac{m}{2}\left(\frac{d\widetilde{u}_{\varepsilon\epsilon}(t)}{dt}\right)^2 - a_\varepsilon(t,\lambda)u_{\nu(\varepsilon,\epsilon)}^2(t) + \right. \right. \right. \right. \right.$$

$$\left. \left. \left. \left. \left. + b_\varepsilon(t,\lambda)\frac{u_{\nu(\varepsilon,\epsilon)}(t)}{\sqrt{\epsilon}} \right] \right) \right] \right| \right]^{1/p_1} \times$$

$$\left[ \overline{\int} dy \left\{ \overline{\int\limits_{u(T)=\frac{x-\lambda}{\sqrt{\epsilon}},\ u(0)=\frac{y-\lambda}{\sqrt{\epsilon}}} \widetilde{D}_{\mathbf{m}}^{(2)}[u(t),\varepsilon]\Psi_0(y)\sqrt{\sqrt{\epsilon}\,|u(T)|}} \times \right. \right.$$

$$\left| \sin\left[ \left( \int\limits_0^T dt \left[ \frac{m}{2}\left(\frac{d\widetilde{u}_{\varepsilon\epsilon}(t)}{dt}\right)^2 - a_\varepsilon(t,\lambda)u_{\nu(\varepsilon,\epsilon)}^2(t) + \right. \right. \right. \right.$$

Therefore from Eqs.(2.2.2.100.a,b) we obtain

$$|\widetilde{v}_{1,\varepsilon}(T,x,x_0,\lambda,\epsilon,\eta)| \leq$$

$$\left[\,\overline{\int} dy \left\{ \overline{\int_{u(T)=\frac{x-\lambda}{\sqrt{\epsilon}},\ u(0)=\frac{y-\lambda}{\sqrt{\epsilon}}}} \widetilde{D}_{\mathbf{m}}^{(1)}[u(t),\varepsilon]\Psi_0(y)\sqrt{\sqrt{\epsilon}\,|u(T)|}\ \times\right.\right.$$

$$\left|\cos\left[\left(\int_0^T dt\left[\frac{m}{2}\left(\frac{d\widetilde{u}_{\varepsilon\epsilon}(t)}{dt}\right)^2 - a_\varepsilon(t,\lambda)u_{\nu(\varepsilon,\epsilon)}^2(t)\ +\right.\right.\right.\right.$$

$$\left.\left.\left.\left. +b_\varepsilon(t,\lambda)\frac{u_{\nu(\varepsilon,\epsilon)}(t)}{\sqrt{\epsilon}}\right]\right)\right]\right|\right]^{1/p_1}\ \times \qquad (2.2.2.101.a)$$

$$\left[\,\overline{\int} dy \left\{ \overline{\int_{u(T)=\frac{x-\lambda}{\sqrt{\epsilon}},\ u(0)=\frac{y-\lambda}{\sqrt{\epsilon}}}} \widetilde{D}_{\mathbf{m}}^{(1)}[u(t),\varepsilon]\Psi_0(y)\sqrt{\sqrt{\epsilon}\,|u(T)|}\ \times\right.\right.$$

$$\left|\cos\left[\left(\int_0^T dt\left[\frac{m}{2}\left(\frac{d\widetilde{u}_{\varepsilon\epsilon}(t)}{dt}\right)^2 - a_\varepsilon(t,\lambda)u_{\nu(\varepsilon,\epsilon)}^2(t)\ +\right.\right.\right.\right.$$

$$\left.\left.\left.\left. +b_\varepsilon(t,\lambda)\frac{u_{\nu(\varepsilon,\epsilon)}(t)}{\sqrt{\epsilon}}\right]\right)\right]\right|\right]^{1/p_2}$$

and

$$|\widetilde{v}_{2,\varepsilon}(T,x,x_0,\lambda,\epsilon,\eta)| \leq$$

$$\left[\,\bar{\int} dy \left\{\overline{\int\limits_{u(T)=\frac{x-\lambda}{\sqrt{\epsilon}},\;u(0)=\frac{y-\lambda}{\sqrt{\epsilon}}} \widetilde{D}_{\mathbf{m}}^{(2)}[u(t),\varepsilon]\Psi_0(y)\,\sqrt{\sqrt{\epsilon}\,|u(T)|}\,\times}\right.\right.$$

$$\left|\sin\left[\left(\int\limits_0^T dt\left[\frac{m}{2}\left(\frac{d\widetilde{u}_{\varepsilon\epsilon}(t)}{dt}\right)^2 - a_\varepsilon(t,\lambda)u_{\nu(\varepsilon,\epsilon)}^2(t)+\right.\right.\right.\right.$$

$$\left.\left.\left.\left.+b_\varepsilon(t,\lambda)\frac{u_{\nu(\varepsilon,\epsilon)}(t)}{\sqrt{\epsilon}}\right]\right)\right]\right|\right]^{1/p_1}\times \qquad (2.2.2.101.b)$$

$$\left[\,\bar{\int} dy \left\{\overline{\int\limits_{u(T)=\frac{x-\lambda}{\sqrt{\epsilon}},\;u(0)=\frac{y-\lambda}{\sqrt{\epsilon}}} \widetilde{D}_{\mathbf{m}}^{(2)}[u(t),\varepsilon]\Psi_0(y)\,\sqrt{\sqrt{\epsilon}\,|u(T)|}\,\times}\right.\right.$$

$$\left|\sin\left[\left(\int\limits_0^T dt\left[\frac{m}{2}\left(\frac{d\widetilde{u}_{\varepsilon\epsilon}(t)}{dt}\right)^2 - a_\varepsilon(t,\lambda)u_{\nu(\varepsilon,\epsilon)}^2(t)+\right.\right.\right.\right.$$

$$\left.\left.\left.\left.+b_\varepsilon(t,\lambda)\frac{u_{\nu(\varepsilon,\epsilon)}(t)}{\sqrt{\epsilon}}\right]\right)\right]\right|\right]^{1/p_2}.$$

From Eqs.(2.2.2.101.a,b) we obtain

$$\widetilde{v}_{1,\varepsilon}(T,x,x_0,\lambda,\epsilon,\eta) \leq \mathbf{J}_{1,\varepsilon}(T,x,x_0,\lambda,\epsilon,\eta) \times \Gamma_{1,\varepsilon}(T,x,x_0,\lambda,\epsilon,\eta,p_2),$$

$$\qquad (2.2.2.102)$$

$$\widetilde{v}_{2,\varepsilon}(T,x,x_0,\lambda,\epsilon,\eta) \leq \mathbf{J}_{2,\varepsilon}(T,x,x_0,\lambda,\epsilon,\eta) \times \Gamma_{2,\varepsilon}(T,x,x_0,\lambda,\epsilon,\eta,p_2).$$

Here

$$\mathbf{J}_{1,\varepsilon}(T,x,x_0,\lambda,\epsilon,\eta) =$$

$$\left[\ \overline{\int}dy\left\{\overline{\int_{u(T)=\frac{x-\lambda}{\sqrt{\epsilon}},\ u(0)=\frac{y-\lambda}{\sqrt{\epsilon}}}}\widetilde{D}_{\mathbf{m}}^{(1)}[u(t),\varepsilon]\Psi_0(y)\sqrt{\sqrt{\epsilon}\,|u(T)|}\ \times\right.\right.$$

$$\left|\cos\left[\left(\int_0^T dt\left[\frac{m}{2}\left(\frac{d\widetilde{u}_{\varepsilon\epsilon}(t)}{dt}\right)^2 - a_\varepsilon(t,\lambda)u_{\nu(\varepsilon,\epsilon)}^2(t) +\right.\right.\right.\right.$$

$$\left.\left.\left.\left.+b_\varepsilon(t,\lambda)\frac{u_{\nu(\varepsilon,\epsilon)}(t)}{\sqrt{\epsilon}}\right]\right)\right]\right|\right]^{1/p_1},$$

$$(2.2.2.103.a)$$

$$\Gamma_{1,\varepsilon}(T,x,x_0,\lambda,\epsilon,\eta,p_2) =$$

$$\left[\ \overline{\int}dy\left\{\overline{\int_{u(T)=\frac{x-\lambda}{\sqrt{\epsilon}},\ u(0)=\frac{y-\lambda}{\sqrt{\epsilon}}}}\widetilde{D}_{\mathbf{m}}^{(1)}[u(t),\varepsilon]\Psi_0(y)\sqrt{\sqrt{\epsilon}\,|u(T)|}\ \times\right.\right.$$

$$\left|\cos\left[\left(\int_0^T dt\left[\frac{m}{2}\left(\frac{d\widetilde{u}_{\varepsilon\epsilon}(t)}{dt}\right)^2 - a_\varepsilon(t,\lambda)u_{\nu(\varepsilon,\epsilon)}^2(t) +\right.\right.\right.\right.$$

$$\left.\left.\left.\left.+b_\varepsilon(t,\lambda)\frac{u_{\nu(\varepsilon,\epsilon)}(t)}{\sqrt{\epsilon}}\right]\right)\right]\right|\right]^{1/p_2}$$

and

$$\mathbf{J}_{2,\varepsilon}(T,x,x_0,\lambda,\epsilon,\eta) =$$

$$\left[\;\overline{\int} dy \left\{ \overline{\int_{u(T)=\frac{x-\lambda}{\sqrt{\epsilon}},\; u(0)=\frac{y-\lambda}{\sqrt{\epsilon}}}} \widetilde{D}^{(2)}_{\mathbf{m}}[u(t),\varepsilon]\Psi_0(y)\; \sqrt{\sqrt{\epsilon}\,|u(T)|}\;\times \right. \right.$$

$$\left|\sin\left[\left(\int\limits_0^T dt\left[\frac{m}{2}\left(\frac{d\widetilde{u}_{\varepsilon\epsilon}(t)}{dt}\right)^2 - a_\varepsilon(t,\lambda)u^2_{\nu(\varepsilon,\epsilon)}(t) + \right.\right.\right.\right.$$

$$\left.\left.\left.\left.+ b_\varepsilon(t,\lambda)\frac{u_{\nu(\varepsilon,\epsilon)}(t)}{\sqrt{\epsilon}}\right]\right)\right]\right|^{1/p_1},$$

$$(2.2.2.103.b)$$

$$\Gamma_{2,\varepsilon}(T,x,x_0,\lambda,\epsilon,\eta,p_2) =$$

$$\left[\;\overline{\int} dy \left\{ \overline{\int_{u(T)=\frac{x-\lambda}{\sqrt{\epsilon}},\; u(0)=\frac{y-\lambda}{\sqrt{\epsilon}}}} \widetilde{D}^{(2)}_{\mathbf{m}}[u(t),\varepsilon]\Psi_0(y)\; \sqrt{\sqrt{\epsilon}\,|u(T)|}\;\times \right. \right.$$

$$\left|\sin\left[\left(\int\limits_0^T dt\left[\frac{m}{2}\left(\frac{d\widetilde{u}_{\varepsilon\epsilon}(t)}{dt}\right)^2 - a_\varepsilon(t,\lambda)u^2_{\nu(\varepsilon,\epsilon)}(t) + \right.\right.\right.\right.$$

$$\left.\left.\left.\left.+ b_\varepsilon(t,\lambda)\frac{u_{\nu(\varepsilon,\epsilon)}(t)}{\sqrt{\epsilon}}\right]\right)\right]\right|^{1/p_2}.$$

From Eq.(2.2.2.103.a,b) by using definition (see Eq.(2.2.2.95)) of the submeasures $\widetilde{D}^{(1)}_{\mathbf{m}}[q(t),\varepsilon,\epsilon], \widetilde{D}^{(2)}_{\mathbf{m}}[q(t),\varepsilon,\epsilon]$ we obtain

$$\mathbf{J}_{1,\varepsilon}(T,x,x_0,\lambda,\epsilon,\eta) \leq$$

$$\left[\,\overline{\int} dy \left\{ \overline{\int\limits_{u(T)=\frac{x-\lambda}{\sqrt{\epsilon}},\ u(0)=\frac{y-\lambda}{\sqrt{\epsilon}}}} \widetilde{D}_{\mathbf{sm}}^{(1)}[u(t),\varepsilon]\Psi_0(y)\,\sqrt{\sqrt{\epsilon}\,|u(T)|}\,\times \right.\right.$$

$$\left|\cos\left[\left(\int\limits_0^T dt\left[\,\frac{m}{2}\left(\frac{d\widetilde{u}_{\varepsilon\epsilon}(t)}{dt}\right)^2 - a_\varepsilon(t,\lambda)u_{\nu(\varepsilon,\epsilon)}^2(t) + \right.\right.\right.\right.$$

$$\left.\left.\left.\left. + b_\varepsilon(t,\lambda)\frac{u_{\nu(\varepsilon,\epsilon)}(t)}{\sqrt{\epsilon}}\,\right]\right)\right]\right| \right] +$$

$$(2.2.2.104.a)$$

$$\left[\,\overline{\int} dy \left\{ \overline{\int\limits_{u(T)=\frac{x-\lambda}{\sqrt{\epsilon}},\ u(0)=\frac{y-\lambda}{\sqrt{\epsilon}}}} \widetilde{D}_{\eta}^+[u(t),\varepsilon]\Psi_0(y)\,\sqrt{\sqrt{\epsilon}\,|u(T)|}\,\times \right.\right.$$

$$\left|\cos\left[\left(\int\limits_0^T dt\left[\,\frac{m}{2}\left(\frac{d\widetilde{u}_{\varepsilon\epsilon}(t)}{dt}\right)^2 - a_\varepsilon(t,\lambda)u_{\nu(\varepsilon,\epsilon)}^2(t) + \right.\right.\right.\right.$$

$$\left.\left.\left.\left. + b_\varepsilon(t,\lambda)\frac{u_{\nu(\varepsilon,\epsilon)}(t)}{\sqrt{\epsilon}}\,\right]\right)\right]\right| \right] =$$

$$\Theta_{1,\varepsilon}(T,x,x_0,\lambda,\epsilon,\eta) + \Lambda_{1,\varepsilon}(T,x,x_0,\lambda,\epsilon,\eta)$$

and

$$\mathbf{J}_{2,\varepsilon}(T,x,x_0,\lambda,\epsilon,\eta) \leq$$

$$\left[ \overline{\overline{\int}} dy \left\{ \overline{\int_{u(T)=\frac{x-\lambda}{\sqrt{\epsilon}},\, u(0)=\frac{y-\lambda}{\sqrt{\epsilon}}}} \widetilde{D}_{\mathbf{sm}}^{(2)}[u(t),\varepsilon]\Psi_0(y)\sqrt{\sqrt{\epsilon}\,|u(T)|} \times \right.\right.$$

$$\left. \left| \sin\left[ \left( \int_0^T dt \left[ \frac{m}{2}\left(\frac{d\widetilde{u}_{\varepsilon\epsilon}(t)}{dt}\right)^2 - a_\varepsilon(t,\lambda)u_{\nu(\varepsilon,\epsilon)}^2(t) + \right.\right.\right.\right.\right.$$

$$\left.\left.\left.\left.\left. + b_\varepsilon(t,\lambda)\frac{u_{\nu(\varepsilon,\epsilon)}(t)}{\sqrt{\epsilon}} \right] \right) \right] \right| \right] +$$

$$(2.2.2.104.b)$$

$$\left[ \overline{\overline{\int}} dy \left\{ \overline{\int_{u(T)=\frac{x-\lambda}{\sqrt{\epsilon}},\, u(0)=\frac{y-\lambda}{\sqrt{\epsilon}}}} \widetilde{D}_{\mathbf{\eta}}^{+}[u(t),\varepsilon]\Psi_0(y)\sqrt{\sqrt{\epsilon}\,|u(T)|} \times \right.\right.$$

$$\left. \left| \sin\left[ \left( \int_0^T dt \left[ \frac{m}{2}\left(\frac{d\widetilde{u}_{\varepsilon\epsilon}(t)}{dt}\right)^2 - a_\varepsilon(t,\lambda)u_{\nu(\varepsilon,\epsilon)}^2(t) + \right.\right.\right.\right.\right.$$

$$\left.\left.\left.\left.\left. + b_\varepsilon(t,\lambda)\frac{u_{\nu(\varepsilon,\epsilon)}(t)}{\sqrt{\epsilon}} \right] \right) \right] \right| \right] =$$

$$\Theta_{2,\varepsilon}(T,x,x_0,\lambda,\epsilon,\eta) + \Lambda_{2,\varepsilon}(T,x,x_0,\lambda,\epsilon,\eta).$$

Here

$$\Theta_{1,\varepsilon}(T,x,x_0,\lambda,\epsilon,\eta) =$$

$$\left[\,\overline{\int} dy \left\{ \overline{\int\limits_{u(T)=\frac{x-\lambda}{\sqrt{\epsilon}},\, u(0)=\frac{y-\lambda}{\sqrt{\epsilon}}}} \widetilde{D}_{\mathbf{sm}}^{(1)}[u(t),\varepsilon] \Psi_0(y) \sqrt{\sqrt{\epsilon}\,|u(T)|} \,\times \right.\right.$$

$$\left| \cos\left[ \left( \int\limits_0^T dt \left[ \frac{m}{2}\left(\frac{d\widetilde{u}_{\varepsilon\epsilon}(t)}{dt}\right)^2 - a_\varepsilon(t,\lambda) u_{\nu(\varepsilon,\epsilon)}^2(t) + \right.\right.\right.\right.$$

$$\left.\left.\left.\left. +\, b_\varepsilon(t,\lambda) \frac{u_{\nu(\varepsilon,\epsilon)}(t)}{\sqrt{\epsilon}} \,\right] \right)\right] \right|\, \right] = \qquad (2.2.2.105.a)$$

$$\left[\,\overline{\int} dy \left\{ \overline{\int\limits_{u(T)=\frac{x-\lambda}{\sqrt{\epsilon}},\, u(0)=\frac{y-\lambda}{\sqrt{\epsilon}}}} \widetilde{D}_{\eta}^{+}[u(t),\varepsilon] \Psi_0(y) \sqrt{\sqrt{\epsilon}\,|u(T)|} \,\times \right.\right.$$

$$\cos\left[ \left( \int\limits_0^T dt \left[ \frac{m}{2}\left(\frac{d\widetilde{u}_{\varepsilon\epsilon}(t)}{dt}\right)^2 - a_\varepsilon(t,\lambda) u_{\nu(\varepsilon,\epsilon)}^2(t) + \right.\right.\right.$$

$$\left.\left.\left.\left. +\, b_\varepsilon(t,\lambda) \frac{u_{\nu(\varepsilon,\epsilon)}(t)}{\sqrt{\epsilon}} \,\right] \right)\right]\right],$$

$$\Theta_{2,\varepsilon}(T,x,x_0,\lambda,\epsilon,\eta) =$$

$$\left[ \overline{\int} dy \left\{ \overline{\int\limits_{u(T)=\frac{x-\lambda}{\sqrt{\epsilon}},\ u(0)=\frac{y-\lambda}{\sqrt{\epsilon}}} \widetilde{D}_{\mathbf{sm}}^{(2)}[u(t),\varepsilon] \Psi_0(y) \sqrt{\sqrt{\epsilon}\,|u(T)|}} \times \right. \right.$$

$$\left| \sin\left[ \left( \int\limits_0^T dt \left[ \frac{m}{2}\left( \frac{d\widetilde{u}_{\varepsilon\epsilon}(t)}{dt} \right)^2 - a_\varepsilon(t,\lambda)u_{\nu(\varepsilon,\epsilon)}^2(t) + \right. \right. \right. \right.$$

$$\left. \left. \left. \left. + b_\varepsilon(t,\lambda)\frac{u_{\nu(\varepsilon,\epsilon)}(t)}{\sqrt{\epsilon}} \right] \right) \right] \right| \right] = \qquad (2.2.2.105.b)$$

$$\left[ \overline{\int} dy \left\{ \overline{\int\limits_{u(T)=\frac{x-\lambda}{\sqrt{\epsilon}},\ u(0)=\frac{y-\lambda}{\sqrt{\epsilon}}} \widetilde{D}_{\mathfrak{\eta}}^{+}[u(t),\varepsilon] \Psi_0(y) \sqrt{\sqrt{\epsilon}\,|u(T)|}} \times \right. \right.$$

$$\sin\left[ \left( \int\limits_0^T dt \left[ \frac{m}{2}\left( \frac{d\widetilde{u}_{\varepsilon\epsilon}(t)}{dt} \right)^2 - a_\varepsilon(t,\lambda)u_{\nu(\varepsilon,\epsilon)}^2(t) + \right. \right. \right.$$

$$\left. \left. \left. + b_\varepsilon(t,\lambda)\frac{u_{\nu(\varepsilon,\epsilon)}(t)}{\sqrt{\epsilon}} \right] \right) \right] \right]$$

and

$$\Lambda_{1,\varepsilon}(T,x,x_0,\lambda,\epsilon,\eta) =$$

$$\left[\,\overline{\int} dy \left\{ \overline{\int_{u(T)=\frac{x-\lambda}{\sqrt{\epsilon}},\ u(0)=\frac{y-\lambda}{\sqrt{\epsilon}}}} \widetilde{D}^+_\eta[u(t),\varepsilon]\Psi_0(y)\,\sqrt{\sqrt{\epsilon}\,|u(T)|}\,\times \right.\right.$$

$$\left|\cos\left[\left(\int_0^T dt\left[\frac{m}{2}\left(\frac{d\widetilde{u}_{\varepsilon\epsilon}(t)}{dt}\right)^2 - a_\varepsilon(t,\lambda)u^2_{\nu(\varepsilon,\epsilon)}(t) + \right.\right.\right.\right.$$

$$\left.\left.\left.\left. + b_\varepsilon(t,\lambda)\frac{u_{\nu(\varepsilon,\epsilon)}(t)}{\sqrt{\epsilon}}\,\right]\right)\right]\right],$$

$$(2.2.2.106.a)$$

$$\Lambda_{2,\varepsilon}(T,x,x_0,\lambda,\epsilon,\eta) =$$

$$\left[\,\overline{\int} dy \left\{ \overline{\int_{u(T)=\frac{x-\lambda}{\sqrt{\epsilon}},\ u(0)=\frac{y-\lambda}{\sqrt{\epsilon}}}} \widetilde{D}^+_\eta[u(t),\varepsilon]\Psi_0(y)\,\sqrt{\sqrt{\epsilon}\,|u(T)|}\,\times \right.\right.$$

$$\left|\sin\left[\left(\int_0^T dt\left[\frac{m}{2}\left(\frac{d\widetilde{u}_{\varepsilon\epsilon}(t)}{dt}\right)^2 - a_\varepsilon(t,\lambda)u^2_{\nu(\varepsilon,\epsilon)}(t) + \right.\right.\right.\right.$$

$$\left.\left.\left.\left. + b_\varepsilon(t,\lambda)\frac{u_{\nu(\varepsilon,\epsilon)}(t)}{\sqrt{\epsilon}}\,\right]\right)\right]\right],$$

$$(2.2.2.106.b)$$

From (2.2.2.102) and (2.2.2.104.a,b) we obtain the inequalities:

$$|\widetilde{v}_{1,\varepsilon}(T,x,x_0,\lambda,\epsilon,\eta)| \leq \mathbf{J}_{1,\varepsilon}(T,x,x_0,\lambda,\epsilon,\eta) \leq$$

$$\Theta_{1,\varepsilon}(T,x,x_0,\lambda,\epsilon,\eta) + \Lambda_{1,\varepsilon}(T,x,x_0,\lambda,\epsilon,\eta),$$

$$|\widetilde{v}_{2,\varepsilon}(T,x,x_0,\lambda,\epsilon,\eta)| \leq \mathbf{J}_{2,\varepsilon}(T,x,x_0,\lambda,\epsilon,\eta) \leq$$

$$\Theta_{2,\varepsilon}(T,x,x_0,\lambda,\epsilon,\eta) + \Lambda_{2,\varepsilon}(T,x,x_0,\lambda,\epsilon,\eta).$$

$$(2.2.2.107)$$

Therefore

$$v_{1,\varepsilon}(T,x,x_0,\lambda,\epsilon,\eta) = \widetilde{v}_{1,\varepsilon}(T,x,x_0,\lambda,\epsilon,\eta) - v_{1,\varepsilon}^*(T,x,x_0,\lambda,\epsilon,\eta) \leq$$

$$\mathbf{J}_{1,\varepsilon}(T,x,x_0,\lambda,\epsilon,\eta) - v_{1,\varepsilon}^*(T,x,x_0,\lambda,\epsilon,\eta) \leq \qquad (2.2.2.108.a)$$

$$\Theta_{1,\varepsilon}(T,x,x_0,\lambda,\epsilon,\eta) + \Lambda_{1,\varepsilon}(T,x,x_0,\lambda,\epsilon,\eta) - v_{1,\varepsilon}^*(T,x,x_0,\lambda,\epsilon,\eta),$$

and

$$v_{2,\varepsilon}(T,x,x_0,\lambda,\epsilon,\eta) = \widetilde{v}_{2,\varepsilon}(T,x,x_0,\lambda,\epsilon,\eta) - v_{2,\varepsilon}^*(T,x,x_0,\lambda,\epsilon,\eta) \leq$$

$$\mathbf{J}_{2,\varepsilon}(T,x,x_0,\lambda,\epsilon,\eta) - v_{2,\varepsilon}^*(T,x,x_0,\lambda,\epsilon,\eta) \leq \qquad (2.2.2.108.b)$$

$$\Theta_{2,\varepsilon}(T,x,x_0,\lambda,\epsilon,\eta) + \Lambda_{2,\varepsilon}(T,x,x_0,\lambda,\epsilon,\eta) - v_{2,\varepsilon}^*(T,x,x_0,\lambda,\epsilon,\eta).$$

Substitution Eqs.(2.2.2.99.a,b) and Eqs.(2.2.2.106.a,b) into Eqs.(2.2.2.108.a,b) gives

$$v_{1,\varepsilon}(T,x,x_0,\lambda,\epsilon,\eta) \leq$$

$$\Theta_{1,\varepsilon}(T,x,x_0,\lambda,\epsilon,\eta) + \Lambda_{1,\varepsilon}(T,x,x_0,\lambda,\epsilon,\eta) -$$

$$-\overline{\int\limits_{u(T)=\frac{x-\lambda}{\sqrt{\epsilon}},\,u(0)=\frac{y-\lambda}{\sqrt{\epsilon}}} dy D^+[u(t),\varepsilon]\Psi_0(y)\sqrt{\sqrt{\epsilon}\,|u(T)|}\,\Xi_1(q(t),\lambda,\varepsilon,\epsilon)\times}$$

$$\exp\left\{-\frac{1}{\epsilon}\int\limits_0^T dt\Big[\big\{\sqrt{\epsilon}\,u(t)+\lambda\big\}_\eta - \eta g(t,\sqrt{\epsilon}\,u(T)+\lambda,\sqrt{\epsilon}\,u(0)+\lambda)\Big]^2\right\} +$$

$$\sum_{n=1}^{M}\frac{(-1)^n\epsilon^n}{(2n)!}\times \tag{2.2.2.109.a}$$

$$\Bigg[\overline{\int\limits_{u(T)=\frac{x-\lambda}{\sqrt{\epsilon}},\,u(0)=\frac{y-\lambda}{\sqrt{\epsilon}}} dy D^+[u(t),\varepsilon]\Psi_0(y)\sqrt{\sqrt{\epsilon}\,|u(T)|}\,\Xi_1(q(t),\lambda,\varepsilon,\epsilon)\times}$$

$$\left(\int\limits_0^T \widetilde{V}_1(q_{v(\varepsilon,\epsilon)}(t))dt\right)^{2n}\times$$

$$\exp\left\{-\frac{1}{\epsilon}\int\limits_0^T dt\Big[\big\{\sqrt{\epsilon}\,u(t)+\lambda\big\}_\eta - \eta g(t,\sqrt{\epsilon}\,u(T)+\lambda,\sqrt{\epsilon}\,u(0)+\lambda)\Big]^2\right\}\Bigg]$$

and

$$\nu_{2,\varepsilon}(T, x, x_0, \lambda, \epsilon, \eta) \leq$$

$$\Theta_{2,\varepsilon}(T, x, x_0, \lambda, \epsilon, \eta) + \Lambda_{2,\varepsilon}(T, x, x_0, \lambda, \epsilon, \eta) -$$

$$-\overline{\int\limits_{u(T)=\frac{x-\lambda}{\sqrt{\epsilon}},\, u(0)=\frac{y-\lambda}{\sqrt{\epsilon}}} dy\, D^+[u(t), \varepsilon] \Psi_0(y)\, \sqrt{\sqrt{\epsilon}\, |u(T)|}\, \Xi_2(q(t), \lambda, \varepsilon, \epsilon) \times}$$

$$\exp\left\{ -\frac{1}{\epsilon} \int\limits_0^T dt \Big[ \left\{ \sqrt{\epsilon}\, u(t) + \lambda \right\}_\eta - \eta g(t, \sqrt{\epsilon}\, u(T) + \lambda, \sqrt{\epsilon}\, u(0) + \lambda) \Big]^2 \right\} \times$$

$$+ \sum_{n=1}^M \frac{(-1)^n \epsilon^{n+1/2}}{(2n+1)!} \times \qquad\qquad (2.2.2.109.b)$$

$$\left[ \overline{\int\limits_{u(T)=\frac{x-\lambda}{\sqrt{\epsilon}},\, u(0)=\frac{y-\lambda}{\sqrt{\epsilon}}} dy\, D^+[u(t), \varepsilon] \Psi_0(y)\, \sqrt{\sqrt{\epsilon}\, |u(T)|}\, \Xi_2(q(t), \lambda, \varepsilon, \epsilon) \times} \right.$$

$$\left( \int\limits_0^T \widetilde{V}_1(q_{\nu(\varepsilon,\epsilon)}(t)) dt \right)^{2n+1} \times$$

$$\left. \exp\left\{ -\frac{1}{\epsilon} \int\limits_0^T dt \Big[ \left\{ \sqrt{\epsilon}\, u(t) + \lambda \right\}_\eta - \eta g(t, \sqrt{\epsilon}\, u(T) + \lambda, \sqrt{\epsilon}\, u(0) + \lambda) \Big]^2 \right\} \right].$$

Therefore from Eqs.(2.2.2.109.a,b) we obtain

$$v_{1,\varepsilon}(T,x,x_0,\lambda,\epsilon,\eta) \leq$$

$$\Theta_{1,\varepsilon}(T,x,x_0,\lambda,\epsilon,\eta) -$$

$$-\sum_{n=1}^{M} \frac{(-1)^n \epsilon^n}{(2n)!} \times$$

$$\left[ \overline{\int\limits_{u(T)=\frac{x-\lambda}{\sqrt{\epsilon}},\, u(0)=\frac{y-\lambda}{\sqrt{\epsilon}}} dy\, D^+[u(t),\varepsilon]\Psi_0(y) \sqrt{\sqrt{\epsilon}\,|u(T)|}\, \Xi_1(q(t),\lambda,\varepsilon,\epsilon) \times } \right. \qquad (2.2.2.110.a)$$

$$\left( \int\limits_0^T \widetilde{\mathcal{V}}_1(q_{v(\varepsilon,\epsilon)}(t))dt \right)^{2n} \times$$

$$\left. \exp\left\{ -\frac{1}{\epsilon} \int\limits_0^T dt \Big[ \left\{ \sqrt{\epsilon}\,u(t)+\lambda \right\}_\eta - \eta g(t,\sqrt{\epsilon}\,u(T)+\lambda, \sqrt{\epsilon}\,u(0)+\lambda) \Big]^2 \right\} \right]$$

and

$$v_{2,\varepsilon}(T, x, x_0, \lambda, \epsilon, \eta) \leq$$

$$\Theta_{2,\varepsilon}(T, x, x_0, \lambda, \epsilon, \eta) -$$

$$-\sum_{n=1}^{M} \frac{(-1)^n \epsilon^{n+1/2}}{(2n+1)!} \times$$

$$\left[ \overline{\int_{u(T)=\frac{x-\lambda}{\sqrt{\epsilon}}, \, u(0)=\frac{y-\lambda}{\sqrt{\epsilon}}} dy \, D^+[u(t), \varepsilon] \Psi_0(y) \sqrt{\sqrt{\epsilon}\, |u(T)|} \, \Xi_2(q(t), \lambda, \varepsilon, \epsilon)} \times \right. \tag{2.2.2.110.b}$$

$$\left( \int_0^T \widetilde{V}_1(q_{v(\varepsilon,\epsilon)}(t)) dt \right)^{2n+1} \times$$

$$\left. \exp\left\{ -\frac{1}{\epsilon} \int_0^T dt \Big[ \left\{ \sqrt{\epsilon}\, u(t) + \lambda \right\}_\eta - \eta g(t, \sqrt{\epsilon}\, u(T) + \lambda, \sqrt{\epsilon}\, u(0) + \lambda) \Big]^2 \right\} \right].$$

Let us calculate now oscillatory path integral $\Theta_{1,\varepsilon}(T, x, x_0, \lambda, \epsilon, \eta)$ :

$$\Theta_{1,\varepsilon}(T, x, x_0, \lambda, \epsilon, \eta) =$$

$$\left[ \overline{\int} dy \left\{ \overline{\int_{u(T)=\frac{x-\lambda}{\sqrt{\epsilon}},\, u(0)=\frac{y-\lambda}{\sqrt{\epsilon}}}} \widetilde{D}_{\eta}^{+}[u(t), \varepsilon] \Psi_0(y) \sqrt{\sqrt{\epsilon}\, |u(T)|} \times \right. \right.$$

$$\cos\left[ \left( \int_0^T dt \left[ \frac{m}{2} \left( \frac{d\widetilde{u}_{\varepsilon\epsilon}(t)}{dt} \right)^2 - a_{\varepsilon}(t, \lambda) u^2_{\nu(\varepsilon,\epsilon)}(t) + \right. \right. \right.$$

$$\left. \left. \left. \left. \left. + b_{\varepsilon}(t, \lambda) \frac{u_{\nu(\varepsilon,\epsilon)}(t)}{\sqrt{\epsilon}} \right] \right) \right] \right] ,$$

(2.2.2.111.a)

and oscillatory path integral $\Theta_{2,\varepsilon}(T, x, x_0, \lambda, \epsilon, \eta)$.

$$\Theta_{2,\varepsilon}(T, x, x_0, \lambda, \epsilon, \eta) =$$

$$\left[ \overline{\int} dy \left\{ \overline{\int_{u(T)=\frac{x-\lambda}{\sqrt{\epsilon}},\, u(0)=\frac{y-\lambda}{\sqrt{\epsilon}}}} \widetilde{D}_{\eta}^{+}[u(t), \varepsilon] \Psi_0(y) \sqrt{\sqrt{\epsilon}\, |u(T)|} \times \right. \right.$$

$$\sin\left[ \left( \int_0^T dt \left[ \frac{m}{2} \left( \frac{d\widetilde{u}_{\varepsilon\epsilon}(t)}{dt} \right)^2 - a_{\varepsilon}(t, \lambda) u^2_{\nu(\varepsilon,\epsilon)}(t) + \right. \right. \right.$$

$$\left. \left. \left. \left. \left. + b_{\varepsilon}(t, \lambda) \frac{u_{\nu(\varepsilon,\epsilon)}(t)}{\sqrt{\epsilon}} \right] \right) \right] \right] .$$

(2.2.2.111.b)

From Eqs.(2.2.2.110. a,b) we obtain

$$\Theta_{1,\varepsilon}(T,x,x_0,\lambda,\epsilon,\eta) =$$

$$\frac{1}{2}\overline{\int_{u(T)=\frac{x-\lambda}{\sqrt{\epsilon}},\; u(0)=\frac{y-\lambda}{\sqrt{\epsilon}}} dy\,\widetilde{D}_{\eta}^{+}[u(t),\varepsilon]\,\Psi_0(y)\,\sqrt{\sqrt{\epsilon}\,|u(T)|}\;\times}$$

$$\left\{\left[\exp\left[i\widetilde{S}_0(T,x,x_0,\lambda,u(t),\varepsilon,\epsilon)\right]+\exp\left[-i\widetilde{S}_0(T,x,x_0,\lambda,u(t),\varepsilon,\epsilon)\right]\right]\right\} =$$

$$\frac{1}{2}\widetilde{\Theta}_{1,\varepsilon}(T,x,x_0,\lambda,\epsilon,\eta)+\frac{1}{2}\widetilde{\Theta}_{2,\varepsilon}(T,x,x_0,\lambda,\epsilon,\eta)$$

and

$$\Theta_{2,\varepsilon}(T,x,x_0,\lambda,\epsilon,\eta) =$$

$$\frac{1}{2i}\overline{\int_{u(T)=\frac{x-\lambda}{\sqrt{\epsilon}},\; u(0)=\frac{y-\lambda}{\sqrt{\epsilon}}} dy\,\widetilde{D}_{\eta}^{+}[u(t),\varepsilon]\,\Psi_0(y)\,\sqrt{\sqrt{\epsilon}\,|u(T)|}\;\times}$$

$$\left\{\left[\exp\left[i\widetilde{S}_0(T,x,x_0,\lambda,u(t),\varepsilon,\epsilon)\right]-\exp\left[-i\widetilde{S}_0(T,x,x_0,\lambda,u(t),\varepsilon,\epsilon)\right]\right]\right\} =$$

$$\frac{1}{2i}\widetilde{\Theta}_{1,\varepsilon}(T,x,x_0,\lambda,\epsilon,\eta)-\frac{1}{2i}\widetilde{\Theta}_{2,\varepsilon}(T,x,x_0,\lambda,\epsilon,\eta).$$

Here

$$\widetilde{S}_0(T,x,x_0,\lambda,q(t),\varepsilon,\epsilon) =$$

$$\int_0^T dt\left[\frac{m}{2}\left(\frac{d\widetilde{u}_{\varepsilon\epsilon}(t)}{dt}\right)^2 - a_\varepsilon(t,\lambda)u_{\nu(\varepsilon,\epsilon)}^2(t)+b_\varepsilon(t,\lambda)\frac{u_{\nu(\varepsilon,\epsilon)}(t)}{\sqrt{\epsilon}}\right]$$

and

$$\Theta_{1,\varepsilon}(T,x,x_0,\lambda,\epsilon,\eta) =$$

$$\overline{\int\limits_{u(T)=\frac{x-\lambda}{\sqrt{\epsilon}},\ u(0)=\frac{y-\lambda}{\sqrt{\epsilon}}}} dy\, \widetilde{D}_\eta^+[u(t),\varepsilon]\Psi_0(y)\sqrt{\sqrt{\epsilon}\,|u(T)|}\ \exp\Big[\,i\widetilde{S}_0(T,x,x_0,\lambda,u(t),\varepsilon,\epsilon)\,\Big], \qquad (2.2.2.114.a)$$

$$\Theta_{2,\varepsilon}(T,x,x_0,\lambda,\epsilon,\eta) =$$

$$\overline{\int\limits_{u(T)=\frac{x-\lambda}{\sqrt{\epsilon}},\ u(0)=\frac{y-\lambda}{\sqrt{\epsilon}}}} dy\, \widetilde{D}_\eta^+[u(t),\varepsilon]\Psi_0(y)\sqrt{\sqrt{\epsilon}\,|u(T)|}\ \exp\Big[\,-i\widetilde{S}_0(T,x,x_0,\lambda,u(t),\varepsilon,\epsilon)\,\Big]. \qquad (2.2.2.114.b)$$

Let's calculate now path integral $\widetilde{\Theta}_{1,\varepsilon}(T,x,x_0,\lambda,\epsilon,\eta)$ by using stationary phase method.

By using replacement $q(t):=\dfrac{q(t)}{\sqrt{\epsilon}}$ we rewrite (2.2.2.114.a) in the next equivalent form

$$\Theta_{1,\varepsilon}(T,x,x_0,\lambda,\epsilon,\eta) =$$

$$\overline{\int\limits_{u(T)=x-\lambda,\ u(0)=y-\lambda}} dy\, \widetilde{D}_\eta^+[u(t),\varepsilon]\Psi_0(y)\sqrt{|u(T)|}\ \times$$

$$\exp\left\{\frac{i}{\epsilon}\int\limits_0^T dt\left[\frac{m}{2}\left(\frac{d\widetilde{u}_\varepsilon(t)}{dt}\right)^2 - a_\varepsilon(t,\lambda)u_\varepsilon^2(t) + b_\varepsilon(t,\lambda)u_\varepsilon(t)\right]\right\}. \qquad (2.2.2.115)$$

We assume now that $a_\varepsilon(t, \lambda) = \dfrac{m\omega^2}{2} + \widetilde{a}_\varepsilon(\lambda)$. Then the corresponding master Lagrangian $\widetilde{\mathcal{L}}_{\varepsilon \approx 0}(\dot{u}, u) = \widetilde{\mathcal{L}}(\dot{u}, u) + o(\varepsilon)$ is

$$\widetilde{\mathcal{L}}(\dot{u}, u) = \frac{m}{2}\dot{u}^2 - \left(\frac{m\omega^2}{2} + \widetilde{a}_\varepsilon(\lambda)\right)u^2 + b(t, \lambda)u. \tag{2.2.2.116}$$

The Euler-Lagrange equation of $\widetilde{\mathcal{L}}(\dot{u}, u)$ with corresponding boundary conditions is

$$\frac{d}{dt}\left(\frac{\partial \widetilde{\mathcal{L}}}{\partial \dot{u}}\right) - \frac{\partial \widetilde{\mathcal{L}}}{\partial u} = 0,$$

$$u(0) = -\lambda, u(T) = x - \lambda = x. \tag{2.2.2.117}$$

We assume now that $\dfrac{m\varpi^2}{2} = \dfrac{m\omega^2}{2} + \widetilde{a}_\varepsilon(\lambda) = \dfrac{m}{2}\left(\omega^2 + \dfrac{2\widetilde{a}_\varepsilon(\lambda)}{m}\right) \geq 0$ and rewrite Eq.(2.2.2.116) in the canonical form:

$$\widetilde{\mathcal{L}}(\dot{u}, u) = \frac{m}{2}\dot{u}^2 - \frac{m\varpi^2}{2}u^2 + b(t, \lambda)u,$$

$$\varpi(\lambda) = \left|\sqrt{\left(\omega^2 + \frac{2\widetilde{a}_\varepsilon(\lambda)}{m}\right)}\right|. \tag{2.2.2.118}$$

Thus Eq.(2.2.2.117) takes the form

$$m\ddot{u} + m\varpi^2(\lambda)u - b(t,\lambda) = 0,$$

$$u(0) = y - \lambda = \widetilde{y}, \qquad (2.2.2.119)$$

$$u(T) = x - \lambda = \widetilde{x}.$$

Corresponding general solution is

$$u_{\text{el}}(t) = \widetilde{u}(t) + \int_0^t \widetilde{u}(t-s)b(s,\lambda)ds,$$

$$(2.2.2.120)$$

$$\widetilde{u}(t) = c_1 \sin(\varpi t) + c_2 \cos(\varpi t).$$

Substituting Eq.(2.2.2.120) into boundary conditions (2.2.2.119) we obtain

$$\widetilde{u}(0) = c_2 = \widetilde{y},$$

$$\widetilde{u}(T) = \widetilde{x} = c_1 \sin \varpi T + \widetilde{y} \cos \varpi T +$$

$$(2.2.2.121)$$

$$\int_0^T [c_1 \sin(\varpi(t-s)) + \widetilde{y} \cos(\varpi(t-s))]b(s,\lambda)ds.$$

and

$$c_1 = \frac{\tilde{x} - \tilde{y}\cos\varpi T - \tilde{y}\int_0^T \cos(\varpi(t-s))b(s,\lambda)ds}{\sin\varpi T + \int_0^T \sin(\varpi(t-s))b(s,\lambda)ds}. \qquad (2.2.2.122)$$

Corresponding classical action $\tilde{S}_{\text{el}}$ is [13]:

$$\tilde{S}_{\text{el}}(T) =$$

$$\frac{m\varpi}{2\sin\varpi T}\Bigg[(\cos\varpi T)(\lambda^2 + \bar{x}^2) - 2\tilde{x}\tilde{y} + \frac{2\tilde{x}}{m\varpi}\int_0^T g(\lambda,t)\sin(\varpi t)dt +$$

$$+ \frac{2\tilde{y}}{m\varpi}\int_0^T g(\lambda,t)\sin\varpi(T-t)dt - \qquad (2.2.2.123)$$

$$- \frac{2}{m^2\varpi^2}\int_0^T\int_0^t g(\lambda,t)g(\lambda,s)\sin\varpi(T-t)\sin(\varpi s)dsdt\Bigg].$$

Substitution Eq.(2.2.123) into Eq.(2.2.115) gives

$$\Theta_{1,\varepsilon}(T,x,x_0,\lambda,\epsilon,\eta) =$$

$$\int_{-\infty}^{+\infty} dy D\Psi_0(y)\sqrt{|x-\lambda|}\exp\left\{\frac{i}{\epsilon}\int_0^T dt\tilde{\mathcal{L}}(\dot{u}(t),u(t))\right\} = \qquad (2.2.2.124)$$

the equation $d\widetilde{S}_{\mathrm{cl}}/d\widetilde{x} = 0$ takes the form

$$2\widetilde{x}(\cos \varpi T) - 2\widetilde{y} + \frac{2}{m\varpi} \int_0^T b(\lambda, t) \sin(\varpi t)dt = 0. \qquad (2.2.2.12)$$

From Eq.(2.3.2.19) we obtain

$$x = \lambda + \left( \frac{\widetilde{y}}{\cos \varpi T} - \frac{1}{m\varpi \cos \varpi T} \int_0^T g(\lambda, t) \sin(\varpi t)dt \right). \qquad (2.2.2.12)$$

$$(2.2.2.1)$$

$$(2.2.2.1)$$

From Theorem 2.2.2.1 and Eqs.(2.2.2.66) we obtain:
**Theorem 2.2.2.2**.

$$|\langle T, x_0, \lambda; \epsilon, \delta \rangle| \leq \widehat{U}(T, x_0, \lambda)(C + O(\delta)).$$

$$(2.2.2.67)$$

$$\liminf_{0 < h \ll \delta} |\langle T, x_0, \lambda; \epsilon, \delta \rangle| \leq \widehat{U}(T, x_0, \lambda) \times C.$$

Here

$$U(T, x_0, \lambda) = \left| \frac{1}{m\varpi(\lambda)} \int_0^T \widehat{b}(t, \lambda) \sin[\varpi(\lambda)(T - t)]dt + x_0 \right|,$$

$$(2.2.2.68)$$

$$\varpi^2(\lambda) = \frac{m\widehat{a}_\varepsilon(\lambda)}{2}.$$

From enequality we obtain(2.2.2.67) master Equation for the Feynman-Colombeau path integral.

**Theorem 2.2.2.2. (Master Equation for the Feynman-Colombeau path**

**integral**).

$$U(T, x_0, \lambda) = \left| \frac{1}{m\varpi(\lambda)} \int_0^T \widehat{b}(t, \lambda) \sin[\varpi(\lambda)(T - t)]dt + x_0 \right| = 0. \qquad (2.2.2.69)$$

# II.3. Quantum anharmonic oscillator with a cubic potential. Quasiclassical asimptotic of the quantum trajectories.

## II.3.1. Quantum harmonic oscillator. Quasiclassical asimptotic of the quantum trajectories.

As a first example we consider quantum harmonic oscillator (as cartooned in fig.2.3.1),

$$V(x) = m\frac{\omega^2}{2}x^2 - bx .$$  (2.3.1.1)

supplemented by an additive sinusoidal driving.

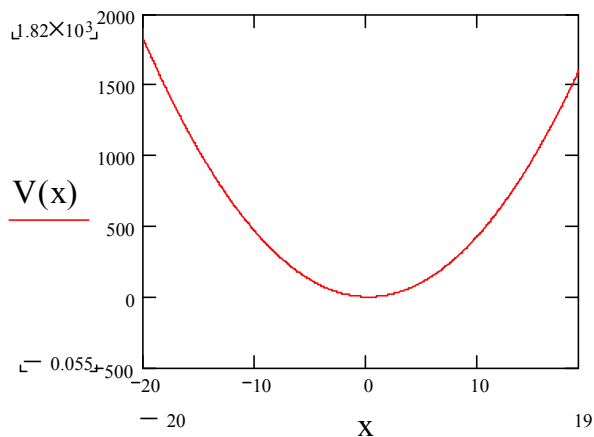

Fig. 2.3.1. $b = 1, \omega = 3, m = 1.$

**Example.1.1.**

$m = 1, \Omega = 10, \omega = 3, b = 1, A = 0.$

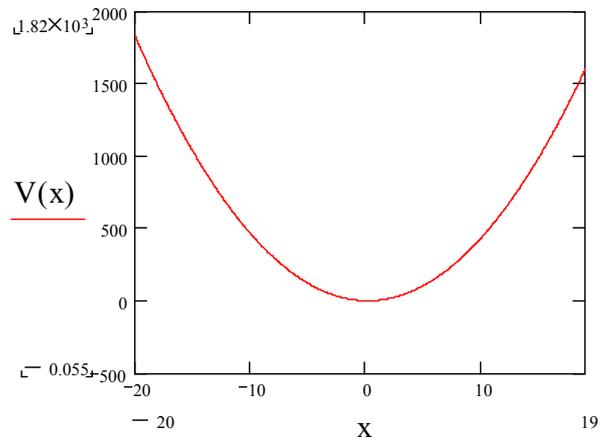

Fig. 1.1.1. $b = 1, \omega = 3, m = 1$.

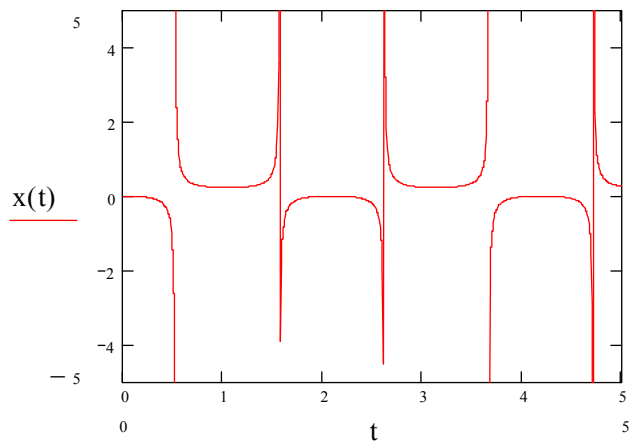

Fig. 1.1.2.

**Example**.**1.2**.

$m = 1, \Omega = 10, \omega = 3, b = 1, A = 10$.

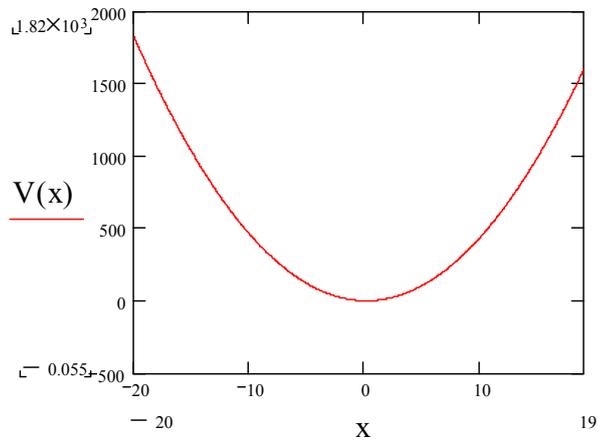

Fig. 1.2.1. $b = 1, \omega = 3, m = 1.$

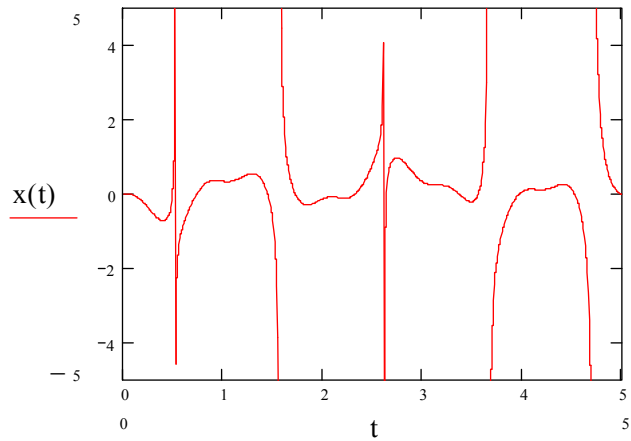

Fig. 1.2.2.

**Example**.**1**.**3**.

$m = 1, \Omega = 20, \omega = 5, b = 1, A = 10.$

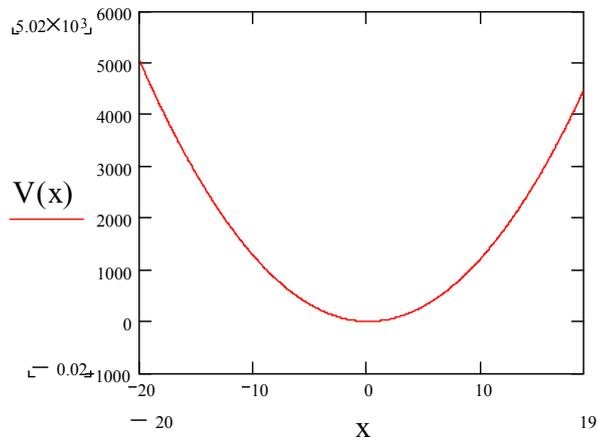

Fig. 1.3.1.

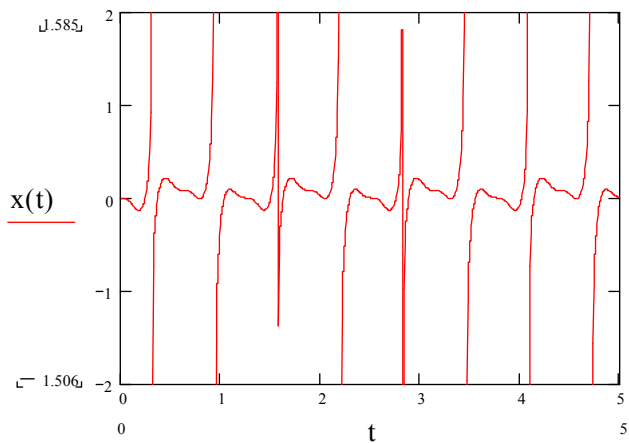

Fig. 1.3.2.

## II.3.2. Quantum anharmonic oscillator with a cubic potential. Quasiclassical asimptotic of the generalized quantum trajectories.

Let's calculate generalized quantum tradjectories (generalized quantum averages) of the one dimensional quantum anharmonic oscillator with a cubic potential, i.e.

$$\langle x_i, t, x_0; \hbar \rangle = \frac{\int\limits_{-L}^{L} x \Psi(x,t,x_0;\hbar)\,dx}{\int\limits_{-L}^{L} \Psi(x,t,x_0;\hbar)\,dx},$$

(2.3.2.0)

$$\Psi(x,0,x_0;\hbar) = \frac{1}{\sqrt{2\pi\varepsilon}} \exp\left(-\frac{x^2}{\varepsilon}\right) \simeq \delta(x-x_0),$$

$$\varepsilon \ll \hbar \ll 1.$$

As a second example we consider quantum anharmonic oscillator with a cubic potential (as cartooned in fig.2.3.2),

$$V(x) = \frac{a}{4}x^3 - \frac{b}{2}x\,.$$

(2.3.2.1)

supplemented by an additive sinusoidal driving.

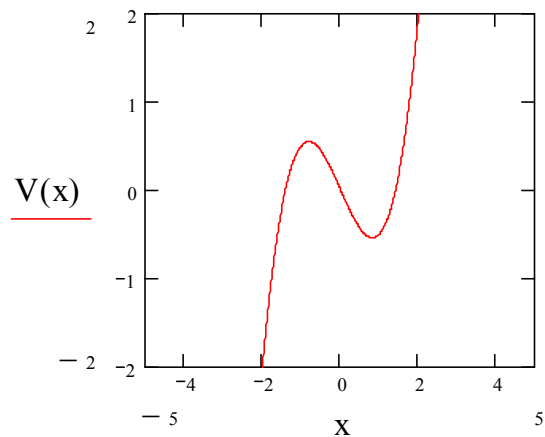

Fig.$2.3.2.\,a = 0.5, b = 1.$
$$V(x) = \frac{a}{4}x^3 - \frac{b}{2}x\,, V(x,t) = \frac{m\omega^2}{2}x^2 + ax^3 - bx - [A\sin(\Omega t)]x.$$

Let us consider transition probability amplitude given via formulae (2.1.1)-(2.1.2) with

$$V(x,t) = \frac{m\omega^2}{2}x^2 + ax^3 - bx - [A\sin(\Omega t)]x. \qquad (2.3.2.2)$$

The corresponding classical Lagrangian is

$$\mathcal{L}(\dot{x},x) = \frac{m}{2}\dot{x}^2 - \frac{m\omega^2}{2}x^2 - ax^3 + bx + [A\sin(\Omega t)]x. \qquad (2.3.2.3)$$

By substitution $x(\tau) = u(\tau) + \lambda$ we obtain

$$\mathcal{L}(\dot{u},u) = \frac{m}{2}\dot{u}^2 - \frac{m\omega^2}{2}(u+\lambda)^2 - a(u+\lambda)^3 +$$

$$b(u+\lambda) + [A\sin(\Omega t)](u+\lambda) =$$

$$\frac{m}{2}\dot{u}^2 - \frac{m\omega^2}{2}u^2 - m\omega^2\lambda u - \frac{m\omega^2\lambda^2}{2} -$$

$$-au^3 - 3au^2\lambda - 3au\lambda^2 - a\lambda^3 + bu + b\lambda +$$

$$[A\sin(\Omega t)]u + [A\sin(\Omega t)]\lambda = \qquad (2.3.2.4)$$

$$\frac{m}{2}\dot{u}^2 - \frac{m\omega^2}{2}u^2 - 3a\lambda u^2 -$$

$$-m\omega^2\lambda u - 3au\lambda^2 + bu + [A\sin(\Omega t)]u - au^3 -$$

$$-\frac{m\omega^2\lambda^2}{2} - a\lambda^3 + b\lambda + [A\sin(\Omega t)]\lambda.$$

The corresponding master Lagrangian $\widetilde{\mathcal{L}}(\dot{u}, u)$ is

$$\widetilde{\mathcal{L}}(\dot{u}, u) = \frac{m}{2}\dot{u}^2 - \left(\frac{m\omega^2}{2} + 3a\lambda\right)u^2 - (m\omega^2\lambda + 3a\lambda^2 - b - A\sin(\Omega t))u =$$

$$\frac{m}{2}\dot{u}^2 - m\left(\frac{\omega^2}{2} + \frac{3a\lambda}{m}\right)u^2 - (m\omega^2\lambda + 3a\lambda^2 - b - A\sin(\Omega t))u.$$

(2.3.2.5)

The Euler-Lagrange equation of $\widetilde{\mathcal{L}}(\dot{u}, u)$ with corresponding boundary conditions is

$$\frac{d}{dt}\left(\frac{\partial\widetilde{\mathcal{L}}}{\partial\dot{u}}\right) - \frac{\partial\widetilde{\mathcal{L}}}{\partial u} = 0,$$

(2.3.2.6)

$$u(0) = -\lambda, u(T) = x - \lambda = \bar{x}.$$

**Remark 2.3.1**. We assume now that $\frac{\omega^2}{2} + \frac{3a\lambda}{m} \geq 0$ and rewrite Eq.(2.3.2.5) in the form

$$\widetilde{\mathcal{L}}(\dot{u}, u) = \frac{m}{2}\dot{u}^2 - \frac{m\varpi^2(\lambda)}{2}u^2 + g(\lambda, t)u,$$

$$\varpi(\lambda) = \left|\sqrt{2\left(\frac{\omega^2}{2} + \frac{3a\lambda}{m}\right)}\right|,$$

(2.3.2.7)

$$g(\lambda, t) = -(m\omega^2\lambda + 3a\lambda^2 - b - A\sin(\Omega t)).$$

Thus Eq.(2.3.2.6) takes the form

$$m\ddot{u} + 2m\varpi^2(\lambda)u + d(\lambda) - A\sin(\Omega t) = 0,$$

$$d(\lambda) = m\omega^2\lambda + 3a\lambda^2 - b.$$

(2.3.2.8)

Plugging in Eq.(2.3.2.8) this trial solution

$$u(t) = C + D\sin(\Omega t)$$

(2.3.2.9)

we obtain

$$-mD\Omega^2\sin(\Omega t) + 2m\varpi^2(\lambda)(C + D\sin(\Omega t)) + d(\lambda) - A\sin(\Omega t) = 0.$$

(2.3.2.10)

Hence

$$2m\varpi^2(\lambda)C + d(\lambda) = 0,$$

$$-mD\Omega^2 + 2m\varpi^2(\lambda)D - A = 0.$$

(2.3.2.11)

Solving Eqs.(2.3.2.11) we obtain

$$C = -\frac{d(\lambda)}{2m\varpi^2(\lambda)},$$

$$D = \frac{A}{(2\varpi^2(\lambda) - \Omega^2)m}.$$

(2.3.2.12)

Thus general solution of the Eq.(2.3.2.8) is

$$u(t, \lambda) = c_1 \sin \varpi t + c_2 \cos \varpi t + D \sin(\Omega t) + C. \qquad (2.3.2.13)$$

Substituting Eq.(2.3.2.13) into boundary conditions

$$u(0) = -\lambda,$$

$$u(T) = x - \lambda = \bar{x} \qquad (2.3.2.14)$$

we obtain

$$u(0) = c_2 + C = -\lambda,$$

$$u(T) = c_1 \sin \varpi T + c_2 \cos \varpi T + D \sin(\Omega T) + C = \bar{x}. \qquad (2.3.2.15)$$

Hence

$$c_2 = -(C + \lambda),$$

$$c_1 = \frac{\bar{x} - C - D \sin(\Omega T) + (C + \lambda) \cos \varpi T}{\sin \varpi T} \qquad (2.3.2.16)$$

and

$$u(t,\lambda,\bar{x}) = \frac{\bar{x} - C - D\sin(\Omega T) + (C+\lambda)\cos\varpi T}{\sin\varpi T}\sin\varpi t - (C+\lambda)\cos\varpi t +$$

$$+D\sin(\Omega t) + C, \tag{2.3.2.17}$$

$$C = -\frac{d(\lambda)}{2m\varpi^2(\lambda)}, D = \frac{A}{(2\varpi^2(\lambda) - \Omega^2)m},$$

Note that

$$u(T,\lambda,\bar{x}) = \frac{\bar{x} - C - D\sin(\Omega T) + (C+\lambda)\cos\varpi T}{\sin\varpi T}\sin\varpi T - (C+\lambda)\cos\varpi T +$$

$$+D\sin(\Omega T) + C = \tag{2.3.17'}$$

$$\bar{x} - C - D\sin(\Omega T) + (C+\lambda)\cos\varpi T - (C+\lambda)\cos\varpi T + D\sin(\Omega T) + C = \bar{x},$$

$$u(T,\lambda,\bar{x}) = \bar{x}.$$

Corresponding master action $\widetilde{S}_{\text{cl}}$ is [13]:

$$\widetilde{S}_{\text{cl}} = \frac{m\varpi}{2\sin\varpi T}\left[(\cos\varpi T)(\lambda^2 + \bar{x}^2) + 2\lambda\bar{x} + \frac{2\bar{x}}{m\varpi}\int_0^T g(\lambda,t)\sin(\varpi t)dt -\right.$$

$$-\frac{2\lambda}{m\varpi}\int_0^T g(\lambda,t)\sin\varpi(T-t)dt - \tag{2.3.2.18}$$

$$\left.-\frac{2}{m^2\varpi^2}\int_0^T\int_0^t g(\lambda,t)g(\lambda,s)\sin\varpi(T-t)\sin(\varpi s)dsdt\right].$$

Thus the equation $d\widetilde{S}_{\text{cl}}/d\bar{x} = 0$ takes the form

$$2\bar{x}(\cos \varpi T) + 2\lambda + \frac{2}{m\varpi} \int_0^T g(\lambda, t) \sin(\varpi t) dt = 0. \qquad (2.3.2.19)$$

From Eq.(2.3.2.19) we obtain

$$\bar{x} = -\left( \frac{\lambda}{\cos \varpi T} + \frac{1}{m\varpi \cos \varpi T} \int_0^T g(\lambda, t) \sin(\varpi t) dt \right) \qquad (2.3.2.20)$$

Substituting Eq.(2.3.2.20) into Eq.(2.3.2.17′) gives master equation of the form

$$u(T, \lambda) = 0,$$

$$u(T, \lambda) = -\frac{1}{\sin \varpi T} \left( \frac{\lambda}{\cos \varpi T} + \frac{1}{m\varpi \cos \varpi T} \int_0^T g(\lambda, t) \sin(\varpi t) dt \right) \sin \varpi T +$$

$$\frac{-C - D\sin(\Omega T) + (C + \lambda)\cos \varpi T}{\sin \varpi T} \sin \varpi T - (C + \lambda)\cos \varpi T +$$

$$+ D\sin(\Omega T) + C =$$

$$-\left( \frac{\lambda}{\cos \varpi T} + \frac{1}{m\varpi \cos \varpi T} \int_0^T g(\lambda, t) \sin(\varpi t) dt \right) = 0, \qquad (2.3.2.22)$$

$$g(\lambda, t) = -(m\omega^2 \lambda + 3a\lambda^2 - b - A\sin(\Omega t)) = -d(\lambda) + A\sin(\Omega t),$$

$$C = -\frac{d(\lambda)}{2m\varpi^2(\lambda)}, D = \frac{A}{(2\varpi^2(\lambda) - \Omega^2)m},$$

$$\varpi^2(\lambda) = 2\left( \frac{\omega^2}{2} + \frac{3a\lambda}{m} \right) \geq 0,$$

$$d(\lambda) = m\omega^2 \lambda + 3a\lambda^2 - b.$$

Finally we obtain master equation of the form

$$\lambda + \frac{1}{m\varpi} \int_0^T g(\lambda, t) \sin(\varpi t) dt = 0,$$

(2.3.2.23)

$$\cos[\varpi(\lambda)T] \neq 0, \varpi(\lambda) \neq 0.$$

By simple calculation we obtain

$$\int_0^T g(\lambda, t) \sin(\varpi t) dt = \int_0^T [-d(\lambda) + A \sin(\Omega t)] \sin(\varpi t) dt =$$

$$-d(\lambda) \int_0^T \sin(\varpi t) dt + A \int_0^T \sin(\varpi t) \sin(\Omega t) dt =$$

$$-d(\lambda) \int_0^T \sin(\varpi t) dt + \frac{A}{2} \left( \int_0^T \cos[(\varpi - \Omega)t] dt - \int_0^T \cos[(\varpi + \Omega)t] dt \right) =$$

$$-d(\lambda) \left( -\frac{\cos(\varpi t)}{\varpi} \Big|_0^T \right) + \frac{A}{2} \left( \frac{\sin[(\varpi - \Omega)t]}{\varpi - \Omega} \Big|_0^T \right) -$$

(2.3.24)

$$\frac{A}{2} \left( \frac{\sin[(\varpi + \Omega)t]}{\varpi + \Omega} \Big|_0^T \right) =$$

$$-d(\lambda) \left( -\frac{\cos(\varpi t)}{\varpi} + \frac{1}{\varpi} \right) + \frac{A}{2} \frac{\sin[(\varpi - \Omega)T]}{\varpi - \Omega} - \frac{A}{2} \frac{\sin[(\varpi + \Omega)T]}{\varpi + \Omega} =$$

$$\left( \frac{d(\lambda) \cos(\varpi t)}{\varpi} - \frac{d(\lambda)}{\varpi} \right) + \frac{A}{2} \frac{\sin[(\varpi - \Omega)T]}{\varpi - \Omega} - \frac{A}{2} \frac{\sin[(\varpi + \Omega)T]}{\varpi + \Omega}.$$

Substitution Eq.(2.3.24) into Eq.(2.3.23) gives

$$\lambda + \frac{1}{m\varpi(\lambda)} \times$$

$$\left[ \frac{d(\lambda)\cos(\varpi T)}{\varpi(\lambda)} - \frac{d(\lambda)}{\varpi(\lambda)} - \frac{A\sin[(\varpi(\lambda)-\Omega)T]}{2(\varpi(\lambda)-\Omega)} + \frac{A\sin[(\varpi(\lambda)+\Omega)T]}{2(\varpi(\lambda)+\Omega)} \right] = 0. \quad (2.3.2.25)$$

Let us consider now quantum harmonic oscillator,i.e.in Eqs.(2.3.22) $a = 0$ and therefore

$$d(\lambda, a = 0) = m\omega^2\lambda - b,$$

$$\varpi(\lambda) = \omega. \quad (2.3.2.26)$$

$$\lambda + \frac{1}{m\omega} \times$$

$$\left[ \frac{(m\omega^2\lambda - b)\cos(\omega T)}{\omega} - \frac{m\omega^2\lambda - b}{\omega} \right] = 0.$$

$$\lambda + \frac{(m\omega^2\lambda - b)\cos(\omega T)}{m\omega^2} - \frac{m\omega^2\lambda - b}{m\omega^2} = 0,$$

$$\lambda + \lambda\cos(\omega T) - \frac{b\cos(\omega T)}{m\omega^2} - \lambda + \frac{b}{m\omega^2} = 0, \qquad (2.3.2.27)$$

$$\frac{\lambda}{\cos\varpi T} + \lambda\cos(\omega T) - \frac{b\cos(\omega T)}{m\omega^2} + \frac{b}{m\omega^2} = 0,$$

$$\lambda = \frac{b}{m\omega^2} - \frac{b}{m\omega^2\cos(\omega T)} =$$

$$\frac{b}{m\omega^2}\left(1 - \frac{1}{\cos(\omega T)}\right).$$

**Example**.**2**.**1**.

$$m = 10, \Omega = 10, \omega = 30, a = -200, b = 1, A = 1$$

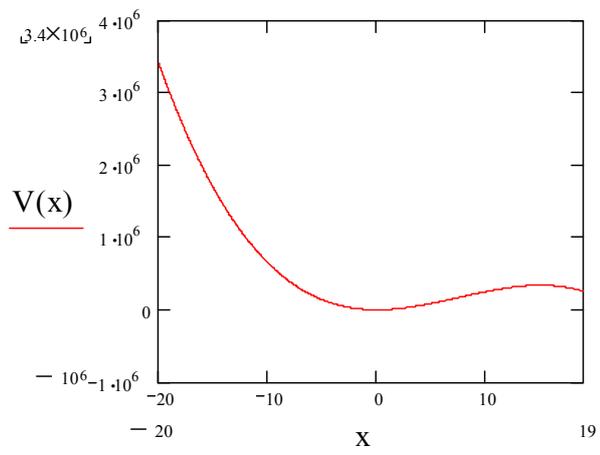

Fig. 2.1.1.

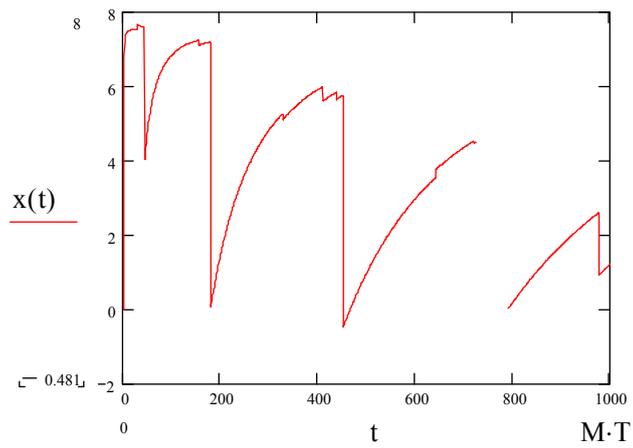

Fig. 2.1.2.a.

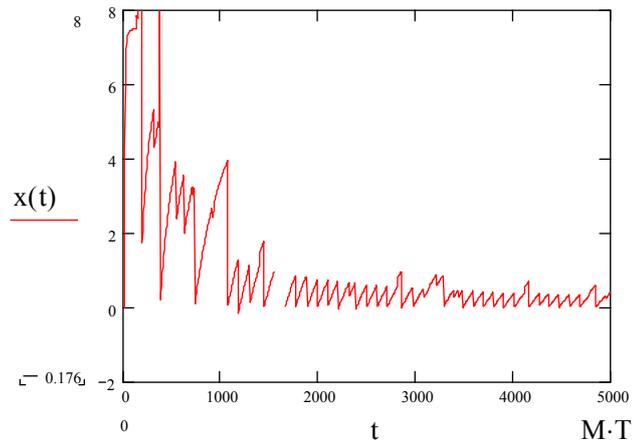

Fig. 2.1.2.b.

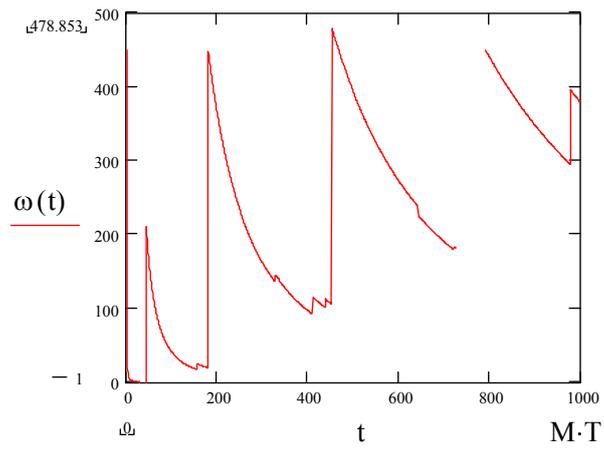

Fig. 2.1.3.a $\omega(t) \triangleq \varpi[\lambda(t)]$.

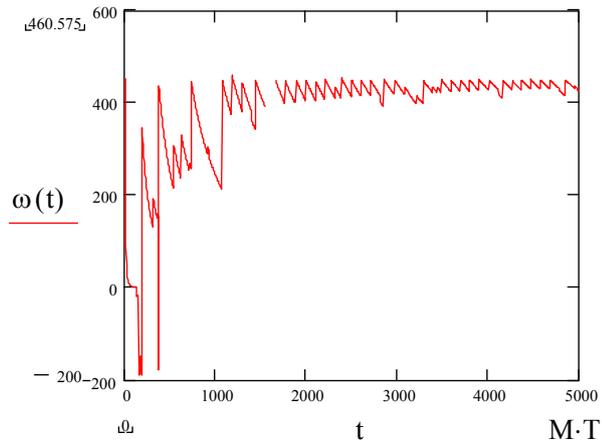

Fig. 2.1.3.b  $\omega(t) \triangleq \varpi[\lambda(t)]$.

**Example**.2.2.

$m = 10, \Omega = 10, \omega = 30, a = -200, b = 1, A = 10$

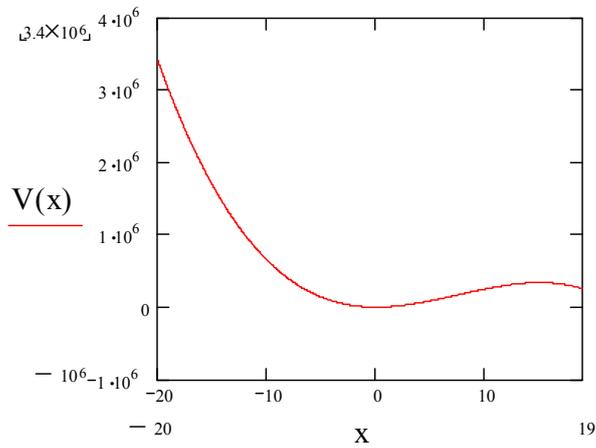

Fig. 2.2.1.

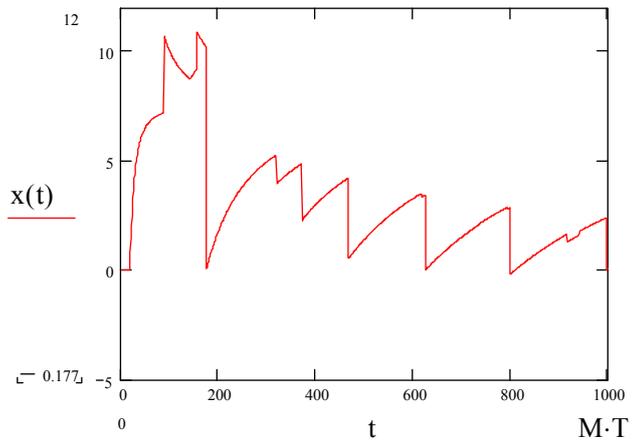

Fig. 2.2.2.a.

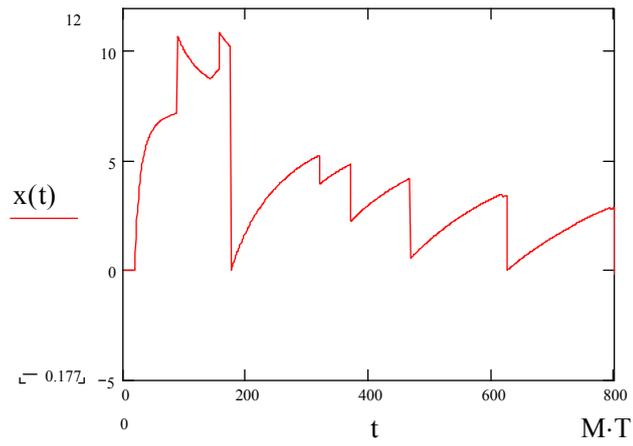

Fig. 2.2.2.b.

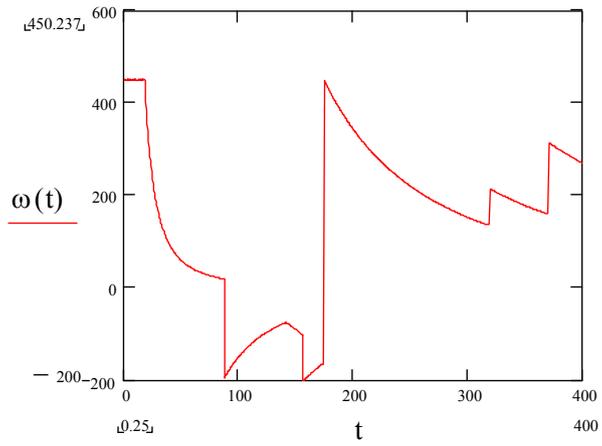

Fig.2.2.3.a.$\omega(t) \triangleq \varpi[\lambda(t)]$.

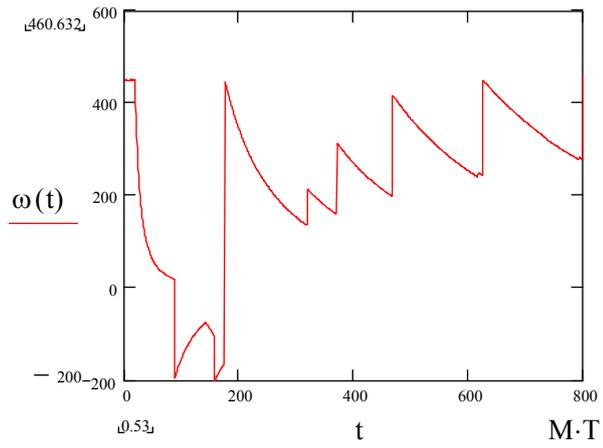

Fig.2.2.3.b.$\omega(t) \triangleq \varpi[\lambda(t)]$.

**Example**.2.3.

$m = 10, \Omega = 10, \omega = 40, a = -200, b = 1, A = 10$

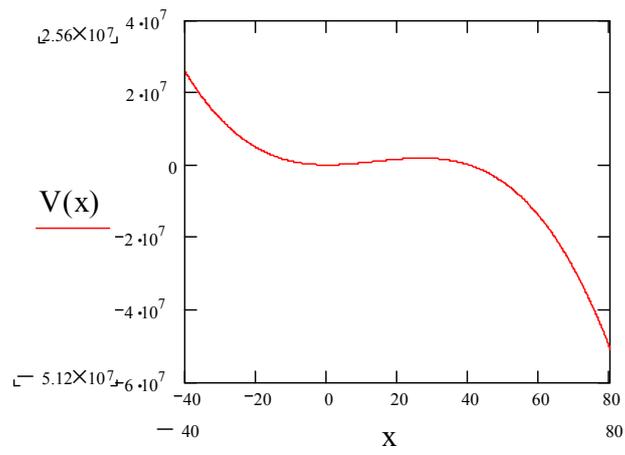

Fig. 2.3.1.

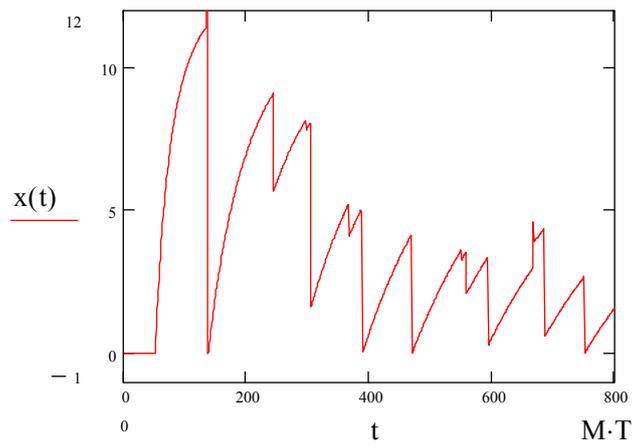

Fig. 2.3.2.

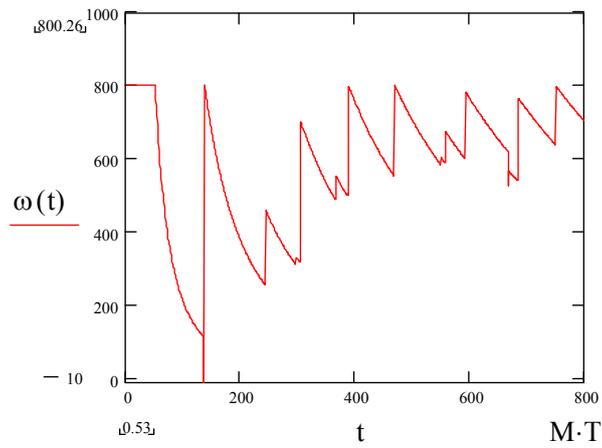

Fig. 2.3.3. $\omega(t) \triangleq \varpi[\lambda(t)]$.

**Example**.2.4.

$$m = 10, \Omega = 20, \omega = 40, a = -200, b = 1, A = 10$$

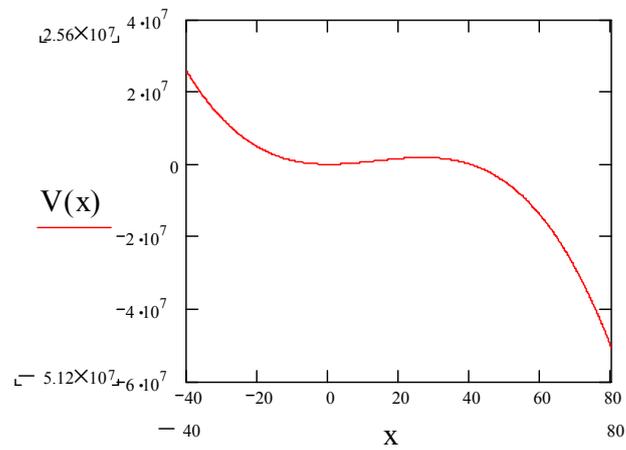

Fig. 2.4.1.

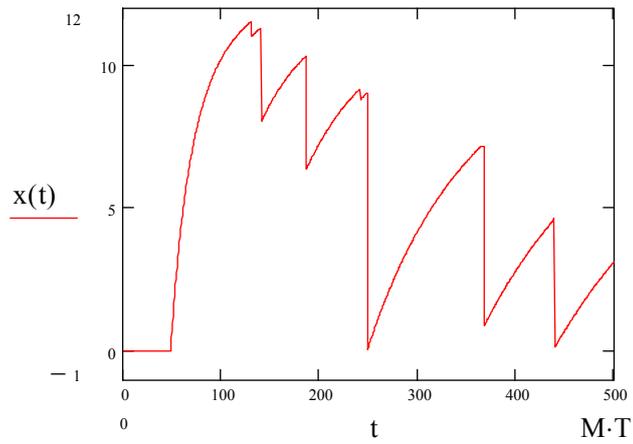

Fig. 2.4.2.

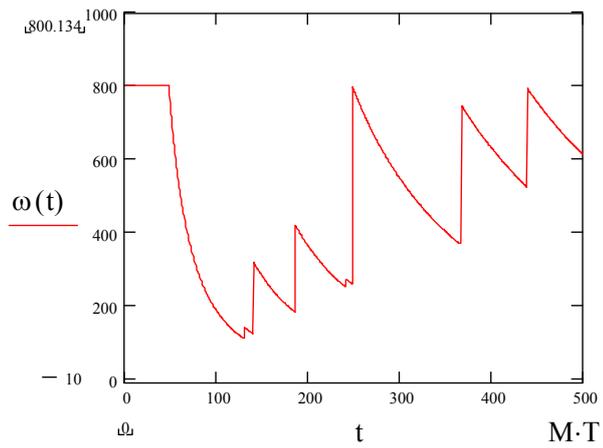

Fig. 2.4.3. $\omega(t) \triangleq \varpi[\lambda(t)]$.

**Example**.2.5.

$m = 10, \Omega = 25, \omega = 40, a = -200, b = 1, A = 10$

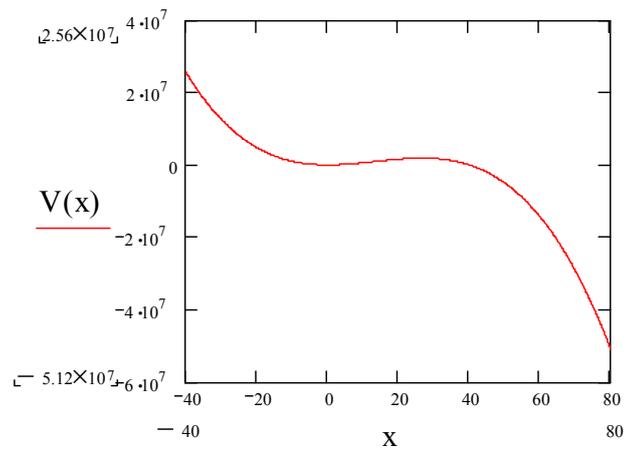

Fig. 2.5.1.

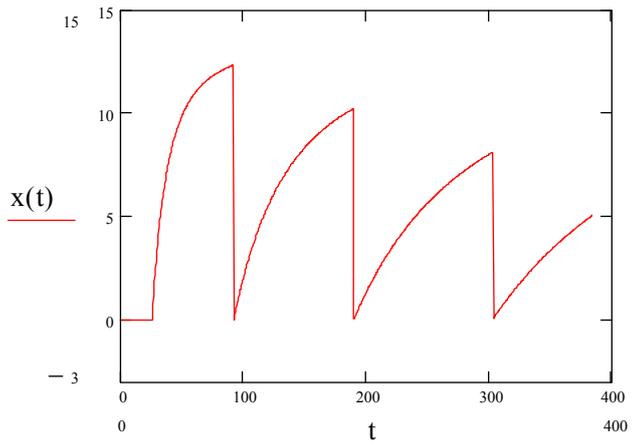

Fig. 2.5.2.

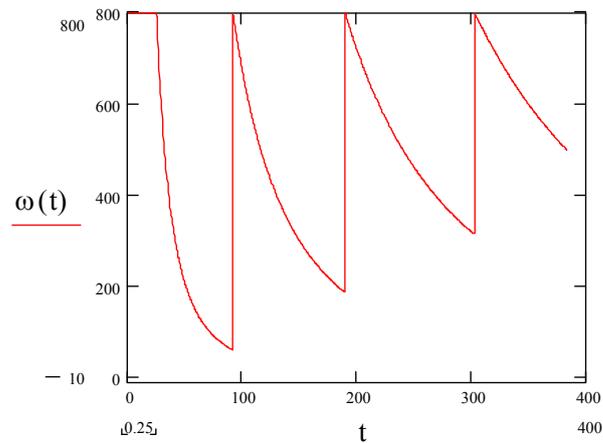

Fig. 2.5.3. $\omega(t) \triangleq \varpi[\lambda(t)]$.

## II.4. Quantum anharmonic oscillator with a duble well potential.

As a a second example we consider quantum anharmonic oscillator with a duble well potential (as cartooned in fig.2.4.1),

$$V(x) = \frac{a}{4} x^4 - \frac{b}{2} x^2 \ , \ \ a, b > 0 \qquad\qquad (2.4.1)$$

supplemented by an additive sinusoidal driving.

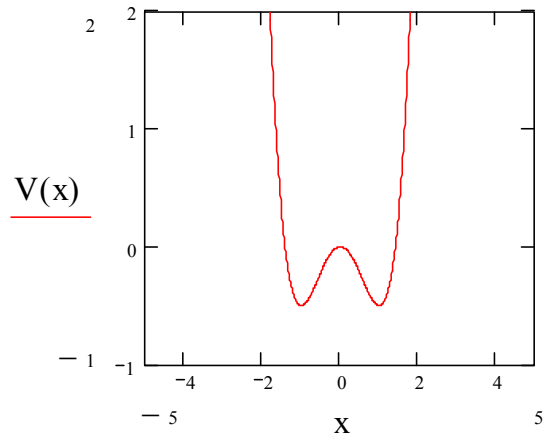

Fig.2.4.1. $a = 0.5, b = 1$

Let us consider transition probability amplitude given via formulae (2.1.1)-

(2.1.2) with

$$V(x) = \frac{m\omega^2}{2}x^2 - ax^4 + bx^3 - [A\sin(\Omega t)]x.$$

$$(2.4.2)$$

The corresponding classical Lagrangian is

$$\mathcal{L}(\dot{x}, x) = \frac{m}{2}\dot{x}^2 - \frac{m\omega^2}{2}x^2 - bx^3 + ax^4 + [A\sin(\Omega t)]x.$$

$$(2.4.3)$$

By substitution $x = u + \lambda$ we obtain

$$\mathcal{L}(\dot{u}, u) = \frac{m}{2}\dot{u}^2 - \frac{m\omega^2}{2}(u+\lambda)^2 - b(u+\lambda)^3 + a(u+\lambda)^4 +$$

$$[A\sin(\Omega t)](u+\lambda) =$$

$$\frac{m}{2}\dot{u}^2 - \frac{m\omega^2}{2}u^2 - m\omega^2\lambda u - 3b\lambda u^2 - 3b\lambda^2 u + 6a\lambda^2 u^2 + 4a\lambda^3 u + \ldots = \qquad (2.4.4)$$

$$\frac{m}{2}\dot{u}^2 - \left(\frac{m\omega^2}{2} + 3b\lambda - 6a\lambda^2\right)u^2 - (m\omega^2\lambda + 3b\lambda^2 - 4a\lambda^3)u +$$

$$[A\sin(\Omega t)]u + au^4 - 4a\lambda u^3 + a\lambda^4 + [A\sin(\Omega t)]\lambda + \ldots$$

The corresponding master Lagrangian $\widetilde{\mathcal{L}}(\dot{u}, u)$ is

$$\widetilde{\mathcal{L}}(\dot{u}, u) = \frac{m}{2}\dot{u}^2 - \left(\frac{m\omega^2}{2} + 3b\lambda - 6a\lambda^2\right)u^2 - (m\omega^2\lambda + 3b\lambda^2 - 4a\lambda^3 - A\sin(\Omega t))u =$$

$$(2.4.5)$$

$$\frac{m}{2}\dot{u}^2 - m\left(\frac{\omega^2}{2} + \frac{3b\lambda}{m} - \frac{6a\lambda^2}{m}\right)u^2 - (m\omega^2\lambda + 3b\lambda^2 - 4a\lambda^3 - A\sin(\Omega t))u.$$

The Euler-Lagrange equation of $\widetilde{\mathcal{L}}(\dot{u}, u)$ with corresponding boundary conditions is

$$\frac{d}{dt}\left(\frac{\partial\widetilde{\mathcal{L}}}{\partial\dot{u}}\right) - \frac{\partial\widetilde{\mathcal{L}}}{\partial u} = 0,$$

$$(2.4.6)$$

$$u(0) = -\lambda, u(T) = x.$$

**Remark 2.4.1.** We assume now that $2\varpi^2(\lambda) = \frac{\omega^2}{2} + \frac{3b\lambda}{m} - \frac{6a\lambda^2}{m} \geq 0$ and rewrite Eq.(2.4.5) of the form

$$\widetilde{\mathcal{L}}(\dot{u}, u) = \frac{m}{2}\dot{u}^2 - \frac{m\varpi^2(\lambda)}{2}u^2 + g(\lambda, t)u,$$

$$\varpi(\lambda) = \sqrt{2\left(\frac{\omega^2}{2} + \frac{3b\lambda}{m} - \frac{6a\lambda^2}{m}\right)}, \qquad (2.4.7)$$

$$g(\lambda, t) = -[m\omega^2\lambda + 3b\lambda^2 - 4a\lambda^3 - A\sin(\Omega t)].$$

Finally we obtain

$$\lambda + \frac{1}{m\varpi(\lambda)} \times$$

$$\left[\frac{d(\lambda)\cos(\varpi T)}{\varpi(\lambda)} - \frac{d(\lambda)}{\varpi(\lambda)} + \frac{A\sin[(\varpi(\lambda) - \Omega)T]}{2(\varpi(\lambda) - \Omega)} - \frac{A\sin[(\varpi(\lambda) + \Omega)T]}{2(\varpi(\lambda) + \Omega)}\right] = 0,$$

$$\varpi(\lambda) = \sqrt{2\left(\frac{\omega^2}{2} + \frac{3b\lambda}{m} - \frac{6a\lambda^2}{m}\right)}, \qquad (2.4.8)$$

$$d(\lambda) = m\omega^2\lambda + 3b\lambda^2 - 4a\lambda^3.$$

**Example**.**2**.**1**.

$$m = 10, \Omega = 20, \omega = 10, a = -20, b = 10, A = 10$$

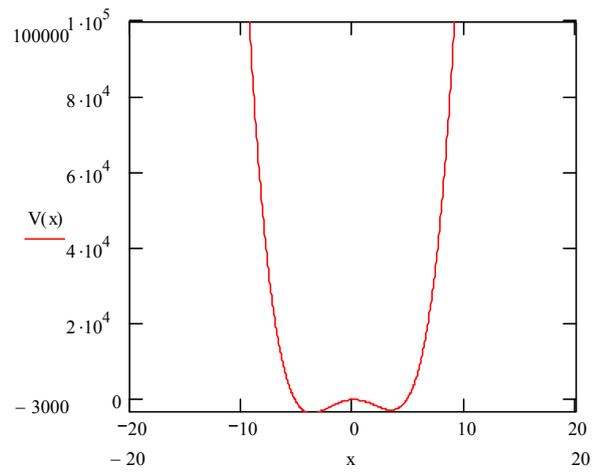

2.1.1.

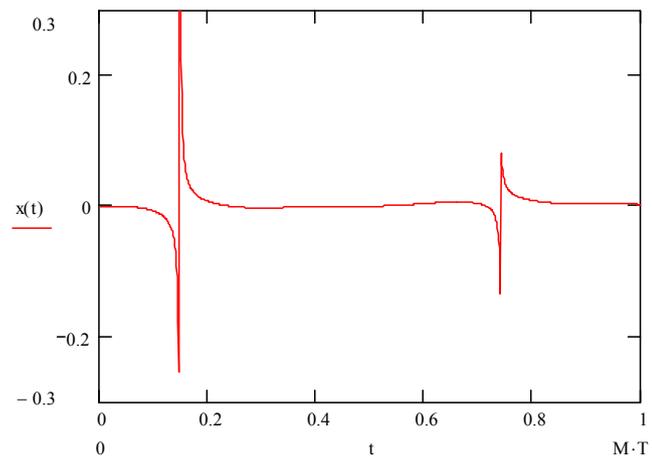

2.1.2.

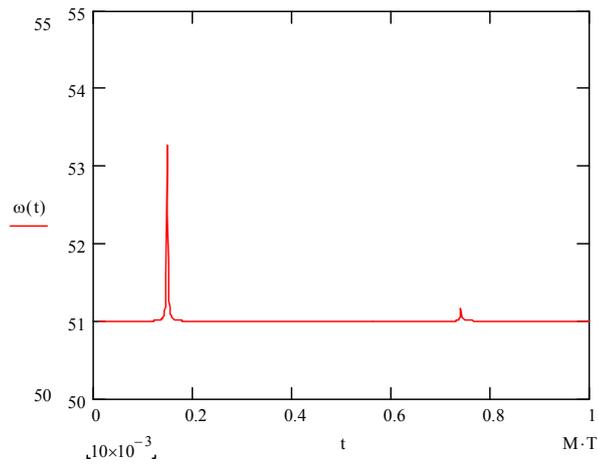

Fig. 2.1.3. $\omega(t) \triangleq \varpi[\lambda(t)]$.

**Example**.2.2.

$m = 100, \Omega = 20, \omega = 10, a = -20, b = 50, A = 10$

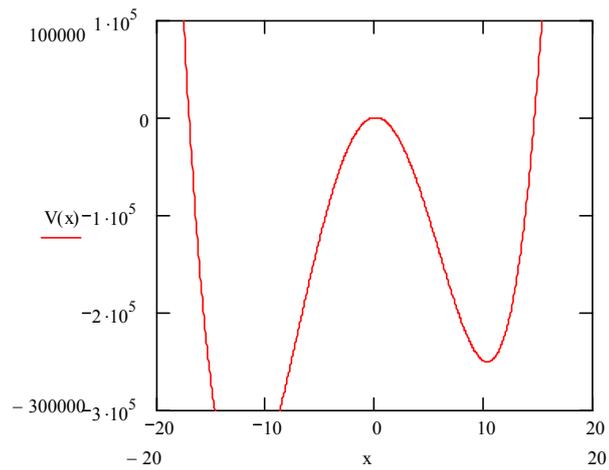

2.1.1.

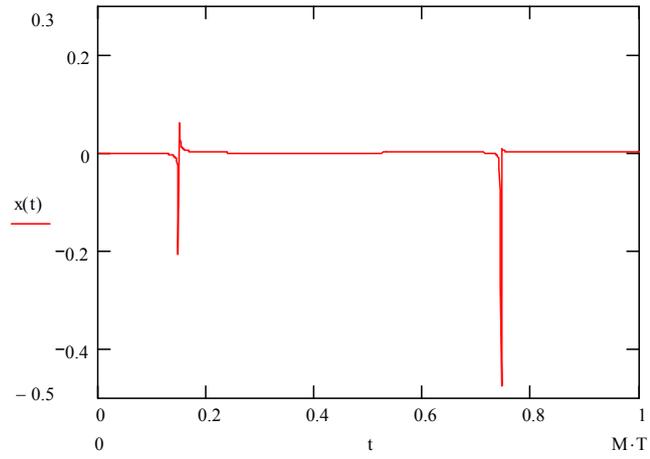

2.1.2.

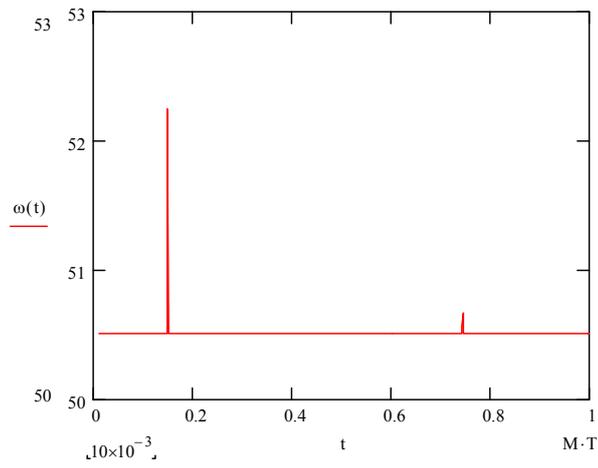

Fig. 2.1.3. $\omega(t) \triangleq \varpi[\lambda(t)]$.

**Example**.2.3.

$m = 100, \Omega = 20, \omega = 10, a = 20, b = 50, A = 100$

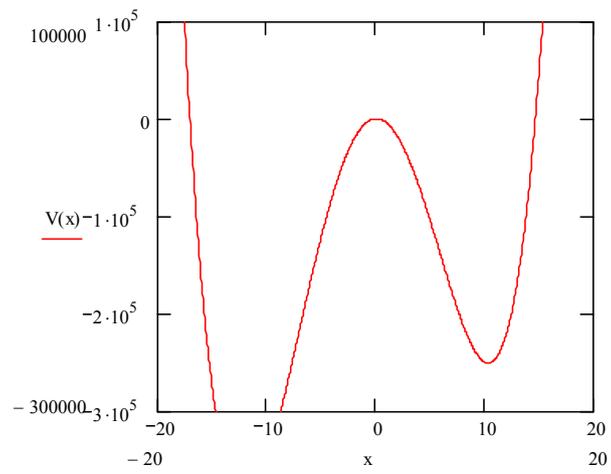

2.1.1.

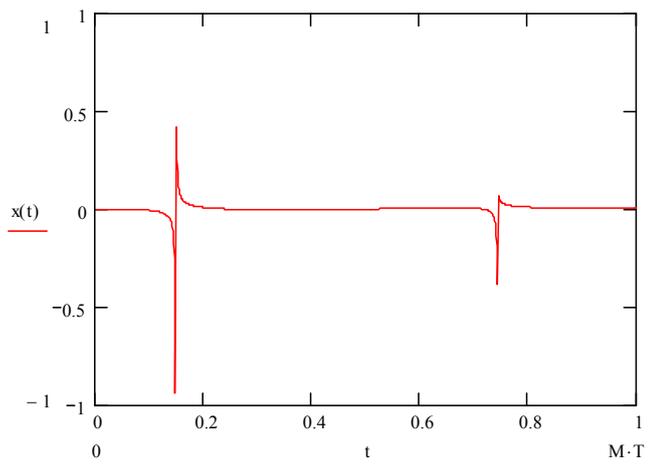

2.1.2.

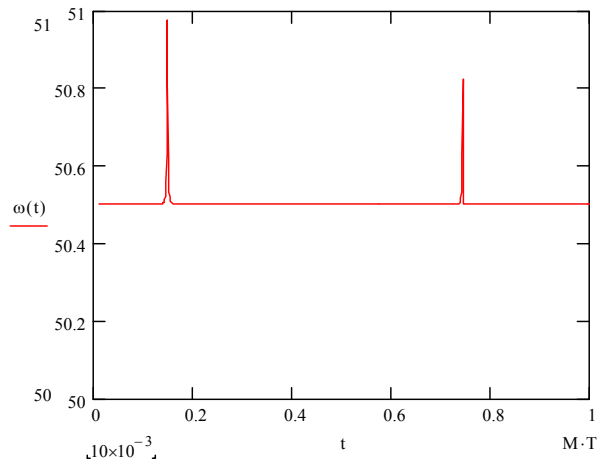

Fig. 2.1.3. $\omega(t) \triangleq \varpi[\lambda(t)]$.

# II.5. Quantum tradjectories of the one dimensional quantum anharmonic oscillator with a cubic potential.

Let's calculate canonical quantum tradjectories (canonical quantum averages) of the one dimensional quantum anharmonic oscillator with a cubic potential,i.e.

$$\langle x_i, t, x_0; \hbar \rangle = \frac{\int\limits_{-\infty}^{\infty} x |\Psi(x, t, x_0; \hbar)|^2 d\,x}{\int\limits_{-\infty}^{\infty} |\Psi(x, t, x_0; \hbar)|^2 dx},$$

$$\Psi(x, 0, x_0; \hbar) = \frac{1}{\sqrt[4]{2\pi\varepsilon}} \exp\left(-\frac{x^2}{\varepsilon}\right),$$

$$\Psi^2(x, 0, x_0; \hbar) \approx \delta(x - x_0),$$

$$\varepsilon \ll \hbar \ll 1.$$

(2.5.1)

Corresponding master action $\widetilde{S}_{\text{cl}}$ is :

$$\widetilde{S}_{\text{cl}}(\bar{x}, y, \lambda, T) = \frac{m\varpi}{2 \sin \varpi T} \bigg[ (\cos \varpi T)(y^2 + \bar{x}^2) - 2\bar{x}y + \frac{2\bar{x}}{m\varpi} \int_0^T g(\lambda, t) \sin(\varpi t) dt - $$

$$+ \frac{2y}{m\varpi} \int_0^T g(\lambda, t) \sin \varpi (T - t) dt - $$

$$- \frac{2}{m^2 \varpi^2} \int_0^T \int_0^t g(\lambda, t) g(\lambda, s) \sin \varpi (T - t) \sin(\varpi s) ds dt \bigg],$$

(2.5.2)

$$\bar{x} = x - \lambda.$$

**Remark**. We remaind that

$$\int\limits_{x(T)=x} (x(\tau) - \lambda) \Psi(x(0)) \exp[\Phi(x(\tau))] Dx(\tau) =$$

$$\left[ \int\limits_{\substack{x(T)=x \\ x(0)=y}} (x(\tau) - \lambda) \Psi(x(0) = y) \exp[\Phi(x(\tau))] Dx(\tau) \right] dy =$$

(2.5.2′)

$$\left[ \int\limits_{\substack{u(T)=x-\lambda \\ u(0)=y-\lambda}} (u(\tau)) \Psi(u(0) + \lambda) \exp[\Phi(u(\tau) + \lambda)] Du(\tau) \right] dy.$$

Here $x(\tau) - \lambda = u(\tau)$, $x(\tau) = u(\tau) + \lambda$ .

Therefore from Eq.(2.5.2$'$) and Maslov formula (0.1.5) we obtain:

$$\Psi(\bar{x}, \lambda, T) = \frac{\sqrt{m\varpi}}{\sqrt[4]{\varepsilon}\,\sqrt{2\pi i\hbar}\,|\sin\varpi T|^{1/2}}\int\limits_{-\infty}^{\infty} dy \exp\left(-\frac{(y+\lambda)^2}{\varepsilon}\right) \times$$

$$\times \exp\left[\frac{i}{\hbar}\widetilde{S}_{\mathrm{el}}(\bar{x}, y, \lambda, T) + \frac{i\pi\gamma}{2}\right] =$$

$$\frac{\sqrt{m\varpi}}{\sqrt[4]{\varepsilon}\,\sqrt{2\pi i\hbar}\,|\sin\varpi T|^{1/2}}\int\limits_{-\infty}^{\infty} dy \exp\left(-\frac{(y+\lambda)^2}{\varepsilon}\right) \times$$

$$\exp\left\{\frac{im\varpi}{2\hbar\sin\varpi T} \times\right.$$

$$\left[(\cos\varpi T)(y^2+\bar{x}^2) - 2\bar{x}y + \frac{2\bar{x}}{m\varpi}\int_0^T g(\lambda, t)\sin(\varpi t)dt +\right.$$

$$\frac{2y}{m\varpi}\int_0^T g(\lambda, t)\sin\varpi(T-t)dt -$$

$$\left.\left.-\frac{2}{m^2\varpi^2}\int_0^T\int_0^t g(\lambda, t)g(\lambda, s)\sin\varpi(T-t)\sin(\varpi s)dsdt\right] + \frac{i\pi\gamma}{2}\right\},$$

(2.5.3)

Let's to calculate integral (2.5.3) using stationary phase method.

$q(t) - \lambda = u(t)$

$q(0) - \lambda = u(0), q(0) = u(0) + \lambda$

$q(T) - \lambda = u(T), q(T) = u(T) + \lambda$

$q(T) = x, u(T) = x - \lambda$

$q(0) = y, u(0) = y - \lambda$

$$\int\limits_{-\infty}^{\infty} dy \int\limits_{q(T)=x} [q(T)-\lambda]\exp\left(-[q(0)]^2\right) = \int\limits_{-\infty}^{\infty} dy \int\limits_{u(T)+\lambda=x} u(T)\exp\left(-[u(0)+\lambda]^2\right) =$$

$$\int\limits_{-\infty}^{\infty} dy \int\limits_{\substack{u(T)=x-\lambda \\ u(0)=y-\lambda}} u(T)\exp(y^2), \qquad \frac{-\left[\sqrt{\hbar}\,u'(0)+\lambda\right]^2}{\hbar} = -\left[u'(0)+\frac{\lambda}{\sqrt{\hbar}}\right]^2$$

$$\frac{a+b}{2} - a = \frac{b}{2} - \frac{a}{2}$$

$$(x-a)^2 + (x-b)^2$$

$$2x - a - b = 0$$

$$\frac{\partial \widetilde{S}_{\mathbf{el}}}{\partial y} = 2y \cos \varpi T - 2\bar{x} + \frac{2}{m\varpi} \int_0^T g(\lambda, t) \sin \varpi (T-t) dt = 0 \qquad (2.5.4)$$

Therefore

$$y_{\mathbf{cr}}(\bar{x}, \lambda, T) = \frac{\bar{x}}{\cos \varpi T} - \frac{1}{m\varpi \cos \varpi T} \int_0^T g(\lambda, t) \sin \varpi (T-t) dt =$$

$$\frac{1}{\cos \varpi T} \left( \bar{x} - \frac{1}{m\varpi} \int_0^T g(\lambda, t) \sin \varpi (T-t) dt \right). \qquad (2.5.5)$$

Simple calculation gives

$$\int_0^T g(\lambda, t) \sin[\varpi(T-t)] dt = \int_0^T [-d(\lambda) + A\sin(\Omega t)] \sin[\varpi(T-t)] dt =$$

$$-d(\lambda) \int_0^T \sin[\varpi(T-t)] dt + A \int_0^T \sin[\varpi(T-t)] \sin(\Omega t) dt =$$

$$d(\lambda) \int_0^T \sin[\varpi(T-t)] d(T-t) +$$

$$\frac{A}{2} \left( \int_0^T \cos[\varpi(T-t) - \Omega t] dt - \int_0^T \cos[\varpi(T-t) + \Omega t] dt \right) =$$

$$d(\lambda) \left( -\frac{\cos[\varpi(T-t)]}{\varpi} \Big|_0^T \right) - \frac{A}{2} \left( \frac{\sin[\varpi(T-t) - \Omega t]}{\varpi + \Omega} \Big|_0^T \right) +$$

$$+ \frac{A}{2} \left( \frac{\sin[\varpi(T-t) + \Omega t]}{\varpi - \Omega} \Big|_0^T \right) =$$

$$d(\lambda) \left( \frac{\cos(\varpi T)}{\varpi} - \frac{1}{\varpi} \right) - \frac{A}{2} \left( \frac{\sin(-\Omega T)}{\varpi + \Omega} - \frac{\sin(\varpi T)}{\varpi + \Omega} \right) + \qquad (2.5.6)$$

$$+ \frac{A}{2} \left( \frac{\sin(\Omega T)}{\varpi - \Omega} - \frac{\sin(\varpi T)}{\varpi - \Omega} \right) =$$

$$d(\lambda) \left( \frac{\cos(\varpi T)}{\varpi} - \frac{1}{\varpi} \right) +$$

$$\frac{A}{2} \frac{\sin(\Omega T) + \sin(\varpi T)}{\varpi + \Omega} + \frac{A}{2} \frac{\sin(\Omega T) - \sin(\varpi T)}{\varpi - \Omega} =$$

$$d(\lambda) \left( \frac{\cos(\varpi T)}{\varpi} - \frac{1}{\varpi} \right) +$$

$$\frac{A}{2} \frac{[\sin(\Omega T) + \sin(\varpi T)](\varpi - \Omega) + [\sin(\Omega T) - \sin(\varpi T)](\varpi + \Omega)}{\varpi^2 - \Omega^2} =$$

$$d(\lambda) \left( \frac{\cos(\varpi T)}{\varpi} - \frac{1}{\varpi} \right) + A \frac{\varpi \sin(\Omega T) - \Omega \sin(\varpi T)}{\varpi^2 - \Omega^2}.$$

Using stationary phase method, from Eqs.(2.5.3)-(2.5.5) we obtain

$$\Psi(\bar{x}, \lambda, T) = \frac{1}{\sqrt[4]{2\pi\varepsilon} \, |\cos \varpi T|^{1/2}} \exp\left[ -\frac{(y_{\mathrm{cr}}(\bar{x}, \lambda, T) + \lambda)^2}{\varepsilon} \right] \times$$

$$\exp\left[ \frac{i}{h} \widetilde{S}_{\mathrm{cl}}(\bar{x}, y_{\mathrm{cr}}, \lambda, T) + \frac{i\pi\gamma}{2} + \frac{i\pi\gamma'}{4} \right] + O(h) =$$

$$\frac{1}{\sqrt[4]{2\pi\varepsilon} \, |\cos \varpi T|^{1/2}} \times \tag{2.5.7}$$

$$\times \exp\left[ -\frac{1}{\varepsilon} \left( \frac{\bar{x}}{\cos \varpi T} - \frac{1}{m\varpi \cos \varpi T} \int_0^T g(\lambda, t) \sin \varpi(T - t) dt + \lambda \right)^2 \right] \times$$

$$\exp\left[ \frac{i}{h} \widetilde{S}_{\mathrm{cl}}(\bar{x}, y_{\mathrm{cr}}, \lambda, T) + \frac{i\pi\gamma}{2} + \frac{i\pi\gamma'}{4} \right] + O(h).$$

Here

$$\gamma' = \mathbf{sgn}\left( \frac{\partial^2 \widetilde{S}_{\mathrm{cl}}}{\partial y^2} \right). \tag{2.5.8}$$

From Eq.(2.5.7) we obtain

$$|\Psi(\bar{x}, \lambda, T)|^2 = \frac{1}{\sqrt{2\pi\varepsilon} \, |\cos \varpi T|} \exp\left[ -\frac{2[y_{\mathrm{cr}}(\bar{x}, \lambda, T) + \lambda]^2}{\varepsilon} \right] \tag{2.5.9}$$

Using now Laplace method, from Eq.(2.5.5) and Eq.(2.5.9) we obtain

$$\int\limits_{-\infty}^{\infty} d\bar{x}\bar{x} |\Psi(\bar{x}, \lambda, T)|^2 =$$

$$\frac{1}{\sqrt{2\pi\varepsilon} \, |\cos \varpi T|} \int\limits_{-\infty}^{\infty} d\bar{x}\bar{x} \exp\left(-\frac{2[y_{\text{cr}}(\bar{x}, \lambda, T) + \lambda]^2}{\varepsilon}\right) =$$

(2.5.10)

$$\frac{1}{\sqrt{\varepsilon} \, |\cos \varpi T|} \times$$

$$\times \int\limits_{-\infty}^{\infty} d\bar{x}\bar{x} \exp\left\{-\frac{2}{\varepsilon}\left[\left(\frac{\bar{x}}{\cos \varpi T} - \frac{1}{m\varpi \cos \varpi T}\int_0^T g(\lambda, t)\sin \varpi(T-t)dt + \lambda\right)^2\right]\right\}.$$

Therefore

$$\bar{x}_{\text{cr}}(\lambda, T) = \frac{1}{m\varpi}\int_0^T g(\lambda, t)\sin \varpi(T-t)dt - \lambda\cos \varpi T,$$

(2.5.11)

and master equation is:

$$\bar{x}_{\text{cr}}(\lambda, T) = \frac{1}{m\varpi}\int_0^T g(\lambda, t)\sin \varpi(T-t)dt - \lambda\cos \varpi T = 0.$$

(2.5.12)

Using Eq.(2.5.6) finally we obtain master equation of the form

$$d(\lambda)\left(\frac{\cos(\varpi T)}{\varpi} - \frac{1}{\varpi}\right) + A\frac{\varpi\sin(\Omega T) - \Omega\sin(\varpi T)}{\varpi^2 - \Omega^2} - \lambda m\varpi\cos(\varpi T) = 0.$$

(2.5.13)

Let us consider now quantum harmonic oscillator, i.e. $a = 0$ and therefore

$$d(\lambda) = m\omega^2\lambda - b,$$

$$\varpi(\lambda) = \omega. \tag{2.5.14}$$

Substitution Eq.(2.5.14) into Eq.(2.5.13) gives

$$(m\omega^2\lambda - b)\left(\frac{\cos(\omega t)}{\omega} - \frac{1}{\omega}\right) + A\frac{\omega\sin(\Omega T) - \Omega\sin(\omega T)}{\omega^2 - \Omega^2} - \lambda m\omega\cos(\omega T) = 0,$$

$$m\omega\lambda\cos(\omega t) - m\omega\lambda - b\left(\frac{\cos(\omega t)}{\omega} - \frac{1}{\omega}\right) +$$

$$A\frac{\omega\sin(\Omega T) - \Omega\sin(\omega T)}{\omega^2 - \Omega^2} - \lambda m\omega\cos(\omega T) = 0, \tag{2.5.15}$$

$$-m\omega\lambda - b\left(\frac{\cos(\omega t)}{\omega} - \frac{1}{\omega}\right) + A\frac{\omega\sin(\Omega T) - \Omega\sin(\omega T)}{\omega^2 - \Omega^2} = 0.$$

Thus

$$\lambda = -\frac{b}{m\omega}\left(\frac{\cos(\omega t)}{\omega} - \frac{1}{\omega}\right) + A\frac{\omega\sin(\Omega T) - \Omega\sin(\omega T)}{m\omega(\omega^2 - \Omega^2)} \tag{2.5.16}$$

**Example**.**2**.**5**.**1**.(Linear Quantum harmonic oscillator, $a = 0$). **No quantum jumps**.

$$m = 1, \Omega = 0, \omega = 9, a = 0, b = 10, A = 0.$$

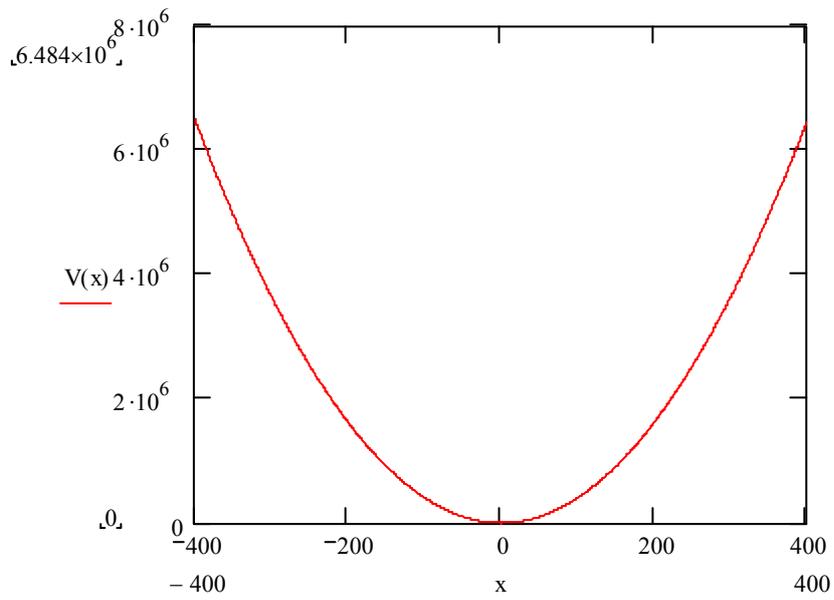

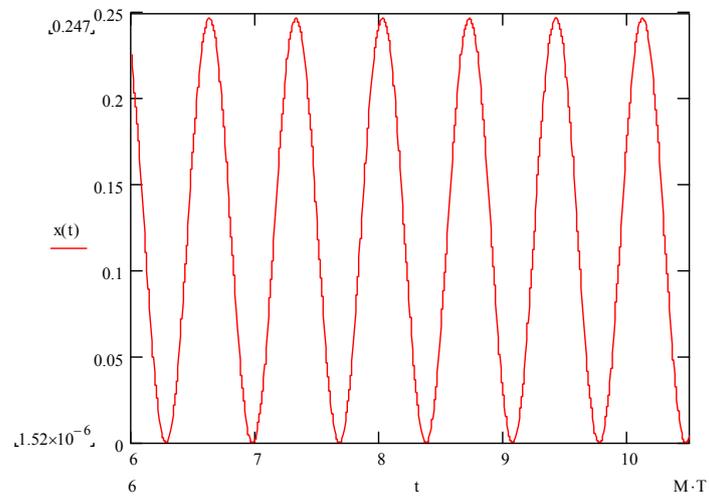

**Example**.2.5.2.

$m = 1, \Omega = 0, \omega = 9, a = 3, b = 10, A = 0.$

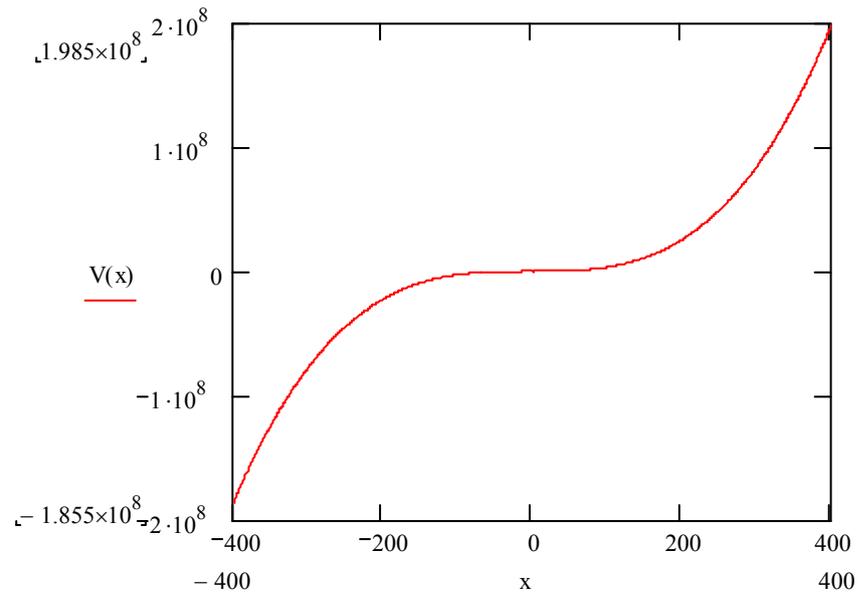

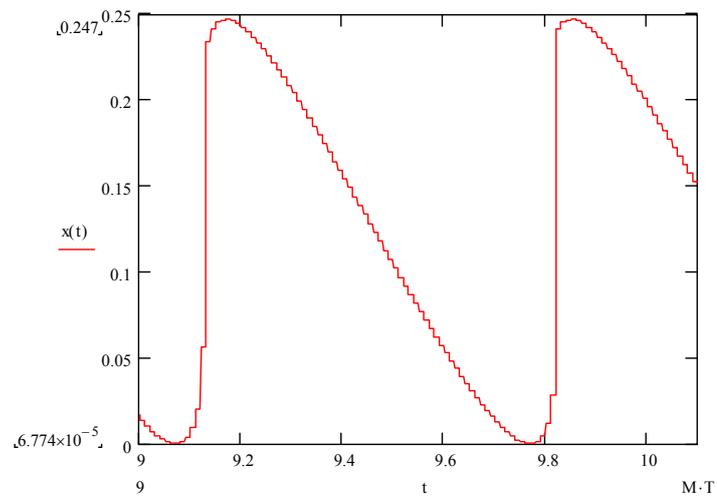

$m = 1, \Omega = 0, \omega = 9, a = 3, b = 10, A = 0.$

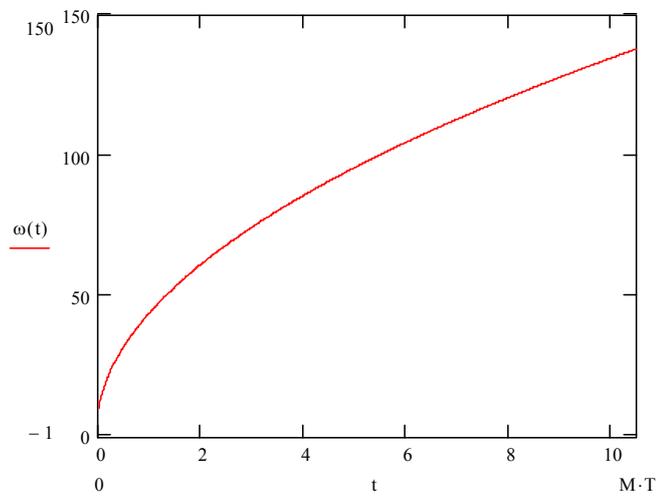

**Example**.2.5.3.

$m = 1, \Omega = 3, \omega = 9, a = 3, b = 10, A = 1.$

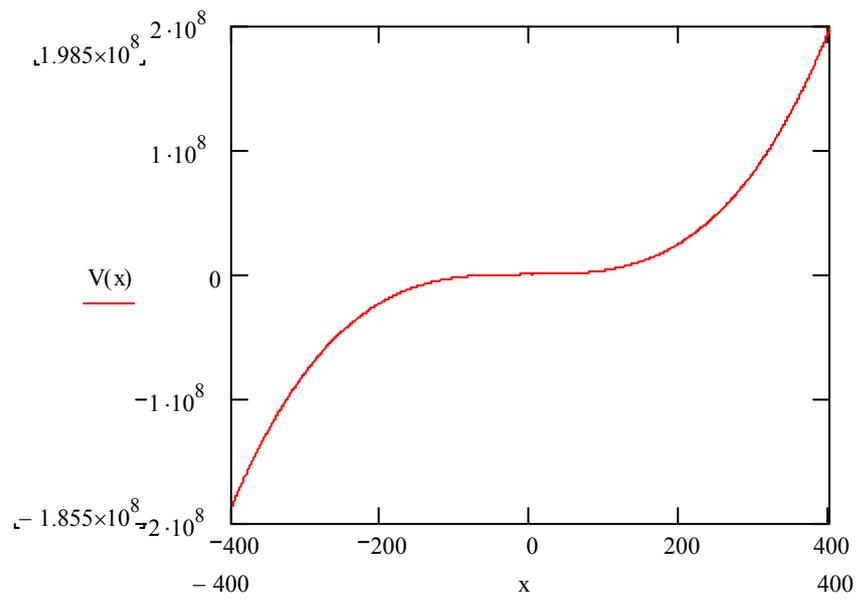

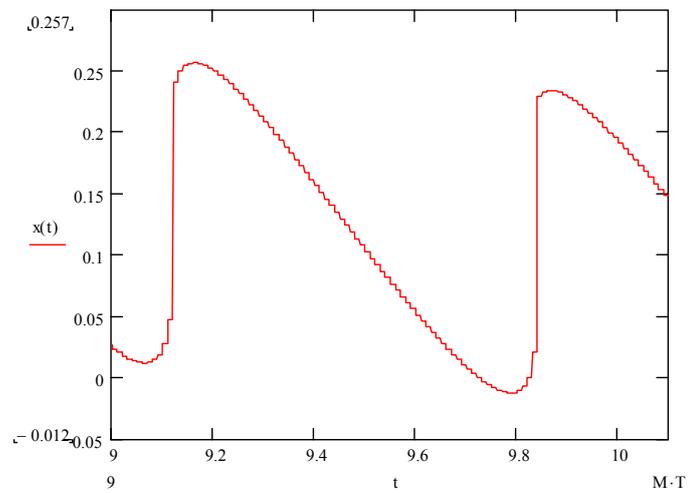

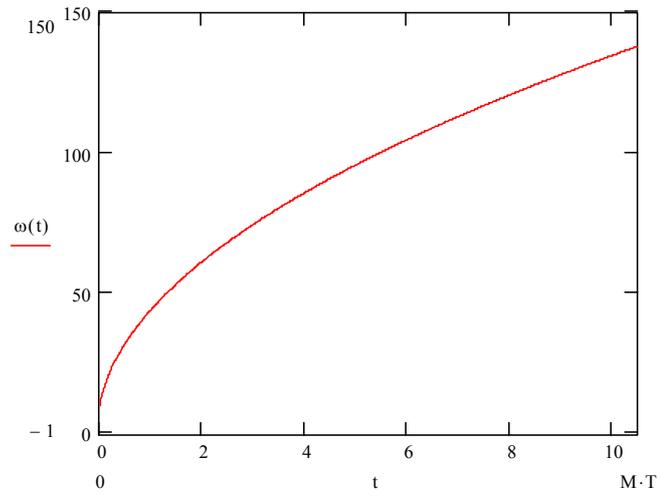

**Example**.**2.5.4**.

$m = 1, \Omega = 3, \omega = 15, a = 3, b = 10, A = 1.$

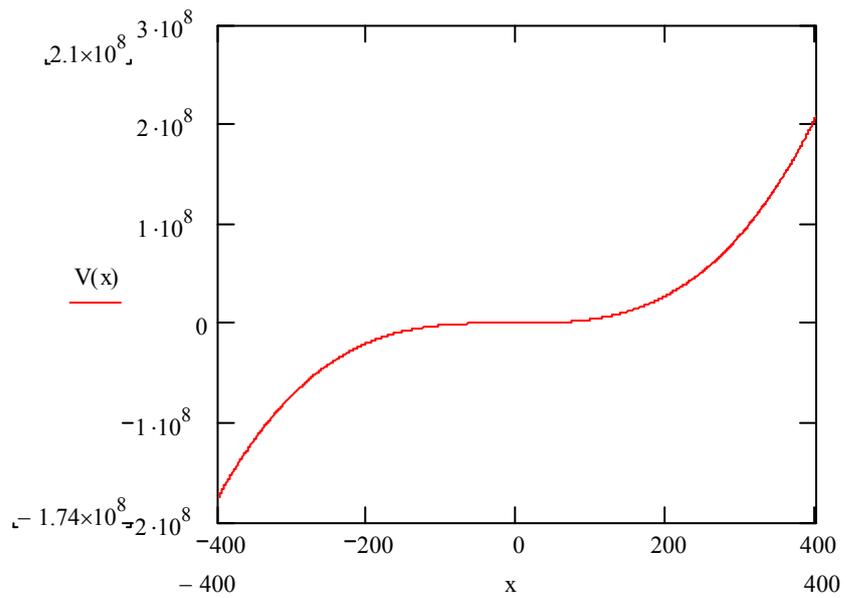

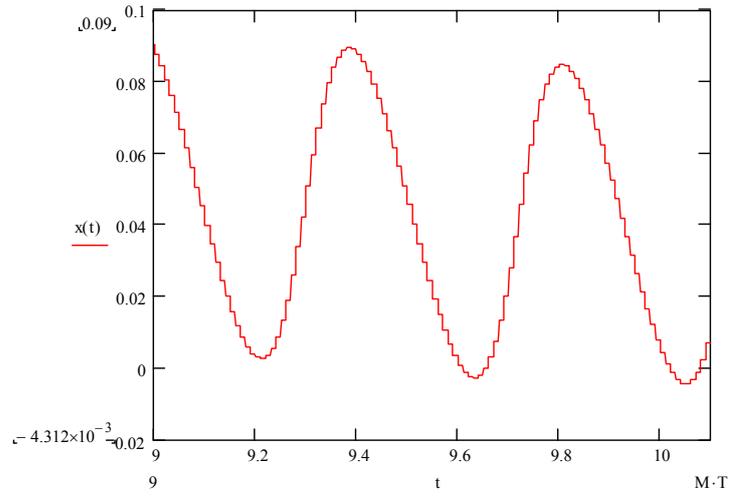

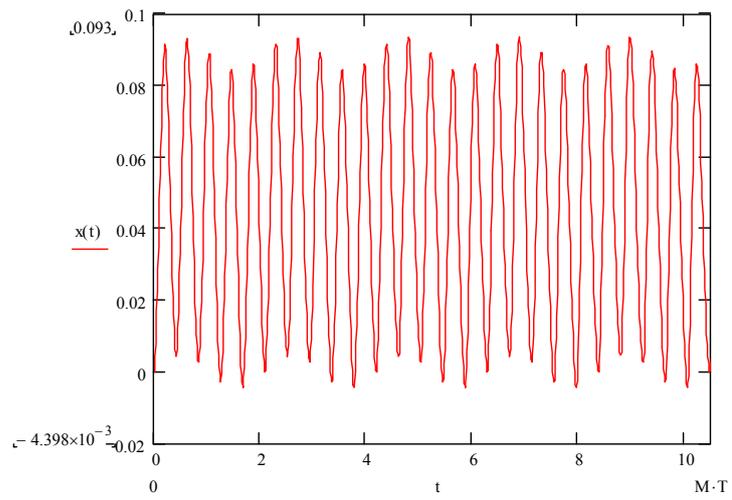

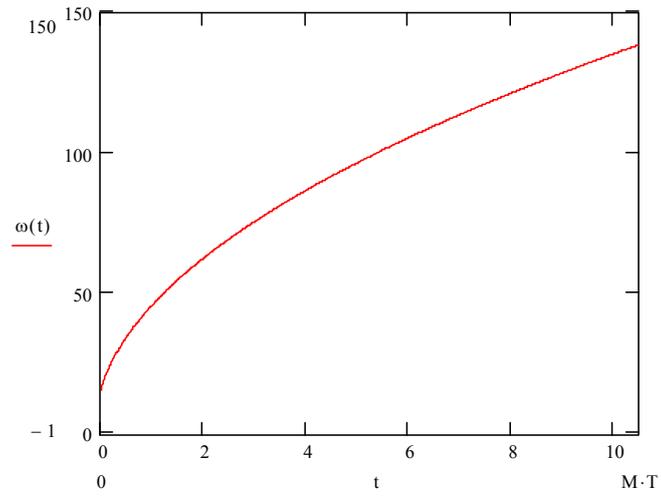

**Example**.**2.5.5**.

$m = 1, \Omega = 3, \omega = 10, a = 3, b = 10, A = 3.$

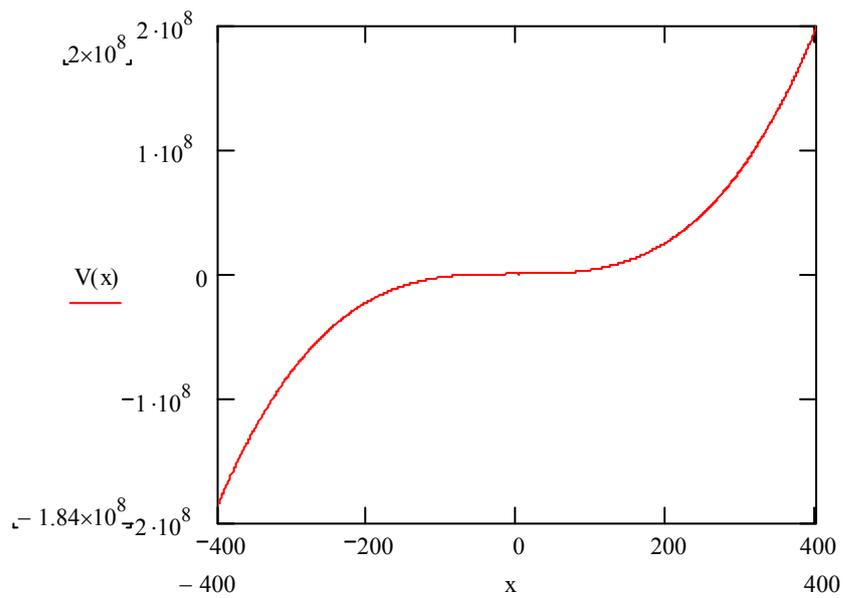

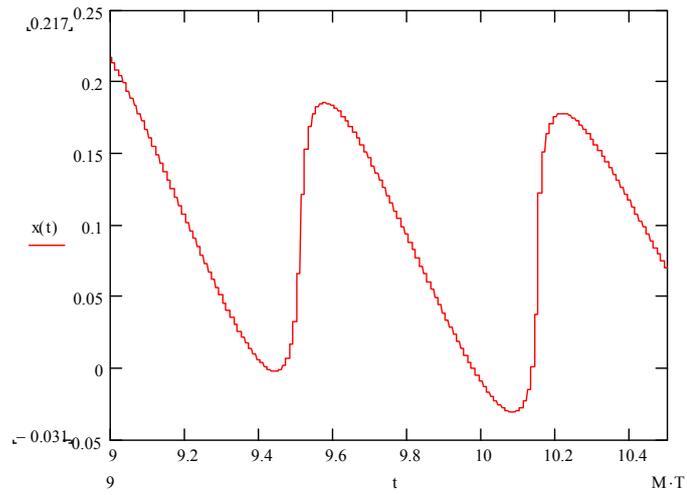

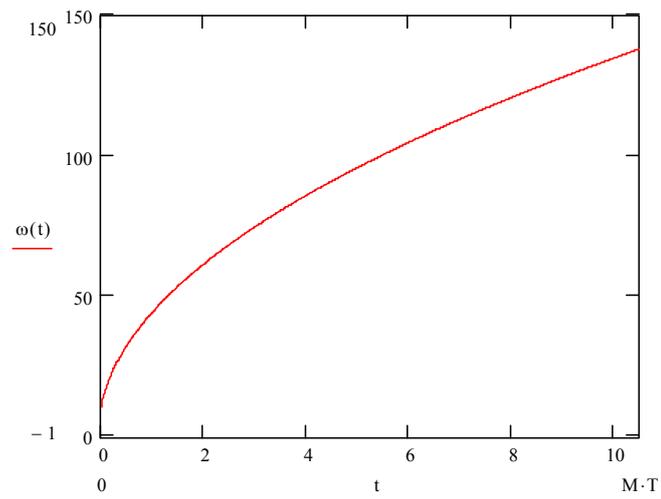

**Example**.2.5.6.

$m = 1, \Omega = 12, \omega = 9, a = 1, b = 15, A = 2.$

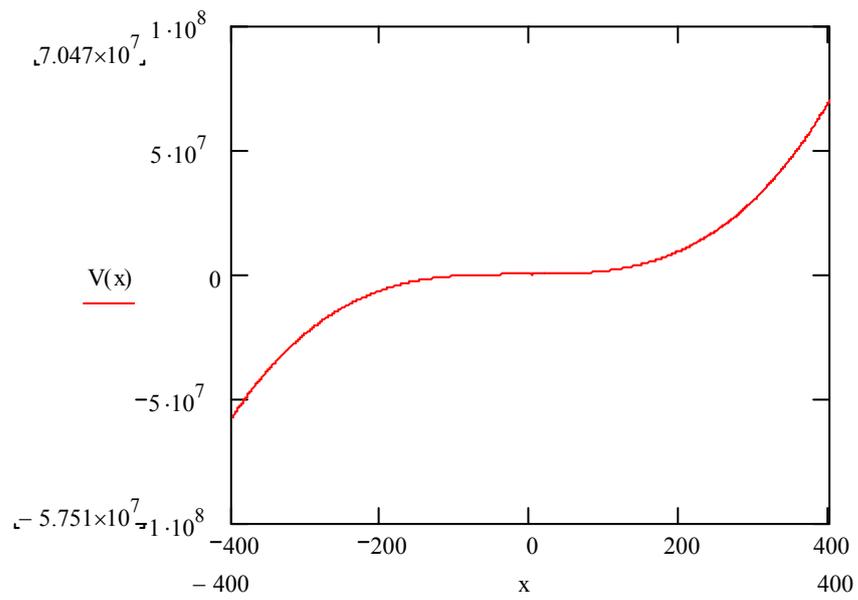

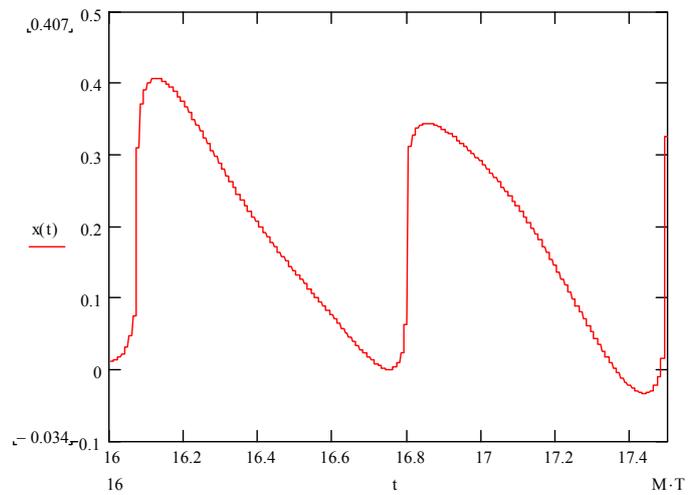

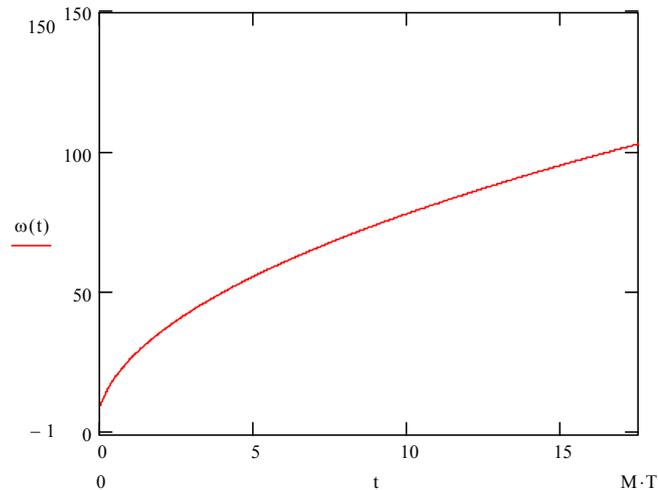

# III.Comparizon with a perturbation theory.

## III.1.Path-integrals calculation by using stationary-point approximation.

Let us consider classical dynamics of the form

$$\dot{x}(t,x_0) = F(x(t,x_0),t), x(0) = x_0,$$

$$F(x,t) = -V'(x) + A\sin(\Omega t),$$
(3.1.1)

$$\Omega = \frac{2\pi}{T}.$$

Where $V(x)$ is an static potential. Let us consider now corresponding quantum dynamics. Note that the transition probability amplitude $K(x,t|x_0,0) = K(x,t)$ is

governed by the Schrödinger type equation:

$$-i\epsilon \, \frac{\partial}{\partial t} \, K(x,t) \,=\, \frac{\partial}{\partial x} \left\{ \left[ -F(x,t) + i\epsilon^2 \frac{\partial}{\partial x} \right] K(x,t) \right\},$$

$$(3.1.3)$$

$$K(x,0) \,=\, \delta(x - x_0).$$

Let us define the conditional averadge $\langle x^2(t) \rangle$ via formula $\langle x^2(t) \rangle \triangleq \int_{-\infty}^{\infty} x_f^2 K(x_f,t) dx_f$. Hence for the conditional averadge $\langle x^2(t) \rangle$ we have the canonical path integral representation:

$$\langle x^2(t) \rangle \,=\, \int_{-\infty}^{\infty} x_f^2 K(x_f,t) dx_f \,=$$

$$\int_{-\infty}^{\infty} \int_{x(t_0)=x_0}^{x(t)=x_f} x_f^2 \exp\left[ -\frac{i}{\epsilon} S[\dot{x}(t), x(t)] \right] dx_f [Dx(t)],$$

$$(3.1.4)$$

$$S[\dot{x}(t), x(t)] \,=\, \frac{1}{4} \int_0^t [\dot{x}(t') - F(x(t'), t')] dt'.$$

Note that for propagator $K_N(x_0, t_0 | x_f, t_f)$ the time discretized path-integral representation is :

$$K_N(x_f, t_f | x_0, t_0) = \underbrace{\int \ldots \int}_{N-1} \exp\left\{ -\frac{iS_N(x_0, \ldots, x_N)}{\epsilon} \right\} \frac{dx_1 \cdots dx_{N-1}}{(4\pi\epsilon \, \Delta t)^{N/2}},$$

$$(3.1.5)$$

with

$$S_N(x_0, \ldots, x_N) \,=\, \sum_{n=0}^{N-1} \frac{\Delta t}{4} \left[ \frac{x_{n+1} - x_n}{\Delta t} - F(x_n, t_n) \right]^2 \,=$$

$$(3.1.6)$$

where the initial-$x_0$ and end-points $x_N$ are fixed by the prescribed $x_0$ and by the additional constraint $x_N = x_f$. From Eqs.(3.1.4)-(3.1.6) we obtain

$$|\langle x^2(t)\rangle| = \lim_{\Delta t \to 0} |\langle x_N^2(t)\rangle| =$$

$$\left| \int_{-\infty}^{\infty} x_N^2 K_N(x_0, t_0 | x_N, t) dx_N \right| =$$

(3.1.7)

$$\left| \underbrace{\int \cdots \int}_{N} x_N^2 \exp\left\{ -\frac{iS_N(x_0, \ldots, x_N)}{\epsilon} \right\} \frac{dx_1 \cdots dx_{N-1} dx_N}{(4\pi\epsilon\Delta t)^{N/2}} \right|.$$

Denoting an critical point of the discrete-time action $S_N(x_0, \ldots, x_N)$ by $\mathbf{x}_k^\# = (x_{0,k}, x_{1,k}^\#, \ldots, x_{n,k}^\#, \ldots, x_{N,k}^\#)$ it follows that $x^\#$ satisfies the critical point conditions

$$\frac{\partial S_N(\mathbf{x}_k^\#)}{\partial x_{n,k}^\#} = 0,$$

(3.1.8)

for $n = 1, \ldots, N$, supplemented by the prescribed boundary conditions for $n = 0, n = N$:

$$x_{0,k}^\# = x_0, x_{N,k}^\# = x_f$$

(3.1.9)

From Eq.(3.1.7) in the limit $\epsilon \to 0$ we obtain

$$\langle x_N^2(t)\rangle =$$

$$\int_{-\infty}^{\infty} x_N^2 K_N(x_0, t_0 | x_N, t) \, dx_N =$$

$$\underbrace{\int \ldots \int}_{N} x_N^2 \exp\left\{-\frac{iS_N(x_0,\ldots,x_N)}{\epsilon}\right\} \frac{dx_1 \cdots dx_{N-1} \, dx_N}{(4\pi\epsilon\Delta t)^{N/2}} = \qquad (3.1.10)$$

$$\mathbf{Z}_N(\mathbf{x}^\#) x_N^{\#2} \exp\left[-\frac{i}{\epsilon} S(\mathbf{x}^\#)\right] + o(\epsilon),$$

where the prefactor $\mathbf{Z}_N(\mathbf{x}^\#)$ is given via $N$-dimensional Gaussian integral of the canonical form as

$$\mathbf{Z}_N(\mathbf{x}^\#) = \underbrace{\int \ldots \int}_{N} \exp\left\{-\frac{i}{2\epsilon} \sum_{n,m=1}^{N} y_n \frac{\partial^2 S(\mathbf{x}^\#)}{\partial x_n^\# \partial x_m^\#} y_m\right\} \frac{dy_1 \cdots dy_{N-1} \, dy_N}{(4\pi D \Delta t)^{N/2}}. \qquad (3.1.11)$$

The Gaussian integral in (3.1.11) is given via formula

$$Z_N(\mathbf{x}^*) = \left[2\Delta t \det\left(2\Delta t \frac{\partial^2 S(\mathbf{x}^\#)}{\partial x_n^\# \partial x_m^\#}\right)\right]^{-\frac{1}{2}} \qquad (3.1.12)$$

$$n = 1,\ldots,N, m = 1,\ldots,N,$$

Equation (3.1.12) we rewrite in the following form:

$$Z_N(\mathbf{x}^*) = [2D Q_N^\#]^{-1/2}. \qquad (3.1.13)$$

Quantity $Q_N^\#$ in (3.1.13) can be calculated by using a second order linear recursion procedure given via formulae:

$$\frac{Q_{n+1}^* - 2Q_n^* - Q_{n-1}^*}{\Delta t^2} = \frac{Q_{n+1}^* - 2Q_n^* - Q_{n-1}^*}{\Delta t^2} -$$

$$Q_n^* \left[ \frac{x_{n+1}^* - x_n^*}{\Delta t} - F(x_n^*, t_n) \right] F''(x_n^*, t_n) +$$

$$Q_n^* F'(x_n^*, t_n)^2 - Q_{n-1}^* F'(x_{n-1}^*, t_{n-1})^2, \tag{3.1.14}$$

$$Q_1^* = \Delta t,$$

$$\frac{Q_2^* - Q_1^*}{\Delta t} = 1 + o(\Delta t).$$

As well known shall see later, we have to leave room for the possibility that even for small noise-strengths $\epsilon \to 0$ more than one critical point of the action (3.1.6) notably contributes to the path-integral expression (3.1.7). We label those various critical points $x_k^\#$ by the discrete index $k$. Thus from Eq.(3.1.10) and Eq.(3.1.13) we obtain

$$|\langle x_N^2(t_N, \omega) \rangle| = \left| \int_{-\infty}^{\infty} x_N^2 K_N(x_0, t_0 | x_N, t) dx_N \right| =$$

$$= \left| \sum_{k=0}^{M} \frac{(x_{N,k}^\#)^2 \exp\left[ \frac{iS_N(\mathbf{x}_k^\#)}{\epsilon} \right]}{\sqrt{2Q_{N,k}^\#}} [1 + o(\epsilon))] \right|. \tag{3.1.15}$$

The continuous-time limit of the action $S_N$ is

$$S[\dot{x}(t), x(t)] = \frac{1}{4} \int_{t_0}^{t_f} [\dot{x}^2 - F(x, t)] dt \tag{3.1.16}$$

The extremality conditions for the critical paths $x_k^\#(t)$ in the continuous-time limit are obtained from (3.1.2) and (3.1.4) by letting $\Delta t \to 0$ is:

(I)

$$\dot{x}^{\#}(t,x_0) = F(x^{\#}(t,x_0),t), x^{\#}(0) = x_0. \qquad (3.1.17)$$

**(II)**

$$\ddot{x}_k^{\#}(t,x_k^{\#}(t_f)) = \dot{F}(x_k^{\#}(t,x_k^{\#}(t_f)),t) + F(x_k^{\#}(t,x_k^{\#}(t_f)),t) F'(x_k^{\#}(t,x_k^{\#}(t_f)),t),$$

$$x^{\#}(0,x_0,\widetilde{x}_f) = x_0, x^{\#}(t_f,x_0,\widetilde{x}_f) = \widetilde{x}_f,$$

$$\left. \frac{dS[\dot{x}_k^{\#}(t,x_0,x_k^{\#}(t_f)),x_k^{\#}(t,x_0,x_k^{\#}(t_f))]}{d(x_k^{\#}(t_f))} \right|_{x_k^{\#}(t_f)=\widetilde{x}_{f,k}} = 0.$$

$$(3.1.18)$$

Here we use the canonical definitions $F'(x,t) = \dfrac{\partial F(x,t)}{\partial x}, \dot{F}(x,t) = \dfrac{\partial F(x,t)}{\partial t}$.
Hamiltonian counterpart of the Lagrangian dynamics given by (3.1.18) is

$$H(x,p,t) = p\dot{x} - L = p^2 + p F(x,t),$$

$$\dot{p}_k^{\#}(t) = -p_k^{\#}(t) F'(x_k^{\#}(t),t), \qquad (3.1.19)$$

$$\dot{x}_k^{\#}(t) = 2p_k^{\#}(t) + F(x_k^{\#}(t),t).$$

From the last equation of the Eqs.(3.1.19) one obtain the momentum $p_k^{\#}(t)$

$$p_k^{\#}(t) = \frac{1}{2}[\dot{x}_k^{\#}(t) - F(x_k^{\#}(t),t)]. \qquad (3.1.20)$$

From Eq.(3.1.16) and Eq.(3.1.20) one obtain:

$$\Phi_k(x_f,t_f) = S[x_k^{\#}(t)] = \int_{t_0}^{t_f} p_k^{\#}(t)^2 \, dt. \qquad (3.1.20)$$

In the limit $\Delta t \to 0$ from the extremality conditions (3.1.17) for the critical paths $x^{\#}(t,x_0,\widetilde{x}_f)$ and a second order linear recursion procedure given via formulae (3.1.10) we obtain the second order homogeneous linear differential equation

$$\frac{1}{2}\ddot{Q}_k^\#(t,x_0,\widetilde{x}_f) - \frac{d}{dt}\left[Q_k^\#(t,\widetilde{x}_f)\,F'(x^\#(t,x_0,\widetilde{x}_f),t)\right] \; +$$

$$Q_k^\#(t)[x^\#(t,x_0,\widetilde{x}_f)] \tag{3.1.1}$$

$$Q_0^\#(t_0,x_0,\widetilde{x}_f) = 0 \;, \quad \dot{Q}_0^\#(t_0,x_0,\widetilde{x}_f) = 1.$$

|

# IV. Wave function collapse and weak Colombeau solutions of the Colombeau-Schrödinger equation.

In quantum mechanics, wave function collapse is the phenomenon in which a wave function-initially in a superposition of several eigenstates-appears to reduce to a single eigenstate after interaction with a measuring apparatus.[57] It is the essence of measurement in quantum mechanics, and connects the wave function with classical observables like position and momentum. Collapse is one of two processes by which quantum systems evolve in time; the other is continuous evolution via the Schrödinger equation [58] However in this role, collapse is merely a black box for thermodynamically irreversible interaction with a classical environment.[59] Calculations of quantum decoherence predict apparent wave function collapse when a superposition forms between the quantum system's states and the environment's states. Significantly, the combined wave function of the system and environment continue to obey the Schrödinger equation [60].

When the Copenhagen interpretation was first expressed, Niels Bohr postulated wave function collapse to cut the quantum world from the classical.[60] This tactical move allowed quantum theory to develop without distractions from interpretational worries. Nevertheless it was debated, for if collapse were a fundamental physical phenomenon, rather than just the epiphenomenon of some other process, it would mean nature was fundamentally stochastic, i.e. nondeterministic, an undesirable property for a theory [61] This issue remained until quantum decoherence entered mainstream opinion after its reformulation in the 1980s.[62] Decoherence explains the perception of wave function collapse in terms of interacting large- and small-scale quantum systems, and is commonly taught at the graduate level (e.g. the Cohen-Tannoudji textbook).[63] The quantum filtering approach [64],[65] and the introduction of quantum causality non-demolition principle [66] allows for a classical-environment derivation of wave function collapse from the stochastic Schrödinger equation.

# IV.1. Weak Colombeau solutions of the Colombeau-Schrödinger equation.

Let us consider Colombeau-Schrödinger equation:

$$(i\hbar_\varepsilon \partial \Psi_\varepsilon(\mathbf{x},t)/\partial t)_\varepsilon = \left(\widehat{\mathbf{H}}_\varepsilon \Psi_\varepsilon(\mathbf{x},t), \right)_\varepsilon, (\Psi_\varepsilon(\mathbf{x},t_0))_\varepsilon = (\Psi_\varepsilon(\mathbf{x}))_\varepsilon,$$

$$\mathbf{x} \in \mathbb{R}^d, \mathbf{x} = (x_1, x_2, \ldots, x_d), (\hbar_\varepsilon) \approx 0,$$

where $[(\Psi_\varepsilon^{\#}(\mathbf{x}))_\varepsilon] \in G(\mathbb{R}^d)$.

**Definition 4.1.1.** Assume that: **(i)** $[(\Psi_\varepsilon^{\#}(\mathbf{x}))_\varepsilon] \in G(\mathbb{R}^d)$, **(ii)** $[(\Psi_\varepsilon^{\#}(\mathbf{x},t))_\varepsilon] \in G(\mathbb{R}^d \times \mathbb{R}_+)$.

# IV.2. Gamow theory of the alpha decay via tunneling. Weak Colombeau solutions of the Colombeau-Schrödinger equation coresponding to alpha decay via tunneling.

By 1928, George Gamow had solved the theory of the alpha decay via tunneling []. The alpha particle is trapped in a potential well by the nucleus. Classically, it is forbidden to escape, but according to the (then) newly discovered principles of quantum mechanics, it has a tiny (but non-zero) probability of "tunneling" through the barrier and appearing on the other side to escape the nucleus. Gamow solved a model potential for the nucleus and derived, from first principles, a relationship between the half-life of the decay, and the energy of the emission.

The particle has total energy $E$ and is incident on the barrier from the left.

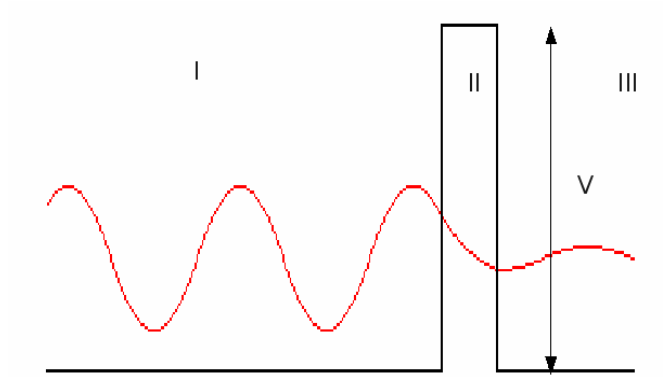

Pic.4.2.1.The particle has total energy $E$ and
is incident on the barrier $V(x)$ from the left.

The Schrodinger equation in each of regions $\mathbf{I} = \{x | x < 0\}, \mathbf{II} = \{x | 0 \leq x \leq l\}$ and
$\mathbf{III} = \{x | x > l\}$ takes the folloving form

$$\left( \frac{\partial^2 \Psi_\varepsilon(x)}{\partial x^2} + \frac{2m}{h_\varepsilon^2}[E - U(x)]\Psi_\varepsilon(x) \right)_\varepsilon = 0, \qquad (\mathbf{4}.2.1)$$

where

$$U(x) = \begin{cases} 0, \text{при } x < 0 \\ U_0, \text{при } 0 \leq x \leq l \\ 0, \text{при } x > l \end{cases} \qquad (\mathbf{4}.2.2)$$

The solutions is

$$(\Psi_{\varepsilon,\mathbf{I}}(x))_\varepsilon = (C_{\varepsilon,+} \exp(ik_\varepsilon x))_\varepsilon + (C_{\varepsilon,-} \exp(-ik_\varepsilon x))_\varepsilon$$

$$(\Psi_{\varepsilon,\mathbf{II}}(x))_\varepsilon = (B_{\varepsilon,+} \exp(k'_\varepsilon x))_\varepsilon + (B_{\varepsilon,-} \exp(-k'_\varepsilon x))_\varepsilon$$

$$(\Psi_{\varepsilon,\mathbf{III}}(x))_\varepsilon = (A_\varepsilon \exp(ik_\varepsilon x))_\varepsilon, \qquad (.93)$$

where

$$(k_\varepsilon)_\varepsilon = \frac{2\pi}{(\hbar_\varepsilon)_\varepsilon} \sqrt{2mE},$$

$$(VI.94)$$

$$(k'_\varepsilon)_\varepsilon = \frac{2\pi}{(\hbar_\varepsilon)_\varepsilon} \sqrt{2m(U_0 - E)}.$$

# IV. 3. Quantum trajectories induced from a strong Colombeau solutions of the Colombeau-Schrödinger equation.

Let us consider Colombeau-Schrödinger equation:

$$\left(i\hbar_\varepsilon \partial \Psi_\varepsilon(\mathbf{x},t)/\partial t\right)_\varepsilon = \left(\widehat{\mathbf{H}}_\varepsilon \Psi_\varepsilon(\mathbf{x},t),\right)_\varepsilon,$$

$$\left(\Psi_\varepsilon(\mathbf{x},t_0)\right)_\varepsilon = \left(\Psi_\varepsilon(\mathbf{x})\right)_\varepsilon \in G(\mathbb{R}^d), \qquad (4.3.1)$$

$$\mathbf{x} \in \mathbb{R}^d, \mathbf{x} = (x_1, x_2, \ldots, x_d), (\hbar_\varepsilon)_\varepsilon \approx 0,$$

where initial data $\left[\left(\Psi_\varepsilon(\mathbf{x})\right)_\varepsilon\right]$ defined by formula

$$(\Psi_\varepsilon(\mathbf{x}))_\varepsilon = (\Psi_\varepsilon(\mathbf{x}, \mathbf{x}_0))_\varepsilon =$$

$$= \left( \frac{\mu(\mathbf{x}_0)}{\sqrt[4]{2\pi\eta_\varepsilon}} \Psi_\varepsilon(\mathbf{x}) \exp\left[ \frac{\langle \mathbf{x} - \mathbf{x}_0, \mathbf{x} - \mathbf{x}_0 \rangle}{2\eta_\varepsilon} \right] \right)_\varepsilon =$$

$$= \left( \frac{1}{\left(\sqrt[4]{2\pi\eta_\varepsilon}\right)^d} \exp\left[ \frac{\langle \mathbf{x} - \mathbf{x}_0, \mathbf{x} - \mathbf{x}_0 \rangle}{2\eta_\varepsilon} \right] \exp\left[ \frac{i}{\hbar_\varepsilon} S(\mathbf{x}) \right] \right)_\varepsilon, \qquad (4.3.2)$$

$$\text{where}$$

$$(\eta_\varepsilon)_\varepsilon \approx 0 \text{ and } \left( \frac{\hbar_\varepsilon}{\eta_\varepsilon} \right)_\varepsilon \approx 0,$$

**Definition 4.3.1.** Quantum trajectory $\{(x_{\varepsilon,i}(t, \mathbf{x}_0; \hbar))_\varepsilon\}_{i=1}^d$ with initial data $\{(x_{\varepsilon,i}(0, \mathbf{x}_0; \hbar))_\varepsilon\}_{i=1}^d = \mathbf{x}_0$ is defined by following formula

$$(x_{\varepsilon,i}(t, \mathbf{x}_0; \hbar))_\varepsilon = \left( \int_{\mathbb{R}^d} x_i |\Psi_\varepsilon(\mathbf{x}, t, \mathbf{x}_0; \hbar)|^2 d^d x \right)_\varepsilon. \qquad (4.3.3)$$

# V. Schrödinger's cat paradox resolution using GRW collapse model.

In this subsection an possible solution of the Schrödinger's cat paradox is considered. We pointed out that: the collapsed state of the cat always shows definite and predictable measurement outcomes even if Schrödinger's cat consists of a superposition: $|\text{cat}\rangle = c_1 |\text{live cat}\rangle + c_2 |\text{death cat}\rangle$.

Contrary to van Kampen's [61] and some others' opinions, "looking" at the outcome changes nothing, beyond informing the observer of what has already happened. We remain: there are widespread claims that Schrödinger's cat is not in a definite alive or dead state but is, instead, in a superposition of the two. van Kampen, for example, writes "The whole system is in a superposition of two states: one in which no decay has occurred and. . .one in which it has occurred. Hence, the state of the cat also consists of a superposition: $|\text{cat}\rangle = c_1|\text{live cat}\rangle + c_2|\text{death cat}\rangle$. The state remains a superposition until an observer looks at the cat" [61].The canonical formulation [62]:

$$|\text{cat}\rangle = c_1|\text{live cat}\rangle|\text{undecayed nucleus }\rangle + c_2|\text{death cat}\rangle|\text{decayed nucleus }\rangle$$

completely obscures the unitary Schrödinger evolution which by using GRW collapse model, predicts specific nonlocal entanglement. The cat state must be

written as:

$$|\text{cat}(t)\rangle = c_1|\text{live cat}(t)\rangle|\text{undecayed nucleus }(t)\rangle +$$
$$+c_2|\text{death cat}(t)\rangle|\text{decayed nucleus}(t)\rangle.$$

This entangled state actually is the collapsed state of both the cat and the nucleus, showing definite outcomes at each instant $t \geq T_{\text{col}}$.

# V.I.Schrödinger's cat paradox

As Weinberg recently reminded us [63], the measurement problem remains a fundamental conundrum. During measurement the state vector of the microscopic system collapses in a probabilistic way to one of a number of classical states, in a way that is unexplained, and cannot be described by the time-dependent Schrödinger equation.To review the essentials, it is sufficient to consider two-state systems. Suppose a nucleus $\mathbf{n}$, whose Hilbert space is spanned by orthonormal states $|s_i(t)\rangle$, $i = 1,2$, where $|s_1(t)\rangle = |\text{undecayed nucleus at instant } t\rangle$ and $|s_2(t)\rangle = |\text{decayed nucleus at instant } t\rangle$ is in the superposition state,

$$|\Psi_t\rangle_{\mathbf{n}} = c_1|s_1(t)\rangle + c_2|s_2(t)\rangle, |c_1|^2 + |c_2|^2 = 1. \qquad (5.1.1)$$

An measurement apparatus $A$, which may be microscopic or macroscopic, is designed to distinguish between states $|s_i(t)\rangle$ by transitioning at each instant $t$ into state $|a_i(t)\rangle$ if it finds $\mathbf{n}$ is in $|s_i(t)\rangle$, $i = 1,2$. Assume the detector is reliable, implying the $|a_1(t)\rangle$ and $|a_2(t)\rangle$ are orthonormal at each instant $t$-i.e., $\langle a_1(t)\|a_2(t)\rangle = 0$ and that the measurement interaction does not disturb states $|s_i\rangle$ -i.e., the measurement is "ideal". When $A$ measures $|\Psi_t\rangle_{\mathbf{n}}$, the Schrödinger equation's unitary time evolution then leads to the "measurement state" (MS) $|\Psi_t\rangle_{\mathbf{n}A}$ :

$$|\Psi_t\rangle_{\mathbf{n}A} = c_1|s_1(t)\rangle|a_1(t)\rangle + c_2|s_2(t)\rangle|a_2(t)\rangle, |c_1|^2 + |c_2|^2 = 1. \qquad (5.1.2)$$

of the composite system $\mathbf{n}A$ following the measurement.

Standard formalism of continuous quantum measurements [62]-[68] leads to a definite but unpredictable measurement outcome, either $|a_1(t)\rangle$ or $|a_2(t)\rangle$ and that $|\Psi_t\rangle_{\mathbf{n}}$ suddenly "collapses" at instant $t'$ into the corresponding state $|s_i(t')\rangle$. But unfortunately equation (1.2) does not appear to resemble such a collapsed state at instant $t'$?.

The measurement problem is as follows [62]-[68]:

(**I**) How do we reconcile canonical collapse models postulate's

(**II**) How do we reconcile the measurement postulate's definite outcomes with the "measurement state" $|\Psi_t\rangle_{\mathbf{n},A}$ at each instant $t$ and

(**III**) how does the outcome become irreversibly recorded in light of the Schrödinger equation's unitary and, hence, reversible evolution?

This paper deals with only the special case of the measurement problem, known as Schrödinger's cat paradox. For a good and complete explanation of this paradox see

Hobson [62] and Leggett [68].

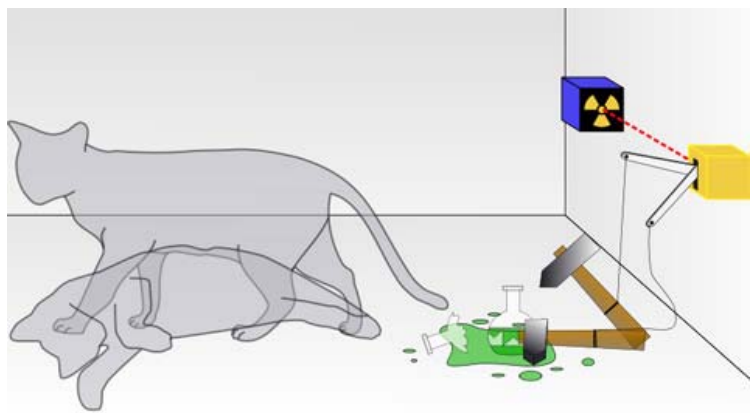

Pic.5.1.1.Schrödinger's cat.

**Schrödinger's cat**: a cat, a flask of poison, and a radioactive source are placed in a sealed box. If an internal monitor detects radioactivity (i.e. a single atom decaying), the flask is shattered, releasing the poison that kills the cat. The Copenhagen interpretation of quantum mechanics implies that after a while, the cat is simultaneously alive and dead. Yet, when one looks in the box, one sees the cat either alive or dead, not both alive and dead. This poses the question of when exactly quantum superposition ends and reality collapses into one possibility or the other.This subsection presents an theoretical approach of the MS that resolves the problem of definite outcomes for Schrödinger's "cat". It shows that the MS actually is the collapsed state of both Schrödinger's "cat" and nucleus, even though it evolved purely unitarily.

**The canonical collapse models**.

In order to appreciate how canonical collapse models work, and what they are able to achieve, we briefly review the GRW model. Let us consider a system of $n$ particles which, only for the sake of simplicity, we take to be scalar and spinless; the GRW model is defined by the following postulates: (**1**) The state of the system is represented by a wave function $\psi_t(\mathbf{x}_1, \mathbf{x}_2, \ldots, \mathbf{x}_n)$ belonging to the Hilbert space $\mathcal{L}_2(\mathbb{R}^{3n})$. **(2)** At random times, the wave function experiences a sudden jump of the form: $\psi_t(\mathbf{x}_1, \mathbf{x}_2, \ldots, \mathbf{x}_n) \rightarrow \psi_t(\mathbf{x}_1, \mathbf{x}_2, \ldots, \mathbf{x}_n; \widetilde{\mathbf{x}}_m)$, where

$$\psi_t(\mathbf{x}_1, \mathbf{x}_2, \ldots, \mathbf{x}_n; \widetilde{\mathbf{x}}_m) = \frac{\mathfrak{R}_m(\widetilde{\mathbf{x}}_m)\psi_t(\mathbf{x}_1, \mathbf{x}_2, \ldots, \mathbf{x}_n)}{\| \mathfrak{R}_m(\widetilde{\mathbf{x}}_m)\psi_t(\mathbf{x}_1, \mathbf{x}_2, \ldots, \mathbf{x}_n) \|_2}, \tag{5.1.3}$$

where $\psi_t(\mathbf{x}_1, \mathbf{x}_2, \ldots, \mathbf{x}_n)$ is the state vector of the whole system at time $t$, immediately prior to the jump process and $\mathfrak{R}_n(\widetilde{\mathbf{x}}_m)$ is a linear operator which is conventionally chosen equal to:

$$\mathfrak{R}_m(\widetilde{\mathbf{x}}_m) = (\pi r_c^2)^{-3/4} \exp\left[ -\frac{(\widehat{\mathbf{x}}_m - \widetilde{\mathbf{x}}_m)^2}{2r_c^2} \right], \tag{5.1.4}$$

where $r_c$ is a new parameter of the model which sets the width of the localization process, and $\widehat{\mathbf{x}}_m$ is the position operator associated to the $m$-th particle of the system and the random variable $\widetilde{\mathbf{x}}_m$ corresponds to the place where the jump occurs. **(3)** It is assumed that the jumps are distributed in time like a Poissonian process with frequency $\lambda = \lambda_{GRW}$ this is the second new parameter of the model. (4) Between two consecutive jumps, the state vector evolves according to the standard Schrödinger equation.

The 1-particle master equation of the GRW model takes the form

$$\frac{d}{dt}\rho(t) = -\frac{i}{\hbar}\left[ \widehat{\mathbf{H}}, \rho(t) \right] - T[\rho(t)]. \tag{5.1.5}$$

Here $\widehat{\mathbf{H}}$ is the standard quantum Hamiltonian of the particle, and $T[\cdot]$ represents the effect of the spontaneous collapses on the particle's wave function. In the position representation, this operator becomes:

$$\langle \mathbf{x}|T[\rho(t)]|\mathbf{y}\rangle = \lambda\left\{ 1 - \exp\left[ -\frac{(\mathbf{x} - \mathbf{y})^2}{4r_c^2} \right] \right\}\langle \mathbf{x}|\rho(t)|\mathbf{y}\rangle. \tag{5.1.6}$$

**Remark 1**.**1**. We note that GRW collapse model follows from the more general

S. Weinberg formalism [1].

Another modern approach to stochastic reduction is to describe it using a stochastic nonlinear Schrödinger equation, an elegant simple example of which is the following one particle case known as Quantum Mechanics with Universal Position Localization [QMUPL]:

$$d|\psi_t(x)\rangle = \left[ -\frac{i}{\hbar}\widehat{\mathbf{H}} - k(\widehat{x} - \langle x_t\rangle)^2 dt \right]|\psi_t(x)\rangle dt + \sqrt{2k}\,(\widehat{q} - \langle q_t\rangle)dW_t|\psi_t(x)\rangle. \tag{5.1.7}$$

Here $\hat{x}$ is the position operator, $\langle x_t \rangle = \langle \psi_t | \hat{x} | \psi_t \rangle$ it is its expectation value, and $k$ is a constant, characteristic of the model, which sets the strength of the collapse mechanics, and it is chosen proportional to the mass $m$ of the particle according to the formula: $k = (m/m_0)\lambda_0$, where $m_0$ is the nucleon's mass and $\lambda_0$ measures the collapse strength. It is easy to see that Eq.(5.1.5) contains both non-linear and stochastic terms, which are necessary to induce the collapse of the wave function.

For an example let us consider a free particle $(\widehat{\mathbf{H}} = p^2/2m)$, and a Gaussian state:

$$\psi_t(x) = \exp\left\{-a_t(x - \bar{x}_t)^2 + i\bar{k}_t x\right\}. (5.1.8)$$

It is easy to see that $\psi_t(x)$ given by Eq.(5.1.6) is solution of Eq.(5.1.5), where

$$\frac{da_t}{dt} = k - \frac{2i\hbar}{m}a_t^2, \frac{d\bar{x}_t}{dt} = \frac{\hbar}{m}\bar{k}_t + \frac{\sqrt{k}}{2\operatorname{Re}(a_t)}\dot{W}_t, \frac{d\bar{k}_t}{dt} = -\sqrt{k}\frac{\operatorname{Im}(a_t)}{\operatorname{Re}(a_t)}\dot{W}_t. (5.1.9)$$

The CSL model is defined by the following stochastic differential equation in the Fock space:

$$d|\psi_t(\mathbf{x})\rangle = \left[-\frac{i}{\hbar}\widehat{\mathbf{H}} - k\left(\widehat{M}(\mathbf{x}) - \langle M_t(\mathbf{x})\rangle\right)^2 dt\right]|\psi_t(\mathbf{x})\rangle dt +$$

$$+ \sqrt{2k}\left(\widehat{M}(\mathbf{x}) - \langle M_t(\mathbf{x})\rangle\right)dW_t(\mathbf{x})|\psi_t(\mathbf{x})\rangle. (5.1.10)$$

# V.II. Generalized Gamow theory of the alpha decay via tunneling using GRW collapse model.

By 1928, George Gamow had solved the theory of the alpha decay via tunneling [69]. The alpha particle is trapped in a potential well by the nucleus. Classically, it is forbidden to escape, but according to the (then) newly discovered principles of quantum mechanics, it has a tiny (but non-zero) probability of "tunneling" through the barrier and appearing on the other side to escape the nucleus. Gamow solved a model potential for the nucleus and derived, from first principles, a relationship between the half-life of the decay, and the energy of the emission. The $\alpha$-particle has total energy $E$ and is incident on the barrier from the right to left.

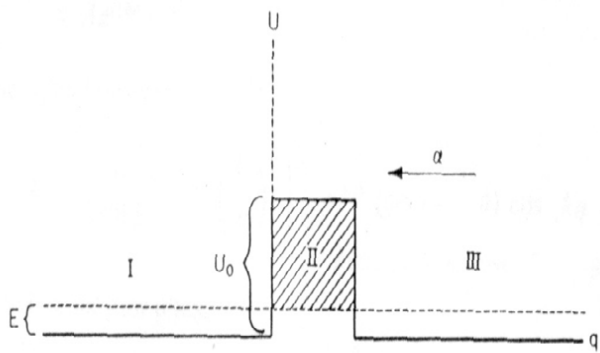

Pic. 5.2.1. The particle has total energy $E$ and
is incident on the barrier $V(x)$ from right to left.

The Schrödinger equation in each of regions $\mathbf{I} = \{x | x < 0\}, \mathbf{II} = \{x | 0 \leq x \leq l\}$ and $\mathbf{III} = \{x | x > l\}$ takes the following form

$$\frac{\partial^2 \Psi(x)}{\partial x^2} + \frac{2m}{\hbar^2}[E - U(x)]\Psi(x) = 0, \tag{5.2.1}$$

where

$$U(x) = \begin{cases} 0 \text{ for } x < 0 \\ U_0 \text{ for } 0 \leq x \leq l \\ 0 \text{ for } x > l \end{cases} \tag{5.2.2}$$

The solutions reads [8]:

$$\Psi_{\mathbf{III}}(x) = C_+ \exp(ikx) + C_- \exp(-ikx),$$
$$\Psi_{\mathbf{II}}(x) = B_+ \exp(k'x) + B_- \exp(-k'x),$$
$$\Psi_{\mathbf{I}}(x) = A\cos(kx) = \frac{A}{2}[\exp(ikx) + \exp(-ikx)], \tag{5.2.3}$$

where

$$k = \frac{2\pi}{h}\sqrt{2mE}, k' = \frac{2\pi}{h}\sqrt{2m(U_0 - E)}. \tag{5.2.4}$$

At the boundary $x = 0$ we have the following boundary conditions:

$$\Psi_{\mathbf{I}}(0) = \Psi_{\mathbf{II}}(0), \frac{\partial \Psi_{\mathbf{I}}(x)}{\partial x}\bigg|_{x=0} = \frac{\partial \Psi_{\mathbf{II}}(x)}{\partial x}\bigg|_{x=0}. \tag{5.2.5}$$

At the boundary $x = l$ we have the following boundary conditions

$$\Psi_{\mathbf{II}}(l) = \Psi_{\mathbf{III}}(l), \frac{\partial \Psi_{\mathbf{II}}(x)}{\partial x}\bigg|_{x=l} = \frac{\partial \Psi_{\mathbf{III}}(x)}{\partial x}\bigg|_{x=l}. \tag{5.2.6}$$

From the boundary conditions (5.2.5)-(5.2.6) one obtains [69]:

$$B_+ = \frac{A}{2}\left(1 + i\frac{k}{k'}\right), B_- = \frac{A}{2}\left(1 - i\frac{k}{k'}\right),$$
$$C_+ = A[ch(k'l) + iDsh(k'l)], C_- = i(ASsh(k'l)\exp(ikl)),$$
$$D = \frac{1}{2}\left(\frac{k}{k'} - \frac{k'}{k}\right), S = \frac{1}{2}\left(\frac{k}{k'} + \frac{k'}{k}\right). \tag{5.2.7}$$

From (5.2.7) one obtain the conservation law

$$|A|^2 = |C_+|^2 - |C_-|^2.$$

Let us introduce now a function $E_{II}(x,l) = \theta_2(x,l)E_2(x,l)$ where

$$E_2(x,l) = \begin{cases} (\pi r_c^2)^{-1/4}\exp\left(-\dfrac{x^2}{2r_c^2}\right) \text{ for } -\infty < x < \dfrac{l}{2} \\[2mm] (\pi r_c^2)^{-1/4}\exp\left(-\dfrac{(x-l)^2}{2r_c^2}\right) \text{ for } \dfrac{l}{2} \le x < \infty \end{cases} \qquad (5.2.8)$$

$$\theta_2(x,l) = \begin{cases} 1 \text{ for } x \in [0,l] \\ 0 \text{ for } x \notin [0,l] \end{cases}$$

**Assumption 5.2.1.** We assume now that:
(i) at instant $t = 0$ the wave function $\Psi_I(x)$ experiences a sudden jump of the form

$$\Psi_I(x) \rightarrow \Psi_I^\#(x) = \frac{\mathfrak{R}_I(\hat{x})\Psi_I(x)}{\|\mathfrak{R}_I(\hat{x})\Psi_I(x)\|_2}, \qquad (5.2.9)$$

where $\mathfrak{R}_I(\hat{x})$ is a linear operator which is chosen equal to:

$$\mathfrak{R}_I(\hat{x}) = (\pi r_c^2)^{-1/4}\theta_1(\hat{x},l)\exp\left[-\frac{\hat{x}^2}{2r_c^2}\right]; \qquad (5.2.10)$$

where

$$\theta_1(x,l) = \begin{cases} 1 \text{ for } x \in [-l,0], \\ 0 \text{ for } x \notin [-l,0]. \end{cases}$$

**Remark 5.2.1.** Note that: $\text{supp}(\Psi_I^\#(x)) \subseteq [-l,0]$
(ii) at instant $t = 0$ the wave function $\Psi_{II}(x)$ experiences a sudden jump of the form

$$\Psi_{II}(x) \rightarrow \Psi_{II}^\#(x) = \frac{\mathfrak{R}_{II}(\hat{x})\Psi_{II}(x)}{\|\mathfrak{R}_{II}(\hat{x})\Psi_{II}(x)\|_2}, \qquad (5.2.11)$$

where $\mathfrak{R}_{II}(\hat{x})$ is a linear operator which is chosen equal to:

$$\mathfrak{R}_{II}(\hat{x}) = E_{II}(\hat{x},l); \qquad (5.2.12)$$

**Remark 5.2.2.** Note that: $\text{supp}(\Psi_{II}^\#(x)) \subseteq [0,l]$.
(iii) at instant $t = 0$ the wave function $\Psi_{III}(x)$ experiences a sudden jump of the form

$$\Psi_{\mathbf{III}}(x) \to \Psi_{\mathbf{III}}^{\#}(x) = \frac{\Re_{\mathbf{III}}(\hat{x})\Psi_{\mathbf{III}}(x)}{\|\Re_{\mathbf{III}}(\hat{x})\Psi_{\mathbf{III}}(x)\|_2}, \tag{5.2.13}$$

where $\Re_{\mathbf{III}}(\hat{x})$ is a linear operator which is chosen equal to:

$$\Re_{\mathbf{III}}(\hat{x}) = (\pi r_c^2)^{-1/4} \exp\left[ -\frac{(\hat{x}-l)^2}{2r_c^2} \right]. \tag{5.2.14}$$

**Remark 5.2.3.** Note that. We have choose operators (5.2.10),(5.2.12) and (5.2.14) such that the boundary conditions (5.2.5),(5.2.6) is satisfied.

**Definition 5.2.1.** Let $\Psi(x)$ be an solution of the Schrödinger equation (5.2.1). The stationary Schrödinger equation (5.2.1) is a weakly well preserved in region $\Gamma \subseteq \mathbb{R}$ by collapsed wave function $\Psi^{\#}(x)$ if there exist an wave function $\Psi(x)$ such that the estimate

$$\int_{\Gamma} \left\{ \frac{\partial^2 \Psi^{\#}(x)}{\partial x^2} + \frac{2m}{h^2}[E - U(x)]\Psi^{\#}(x) \right\} dx = O(h^{2+\alpha}), \tag{5.2.15}$$

where the inequality $\alpha \geq 1$, is satisfied.

**Proposition 5.2.1.** The Schrödinger equation in each of regions $\mathbf{I}, \mathbf{II}, \mathbf{III}$ is a weakly well preserved by collapsed wave function $\Psi_{\mathbf{I}}^{\#}(x), \Psi_{\mathbf{II}}^{\#}(x)$ and $\Psi_{\mathbf{III}}^{\#}(x)$ correspondingly.

**Proof.** See Appendix CAT-B.

**Definition 5.2.2.** Let us consider the time-dependent Schrödinger equation:

$$ih\frac{\partial \Psi(\mathbf{x},t)}{\partial t} = \widehat{\mathbf{H}}\Psi(\mathbf{x},t), \tag{5.2.16}$$

$$t \in [0,T], \mathbf{x} \in \mathbb{R}^{3n}.$$

The time-dependent Schrödinger equation (5.2.16) is a weakly well preserved by corresponding to $\Psi(\mathbf{x},t)$ collapsed wave function $\Psi^{\#}(\mathbf{x},t)$

$$\Psi^{\#}(\mathbf{x}_1, \mathbf{x}_2, \ldots, \mathbf{x}_n, t) =$$

$$\Psi(\mathbf{x}_1, \mathbf{x}_2, \ldots, \mathbf{x}_n, t; \widetilde{\mathbf{x}}_{m_1}, \ldots, \widetilde{\mathbf{x}}_{m_k}) =$$

$$= \frac{\Re_{m_1, \ldots, m_k}(\widetilde{\mathbf{x}}_{m_1}, \ldots, \widetilde{\mathbf{x}}_{m_k}) \Psi(\mathbf{x}_1, \mathbf{x}_2, \ldots, \mathbf{x}_n, t)}{\| \Re_{m_1, \ldots, m_k}(\widetilde{\mathbf{x}}_{m_1}, \ldots, \widetilde{\mathbf{x}}_{m_k}) \Psi(\mathbf{x}_1, \mathbf{x}_2, \ldots, \mathbf{x}_n, t) \|_2},$$

$$\Re_{m_1, \ldots, m_k}(\widetilde{\mathbf{x}}_{m_1}, \ldots, \widetilde{\mathbf{x}}_{m_k}) = \prod_{i=1}^{k} \Re_{m_i}(\widetilde{\mathbf{x}}_{m_i})$$

in region $\Gamma \subseteq \mathbb{R}^{3d}$ if there exist an wave function $\Psi(\mathbf{x}, t)$ such that the estimate

$$\int_{\Gamma} \left\{ i\hbar \frac{\partial \Psi^{\#}(\mathbf{x}, t)}{\partial t} - \widehat{\mathbf{H}} \Psi^{\#}(\mathbf{x}, t) \right\} d^{3d}x = O(\hbar^{\alpha}), \tag{5.2.17}$$

$$t \in [0, T], \mathbf{x} \in \mathbb{R}^{3d},$$

where the inequality $\alpha \geq 1$, is satisfied.

**Definition 5.2.3.** Let $\Psi^{\#}(\mathbf{x}, t) = \Psi^{\#}(\mathbf{x}_1, \mathbf{x}_2, \ldots, \mathbf{x}_d, t)$ be a function $\Psi(\mathbf{x}_1, \mathbf{x}_2, \ldots, \mathbf{x}_d, t; \widetilde{\mathbf{x}}_1, \ldots, \widetilde{\mathbf{x}}_d)$. Let us consider the Probability Current Law

$$\frac{\partial}{\partial t} P(\Gamma, t) + \int_{\partial \Gamma} \mathbf{J}(\mathbf{x}_1, \mathbf{x}_2, \ldots, \mathbf{x}_d, t) \cdot \mathbf{n} d^{2d}x = O(\hbar^{\alpha}),$$

$$\mathbf{J}(\mathbf{x}_1, \mathbf{x}_2, \ldots, \mathbf{x}_d, t) = \Psi(\mathbf{x}, t)\nabla\overline{\Psi(\mathbf{x}, t)} - \overline{\Psi(\mathbf{x}, t)}\nabla\Psi(\mathbf{x}, t), \tag{5.2.18}$$

$$t \in [0, T], \mathbf{x} \in \mathbb{R}^{3d},$$

corresponding to Schrödinger equation (5.2.16). Probability Current Law (5.2.18) is a weakly well preserved by corresponding to $\Psi(\mathbf{x}, t)$ collapsed wave function $\Psi^{\#}(\mathbf{x}, t)$ in region $\Gamma \subseteq \mathbb{R}^{3d}$ if there exist an wave function $\Psi(\mathbf{x}, t)$ such that the estimate

$$\frac{\partial}{\partial t} P(\Gamma, t) + \int_{\partial \Gamma} \mathbf{J}^{\#}(\mathbf{x}_1, \mathbf{x}_2, \ldots, \mathbf{x}_d, t) \cdot \mathbf{n} d^{2d}x = O(\hbar^{\alpha}),$$

$$\mathbf{J}^{\#}(\mathbf{x}_1, \mathbf{x}_2, \ldots, \mathbf{x}_d, t) = \Psi^{\#}(\mathbf{x}, t)\nabla\overline{\Psi^{\#}(\mathbf{x}, t)} - \overline{\Psi^{\#}(\mathbf{x}, t)}\nabla\Psi^{\#}(\mathbf{x}, t)$$

$$= O(\hbar^{\alpha}), \tag{5.2.19}$$

$$t \in [0, T], \mathbf{x} \in \mathbb{R}^{3d},$$

where the inequality $\alpha \geq 1$, is satisfied.

**Proposition 5.2.2.** Assume that there exist an wave function $\Psi(\mathbf{x}, t)$ such that the estimate (5.2.17) is satisfied. Then Probability Current Law (2.18) is a weakly well preserved by corresponding to $\Psi(\mathbf{x}, t)$ collapsed wave function $\Psi^{\#}(\mathbf{x}, t)$ in region $\Gamma \subseteq \mathbb{R}^{3d}$, i.e. the estimate (5.2.19) is satisfied on the wave function $\Psi^{\#}(\mathbf{x}, t)$.

# V.III. Schrödinger's Cat paradox resolution

In this section we shall consider the problem of the collapse of the cat state vector on the basis of two different hypotheses:

(**A**) The canonical postulate of QM is correct in all cases.

(**B**) The canonical interpretation of the wave function $\psi = c_1\psi_1 + c_2\psi_2$ is correct only when the supports the wave functions $\psi_1$ and $\psi_2$ essentially overlap. When the wave functions $\psi_1$, $\psi_2$ have separated supports (as in the case of the experiment that we are considering in section II) we claim that canonical interpretation of the wave function $\psi = c_1\psi_1 + c_2\psi_2$ is no longer valid for a such cat state, for details see Appendix CAT-C.

# V.III.1. Considerationtion of the Schrödinger's cat paradox using canonical von Neumann postulate

Let $|s_1(t)\rangle$ and $|s_2(t)\rangle$ be

$$|s_1(t)\rangle = \Big| \text{undecayed nucleus at instant } t \Big\rangle,$$

(5.3.1)

$$|s_2(t)\rangle = \Big| \text{decayed nucleus at instant } t \Big\rangle.$$

In a good approximation we assume now that

$$|s_1(0)\rangle = \int_{-\infty}^{+\infty} \Psi_{\mathbf{II}}^{\#}(x)|x\rangle dx$$

(5.3.2)

and

$$|s_2(0)\rangle = \int_{-\infty}^{+\infty} \Psi_{\mathbf{I}}^{\#}(x)|x\rangle dx.$$

(5.3.3)

**Remark 5.3.1**. Note that: (i) $|s_2(0)\rangle = \Big| \text{decayed nucleus at instant } 0 \Big\rangle = \Big| \text{free } \alpha\text{-particle at instant } 0 \Big\rangle$. (ii) Feynman propagator of a free $\alpha$-particle are [70]:

$$K_2(x, t, x_0) = \left( \frac{m}{2\pi i \hbar t} \right)^{1/2} \exp\left\{ \frac{i}{\hbar} \left[ \frac{m(x - x_0)^2}{2t} \right] \right\}.$$

(5.3.4)

Therefore from Eq.(5.3.3),Eq.(5.2.9) and Eq.(5.3.4) we obtain

$$|s_2(t)\rangle = \int_{-\infty}^{+\infty} \Psi_{\mathbf{I}}^{\#}(x,t)|x\rangle dx,$$

$$\Psi_{\mathbf{I}}^{\#}(x,t) = \int_{-\infty}^{0} \Psi_{\mathbf{I}}^{\#}(x_0) K_2(x,t,x_0) dx_0 =$$

$$(\pi r_c^2)^{-1/4} \times \left(\frac{m}{2\pi i h t}\right)^{1/2} \times \int_{-\infty}^{0} \theta_1(x_0,l) \exp\left(-\frac{x_0^2}{2r_c^2}\right) \exp\left(-i\frac{2\pi}{h}\sqrt{2mE}\,x_0\right) \times$$

$$\times \exp\left\{\frac{i}{h}\left[\frac{m(x-x_0)^2}{2t}\right]\right\} dx_0 =$$

<div style="text-align:right">(5.3.5)</div>

$$(\pi r_c^2)^{-1/4} \times \left(\frac{m}{2\pi i h_\varepsilon t}\right)^{1/2} \times \int_{-l}^{0} \theta_1(x_0,l) \exp\left(-\frac{x_0^2}{2r_c^2}\right) \times$$

$$\times \exp\left\{\frac{i}{h}\left[\frac{m(x-x_0)^2}{2t} - \pi\sqrt{4mE}\,x_0\right]\right\} dx_0 =$$

$$(\pi r_c^2)^{-1/4} \times \left(\frac{m}{2\pi i h t}\right)^{1/2} \times \int_{-l}^{0} \theta_1(x_0,l) \exp\left(-\frac{x_0^2}{2r_c^2}\right) \times \exp\left\{\frac{i}{h}\left[S(t,x,x_0)\right]\right\} dx_0,$$

where

<div style="text-align:right">(5.3.6)</div>

$$S(t,x,x_0) = \frac{m(x-x_0)^2}{2t} - \pi\sqrt{8mE}\,x_0.$$

We assume now that

$$h \ll 2r_c^2 \ll l^2 < 1. \qquad (5.3.7)$$

Oscillatory integral in RHS of Eq.(5.3.5) is calculated now directly using stationary phase approximation. The phase term $S(x,x_0)$ given by Eq.(5.3.6) is stationary when

$$\frac{\partial S(t,x,x_0)}{\partial x_0} = -\frac{m(x-x_0)}{t} - \pi\sqrt{8mE} = 0. \qquad (5.3.8)$$

Therefore

$$-\frac{m(x-x_0)}{t} - \pi\sqrt{8mE} = 0,$$
$$-(x-x_0) = \pi t\sqrt{8E/m}\,, \qquad (5.3.9)$$

and thus stationary point $x_0(t,x)$ are

$$x_0(t,x) = \pi t \sqrt{8E/m} + x. \tag{5.3.10}$$

Thus from Eq.(5.3.5) and Eq.(5.3.10) using stationary phase approximation we obtain

$$|s_2(t)\rangle = |s_2(t)\rangle = \int_{-\infty}^{+\infty} \Psi_1^\#(x,t)|x\rangle dx,$$

$$\Psi_1^\#(x,t) = \tag{5.3.11}$$

$$(\pi r_c^2)^{-1/4} \times \theta_1(x_0(t,x),l) \exp\left[-\frac{x_0^2(t,x)}{2r_c^2}\right] \times \exp\left\{\frac{i}{\hbar}\left[S(t,x,x_0(t,x))\right]\right\} + O(\hbar),$$

where

$$S(x,x_0(t,x)) = \frac{m(x-x_0(t,x))^2}{2t} - \pi\sqrt{8mE}\, x_0(t,x). \tag{5.3.12}$$

From Eq.(5.3.11) we obtain

$$\langle s_2(t)\|s_2(t)\rangle \simeq (\pi r_c^2)^{-1/2} \times \theta_1\left(x + \pi t\sqrt{8E/m}, l\right) \exp\left[-\frac{\left(x + \pi t\sqrt{8E/m}\right)^2}{r_c^2}\right]. \tag{5.3.13}$$

**Remark 5.3.2.** From the inequality (5.3.7) and Eq.(5.3.13) follows that $\alpha$-particle at each instant $t \geq 0$ moves quasiclassically from right to left by the law

$$x(t) = -\pi t\sqrt{8E/m}, \tag{5.3.14}$$

i.e. i.e.,estimating the position $x(t,x_0,t_0;\hbar)$ at each instant $t \geq 0$ with final error $r_c$ gives $|\langle x\rangle(t) - x(t)| \leq r_c, i = 1,\ldots,d$ with a probability

$$\mathbf{P}\{|\langle x\rangle(t,0,0;\hbar) - x(t)| \leq r_c\} = 1.$$

**Remark 5.3.3.** We assume now that a distance between radioactive source and internal monitor which detects a single atom decaying (see Pic.1) is equal to $L$.

**Proposition 5.3.1.** After $\alpha$-decay at instant $t = 0$ the collaps: $|\text{live cat}\rangle \rightarrow |\text{death cat}\rangle$ arises at instant

$$T = \frac{L}{\pi\sqrt{8E/m}} \tag{5.3.15}$$

with a probability $\mathbf{P}_T\big(|\text{death cat}\rangle\big)$ to observe a state $|\text{death cat}\rangle$ at instant $T$ is $\mathbf{P}_T\big(|\text{death cat}\rangle\big) = 1.$

**Proof.** Note that. In this case Schrödinger's cat in fact perrmorm the single measurement of $\alpha$-particle position with accuracy of $\delta x = l$ at instant $t = T$ (given by Eq.(5.3.15)) by internal monitor (see Pic.5.1.1). The probability of getting a result $L$

with accuracy of $\delta x = l$ given by

$$\int_{|L-x|\leq l/2} |\langle x|\|c_2|^2 s_2(T)\rangle|^2 dx = 1. \tag{5.3.16}$$

Therefore at instant $T$ the $\alpha$-particle kills Schrödinger's cat with a probability $\mathbf{P}_T\big(\big|\text{death cat}\big\rangle\big) = 1$.

**Remark 5.3.4.** Note that. When Schrödinger's cat has permormed this measurement the immediate post measurement state of $\alpha$-particle (by von Neumann postulate C.4, see Appendix CAT-C) will end up in the state

$$|\Psi_T\rangle = \frac{\displaystyle\int_{|L-x|\leq l/2} |x\rangle\langle x|\|s_2(T)\rangle dx}{\sqrt{\displaystyle\int_{|L-x|\leq l/2} |\langle x|\|s_2(T)\rangle|^2 dx}} \in \mathbf{S}_\Theta, \Theta = \{x|\|L-x| \leq l/2\} \tag{5.3.17}$$

From Eq.(5.3.17) one obtains

$$\langle x'\|\Psi_T\rangle = \int_{|L-x|\leq l/2} \langle x'\|x\rangle\langle x|\|s_2(T)\rangle dx = \int_{|L-x|\leq l/2} \delta(x'-x)\langle x|\|s_2(T)\rangle dx = \Psi_{\mathbf{I}}^{\#}(x',t). \tag{5.3.18}$$

Therefore the state $|\Psi_T\rangle$ again kills Schrödinger's cat with a probability $\mathbf{P}_T\big(\big|\text{death cat}\big\rangle\big) = 1$.

Suppose now that a nucleus $\mathbf{n}$, whose Hilbert space is spanned by orthonormal states $|s_i(t)\rangle$, $i = 1, 2$, where $|s_1(t)\rangle = \big|\text{undecayed nucleus at instant } t\big\rangle$ and $|s_2(t)\rangle = \big|\text{decayed nucleus at instant } t\big\rangle$ is in the superposition state

$$|\Psi_t\rangle_{\mathbf{n}} = c_1|s_1(t)\rangle + c_2|s_2(t)\rangle, |c_1|^2 + |c_2|^2 = 1. \tag{5.3.19}$$

**Remark 5.3.5.** Note that: (i) $|s_1(0)\rangle = \big|\text{undecayed nucleus at instant } t = 0\big\rangle =$

$= \big|\alpha\text{-particle iside region } (0, l] \text{ at instant } t = 0\big\rangle$. (ii) Feynman propagator of $\alpha$-particle inside region $(0, l]$ are [70]:

$$K_2(x, t, x_0) = \left(\frac{m}{2\pi i\hbar t}\right)^{1/2} \exp\left\{\frac{i}{\hbar}[S(t, x, x_0)]\right\}, \tag{5.3.20}$$

where

$$S(t, x, x_0) = \frac{m(x-x_0)^2}{2t} + mt(U_0 - E). \tag{5.3.21}$$

Therefore from Eq.(5.2.11)-Eq.(5.2.12) and Eq.(5.3.20)-Eq.(5.3.21) we obtain

$$|s_1(t)\rangle = \int_{-\infty}^{+\infty} \Psi_{\mathbf{II}}^{\#}(x, t)|x\rangle dx,$$

$$\Psi_{\mathbf{II}}^{\#}(x, t) = \int_0^l \Psi_{\mathbf{II}}^{\#}(x_0) K_2(x, t, x_0) dx_0 =$$

$$\left(\frac{m}{2\pi i\hbar t}\right)^{1/2} \int_0^l E(x_0, l)\Psi_{\mathbf{II}}(x_0)\theta_l(x_0) \exp\left\{\frac{i}{\hbar}[S(t, x, x_0)]\right\} dx_0, \tag{5.3.22}$$

where

$$\theta_l(x) = \begin{cases} 1 \text{ for } x \in [0,l] \\ 0 \text{ for } x \notin [0,l] \end{cases}$$

**Remark 5.3.6.** We assume for simplification now that

$$k' \leq 1. \tag{5.3.23}$$

Therefore oscillatory integral in RHS of Eq.(5.3.22) is calculated now directly using stationary phase approximation. The phase term $S(x, x_0)$ given by Eq.(5.3.21) is stationary when

$$\frac{\partial S(t, x, x_0)}{\partial x_0} = -\frac{m(x - x_0)}{t} = 0. \tag{5.3.24}$$

and thus stationary point $x_0(t, x)$ are

$$-x + x_0 = 0$$
$$x_0(t, x) = x. \tag{5.3.25}$$

Thus from Eq.(5.3.22) and Eq.(5.3.25) using stationary phase approximation we obtain

$$\Psi_{\mathbf{II}}^{\#}(x, t) =$$
$$E(x_0(t,x), l)\Psi_{\mathbf{II}}(x_0(t,x))\theta_l(x_0(t,x)) \exp\left\{ \frac{i}{\hbar}\left[ S(t, x, x_0(t,x)) \right] \right\} + O(h) =$$

$$= E(x, l)\Psi_{\mathbf{II}}(x)\theta_l(x) \exp\left\{ \frac{i}{\hbar}\left[ mt(U_0 - E) \right] \right\} + O(h) = \tag{5.3.26}$$

$$E(x, l)\theta_l(x)O(1) \exp\left\{ \frac{i}{\hbar}\left[ mt(U_0 - E) \right] \right\} + O(h).$$

Therefore from Eq.(3.22) and Eq.(3.26) we obtain

$$\langle s_1(t) | s_1(t) \rangle = |\Psi_{\mathbf{II}}^{\#}(x, t)|^2 = E^2(x, l)\theta_l(x)O(1) + O(h). \tag{3.27}$$

**Remark 5.3.7.** Note that for each instant $t > 0$ :
$\mathrm{supp}(\Psi_{\mathbf{II}}^{\#}(x, t)) \cap \mathrm{supp}(\Psi_{\mathbf{I}}^{\#}(x, t)) = \varnothing$.

**Remark 5.3.8.** Note that. From Eq.(5.3.11), Eq.(5.3.13), Eq.(5.3.19), Eq.(5.3.22)-Eq.(5.3.27) and Eq.(A.13), see Appendix CAT-A, by Remark 5.3.7 we obtain

$$_\mathbf{n}\langle\Psi_t|\widehat{x}|\Psi_t\rangle_\mathbf{n} = |c_1|^2\langle s_1(t)|\widehat{x}|s_1(t)\rangle + |c_2|^2\langle s_2(t)|\widehat{x}|s_2(t)\rangle +$$

$$c_1c_2^*\langle s_2(t)|\widehat{x}|s_1(t)\rangle + c_1^*c_2\langle s_2(t)|\widehat{x}|s_1(t)\rangle^* = \tag{5.3.28}$$

$$|c_1|^2\langle s_1(t)|\widehat{x}|s_1(t)\rangle + |c_2|^2\langle s_2(t)|\widehat{x}|s_2(t)\rangle = |c_1|^2 l + |c_2|^2 T\pi\sqrt{8E/m}\,.$$

**Proposition 5.3.2. (i)** Suppose that a nucleus $\mathbf{n}$ is in the superposition state $|\Psi_t\rangle_\mathbf{n}$ ($|\Psi_t\rangle_\mathbf{n}$-particle) given by Eq.(5.3.19). Then the collaps: $|$live cat$\rangle \to |$death cat$\rangle$ arises at instant

$$T_{\mathrm{col}} \approx \frac{L \pm l}{|c_2|^2\sqrt{8\pi^2 E/m}}\,. \tag{5.3.29}$$

with a probability $\mathbf{P}_{T_{\mathrm{col}}}\big(|$death cat$\rangle\big)$ to observe a state $|$death cat$\rangle$ at instant $T_{\mathrm{col}}$ is $\mathbf{P}_{T_{\mathrm{col}}}\big(|$death cat$\rangle\big) = |c_2|^2$.

**(ii)** Assume now a Schrödinger's cat has performed the single measurement of $|\Psi_t\rangle_\mathbf{n}$-particle position with accuracy of $\delta x = l$ at instant $T = T_{\mathrm{col}}$ (given by Eq.(5.3.29)) by internal monitor (see Pic.5.1.1) and the result $x \approx L \pm l$ is not observed by Schrödinger's cat. Then the collaps: $|$live cat$\rangle \to |$death cat$\rangle$ never arises at any instant $T > T_{\mathrm{col}}$ and a probability $\mathbf{P}_{T>T_{\mathrm{col}}}\big(|$death cat$\rangle\big)$ to observe a state $|$death cat$\rangle$ at instant $T > T_{\mathrm{col}}$ is $\mathbf{P}_{T>T_{\mathrm{col}}}\big(|$death cat$\rangle\big) = 0$.

**Proof. (i)** Note that for $t > 0$ the marginal density matrix $\rho(t)$ is diagonal

$$\rho(t) = \begin{pmatrix} |c_1|^2\int|\Psi_\mathbf{II}^{\#}(x,t)|^2dx & 0 \\ 0 & |c_2|^2\int|\Psi_\mathbf{I}^{\#}(x,t)|^2dx \end{pmatrix}$$

In this case a Schrödinger's cat in fact perform the single measurement of $|\Psi_t\rangle_\mathbf{n}$-particle position with accuracy of $\delta x = l$ at instant $t = T_{\mathrm{col}}$ (given by Eq.(5.3.29)) by internal monitor (see Pic.5.1.1). The probability of getting a result $L$

at instant $T \approx T_{\mathrm{col}}$ with accuracy of $\delta x = l$ given by

$$\int_{|L-x|\leq l/2}|\langle x||\Psi_t\rangle_\mathbf{n}|^2dx = \int_{|L-x|\leq l/2}|\langle x|c_1|s_1(T)\rangle + c_2|s_2(T)\rangle|^2dx =$$

$$\int_{|L-x|\leq l/2}|c_1\langle x||s_1(T)\rangle + c_2\langle x||s_2(T)\rangle|^2dx = \tag{5.3.30}$$

$$\int_{|L-x|\leq l/2}|c_1^2\Psi_\mathbf{II}^{\#2}(x,T) + c_2^2\Psi_\mathbf{I}^{\#2}(x,T) + 2c_1c_2\Psi_\mathbf{I}^{\#}(x,T)\Psi_\mathbf{II}^{\#}(x,T)|dx.$$

From Eq.(5.3.30) by Remark 5.3.7 and Eq.(5.3.13) one obtains

$$\int_{|L-x|\leq l/2}|\langle x||\Psi_T\rangle_\mathbf{n}|^2dx = \int_{|L-x|\leq l/2}|c_2^2\Psi_\mathbf{I}^{\#2}(x,T)|dx = |c_2|^2\int_{|L-x|\leq l/2}|\Psi_\mathbf{I}^{\#}(x,T)|^2dx = |c_2|^2. \tag{5.3.31}$$

Note that. When Schrödinger's cat has permormed this measurement and the result $x \approx L \pm l$ is observed, then the immediate post measurement state of $\alpha$-particle (by von Neumann postulate C.4,see Appendix CAT-C) will end up in the state

$$|\Psi_{T_{\text{col}}}\rangle_{\mathbf{n}} = \frac{\int_{|L-x|\leq l/2}|x\rangle\langle x||\Psi_{T_{\text{col}}}\rangle_{\mathbf{n}}dx}{\sqrt{\int_{|L-x|\leq l/2}|\langle x||\Psi_{T_{\text{col}}}\rangle_{\mathbf{n}}|^2 dx}} =$$

$$\frac{\int_{|L-x|\leq l/2}|x\rangle\langle x|(c_1|s_1(T_{\text{col}})\rangle + c_2|s_2(T_{\text{col}})\rangle)dx}{\sqrt{\int_{|L-x|\leq l/2}|\langle x||\Psi_{T_{\text{col}}}\rangle_{\mathbf{n}}|^2 dx}} = \qquad (5.3.32)$$

$$\frac{c_1\int_{|L-x|\leq l/2}|x\rangle\langle x||s_1(T_{\text{col}})\rangle + c_2\int_{|L-x|\leq l/2}|x\rangle\langle x||s_2(T_{\text{col}})\rangle dx}{\sqrt{\int_{|L-x|\leq l/2}|\langle x||\Psi_{T_{\text{col}}}\rangle_{\mathbf{n}}|^2 dx}} \in \mathbf{S}_{\Theta}, \Theta = \{x||L-x|\leq l/2\}.$$

From Eq.(5.3.32) by Eq.(5.3.31) and by Remark 5.3.7 one obtains

$$|\Psi_{T_{\text{col}}}\rangle_{\mathbf{n}} = \frac{\int_{|L-x|\leq l/2}|x\rangle\langle x||\Psi_{T_{\text{col}}}\rangle_{\mathbf{n}}dx}{\sqrt{\int_{|L-x|\leq l/2}|\langle x||\Psi_{T_{\text{col}}}\rangle_{\mathbf{n}}|^2 dx}} =$$

$$\frac{\int_{|L-x|\leq l/2}|x\rangle\langle x|(c_1|s_1(T_{\text{col}})\rangle + c_2|s_2(T_{\text{col}})\rangle)dx}{\sqrt{\int_{|L-x|\leq l/2}|\langle x||\Psi_{T_{\text{col}}}\rangle_{\mathbf{n}}|^2 dx}} = \qquad (5.3.32)$$

$$= \frac{c_2}{|c_2|}\int_{|L-x|\leq l/2}|x\rangle\langle x||s_2(T_{\text{col}})\rangle dx.$$

Obviously by Remark 5.3.4 the staite $|\Psi_{T_{\text{col}}}\rangle_{\mathbf{n}}$ kills Schrödinger's cat with a probability $\mathbf{P}_{T_{\text{col}}}\big(\big|\text{death cat}\big\rangle\big) = 1.$

**Proof**. **(ii)** The probability of getting a result $L$ at any instant $T > T_{\text{col}}$ with accuracy of $\delta x = l$ by Eq.(5.3.31) and Eq.(5.3.13) given by

$$\int_{|L-x|\leq l/2}|\langle x||\Psi_T\rangle_{\mathbf{n}}|^2 dx = \int_{|L-x|\leq l/2}|c_2^2\Psi_{\mathbf{I}}^{\#2}(x,T)|dx = |c_2|^2\int_{|L-x|\leq l/2}|\Psi_{\mathbf{I}}^{\#}(x,T)|^2 dx =$$

$$\simeq (\pi r_c^2)^{-1/2}\times\theta_1\Big(x+\pi T\sqrt{8E/m},l\Big)\exp\left[-\frac{\left(x+\pi T\sqrt{8E/m}\right)^2}{r_c^2}\right] = 0.$$

Thus standard formalism of continuous quantum measurements [65],[67] leads to a definite but unpredictable measurement outcomes, either $|s_1(t)\rangle$ or $|s_2(t)\rangle$ and thus $|\Psi_t\rangle_{\mathbf{n}}$ suddenly "collapses" at unpredictable instant $t'$ into the corresponding state $|s_i(t')\rangle, i = 1,2.$

# V.III.2.Resolution of the Schrödinger's cat paradox using generalized von Neumann postulate.

**Proposition 5.3.3.** Suppose that a nucleus $\mathbf{n}$ is in the superposition state given by Eq.(5.3.19). the collaps: $|\text{live cat}\rangle \rightarrow |\text{death cat}\rangle$ arises at instant

$$T_{\text{col}} = \frac{L}{|c_2|^2 \sqrt{8\pi^2 E/m}}. \tag{5.3.33}$$

with a probability $\mathbf{P}_{T_{\text{col}}}\big(|\text{death cat}\rangle\big)$ to observe a state $|\text{death cat}\rangle$ at instant $T_{\text{col}}$ is $\mathbf{P}_{T_{\text{col}}}\big(|\text{death cat}\rangle\big) = 1$.

**Proof.** Let us consider now a state $|\Psi_t\rangle_{\mathbf{n}}$ given by Eq.(5.3.19). This state consists of a summ of two wave packets $c_1 \Psi_{\mathbf{II}}^{\#}(x,t)$ and $c_2 \Psi_{\mathbf{I}}^{\#}(x,t)$. Wave packet $c_1 \Psi_{\mathbf{II}}^{\#}(x,t)$ present an $\alpha_{\mathbf{II}}$-particle which lives in region $\mathbf{II}$ with a probability $|c_1|^2$ (see Pic. 5.2.1). Wave packet $c_2 \Psi_{\mathbf{I}}^{\#}(x,t)$ present an $\alpha_{\mathbf{I}}$-particle which lives in region $\mathbf{I}$ with a probability $|c_2|^2$ (see Pic. 5.2.1) and moves from the right to the left. Note that $\mathbf{I} \cap \mathbf{II} = \varnothing$. From Eq.(5.3.28) follows that $\alpha_{\mathbf{I}}$-particle at each instant $t \geq 0$ moves quasiclassically from right to left by the law

$$x(t) = -|c_2|^2 \pi t \sqrt{8E/m}, \tag{5.3.34}$$

From Eq.(5.3.34) one obtains

$$T = T_{\text{col}} \simeq \frac{L}{|c_2| \sqrt{8\pi^2 E/m}}. \tag{5.3.35}$$

Note that. In this case Schrödinger's cat in fact permorm a single measurement of $|\Psi_t\rangle_{\mathbf{n}}$-particle position with accuracy of $\delta x = l$ at instant $t = T = T_{\text{col}}$ (given by Eq.(5.3.35)) by internal monitor (see Pic.5.1.1). The probability of getting the result

$L$ at instant $t = T_{\text{col}}$ with accuracy of $\delta x = l$ by Remark 5.3.7 and by postulate C.V.2 (see Appendix CAT-C) and by postulate C.IV.3 (see apendix CAT-C) given by

$$\int_{|L-x|\leq l/2} \Big[ |\langle x|c_1|s_1(T_{\text{col}})\rangle|^2 * |\langle x|c_2|s_2(T_{\text{col}})\rangle|^2 \Big] dx =$$

$$\int_{|L-x|\leq l/2} |c_2|^{-2}|c_1|^{-2} \Big[ |\langle x|c_1|^{-2}||s_1(T_{\text{col}})\rangle|^2 * |\langle x|c_2|^{-2}||s_2(T_{\text{col}})\rangle|^2 \Big] dx = \tag{5.3.36}$$

$$\int_{|L-x|\leq l/2} |c_2|^{-2}|c_1|^{-2} \Big[ |\Psi_{\mathbf{I}}^{\#}(x|c_2|^{-2}, T_{\text{col}})|^2 * |\Psi_{\mathbf{II}}^{\#}(x|c_1|^{-2}, T_{\text{col}})|^2 \Big] dx = 1.$$

Note that. When Schrödinger's cat has permormed this measurement and the result $x \approx L \pm l$ is observed, then the immediate post measurement state of $\alpha$-particle (by generalized von Neumann postulate C.V.3, see Appendix CAT-C) will end up in the state

$$|\Psi_{T_{\text{col}}}\rangle_{\mathbf{n}} = \frac{\int_{|L-x|\leq l/2} |x\rangle\langle x||\Psi_{T_{\text{col}}}\rangle_{\mathbf{n}} dx}{\sqrt{\int_{|L-x|\leq l/2} \Big[ |\langle x|s_1(T_{\text{col}})\rangle|^2 + |\langle x|s_2(T_{\text{col}})\rangle|^2 \Big] dx}} =$$

$$\frac{\int_{|L-x|\leq l/2} |x\rangle\langle x||s_2^{c_2}(T_{\text{col}})\rangle dx}{\sqrt{\int_{|L-x|\leq l/2} \Big[ |\langle x|s_2(T_{\text{col}})\rangle|^2 \Big] dx}} \in \mathbf{H}_\Theta, \Theta = \{x||L-x| \leq l/2\}. \tag{5.3.37}$$

The staite $|\Psi_{T_{\text{col}}}\rangle_{\mathbf{n}}$ again kills Schrödinger's cat with a probability $\mathbf{P}_{T_{\text{col}}}\big(|\text{death cat}\rangle\big) = 1$.

**Thus is the collapsed state of the cat always shows definite and predictable outcomes even if cat also consists of a superposition**:

$$\big|\text{cat}\big\rangle = c_1\big|\text{live cat}\big\rangle + c_2\big|\text{death cat}\big\rangle.$$

Contrary to van Kampen's [61] and some others' opinions, "looking" at the outcome changes nothing, beyond informing the observer of what has already happened.

We remain: there are widespread claims that Schrödinger's cat is not in a definite alive or dead state but is, instead, in a superposition of the two. van Kampen, for example, writes "The whole system is in a superposition of two states: one in which no decay has occurred and one in which it has occurred. Hence, the state of the cat also consists of a superposition: $\big|\text{cat}\big\rangle = c_1\big|\text{live cat}\big\rangle + c_2\big|\text{death cat}\big\rangle$. The state remains a superposition until an observer looks at the cat" [61].

# Appendix. CAT-A.

The time-dependent Schrodinger equation governs the time evolution of a quantum mechanical system:

$$i\hbar\frac{\partial\Psi(\mathbf{x},t)}{\partial t} = \widehat{\mathbf{H}}\Psi(\mathbf{x},t). \tag{$A.1$}$$

The average, or expectation, value $\langle x_i\rangle$ of an observable $x_i$ corresponding to a quantum mechanical operator $\hat{x}_i$ is given by:

$$\langle x_i\rangle(t,\mathbf{x}_0,t_0;\hbar) = \frac{\int\limits_{\mathbb{R}^d} x_i|\Psi(\mathbf{x},t,\mathbf{x}_0,t_0;\hbar)|^2 d^dx}{\int\limits_{\mathbb{R}^d}|\Psi(\mathbf{x},t,\mathbf{x}_0,t_0;\hbar)|^2 d^dx}. \tag{$A.2$}$$

$$i = 1,\ldots,d.$$

**Remark A.1.** We assume now that: the solution $\Psi(\mathbf{x},t,\mathbf{x}_0,t_0;\hbar)$ of the time-dependent Schrödinger equation (A.1) has a good approximation by a delta function such that

$$|\Psi(\mathbf{x},t,\mathbf{x}_0,t_0;\hbar)|^2 \simeq \prod_{i=1}^{d}\delta(x_i - x_i(t,\mathbf{x}_0,t_0)), \tag{$A.3$}$$

$$x_i(t,\mathbf{x}_0,t_0) = x_{i,0},$$

$$i = 1,\ldots,d.$$

**Remark A.2.** Note that under conditions given by Eq.(A.3) QM-system which governed by Schrödinger equation Eq.(A.1) completely evolve quasiclassically i.e. estimating the position $\{x_i(t,\mathbf{x}_0,t_0;\hbar)\}_{i=1}^{d}$ at each instant $t$ with final error $\delta$ gives $|\langle x_i\rangle(t,\mathbf{x}_0,t_0;\hbar) - x_i(t,\mathbf{x}_0,t_0)| \leq \delta, i = 1,\ldots,d$ with a probability

$$\mathbf{P}\{|\langle x_i\rangle(t,\mathbf{x}_0,t_0;\hbar) - x_i(t,\mathbf{x}_0,t_0)| \le \delta\} \simeq 1.$$

Thus from Eq.(A.2) and Eq.(A.3) we obtain

$$\langle x_i\rangle(t,\mathbf{x}_0,t_0;\hbar) \simeq$$

$$\simeq \frac{\int_{\mathbb{R}^d} x_i \prod_{i=1}^{d=1} \delta(x_i - x_i(t,\mathbf{x}_0,t_0)) d^d x}{\int_{\mathbb{R}^d} \prod_{i=1}^{d=1} \delta(x_i - x_i(t,\mathbf{x}_0,t_0)) d^d x} = x_i(t,\mathbf{x}_0,t_0). \qquad (A.4)$$

$$i = 1,\ldots,d.$$

Thus under condition given by Eq.(A.3) one obtain

$$\langle x_{i,t}\rangle(t,\mathbf{x}_0,t_0;\hbar) \simeq x_i(t,\mathbf{x}_0,t_0), \qquad (A.5)$$
$$i = 1,\ldots,d.$$

**Remark A.3**. Let $\Psi_i(\mathbf{x},t,\mathbf{x}_0,t_0), i = 1,2$ be the solutions of the time-dependent Schrödinger equation (A.1). We assume now that $\Phi(\mathbf{x},t,\mathbf{x}_0,\mathbf{y}_0,t_0)$ is a linear superposition such that

$$\Phi(\mathbf{x},t,\mathbf{x}_0,\mathbf{y}_0,t_0) = c_1\Psi_1(\mathbf{x},t,\mathbf{x}_0,t_0) + c_2\Psi_2(\mathbf{x},t,\mathbf{y}_0,t_0). \qquad (A.6)$$
$$|c_1|^2 + |c_2|^2 = 1.$$

Then we obtain

$$|\Phi(\mathbf{x},t,\mathbf{x}_0,\mathbf{y}_0,t_0)|^2 = (\Phi(\mathbf{x},t,\mathbf{x}_0,\mathbf{y}_0,t_0)\Phi^*(\mathbf{x},t,\mathbf{x}_0,\mathbf{y}_0,t_0)) =$$

$$= ([c_1\Psi_1(\mathbf{x},t,\mathbf{x}_0,t_0) + c_2\Psi_2(\mathbf{x},t,\mathbf{y}_0,t_0)]) \times$$

$$\times([c_1^*\Psi_1^*(\mathbf{x},t,\mathbf{x}_0,t_0) + c_2^*\Psi_2^*(\mathbf{x},t,\mathbf{x}_0,\mathbf{y}_0,t_0)]) = \qquad (A.7)$$

$$= |c_1|^2\left(|\Psi_1(\mathbf{x},t,\mathbf{x}_0,t_0)|^2\right) + c_1^*c_2(\Psi_1^*(\mathbf{x},t,\mathbf{x}_0)\Psi_2(\mathbf{x},t,\mathbf{y}_0,t_0)) +$$

$$|c_2|^2\left(|\Psi_2(\mathbf{x},t,\mathbf{y}_0,t_0)|^2\right) + c_1c_2^*(\Psi_1(\mathbf{x},t,\mathbf{x}_0)\Psi_2^*(\mathbf{x},t,\mathbf{y}_0,t_0)).$$

**Definition A.1**. Let $\langle\mathbf{x}\rangle(t,\mathbf{x}_0,\mathbf{y}_0,t_0)$ be a vector-function

$$\langle \mathbf{x} \rangle (t, \mathbf{x}_0, \mathbf{y}_0, t_0) \,:\, [0, T] \times \mathbb{R}^d \times \mathbb{R}^d \times [0, T] \to \mathbb{R}^d$$

$$\langle \mathbf{x} \rangle (t, \mathbf{x}_0, \mathbf{y}_0, t_0) = \{ \langle x_1 \rangle (t, \mathbf{x}_0, \mathbf{y}_0, t_0), \ldots, \langle x_d \rangle (t, \mathbf{x}_0, \mathbf{y}_0, t_0) \}, \qquad (A.8)$$

where

$$\langle x_i \rangle (t, \mathbf{x}_0, \mathbf{y}_0, t_0) = \int\limits_{\mathbb{R}^d} x_i |\Phi(\mathbf{x}, t, \mathbf{x}_0, \mathbf{y}_0, t_0)|^2 d^d x =$$

$$= |c_1|^2 \int\limits_{\mathbb{R}^d} x_i |\Psi_1(\mathbf{x}, t, \mathbf{x}_0, t_0)|^2 d^d x +$$

$$+ |c_2|^2 \int\limits_{\mathbb{R}^d} x_i |\Psi_2(\mathbf{x}, t, \mathbf{y}_0, t_0)|^2 d^d x + \qquad (A.9)$$

$$+ c_1^* c_2 \int\limits_{\mathbb{R}^d} x_i \Psi_1^*(\mathbf{x}, t, \mathbf{x}_0, t_0) \Psi_2(\mathbf{x}, t, \mathbf{y}_0, t_0) d^d x +$$

$$+ c_1 c_2^* \int\limits_{\mathbb{R}^d} x_i \Psi_1(\mathbf{x}, t, \mathbf{x}_0, t_0) \Psi_2^*(\mathbf{x}, t, \mathbf{y}_0, t_0) d^d x.$$

**Definition A.2**. Let $\Delta(t, \mathbf{x}_0, \mathbf{y}_0, t_0)$ be a vector-function

$$\Delta(t, \mathbf{x}_0, \mathbf{y}_0, t_0) \,:\, [0, T] \times \mathbb{R}^d \times \mathbb{R}^d \to \mathbb{R}^d$$

$$(\Delta(t, \mathbf{x}_0, \mathbf{y}_0, t_0)) = \{ \delta_1(t, \mathbf{x}_0, \mathbf{y}_0, t_0), \ldots, \delta_d(t, \mathbf{x}_0, \mathbf{y}_0, t_0) \}, \qquad (A.10)$$

where

$$\delta_i(t, \mathbf{x}_0, \mathbf{y}_0, t_0) = \delta[x_i(t, \mathbf{x}_0, \mathbf{y}_0, t_0)] =$$

$$= c_1^* c_2 \int_{\mathbb{R}^d} x_i \Psi_1^*(\mathbf{x}, t, \mathbf{x}_0, t_0) \Psi_2(\mathbf{x}, t, \mathbf{y}_0, t_0) d^d x +$$

$$+ c_1 c_2^* \int_{\mathbb{R}^d} x_i \Psi_1(\mathbf{x}, t, \mathbf{x}_0, t_0) \Psi_2^*(\mathbf{x}, t, \mathbf{y}_0, t_0) d^d x.$$

$(A.11)$

Substituting Eqs.(A.11) into Eqs.(A.9) gives

$$\langle x_i \rangle(t, \mathbf{x}_0, \mathbf{y}_0, t_0) = \int_{\mathbb{R}^d} x_i |\Phi(\mathbf{x}, t, \mathbf{x}_0, \mathbf{y}_0, t_0)|^2 d^d x =$$

$$= |c_1|^2 \int_{\mathbb{R}^d} x_i |\Psi_1(\mathbf{x}, t, \mathbf{x}_0, t_0)|^2 d^d x +$$

$$+ |c_2|^2 \int_{\mathbb{R}^d} x_i |\Psi_2(\mathbf{x}, t, \mathbf{y}_0, t_0)|^2 d^d x + \delta_i(t, \mathbf{x}_0, \mathbf{y}_0, t_0) =$$

$(A.12)$

$$= |c_1|^2 \langle x_i \rangle(t, \mathbf{x}_0, t_0) + |c_2|^2 \langle x_i \rangle(t, \mathbf{y}_0, t_0) + \delta_i(t, \mathbf{x}_0, \mathbf{y}_0, t_0).$$

Substitution equations (A.5) into equations (A.12) gives

$$\langle x_i \rangle(t, \mathbf{x}_0, \mathbf{y}_0, t_0) = \int_{\mathbb{R}^d} x_i |\Phi(\mathbf{x}, t, \mathbf{x}_0, \mathbf{y}_0, t_0)|^2 d^d x =$$

$$= |c_1|^2 \langle x_i \rangle(t, \mathbf{x}_0, t_0) + |c_2|^2 \langle x_i \rangle(t, \mathbf{y}_0, t_0) + \delta_i(t, \mathbf{x}_0, \mathbf{y}_0, t_0)$$

$(A.13)$

$$\simeq |c_1|^2 x_i(t, \mathbf{x}_0, t_0) + |c_2|^2 x_i(t, \mathbf{y}_0, t_0) + \delta_i(t, \mathbf{x}_0, \mathbf{y}_0, t_0).$$

# Appendix. CAT-B.

The Schrödinger equation (2.1) in region $\mathbf{I} = \{x | x < 0\}$ has the following form

$$\hbar^2 \frac{\partial^2 \Psi_{\mathbf{I}}(x)}{\partial x^2} + 2mE\Psi_{\mathbf{I}}(x) = 0. \qquad (B.1)$$

From Schrödinger equation (B.1) follows

$$\hbar^2 \int_{-\infty}^{0} \frac{\partial^2 \Psi_{\mathbf{I}}(x)}{\partial x^2} dx + 2mE \int_{-\infty}^{0} \Psi_{\mathbf{I}}(x) dx = 0. \qquad (B.2)$$

Let $\Psi_{\mathbf{I}}^{\#}(x)$ be a function

$$\Psi_{\mathbf{I}}^{\#}(x) = \phi(x)\Psi_{\mathbf{I}}(x), \qquad (B.3)$$

where

$$\phi(x) = (\pi r_c^2)^{-1/4} \exp\left(\frac{x^2}{2r_c^2}\right) \qquad (B.4)$$

see Eq.(2.9). Note that

$$\frac{\partial^2 [\phi(x)\Psi_{\mathbf{I}}(x)]}{\partial x^2} = \frac{\partial}{\partial x} \left[ \Psi_{\mathbf{I}}(x)\frac{\partial \phi(x)}{\partial x} + \phi(x)\frac{\partial \Psi_{\mathbf{I}}(x)}{\partial x} \right] =$$

$$\qquad (B.5)$$

$$2\frac{\partial \Psi_{\mathbf{I}}(x)}{\partial x}\frac{\partial \phi(x)}{\partial x} + \Psi_{\mathbf{I}}(x)\frac{\partial^2 \phi(x)}{\partial x^2} + \phi(x)\frac{\partial^2 \Psi_{\mathbf{I}}(x)}{\partial x^2}.$$

Therefore substitution (B.2) into LHS of the Schrödinger equation (B.1) gives

$$\hbar^2 \int\limits_{-\infty}^{0} \frac{\partial^2 \Psi_{\mathbf{I}}^{\#}(x)}{\partial x^2} dx + 2mE \int\limits_{-\infty}^{0} \Psi_{\mathbf{I}}^{\#}(x) dx =$$

$$\hbar^2 \int\limits_{-\infty}^{0} \frac{\partial^2 \phi(x)\Psi_{\mathbf{I}}(x)}{\partial x^2} dx + 2Em \int\limits_{-\infty}^{0} \phi(x)\Psi_{\mathbf{I}}(x) dx =$$

$$2\hbar^2 \int\limits_{-\infty}^{0} \frac{\partial \Psi_{\mathbf{I}}(x)}{\partial x} \frac{\partial \phi(x)}{\partial x} dx + \hbar^2 \int\limits_{-\infty}^{0} \Psi_{\mathbf{I}}(x) \frac{\partial^2 \phi(x)}{\partial x^2} dx +$$

$$+ \int\limits_{-\infty}^{0} \phi(x) \left\{ \hbar^2 \frac{\partial^2 \Psi_{\mathbf{I}}(x)}{\partial x^2} + 2Em \int\limits_{-\infty}^{0} \Psi_{\mathbf{I}}(x) \right\} dx.$$

$(B.6)$

Note that

$$\int\limits_{-\infty}^{0} \phi(x) \left\{ \hbar^2 \frac{\partial^2 \Psi_{\mathbf{I}}(x)}{\partial x^2} + 2Em \int\limits_{-\infty}^{0} \Psi_{\mathbf{I}}(x) \right\} dx = 0. \qquad \left( B.7 \right)$$

Therefore from Eq.(B.6) and Eq.(2.3)-Eq.(2.4) one obtain

$$\hbar^2 \int\limits_{-\infty}^{0} \frac{\partial^2 \Psi_{\mathbf{I}}^{\#}(x)}{\partial x^2} dx + 2mE \int\limits_{-\infty}^{0} \Psi_{\mathbf{I}}^{\#}(x) dx =$$

$$\hbar^2 \int\limits_{-\infty}^{0} \frac{\partial^2 \phi(x)\Psi_{\mathbf{I}}(x)}{\partial x^2} dx + 2Em \int\limits_{-\infty}^{0} \phi(x)\Psi_{\mathbf{I}}(x) dx =$$

$(B.8)$

$$= 2\hbar^2 \int\limits_{l}^{\infty} \frac{\partial \Psi_{\mathbf{I}}(x)}{\partial x} \frac{\partial \phi(x)}{\partial x} dx + \hbar^2 \int\limits_{l}^{\infty} \Psi_{\mathbf{I}}(x) \frac{\partial^2 \phi(x)}{\partial x^2} dx.$$

From Eq.(B.6) one obtain

$$\frac{\partial \phi(x)}{\partial x} = (\pi r_c^2)^{-1/4} \frac{\partial}{\partial x} \exp\left[-\frac{x^2}{2r_c^2}\right] = -(\pi r_c^2)^{-1/4} r_c^{-2} x \exp\left[-\frac{x^2}{2r_c^2}\right],$$

$$\frac{\partial^2 \phi(x)}{\partial x^2} = -(\pi r_c^2)^{-1/4} r_c^{-2} \exp\left[-\frac{x^2}{2r_c^2}\right] + \hspace{2cm} (B.9)$$

$$+(\pi r_c^2)^{-1/4} r_c^{-4} x^2 \exp\left[-\frac{x^2}{2r_c^2}\right].$$

From Eq.(B.9) and Eq.(2.3)-Eq.(2.4) one obtain

$$\hbar^2 \int\limits_{-\infty}^{0} \frac{\partial \Psi_{\mathbf{I}}(x)}{\partial x} \frac{\partial \phi(x)}{\partial x} dx =$$

$$-\frac{\hbar^2}{(\pi r_c^2)^{1/4} r_c^2} \int\limits_{-\infty}^{0} \frac{\partial \exp(ikx)}{\partial x} x \exp\left[-\frac{x^2}{2r_c^2}\right] dx =$$

$$-\frac{2\pi\sqrt{2mE}\,\hbar}{(\pi r_c^2)^{1/4} r_c^2} \int\limits_{-\infty}^{0} x \exp\left(i\frac{2\pi\sqrt{2mE}}{\hbar}x\right) \exp\left[-\frac{x^2}{2r_c^2}\right] dx, \hspace{1cm} (B.10)$$

$$k = \frac{2\pi}{\hbar}\sqrt{2mE}.$$

and

$$\hbar^2 \int\limits_{-\infty}^{0} \Psi_{\mathbf{I}}(x) \frac{\partial^2 \phi(x)}{\partial x^2} dx = -\frac{\hbar^2}{(\pi r_c^2)^{3/4} r_c^2} \int\limits_{-\infty}^{0} \exp(ikx) \exp\left[-\frac{x^2}{2r_c^2}\right] dx +$$

$$+\frac{\hbar^2}{(\pi r_c^2)^{1/4} r_c^2} \int\limits_{-\infty}^{0} x^2 \exp(ikx) \exp\left[-\frac{x^2}{2r_c^2}\right] dx. \hspace{1cm} (B.11)$$

# Appendix. CAT-C. Generalized Postulates for Continuous Valued Observables.

Suppose we have an $n$-dimensional quantum system.

**I.Then we claim the following**:

**C.I**. Any given $n$-dimensional quantum system is identified by a set $\Re$
$\Re \triangleq \langle \mathbf{H}, \Im, \Re, \mathscr{L}_{2,1}, \mathbf{G}, |\psi_t\rangle \rangle$ where:

**(i)** $\mathbf{H}$ that is some infinite-dimensional complex Hilbert space, **(ii)** $\Im = (\Omega, \mathscr{F}, \mathbf{P})$ that is complete probability space,**(iii)** $\Re = (\mathbb{R}^n, \Sigma)$ that is measurable space ,**(iv)** $\mathscr{L}_{2,1}(\Omega)$ that is complete space of random variables $X : \Omega \to \mathbb{R}^n$ such that $\int_\Omega \|X(\omega)\| d\mathbf{P} < \infty, \int_\Omega \|X(\omega)\|^2 d\mathbf{P} < \infty$ and **(v)** $\mathbf{G} : \mathbf{H} \to \mathscr{L}_{2,1}(\Omega)$ that is one to one correspondence such that

$$\left| \langle \overline{\psi} | \widehat{Q} | \psi \rangle \right| = \int_\Omega \left( \mathbf{G}\left[ \widehat{Q} | \psi \rangle \right](\omega) \right) d\mathbf{P} = \mathbf{E}_\Omega \left( \mathbf{G}\left[ \widehat{Q} | \psi \rangle \right](\omega) \right) \qquad (C.1)$$

for any $|\psi\rangle \in \mathbf{H}$ and for any Hermitian operator $\widehat{Q} : \mathbf{H} \to \mathbf{H}$,

**(vi)** $|\psi_t\rangle$ is an continuous vector function $|\psi_t\rangle : \mathbb{R}_+ \to \mathbf{H}$ which representedthe evolution of the quantum system $\Re$.

**C.I.2.** For any $|\psi_1\rangle, |\psi_2\rangle \in \mathbf{H}$ and for any Hermitian operator $\widehat{Q} : \mathbf{H} \to \mathbf{H}$ such that $\left\langle \overline{\psi}_1 \left| \widehat{Q} \right| \psi_2 \right\rangle = \left\langle \psi_2 \left| \widehat{Q} \right| \psi_1 \right\rangle = 0$ :

$$\mathbf{G}\left[ \widehat{Q}(|\psi_1\rangle + |\psi_2\rangle) \right](\omega) = \mathbf{G}\left[ \widehat{Q} | \psi_1 \rangle \right](\omega) + \mathbf{G}\left[ \widehat{Q} | \psi_2 \rangle \right](\omega). \qquad (C.2)$$

**C.I.3.** Suppose that the evolution of the quantum system is represented by continuous vector function $|\psi_t\rangle : \mathbb{R}_+ \to \mathbf{H}$. Then any process of continuous measurements on measuring observable $\widehat{Q}$ for the system in state $|\psi_t\rangle$ one can to describe by an continuous $\mathbb{R}^n$-valued stochastic processes $X_t(\omega) = X_t(\omega; |\psi_t\rangle)$ given on probability space $(\Omega, \mathscr{F}, \mathbf{P})$ and a measurable space $(\mathbb{R}^n, \Sigma)$.

**Remark C.1**.We assume now for short but without loss of generality that $n = 1$.

**Remark C.2**. Let $X(\omega)$ be random variable $X(\omega) \in \mathscr{L}_{2,1}(\Omega)$ such that $X(\omega) = \mathbf{G}[|\psi\rangle](\omega)$, then we denote such random variable by $X_{|\psi\rangle}(\omega)$. The probability density of random variable $X_{|\psi\rangle}(\omega)$ we denote by $p_{|\psi\rangle}(q), q \in \mathbb{R}$.

**Definition C.1**. The *classical pure states* correspond to vectors $\mathbf{v} \in \mathbf{H}$ of norm $\|\mathbf{v}\| = 1$. Thus the set of all classical pure states corresponds to the unit sphere $\mathbf{S}^\infty \subset \mathbf{H}$ in a Hilbert space $\mathbf{H}$.

**Definition C.2**.The projective Hilbert space $P(\mathbf{H})$ of a complex Hilbert space $\mathbf{H}$ is the set of equivalence classes $[\mathbf{v}]$ of vectors $\mathbf{v}$ in $\mathbf{H}$, with $\mathbf{v} \neq \mathbf{0}$, for the equivalence relation given by $\mathbf{v} \sim_P \mathbf{w} \Leftrightarrow \mathbf{v} = \lambda \mathbf{w}$  for some non-zero complex number $\lambda \in \mathbb{C}$. The equivalence classes for the relation $\sim_P$  are also called rays or projective rays.

**Remark C.3**.The physical significance of the projective Hilbert space $P(\mathbf{H})$ is that in canonical quantum theory, the states $|\psi\rangle$ and $\lambda|\psi\rangle$ represent the same physical state of the quantum system, for any $\lambda \neq 0$. It is conventional to choose a state $|\psi\rangle$ from the ray $[|\psi\rangle]$ so that it has unit norm $\sqrt{\langle \psi | \psi \rangle} = 1$.

**Remark C.4**. In contrast with canonical quantum theory we have used instead contrary to $\sim_P$ equivalence relation $\sim_Q$, a Hilbert space $\mathbf{H}$, see Definition C.7.

**Definition C.3**.The *non-classical pure states* correspond to the vectors $\mathbf{v} \in \mathbf{H}$ of a norm $\|\mathbf{v}\| \neq 1$. Thus the set of all non-classical pure states corresponds to the set $\mathbf{H} \backslash \mathbf{S}^\infty \subset \mathbf{H}$ in the Hilbert space $\mathbf{H}$.

Suppose we have an observable $Q$ of a quantum system that is found through an exhaustive series of measurements, to have a set $\mathfrak{I}$ of values $q \in \mathfrak{I}$ such that $\mathfrak{I} = \cup_{i=1}^{m} (\theta_1^i, \theta_2^i), m \geq 2, (\theta_1^i, \theta_2^i) \cap (\theta_1^j, \theta_2^j) = \varnothing, i \neq j$. Note that in practice any observable $Q$ is measured to an accuracy $\delta q$ determined by the measuring device. We represent now by $|q\rangle$ the idealized state of the system in the limit $\delta q \to 0$, for which the observable definitely has the value $q$.

**II**.**Then we claim the following**:

**C**.**II**.**1**.The states $\{|q\rangle : q \in \mathfrak{I}\}$ form a complete set of $\delta$-function normalized basis states for the state space $\mathbf{H}_\mathfrak{I}$ of the system.

That the states $\{|q\rangle : q \in \mathfrak{I}\}$ form a complete set of basis states means that any state $|\psi[\mathfrak{I}]\rangle \in \mathbf{H}_\mathfrak{I}$ of the system can be expressed as: $|\psi[\mathfrak{I}]\rangle = \int_\mathfrak{I} c_{\psi[\mathfrak{I}]}(q) dq$, where $\operatorname{supp}(c_{\psi[\mathfrak{I}]}(q)) \subseteq \mathfrak{I}$ and while $\delta$-function normalized means that $\langle q|q'\rangle = \delta(q - q')$ from which follows $c_{\psi[\mathfrak{I}]}(q) = \langle q|\psi[\mathfrak{I}]\rangle$ so that $|\psi[\mathfrak{I}]\rangle = \int_\mathfrak{I} |q\rangle\langle q|\psi[\mathfrak{I}]\rangle dq$.

The completeness condition can then be written as $\int_\mathfrak{I} |q\rangle\langle q| dq = \widehat{\mathbf{1}}_{\mathbf{H}_\mathfrak{I}}$.

**C**.**II**.**2**.For the system in state $|\psi[\mathfrak{I}]\rangle$ the probability $P(q, q + dq; |\psi[\mathfrak{I}]\rangle)$ of obtaining the result $q \in \mathfrak{I}$ lying in the range $(q, q + dq) \subset \mathfrak{I}$ on measuring observable $Q$ is given by

$$P(q, q + dq; |\psi[\mathfrak{I}]\rangle) = p_{|\psi[\mathfrak{I}]\rangle}(q) dq \qquad (C.3)$$

for any $|\psi[\mathfrak{I}]\rangle \in \mathbf{H}_\mathfrak{I}$.

**Remark C**.**5**. Note that in general case $p_{|\psi[\mathfrak{I}]\rangle}(q) \neq |c_{\psi[\mathfrak{I}]}(q)|^2$.

**C**.**II**.**3**.The observable $Q_\mathfrak{I}$ is represented by a Hermitian operator $\widehat{Q}_\mathfrak{I} : \mathbf{H}_\mathfrak{I} \to \mathbf{H}_\mathfrak{I}$ whose eigenvalues are the possible results $\{q : q \in \mathfrak{I}\}$, of a measurement of $Q_\mathfrak{I}$, and the associated eigenstates are the states $\{|q\rangle : q \in \mathfrak{I}\}$, i.e. $\widehat{Q}_\mathfrak{I}|q\rangle = q|q\rangle, q \in \mathfrak{I}$.

**Remark C**.**6**. Note that the spectral decomposition of the operator $\widehat{Q}_\mathfrak{I}$ is then

$$\widehat{Q}_\mathfrak{I} = \int_\mathfrak{I} q|q\rangle\langle q| dq. \qquad (C.3)$$

**Definition C**.**4**. A connected set in $\mathbb{R}$ is a set $X \subset \mathbb{R}$ that cannot be partitioned into two nonempty subsets which are open in the relative topology induced on the set. Equivalently, it is a set which cannot be partitioned into two nonempty subsets such that each subset has no points in common with the set closure of the other.

**Definition C**.**5**. The *well localized pure states* $|\psi[\Theta]\rangle$ *with a support* $\Theta = (\theta_1, \theta_2)$ correspond to vectors of norm 1 and such that: $\operatorname{supp}(c_{\psi[\Theta]}(q)) = \Theta$ is a connected set in $\mathbb{R}$ Thus the set of all well localized pure states corresponds to the unit sphere $\mathbf{S}_\Theta^\infty \subsetneqq \mathbf{S}^\infty \subset \mathbf{H}$ in the Hilbert space $\mathbf{H}_\Theta \subsetneqq \mathbf{H}$.

Suppose we have an observable $Q_\Theta$ of a system that is found through an exhaustive series of measurements, to have a continuous range of values $q$ : $\theta_1 < q < \theta_2$.

**III**.**Then we claim the following**:

**C**.**III**.**1**.For the system in well localized pure statestate $|\psi[\Theta]\rangle$ such that:

(i) $|\psi[\Theta]\rangle \in \mathbf{S}_\Theta^\infty$ and

(ii) supp$(c_{\psi[\Theta]}(q)) \triangleq \{q|c_{\psi[\Theta]}(q) \neq 0\}$ is a connected set in $\mathbb{R}$, then the probability $P(q, q + dq; |\psi[\Theta]\rangle)$ of obtaining the result $q$ lying in the range $(q, q + dq)$ on measuring

observable $Q_\Theta$ is given by

$$P(q, q + dq; |\psi[\Theta]\rangle) = |\langle q|\psi[\Theta]\rangle|^2 dq = |c_{\psi[\Theta]}(q)|^2 dq. \qquad (C.4)$$

**C.III.2.** $p_{|\psi[\Theta]\rangle}(q)dq = |\langle q|\psi[\Theta]\rangle|^2 dq = |c_{\psi[\Theta]}(q)|^2 dq.$

**C.III.3.** Let $|\psi[\Theta_1]\rangle$ and $|\psi[\Theta_2]\rangle$ be well localized pure states with $\Theta_1 = (\theta_1^1, \theta_2^1)$ and $\Theta_2 = (\theta_1^2, \theta_2^2)$ correspondingly. Let $X_1(\omega) = X_{|\psi[\Theta_1]\rangle}(\omega)$ and $X_2(\omega) = X_{|\psi[\Theta_2]\rangle}(\omega)$ correspondingly. Assume that $\overline{\Theta}_1 \cap \overline{\Theta}_2 = \varnothing$ (here the closure of $\Theta_i, i = 1, 2$ is denoted by

$\overline{\Theta}_i, i = 1, 2$) then random variables $X_1(\omega)$ and $X_2(\omega)$ are independent.

**C.III.4.** If the system is in well localized pure state $|\psi[\Theta]\rangle$ the state $|\psi[\Theta]\rangle$ described by a wave function $\psi(q, \Theta) = \langle q||\psi[\Theta]\rangle$ and the value of observable $Q_\Theta$ is measured once each

on many identically prepared system, the average value of all the measurements will be

$$\langle Q_\Theta \rangle = \frac{\int\limits_\Theta q|\psi(q, \Theta)|^2 dq}{\int\limits_\Theta |\psi(q, \Theta)|^2 dq}. \qquad (C.5)$$

The completeness condition can then be written as $\int_\Theta |q\rangle\langle q| dq = \widehat{\mathbf{1}}_{\mathbf{H}_\Theta}$.

Completeness means that for any state $|\psi[\Theta]\rangle \in \mathbf{S}_\Theta^\infty$ it must be the case that $\int_\Theta |\langle q|\psi[\Theta]\rangle|^2 dq \neq 0$, i.e. there must be a non-zero probability to get some result on measuring observable $Q_\Theta$.

**C.III.5.** (**von Neumann measurement postulate**) Assume that

(i) $|\psi\rangle \in \mathbf{S}_\Theta^\infty$ and (ii) supp$(c_\psi(q)) = \Theta$ is a connected set in $\mathbb{R}$. Then if on performing a measurement of $Q_\Theta$ with an accuracy $\delta q$, the result is obtained in the range $(q - \frac{1}{2}\delta q, q + \frac{1}{2}\delta q)$, then the system will end up in the state

$$\frac{\widehat{P}(q, \delta q)|\psi[\Theta]\rangle}{\sqrt{\langle\psi|\widehat{P}(q, \delta q)|\psi[\Theta]\rangle}} = \frac{\int_{|q-q'|\leq\delta q/2} |q'\rangle\langle q'||\psi[\Theta]\rangle dq'}{\sqrt{\int_{|q-q'|\leq\delta q/2} |\langle q'||\psi[\Theta]\rangle|^2 dq'}}. \qquad (C.6)$$

**IV.** **We claim the following**:

**C.IV.1** For the system in state $|\psi^a[\Theta]\rangle = a|\psi[\Theta]\rangle \in \mathbf{H}_\Theta$, where: (i) $|\psi[\Theta]\rangle \in \mathbf{S}_\Theta^\infty, |a| \neq 1$,

(ii) supp$(c_{\psi[\Theta]}(q))$ is a connected set in $\mathbb{R}$ and (iii) $|\psi[\Theta]\rangle = \int_{\theta_1}^{\theta_2} c_{\psi[\Theta]}(q)|q\rangle dq$

$$\mathbf{G}\left[\widehat{Q}_\Theta|\psi^a[\Theta]\rangle\right] = |a|^2 \mathbf{G}\left[\widehat{Q}_\Theta|\psi[\Theta]\rangle\right]. \qquad (C.6)$$

**C.IV.2**. Assume that the system in state $|\psi^a[\Theta]\rangle = a|\psi[\Theta]\rangle \in \mathbf{H}_\Theta$, where (i) $|\psi[\Theta]\rangle \in \mathbf{S}_\Theta^\infty$, $|a| \neq 1$, (ii) $\text{supp}(c_{\psi[\Theta]}(q))$ is a connected set in $\mathbb{R}$ and (iii) $|\psi[\Theta]\rangle = \int_{\theta_1}^{\theta_2} c_{\psi[\Theta]}(q)|q\rangle dq$.

Then if the system is in state $|\psi^a[\Theta]\rangle$ described by a wave function $\psi^a(q;\Theta) = \langle q||\psi^a[\Theta]\rangle$ and the value of observable $Q_\Theta$ is measured once each on many identically prepared system, the average value of all the measurements will be

$$\langle Q_\Theta \rangle = \int_\Theta q|\psi^a(q;\Theta)|^2 dq. \qquad (C.7)$$

**C.IV.3**. The probability $P(q, q+dq; |\psi^a[\Theta]\rangle)dq$ of obtaining the result $q$ lying in the range

$(q, q+dq)$ on measuring $Q_\Theta$ is

$$P(q, q+dq; |\psi^a[\Theta]\rangle)dq = |a|^{-2}|c_{\psi[\Theta]}(q|a|^{-2})|^2 dq. \qquad (C.8)$$

**Remark C.7**. Note that C.IV.3 immediately folows from C.IV.1 and C.III.2.

**C.IV.4**. (**Generalized von Neumann measurement postulate**) If on performing a measurement of observable $Q_\Theta$ with an accuracy $\delta q$, the result is obtained in the range $(q - \frac{1}{2}\delta q, q + \frac{1}{2}\delta q)$, then the system immediately after measurement will end up in the state

$$\frac{\widehat{P}(q, \delta q)|\psi^a[\Theta]\rangle}{\sqrt{\langle \psi|\widehat{P}(q, \delta q)|\psi[\Theta]\rangle}} = \frac{\int_{|q-q'|\leq \delta q/2} |q'\rangle\langle q'||\psi^a[\Theta]\rangle dq'}{\sqrt{\int_{|q-q'|\leq \delta q/2} |\langle q'||\psi[\Theta]\rangle|^2 dq'}} =$$

$$\frac{a\int_{|q-q'|\leq \delta q/2} |q'\rangle\langle q'||\psi[\Theta]\rangle dq'}{\sqrt{\int_{|q-q'|\leq \delta q/2} |\langle q'||\psi[\Theta]\rangle|^2 dq'}} \in \mathbf{H}_\Theta. \qquad (C.9)$$

**C.V.1**. Let $|\Psi^{a_1,a_2}[\Theta_1, \Theta_2]\rangle = |\psi_1^{a_1}[\Theta_1]\rangle + |\psi_2^{a_2}[\Theta_2]\rangle \in \mathbf{H}_{1,2} \triangleq \mathbf{H}_{\Theta_1} \oplus \mathbf{H}_{\Theta_2} \subsetneqq \mathbf{H}$, where (i) $|\psi_i^{a_i}[\Theta_i]\rangle = a_i|\psi_i[\Theta_i]\rangle \in \mathbf{H}_{\Theta_i}, |\psi_i\rangle = |\psi_i[\Theta_i]\rangle \in \mathbf{S}_{\Theta_i}^\infty, |a_i| \neq 1, i = 1, 2$;
(ii) $\text{supp}(c_{\psi_i[\Theta_i]}(q)), i = 1, 2$ is a connected sets in $\mathbb{R}$;
(iii) $\left(\text{supp}(c_{\psi_1[\Theta_1]}(q))\right) \cap \left(\text{supp}(c_{\psi_2[\Theta_2]}(q))\right) = \varnothing$ and
(iv) $|\psi_i[\Theta_i]\rangle = \int_{\theta_1}^{\theta_2} c_{\psi_i[\Theta_i]}(q)|q\rangle dq, i = 1, 2$.

Then if the system is in a state $|\Psi^{a_1,a_2}[\Theta_1, \Theta_2]\rangle$ described by a wave function $\Psi^{a_1,a_2}(q; \Theta_1, \Theta_2) = \langle q||\Psi^{a_1,a_2}[\Theta_1, \Theta_2]\rangle, q \in \Theta_1 \cup \Theta_2$ and the value of observable $Q_{\Theta_1,\Theta_2}$ is measured once each on many identically prepared system, the average value of all the measurements will be

$$\langle Q_{\Theta_1, \Theta_2} \rangle = \int_{\Theta_1 \cup \Theta_2} q |\Psi^{a_1, a_2}(q; \Theta_1, \Theta_2)|^2 dq. \tag{C.10}$$

**C.V. 2.** The probability of getting a result $q$ with an accuracy $\delta q$ such that $(q - \frac{1}{2}\delta q, q + \frac{1}{2}\delta q) \in supp(c_{\psi_1}(q))$ or $(q - \frac{1}{2}\delta q, q + \frac{1}{2}\delta q) \in supp(c_{\psi_2}(q))$ given by

$$\int_{|q - q'| \leq \delta q/2} \Big[ \left( |\langle q' \| \psi_1^{a_1}[\Theta_1] \rangle|^2 \right) * \left( |\langle q' \| \psi_2^{a_2}[\Theta_2] \rangle|^2 \right) \Big] dq'. \tag{C.11}$$

**Remark C.8.** Note that C.IV.3 immediately folows from C.III.3.

**C.Ψ. 3.** Assume that the system is initially in the state $|\Psi^{a_1, a_2}[\Theta_1, \Theta_2]\rangle$. If on performing a measurement of $Q_{\Theta_1, \Theta_2}$ with an accuracy $\delta q$, the result is obtained in the range $(q - \frac{1}{2}\delta q, q + \frac{1}{2}\delta q)$, then the state of the system immediately after measurement given by

$$\frac{\widehat{P}(q_i, \delta q)|\Psi^{a_1, a_2}[\Theta_1, \Theta_2]\rangle}{\sqrt{\langle \psi |\widehat{P}(q_i, \delta q)|\psi \rangle}} =$$

$$\frac{\int_{|q_i - q'| \leq \delta q/2} (|q'\rangle\langle q' \| \psi_1^{a_1}[\Theta_1]\rangle + |q'\rangle\langle q' \| \psi_2^{a_2}[\Theta_2]\rangle) dq'}{\sqrt{\int_{|q_i - q'| \leq \delta q/2} \Big[ |\langle q' \| \psi_1[\Theta_1]\rangle|^2 + |\langle q' \| \psi_2[\Theta_2]\rangle|^2 \Big] dq'}} = \tag{C.12}$$

$$\frac{\int_{|q_i - q'| \leq \delta q/2} (a_1 |q'\rangle\langle q' \| \psi_1[\Theta_1]\rangle + a_2 |q'\rangle\langle q' \| \psi_2[\Theta_2]\rangle) dq'}{\sqrt{\int_{|q_i - q'| \leq \delta q/2} \Big[ |\langle q' \| \psi_1[\Theta_1]\rangle|^2 + |\langle q' \| \psi_2[\Theta_2]\rangle|^2 \Big] dq'}} \in \mathbf{H}_{\Theta_i},$$

$$q_i \in \Theta_i, i = 1, 2.$$

**Definition C.6.** Let $\mathbf{H}_{1,2}$ be $\mathbf{H}_{1,2} \triangleq \mathbf{H}_{\Theta_1} \oplus \mathbf{H}_{\Theta_2}$.

**Definition C.7.** Let $|\psi^a\rangle$ be a state $|\psi^a\rangle = a|\psi\rangle$, where $|\psi\rangle \in \mathbf{S}^\infty, |a| \neq 1$ and $|\psi\rangle = \int_{\theta_1}^{\theta_2} c_\psi(q)|q\rangle dq$. Let $|\psi_a\rangle$ be an state such that $|\psi_a\rangle \in \mathbf{S}^\infty$. States $|\psi^a\rangle$ and $|\psi_a\rangle$ is a $Q$-equivalent: $|\psi^a\rangle \sim_Q |\psi_a\rangle$ iff

$$P(q, q + dq; |\psi^a\rangle) = |a|^{-2} |c_\psi(q|a|^{-2})|^2 dq = P(qq + dq; |\psi_a\rangle) dq \tag{C.13}$$

**C.V.** For any state $|\psi^a\rangle = a|\psi\rangle$, where $|\psi\rangle \in \mathbf{S}^\infty, |a| \neq 1$ and $|\psi\rangle = \int_{\theta_1}^{\theta_2} c_\psi(q)|q\rangle dq$ there exist an state $|\psi_a\rangle \in \mathbf{S}^\infty$ such that: $|\psi^a\rangle \sim_Q |\psi_a\rangle$.

**Definition C.8.** Let $|\psi^a\rangle$ be a state $|\psi^a\rangle = a|\psi\rangle$, where $|\psi\rangle \in \mathbf{S}^\infty, |a| \neq 1$ and $|\psi\rangle = \int_{\theta_1}^{\theta_2} c_\psi(q)|q\rangle dq$. Let $|\psi_a\rangle$ be an state such that $|\psi_a\rangle \in \mathbf{S}^\infty$. States $|\psi^a\rangle$ and $|\psi_a\rangle$ is a $\widehat{Q}$-equivalent ($|\psi^a\rangle \sim_{\widehat{Q}} |\psi_a\rangle$) iff: $\langle \psi^a| \widehat{Q} \psi^a\rangle = \langle \psi_a| \widehat{Q} \psi_a\rangle$.

**C.VI.** For any state $|\psi^a\rangle = a|\psi\rangle$, where $|\psi\rangle \in \mathbf{S}^\infty, |a| \neq 1$ and $|\psi\rangle = \int_{\theta_1}^{\theta_2} c_\psi(q)|q\rangle dq$ there exist an state $|\psi_a\rangle \in \mathbf{S}^\infty$ such that: $|\psi^a\rangle \sim_{\widehat{Q}} |\psi_a\rangle$.

# Apendix I.1.Derivation of the Path Integral for the non-Hermitian Hamiltonian.

Here we derive Path Integral for the case of a generalized non-Hermitian Hamiltonian

$$(H_\varepsilon(p,x))_\varepsilon = \frac{p^2}{2m} - \frac{im\omega x^2}{2} + (V_\varepsilon(x))_\varepsilon. \qquad (1.1.1)$$

Here $\mathbf{cl}[(V_\varepsilon(x))_\varepsilon] \in G(\mathbb{R}^n)$.

The quantity we want to calculate is

$$(K_\varepsilon(x_1, t + \Delta t | x_0, t))_\varepsilon = (\langle x_1, t + \Delta t | x_0, t \rangle_\varepsilon)_\varepsilon =$$

$$\left( \left\langle x_1 \left| \exp\left(-\frac{iH_\varepsilon \Delta t}{\epsilon}\right) \right| x_0, t \right\rangle \right)_\varepsilon = \left( \int dp \langle x_1 | p \rangle \left\langle p \left| \exp\left(-\frac{iH_\varepsilon \Delta t}{\epsilon}\right) \right| x_0 \right\rangle \right)_\varepsilon \qquad (1.1.2)$$

the last factor in Eq.(1.1.2) can be estimated as

$$\left( \left\langle p \left| \exp\left(-\frac{iH_\varepsilon \Delta t}{\epsilon}\right) \right| x_0 \right\rangle \right)_\varepsilon = \left( \left\langle p \left| 1 - \frac{iH_\varepsilon \Delta t}{\epsilon} + O(\Delta t)^2 \right| x_0 \right\rangle \right)_\varepsilon =$$

$$\left( 1 - \frac{i}{\epsilon} \frac{p^2 \Delta t}{2m} - \frac{1}{\epsilon} \frac{x_0^2 m\omega \Delta t}{2} - \frac{i}{\epsilon} (V_\varepsilon(x_0, t))_\varepsilon \Delta t + O(\Delta t)^2 \right) \frac{\exp(-ipx_0/\epsilon)}{\sqrt{2\pi\epsilon}} = \qquad (1.1.3)$$

$$\frac{1}{\sqrt{2\pi\epsilon}} \exp\left[ -\frac{i}{\epsilon} \left( px_0 + \frac{p^2}{2m}\Delta t - \frac{x_0^2 m\omega}{2}\Delta t + (V_\varepsilon(x_0, t))_\varepsilon \Delta t + O(\Delta t)^2 \right) \right].$$

Substitution Eq.(1.1.3) into Eq.(1.1.2) gives

$$\langle x_1, t + \Delta t | x_0, t \rangle =$$

$$\frac{1}{\sqrt{2\pi\epsilon}} \left( \int dp \exp\left( -\frac{ipx_1}{\epsilon} \right) \times \right.$$

$$\left. \exp\left[ -\frac{i}{\epsilon} \left( px_0 + \frac{p^2}{2m}\Delta t - \frac{i\omega x_0^2}{2}\Delta t + (V_\varepsilon(x_0, t))_\varepsilon \Delta t + O(\Delta t)^2 \right) \right] \right)_\varepsilon = \qquad (1.1.4)$$

$$\sqrt{\frac{m}{2\pi i\epsilon\Delta t}} \exp\left[ \frac{i}{\epsilon} \left( \frac{m}{2} \frac{(x_1 - x_0)^2}{\Delta t} + \frac{im\omega x_0^2}{2}\Delta t - (V_\varepsilon(x_0, t))_\varepsilon \Delta t + O(\Delta t)^2 \right) \right]$$

By formally identifying $\dot{x}^2(t_1) = (x_1 - x_0)^2/(\Delta t)^2$ from Eq.(1.1.4) we obtain

$$\langle x_1, t + \Delta t | x_0, t \rangle = \sqrt{\frac{m}{2\pi i\epsilon\Delta t}} \left( \exp\left[ \frac{i}{\epsilon} (\mathcal{L}_\varepsilon(\dot{x}(t), x(t), t)\Delta t) \right] \right)_\varepsilon, \qquad (1.1.5)$$

Here

$$\mathcal{L}_\varepsilon(\dot{x}(t), x(t), t) = \frac{m}{2}\dot{x}^2(t) + \frac{im\omega x^2}{2} - (V_\varepsilon(x, t))_\varepsilon \qquad (1.1.6)$$

is the complex classical Lagrangian.Using the completeness relation $N$ times one obtain

$$\left( \langle x_f, t_f | x_i, t_i \rangle_\varepsilon \right)_\varepsilon = \int ((\langle x_f, t_f | x_{N-1}, t_{N-1} \rangle)_\varepsilon) dx ((_{N-1}\langle x_{N-1}, t_{N-1} | x_{N-2}, t_{N-2} \rangle)_\varepsilon) dx_{N-2} \ldots$$
$$\ldots dx_2 ((\langle x_2, t_2 | x_1, t_1 \rangle)_\varepsilon) dx_1 ((\langle x_1, t_1 | x_i, t_i \rangle)_\varepsilon). \qquad (1.1.7)$$

The time interval for each factor is $(\Delta t_\varepsilon)_\varepsilon = (t_f - t_i)/(N_\varepsilon)_\varepsilon$. Substitution Eq.(1.1.5) into Eq.(1.1.6) by taking the limit $N_\varepsilon \to \infty, \varepsilon \in (0, 1]$ gives:

$$( \langle x_f, t_f | x_i, t_i \rangle )_\varepsilon = \left( \lim_{N_\varepsilon \to \infty} \int \prod_{i=1}^{N_\varepsilon - 1} dx_i \exp\left[ \frac{i}{\epsilon} \Delta t_\varepsilon \sum_{j=0}^{N_\varepsilon - 1} \mathcal{L}_\varepsilon(\dot{x}(t_j), x(t_j), t_j) \right] \right)_\varepsilon . \qquad (1.1.8)$$

Here $t_0 = t_i$. The exponent becomes to Colombeau time-integral of the generalized Lagrangian, namely the generalized action $(S_\varepsilon[x(t), t_i, t_f])_\varepsilon$ for each path. A change of variable: $x(t) \to x(t) + \delta x(t), \delta x(t_f) = \delta x(t_i) = 0$ does not change the result of the Feynman-Colombeau path integral, and we find

$$\left( \int\limits_{\substack{x(t_f) = x_f \\ x(t_i) = x_i}} D[x(t)] \exp\left\{ \frac{i}{\hbar} S_\varepsilon[x(t) + \delta x(t), t_i, t_f] \right\} \right)_\varepsilon =$$

$$\int\limits_{\substack{x(t_f) = x_f \\ x(t_i) = x_i}} D[x(t)] \exp\left\{ \frac{i}{\hbar} \left( S_\varepsilon[x(t) + \delta x(t), t_i, t_f] \right)_\varepsilon \right\} =$$

$$\left( \int\limits_{\substack{x(t_f) = x_f \\ x(t_i) = x_i}} D[x(t)] \exp\left\{ \frac{i}{\hbar} S_\varepsilon[x(t), t_i, t_f] \right\} \right)_\varepsilon =$$

$$\int\limits_{\substack{x(t_f) = x_f \\ x(t_i) = x_i}} D[x(t)] \exp\left\{ \frac{i}{\hbar} \left( S_\varepsilon[x(t), t_i, t_f] \right)_\varepsilon \right\} .$$

$(1.1.9)$

Therefore

$$\int\limits_{\substack{x(t_f)=x_f \\ x(t_i)=x_i}} D[x(t)] \exp\left\{\frac{i}{\hbar}\left(S_\varepsilon[x(t)+\delta x(t), t_i, t_f]\right)_\varepsilon\right\} -$$

$$\int\limits_{\substack{x(t_f)=x_f \\ x(t_i)=x_i}} D[x(t)] \exp\left\{\frac{i}{\hbar}\left(S_\varepsilon[x(t), t_i, t_f]\right)_\varepsilon\right\} = \tag{1.1.10}$$

$$\int\limits_{\substack{x(t_f)=x_f \\ x(t_i)=x_i}} D[x(t)] \frac{i}{\hbar}\left((\delta S_\varepsilon[x(t)])_\varepsilon\right) \exp\left\{\frac{i}{\hbar}\left(S_\varepsilon[x(t)]\right)_\varepsilon\right\} = 0.$$

From Eq.(1.1.10) we obtai

$$\int\limits_{\substack{x(t_f)=x_f \\ x(t_i)=x_i}} D[x(t)] \frac{i}{\hbar}\left((\delta S_\varepsilon[x(t), t_i, t_f])_\varepsilon\right) \times$$

$$\exp\left\{\frac{i}{\hbar}\left(S_\varepsilon[x(t), t_i, t_f]\right)_\varepsilon\right\} = \tag{1.1.11}$$

$$\int\limits_{\substack{x(t_f)=x_f \\ x(t_i)=x_i}} D[x(t)] \frac{i}{\hbar}\left[\int\limits_{t_i}^{t_f}\left(\frac{\partial \mathcal{L}_\varepsilon}{\partial x} - \frac{d}{dt}\frac{\partial \mathcal{L}_\varepsilon}{\partial \dot{x}}\right)_\varepsilon \delta x(t) dt\right] \times$$

$$\exp\left\{\frac{i}{\hbar}\left(S_\varepsilon[x(t), t_i, t_f]\right)_\varepsilon\right\} = 0.$$

Because $\delta x(t)$ is an arbitrary change of variable, the expression must be zero at all $t$ independently:

$$\int_{\substack{x(t_f)=x_f \\ x(t_i)=x_i}} D[x(t)] \left[ \left( \frac{\partial \mathcal{L}_\varepsilon}{\partial x} - \frac{d}{dt} \frac{\partial \mathcal{L}_\varepsilon}{\partial \dot{x}} \right)_\varepsilon \right] \exp \left\{ \frac{i}{\hbar} \left( S_\varepsilon[x(t), t_i, t_f] \right)_\varepsilon \right\} =$$

$$\int_{\substack{x(t_f)=x_f \\ x(t_i)=x_i}} D[x(t)] \left( \frac{\partial \mathcal{L}_\varepsilon}{\partial x} \right)_\varepsilon \exp \left\{ \frac{i}{\hbar} \left( S_\varepsilon[x(t), t_i, t_f] \right)_\varepsilon \right\} -$$

$$\int_{\substack{x(t_f)=x_f \\ x(t_i)=x_i}} D[x(t)] \left( \left( \frac{d}{dt} \frac{\partial \mathcal{L}_\varepsilon}{\partial \dot{x}} \right)_\varepsilon \right) \exp \left\{ \frac{i}{\hbar} \left( S_\varepsilon[x(t), t_i, t_f] \right)_\varepsilon \right\} = \qquad (1.1.12)$$

$$\left( \left\langle \frac{\partial \mathcal{L}_\varepsilon}{\partial x} \right\rangle \right)_\varepsilon - \left( \left\langle \frac{d}{dt} \frac{\partial \mathcal{L}_\varepsilon}{\partial \dot{x}} \right\rangle \right)_\varepsilon =$$

$$\left( \left\langle \frac{\partial \mathcal{L}_\varepsilon}{\partial x} \right\rangle \right)_\varepsilon - \frac{d}{dt} \left( \left\langle \frac{\partial \mathcal{L}_\varepsilon}{\partial \dot{x}} \right\rangle \right)_\varepsilon = 0.$$

Therefore, the Euler–Lagrange equation must hold as an generalized expectation value,nothing but the Ehrenfest's theorem.

## Apendix I.2.Colombeau-Schrödinger Equation from Feynman-Colombeau Path Integral.

Here we derive Colombeau-Schrödinger equation from the Feynman-Colombeau path integral.Let us first see that the momentum is given by a derivative. Starting from the Feynman-Colombeau path integral

$$\left( \int\limits_{\substack{x(t_f)=x_f \\ x(t_i)=x_i}} D[x(t)] \exp\left\{ \frac{i}{\hbar} S_\varepsilon[x(t), t_i, t_f] \right\} \right)_\varepsilon =$$

$$\int\limits_{\substack{x(t_f)=x_f \\ x(t_i)=x_i}} D[x(t)] \exp\left\{ \frac{i}{\hbar} \left( S_\varepsilon[x(t), t_i, t_f] \right)_\varepsilon \right\}.$$

(1.2.1)

we shift the trajectory $x(t)$ by a small amount $x(t) + \delta x(t)$ with the boundary condition that $x_i$ is held fixed ( $\delta x(t_i) = 0$) while $x_f$ is varied ( $x(t_f) \neq 0$). Under this variation, the propagator changes by

$$\left( \langle x_f + \delta x(t_f), t_f | x_i, t_i \rangle_\varepsilon \right)_\varepsilon - \left( \langle x_f, t_f | x_i, t_i \rangle_\varepsilon \right)_\varepsilon = \left( \left( \frac{\partial}{\partial x_f} \langle x_f, t_f | x_i, t_i \rangle_\varepsilon \right)_\varepsilon \right) (\delta x(t_f))$$

(1.2.2)

On the other hand, the Feynman-Colombeau path integral changes by

$$\int\limits_{\substack{x(t_f)=x_f+\delta x(t_f) \\ x(t_i)=x_i}} D[x(t)] \exp\left\{ \frac{i}{\hbar} \left( S_\varepsilon[x(t) + \delta x(t), t_i, t_f + \delta x(t_f)] \right)_\varepsilon \right\} -$$

$$\int\limits_{\substack{x(t_f)=x_f \\ x(t_i)=x_i}} D[x(t)] \exp\left\{ \frac{i}{\hbar} \left( S_\varepsilon[x(t), t_i, t_f] \right)_\varepsilon \right\} =$$

(1.2.3)

$$\frac{i}{\hbar} \int\limits_{\substack{x(t_f)=x_f \\ x(t_i)=x_i}} D[x(t)] \left( (\delta S_\varepsilon[x(t), t_i, t_f])_\varepsilon \right) \exp\left\{ \frac{i}{\hbar} \left( S_\varepsilon[x(t), t_i, t_f] \right)_\varepsilon \right\}$$

and the action changes by

$$(\delta S_\varepsilon[x(t), t_i, t_f])_\varepsilon = (S_\varepsilon[x(t) + \delta x(t), t_i, t_f])_\varepsilon - (S_\varepsilon[x(t), t_i, t_f])_\varepsilon =$$

$$\int_{t_i}^{t_f} dt \left( \frac{\partial \mathcal{L}_\varepsilon}{\partial x} \delta x(t) + \frac{\partial \mathcal{L}_\varepsilon}{\partial \dot{x}} \delta \dot{x}(t) \right)_\varepsilon =$$

$$\left( \frac{\partial \mathcal{L}_\varepsilon}{\partial \dot{x}} \delta x(t) \right)_\varepsilon \bigg|_{t_i}^{t_f} + \int_{t_i}^{t_f} \left( \left( \frac{\partial \mathcal{L}_\varepsilon}{\partial x} - \frac{\partial \mathcal{L}_\varepsilon}{\partial \dot{x}} \right)_\varepsilon \right) \delta x(t) dt.$$

(1.2.4)

The last terms vanishes because of the equation of motion (which holds as an expectation value, as we saw in the previous **Apendix I.1.1**), and we are left with

$$(\delta S_\varepsilon[x(t), t_i, t_f])_\varepsilon = \left( \left( \frac{\partial \mathcal{L}_\varepsilon}{\partial \dot{x}} \right)_\varepsilon \right) \delta x(t_f) = ((p_\varepsilon(t_f))_\varepsilon) \delta x(t_f).$$

(1.2.5)

By putting them together, and dropping $\delta x(t_f)$, we obtain

$$\left( \frac{\partial}{\partial x_f} \langle x_f, t_f | x_i, t_i \rangle_\varepsilon \right)_\varepsilon = \frac{i}{\hbar} \int_{\substack{x(t_f)=x_f \\ x(t_i)=x_i}} D[x(t)] ((p_\varepsilon(t_f))_\varepsilon) \exp\left\{ \frac{i}{\hbar} (S_\varepsilon[x(t), t_i, t_f])_\varepsilon \right\}$$

(1.2.6)

This is precisely how the momentum operator is represented in the position space. Now the Colombeau-Schrödinger equation can be derived by taking a variation with respect to $t_f$. Note that

$$\left( \frac{\partial S_\varepsilon}{\partial t_f} \right)_\varepsilon = -(H_\varepsilon(t_f))_\varepsilon$$

(1.2.7)

by using the equation of motion. Therefore,

$$\left( \frac{\partial}{\partial t_f} \langle x_f, t_f | x_i, t_i \rangle_\varepsilon \right)_\varepsilon = -\frac{i}{\hbar} \int\limits_{\substack{x(t_f)=x_f \\ x(t_i)=x_i}} D[x(t)] ((H_\varepsilon(t_f))_\varepsilon) \exp\left\{ \frac{i}{\hbar} (S_\varepsilon[x(t), t_i, t_f])_\varepsilon \right\} \qquad (1.2.8)$$

If

$$(H_\varepsilon)_\varepsilon = \frac{p^2}{2m} + (V_\varepsilon(x, t))_\varepsilon \qquad (1.2.9)$$

the momentum can be rewritten using Eq. (1.2.6), and we recover the Colombeau-Schrödinger equation:

$$i\hbar \frac{\partial}{\partial t} (\langle x, t | x_i, t_i \rangle_\varepsilon) = \left[ \frac{1}{2m} \left( \frac{\hbar}{i} \frac{\partial}{\partial x} \right)^2 + (V_\varepsilon(x, t))_\varepsilon \right] (\langle x, t | x_i, t_i \rangle_\varepsilon). \qquad (1.2.10)$$

# Apendix I.3.Colombeau-Feynman oscillatory path integral.

Let us calculate a simple Colombeau-Feynman oscillatory path integral:

$$\left(\langle x_f, t_f | x_i, t_i \rangle_\varepsilon\right)_\varepsilon \triangleq \left( \int\limits_{\substack{x(t_f)=x_f \\ x(t_i)=x_i}} D[x(t),\varepsilon] \exp\left[ \frac{i}{\hbar_\varepsilon} \int_{t_i}^{t_f} \frac{\mathbf{m}}{2} \left[ \frac{d_{\widetilde{\mathbb{R}},\varepsilon}}{d_{\widetilde{\mathbb{R}},\varepsilon} t} x(t) \right]^2 dt \right] \right)_{\varepsilon\in(0,1]} =$$

$$\left( \int\limits_{\substack{x(t_f)=x_f \\ x(t_i)=x_i}} D[x(t),\varepsilon] \exp\left[ \frac{i}{\hbar_\varepsilon} \int_{t_i}^{t_f} \frac{\mathbf{m}}{2} \left[ \dot{x}_{\widetilde{\mathbb{R}},\varepsilon}(t) \right]^2 dt \right] \right)_{\varepsilon\in(0,1]} =$$

$$\left( \int\limits_{\substack{x(t_f)=x_f \\ x(t_i)=x_i}} D[x(t),\varepsilon] \exp\left[ \frac{i}{\hbar_\varepsilon} \int_{t_i}^{t_f} \frac{\mathbf{m}}{2} \dot{x}_\varepsilon^2(t) dt \right] \right)_{\varepsilon\in(0,1]} = \tag{1.3.1}$$

$$\left( \int\limits_{\substack{x(t_f)=x_f \\ x(t_i)=x_i}} D[x(t),\varepsilon] \exp\left[ \frac{i}{\hbar_\varepsilon} \int_{t_i}^{t_f} \frac{\mathbf{m}}{2} \dot{x}^2(t) dt \right] \right)_{\varepsilon\in(0,1]},$$

$$\mathbf{m} = m_1 + im_2, m_2 \geq 0,$$

$$x(t) \in \mathbf{C}[t_i, t_f],$$

$$\left(x_\varepsilon(t)\right)_\varepsilon \in G(\mathbb{R}).$$

Here operator $\dfrac{d_{\widetilde{\mathbb{R}},\varepsilon}}{d_{\widetilde{\mathbb{R}},\varepsilon} t} : \mathbf{C}[t_i, t_f] \mapsto G(\mathbb{R})$ and Colombeau-Feynman measure $\left(D[x(t),\varepsilon]\right)_\varepsilon$ will be fixed later.

**Remark**.**1**.**3**.**1**. We let for short abbreviation $\dfrac{d_{\widetilde{\mathbb{R}},\varepsilon}}{d_{\widetilde{\mathbb{R}},\varepsilon} t} x(t) = \dot{x}_{\widetilde{\mathbb{R}},\varepsilon}(t) = \dot{x}_\varepsilon(t)$ or even simply $\dot{x}(t)$.

The classical path is

$$x_{\mathbf{cl}}(t) = x_i + \frac{x_f - x_i}{t_f - t_i}(t - t_i).$$

(1.3.2)

Assume that $x(t) \in \mathbf{C}[t_i, t_f]$. Hence one can write

$$x(t) =_{\widetilde{\mathbb{R}}} x_{\mathbf{cl}}(t) + \big(\delta_\varepsilon(t)\big)_\varepsilon, \varepsilon \in (0, 1],$$

(1.3.3)

where the left-over piece must vanish at the initial and the final time. Therefore, we can expand $\big(\delta_\varepsilon(t)\big)_\varepsilon$ without loss of generality, in Colombeau-Fourier series such that

$$\big(\delta_\varepsilon(t)\big)_\varepsilon = \left( \sum_{n=1}^{m(\varepsilon)} a_{n,\varepsilon} \sin\!\left( \frac{n\pi(t - t_i)}{t_f - t_i} \right) \right)_{\varepsilon \in (0,1]}$$

(1.3.4)

with arbitraly function $m(\varepsilon)$ such that $m(\varepsilon) \to \infty$ iff $\varepsilon \to 0$. Then the integral over all paths $x(t) \in \mathbf{C}[t_i, t_f]$ such that $a_{n,\varepsilon} = 0$ for $n > m(\varepsilon)$, can then be viewed as integrals over all $a_{n,\varepsilon}$, with $n \le m(\varepsilon)$, i.e.

$$\int D[x(t), \varepsilon](\cdot) = \int D[\{a_{n,\varepsilon}\}, \varepsilon](\cdot) = N_\varepsilon \int_{-\infty}^{\infty} \prod_{n=1}^{m(\varepsilon)} da_n(\cdot)$$

(1.3.5)

and therefore the integral over all paths $x(t) \in \mathbf{C}[t_i, t_f]$ can be viewed as Colombeau over all $a_{n,\varepsilon \in (0,1]}$, i.e.

$$\left( \int D[x(t), \varepsilon](\cdot) \right)_\varepsilon = \left( \int D[\{a_n\}, \varepsilon](\cdot) \right)_\varepsilon = \big( (N_\varepsilon)_\varepsilon \big) \left( \int_{-\infty}^{\infty} \prod_{n=1}^{m(\varepsilon)} da_n(\cdot) \right)_\varepsilon$$

(1.3.6)

The overall normalization factor $(N_\varepsilon)_\varepsilon$ that depends on **m** and $t_f - t_i$ will be fixed later.

The action, on the other hand, can be calculated. The starting point is

$$\dot{x}_\varepsilon(t) = \dot{x}_{\mathbf{cl}}(t) + \dot{\tilde{\delta}}_{\widetilde{\mathbb{R}},\varepsilon}(t) = \frac{x_f - x_i}{t_f - t_i} + \sum_{n=1}^{m(\varepsilon)} a_{n,\varepsilon} \frac{n\pi}{t_f - t_i} \cos\left(\frac{n\pi(t - t_i)}{t_f - t_i}\right). \qquad (1.3.7)$$

Because different modes are orthogonal upon $t$-integral, the action $(S_\varepsilon)_\varepsilon$ in coordinates $\{a_{n,\varepsilon}\}_{n=1}^{m(\varepsilon)}$ is

$$(S_\varepsilon)_\varepsilon = \frac{\mathbf{m}}{2}\left(\int_{t_i}^{t_f} \dot{x}_{\widetilde{\mathbb{R}},\varepsilon}^2(t)\,dt\right)_\varepsilon = \frac{\mathbf{m}}{2}\frac{(x_f - x_i)^2}{t_f - t_i} + \frac{\mathbf{m}}{2}\left(\sum_{n=1}^{m(\varepsilon)}\frac{1}{2}\frac{n^2\pi^2}{t_f - t_i}a_{n,\varepsilon}^2\right)_\varepsilon. \qquad (1.3.8)$$

The first term is nothing but the classical action and thus the path integral (1.3.1) reduces to Colombeau oscillatory integral:

$$\left(\int_{\substack{x(t_f)=x_f \\ x(t_i)=x_i}} D[x(t),\varepsilon]\exp\left[\frac{i}{\hbar_\varepsilon}\int_{t_i}^{t_f}\frac{\mathbf{m}}{2}\dot{x}_{\widetilde{\mathbb{R}},\varepsilon}^2(t)\,dt\right]\right)_\varepsilon =$$

$$((N_\varepsilon)_\varepsilon)\exp\left[\frac{i}{(\hbar_\varepsilon)_\varepsilon}\frac{\mathbf{m}}{2}\frac{(x_f - x_i)^2}{t_f - t_i}\right]\times$$

$$\left(\int_{-\infty}^{\infty}\prod_{n=1}^{m(\varepsilon)} da_{n,\varepsilon}\exp\left[\left(\frac{i}{\hbar_\varepsilon}\frac{\mathbf{m}}{4}\sum_{n=1}^{m(\varepsilon)}\frac{1}{2}\frac{n^2\pi^2}{t_f - t_i}\right)a_{n,\varepsilon}^2\right]\right)_\varepsilon.$$

$$(1.3.7)$$

Thus we obtain

$$\left( \langle x_f, t_f | x_i, t_i \rangle_\varepsilon \right)_\varepsilon = \left( (N_\varepsilon)_\varepsilon \right) \exp\left[ \frac{i}{(\hbar_\varepsilon)_\varepsilon} \frac{\mathbf{m}}{2} \frac{(x_f - x_i)^2}{t_f - t_i} \right] \times$$

$$\left( \prod_{n=1}^{m(\varepsilon)} \left( -\frac{i}{\pi \hbar_\varepsilon} \frac{\mathbf{m}}{2} \frac{1}{2} \frac{n^2 \pi^2}{t_f - t_i} \right)^{-1/2} \right)_\varepsilon . \tag{1.3.9}$$

Therefore, the result is simply

$$\left( \langle x_f, t_f | x_i, t_i \rangle_\varepsilon \right)_\varepsilon = (f_\varepsilon(t_f - t_i))_\varepsilon \exp\left[ \frac{i}{(\hbar_\varepsilon)_\varepsilon} \frac{\mathbf{m}}{2} \frac{(x_f - x_i)^2}{t_f - t_i} \right] \tag{1.3.10}$$

The normalization Colombeau constant $(f_\varepsilon(t_f - t_i))_\varepsilon$ can depend only on the time interval $t_f - t_i$ and is determined by the requirement that

$$\left( \int dx (\langle x_f, t_f | x, t \rangle_\varepsilon) \left( \langle x, t | x_i, t_i \rangle_\varepsilon \right) \right)_\varepsilon = (\langle x_f, t_f | x_i, t_i \rangle_\varepsilon)_\varepsilon \tag{1.3.11}$$

And hence one obtain

$$((f_\varepsilon(t_f - t))_\varepsilon)((f_\varepsilon(t - t_i))_\varepsilon) \sqrt{\frac{2\pi i (\hbar_\varepsilon)_\varepsilon (t_f - t)(t - t_i)}{\mathbf{m}(t_f - t_i)}} = (f_\varepsilon(t_f - t_i))_\varepsilon. \tag{1.3.12}$$

Therefore we find

$$(f_\varepsilon(t))_\varepsilon = \sqrt{\frac{\mathbf{m}}{2\pi i ((\hbar_\varepsilon)_\varepsilon) t}} . \tag{1.3.13}$$

Finally we obtain

$$\left( \langle x_f, t_f | x_i, t_i \rangle_\varepsilon \right)_\varepsilon = \sqrt{\frac{\mathbf{m}}{2\pi i ((\hbar_\varepsilon)_\varepsilon)(t_f - t_i)}} \, \exp\left[ \frac{i}{((\hbar_\varepsilon)_\varepsilon)} \, \frac{\mathbf{m}}{2} \, \frac{(x_f - x_i)^2}{t_f - t_i} \right]. \quad (1.3.14)$$

In order to obtain this result directly from the Colombeau-Feynman path integral, we should have chosen the normalization of the complex Colombeau-Feynman measure $(D[x(t), \varepsilon])_\varepsilon$ to be

$$(D[x(t), \varepsilon])_\varepsilon =$$

$$\sqrt{\frac{\mathbf{m}}{2\pi T i ((\hbar_\varepsilon)_\varepsilon)}} \left( \left( \prod_{n=1}^{m(\varepsilon)} \left( \frac{1}{i\pi \hbar_\varepsilon} \, \frac{\mathbf{m}}{2} \, \frac{1}{2} \, \frac{n^2 \pi^2}{T} \right)^{1/2} \right)_\varepsilon \right) \left( \left( \int_{-\infty}^{\infty} \prod_{n=1}^{m(\varepsilon)} da_{n,\varepsilon} \right)_\varepsilon \right), \quad (1.3.15)$$

$$T = t_f - t_i.$$

**Definition 1.3.1. (i)** If $m_2 = 0$, we define positive Colombeau-Feynman measure $(D^+[x(t), \varepsilon])_\varepsilon$ by formula

$$(D^+[x(t), \varepsilon])_\varepsilon =$$

$$\left( \left( \prod_{n=1}^{8m(\varepsilon)} \left( \frac{1}{i\pi \hbar_\varepsilon} \, \frac{m_1}{2} \, \frac{1}{2} \, \frac{n^2 \pi^2}{T} \right)^{1/2} \right)_\varepsilon \right) \left( \left( \int_{-\infty}^{\infty} \prod_{n=1}^{8m(\varepsilon)} da_{n,\varepsilon}(\bullet) \right)_\varepsilon \right) = \quad (1.3.16.a)$$

$$\left( \left( \prod_{n=1}^{8m(\varepsilon)} \left( \frac{1}{\pi \hbar_\varepsilon} \, \frac{m_1}{2} \, \frac{1}{2} \, \frac{n^2 \pi^2}{T} \right)^{1/2} \right)_\varepsilon \right) \left( \left( \int_{-\infty}^{\infty} \prod_{n=1}^{8m(\varepsilon)} da_{n,\varepsilon}(\bullet) \right)_\varepsilon \right)$$

**(ii)** If $m_2 > 0$ we define modulus $(|D[x(t), \varepsilon]|)_\varepsilon$ of the complex Colombeau-Feynman measure $(D[x(t), \varepsilon])_\varepsilon$ by formula

$$\left(\left|D[x(t),\varepsilon]\right|\right)_\varepsilon =$$

$$\left(\left(\left|\prod_{n=1}^{m(\varepsilon)}\left(\frac{1}{i\pi\hbar_\varepsilon}\frac{\mathbf{m}}{2}\frac{1}{2}\frac{n^2\pi^2}{T}\right)^{1/2}\right|\right)_\varepsilon\right)\left(\left(\int_{-\infty}^{\infty}\prod_{n=1}^{m(\varepsilon)}da_{n,\varepsilon}(\bullet)\right)_\varepsilon\right) = \qquad (1.3.16.b)$$

$$\left(\left(\prod_{n=1}^{m(\varepsilon)}\left(\frac{1}{\pi\hbar_\varepsilon}\frac{\sqrt{m_1^2+m_2^2}}{2}\frac{1}{2}\frac{n^2\pi^2}{T}\right)^{1/2}\right)_\varepsilon\right)\left(\left(\int_{-\infty}^{\infty}\prod_{n=1}^{m(\varepsilon)}da_{n,\varepsilon}(\bullet)\right)_\varepsilon\right)$$

Let us calculate now Colombeau oscillatory path integral

$$\left(\langle x_f, t_f | x_i, t_i \rangle_\varepsilon\right)_\varepsilon = \left( \int\limits_{\substack{x(T)=x_f \\ x(0)=x_i}} D[x(t), \varepsilon] \times \right.$$

$$\exp\left\{ \frac{i}{\hbar_\varepsilon} \left( \int\limits_{t_i}^{t_f} dt \left[ \frac{\mathbf{m}}{2} \left( \frac{d_{\widetilde{\mathbb{R}}, \varepsilon} x(t)}{d_{\widetilde{\mathbb{R}}, \varepsilon} t} \right)_\varepsilon^2 - \frac{\omega_1^2}{\hbar_\varepsilon} \frac{\mathbf{m}}{2} x^2(t) \right] \right) - \right.$$

$$\left. - \frac{\omega_2^2}{\hbar_\varepsilon} \int\limits_0^T dt \left[ \frac{\mathbf{m}}{2} x^2(t) \right] \right)_\varepsilon = \tag{1.3.17}$$

$$\left( \int\limits_{\substack{x(T)=x_f \\ x(0)=x_i}} D[x(t), \varepsilon] \times \right.$$

$$\left. \exp\left\{ \frac{i}{\hbar_\varepsilon} \left( \int\limits_{t_i}^{t_f} dt \left[ \frac{\mathbf{m}}{2} \left( \frac{d_{\widetilde{\mathbb{R}}, \varepsilon} x(t)}{d_{\widetilde{\mathbb{R}}, \varepsilon} t} \right)_\varepsilon^2 - \frac{\varpi^2}{\hbar_\varepsilon} \frac{\mathbf{m}}{2} x^2(t) \right] \right) \right\} \right)_\varepsilon,$$

$$\varpi^2 = \omega_1^2 - i\omega_2^2.$$

The classical path is

$$x_{\mathbf{cl}}(t) = x_i \frac{\sin[\varpi(t_f - t)]}{\sin[\varpi(t_f - t_i)]} + x_f \frac{\sin[\varpi(t - t_i)]}{\sin[\varpi(t_f - t_i)]}. \tag{1.3.18}$$

The action $S_{\mathbf{cl}}$ along the classical path is

$$S_{\mathbf{cl}} = \int\limits_{t_i}^{t_f} dt \left[ \frac{\mathbf{m}}{2} \left( \frac{dx_{\mathbf{cl}}(t)}{dt} \right)^2 - \varpi^2 \frac{\mathbf{m}}{2} x_{\mathbf{cl}}^2(t) \right] =$$

(1.3.19)

$$\frac{\mathbf{m}\varpi}{2} \frac{(x_i^2 + x_f^2)\cos[\varpi(t_f - t_i)] - 2x_i x_f}{\sin[\varpi(t_f - t_i)]}$$

One can write

$$x(t) =_{\widetilde{\mathbb{R}}} x_{\mathbf{cl}}(t) + (\delta_\varepsilon(t))_\varepsilon,$$

(1.3.20)

where the left-over piece must vanish at the initial and the final time. Therefore, we can expand $(\delta_\varepsilon(t))_\varepsilon$ without loss of generality, in Fourier-Colombeau series such that

$$(\delta_\varepsilon(t))_\varepsilon = \left( \sum_{n=1}^{m(\varepsilon)} a_{n,\varepsilon} \sin\left( \frac{n\pi(t - t_i)}{t_f - t_i} \right) \right)_\varepsilon$$

(1.3.21)

with arbitraly function $m(\varepsilon)$ such that $m(\varepsilon) \to \infty$ iff $\varepsilon \to 0$. Then the integral over all paths $x(t) \in \mathbf{C}[t_i, t_f]$ such that $a_{n,\varepsilon} = 0$ for $n > m(\varepsilon)$, can then be viewed as integrals over all $a_{n,\varepsilon}$, with $n \leq m(\varepsilon)$, i.e.

$$\int D[x(t), \varepsilon](\cdot) = \int D[\{a_{n_\varepsilon}\}, \widetilde{\varepsilon}](\cdot) = N_\varepsilon \int\limits_{-\infty}^{\infty} \prod_{n=1}^{m(\varepsilon)} da_{n,\varepsilon}(\cdot)$$

(1.3.22)

and therefore the integral over all paths $x(t) \in \mathbf{C}[t_i, t_f]$ can be viewed as Colombeau over all $a_{n,\varepsilon}$, i.e.

$$\left( \int D[x(t), \varepsilon](\cdot) \right)_{\varepsilon} = \left( \int D[\{a_{n,\varepsilon}\}, \varepsilon](\cdot) \right)_{\varepsilon} = \left( (N_{\varepsilon})_{\varepsilon} \left( \int\limits_{-\infty}^{\infty} \prod_{n=1}^{m(\varepsilon)} da_{n,\varepsilon}(\cdot) \right) \right)_{\varepsilon} \quad (1.3.23)$$

The overall normalization factor $(f_{\varepsilon} t_f - t_i)_{\varepsilon}$ that depends on $\mathbf{m}$ and $t_f - t_i$ will be fixed later. The action, on the other hand, can be calculated. The starting point is

$$\left( \dot{x}_{\widetilde{\mathbb{R}},\varepsilon}(t) \right)_{\varepsilon} = \dot{x}_{\mathbf{cl}}(t) + \dot{\delta}_{\varepsilon}(t) = \dot{x}_{\mathbf{cl}}(t) + \left( \sum_{n=1}^{m(\varepsilon)} a_{n,\varepsilon} \frac{n\pi}{t_f - t_i} \cos\left( \frac{n\pi(t - t_i)}{t_f - t_i} \right) \right)_{\varepsilon}. \quad (1.3.24)$$

Because different modes are orthogonal upon $t$-integral, the action $(S_{\varepsilon})_{\varepsilon}$ in coordinates $\{a_{n,\varepsilon}\}_{n=1}^{m(\varepsilon)}$ is

$$(S_{\varepsilon})_{\varepsilon} = S_{\mathbf{cl}} + \frac{\mathbf{m}}{4} \left( \sum_{n=1}^{m(\varepsilon)} \left[ \frac{n^2 \pi^2}{t_f - t_i} - \varpi^2 (t_f - t_i) \right] a_{n,\varepsilon}^2 \right)_{\varepsilon}. \quad (1.3.25)$$

Therefore the Colombeau-Feynman path integral (1.3.17) is Colombeau-Fresnel integral over $\{a_{n,\varepsilon}\}_{n=1}^{\infty}$. Using the Colombeau-Feynman measure given by Eq.(1.3.15)

we obtain

$$\left( \langle x_f, t_f | x_i, t_i \rangle_{\varepsilon} \right)_{\varepsilon} =$$

$$\left[ \exp\left( \frac{i S_{\mathbf{cl}}}{(\hbar_{\varepsilon})_{\varepsilon}} \right) \right] \times \sqrt{\frac{\mathbf{m}}{2\pi T i (\hbar_{\varepsilon})_{\varepsilon}}} \left( \left( \prod_{n=1}^{m(\varepsilon)} \left( \frac{1}{i\pi\hbar_{\varepsilon}} \frac{\mathbf{m}}{2} \frac{1}{2} \frac{n^2\pi^2}{T} \right)^{1/2} \right)_{\varepsilon} \right) \times \quad (1.3.26)$$

$$\left( \int\limits_{-\infty}^{\infty} da_{n,\varepsilon} \exp\left[ \frac{i}{\hbar_{\varepsilon}} \frac{\mathbf{m}}{2} \left( \frac{n^2\pi^2}{t_f - t_i} - \varpi^2(t_f - t_i) \right) \frac{1}{2} a_{n,\varepsilon}^2 \right] \right)_{\varepsilon}.$$

Now we resort to the following Colombeau representation of the sine function

$$\left( \prod_{n=1}^{m(\varepsilon)} \left( 1 - \frac{x^2}{n^2} \right) \right)_{\varepsilon} =_{\widetilde{\mathbb{R}}} \frac{\sin \pi x}{\pi x}. \qquad (1.3.27)$$

Thus we find

$$\left( \langle x_f, t_f | x_i, t_i \rangle_{\varepsilon} \right) =_{\widetilde{\mathbb{C}}} \left[ \exp \left( \frac{iS_{\mathbf{cl}}}{(\hbar_{\varepsilon})} \right) \right] \times \sqrt{\frac{m\varpi}{2\pi i (\hbar_{\varepsilon}) \sin \varpi (t_f - t_i)}}. \qquad (1.3.28)$$

This agrees with the result from the more conventional method as given in literature.

**Definition 1.3.2.** Assume that $(V_{\varepsilon}(x))_{\varepsilon} \in G(\mathbb{R}^n)$. We define Colombeau-Feinman propagator $\left( \langle x_f, t_f | x_i, t_i \rangle_{\varepsilon} \right)_{\varepsilon \in (0,1]}$ by formula

$$\left( \langle x_f, t_f | x_i, t_i \rangle_{\varepsilon} \right)_{\varepsilon \in (0,1]} =$$

$$\left( \int_{\substack{q(t_f)=x_f \\ q(t_i)=x_i}} D[x(t), \varepsilon] \times \exp \left\{ \frac{i}{\hbar_{\varepsilon}} \left( \int_{t_i}^{t_f} dt \left[ \frac{\mathbf{m}}{2} \left( \frac{d_{\widetilde{\mathbb{R}}, \varepsilon} x(t)}{d_{\widetilde{\mathbb{R}}, \varepsilon} t} \right)^2 - V_{\delta}(x(t)) \right] \right) \right\} \right)_{\varepsilon \in (0,1]}. \qquad (1.3.29)$$

And in a more general case we define Colombeau-Feinmanpropagator $\left( \langle x_f, t_f | x_i, t_i \rangle_{\varepsilon, \delta} \right)_{\varepsilon, \delta \in (0,1]}$ by formula

$$\left( \langle x_f, t_f | x_i, t_i \rangle_{\varepsilon, \delta} \right)_{\varepsilon, \delta \in (0,1]} =$$

$$\left( \int_{\substack{q(t_f)=x_f \\ q(t_i)=x_i}} D[x(t), \varepsilon] \times \exp \left\{ \frac{i}{\hbar_{\varepsilon, \delta}} \left( \int_{t_i}^{t_f} dt \left[ \frac{\mathbf{m}}{2} \left( \frac{d_{\widetilde{\mathbb{R}}, \varepsilon} x(t)}{d_{\widetilde{\mathbb{R}}, \varepsilon} t} \right)^2 - V_\delta(x(t)) \right] \right) \right\} \right)_{\varepsilon, \delta \in (0,1]} . \tag{1.3.30}$$

**Remark**.1.3.2. We let for short the abbreviations

$$\left( \langle x_f, t_f | x_i, t_i \rangle_{\varepsilon} \right)_{\varepsilon (0,1]} =$$

$$\left( \int_{\substack{q(t_f)=x_f \\ q(t_i)=x_i}} D[x(t), \varepsilon] \times \exp \left\{ \frac{i}{\hbar_{\varepsilon}} \left( \int_{t_i}^{t_f} dt \left[ \frac{\mathbf{m}}{2} \left( \frac{d_{\widetilde{\mathbb{R}}} x(t)}{d_{\widetilde{\mathbb{R}}} t} \right)^2 - V_\varepsilon(x(t)) \right] \right) \right\} \right)_{\varepsilon \in (0,1]} \tag{1.3.31}$$

and

$$\left( \langle x_f, t_f | x_i, t_i \rangle_{\varepsilon, \delta} \right)_{\varepsilon, \delta \in (0,1]} =$$

$$\left( \int_{\substack{q(t_f)=x_f \\ q(t_i)=x_i}} D[x(t), \varepsilon] \times \exp \left\{ \frac{i}{\hbar_{\varepsilon}} \left( \int_{t_i}^{t_f} dt \left[ \frac{\mathbf{m}}{2} \left( \frac{d_{\widetilde{\mathbb{R}}} x(t)}{d_{\widetilde{\mathbb{R}}} t} \right)^2 - V_\delta(x(t)) \right] \right) \right\} \right)_{\varepsilon, \delta \in (0,1]} \tag{1.3.32}$$

or even simply

$$\left(\langle x_f, t_f | x_i, t_i \rangle_\varepsilon\right)_{\varepsilon \in (0,1]} =$$

$$\left(\int\limits_{\substack{q(t_f)=x_f \\ q(t_i)=x_i}} D[x(t), \varepsilon] \times \exp\left\{\frac{i}{\hbar_\varepsilon}\left(\int\limits_{t_i}^{t_f} dt\left[\frac{\mathbf{m}}{2}\left(\frac{dx(t)}{dt}\right)^2 - V_\varepsilon(x(t))\right]\right)\right\}\right)_{\varepsilon \in (0,1]} \qquad (1.3.33)$$

and

$$\left(\langle x_f, t_f | x_i, t_i \rangle_{\varepsilon,\delta}\right)_{\varepsilon,\delta \in (0,1]} =$$

$$\left(\int\limits_{\substack{q(t_f)=x_f \\ q(t_i)=x_i}} D[x(t), \varepsilon] \times \exp\left\{\frac{i}{\hbar_\varepsilon}\left(\int\limits_{t_i}^{t_f} dt\left[\frac{\mathbf{m}}{2}\left(\frac{dx(t)}{dt}\right)^2 - V_\delta(x(t))\right]\right)\right\}\right)_{\varepsilon,\delta \in (0,1]} \qquad (1.3.34)$$

.

# Apendix II.1.1.Upper and Lower Infinite-Dimensional Subintegrals.Colombeau-Feynman Measure.

Let $\mathbb{R}^\infty$ be the space

$$\mathbb{R}^\infty = \bigoplus_{n=1}^{\infty} \mathbb{R}^n \qquad (1.1.1)$$

**Definition 1.1.1.** For $f \in (\mathbb{R}^\infty)^*$ we define the infinite dimensional normed by sequence $\{c_n\}$, Feynman integral $\mathfrak{I}_\infty(f)$ as limit

$$\mathfrak{I}_\infty(f) = \int_{\mathbb{R}^\infty} f(\mathbf{x}) D[\mathbf{x}] = \lim_{n \to \infty} \int_{\mathbb{R}^n} f_n(x) c_n d^n x =$$

$$\lim_{n \to \infty} \int \ldots \int f_n(x_1, x_2, \ldots, x_n) c_n dx_1 dx_2 \ldots dx_n$$

$$(1.1.2)$$

iff limit in RHS of Eq.(1.1.2) exists and $f_n(x) = f(\mathbf{x}) \upharpoonright \mathbb{R}^n$.

**Remark 1.1.1.** We note that one can write Definition 1.1.1 using completely Colombeau form.

**Definition 1.1.2.** For $f \in (\mathbb{R}^\infty)^*$ we define the infinite dimensional integral $\mathfrak{I}_\infty(f)$ as Colombeau quantity such that

$$\Im_\infty(f) = \int_{\mathbb{R}^\infty} f(\mathbf{x}) D[\mathbf{x}] =_{\widetilde{\mathbb{C}}} \int_{\mathbb{R}^\infty} f(\mathbf{x}) D[\mathbf{x}, \varepsilon] =_{\widetilde{\mathbb{C}}}$$

$$\left( \int_{\mathbb{R}^{n_\varepsilon}} f_{n_\varepsilon}(x) c_{n_\varepsilon} d^{n_\varepsilon} x \right)_\varepsilon =_{\widetilde{\mathbb{C}}}$$

$$\left( \int \ldots \int f(x_1, x_2, \ldots, x_{n_\varepsilon}) c_{n_\varepsilon} dx_1 dx_2 \ldots dx_{n_\varepsilon} \right)_\varepsilon =_{\widetilde{\mathbb{C}}} \qquad (1.1.3)$$

$$((c_{n_\varepsilon})_\varepsilon) \left( \int_{\mathbb{R}^{n_\varepsilon}} f_{n_\varepsilon}(x) d^{n_\varepsilon} x \right)_\varepsilon =_{\widetilde{\mathbb{C}}}$$

$$((c_{n_\varepsilon})_\varepsilon) \left( \int \ldots \int f(x_1, x_2, \ldots, x_{n_\varepsilon}) dx_1 dx_2 \ldots dx_{n_\varepsilon} \right)_\varepsilon,$$

where $n_\varepsilon \to \infty$ if $\varepsilon \to 0$.

**Remark 1.1.2.** We note that Definition 1.1.2 work if iven limit in RHS of Eq.(1.1.2) dos not exists.

**Definition 1.1.3.** We define Colombeau-Feynman measure $(D[\mathbf{x}, \varepsilon])_\varepsilon$ by formula

$$(D[\mathbf{x}, \varepsilon])_\varepsilon = \left( c_{n_\varepsilon} dx_1 dx_2 \ldots dx_{n_\varepsilon} \right)_\varepsilon =$$

$$((c_{n_\varepsilon})_\varepsilon) \left( dx_1 dx_2 \ldots dx_{n_\varepsilon} \right)_\varepsilon. \qquad (1.1.4)$$

**Definition 1.1.4.** Let $\Sigma(\mathbb{R}^\infty)$ be the measurable space $(\mathbb{R}^\infty, \Sigma)$ and Let $\mathbf{M}(\mathbb{R}^\infty, \Sigma)$ be algebra with operations: (i) $\forall c \in \mathbb{C}, \forall A \in \Sigma : (c\mu)(A) = c\mu(A)$ (ii) $\forall A \in \Sigma : (\mu + \nu)(A) = \mu(A) + \nu(A)$. Set $E(\mathbf{M}(\mathbb{R}^\infty, \Sigma)) = (\mathbf{M}(\mathbb{R}^\infty, \Sigma))^I, I = (0, 1],$

$$E_M(\mathbf{M}(\mathbb{R}^\infty, \Sigma), \{c_{n_\varepsilon}\}) = \left\{ \left( c_{n_\varepsilon} \mu_{n_\varepsilon}(\bullet) \right)_\varepsilon \in E(\mathbf{M}(\mathbb{R}^\infty, \Sigma)) | (\forall n \in \mathbb{N})((\forall A \in \Sigma(\mathbb{R}^\infty))) \right.$$

$$\left. (\exists p \in \mathbb{N})[|c_{n_\varepsilon} \mu_{n_\varepsilon}(A)| = O(\varepsilon^{-p})] \right\},$$

$$\mu_{n_\varepsilon}(\bullet) \in \mathbf{M}(\mathbb{R}^{n_\varepsilon}, \Sigma_{n_\varepsilon}), \Sigma_{n_\varepsilon} = \mathbb{R}^{n_\varepsilon} \cap \Sigma(\mathbb{R}^\infty),$$

$$|c_{n_\varepsilon}| \to \infty \text{ if } \varepsilon \to 0,$$

(1.1.5)

$$N(\mathbf{M}(\mathbb{R}^\infty, \Sigma), \{c_{n_\varepsilon}\}) = \left\{ \left( c_{n_\varepsilon} \mu_{n_\varepsilon}(\bullet) \right)_\varepsilon \in E(\mathbf{M}(\mathbb{R}^\infty, \Sigma)) | (\forall n \in \mathbb{N})((\forall A \in \Sigma(\mathbb{R}^\infty))) \right.$$

$$\left. (\forall q \in \mathbb{N})[|c_{n_\varepsilon} \mu_{n_\varepsilon}(A)| = O(\varepsilon^q)] \right\},$$

(1.1.6)

$$\mu_{n_\varepsilon}(\bullet) \in \mathbf{M}(\mathbb{R}^{n_\varepsilon}, \Sigma_{n_\varepsilon}), \Sigma_{n_\varepsilon} = \mathbb{R}^{n_\varepsilon} \cap \Sigma(\mathbb{R}^\infty).$$

$$G(\mathbf{M}(\mathbb{R}^\infty, \Sigma), \{c_{n_\varepsilon}\}) = E_M(\mathbf{M}(\mathbb{R}^\infty, \Sigma), \{c_{n_\varepsilon}\})/N_M(\mathbf{M}(\mathbb{R}^\infty, \Sigma), \{c_{n_\varepsilon}\}).$$

(1.1.7)

Elements of $E_M(\mathbf{M}(\mathbb{R}^\infty, \Sigma), \{c_{n_\varepsilon}\})$ and $N_M(\mathbf{M}(\mathbb{R}^\infty, \Sigma), \{c_{n_\varepsilon}\})$ are called moderate, resp.negligible Normed by secuence $\{c_{n_\varepsilon}\}$ Colombeau-Feynman measures or simply **Colombeau-Feynman measures**.

Note that $E_M(\mathbf{M}(\mathbb{R}^\infty, \Sigma), \{c_{n_\varepsilon}\})$ is an algebra with operations (i) $\forall c \in \mathbb{C}, \forall A \in \Sigma$ : $((c\mu_\varepsilon)_\varepsilon)(A) = c(\mu_\varepsilon(A))_\varepsilon$ (ii) $\forall A \in \Sigma$ : $((\mu_\varepsilon + \nu_\varepsilon)_\varepsilon)(A) = (\mu_\varepsilon(A))_\varepsilon + (\nu_\varepsilon(A))_\varepsilon$.

It is the largest subalgebra of $E(\mathbf{M}(\mathbb{R}^\infty, \Sigma))$ in which $N_M(\mathbf{M}(\mathbb{R}^\infty, \Sigma), \{c_{n_\varepsilon}\})$ is ideal. Thus, $G(\mathbf{M}(\mathbb{R}^\infty, \Sigma), \{c_{n_\varepsilon}\})$ is an associative, commutative algebra.

If $(\mu_\varepsilon(\bullet))_\varepsilon \in E_M(\mathbf{M}(\mathbb{R}^\infty, \Sigma), \{c_{n_\varepsilon}\})$ is a representative of $\mu \in G(\mathbf{M}(\mathbb{R}^\infty, \Sigma), \{c_{n_\varepsilon}\})$, we write $\mu = \mathbf{cl}[(\mu_\varepsilon(\bullet))_\varepsilon]$ or simply $\mu = [(\mu_\varepsilon(\bullet))_\varepsilon]$.

**Definition 1.1.5.** (1) For $\mathbb{R}$-valued functional $f \in (\mathbb{R}^\infty)^*$ we define the infinite dimensional upper integral $\overline{\mathfrak{I}}_\infty(f)$ as upper limit

$$\overline{\mathfrak{J}}_\infty(f) = \overline{\int_{\mathbb{R}^\infty} f(\mathbf{x}) D[\mathbf{x}]} = \limsup_{n\to\infty} \int_{\mathbb{R}^n} f_n(x) c_n d^n x =$$

$$\limsup_{n\to\infty} \int \ldots \int f_n(x_1, x_2, \ldots, x_n) c_n dx_1 dx_2 \ldots dx_n.$$

(1.1.8)

(2) For $\mathbb{R}$-valued functional $f \in (\mathbb{R}^\infty)^*$ we define the infinite dimensional lower integral $\underline{\mathfrak{J}}_\infty(f)$ as lower limit

$$\underline{\mathfrak{J}}_\infty(f) = \underline{\int_{\mathbb{R}^\infty} f(\mathbf{x}) D[\mathbf{x}]} = \liminf_{n\to\infty} \int_{\mathbb{R}^\infty} f_n(x) c_n d^n x =$$

$$\liminf_{n\to\infty} \int \ldots \int f_n(x_1, x_2, \ldots, x_n) c_n dx_1 dx_2 \ldots dx_n.$$

(1.1.9)

**Proposition 1.1.1.** Assume that $f \in (\mathbb{R}^\infty)^*$ is an $\mathbb{R}$-valued functional, then $\underline{\mathfrak{J}}_\infty(f) \leq \overline{\mathfrak{J}}_\infty(f)$ and

$$\overline{\mathfrak{J}}_\infty(f) = \overline{\int_{\mathbb{R}^\infty} f(\mathbf{x}) D[\mathbf{x}]} = \lim_{k\to\infty} \left( \sup_{m\geq k} \int_{\mathbb{R}^m} f_m(x) c_m d^m x \right),$$

$$\underline{\mathfrak{J}}_\infty(f) = \underline{\int_{\mathbb{R}^\infty} f(\mathbf{x}) D[\mathbf{x}]} = \lim_{k\to\infty} \left( \inf_{m\geq k} \int_{\mathbb{R}^m} f_m(x) c_m d^m x \right).$$

(1.1.10)

**Proposition 1.1.2.** Suppose that $\overline{\mathfrak{J}}_\infty(f) = \underline{\mathfrak{J}}_\infty(f) = a_f < \infty$. Then $\mathfrak{J}_\infty(f) = a_f$.

**Proposition 1.1.3.** Assume that $f_1, f_2 \in (\mathbb{R}^\infty)^*$ and $f_1(\mathbf{x}) \leq f_2(\mathbf{x})$, then

$$\overline{\int}_{\mathbb{R}^\infty} [f_1(\mathbf{x})]D[\mathbf{x}] \leq \overline{\int}_{\mathbb{R}^\infty} [f_2(\mathbf{x})]D[\mathbf{x}],$$

$$(1.1.11)$$

$$\underline{\int}_{\mathbb{R}^\infty} [f_1(\mathbf{x})]D[\mathbf{x}] \leq \underline{\int}_{\mathbb{R}^\infty} [f_2(\mathbf{x})]D[\mathbf{x}].$$

**Proposition 1.1.4.** Assume that $f_1, f_2 \in (\mathbb{R}^\infty)^*$ is an $\mathbb{R}$-valued functionals then

$$\underline{\int}_{\mathbb{R}^\infty} [f_1(\mathbf{x})]D[\mathbf{x}] + \overline{\int}_{\mathbb{R}^\infty} [f_2(\mathbf{x})]D[\mathbf{x}] \leq \overline{\int}_{\mathbb{R}^\infty} [f_1(\mathbf{x}) + f_2(\mathbf{x})]D[\mathbf{x}] \leq$$

$$\leq \overline{\int}_{\mathbb{R}^\infty} [f_1(\mathbf{x})]D[\mathbf{x}] + \overline{\int}_{\mathbb{R}^\infty} [f_2(\mathbf{x})]D[\mathbf{x}],$$

$$(1.1.12)$$

$$\underline{\int}_{\mathbb{R}^\infty} [f_1(\mathbf{x})]D[\mathbf{x}] + \overline{\int}_{\mathbb{R}^\infty} [f_2(\mathbf{x})]D[\mathbf{x}] \geq \underline{\int}_{\mathbb{R}^\infty} [f_1(\mathbf{x}) + f_2(\mathbf{x})]D[\mathbf{x}] \geq$$

$$\underline{\int}_{\mathbb{R}^\infty} [f_1(\mathbf{x})]D[\mathbf{x}] + \underline{\int}_{\mathbb{R}^\infty} [f_2(\mathbf{x})]D[\mathbf{x}].$$

# Apendix II.1.2. Upper and Lower Infinite-Dimensional Colombeau-Feynman Subintegrals.

**Definition 1.2.1.** Let $\mathbb{R}^\infty$ be the space $\mathbb{R}^\infty = \bigoplus_{n=1}^\infty \mathbb{R}^n$. Set

$$E(\mathbb{R}^\infty) = (\mathbf{C}^\infty(\mathbb{R}^\infty))^I, I = (0,1],$$

$$E_M(\mathbb{R}^\infty) = \{(f_\varepsilon(\mathbf{x}))_\varepsilon \in E(\mathbb{R}^\infty) | (\forall n \in \mathbb{N})((\forall m \in \mathbb{N})$$

$$(\forall K \Subset \mathbb{R}^n)(\forall \alpha \in \mathbb{N}_0^m)(\exists p \in \mathbb{N})\left[\sup_{x \in K} |\partial^\alpha f_{n,\varepsilon}(x)| = O(\varepsilon^{-p})\right]\}, \qquad (1.2.1)$$

$$(f_{n,\varepsilon}(\mathbf{x}))_\varepsilon = (f_\varepsilon(\mathbf{x}) \upharpoonright \mathbb{R}^n)_\varepsilon,$$

$$N_M(\mathbb{R}^\infty) = \{(f_\varepsilon(\mathbf{x}))_\varepsilon \in E(\mathbb{R}^\infty) | (\forall n \in \mathbb{N})((\forall m \in \mathbb{N})$$

$$(\forall K \Subset \mathbb{R}^n)(\forall \alpha \in \mathbb{N}_0^m)(\forall q \in \mathbb{N})\left[\sup_{x \in K} |\partial^\alpha f_{n,\varepsilon}(x)| = O(\varepsilon^q)\right]\}, \qquad (1.2.2)$$

$$G(\mathbb{R}^\infty) = E_M(\mathbb{R}^\infty)/N_M(\mathbb{R}^\infty). \qquad (1.2.3)$$

Elements of $E_M(\mathbb{R}^\infty)$ and $N_M(\mathbb{R}^\infty)$ are called moderate, resp. negligible functionals. Note that $E_M(\mathbb{R}^\infty)$ is a differential algebra with pointwise operations. It is the largest differential subalgebra of $E(\mathbb{R}^\infty)$ in which $N_M(\mathbb{R}^\infty)$ is a differential ideal. Thus, $G(\mathbb{R}^\infty)$ is an associative, commutative differential algebra. If $(f_\varepsilon(\mathbf{x}))_\varepsilon \in E_M(\mathbb{R}^\infty)$ is a representative of $f \in G(\mathbb{R}^\infty)$, we write $f = \mathbf{cl}[(f_\varepsilon(\mathbf{x}))_\varepsilon]$ or simply $f = [(f_\varepsilon(\mathbf{x}))_\varepsilon]$.

**Definition 1**.**2**.**2**. For $(f_\varepsilon)_\varepsilon \in G_{\mathbf{C}^\infty(\mathbb{R}^\infty)} = G(\mathbb{R}^\infty)$ we define (by using Definition 1.1.1) the infinite dimensional Colombeau-Feinman integral $\mathfrak{I}_\infty(\mathbf{cl}[(f_\varepsilon)_\varepsilon]) \triangleq \mathbf{cl}[(\mathfrak{I}_\infty(f_\varepsilon))_\varepsilon]$ by formula

$$\left(\Im_\infty(f_\varepsilon)\right)_\varepsilon = \left(\int\limits_{\mathbb{R}^\infty} f_\varepsilon(\mathbf{x})D[\mathbf{x}]\right)_\varepsilon =$$

$$\left(\lim_{n\to\infty}\int\limits_{\mathbb{R}^n} f_{n,\varepsilon}(x)c_n d^n x\right)_\varepsilon = \qquad (1.2.4)$$

$$\left(\lim_{n\to\infty}\int\ldots\int f_{n,\varepsilon}(x_1,x_2,\ldots,x_n)c_n dx_1 dx_2\ldots dx_n\right)_\varepsilon$$

iff $\forall\varepsilon\in(0,1]$ limit in RHS of Eq.(1.2.4) exists.

**Remark 1.2.1.** We note that one can rewrite Definition 1.2.2 using completely Colombeau form.

**Definition 1.2.3.** For $f\in G(\mathbb{R}^\infty)$ we define the infinite dimensional Colombeau-Feinman integral $\Im_\infty((f_\delta)_\delta),\delta\in(0,1]$ as Colombeau quantity such that

$$\left(\Im_\infty(f_\delta)\right)_\delta = \left(\int\limits_{\mathbb{R}^\infty} f_\delta(\mathbf{x})D[\mathbf{x}]\right)_{\delta\in(0,1]} =_{\widetilde{\mathbb{C}}} \left(\left(\int\limits_{\mathbb{R}^\infty} f_\delta(\mathbf{x})D[\mathbf{x},\varepsilon]\right)_{\delta\in(0,1]}\right)_{\varepsilon\in(0,1]} =_{\widetilde{\mathbb{C}}}$$

$$\left(\left(\int\limits_{\mathbb{R}^{n_\varepsilon}} f_{n_\varepsilon,\delta}(x)c_{n_\varepsilon}d^{n_\varepsilon}x\right)_{\delta\in(0,1]}\right)_{\varepsilon\in(0,1]} =_{\widetilde{\mathbb{C}}}$$

$$(1.2.5)$$

$$((c_{n_\varepsilon})_\varepsilon)\left(\left(\int\limits_{\mathbb{R}^{n_\varepsilon}} f_{n_\varepsilon,\delta}(x)d^{n_\varepsilon}x\right)_{\delta\in(0,1]}\right)_{\varepsilon\in(0,1]} =_{\widetilde{\mathbb{C}}}$$

$$((c_{n_\varepsilon})_\varepsilon)\left(\left(\int\ldots\int f_{n_\varepsilon,\delta}(x_1,x_2,\ldots,x_{n_\delta})dx_1 dx_2\ldots dx_{n_\varepsilon}\right)_{\delta\in(0,1]}\right)_{\varepsilon\in(0,1]}.$$

**Definition 1.2.4.** (1) For $\widetilde{\mathbb{R}}$-valued functional $f \in G(\mathbb{R}^\infty)$ we define the infinite dimensional Colombeau-Feinman upper integral $\overline{\mathfrak{J}}_\infty((f_\varepsilon)_\varepsilon)$ by formula

$$\overline{\mathfrak{J}}_\infty((f_\varepsilon)_\varepsilon) = \left( \overline{\int_{\mathbb{R}^\infty}} f_\varepsilon(\mathbf{x}) D[\mathbf{x}] \right)_\varepsilon = \left( \limsup_{n \to \infty} \int_{\mathbb{R}^n} f_{n,\varepsilon}(x) c_n d^n x \right)_\varepsilon =$$

$$\left( \limsup_{n \to \infty} \int \ldots \int f_{n,\varepsilon}(x_1, x_2, \ldots, x_n) c_n dx_1 dx_2 \ldots dx_n \right)_\varepsilon .$$

(1.2.6)

(2) For $\widetilde{\mathbb{R}}$-valued Colombeau generalized functional $f \in (\mathbb{R}^\infty)^*$ we define the infinite dimensional lower integral $\underline{\mathfrak{J}}_\infty((f_\varepsilon)_\varepsilon)$ by formula

$$\underline{\mathfrak{J}}_\infty((f_\varepsilon)_\varepsilon) = \left( \underline{\int_{\mathbb{R}^\infty}} f_\varepsilon(\mathbf{x}) D[\mathbf{x}] \right)_\varepsilon = \left( \liminf_{n \to \infty} \int_{\mathbb{R}^n} f_{n,\varepsilon}(x) c_n d^n x \right)_\varepsilon =$$

$$\left( \liminf_{n \to \infty} \int \ldots \int f_{n,\varepsilon}(x_1, x_2, \ldots, x_n) c_n dx_1 dx_2 \ldots dx_n \right)_\varepsilon .$$

(1.2.7)

**Proposition 1.2.1.** Assume that $f \in G(\mathbb{R}^\infty)$ is $\widetilde{\mathbb{R}}$-valued Colombeau generalized functional, then $\underline{\mathfrak{J}}_\infty((f_\varepsilon)_\varepsilon) \leq_{\widetilde{\mathbb{R}}} \overline{\mathfrak{J}}_\infty((f_\varepsilon)_\varepsilon)$ and

$$\overline{\mathfrak{J}}_\infty((f_\varepsilon)_\varepsilon) = \left( \overline{\int_{\mathbb{R}^\infty}} f_\varepsilon(\mathbf{x}) D[\mathbf{x}] \right)_\varepsilon =_{\widetilde{\mathbb{R}}} \left( \lim_{k \to \infty} \left( \sup_{m \geq k} \int_{\mathbb{R}^m} f_{m,\varepsilon}(x) c_m d^m x \right) \right)_\varepsilon ,$$

$$\underline{\mathfrak{J}}_\infty((f_\varepsilon)_\varepsilon) = \left( \underline{\int_{\mathbb{R}^\infty}} f_\varepsilon(\mathbf{x}) D[\mathbf{x}] \right)_\varepsilon =_{\widetilde{\mathbb{R}}} \left( \lim_{k \to \infty} \left( \inf_{m \geq k} \int_{\mathbb{R}^m} f_{m,\varepsilon}(x) c_m d^m x \right) \right)_\varepsilon .$$

(1.2.8)

**Proposition 1.2.2.** Suppose that $\overline{\mathfrak{J}}_\infty((f_\varepsilon)_\varepsilon) =_{\widetilde{\mathbb{R}}} \underline{\mathfrak{J}}_\infty((f_\varepsilon)_\varepsilon) =_{\widetilde{\mathbb{R}}} a_{((f_\varepsilon)_\varepsilon)} < \infty, a_f \in \widetilde{\mathbb{R}}$.

Then $\mathfrak{J}_\infty((f_\varepsilon)_\varepsilon) =_{\widetilde{\mathbb{R}}} a_{((f_\varepsilon)_\varepsilon)}$.

**Proposition 1.2.3.** Assume that $f_1, f_2 \in G(\mathbb{R}^\infty)$ is $\widetilde{\mathbb{R}}$-valued Colombeau generalized

functionals and $f_1(\mathbf{x}) \leq_{\widetilde{\mathbb{R}}} f_2(\mathbf{x})$, then

$$\left(\overline{\int_{\mathbb{R}^{\infty}}}[f_{1,\varepsilon}(\mathbf{x})]D[\mathbf{x}]\right)_{\varepsilon} \leq_{\widetilde{\mathbb{R}}} \left(\overline{\int_{\mathbb{R}^{\infty}}}[f_2(\mathbf{x})]D[\mathbf{x}]\right)_{\varepsilon},$$

$$\left(\underline{\int_{\mathbb{R}^{\infty}}}[f_{1,\varepsilon}(\mathbf{x})]D[\mathbf{x}]\right)_{\varepsilon} \leq_{\widetilde{\mathbb{R}}} \left(\underline{\int_{\mathbb{R}^{\infty}}}[f_{2,\varepsilon}(\mathbf{x})]D[\mathbf{x}]\right)_{\varepsilon}.$$

$(1.2.9)$

**Proposition 1.2.4.** Assume that $f_1, f_2 \in G(\mathbb{R}^{\infty})$ is an $\widetilde{\mathbb{R}}$-valued Colombeau generalized functionals then

$$\left(\underline{\int_{\mathbb{R}^{\infty}}}[f_{1,\varepsilon}(\mathbf{x})]D[\mathbf{x}]\right)_{\varepsilon} + \left(\overline{\int_{\mathbb{R}^{\infty}}}[f_{2,\varepsilon}(\mathbf{x})]D[\mathbf{x}]\right)_{\varepsilon} \leq_{\widetilde{\mathbb{R}}} \left(\overline{\int_{\mathbb{R}^{\infty}}}[f_{1,\varepsilon}(\mathbf{x}) + f_{2,\varepsilon}(\mathbf{x})]D[\mathbf{x}]\right)_{\varepsilon}$$

$$\leq_{\widetilde{\mathbb{R}}} \left(\overline{\int_{\mathbb{R}^{\infty}}}[f_{1,\varepsilon}(\mathbf{x})]D[\mathbf{x}]\right)_{\varepsilon} + \left(\overline{\int_{\mathbb{R}^{\infty}}}[f_{2,\varepsilon}(\mathbf{x})]D[\mathbf{x}]\right)_{\varepsilon},$$

$$\left(\underline{\int_{\mathbb{R}^{\infty}}}[f_{1,\varepsilon}(\mathbf{x})]D[\mathbf{x}]\right)_{\varepsilon} + \left(\overline{\int_{\mathbb{R}^{\infty}}}[f_{2,\varepsilon}(\mathbf{x})]D[\mathbf{x}]\right)_{\varepsilon} \widetilde{\mathbb{R}} \geq$$

$$\left(\underline{\int_{\mathbb{R}^{\infty}}}[f_{1,\varepsilon}(\mathbf{x}) + f_{2,\varepsilon}(\mathbf{x})]D[\mathbf{x}]\right)_{\varepsilon} \widetilde{\mathbb{R}} \geq$$

$$\left(\underline{\int_{\mathbb{R}^{\infty}}}[f_{1,\varepsilon}(\mathbf{x})]D[\mathbf{x}]\right)_{\varepsilon} + \left(\underline{\int_{\mathbb{R}^{\infty}}}[f_{2,\varepsilon}(\mathbf{x})]D[\mathbf{x}]\right)_{\varepsilon}.$$

$(1.2.10)$

# Apendix II.2.1. Colombeau Signed Measures

In this section we introduced Colombeau type generalization of the notion of the signed measures. We also introduced Colombeau type generalization of the notion of the submeasures.

**Definition 2.1.1.** Suppose $\Sigma_X$ is an $\sigma$-algebra on a set $X$. A function $\mu : \to [-\infty, \infty]$ is called a signed measure on $\Sigma_X$, if it has the properties below:

$$(1)\ \mu(A) \leq \infty, \forall A \in \Sigma_X,$$

$$(2)\ \mu(A) \geq -\infty, \forall A \in \Sigma_X, \tag{2.1.1}$$

$$(3)\ \mu(\varnothing) = 0,$$

(4) For any pairwise disjoint sequence $\{A_n\}_{n=1}^{\infty}$, $A_n \subset \Sigma_X$ one has the equality

$$(4)\ \mu(\cup_{n=1}^{\infty} A_n) = \sum_{n=1}^{\infty} A_n. \tag{2.1.2}$$

Let us to use the term "honest" measure, for a measure in the usual sense.

**Definition 2.1.2.** Let $X$ be a non-empty set, let $\Sigma_X$ be a $\sigma$-algebra on $X$, and let $\mu$ be a signed measure on $\Sigma_X$ such that $\forall A \in \Sigma_X$ :

$$-\infty < \mu(A) < \infty \tag{2.1.3}$$

A signed measure $\mu(A)$ with property (2.2.1.3) is called finite.

**Theorem 2.1.1.** Let $X$ be a non-empty set, let $\Sigma_X$ be a $\sigma$-algebra on $X$, and let $\mu$ be a signed measure on $\Sigma_X$. Then there exist subsets $X^+, X^- \in \Sigma_X$, with the following properties:

**(i)** $X^+ \cap X^- = \varnothing$, $X^+ \cup X^- = X$;

**(ii)** the maps $\mu^{\pm} : \to [-\infty, \infty]$, defined by

$$\mu^{\pm}(A) = \pm\mu(A \cap X^{\pm}), \forall A \in \Sigma_X, \tag{2.1.4}$$

are "honest" measures on $\Sigma_X$;

(iii) one of the measures $\mu^{\pm}(A)$ is finite, and one has the equality

$$\mu = \mu^+ - \mu^-. \tag{2.1.5}$$

**Definition 2.1.3.** Let $X$ be a non-empty set, let $\Sigma_X$ be a $\sigma$-algebra on $X$. A function $\mu : \Sigma_X \to \mathbb{C}$ is called a complex measure on $\Sigma_X$, if it is $\sigma$-additive in the sense that for any pairwise disjoint sequence $\{A_n\}_{n=1}^{\infty}, \forall n (A_n \in \Sigma_X)$, one has the equality

$$(4)\ \mu(\cup_{n=1}^{\infty} A_n) = \sum_{n=1}^{\infty} A_n. \tag{2.1.6}$$

**Theorem 2.1.2.** Let $X$ be a non-empty set, let $\Sigma_X$ be a $\sigma$-algebra on $X$ and let $\mu$ be either a signed measure, or a complex measure on $\Sigma_X$. For every $A \in \Sigma_X$, we define

$$|\mu|(A) = \sup_{\cup_{k=1}^{\infty} A_k = X} \left\{ \sum_{k=1}^{\infty} |\mu(A_k)| : \{A_n\}_{n=1}^{\infty} \text{pairwise disjoint} \right\}. \tag{2.1.7}$$

Then the map $|\mu|(A) : \Sigma_X \to [0, \infty]$ is an "honest" measure on $\Sigma_X$.

**Definition 2.1.4.** The "honest" measure $|\mu|(A)$, defined by (2.2.1.7), is called the variation measure of $\mu(A)$. By construction, we have the inequality

$$|\mu(A)| \leq |\mu|(A). \tag{2.1.8}$$

In the case of signed measures, the variation measure is also given by the following.

**Theorem 2.1.3.** Let $X$ be a non-empty set, let $\Sigma_X$ be a $\sigma$-algebra on $X$ and let $\mu$ be a signed measure on $\Sigma_X$. Then one has the equality

$$|\mu| = \mu^+ + \mu^-, \tag{2.1.9}$$

where $\mu = \mu^+ - \mu^-$ is the Hahn-Jordan decomposition of $\mu$.

**Definition 2.2.1.5.** Let $X$ be a non-empty set, let $\Sigma_X$ be a $\sigma$-algebra on $X$. Let

$\mathbf{M}_{\mathbb{R}}(\Sigma_X)$ be a set of the all finite signed measures on $\Sigma_X$. Set
$E_{\mathbb{R}}(\Sigma_X) = (\mathbf{M}_{\mathbb{R}}(\Sigma_X))^I, I = (0,1],$

$$E_M^{\mathbb{R}}(\Sigma_X) =$$

$$\left\{ (\mu_\varepsilon)_\varepsilon \in E_{\mathbb{R}}(\Sigma_X) | (\forall n \in \mathbb{N})((\forall m \in \mathbb{N}))(\exists p \in \mathbb{N}) \left[ \sup_{A \in \Sigma_X} |\mu_\varepsilon(A)| = O(\varepsilon^{-p}) \right] \right\}, \qquad (2.1.10)$$

$$N_M^{\mathbb{R}}(\Sigma_X) =$$

$$\left\{ (\mu_\varepsilon)_\varepsilon \in E_{\mathbb{R}}(\Sigma_X) | (\forall n \in \mathbb{N})((\forall m \in \mathbb{N}))(\forall q \in \mathbb{N}) \left[ \sup_{A \in \Sigma_X} |\mu_\varepsilon(A)| = O(\varepsilon^{q}) \right] \right\}, \qquad (2.1.11)$$

$$G_{\mathbb{R}}(\Sigma_X) = E_M^{\mathbb{R}}(\Sigma_X)/N_M^{\mathbb{R}}(\Sigma_X). \qquad (2.1.12)$$

Elements of $E_M^{\mathbb{R}}(\Sigma_X)$ and $N_M^{\mathbb{R}}(\Sigma_X)$ are called moderate, resp. negligible signed measures. Note that $G_{\mathbb{R}}(\Sigma_X)$ is algebra with pointwise operations. It is the largest subalgebra of $E_M^{\mathbb{R}}(\Sigma_X)$ in which $N_M^{\mathbb{R}}(\Sigma_X)$ is ideal. Thus, $G_{\mathbb{R}}(\Sigma_X)$ is an associative, commutative algebra. If $(\mu_\varepsilon)_\varepsilon \in E_M^{\mathbb{R}}(\Sigma_X)$ is a representative of $\mu \in G_{\mathbb{R}}(\Sigma_X)$ we write
$\mu = \mathbf{cl}[(\mu_\varepsilon)_\varepsilon]$ or simply
$\mu = [(\mu_\varepsilon)_\varepsilon]$.

**Definition 2.1.6.** An element $\mu \in G_{\mathbb{R}}(\Sigma_X)$ and any representative $(\mu_\varepsilon)_\varepsilon \in E_M^{\mathbb{R}}(\Sigma_X)$ of $\mu$ is called Colombeau signed measure on $\Sigma_X$.

**Definition 2.1.7.** Let $X$ be a non-empty set, let $\Sigma_X$ be a $\sigma$-algebra on $X$. Let $\mathbf{M}_{\mathbb{C}}(\Sigma_X)$ be a set of the all complex finite measures on $\Sigma_X$. Set
$E_{\mathbb{C}}(\Sigma_X) = (\mathbf{M}_{\mathbb{C}}(\Sigma_X))^I, I = (0,1],$

$$E_M^{\mathbb{C}}(\Sigma_X) =$$

$$\left\{ (\mu_\varepsilon)_\varepsilon \in E_{\mathbb{C}}(\Sigma_X) | (\forall n \in \mathbb{N})((\forall m \in \mathbb{N}))(\exists p \in \mathbb{N}) \left[ \sup_{A \in \Sigma_X} |\mu_\varepsilon(A)| = O(\varepsilon^{-p}) \right] \right\}, \qquad (2.1.13)$$

$$N_M^{\mathbb{C}}(\Sigma_X) =$$

$$\left\{ (\mu_\varepsilon)_\varepsilon \in E_{\mathbb{C}}(\Sigma_X) | (\forall n \in \mathbb{N})((\forall m \in \mathbb{N}))(\forall q \in \mathbb{N}) \left[ \sup_{A \in \Sigma_X} |\mu_\varepsilon(A)| = O(\varepsilon^q) \right] \right\}, \qquad (2.1.14)$$

$$G_{\mathbb{C}}(\Sigma_X) = E_M^{\mathbb{C}}(\Sigma_X)/N_M^{\mathbb{C}}(\Sigma_X). \qquad (2.1.15)$$

Elements of $E_M^{\mathbb{C}}(\Sigma_X)$ and $N_M^{\mathbb{C}}(\Sigma_X)$ are called moderate, resp. negligible complex measures.

Note that $G_{\mathbb{R}}(\Sigma_X)$ is algebra with pointwise operations. It is the largest subalgebra of $E_M^{\mathbb{R}}(\Sigma_X)$ in which $N_M^{\mathbb{R}}(\Sigma_X)$ is ideal. Thus, $G_{\mathbb{R}}(\Sigma_X)$ is an associative, commutative algebra. If $(\mu_\varepsilon)_\varepsilon \in E_M^{\mathbb{R}}(\Sigma_X)$ is a representative of $\mathbf{\mu} \in G_{\mathbb{R}}(\Sigma_X)$ we write $\mathbf{\mu} = \mathbf{cl}[(\mu_\varepsilon)_\varepsilon]$ or simply $\mathbf{\mu} = [(\mu_\varepsilon)_\varepsilon]$.

**Definition 2.1.8**. An element $\mathbf{\mu} \in G_{\mathbb{R}}(\Sigma_X)$ and any representative $(\mu_\varepsilon)_\varepsilon \in E_M^{\mathbb{R}}(\Sigma_X)$ of $\mathbf{\mu}$ is called Colombeau complex measure on $\Sigma_X$.

**Theorem 2.1.4**.Let $X$ be a non-empty set, let $\Sigma_X$ be a $\sigma$-algebra on $X$, and let $(\mu_\varepsilon)_\varepsilon$ be Colombeau signed measure on $\Sigma_X$. Then $\forall \varepsilon \in (0, 1]$ there exist subsets $X_\varepsilon^+, X_\varepsilon^- \in \Sigma_X$, with the following properties:

**(i)** $\forall \varepsilon (\varepsilon \in (0,1]) \left[ X_\varepsilon^+ \cap X_\varepsilon^- = \varnothing \right]$, $\forall \varepsilon (\varepsilon \in (0,1]) \left[ X_\varepsilon^+ \cup X_\varepsilon^- = X \right]$;

**(ii)** the maps $(\mu_\varepsilon^\pm)_\varepsilon \to [-\infty, +\infty]^I$ and $\mathbf{\mu}^\pm : \to \widetilde{\mathbb{R}}$, defined by

$$\mu_\varepsilon^\pm(A) = \pm\mu_\varepsilon(A \cap X_\varepsilon^\pm), \forall A \in \Sigma_X, \forall \varepsilon(\varepsilon \in (0,1])$$

$$\qquad (2.1.16)$$

$$\mathbf{\mu}^\pm(A) = \mathbf{cl}[(\mu_\varepsilon^\pm)_\varepsilon(A)]$$

are "honest" Colombeau measures on $\Sigma_X$;

**(iii)** $\forall \varepsilon \in (0, 1]$ measures $\mu_\varepsilon^\pm(A)$ is finite, and one has the equality

$$(\mu_\varepsilon)_\varepsilon = (\mu_\varepsilon^+)_\varepsilon - (\mu_\varepsilon^-)_\varepsilon. \qquad (2.1.17)$$

**Theorem 2.1.5.**Let $X$ be a non-empty set, let $\Sigma_X$ be a $\sigma$-algebra on $X$ and let $(\mu_\varepsilon)_\varepsilon$ be either a signed Colombeau measure, or a complex Colombeau measure on $\Sigma_X$. For every $A \in \Sigma_X$, we define

$$|(\mu_\varepsilon)_\varepsilon|(A) = \left( \sup_{\cup_{k=1}^\infty A_k = X} \left\{ \sum_{k=1}^\infty |\mu_\varepsilon(A_k)| : \{A_n\}_{n=1}^\infty \text{pairwise disjoint} \right\} \right)_\varepsilon. \quad (2.1.18)$$

Then the map $|(\mu_\varepsilon)_\varepsilon|(A) : \Sigma_X \to [0,\infty]^I$ is an "honest" Colombeau measure on $\Sigma_X$.

**Definition 2.1.9.**The "honest" Colombeau measure $|(\mu_\varepsilon)_\varepsilon|(A)$,defined by (2.2.1.18), is called the variation of Colombeau measure $(\mu_\varepsilon)_\varepsilon(A)$. By construction, we have the inequality

$$|(\mu_\varepsilon(A))_\varepsilon| \le (|\mu_\varepsilon|(A))_\varepsilon. \quad (2.1.19)$$

In the case of signed Colombeau measures, the variation measure is also given by the following.

**Theorem 2.1.6.**Let $X$ be a non-empty set, let $\Sigma_X$ be a $\sigma$-algebra on $X$ and let $(\mu_\varepsilon)_\varepsilon$ be a signed Colombeau measure on $\Sigma_X$. Then $\forall A \in \Sigma_X$ one has the equality

$$|(\mu_\varepsilon)_\varepsilon|(A) = ((\mu_\varepsilon^+)_\varepsilon)(A) + ((\mu_\varepsilon^-)_\varepsilon)(A), \quad (2.1.20)$$

where $\varepsilon \in (0,1]$ and $\mu_\varepsilon = \mu_\varepsilon^+ - \mu_\varepsilon^-$ is the Hahn-Jordan decomposition of $\mu_\varepsilon$. We remind now that definitions.

**Definition 2.1.10.** The signum function $\mathbf{sgn}(x)$ of a real number $x \in \mathbb{R}$ is defined as follows:

$$\mathbf{sgn}(x) = \begin{cases} -1 \text{ if } x < 0 \\ \phantom{-}0 \text{ if } x = 0 \\ \phantom{-}1 \text{ if } x > 0 \end{cases} \quad (2.1.21)$$

We note that any real number can be expressed as the product of its absolute value and its sign function:

$$x = |x|\mathbf{sgn}(x) \tag{2.1.22}$$

**Definition 2.1.11.** Let $X$ be a non-empty set. (1) The positive part $f^+(x)$ of a real-valued function $f : X \to \mathbb{R}$ is defined via formula:

$$f^+(x) = \max\{f(x), 0\} =$$

$$\begin{cases} f(x) \text{ if } f(x) > 0 \\ 0 \text{ otherwise} \end{cases} \tag{2.1.23}$$

(2) The negative part $f^-(x)$ of a real-valued function $f : X \to \mathbb{R}$ is defined via formula

$$f^-(x) = \max\{-f(x), 0\} =$$

$$-\min\{f(x), 0\} =$$

$$\begin{cases} -f(x) \text{ if } f(x) < 0 \\ 0 \text{ otherwise} \end{cases} \tag{2.1.24}$$

**Remark.2.1.1.**Note that both $f^+(x)$ and $f^-(x)$ are non-negative functions.The function $f(x)$ can be expressed in terms of $f^+(x)$ and $f^-(x)$ as

$$f(x) = f^+(x) - f^-(x) \tag{2.1.25}$$

Also note that

$$f^+(x) = \big| f(x) \big| \mathbf{sgn}^+\big(f(x)\big),$$

$$f^-(x) = \big| f(x) \big| \mathbf{sgn}^-\big(f(x)\big).$$

(2.1.26)

Therefore from Eqs.(2.1.25)-(2.1.26) one obtain

$$f(x) = \Big[ \mathbf{sgn}^+\big(f(x)\big) - \mathbf{sgn}^-\big(f(x)\big) \Big] \big| f(x) \big|.$$

(2.1.27)

**Remark**.**2**.**1**.**2**. Consider a nonnegative finite measure $v : \Sigma_X \to \mathbb{R}_+$ on the space $(X, \Sigma_X)$ and a measurable function $f : X \to \mathbb{R}$ such that $\int_A \big| f(x) \big| dv(x) < \infty.$ Then, a finite signed measure $\mu_{\mathbf{sg}}(A)$ is given by formula:

$$\mu_{\mathbf{sg}}(A) = \int_A f(x) dv(x) =$$

$$\int_A \Big[ f^+(x) - f^-(x) \Big] dv(x) =$$

(2.1.28)

$$\int_A f^+(x) dv(x) - \int_A f^-(x) dv(x) =$$

$$\mu_{\mathbf{sg}}^+(A) - \mu_{\mathbf{sg}}^-(A).$$

for all $A \in \Sigma_X.$

**Remark**.**2**.**1**.**3**. Substitution Eq.(2.1.27) into Eq.(2.1.28) gives

$$\mu_{\mathbf{sg}}(A) = \int_A f(x) d\nu(x) =$$

$$\int_A \Big[ \mathbf{sgn}^+\big(f(x)\big) - \mathbf{sgn}^-\big(f(x)\big) \Big] \big| f(x) \big| d\nu(x) = \qquad (2.1.29)$$

$$\int_A \big| f(x) \big| d\widetilde{\nu}_{\mathbf{sg}}(x),$$

where $d\widetilde{\nu}_{\mathbf{sg}}(x) = \Big[ \mathbf{sgn}^+\big(f(x)\big) - \mathbf{sgn}^-\big(f(x)\big) \Big] d\nu(x)$.

**Example 2.1.1.** Freshnel signed measure $\mu_{Fr}$ on $\Sigma_X$. $X = \mathbb{R}^n$, $f(x) = \sin(x^2)$, $d\nu(x) = d^n x$.

$$\mu_{Fr}(A) = \int_A \sin(x^2) d^n x =$$

$$\int_A \Big[ \mathbf{sgn}^+\big(\sin(x^2)\big) - \mathbf{sgn}^-\big(\sin(x^2)\big) \Big] \big| \sin(x^2) \big| d^n x = \qquad (2.1.30)$$

$$\int_A \big| \sin(x^2) \big| d\widetilde{\nu}_{\mathbf{sg}}(x),$$

$$d\widetilde{\nu}_{\mathbf{sg}}(x) = \Big[ \mathbf{sgn}^+\big(\sin(x^2)\big) - \mathbf{sgn}^-\big(\sin(x^2)\big) \Big] d^n x.$$

**Remark 2.1.4.** Let us rewrite (2.1.29) in equivalent form

$$\mu_{\mathbf{sg}}(A) = \int\limits_A f(x) dv(x) =$$

$$\int\limits_A \left[ \mathbf{sgn}^+\big(f(x)\big) - \mathbf{sgn}^-\big(f(x)\big) \right] \big| f(x) \big| dv(x) =$$

$$\int\limits_A \left\{ \left[ \mathbf{sgn}^+\big(f(x)\big) - \mathbf{sgn}^-\big(f(x)\big) \right] + 1 - 1 \right\} \big| f(x) \big| dv(x) = \qquad (2.1.31)$$

$$\int\limits_A \left\{ \left[ \mathbf{sgn}^+\big(f(x)\big) - \mathbf{sgn}^-\big(f(x)\big) \right] + 1 \right\} \big| f(x) \big| dv(x) - \int\limits_A \big| f(x) \big| dv(x) =$$

$$\int\limits_A \big| f(x) \big| d\widetilde{v}_+(x) - \int\limits_A \big| f(x) \big| dv(x),$$

where positive measure $d\widetilde{v}_+(x)$ is given via formula

$$d\widetilde{v}_+(x) = \left\{ \left[ \mathbf{sgn}^+\big(f(x)\big) - \mathbf{sgn}^-\big(f(x)\big) \right] + 1 \right\} dv(x). \qquad (2.1.32)$$

**Example 2.1.2.** Freshnel signed measure $\mu_{Fr}$ on $\Sigma_X$. $X = \mathbb{R}^n$, $f(x) = \sin(x^2)$, $dv_+(x) = d^n x$.

$$\mu_{Fr}(A) = \int\limits_{A} \sin(x^2) d^{\,n}x =$$

$$\int\limits_{A} \left\{ \left[ \mathbf{sgn}^+\left(\sin(x^2)\right) - \mathbf{sgn}^-\left(\sin(x^2)\right)\right] + 1 - 1 \right\} \left| \sin(x^2)\right| d^n x =$$

$$\int\limits_{A} \left\{ \left[ \mathbf{sgn}^+\left(\sin(x^2)\right) - \mathbf{sgn}^-\left(\sin(x^2)\right)\right] + 1 \right\} \left| \sin(x^2)\right| d^n x - \int\limits_{A} \left| \sin(x^2)\right| d^n x = \qquad (2.1.33)$$

$$\int\limits_{A} \left| \sin(x^2)\right| d\widetilde{v}_+(x) - \int\limits_{A} \left| \sin(x^2)\right| d^n x,$$

$$d\widetilde{v}_+(x) = \left\{ \left[ \mathbf{sgn}^+\left(\sin(x^2)\right) - \mathbf{sgn}^-\left(\sin(x^2)\right)\right] + 1 \right\} d^{\,n}x.$$

**Definition 2.1.12.** Let $X$ be a non-empty set.

(1) The positive part $(f_\varepsilon(x))_\varepsilon^+$ of a real-valued generalized function $(f_\varepsilon)_\varepsilon : X \to \widetilde{\mathbb{R}}$ is defined via formula:

$$(f_\varepsilon(x))_\varepsilon^+ = \left(f_\varepsilon^+(x)\right)_\varepsilon; \qquad (2.1.34)$$

(2) The negative part $(f_\varepsilon(x))_\varepsilon^-$ of a real-valued generalized function $(f_\varepsilon)_\varepsilon : X \to \widetilde{\mathbb{R}}$ is defined via formula

$$(f_\varepsilon(x))_\varepsilon^- = \left(f_\varepsilon^-(x)\right)_\varepsilon. \qquad (2.1.35)$$

**Remark.2.1.5.** We note that both $\left(f_\varepsilon^+(x)\right)_\varepsilon$ and $\left(f_\varepsilon^-(x)\right)_\varepsilon$ are non-negative Colombeau generalized functions. The any Colombeau generalized function $(f_\varepsilon(x))_\varepsilon$ can be expressed in terms of $(f_\varepsilon(x))_\varepsilon^+$ and $(f_\varepsilon(x))_\varepsilon^-$ as

$$(f_\varepsilon(x))_\varepsilon = (f_\varepsilon(x))_\varepsilon^+ - (f_\varepsilon(x))_\varepsilon^-. \qquad (2.1.36)$$

Also we note that

$$(f_\varepsilon(x))_\varepsilon^+ = \big| \, (f_\varepsilon(x))_\varepsilon \, \big| \mathbf{sgn}^+ \big[ \big( \big( f_\varepsilon(x) \big) \big)_\varepsilon \big],$$

$$(2.1.37)$$

$$(f_\varepsilon(x))_\varepsilon^- = \big| (f_\varepsilon(x))_\varepsilon \, \big| \mathbf{sgn}^- \big[ \big( \big( f_\varepsilon(x) \big) \big)_\varepsilon \big].$$

Therefore from Eqs.(2.1.36)-(2.1.37) one obtain

$$(f_\varepsilon(x))_\varepsilon = \Big[ \mathbf{sgn}^+ \big[ \big( \big( f_\varepsilon(x) \big) \big)_\varepsilon \big] - \mathbf{sgn}^- \big[ \big( \big( f_\varepsilon(x) \big) \big)_\varepsilon \big] \Big] \big| \, (f_\varepsilon(x))_\varepsilon \, \big|. \qquad (2.1.38)$$

**Remark**.**2**.**1**.**6**.Consider a nonnegative Colombeau measure $(v_\varepsilon)_\varepsilon : \Sigma_X \to \widetilde{\mathbb{R}}_+$ on the space $(X, \Sigma_X)$ and a measurable Colombeau generalized function $(f_\varepsilon)_\varepsilon : X \to \widetilde{\mathbb{R}}$ such that $\left( \int\limits_A \big| f_\varepsilon(x) \big| dv_\varepsilon(x) \right)_\varepsilon < \infty.$ Then, a finite signed Colombeau measure $(\mu_{\mathbf{sg}}(A, \varepsilon))_\varepsilon$ is given by formula:

$$\left(\mu_{\mathbf{sg}}(A,\varepsilon)\right)_{\varepsilon} = \left(\int_A f_{\varepsilon}(x)dv_{\varepsilon}(x)\right)_{\varepsilon} =$$

$$\left(\int_A \left[f_{\varepsilon}^{+}(x) - f_{\varepsilon}^{-}(x)\right]dv_{\varepsilon}(x)\right)_{\varepsilon} =$$

$$\left(\int_A f_{\varepsilon}^{+}(x)dv_{\varepsilon}(x)\right)_{\varepsilon} - \left(\int_A f_{\varepsilon}^{-}(x)dv_{\varepsilon}(x)\right)_{\varepsilon} = \qquad (2.1.39)$$

$$\left(\mu_{\mathbf{sg}}^{+}(A,\varepsilon)\right)_{\varepsilon} - \left(\mu_{\mathbf{sg}}^{-}(A,\varepsilon)\right)_{\varepsilon}.$$

for all $A \in \Sigma_X$.

**Remark.2.1.7.** Substitution Eq.(2.1.38) into Eq.(2.1.39) gives

$$\mu_{\mathbf{sg}}(A,\varepsilon) = \int_A f_{\varepsilon}(x)dv_{\varepsilon}(x) =$$

$$\left(\int_A \left[\mathbf{sgn}^{+}\left(f_{\varepsilon}(x)\right) - \mathbf{sgn}^{-}\left(f_{\varepsilon}(x)\right)\right]\left|f_{\varepsilon}(x)\right|dv_{\varepsilon}(x)\right)_{\varepsilon} = \qquad (2.1.40)$$

$$\left(\int_A \left|f_{\varepsilon}(x)\right|d\widetilde{v}_{\mathbf{sg}}(x,\varepsilon)\right)_{\varepsilon},$$

where

$$\left(d\widetilde{v}_{\mathbf{sg}}(x,\varepsilon)\right)_{\varepsilon} = \left[\mathbf{sgn}^{+}\left(f(x)\right) - \mathbf{sgn}^{-}\left(f(x)\right)\right]dv(x). \qquad (2.1.41)$$

**Remark**.2.1.8. Let us rewrite (2.1.40) in equivalent form

$$\left(\mu_{\mathbf{sg}}(A,\varepsilon)\right)_{\varepsilon} = \left(\int\limits_{A} f_{\varepsilon}(x)dv_{\varepsilon}(x)\right)_{\varepsilon} =$$

$$\left(\int\limits_{A}\left[\mathbf{sgn}^{+}\big(f_{\varepsilon}(x)\big) - \mathbf{sgn}^{-}\big(f_{\varepsilon}(x)\big)\right]\big|f_{\varepsilon}(x)\big|dv_{\varepsilon}(x)\right)_{\varepsilon} =$$

$$\left(\int\limits_{A}\left\{\left[\mathbf{sgn}^{+}\big(f_{\varepsilon}(x)\big) - \mathbf{sgn}^{-}\big(f_{\varepsilon}(x)\big)\right] + 1 - 1\right\}\big|f_{\varepsilon}(x)\big|dv_{\varepsilon}(x)\right)_{\varepsilon} =$$

(2.1.42)

$$\left(\int\limits_{A}\left\{\left[\mathbf{sgn}^{+}\big(f_{\varepsilon}(x)\big) - \mathbf{sgn}^{-}\big(f_{\varepsilon}(x)\big)\right] + 1\right\}\big|f_{\varepsilon}(x)\big|dv_{\varepsilon}(x)\right)_{\varepsilon} -$$

$$\left(\int\limits_{A}\big|f_{\varepsilon}(x)\big|dv_{\varepsilon}(x)\right)_{\varepsilon} =$$

$$\left(\int\limits_{A}\big|f_{\varepsilon}(x)\big|d\tilde{v}_{\varepsilon}^{+}(x)\right)_{\varepsilon} - \left(\int\limits_{A}\big|f_{\varepsilon}(x)\big|dv_{\varepsilon}(x)\right)_{\varepsilon},$$

where positive Colombeau measure $(d\tilde{v}_{\varepsilon}^{+}(x))_{\varepsilon}$ is given via formula

$$(d\tilde{v}_{\varepsilon}^{+}(x))_{\varepsilon} = \left(\left\{\left[\mathbf{sgn}^{+}\big(f_{\varepsilon}(x)\big) - \mathbf{sgn}^{-}\big(f_{\varepsilon}(x)\big)\right] + 1\right\}dv_{\varepsilon}(x)\right)_{\varepsilon}.$$

(2.1.43)

**Definition 2.1.13.** Let $(\lambda_{\varepsilon})_{\varepsilon}$ and $(v_{\varepsilon})_{\varepsilon}$ be signed Colombeau measures on $(X, \Sigma_{X})$ then $(\lambda_{\varepsilon})_{\varepsilon}$ is called absolutely continuous with respect to $(v_{\varepsilon})_{\varepsilon}$ $((\lambda_{\varepsilon})_{\varepsilon} \ll (v_{\varepsilon})_{\varepsilon})$ if $(|v_{\varepsilon}|(A))_{\varepsilon} = 0$ imply $(\lambda_{\varepsilon}(A))_{\varepsilon} = 0$.

**Theorem 2.1.8.**Let $(\mu_{\varepsilon})_{\varepsilon}$ be $\sigma$-additive Colombeau measure on the space

$(X, \Sigma_X), (\nu_{\mathbf{sg}}(\cdot, \varepsilon))_\varepsilon$ be a signed Colombeau measure on $(X, \Sigma_X)$ and $(\nu_{\mathbf{sg}}(\cdot, \varepsilon))_\varepsilon$ is absolutely continuous with respect to $(\mu_\varepsilon)_\varepsilon$. Then there exists unique $(f_\varepsilon)_\varepsilon \in G_{L_1(X,\mu)}$ such that

$$(\nu_{\mathbf{sg}}(A, \varepsilon))_\varepsilon = \left( \int\limits_A f_\varepsilon(x) d\mu_\varepsilon \right)_\varepsilon \tag{2.1.44}$$

for any $A \in \Sigma_X$.

# Apendix II.2.2.Colombeau Submeasures.

**Definition 2.2.1.** Let $X$ be a non-empty set, let $\Sigma_X$ be a $\sigma$-algebra on $X$. A submeasure $\mu$ on $\Sigma_X$ is a function $\mu : \Sigma_X \to [0, +\infty]$ which has the properties $\forall A, B \in \Sigma_X$ :

$$\mu(\varnothing) = 0,$$

$$A \subseteq B \Rightarrow \mu(A) \le \mu(B), \tag{2.2.1}$$

$$\mu(A \cup B) \le \mu(A) + \mu(B).$$

**Definition 2.2.2.** Let $X$ be a non-empty set, let $\Sigma_X$ be a $\sigma$-algebra on $X$. Let $\mathbf{SM}(\Sigma_X)$ be a set of the all finite submeasures on $\Sigma_X$. Set $SE(\Sigma_X) = (\mathbf{SM}(\Sigma_X))^I, I = (0, 1]$,

$$SE_M(\Sigma_X) =$$

$$\left\{ (\mu_\varepsilon)_\varepsilon \in SE(\Sigma_X) | (\forall n \in \mathbb{N})((\forall m \in \mathbb{N}))(\exists p \in \mathbb{N}) \left[ \sup_{A \in \Sigma_X} |\mu_\varepsilon(A)| = O(\varepsilon^{-p}) \right] \right\}, \tag{2.2.2}$$

$$SN_M(\Sigma_X) =$$

$$\left\{ (\mu_\varepsilon)_\varepsilon \in SE(\Sigma_X) | (\forall n \in \mathbb{N})((\forall m \in \mathbb{N}))(\forall q \in \mathbb{N}) \left[ \sup_{A \in \Sigma_X} |\mu_\varepsilon(A)| = O(\varepsilon^q) \right] \right\}, \qquad (2.2.3)$$

$$SG(\Sigma_X) = SE_M(\Sigma_X)/SN_M(\Sigma_X). \qquad (2.2.4)$$

Elements of $SE_M(\Sigma_X)$ and $SN_M(\Sigma_X)$ are called moderate, resp. negligible submeasures.

Note that $SG(\Sigma_X)$ is algebra with pointwise operations. It is the largest subalgebra of $SE_M(\Sigma_X)$ in which $SN_M(\Sigma_X)$ is ideal. Thus,$SG(\Sigma_X)$ is an associative, commutative algebra.If $(\mu_\varepsilon)_\varepsilon \in SE_M(\Sigma_X)$ is a representative of $\mu \in SG(\Sigma_X)$ we write $\mu = \mathbf{cl}[(\mu_\varepsilon)_\varepsilon]$ or simply $\mu = [(\mu_\varepsilon)_\varepsilon]$.

**Definition 2.2.3**. An element $\mu \in SG(\Sigma_X)$ and any representative $(\mu_\varepsilon)_\varepsilon \in SE_M(\Sigma_X)$ of $\mu$ is called Colombeau submeasure on $\Sigma_X$.

**Definition 2.2.4**.Let $\Re$ be a ring of subsets of a set $\Im \neq \varnothing$.

A set function $\mu : \Re \to [0,\infty)$ is said to be a Dobrakov submeasure or D-submeasure, if it is

(**1**) monotone: if $E, F \in \Re$ such that $E \subset F$, then $\mu(E) \leq \mu(F)$;

(**2**) subadditively continuous: for every $F \in \Re$ and $\varepsilon > 0$ there exists a $\delta > 0$ such that for every $E \in \Re$ with $\mu(E) < \delta$

$$(i)\mu(E \cup F) \leq \mu(F) + \varepsilon,$$

$$(2.2.5)$$

$$(ii)\mu(F) \leq \mu(F \backslash E) + \varepsilon;$$

# Apendix II.3.1. Colombeau Signed Measures and Submeasures on Infinite-Dimensional Linear

# Spaces.

In this section we introduced a generalization of the notion of the positive Feynman-Colombeau "measure" to the case where the values are allowed to be outside set $[0, \infty]$. We also introduced an infinite-dimensional generalization of the notion of the submeasures.

**Definition 2.2.2.1.** Suppose $\mathfrak{R}^\infty$ is an $\sigma$-algebra on a set $\mathbb{R}^\infty = \bigoplus_{n=1}^{\infty} \mathbb{R}^n$. A function $\mu : \to [-\infty, \infty]$ is called a signed measure on $\mathfrak{R}^\infty$, if it has the properties below:

$$(1) \ \mu(A) \leq \infty, \forall A \in \mathfrak{R}^\infty,$$

$$(2) \ \mu(A) \geq -\infty, \forall A \in \mathfrak{R}^\infty, \qquad (2.2.1)$$

$$(3) \ \mu(\varnothing) = 0,$$

(4) For any pairwise disjoint sequence $\{A_n\}_{n=1}^{\infty}$, $A_n \subset \mathfrak{R}^\infty$ one has the equality

$$(4) \ \mu(\cup_{n=1}^{\infty} A_n) = \sum_{n=1}^{\infty} A_n. \qquad (2.2.2)$$

**Examples 2.2.1.** Let us agree to use the term "honest" measure, for a measure in the usual sense.
1. Any "honest" measure on $\mathfrak{R}^\infty$ is a signed measure on $\mathfrak{R}^\infty$.
2. If $\mu$ is a signed measure on $\mathfrak{R}^\infty$, then $-\mu$ is again a signed measure on $\mathfrak{R}^\infty$.
3. If $\mu_1$ and $\mu_2$ are "honest" measures on $\mathfrak{R}^\infty$, one of which is finite, then $\mu_1 - \mu_2$ is a signed measure on $\mathfrak{R}^\infty$.

**Definition 2.2.2.** Let $\mathbf{M}(\mathfrak{R}^\infty)$ be a set of the all signed measures on $\mathfrak{R}^\infty$.
Set $E(\mathfrak{R}^\infty) = (\mathbf{M}(\mathfrak{R}^\infty))^I, I = (0, 1],$

**Definition 2.** (1) A submeasure $\mu$ on $\mathfrak{R}^\infty$ is a setfunction which has the properties

$$\mu(\varnothing) = 0,$$

$$A \subseteq B \Rightarrow \mu(A) \leq \mu(B), \qquad (2.2.)$$

$$\mu(A \cup B) \leq \mu(A) + \mu(B).$$

# Apendix III. Colombeau-Albeverio Path Integral.

Let **H** be a real separable Hilbert space, considered as a measurable Hilbert space with the o-algebra generated by its open subsets. Let $\|\cdot\|$ ( be the norm in **H**, We want to define a normalized subintegral on **H** of the following font

$$(\mathbf{J}[\Phi_\varepsilon])_\varepsilon = \left( \widetilde{\int_{\mathbf{H}}} \exp\left[ \frac{i}{2} \|\gamma\|^2 \right] \Phi_\varepsilon(\gamma) D[\gamma] \right)_{\varepsilon \in (0,1]} \qquad (3.1)$$

where $i$ is the imaginary unit, $(\Phi_\varepsilon)_\varepsilon$ belongs to a suitable Colombeau class of "integrable functions" on **H** and $\sim$ above the integral reminds us to the normalization. Normalization should mean that

$$(\mathbf{J}[1_\varepsilon])_{\varepsilon \in (0,1]} = \left( \widetilde{\int_{\mathbf{H}}} \exp\left[ \frac{i}{2} \|\gamma\|^2 \right] D[\gamma] \right)_{\varepsilon \in (0,1]} = \mathbf{cl}\left[ (1_\varepsilon)_{\varepsilon \in (0,1]} \right]. \qquad (3.2)$$

To see what tat this implies let us Took first to the case where **H** - is finite i.e. $\mathbf{H} = \mathbb{R}^d$.
Thus one obtain

$$(2\pi i)^{-\frac{d}{2}} \int\limits_{\mathbb{R}^d} \exp\left[\frac{i}{2}\|\gamma\|^2\right] D[\gamma] = 1,$$

(3.3)

$$i^{-\frac{d}{2}} = \exp\left(\frac{i\pi d}{4}\right).$$

Hence from Eq.(3.3),Eq.(3.2) we see that

$$(\mathbf{J}[1_\varepsilon])_{\varepsilon\in(0,1]} = \left(\widetilde{\int\limits_{\mathbf{H}}} \exp\left[\frac{i}{2}\|\gamma\|^2\right] D[\gamma]\right)_{\varepsilon\in(0,1]} =$$

(3.4)

$$(2\pi i)^{-\frac{d}{2}} \int\limits_{\mathbb{R}^d} \exp\left[\frac{i}{2}\|\gamma\|^2\right] d^d\gamma.$$

We also want that the normalized subintegral given via Eq.(3.1) be translation Invariant, which is formally expressed by the writing $D[\gamma]$ in Eq.(3.4) as a Lebesgue measure $d^d\gamma$; in fact when $\mathbf{H} = \mathbb{R}^d$ the translation Invariance of the Lebesgue measure $d^d\gamma$ implies the translation invariance of the usual integral $(2\pi i)^{-\frac{d}{2}} \int\limits_{\mathbb{R}^d} \exp\left[\frac{i}{2}\|\gamma\|^2\right] d^d\gamma.$

The condition for translation fnvariance is

$$\left(\widetilde{\int\limits_{\mathbf{H}}} \exp\left[\frac{i}{2}\|\gamma + \alpha_\varepsilon\|^2\right] D[\gamma]\right)_\varepsilon = (1)_\varepsilon.$$

(3.5)

for all $[(\alpha_\varepsilon)_\varepsilon] \in G_{\mathbf{H}}$ i.e. $(\alpha_\varepsilon)_\varepsilon \in E_M(\mathbf{H})$. Let now $\langle\cdot,\cdot\rangle$ be the scaler product in $\mathbf{H}$. Then $\|\gamma + \alpha\|^2 = \|\gamma\|^2 + + 2\langle\gamma,\alpha\rangle + \|\alpha\|^2, \varepsilon \in (0,1]$ and hence Eq.(3.5) is equivalent with

$$\left(\widetilde{\int\limits_{\mathbf{H}}} \exp\left[\frac{i}{2}\|\gamma\|^2\right] \exp[i\langle\gamma,\alpha\rangle] \exp\left[\frac{i}{2}\|\alpha\|^2\right] D[\gamma]\right)_\varepsilon = (1)_\varepsilon$$

(3.6)

i.e. observing that $\exp\left[\frac{i}{2}\|\alpha\|^2\right]$ is independent of the variable $\gamma$ we get, with the obvious postulate that the normalized subintegral be $\widetilde{\mathbb{C}}$- linear:

$$\left(\int\limits_{\mathbf{H}}^{\widetilde{}}\exp\left[\frac{i}{2}\|\gamma\|^2\right]\exp[i\langle\gamma,\alpha_\varepsilon\rangle]D[\gamma]\right)_\varepsilon = \left(\exp\left[-\frac{i}{2}\|\alpha\|^2\right]\right)_\varepsilon. \qquad (3.7)$$

Set now

$$\left(\Phi_{\alpha_\varepsilon}(\gamma)\right)_\varepsilon = \left(\exp[i\langle\gamma,\alpha_\varepsilon\rangle]\right)_\varepsilon \qquad (3.8)$$

and

$$\left(\mathbf{J}[\Phi_{\alpha_\varepsilon}]\right)_\varepsilon = \left(\int\limits_{\mathbf{H}}^{\widetilde{}}\exp\left[\frac{i}{2}\|\gamma\|^2\right]\Phi_{\alpha_\varepsilon}(\gamma)D[\gamma]\right)_\varepsilon. \qquad (3.9)$$

Then from Eq.(3.7) - Eq.(3.9) we get

$$\left(\mathbf{J}[\Phi_{\alpha_\varepsilon}]\right)_\varepsilon = \left(\int\limits_{\mathbf{H}}^{\widetilde{}}\exp\left[\frac{i}{2}\|\gamma\|^2\right]\Phi_{\alpha_\varepsilon}(\gamma)D[\gamma]\right)_\varepsilon$$

$$= \left(\exp\left[-\frac{i}{2}\|\alpha_\varepsilon\|^2\right]\right)_\varepsilon. \qquad (3.10)$$

Eq.(3.10) gives us the evaluation of the functional $\left(\mathbf{J}[\Phi_{\alpha_\varepsilon}]\right)_\varepsilon$ for the value of the argument $\Phi_\varepsilon^\alpha$. Extending now naturally the definition of $\left(\mathbf{J}[\Phi_{\alpha_\varepsilon}]\right)_\varepsilon$ to linear combinations of functions (3.8), where $\alpha_\varepsilon$ varies over $E_M(\mathbf{H})$, and observing that when $\alpha_\varepsilon$ is the zero vector in $\mathbf{H}$ for all $\varepsilon \in (0,1]$ we have $\left(\Phi_{\alpha_\varepsilon}\right)_\varepsilon = (1)_\varepsilon$ and $\left(\mathbf{J}[1]\right)_\varepsilon = (1)_\varepsilon$, by Eq.(3.10) (which is in this case precisely the normalization condition (3.2) ),we see that $\left(\mathbf{J}[\Phi_{\alpha_\varepsilon}]\right)_\varepsilon$ becomes a normalized $\widetilde{\mathbb{C}}$-linear functional on the linear vector space over $\widetilde{\mathbb{C}}$ generated by the functions (3.8) (which is actually identical with the algebra generated by such functions). We have then

$$\left( \int\limits_{\mathbf{H}}^{\widetilde{\phantom{x}}} \exp\left[ \frac{i}{2} \|\gamma\|^2 \right] \sum_n \Phi_{\alpha_{n,\varepsilon}}(\gamma) D[\gamma] \right)_\varepsilon = \left( \mathbf{J}\left[ \sum_n \Phi_{\alpha_{n,\varepsilon}}(\gamma) \right] \right)_\varepsilon =$$

(3.11)

$$\sum_n c_n \left( \exp\left[ -\frac{i}{2} \|\alpha_{n,\varepsilon}\|^2 \right] \right)_\varepsilon = \left( \int\limits_{\mathbf{H}} \exp\left[ \frac{i}{2} \|\alpha_\varepsilon\|^2 \right] \sum_n \delta(\alpha_\varepsilon - \alpha_{n,\varepsilon}) D[\gamma] \right)_\varepsilon,$$

where $\delta(\alpha_\varepsilon - \alpha_{n,\varepsilon})$ is the measure concentrated at the origin in $\mathbf{H}$. Observing that

$$(\Phi_\varepsilon(\gamma))_\varepsilon = \sum_n (\Phi_{\alpha_{n,\varepsilon}}(\gamma))_\varepsilon$$

(3.11′)

is the Fourier-Colombeau transform on $\mathbf{H}$ of the Colombeau measure $(d\mu_{\Phi_\varepsilon}(\gamma))_\varepsilon$ i.e.

$$(\Phi_\varepsilon(\gamma))_\varepsilon = \left( \int\limits_{\mathbf{H}} \exp[i\langle \gamma, \alpha_\varepsilon \rangle] d\mu_{\Phi_\varepsilon}(\gamma) \right)_\varepsilon$$

(3.12)

we can rewrite Eq.(3.11) in the form

$$(\mathbf{J}[\Phi_\varepsilon(\gamma)])_\varepsilon =$$

$$\left( \int\limits_{\mathbf{H}}^{\widetilde{\phantom{x}}} \exp\left[ \frac{i}{2} \|\gamma\|^2 \right] \Phi_\varepsilon(\gamma) D[\gamma] \right)_\varepsilon = \left( \int\limits_{\mathbf{H}} \exp\left[ -\frac{i}{2} \|\alpha_\varepsilon\|^2 \right] d\mu_{\Phi_\varepsilon}(\gamma) \right)_\varepsilon.$$

(3.13)

We now observe that the right hand side of Eq.(3.13) is well defined also when $\mu_{\Phi_\varepsilon}$ is replaced by any bounded for all $\varepsilon \in (0,1]$ Colombeau complex measure on the measurable space $\mathbf{H}$, $\exp\left[ -\frac{i}{2} \|\alpha_\varepsilon\|^2 \right]$ being a bounded continuous function. Since any such Colombeau measure $(\mu_\varepsilon)_\varepsilon$ has a Colombeau-Fourier transform on $\mathbf{H}$, namely a bounded uniformly continuous function $(\Phi_{\mu_\varepsilon})_\varepsilon$ for all $\varepsilon \in (0,1]$ such that

$$(\Phi_{\mu_\varepsilon}(\gamma))_\varepsilon = \left( \int_{\mathbf{H}} \exp[i\langle \gamma, \alpha_\varepsilon \rangle] d\mu_\varepsilon(\alpha_\varepsilon) \right)_\varepsilon, \qquad (3.14)$$

it is natural to extend the definition of $(\mathbf{J}[\Phi_\varepsilon(\gamma)])_\varepsilon$ to all such functions $(\Phi_\varepsilon(\gamma))_\varepsilon$. We then have

$$(\mathbf{J}[\Phi_\varepsilon])_\varepsilon = \left( \int_{\mathbf{H}} \exp\left[-\frac{i}{2} \|\alpha_\varepsilon\|^2\right] d\mu_\varepsilon(\alpha_\varepsilon) \right)_\varepsilon \qquad (3.15)$$

for all functions $\Phi_\varepsilon(\gamma)$ on $\mathbf{H}$ which are Colombeau-Fourier transforms of bounded complex Colombeau measures $(\mu_\varepsilon)_\varepsilon$ i.e. such that

$$(\Phi_\varepsilon(\gamma))_\varepsilon = \left( \int_{\mathbf{H}} \exp[i\langle \gamma, \alpha_\varepsilon \rangle] d\mu_\varepsilon(\alpha_\varepsilon) \right)_\varepsilon. \qquad (3.16)$$

Because of (3.13) It is natural to use for $(\mathbf{J}[\Phi_\varepsilon])_\varepsilon$ also the notation

$$(\mathbf{J}[\Phi_\varepsilon])_\varepsilon = \left( \widetilde{\int_{\mathbf{H}}} \exp\left[-\frac{i}{2} \|\gamma\|^2\right] \Phi_\varepsilon(\gamma) D[\gamma] \right)_\varepsilon, \qquad (3.17)$$

and to take this, together with Eq.(3.15), as the definition of the normalized Integral (2.1). We have then that the normalized integral is a linear normalized $\widetilde{\mathbb{C}}$- valued functional defined at least for all nets $(\Phi_\varepsilon)_\varepsilon$ in the linear vector space $E_M[\mathscr{F}(\mathbf{H})]$ of Colombeau-Fourier transforms of bounded complex Colombeau measures on $\mathbf{H}$, where $\mathscr{F}(\mathbf{H})$ is the linear vector space of Fourier transforms of bounded complex measures on $\mathbf{H}$. The structure of $\mathscr{F}(\mathbf{H})$ is well studied. It is easy to see that $\mathscr{F}(\mathbf{H})$ the linear vector space is a Banach space equipped with the norm $\|\Phi_\varepsilon\|_{\mathscr{F}} = \|\mu_\varepsilon\|$ for all $\varepsilon \in (0,1]$, where $\|\mu_\varepsilon\|$ is the total variation norm for the corresponding measure $\mu_\varepsilon$, in fact $\mathscr{F}(\mathbf{H})$ is a Banach function algebra with respect to the pointwise multiplication of functions. $\mathscr{F}(\mathbf{H})$ is bijectively isometric with $M(\mathbf{H})$, where $M(\mathbf{H})$ is the space of bounded complex measures on $\mathbf{H}$ equipped with the total variation norm, the bisection being given by the Fourier transform

$$\Phi_\varepsilon(\gamma) = \int_{\mathbf{H}} \exp[i\langle\gamma,\alpha_\varepsilon\rangle]d\mu_\varepsilon(\alpha_\varepsilon).$$ An Intrinsic charac- terization of $\mathscr{F}(\mathbf{H})$ is e.g. as the complex linear hull of those positive definite functions $\Phi_{\varepsilon\in(0,1]}$ which are continuous for all $\varepsilon \in (0,1]$ in the so-called Minlos-Sazonov-Gross topology i.e. those continuous positive definite functions $\Phi_{\varepsilon\in(0,1]}$ on $\mathbf{H}$ for which for any $\delta > 0$ there exists a nuclear operator (positive (symnetric) trace class) $\mathbf{N}_\delta$ such that $\mathrm{Re}(\Phi_\varepsilon(0) - \Phi_\varepsilon(\gamma)) < \delta$ whenever $\langle\gamma,\mathbf{N}_\delta\gamma\rangle < 1$. We call $\mathscr{F}(\mathbf{H})$ the space of (Fresnel) inteqrable functions [38],[39] and we call $(\mathbf{J}[\Phi_\varepsilon])_\varepsilon$ the (Colombeau-Fresnel) integral of $(\Phi_\varepsilon)_\varepsilon$. $(\mathbf{J}[\Phi_\varepsilon])_\varepsilon$ is the normalized integral of $\left(\Phi_\varepsilon(\cdot)\exp\left[\frac{i}{2}\|\cdot\|^2\right]\right)_\varepsilon$ illustrated above. The properties of this Colombeau integral follows directly from the properties of the Fresnel integral which are studied in [39].

▲ **P**. **1**. $(\mathbf{J}[\Phi_\varepsilon])_\varepsilon$ is a $\widetilde{\mathbb{C}}$- linear normalized bounded functional on the Colombeau-Banach module $G_{\mathscr{F}(\mathbf{H})}$. Thus in particular we have

$$\left|\left(\mathbf{J}\left[\prod_{j=1}^n \Phi_{j,\varepsilon}\right]\right)_{\varepsilon\in(0,1]}\right| \leq \prod_{j=1}^n (\|\Phi_{j,\varepsilon}\|_{\mathscr{F}})_{\varepsilon\in(0,1]} \tag{3.18}$$

where $|\cdot|$ on the left hand side is the absolute value in the field $\mathbb{C}$ of complex numbers and we recall that $\|\cdot\|$ is the norm in $\mathscr{F}(\mathbf{H})$.

▲ **P.2**. When $\mathbf{H}$ is finite dimensional i.e. $\mathbf{H} = \mathbb{R}^d$ then

$$(\mathbf{J}[\Phi_\varepsilon])_\varepsilon = \left(\widetilde{\int_{\mathbf{H}}}\exp\left[-\frac{i}{2}\|\gamma\|^2\right]\Phi_\varepsilon(\gamma)D[\gamma]\right)_\varepsilon =$$

$$(2\pi i)^{-d/2}\left(\int_{\mathbb{R}^d}\exp\left[-\frac{i}{2}\|\gamma\|^2\right]\Phi_\varepsilon(\gamma)d\gamma\right)_\varepsilon \tag{3.19}$$

for all $(\Phi_\varepsilon)_\varepsilon$ in $E_M[\mathscr{F}(\mathbb{R}^d)]$ for which the Colombeau-Lebesgue integral exists. Another useful property is the linearity of $(\mathbf{J}[\Phi_\varepsilon])_\varepsilon$ under smooth partitions $\{\varphi_k\}$ of the unit i.e.

$$(\mathbf{J}[\Phi_\varepsilon])_\varepsilon = \left( \widetilde{\int}_{\mathbb{R}^d} \exp\left[ -\frac{i}{2} \|\gamma\|^2 \right] \Phi_\varepsilon(\gamma) d\gamma \right)_\varepsilon =$$

$$\sum_k \left( \widetilde{\int}_{\mathbb{R}^d} \exp\left[ -\frac{i}{2} \|\gamma\|^2 \right] \Phi_\varepsilon(\gamma) \varphi_k(\gamma) d\gamma \right)_\varepsilon = \qquad (3.20)$$

$$\sum_k (\mathbf{J}[\Phi_\varepsilon \varphi_k])_\varepsilon$$

with $\varphi_k(\gamma)$ of compact support (e.g. $\varphi_k(\gamma) = \psi_k(\gamma) - \psi_{k-1}(\gamma), k \in \mathbb{N}, \psi_0(\gamma) = 0$, $\psi_k(\gamma) = \psi(\gamma/k), \psi(\gamma)$ being a $C^\infty$ function equal to $0$ for $\|\gamma\| > 2$ and equal to $1$ for $\|\gamma\| < 1$.

▲ **P.3**. Let $\mathscr{F}_P(\mathbf{H})$ be the subalgebra of $\mathscr{F}(\mathbf{H})$ consisting of functions $f$ in $\mathscr{F}(\mathbf{H})$ with the property that $f(\gamma) = f(P_f\gamma)$ for some individual finite dimensional (i.e. with finite rank=$n$) projection on $\mathbf{H}$. i.e. the functions on $\mathscr{F}_P(\mathbf{H})$ are the cylinder or tame functions on $\mathbf{H}$. Let for $\left( \Phi_\varepsilon \right)_\varepsilon \in E_M(\mathscr{F}_P(\mathbf{H}))$ be $(\mathbf{J}_{P_{\Phi_\varepsilon}\mathbf{H}}[\cdot])_\varepsilon$ the normalized integral on the finite dimensional space $E_M(P_{\Phi_\varepsilon}\mathbf{H})$. Then $(\mathbf{J}_\mathbf{H}[\Phi_\varepsilon])_\varepsilon = \left( \mathbf{J}_{P_{\Phi_\varepsilon}\mathbf{H}}\left[ \widetilde{\Phi}_\varepsilon \right] \right)_\varepsilon$, where $\widetilde{\Phi}_\varepsilon, \varepsilon \in (0,1]$ is the function on $P_{\Phi_\varepsilon}\mathbf{H}$ such that $\left( \widetilde{\Phi}_\varepsilon(P_{\Phi_\varepsilon}\gamma) \right)_\varepsilon = ((\Phi_\varepsilon \circ P_{\Phi_\varepsilon})(\gamma))_\varepsilon$. In particular for $\left( \Phi_\varepsilon \right)_\varepsilon \in E_M(\mathscr{F}_P(\mathbf{H}))$

$$(\mathbf{J}_\mathbf{H}[\Phi_\varepsilon])_\varepsilon = (2\pi i)^{-(\dim P_{\Phi_\varepsilon}\mathbf{H})} \int_{\mathbb{R}^d} \widetilde{\Phi}_\varepsilon(\gamma) \exp\left[ \frac{i}{2} \|\gamma\|^2 \right] d\gamma. \qquad (3.21)$$

whenever the Riemann or Lebesgue integral on the right hand side exists.

▲ **P.4**. The translation invariance of the normalized integral has heen built in from the beginning in our presentation above. One has thus

$$\left( \widetilde{\int}_\mathbf{H} \Phi_\varepsilon(\gamma + \alpha_\varepsilon) \exp\left[ \frac{i}{2} \|\gamma + \alpha_\varepsilon\|^2 \right] D[\gamma] \right)_\varepsilon = \left( \widetilde{\int}_\mathbf{H} \Phi_\varepsilon(\gamma) \exp\left[ \frac{i}{2} \|\gamma\|^2 \right] D[\gamma] \right)_\varepsilon \qquad (3.22)$$

for all $(\alpha_\varepsilon)_\varepsilon \in E_M(\mathscr{F}_P(\mathbf{H}))$ Moreover for any bounded operator $\mathbf{A} : \mathbf{H} \to \mathbf{H}$ with everywhere define bounded inverse we have

$$\left( \int_{\mathbf{H}}^{\sim} \Phi_\varepsilon(\gamma) \exp\left[ \frac{i}{2} \|\gamma\|^2 \right] D[\gamma] \right)_\varepsilon = \left( \int_{\mathbf{H}}^{\sim} \Phi_\varepsilon(\mathbf{A}\gamma) \exp\left[ \frac{i}{2} \|\mathbf{A}\gamma\|^2 \right] D[\gamma] \right)_\varepsilon \qquad (3.23)$$

In particular one has the Euclidean invariance of the normalized integral.

▲ **P**. **5**. For the nornalized Colombeau integral a "Fubini Theorem" concerning iterated integration

$$\left( \overbrace{\int_{\mathbf{H}_1 \times \mathbf{H}_2}^{\sim} \Phi_\varepsilon(\gamma_1, \gamma_2) \exp\left[ \frac{i}{2} \|\gamma\|^2 \right] D[\gamma_1] D[\gamma_2]} \right)_\varepsilon =$$

$$\left( \int_{\mathbf{H}_2}^{\sim} \exp\left[ \frac{i}{2} \|\gamma_2\|^2 \right] \left\{ \int_{\mathbf{H}_1}^{\sim} \Phi_\varepsilon(\gamma_1, \gamma_2) \exp\left[ \frac{i}{2} \|\gamma_1\|^2 \right] D[\gamma_1] \right\} D[\gamma_2] \right)_\varepsilon = \qquad (3.24)$$

$$\left( \int_{\mathbf{H}_2}^{\sim} \exp\left[ \frac{i}{2} \|\gamma_1\|^2 \right] \left\{ \int_{\mathbf{H}_1}^{\sim} \Phi_\varepsilon(\gamma_1, \gamma_2) \exp\left[ \frac{i}{2} \|\gamma_2\|^2 \right] D[\gamma_2] \right\} D[\gamma_1] \right)_\varepsilon .$$

# The nofinalized Feynman-Colombeau path integral with respect to a quadratic form.

In order to treat on the same footing quantum mechanical systems in which a quadratic term in the coordinates arises in the action (and in the Hamiltonlan) it is useful to introduce the concept of an Feynman-Colombeau Integral normalized with respect to a (non necessarily positive) quadratic form.

Let us think as an example to a system in $\mathbb{R}^d$ with classical action $(S_T(\gamma))_\varepsilon, \varepsilon \in (0, 1]$ for the time Interval $[0, t]$

$$(S_{T,\varepsilon}(\gamma))_\varepsilon = (S^0_{T,\varepsilon}(\gamma))_\varepsilon + \left( \int_0^T V_\varepsilon(\gamma(t)) dt \right)_\varepsilon, \qquad (3.25)$$

where

$$(S^0_{T,\varepsilon}(\gamma))_\varepsilon = \frac{1}{2} \int_0^T \dot{\gamma}^2(t) dt - \frac{1}{2} \left( \int_0^T (\gamma(t), A_\varepsilon \gamma(t)) dt \right)_\varepsilon \qquad (3.26)$$

Here $A_\varepsilon$ is a $d \times d$ matrix . Since quadratic functions are not Fourier transform of bounded complex measures, it is natural to let the quadratic form $S^0_{T,\varepsilon}(\gamma)$ take the role of the kinetic energy term in the previous definition of normalized integral (Feynman path integrals for the case under consideration}. Mathematically it turns out that this can be done as a particular case of the following situation. Let **H** be a real separable Hilbert space. Let $(B_\varepsilon)_\varepsilon$ be an net a densely defined symmetric operators in **H**, then $\langle \gamma, B_\varepsilon \gamma \rangle_{\mathbf{H}}$ is a real quadratic form on $\mathbf{D}(B_\varepsilon)$. He would like to define, somewhat in analogy with above construction a normalized Colombeau integral

$$\left( \int_{\mathbf{H}}^{\sim} \Phi_\varepsilon(\gamma) \exp\left[ \frac{i}{2} \langle \gamma, B_\varepsilon \gamma \rangle_{\mathbf{H}} \right] D[\gamma] \right)_\varepsilon, \qquad (3.27)$$

where $\langle \cdot, \cdot \rangle_{\mathbf{H}}$ is the scalar product in **H**. The normalization should be such that

$$\mathbf{cl}\left[ \left( \int_{\mathbf{H}}^{\sim} \exp\left[ \frac{i}{2} \langle \gamma, B_\varepsilon \gamma \rangle_{\mathbf{H}} \right] D[\gamma] \right)_\varepsilon \right] = 1. \qquad (3.28)$$

From the care $\mathbf{H} = \mathbb{R}^d$ see that it is natural to assume that all $B_\varepsilon$, $\varepsilon \in (0, 1]$ it non degenerate, namely in order to define (3.27) In analogy with Eq.(3.13),i.e.by a Parseval relation involving the Colombeau-Fourier transfonws of $\exp[\frac{i}{2}\langle \gamma, B_\varepsilon \gamma \rangle]$ Let us
look first at the situation where $B_\varepsilon$ is strictly positive and bounded everywhere in **H**. Then, with $\langle \gamma, \gamma \rangle_{B_\varepsilon} = \langle \gamma, B_\varepsilon \gamma \rangle_{\mathbf{H}}$ the closure of **H** in the $\| \cdot \|_{B_\varepsilon}$-norm, we have by Eq.(3.13)

$$\left(\widetilde{\int_{\mathbf{H}_{B_\varepsilon}}} \Phi_\varepsilon(\gamma) \exp\left[\frac{i}{2}\langle\gamma,\gamma\rangle_{B_\varepsilon}\right] d\gamma\right)_\varepsilon = \left(\int_{\mathbf{H}_{B_\varepsilon}} \Phi_\varepsilon(\alpha) \exp\left[-\frac{i}{2}\langle\alpha,\alpha\rangle_{B_\varepsilon}\right] d\mu_{\Phi_\varepsilon}(\alpha,B_\varepsilon)\right)_\varepsilon \quad (3.29)$$

for any $\Phi_\varepsilon(\gamma)$ such that $\Phi_\varepsilon(\gamma) \in \mathscr{F}(\mathbf{H}_{B_\varepsilon}), \varepsilon \in (0,1]$ and where $\mu_{\Phi_\varepsilon}(\alpha,B_\varepsilon)$ is the Colombeau measure whose Colombeau-Fourier transfor is $(\Phi_\varepsilon(\gamma))_\varepsilon$ i.e.

$$(\Phi_\varepsilon(\gamma))_\varepsilon = \left(\int_{\mathbf{H}} \exp\left[i\langle\gamma,\alpha\rangle_{B_\varepsilon}\right] d\mu_{\Phi_\varepsilon}(\alpha,B_\varepsilon)\right)_\varepsilon =$$

$$\left(\int_{\mathbf{H}} \exp[i\langle\gamma,\alpha\rangle_{\mathbf{H}}] d\mu_{\Phi_\varepsilon}(\alpha)\right)_\varepsilon, \quad (3.30)$$

where

$$d\mu_{\Phi_\varepsilon}(\alpha) = d\mu_{\Phi_\varepsilon}(B_\varepsilon^{-1}\alpha, B_\varepsilon). \quad (3.31)$$

From (3.29) - (3.31) we get

$$\left( \widetilde{\int_{\mathbf{H}_{B_\varepsilon}}} \Phi_\varepsilon(\gamma) \exp\left[ \frac{i}{2} \langle \gamma, \gamma \rangle_{B_\varepsilon} \right] d\gamma \right)_\varepsilon =$$

$$\left( \int_{\mathbf{H}_{B_\varepsilon}} \Phi_\varepsilon(\alpha) \exp\left[ -\frac{i}{2} \langle \alpha, B_\varepsilon \alpha \rangle_{\mathbf{H}} \right] d\mu_{\Phi_\varepsilon}(\alpha, B_\varepsilon) \right)_\varepsilon =$$

(3.32)

$$\left( \int_{\mathbf{R}(B_\varepsilon)} \Phi_\varepsilon(\beta) \exp\left[ -\frac{i}{2} \langle B_\varepsilon^{-1}\beta, \beta \rangle_{\mathbf{H}} \right] d\mu_{\Phi_\varepsilon}(\beta) \right)_\varepsilon =$$

$$\left( \int_{\mathbf{R}(B_\varepsilon)} \Phi_\varepsilon(\beta) \exp\left[ -\frac{i}{2} \Delta_\varepsilon(\beta, \beta) \right] d\mu_{\Phi_\varepsilon}(\beta) \right)_\varepsilon$$

with $\Delta_\varepsilon(\beta, \beta)$ the bilinear form such that $\Delta_\varepsilon(\beta, \beta) = \langle B_\varepsilon^{-1}\beta, \beta \rangle_{\mathbf{H}}$ and $\mathbf{R}(B_\varepsilon)$ means the range of $B_\varepsilon$. Hence the normalized integral on $\mathbf{H}_{B_\varepsilon}$ is given by a "Parseval relation" involving an integral of $\exp[-\frac{i}{2}\Delta_\varepsilon(\beta, \beta)]$ against a bounded complex measure. This suggests to define in general (3.27) by

$$\left( \int_{\mathbf{D}_\varepsilon} \Phi_\varepsilon(\beta) \exp\left[ -\frac{i}{2} \Delta_\varepsilon(\beta, \beta) \right] d\mu_{\Phi_\varepsilon}(\beta) \right)_\varepsilon$$

(3.33)

where $\mathbf{D}_\varepsilon$ is a suitable linear subset of $\mathbf{H}$ and $\Delta_\varepsilon$ is a suitable "inverse" of $B_\varepsilon$. By the above considerations we see, that it is natural to require that $\Delta_\varepsilon$ be a symmetric bilinear form on $\mathbf{D}_\varepsilon \times \mathbf{D}_\varepsilon$, such that $\Delta_\varepsilon(\alpha, B_\varepsilon\beta) = \langle \alpha, \beta \rangle_{\mathbf{H}}$,
$\forall \alpha \in \mathbf{H}, \beta \in \mathbf{D}_\varepsilon = \mathbf{D}(B_\varepsilon)$, assuming that the range of $B_\varepsilon$ is contained in $\mathbf{D}_\varepsilon$. For applications to the case of QM measurments it is useful to allow $(\Delta_\varepsilon)_\varepsilon$ to be a $\widetilde{\mathbb{C}}$-valued and we see that the simplest assumption to define (3.33) is to require that the integration be taken on some nice space, e.g. $\mathbf{D}_\varepsilon$ a separable Banach space densely contained in $\mathbf{H}$ and with norm $\|\cdot\|_{\mathbf{D}_\varepsilon}$ stronger than the one $\|\cdot\|_{\mathbf{H}}$, in $\mathbf{H}$ (i.e. $\|\cdot\|_{\mathbf{H}} < \delta_\varepsilon \|\cdot\|_{\mathbf{D}_\varepsilon}$ for some $\delta_\varepsilon$), and with $\Delta_\varepsilon$ continuous on $\mathbf{D}_\varepsilon \times \mathbf{D}_\varepsilon$, with $\mathrm{Im}\,\Delta_\varepsilon(\gamma, \gamma) \leq 0$. The definition of (3.27) by (3.33) clearly depends on the choice of $\Delta_\varepsilon$, hence it is useful to write

$$(\mathbf{J}_{\Delta_\varepsilon}[\Phi_\varepsilon])_\varepsilon = \left( \int_{\mathbf{D}_\varepsilon} \Phi_\varepsilon(\gamma) \exp\left[ \frac{i}{2} \Delta_\varepsilon(\gamma, B_\varepsilon \gamma) \right] d\mu_{\Phi_\varepsilon}(\gamma) \right)_\varepsilon \qquad (3.34)$$

instead of (3.27). Hence we come to the definition:

**Definition.** Let $\left( \mathbf{H}, (\mathbf{D}_\varepsilon)_\varepsilon, (\|\cdot\|_\varepsilon)_\varepsilon, (B_\varepsilon)_\varepsilon, \Delta_\varepsilon, \right), \varepsilon \in (0,1]$, be a "Fresnel-Colombeau fourtuple" as described above, i.e. $\mathbf{H}$ a real separable Hilbert space, $\mathbf{D}_\varepsilon$ a real separable Banach space densely contained in $\mathbf{H}$ for all $\varepsilon \in (0,1]$ and with norm $\|\cdot\|_\varepsilon$ dominating the one in $\mathbf{H}, (B_\varepsilon)_{\varepsilon \in (0,1]}$ a net of a symmetric densely defined operators (not necessarily bounded) with range containing $\mathbf{D}_\varepsilon$ and $(\Delta_\varepsilon)_\varepsilon$ a $\widetilde{\mathbb{C}}$-valued symmetric bilinear form with non positive imaginary part and satisfying $\Delta_\varepsilon(\alpha, B_\varepsilon \beta) = \langle \alpha, \beta \rangle_{\mathbf{H}}$. Let $\mathcal{F}(\mathbf{D}_\varepsilon^*)$ be the Banach algebra of functions $\Phi_\varepsilon$ on the topological dual $D_\varepsilon^*$ of $\mathbf{D}_\varepsilon$ which are Fourier transforms of bounded complex measures $\mu_{\Phi_\varepsilon}$ on $\mathbf{D}_\varepsilon$ i.e. $\Phi_\varepsilon(\gamma) = \int \exp[i\langle \gamma, \alpha \rangle] d\mu_{\Phi_\varepsilon}(\alpha)$, where $\langle \cdot, \cdot \rangle$ is the dualization between $\mathbf{D}_\varepsilon$ and $\mathbf{D}_\varepsilon^*$. Multiplication of functions is, as in $\mathcal{F}(\mathbf{H})$, the pointwise one,and the norm $\|\cdot\|_{\mathcal{F}(\mathbf{D}_\varepsilon^*)}$ in $\mathcal{F}(\mathbf{D}_\varepsilon^*)$ is defined as the total variation $\|\mu_{\Phi_\varepsilon}\|$ in the corresponding Banach algebra $M(\mathbf{D}_\varepsilon^*)$ of bounded complex measures $\mu_{\Phi_\varepsilon}$ on $\mathbf{D}_\varepsilon^*$. We call $(\mathbf{J}_{\Delta_\varepsilon}[\Phi_\varepsilon])_\varepsilon$ as normalized Colombeau integral with respect to the quadratic form $(\Delta_\varepsilon)_\varepsilon$. This integral enjoys properties related to those of the normalized integral defined by (3.17).

▲ **P.1**. $\Phi_\varepsilon \to \mathbf{J}_{\Delta_\varepsilon}[\Phi_\varepsilon]$ for all $\varepsilon \in (0,1]$ is a $\mathbb{C}$-valued, linear, bounded, continuous normalized functional on the Banach function algebra $\mathcal{F}(\mathbf{D}_\varepsilon^*)$. One has

$$\left| \left( \mathbf{J}_{\Delta_\varepsilon} \left[ \prod_{j=1}^n \Phi_{j,\varepsilon} \right] \right)_{\varepsilon \in (0,1]} \right| \leq \prod_{j=1}^n \left( \|\Phi_{j,\varepsilon}\|_{\mathcal{F}(\mathbf{D}_\varepsilon^*)} \right)_{\varepsilon \in (0,1]} \qquad (3.35)$$

with $\|\cdot\|_{\mathcal{F}(\mathbf{D}_\varepsilon^*)} = \|\mu_{\Phi_\varepsilon}\|$ =total variation on $\mathbf{D}_\varepsilon$ of $\mu_{\Phi_\varepsilon}$.

▲ **P.2**. When $B_\varepsilon$ is strictly positive and with range dense in $\mathbf{D}_\varepsilon$ for all $\varepsilon \in (0,1]$ then

$$(\mathbf{J}_{\Delta_\varepsilon}[\Phi_\varepsilon])_\varepsilon = \left( \widetilde{\int_{\mathbf{H}_{B_\varepsilon}}} \Phi_\varepsilon(\gamma) \exp\left[ \frac{i}{2} \langle \gamma, B_\varepsilon \gamma \rangle_{B_\varepsilon} \right] d\mu_{\Phi_\varepsilon}(\gamma) \right)_\varepsilon \qquad (3.36)$$

In particular for $\mathbf{cl}[(B_\varepsilon)_\varepsilon] = 1$

$$(\mathbf{J}_{\Delta_\varepsilon}[\Phi_\varepsilon])_\varepsilon = (\mathbf{J}[\Phi_\varepsilon])_\varepsilon = \left( \widetilde{\int_{\mathbf{H}}} \Phi_\varepsilon(\gamma) \exp\left[ \frac{i}{2} \langle \gamma, \gamma \rangle_{\mathbf{H}} \right] d\mu_{\Phi_\varepsilon}(\gamma) \right)_\varepsilon. \qquad (3.37)$$

Thus the normalized Colombeau integral with respect to a quadratic form $(\Delta_\varepsilon)_\varepsilon$ reduces to the normalized integral when the quadratic form is the unit one on $\mathbf{H}$.

▲ **P.3**. If $B_\varepsilon^{-1}$ Is bounded and everywhere defined for all $\varepsilon \in (0,1]$ then $\mathbf{D}_\varepsilon = \mathbf{H}$ with equivalent norms and $\Delta_\varepsilon(\alpha, \beta) = \langle \alpha, B_\varepsilon^{-1}\beta \rangle_{\mathbf{H}}, \forall \alpha, \beta \in \mathbf{H}$. Moreover $B_\varepsilon$ Is self-adjoint, one has a direct sum decomposition $\mathbf{H} = \mathbf{H}_+^\varepsilon \bigoplus \mathbf{H}_-^\varepsilon$ of $\mathbf{H}$ and one has $B_\varepsilon = B_+^\varepsilon \bigoplus B_-^\varepsilon$ with $B_\pm^\varepsilon$ strictly positive, where $B_\pm^\varepsilon$ are the parts of $B_\varepsilon$ in $\mathbf{H}_\pm^\varepsilon$ resp. One has then the "Fubini theorem"

$$(\mathbf{J}_{\Delta_\varepsilon}[\Phi_\varepsilon])_\varepsilon = \left( \widetilde{\int_{\mathbf{H}}} \Phi_\varepsilon(\gamma) \exp\left[ \frac{i}{2} \langle \gamma, B_\varepsilon \gamma \rangle_{\mathbf{H}} \right] d\mu_{\Phi_\varepsilon}(\gamma) \right)_\varepsilon =$$

$$(3.38)$$

$$\left( \widetilde{\int_{\mathbf{H}_+^\varepsilon}} \exp\left[ \frac{i}{2} \| \gamma_1 \|_{\mathbf{H}_+^\varepsilon}^2 \right] \left\{ \widetilde{\int_{\mathbf{H}_-^\varepsilon}} \exp\left[ -\frac{i}{2} \| \gamma_2 \|_{\mathbf{H}_-^\varepsilon}^2 \right] \overline{\Phi_\varepsilon(\gamma_1 \bigoplus \gamma_2)} D[\gamma_2] \right\}^- D[\gamma_1] \right)_\varepsilon$$

where $\overline{\alpha}$ and $\{\cdot\}^-$ mean complex conjugates.

▲ **P.4**. If both $B_\varepsilon$ and $B_\varepsilon^{-1}$ are bounded everywhere on $\mathbf{H}$ for all $\varepsilon \in (0,1]$ then P.3 applies with $\mathbf{H} = \mathbf{D}_\varepsilon$, $\Delta_\varepsilon = B_\varepsilon$ and in this case

$$(\mathbf{J}_{\Delta_\varepsilon}[\Phi_\varepsilon])_\varepsilon = \left( \widetilde{\int_{\mathbf{H}_{B_\varepsilon}}} \Phi_\varepsilon(\gamma) \exp\left[ \frac{i}{2} \langle \gamma, \gamma \rangle_{B_\varepsilon} \right] D[\gamma] \right)_\varepsilon =$$

$$(3.39)$$

$$\left( \widetilde{\int_{\mathbf{H}}} \Phi_\varepsilon(\gamma) \exp\left[ -\frac{i}{2} \langle \alpha, B_\varepsilon \alpha \rangle_{\mathbf{H}} \right] d\mu_{\Phi_\varepsilon}(\alpha, B_\varepsilon) \right)_\varepsilon.$$

▲ **P.5**. If $\dim \mathbf{H} < \infty$ i.e. $\mathbf{H} = \mathbb{R}^d$ then

$$(\mathbf{J}_{\Delta_\varepsilon}[\Phi_\varepsilon])_\varepsilon = \left( \widetilde{\int_{\mathbf{H}_{B_\varepsilon}}} \Phi_\varepsilon(\gamma) \exp\left[ \frac{i}{2}\langle\gamma,\gamma\rangle_{B_\varepsilon} \right] D[\gamma] \right)_\varepsilon =$$

$$(2\pi)^{-d/2}|\det(B_\varepsilon)_\varepsilon| \left( \exp\left[ -\frac{i\pi}{4}\mathbf{sign}(B_\varepsilon) \right] \right)_\varepsilon \left( \int_{\mathbb{R}^d} \Phi_\varepsilon(\gamma) \exp\left[ \frac{i}{2}\langle\gamma,B_\varepsilon\gamma\rangle \right] d\gamma \right)_\varepsilon ,$$

(3.40)

where $\det B_\varepsilon$ and $\mathbf{sign}(B_\varepsilon)$ are the determinant resp. the signature of $B_\varepsilon$.

▲ **P.6**. Invarlance under translations by elements in $\mathbf{D}(B_\varepsilon)$ :

$$(\mathbf{J}_{\Delta_\varepsilon}[\Phi_\varepsilon])_\varepsilon = \left( \widetilde{\int_{\mathbf{H}}} \Phi_\varepsilon(\gamma) \exp\left[ \frac{i}{2}\langle\gamma,B_\varepsilon\gamma\rangle_{\mathbf{H}} \right] D[\gamma] \right)_\varepsilon =$$

$$\left( \widetilde{\int_{\mathbf{H}}} \Phi_\varepsilon(\gamma+\alpha) \exp\left[ \frac{i}{2}\langle\gamma+\alpha,B_\varepsilon(\gamma+\alpha)\rangle_{\mathbf{H}} \right] D[\gamma] \right)_\varepsilon .$$

(3.41)

# Apendix.IV. The Schrödinger equation with singular potential in $E_M[\mathcal{F}(\mathbb{R}^d)]$.

We consider the case of a Schrödinger operator of the form $-\Delta + (V_\varepsilon(x))_\varepsilon$, where $\Delta$ is the Laplacian on $\mathbb{R}^d$ and $(V_\varepsilon(x))_\varepsilon$ is a function in $E_M[\mathcal{F}(\mathbb{R}^d)]$. According to Feynman's heuristic considerations [13] we expect that the solution of the corresponding Colombeau-Schrödinger equation

$$i\frac{\partial}{\partial t}(\Psi_\varepsilon(t,x))_\varepsilon = \left( [-\Delta + V_\varepsilon(x)]\,\Psi_\varepsilon(t,x) \right)_\varepsilon$$

(4.1)

on the time interval $[0,T]$, with initial condition

$$(\Psi_\varepsilon(0,x))_\varepsilon = (\varphi_\varepsilon(x))_\varepsilon. \tag{4.2}$$

be given by a suitable normalized Colombeau path Integral

$$\widetilde{\int_{\gamma(T)=x}} \varphi(\gamma(0)) \exp\left[\frac{i}{2}(S_\varepsilon^T(\gamma))_\varepsilon\right] D[\gamma] \tag{4.3}$$

where $(S_\varepsilon^T(\gamma))_\varepsilon$ is the action correspondent to the given problem. In our problem the action is

$$(S_\varepsilon^T(\gamma))_\varepsilon = \int_0^T \dot{\gamma}^2(t)dt - \int_0^T (V_\varepsilon(\gamma(t)))_\varepsilon dt \tag{4.4}$$

and the integration should be in a suitable space of paths $\gamma$ ending at the instant $T$ in $x$. In order to be able to actually compute integrals of the form given via Eq.(4.3) It is natural to try to exploit the special feature of the phase function $(S_\varepsilon^T(\gamma))_\varepsilon$, namely of consisting of a positive quadratic part plus the term involving $(V_\varepsilon(\gamma(t)))_\varepsilon$, considered as a perturbation. The most simple way then of taking advantage of the quadratic part as a "known" part is to introduce as an integration space a linear space of paths $\gamma$ for which $\int_0^T \dot{\gamma}^2(t)dt$ is finite i.e. a subset of all absolutely continuous paths. For simplicity it is natural to restrict the paths to those building a Hilbert space **H** in which $\|\gamma\| = \left(\int_0^T \dot{\gamma}^2(t)dt\right)^{1/2}$ is the norm and for which Fourier transforms are particularly easily computed, which is the case when the Fourier transform of $\exp\left[\frac{i}{2}\int_0^T \dot{\gamma}^2(t)dt\right]$ evaluated as $\Gamma(s,\cdot)$ is equal to $\exp\left[\frac{i}{2}\Gamma(s,s)\right]$, where $\Gamma(s,t) = T\text{-max}(s,t)$ is the kernel of $-d^2/dt^2$ with zero boundary condition e.g. at $t = T$ and zero derivative condition at $t = 0$. This then selects **H** to be the Hilbert space of absolutely continuous paths on $[0,T]$ of finite kinetic energy and being at the origin at time $t$ (we note that the space has to be at least as big as to contain $\Gamma(s,\cdot)$ ). Then we have for paths $\gamma \in$ **H** that the expression (4.4) reads, for paths $\widetilde{\gamma} := \gamma + x$ (with the correct property for (4.4) of ending at time $T$ in $x$)

$$(S_\varepsilon^T(\widetilde{\gamma}))_\varepsilon = (S_\varepsilon^T(\gamma+x))_\varepsilon = \|\gamma\|^2 - \int_0^T (V_\varepsilon(\gamma(t)+x))_\varepsilon dt \tag{4.4$'$}$$

where $\|\cdot\|$ is the norm in $\mathbf{H}$. Thus the candidate for (4.3) becomes

$$\widetilde{\int_{\gamma(T)=x}} \varphi(\gamma(0)+x)\exp\left[\frac{i}{2}\|\gamma\|^2\right]\exp\left[-i\int_0^T (V_\varepsilon(\gamma(t)+x))_\varepsilon dt\right]D[\gamma]. \qquad (4.5)$$

**Theorem 4.1** Assume $\varphi_\varepsilon \in E_M[\mathcal{F}(\mathbb{R}^d)]$, $V_\varepsilon \in E_M[\mathcal{F}(\mathbb{R}^d)], \varepsilon \in (0,1]$ and define for $T \geq 0$

$$((\mathfrak{I}_T\varphi_\varepsilon)(x))_\varepsilon =$$

$$\widetilde{\int_{\mathbf{H}}} \varphi(\gamma(0)+x)\exp\left[\frac{i}{2}\|\gamma\|^2\right]\exp\left[-i\int_0^T (V_\varepsilon(\gamma(t)+x))_\varepsilon dt\right]D[\gamma]. \qquad (4.5')$$

Then for all $\varepsilon \in (0,1]$, $(\mathfrak{I}_T\varphi_\varepsilon)(x)$ is a continuous bounded semigroup on $E_M[\mathcal{F}(\mathbb{R}^d)]$ with

$$\|(\mathfrak{I}_T\varphi_\varepsilon)(x)\|_{\mathcal{F}} \leq \|\varphi_\varepsilon\|_{\mathcal{F}}\exp(t\|V_\varepsilon\|_{\mathcal{F}}), \qquad (4.5'')$$

where $\|\cdot\|_{\mathcal{F}}$ is the norm on $E_M[\mathcal{F}(\mathbb{R}^d)]$.
**Proof**. Let $\gamma \in \mathbf{H}$ and write

$$\gamma(t) = \gamma_1(t) + \gamma_2(t) \qquad (4.6)$$

where $\gamma_2(t)$ is constant on the subinterval $[0,t_1]$ $0 < t_1 < t$ and $\gamma_2(t) = \gamma(t)$ on $[t_1,T]$, and where $\gamma_1(t) = 0$ on $[t_1,T]$, $\gamma_1(t) = \gamma_2(t) - \gamma_2(T)$ on $[0,t_1]$. Clearly (4.6) gives a splitting of $\mathbf{H}$ in the direct sum $\mathbf{H} = \mathbf{H}_1 \bigoplus \mathbf{H}_2$ where $\mathbf{H}_1$ consists of the paths $\gamma_1(t)$ and $\mathbf{H}_2$ consists of the paths $\gamma_2(t)$. By the property P5) of the normalized integral (Apendix III "Fubini Theorem") we have

$$((\Im_T \varphi_\varepsilon)(x))_\varepsilon =$$

$$\left( \widetilde{\int_{\mathbf{H}}} \varphi(\gamma(0) + x) \exp\left[ \frac{i}{2} \int_0^T \dot{\gamma}^2(t) dt \right] \exp\left[ -i \int_0^T (V_\varepsilon(\gamma(t) + x))_\varepsilon dt \right] D[\gamma] \right)_\varepsilon =$$

$$\left( \widetilde{\int_{\mathbf{H}_1}} \exp\left[ i \int_{t_1}^T \dot{\gamma}_2^2(t) dt \right] \exp\left[ -i \int_{t_1}^T (V_\varepsilon(\gamma_2(t) + x))_\varepsilon dt \right] \times \right.$$

$$\left\{ \widetilde{\int_{\mathbf{H}_1}} \varphi(\gamma_1(0) + \gamma_2(t_1) + x) \exp\left[ i \int_0^{t_1} \dot{\gamma}_1^2(t) dt \right] \times \right.$$

$$\left. \exp\left[ -i \int_0^{t_1} (V_\varepsilon(\gamma_1(t) + \gamma_2(t) + x))_\varepsilon dt \right] D[\gamma_1] \right\} D[\gamma_2] \right)_\varepsilon =$$

$$((\Im_{T-t_1}(\Im_{t_1} \varphi_\varepsilon))(x))_\varepsilon,$$

which proves that $((\Im_T \varphi_\varepsilon)(x))_\varepsilon$ is a semigroup, The continuity of $(\Im_T \varphi_\varepsilon)(x), \varepsilon \in (0,1]$ and the bound on its norm follow from the xplicit computation of the normalized integral defining $(\Im_T \varphi_\varepsilon)(x)$ which is done by using the property PI (Appendix III) of the normalized integral i.e. that the Fresnel integral is a bounded linear functional on the Banach algebra $\mathcal{F}(\mathbf{H})$ so that

$$(\Im_T \varphi_\varepsilon)(x) =$$

$$\widetilde{\int_{\mathbf{H}}} \varphi_\varepsilon(\gamma(0) + x) \exp\left[ \frac{i}{2} \int_0^T \dot{\gamma}^2(t) dt \right] \exp\left[ -i \int_0^T (V_\varepsilon(\gamma(t) + x))_\varepsilon dt \right] D[\gamma] =$$

$$\sum_{n=0}^\infty \frac{(-i)^n}{n!} \underbrace{\int_0^T \dots \int_0^T}_{n} \widetilde{\int_{\mathbf{H}}} \varphi_\varepsilon(\gamma(0) + x) \left( \exp\left[ \frac{i}{2} \int_0^T \dot{\gamma}^2(t) dt \right] \prod_{j=1}^n V_\varepsilon(\gamma(t_j) + x) dt_j \right) D[\gamma] = \tag{4.8}$$

$$\sum_{n=0}^\infty \frac{(-i)^n}{n!} \underbrace{\int_0^T \dots \int_0^T}_{n} \underbrace{\int_{\mathbb{R}^d} \dots \int_{\mathbb{R}^d}}_{n} \exp\left[ i \sum_{j=0,k=0} \alpha_{j,\varepsilon} \Gamma(t_j, t_k) \alpha_{k,\varepsilon} \right] \exp\left[ i \sum_{j=0} \alpha_{j,\varepsilon} \right] \times$$

$$\prod_{j=1}^n (d\mu_V(\alpha_{j,\varepsilon}) dt_j) d\mu_{\varphi_\varepsilon}(\alpha_{0,\varepsilon}),$$

where $\mu_{V_\varepsilon}$ and $\mu_{\varphi_\varepsilon}$ are the measures on $\mathbb{R}^d$ of which $V_\varepsilon$ and $\varphi_\varepsilon$ are Fourier transforms.

In fact Eq.(4.8) yields the bound

$$\|(\Im_T \varphi_\varepsilon)(x)\|_{\mathscr{F}} \le \sum_{n=0}^\infty \frac{T^n}{n!} \|\mu_{V_\varepsilon}\| \|\mu_{\varphi_\varepsilon}\| = \|\varphi_\varepsilon\|_{\mathscr{F}} \exp[T\|V_\varepsilon\|_{\mathscr{F}}].$$

**Theorem 4.2**. The solution of the Schrödinger equation (4.1)-(4.2) with $V_\varepsilon$ and $\varphi_\varepsilon$, $\varepsilon \in (0,1]$ in $\mathscr{F}(\mathbf{H})$ is given by $(\Psi_\varepsilon(T,x))_\varepsilon = ((\Im_{T,\varepsilon}\varphi_\varepsilon)(x))_\varepsilon$ where $((\Im_T\varphi_\varepsilon)(x))_\varepsilon$ is given in Theorem 4.1 i.e. we have:

$$((\Im_{T,\varepsilon}\varphi_\varepsilon)(x))_\varepsilon =$$

$$\left( \widetilde{\int_{\mathbf{H}}} \varphi_\varepsilon(\gamma(0) + x) \exp\left[ \frac{i}{2} \int_0^T \dot{\gamma}^2(t) dt \right] \exp\left[ -i \int_0^T (V_\varepsilon(\gamma(t) + x))_\varepsilon dt \right] D[\gamma] \right)_\varepsilon. \tag{4.9}$$

In the case where $(V_\varepsilon(\gamma))_\varepsilon$ is real the semigroup $\Im_{T,\varepsilon}\varphi_\varepsilon, t > 0, \varepsilon \in (0,1]$ is a semigroup mapping isometrically and bijectively $\mathscr{F}(\mathbb{R}^d) \cap L_2(\mathbb{R}^d)$ into itself, hence extends naturally to an unitary group $U_{T,\varepsilon}$ on $L_2(\mathbb{R}^d)$. In fact $U_{T,\varepsilon} = \exp[-iTH_\varepsilon]$ where $H_\varepsilon$ is the self-adjoint operator $-\Delta + V_\varepsilon, \varepsilon \in (0,1]$ in $L_2(\mathbb{R}^d)$.

**Proof**: It suffices to compare (4.8) with the computation of the convergent Dyson series for $U_{T,\varepsilon} = \exp[-iTH_\varepsilon]$.

**Remark 4**.**1**.The expression for the solution of Schrödinger's equation (4.1)-(4.2) is the realization of the formal Feymman expression (4.3) where the Feynman path $\gamma$ ending at time $t$ in $x$ is written as $\gamma =: \gamma + x$ where $\gamma$ is a path in $\mathbf{H}$ i.e. ending at time $t$ in $0$. Actually the particular translation $\gamma \to \gamma + x$ could be replaced, by using the translation invariance of the normalized integral (property P4) of Apendix III ), by any translation $\gamma \to \gamma + x + \gamma_0$ with $\gamma_0 \in \mathbf{H}$. The formulae in Theorems 4.1,4.2 would then take a different aspect (but would of course be equi-valent to the previous ones). In this sense the formula (4.3) can be understood as a short notation for a family of equivalent formulae, e.g, for the one given in Theorem 4.2.

**Remark 4**.**2**.The expression (4.8) for the solution of the Schrödinger equation (4.1), (4.2) can be put directly in relation with the one obtained from Varadhan and Maslov-Chebotarev's [27] way of expressing the Feynman-Colombeau integrals as integrals with respect to complex Colombeau-Poisson measures.

Let now $(G_{T,\varepsilon}(x,y))_\varepsilon$ be the Green's function for (4.1), (4.2) such that the solution of Colombeau-Schrödlnger,s equation (4.1),(4,2) is given by

$$(\Psi_\varepsilon(T,x))_\varepsilon = \left( \int_{\mathbb{R}^d} \varphi_\varepsilon(y) G_{T,\varepsilon}(x,y) dy \right)_\varepsilon \qquad (4.9)$$

Using P.5 of Apendix III ("Fubini Theorem" for the normalized Colombeau integral) together with the splitting $\mathbf{H} = \mathbf{H}_0 \bigoplus \mathbf{H}_0^\perp$, with $\mathbf{H}_0 = \left\{ y\eta(t) | y \in \mathbb{R}^d, \eta(t) = 1 - \dfrac{t}{T} \right\}$ we get the following:

**Theorem 4**.**3**. The Green's function for Colombeau-Schrödinger equation (4.1),(4.2) with $\varphi_\varepsilon, V_\varepsilon \in \mathscr{F}(\mathbb{R}^d)$ for all $\varepsilon \in (0,1]$ is given by

$$(G_{T,\varepsilon}(x,y))_\varepsilon = G(T,x,y) \times$$

$$\left( \int_{\mathbf{H}}^{\widetilde{\quad}} \exp\left[ \frac{i}{2} \int_0^T \dot{\gamma}^2(t) dt \right] \exp\left[ -i \int_0^T (V_\varepsilon(\gamma(t) + (y-x)\eta(t) + x))_\varepsilon dt \right] D[\gamma] \right)_\varepsilon . \qquad (4.10)$$

where

$$G(T,x,y) = (2\pi i T)^{-d/2} \exp\left[ i \frac{|x-y|^2}{2T} \right] \qquad (4.11)$$

is the Gireen's function in the case $V_\varepsilon \equiv 0$.

Let us consider a Schrödinger-Colombeau operator of the form $-\Delta + (V_\varepsilon(x))_\varepsilon$, $x \in \mathbb{R}^d$ with $V_\varepsilon \in \mathcal{F}(\mathbb{R}^d)$ for all $\varepsilon \in (0,1]$ and real. Let $H_0 = \Delta$ and $H_\varepsilon = H_0 + V_\varepsilon(x)$, as self- adjoint operators on $L_2(\mathbb{R}^d)$ for all $\varepsilon \in (0,1]$. We consider now an operator $\exp(-itH_\varepsilon)\exp(itH_0)$ for $t \in \mathbb{R}$.

**Theorem 4.4.** Let $V_\varepsilon \in \mathcal{F}(\mathbb{R}^d)$ and $\varphi_\varepsilon \in \mathcal{F}(\mathbb{R}^d) \cap L_2(\mathbb{R}^d)$ for all $\varepsilon \in (0,1]$, then, for all $x \in \mathbb{R}^d$

$$([\exp(-itH_\varepsilon)\exp(itH_0)\varphi_\varepsilon](x))_\varepsilon = \left(\int_{\mathbb{R}^d} \exp(i\beta \cdot x) W_{t,\varepsilon}(x,\beta) d\mu_{\varphi_\varepsilon}(\beta)\right)_\varepsilon \quad (4.12)$$

where $\beta \cdot x$ is the scalar product in $\mathbb{R}^d$ and $\mu_{\varphi_\varepsilon}$ is the generalized Colombeau measure on $\mathbb{R}^d$ such that $(V_\varepsilon(x))_\varepsilon$ is the Colombeau-Fourier transform of $\mu_{\varphi_\varepsilon}$, and the quantity $(W_{t,\varepsilon}(x,\beta))_\varepsilon$ is given by the following normalized Colombeau paths integral:

$$(W_{t,\varepsilon}(x,\beta))_\varepsilon =$$

$$\left(\widetilde{\int}_{H_-} \exp\left[\frac{i}{2}\int_{-\infty}^{0}\dot{\gamma}^2(\tau)d\tau\right]\exp\left[-i\int_{-t}^{0}V_\varepsilon(\gamma(\tau) + \beta\tau + x)d\tau\right]D[\gamma]\right)_\varepsilon, \quad (4.13)$$

where $H_-$ is the Hilbert space of absolutely continuous paths $\gamma$ on the time interval $(-\infty, 0]$, with values in $\mathbb{R}^d$, of finite kinetic energy $\int_{-\infty}^{0}\dot{\gamma}^2(\tau)d\tau < \infty$

($\dot{\gamma}(\tau)$ is defined a.e., $\gamma$ being absolutely continuous) equal to the norm in $H_-$ and ending at time $0$ at the origin.
For potentials such that the wave operators $W_{\pm,\varepsilon}$ exist for all $\varepsilon \in (0,1]$ as strong limits in $L_2(\mathbb{R}^d)$ of $\exp(-itH_\varepsilon)\exp(itH_0)$ as $t \to \pm\infty$ we have:

$$(W_{+,\varepsilon}\varphi_\varepsilon)_\varepsilon = \left( \int_{\mathbb{R}^d} \exp(iv_- \cdot x) W_{+,\varepsilon}(x, v_-) d\mu_{\varphi_\varepsilon}(v_-) \right)_\varepsilon \qquad (4.14)$$

with

$$(W_{+,\varepsilon}(x, v_-))_\varepsilon =$$

$$\left( \widetilde{\int}_{H_-} \exp\left[ \frac{i}{2} \int_{-\infty}^0 \dot{\gamma}^2(\tau) d\tau \right] \exp\left[ -i \int_{-t}^0 V_\varepsilon(\gamma(\tau) + v_-\tau + x) d\tau \right] D[\gamma] \right)_\varepsilon, \qquad (4.15)$$

where the integral is an "improper Colombeau-Fresnel integral" defined as the limit for

$t \to \infty$ of the Colombeau-Fresnel integral of the Colombeau-Fresnel integrable function:

$$\exp\left[ -i \int_{-t}^0 V_\varepsilon(\gamma(\tau) + v_-\tau + x) d\tau \right].$$

Similarly

$$(W_{-,\varepsilon}\varphi_\varepsilon)_\varepsilon = \left( \int_{\mathbb{R}^d} \exp(-iv_+ \cdot x) W_{-,\varepsilon}(x, v_+) d\mu_{\varphi_\varepsilon}(v_+) \right)_\varepsilon \qquad (4.16)$$

with

$$(W_{-,\varepsilon}(x, v_-))_\varepsilon =$$

$$\left( \widetilde{\int}_{H_+} \exp\left[ \frac{i}{2} \int_{-\infty}^0 \dot{\gamma}^2(\tau) d\tau \right] \exp\left[ -i \int_{-t}^0 V_\varepsilon(\gamma(\tau) + v_+\tau + x) d\tau \right] D[\gamma] \right)_\varepsilon, \qquad (4.17)$$

here $H_+$ is defined as $H_-$ with the time interval $(-\infty, 0]$ replaced by $[0, +\infty)$.

**Remark 4.3.**In the case where $V_\varepsilon$ is gentle for all $\varepsilon \in (0,1]$ such that the perturbation series in powers of $V_\varepsilon$ for the wave operators converges for all $\varepsilon \in (0,1]$, the improper Colombeau-Fresnel integrals can also be defined by the termwise limit as $t \to \pm\infty$ of the convergent perturbation series obtained by expanding in powers of $V_\varepsilon$ and interchanging the sum and the Colombeau-Fresnel Integral.

**Remark 4.4.**: $(W_{+,\varepsilon}(x,v_-))_\varepsilon$ has e.g. the interpretation of the wave function at time $0$ of a particle which had prescribed velocity $v_-$ at time $-\infty$.

# Apendix V. Generalized Perturbation of Dissipative Operators.Colombeau-Kato generalized Potentials.

**Definition 5.1.**By a contraction semigroup on a Banach space $\Re$ we shall mean a family of bounded everywhere-defined linear operators $\mathbf{P}^t, 0 \le t < \infty$, mapping $\Re$ into itself such that

$$\mathbf{P}^0 = 1, \mathbf{P}^t\mathbf{P}^s = \mathbf{P}^{t+s}, 0 \le t,s < \infty,$$

$$\|\mathbf{P}^t\| \le 1, 0 \le t < \infty, \tag{5.1}$$

$$\lim_{t\to 0} \mathbf{P}^t\psi = \psi, \psi \in \Re.$$

**Definition 5.2.**The infinitesimal generator $\Im$ of a contraction semigroup $\mathbf{P}^t$ is defined by

$$\Im\psi = \lim_{t\to 0}\left[\frac{1}{t}(\mathbf{P}^t\psi - \psi)\right] \tag{5.2}$$

on the domain $\mathbf{D}(A)$ of all $\psi$ in $\Re$ for which the limit exists.

**Definition 5.3.**Let $\Re$ be an Banach space. By a contraction Colombeau semigroup on a Colombeau module $G_\Re$ we shall mean a family of bounded everywhere-defined generalized linear operators $(\mathbf{P}^t_\varepsilon)_{\varepsilon\in(0,1]}$, $0 \le t < \infty$, mapping Colombeau module $G_\Re$ into itself such that :

$$(\mathbf{P}^0_\varepsilon)_\varepsilon = 1,$$

$$(\mathbf{P}^t_\varepsilon)\mathbf{P}^s_\varepsilon = (\mathbf{P}^t_\varepsilon \mathbf{P}^s_\varepsilon)_\varepsilon = (\mathbf{P}^{t+s})_\varepsilon, 0 \le t, s < \infty,$$

$$(\|\mathbf{P}^t_\varepsilon\|)_\varepsilon \le 1, 0 \le t < \infty,$$

$$\left(\lim_{t\to 0} \mathbf{P}^t_\varepsilon \psi\right)_\varepsilon = \psi, \psi \in \Re,$$

$$\varepsilon \in (0, 1].$$

(5.3)

**Definition 5**.4.The infinitesimal generator $(\Im_\varepsilon)_\varepsilon$ of a contraction semigroup $(\mathbf{P}^t_\varepsilon)_\varepsilon$ is defined by

$$((\Im_\varepsilon)_\varepsilon)\psi = \left(\lim_{t\to 0}\left[\frac{1}{t}(\mathbf{P}^t_\varepsilon \psi - \psi)\right]\right)_\varepsilon$$

(5.4)

on the domain $\mathbf{D}(\Im)$ of all $\psi$ in $\Re$ for which the limit exists.

**Definition 5.5.**We shall call an operator $A$ on a Banach space $\Re$ dissipative in case for all $\psi$ in $\mathbf{D}(A)$, if $\varphi$ is in the space $\Re^*$ of continuous linear functionals on $\Re, \|\varphi\| = 1$ and $(\psi, \varphi) = \|\psi\|$ then $\text{Re}(A\psi, \varphi) \le 0$.

**Definition 5.6.**

**Theorem 5**.1.[41].Let $\Im$ be the infinitesimal generator of a contraction semigroup on the Banach space $\Re$ and let $B$ be a dissipative operator with $\mathbf{D}(B) \supset \mathbf{D}(\Im)$. If there exist constants $a$ and $b$ with $a < 1/2$ such that for all $\psi$ in $\mathbf{D}(\Im)$

$$\|B\psi\| \le a\|\Im\| + b\|\psi\|,$$

(5.5)

then $\Im + B$ is the infinitesimal generator of a contraction semigroup.

**Theorem 5**.2.[41].

$$\pi^{-1}I_v(z) = \int_0^\pi e^{z\cos\theta}\cos(v\theta)d\theta - \sin(v\pi)\int_0^\infty e^{-z(cht)-vt}dt$$

$$I_{\sqrt{2ME}}\left(\frac{M\omega u_a u_b}{i\sin(\omega S)}\right), z = \frac{M\omega u_a u_b}{i\sin(\omega S)}, v = \sqrt{2ME}$$

95–104,105–112,1992.

the Study of Language and Information, Stanford University.